\documentclass{memo-l}
\usepackage[utf8]{inputenc}
\usepackage{mathrsfs}
\usepackage{amssymb}
\usepackage{bm}
\usepackage{tikz-cd}
\usetikzlibrary{calc}
\usepackage{calc}
\usepackage{cite}
\usepackage{url}
\usepackage[symbol]{footmisc}
\usepackage{etoolbox}
\usepackage{hyperref}
\usepackage[refpage]{nomencl}
\usepackage{multicol}
\usepackage[greek,english]{babel}

\newif\ifidx
\idxfalse

\makenomenclature
\makeindex

\renewcommand{\nomname}{List of symbols}

\newtheorem{thm}{Theorem}[section]
\newtheorem{lemma}[thm]{Lemma}
\newtheorem{cor}[thm]{Corollary}
\newtheorem{prop}[thm]{Proposition}
\newtheorem*{thm*}{Theorem}
\newtheorem*{prop*}{Proposition}
\newtheorem{introthm}{Theorem}

\theoremstyle{remark}
\newtheorem{rk}[thm]{Remark}
\newtheorem{constr}[thm]{Construction}
\newtheorem{warn}[thm]{Warning}
\newtheorem{ex}[thm]{Example}

\newtheorem{convention}[thm]{Convention}
\newtheorem{notation}[thm]{Notation}

\newtheorem*{claim*}{Claim}
\AtBeginEnvironment{claim*}{\let\oldqed=\qedsymbol%
\renewcommand\qedsymbol{$\triangle$}}
\AtEndEnvironment{claim*}{\let\qedsymbol=\oldqed}

\theoremstyle{definition}
\newtheorem{defi}[thm]{Definition}

\numberwithin{equation}{section}

\newcommand{\Ho}{\textup{Ho}}
\newcommand{\nerve}{\textup{N}}
\newcommand{\h}{\textup{h}}
\newcommand{\cat}[1]{\textbf{\textup{#1}}}
\newcommand{\op}{{\textup{op}}}
\newcommand{\diag}{\mathop{\textup{diag}}\nolimits}
\newcommand{\blank}{{\textup{--}}}
\newcommand{\pr}{{\textup{pr}}}
\newcommand{\Hom}{{\textup{Hom}}}
\newcommand{\id}{{\textup{id}}}
\newcommand{\ev}{{\textup{ev}}}
\newcommand{\im}{\mathop{\textup{im}}}
\newcommand{\essim}{\mathop{\textup{ess im}}}
\newcommand{\incl}{{\textup{incl}}}
\newcommand{\colim}{\mathop{\textup{colim}}\nolimits}
\newcommand{\forget}{\mathop{\textup{forget}}\nolimits}
\newcommand{\Fun}{\textup{Fun}}
\newcommand{\FUN}{{\textup{Fun}_{\textup{enr}}}}
\newcommand{\Inj}{\textup{Inj}}
\newcommand{\supp}{\mathop{\textup{supp}}\nolimits}

\newcommand{\core}{\mathop{\textup{core}}}

\newcommand{\smashp}{\wedge}
\newcommand{\asm}{\mathord{\textup{asm}}}

\newcommand\und{\mathop{\textup{u}}\nolimits}

\newcommand{\Maps}{\mathord{\textup{maps}}}
\newcommand{\Ob}{\mathop{\textup{Ob}}}
\newcommand{\zigzag}{\mathord{\text{\raisebox{.8pt}{$\scriptstyle\bullet$}\kern-4pt$\to$\kern-1pt\raisebox{.8pt}{$\scriptstyle\bullet$}\kern-1pt$\gets$\kern-4pt\raisebox{.8pt}{$\scriptstyle\bullet$}}}}
\let\smashp=\wedge
\newcommand{\Sing}{\mathop{\textup{Sing}}}
\newcommand{\const}{\mathop{\textup{const}}}

\newcommand{\Ex}{\mathop{\textup{Ex}}\nolimits}
\newcommand{\hq}{\mathord{\textup{/\kern-2.5pt/}}}
\newcommand{\myh}{\mathord{\textup{`$\mskip-.1\thinmuskip h\mskip-.4\thinmuskip$'}}}
\newcommand{\discr}{\mathop{\textup{discr}}}
\newcommand{\Discr}{\mathop{\textup{Discr}}}
\newcommand{\End}{\mathord{\textup{End}}}
\newcommand{\Aut}{\mathord{\textup{Aut}}}
\newcommand{\centralizer}{\mathord{\textup{C}}}
\newcommand{\stabilizer}{\mathord{\textup{Stab}}}
\newcommand{\triv}{\mathord{\textup{triv}}}
\newcommand{\sh}{\mathord{\textup{sh}}}
\newcommand{\ppo}{\mathbin{\raise.25pt\hbox{\scalebox{.67}{$\square$}}}}
\newcommand{\symm}{{\bm{\Sigma}}}
\newcommand{\CMon}{{\textup{CMon}}}
\newcommand{\Mon}{{\textup{Mon}}}
\newcommand{\free}{{\textbf{\textup P}}}
\newcommand{\Sym}{{\textup{Sym}}}
\newcommand{\res}{{\textup{res}}}
\newcommand{\shear}{{\textup{shear}}}
\newcommand{\coshear}{{\textup{coshear}}}
\newcommand{\SP}{{\textup{SP}}}
\renewcommand{\Digamma}{{\textup{\foreignlanguage{greek}{Ϝ}}}}
\newcommand{\Gal}{{\textup{Gal}}}
\renewcommand{\digamma}{{\textit{\bfseries\foreignlanguage{greek}{ϝ}}}}
\newcommand{\qedNOW}{\pushQED{\qed}\qedhere\popQED}

\let\phi=\varphi
\let\setminus=\smallsetminus
\let\del=\partial

\def\twocell[#1]{\arrow[#1, dash, phantom, "\Rightarrow"{scale=1.125, yshift=-.4pt, description, allow upside down, sloped, inner sep=0pt}]}

\begin{document}
\title[$G$-Global Algebraic $K$-Theory]{$\bm G$-Global Homotopy Theory and\\ Algebraic $\bm K$-Theory}
\author{Tobias Lenz}
\address{Mathematisches Institut\newline\null\hskip\parindent Rheinische Friedrich-Wilhelms-Universit\"at Bonn\newline\null\hskip\parindent Endenicher Allee 60, 53115 Bonn, Germany \newline{}\newline\null\hskip\parindent Max-Planck-Institut für Mathematik\newline\null\hskip\parindent Vivatsgasse 7, 53111 Bonn, Germany\newline{}}
\email{lenz@math.uni-bonn.de}
\thanks{While writing this monograph, the author was an associate member of the Hausdorff Center for Mathematics, funded by the Deutsche Forschungsgemeinschaft (DFG, German Research Foundation) under Germany's Excellence Strategy (GZ 2047/1, project ID 390685813). He would moreover like to thank the Max Planck Institute for Mathematics in Bonn for their hospitality and support during that time.}
\subjclass[2020]{Primary 55P91, 
19D23, 
18F25, 
Secondary 55P48
}

\begin{abstract}
We develop the foundations of \emph{$G$-global homotopy theory} as a synthesis of classical equivariant homotopy theory on the one hand and global homotopy theory in the sense of Schwede on the other hand. Using this framework, we then introduce the \emph{$G$-global algebraic $K$-theory} of small symmetric monoidal categories with $G$-action, unifying $G$-equivariant algebraic $K$-theory, as considered for example by Shimakawa, and Schwede's global algebraic $K$-theory.

As an application of the theory, we prove that the $G$-global algebraic $K$-theory functor exhibits the category of small symmetric monoidal categories with $G$-action as a model of connective $G$-global stable homotopy theory, generalizing and strengthening a classical non-equivariant result due to Thomason. This in particular allows us to deduce the corresponding statements for global and equivariant algebraic $K$-theory.
\end{abstract}

\frontmatter
\maketitle
\tableofcontents

\mainmatter
\nomenclature[aSSet]{$\cat{SSet}$}{category of simplicial sets \nomnorefpage}
\nomenclature[aSet]{$\cat{Set}$}{category of sets \nomnorefpage}
\nomenclature[aAb]{$\cat{Ab}$}{category of abelian groups \nomnorefpage}
\nomenclature[aCat]{$\cat{Cat}$}{category of small categories \nomnorefpage}
\nomenclature[aN]{$\nerve$}{nerve of a small category \nomnorefpage}
\nomenclature[aeta]{$\eta$}{generic name for an adjunction unit \nomnorefpage}
\nomenclature[aepsilon]{$\epsilon$}{generic name for an adjunction counit \nomnorefpage}
\nomenclature[aphi]{$\phi,\psi,\dots$}{generic names for group homomorphisms from finite groups to $G$ (when discussing $G$-global weak equivalences) \nomnorefpage}
\nomenclature[aalpha]{$\alpha,\beta,\dots$}{generic names for group homorphisms (when discussing functoriality) \nomnorefpage}
\nomenclature[aG]{$G$}{generic discrete group \nomnorefpage}
\nomenclature[amaps]{$\Maps$}{simplicial mapping space (possibly based, possibly equipped with conjugation action) \nomnorefpage}
\nomenclature[auvw]{$u,v,\dots$}{generic names for elements of $\mathcal M$ \nomnorefpage}
\nomenclature[aDelta]{$\Delta$}{simplex category \nomnorefpage}
\nomenclature[aHo]{$\Ho$}{homotopy ($1$-)category of a category with weak equivalences \nomnorefpage}
\nomenclature[aCo]{$\mathscr C^\circ$}{subcategory of cofibrant-fibrant objects of a model category \nomnorefpage}
\nomenclature[aCf]{$\mathscr C^f$}{subcategory of fibrant objects of a model category \nomnorefpage}
\nomenclature[aCc]{$\mathscr C^c$}{subcategory of cofibrant objects of a model category \nomnorefpage}
\nomenclature[aFun1]{$\Fun$}{functor category \nomnorefpage}
\nomenclature[zo1]{$\otimes$}{generic notation for monoidal product \nomnorefpage}
\nomenclature[aTop]{$\cat{Top}$}{category of `topological spaces' (i.e.~compactly generated spaces) \nomnorefpage}
\nomenclature[aim]{$\im$}{image of a map \nomnorefpage}
\nomenclature[aess im]{$\essim$}{essential image of a functor \nomnorefpage}
\nomenclature[aGH]{$(G:H)$}{index of the subgroup $H$ in $G$ \nomnorefpage}
\nomenclature[aSigmanH]{$\Sigma_n\wr H$}{wreath product \nomnorefpage}

\chapter*{Introduction}
\begingroup\baselineskip=\the\baselineskip plus .05pt
\parskip=\the\parskip plus .5pt
\emph{Equivariant homotopy theory} is concerned with spaces carrying additional `symmetries,' encoded in the action of a suitable fixed group, and their (co)homology theories. One of the early successes of the equivariant point of view was the proof of the \emph{Atiyah-Segal Completion Theorem} \cite{atiyah-segal} using equivariant topological $K$-theory, generalizing and greatly simplifying Atiyah's original argument \cite{atiyah-old} based on non-equivariant $K$-theory. Subsequently, this motivated the \emph{Segal Conjecture} on the stable cohomotopy of classifying spaces, from which much of the original impetus for the development of equivariant stable homotopy theory derived \cite{carlsson-survey}, culminating in Carlsson's proof \cite{carlsson-segal}. At around the same time, the spectrum level foundations of equivariant stable homotopy theory were worked out by Lewis, May, and Steinberger \cite{lms}, and several additional equivalent models have been established by now \cite{mandell-may,hausmann-equivariant}.

Since these early days, equivariant homotopy theory has seen striking applications and connections to other areas of mathematics, ranging from the \emph{Sullivan Conjecture} \cite[Chapter~5]{sullivan} (proven by Carlsson \cite{carlsson-sullivan} and in an important special case by Miller \cite{miller}) on the relation between fixed points and homotopy fixed points, which was motivated by questions about the homotopy types of real algebraic varieties, to the celebrated solution of the Kervaire invariant one problem by Hill, Hopkins, and Ravenel \cite{hhr} in all dimensions apart from $126$.

\emph{Global homotopy theory}, as investigated among others by Schwede \cite{schwede-book}, provides a rigorous framework to talk about `uniform equivariant phenomena,' and in particular provides a natural home for many equivariant (co)homology theories that exist in a compatible way for all suitable groups, like equivariant topological $K$-theory and equivariant stable bordism. The global formalism has in several cases led to clean and conceptual descriptions of such uniform phenomena where direct descriptions for each individual group are much more opaque, for example for equivariant formal group laws \cite{hausmann-formal} or for the zeroth equivariant homotopy groups of symmetric products \cite{schwede-sym-prod}.

However, not all $G$-equivariant cohomology theories come from global ones; in particular, while there is a forgetful functor from the global stable homotopy category to the $G$-equivariant one admitting both adjoints, this is not a (Bousfield) localization unless $G$ is trivial, i.e.~global homotopy theory is not in any straight-forward way a generalization or refinement of $G$-equivariant homotopy theory.

The present monograph studies \emph{$G$-global homotopy theory} as a synthesis of the above two approaches; in particular, $1$-global homotopy theory recovers global homotopy theory, while there exists for every finite group $G$ a forgetful functor from (unstable or stable) $G$-global to $G$-equivariant homotopy theory admitting both a fully faithful left and a fully faithful right adjoint, yielding right and left Bousfield localizations, respectively. In fact, we develop our theory more generally for all discrete groups $G$, which for infinite $G$ refines \emph{proper equivariant homotopy theory}, as developed stably in \cite{proper-equivariant}. The idea of $G$-global homotopy theory has been around for some time---for example, specific point-set models of unstable and stable $\mathbb Z/2$-global homotopy theory appeared in preliminary versions of \cite{schwede-book}, while $\Sigma_n$-global weak equivalences for $n\ge0$ have recently been used in \cite{barrero} for the study of operads in unstable global homotopy theory---but this seems to be the first time the theory is developed systematically.

While there should also be a notion of $G$-global homotopy theory with respect to all compact Lie groups, we restrict our attention to the discrete case here. This in particular allows us to construct models of a more combinatorial nature, which is crucial for our main application:

\subsection*{Equivariant algebraic \texorpdfstring{$\bm K$}{K}-theory}
The \emph{algebraic $K$-theory} of small symmetric mo\-noidal categories was introduced by May \cite{may-permutative} using operadic techniques, and later an equivalent construction using Segal's theory of (special) $\Gamma$-spaces \cite{segal-gamma} was given by May \cite{may-unique} and Shimada and Shimakawa \cite{shimada-shimakawa}. Our study of $G$-global homotopy theory is motivated by refinements of this to equivariant and global contexts:

On the one hand, Shimakawa \cite{shimakawa} generalized Segal's theory to construct the \emph{$G$-equivariant algebraic $K$-theory} (a genuine $G$-spectrum) of a small symmetric monoidal category with the action of a finite group $G$, refining a construction of Fröhlich and Wall \cite{froehlich-wall} of the (low-dimensional) equivariant $K$-groups and extending work of Fiedorowicz, Hauschild, and May \cite{fhm}, who considered the equivariant algebraic $K$-theory of rings (with trivial $G$-action). More recently, equivariant algebraic $K$-theory has been revisited in the work of May, Merling, and their collaborators \cite{guillou-may, merling, may-merling-osorno}.

On the other hand, Schwede \cite{schwede-k-theory} introduced the \emph{global algebraic $K$-theory} of so-called \emph{parsummable categories}; he describes how small symmetric monoidal categories yield parsummable categories, and we have shown in \cite{sym-mon-global} that this accounts for all examples up to homotopy.

These two approaches generalize the original construction into two different directions, and in particular neither is a special case of the other---for example, Schwede's theory takes \emph{less} general inputs, but yields \emph{more} structure on the output. Accordingly, any direct comparison of the two constructions would be rather weak, only capturing a fraction of each of them.

Here we take a different route to clarifying the relation between the two approaches by realizing them as facets of a more general construction: using the framework of $G$-global homotopy theory developed in this monograph, we introduce \emph{$G$-global algebraic $K$-theory}, which for $G=1$ again recovers global algebraic $K$-theory, while for general $G$ it refines the equivariant construction. Together with a basic compatibility property of $G$-global algebraic $K$-theory under change of groups, this in particular shows that Schwede's global algebraic $K$-theory of a small symmetric monoidal category $\mathscr C$ forgets to the $G$-equivariant algebraic $K$-theory of $\mathscr C$ equipped with the trivial action, and we suggest to think of the existence of $G$-global algebraic $K$-theory as a general comparison of the above two approaches.

\subsection*{A \texorpdfstring{$\bm G$}{G}-global Thomason Theorem}
A classical result of Thomason says that all connective stable homotopy types arise as $K$-theory spectra of small symmetric monoidal categories; even stronger, he proved as \cite[Theorem~5.1]{thomason} that the $K$-theory functor expresses the homotopy category of connective spectra as a localization of the $1$-category of small symmetric monoidal categories.\index{Thomason Theorem}

Building on our results in \cite{sym-mon-global}, we prove the following $G$-global generalization of Thomason's theorem as the main result of this monograph:

\begin{introthm}\label{introthm:thomason-general}
For any discrete group $G$, the $G$-global algebraic $K$-theory construction exhibits the quasi-category of connective $G$-global stable homotopy types as a quasi-categorical localization of both the category $\cat{$\bm G$-SymMonCat}$ of small symmetric monoidal categories with $G$-action as well as the category $\cat{$\bm G$-ParSumCat}$ of parsummable categories with $G$-action.
\end{introthm}

This in particular immediately yields the corresponding statements for global and equivariant algebraic $K$-theory:

\begin{introthm}\label{introthm:thomason-global}
Schwede's global algebraic $K$-theory functor exhibits the quasi-category of connective global stable homotopy types as a quasi-categorical localization of both $\cat{SymMonCat}$ and $\cat{ParSumCat}$.
\end{introthm}

\begin{introthm}\label{introthm:thomason-equivariant}
For any finite $G$, Shimakawa's equivariant algebraic $K$-theory construction exhibits the quasi-category of connective $G$-equivariant stable homotopy types as a quasi-categorical localization of $\cat{$\bm G$-SymMonCat}$.
\end{introthm}

Theorem~\ref{introthm:thomason-global} (or more precisely its parsummable case) had been conjectured by Schwede in the original preprint version of \cite{schwede-k-theory}; to the best of our knowledge, also Theorem~\ref{introthm:thomason-equivariant} had not been proven before.

In fact, our methods yield a bit more: Mandell \cite{mandell} strengthened Thomason's result to an equivalence between symmetric monoidal categories (up to weak homotopy equivalences) and special $\Gamma$-spaces; in view of Segal's description of the passage from special $\Gamma$-spaces to connective spectra, Mandell called this a `non-group-completed' version of Thomason's theorem. Conversely, this result then gives a conceptual description of algebraic $K$-theory as a `higher group completion.'

Generalizing his result, we also prove `non-group-completed' versions of Theorems~\ref{introthm:thomason-general}--\ref{introthm:thomason-equivariant}, in particular showing that small symmetric monoidal categories model all \emph{ultra-commutative monoids} in the sense of \cite{schwede-book} (when viewed through the eyes of finite groups), and that small symmetric monoidal categories with $G$-action are equivalent to Shimakawa's special $\Gamma$-$G$-spaces. This way, we also get conceptual descriptions of equivariant, global, and $G$-global algebraic $K$-theory as `higher equivariant (global; $G$-global) group completions,' which we can in particular view as evidence that these are the `correct' generalizations of classical algebraic $K$-theory.

While the study of global versions of $\Gamma$-spaces, which is a key ingredient to our proof of Theorem~\ref{introthm:thomason-global}, naturally leads to $G$-global homotopy theory, it is interesting that ($G$-)global techniques are also central to our proof of the purely equivariant Theorem~\ref{introthm:thomason-equivariant}. In particular, the argument given here makes crucial use of parsummable categories (and their $G$-global generalizations), and we are not aware of a simple way to bypass them to obtain a standalone proof of Theorem~\ref{introthm:thomason-equivariant}.

\subsection*{Setup} While most of our arguments take place in Quillen's framework of model categories \cite{quillen}, we ultimately formulate our results in the language of quasi-categories in the sense of Joyal \cite{joyal-quasi-loc} and Lurie \cite{htt}. The reason for this is twofold: Firstly, our main results themselves are actually not model categorical at this point; for example, it is not clear whether $\cat{$\bm G$-ParSumCat}$ admits a model structure whose weak equivalences are created by the $G$-global algebraic $K$-theory functor. Secondly, we think that even in cases where there are model structures available, the higher categorical language is often cleaner: for example, our models will typically only be related by a zig-zag of Quillen equivalences, and similarly several of our constructions will become much easier to compare to each other once one has passed to quasi-categories.

\subsection*{Outline}
Chapter~1 introduces and compares various models of \emph{unstable $G$-global homotopy theory}, laying the foundations for the results established in later chapters. We also compare our approach to usual proper equivariant homotopy theory (Theorem~\ref{thm:G-global-vs-proper-sset}) and to Schwede's model of unstable global homotopy theory in terms of \emph{orthogonal spaces} (Theorem~\ref{thm:L-vs-I}).

We then study various notions of `$G$-globally coherently commutative monoids' in Chapter~2, in particular $G$-global versions of $\Gamma$-spaces or ultra-commutative monoids. As the main result of this chapter (Theorem~\ref{thm:gamma-vs-uc}) we prove that these models are equivalent, harmonizing Schwede's global approach with the classical equivariant theory.

Chapter~3 is concerned with \emph{stable $G$-global homotopy theory}. We introduce generalizations of Hausmann's global model structure \cite{hausmann-global} to symmetric spectra with $G$-action, whose weak equivalences for finite $G$ refine the usual $G$-equivariant stable weak equivalences \cite{hausmann-equivariant}. Using this, we then prove a $G$-global strengthening of a classical result due to Segal (Theorem~\ref{thm:group-completion}): any $G$-global $\Gamma$-space can be delooped to a connective $G$-global spectrum, and this provides an equivalence between so-called \emph{very special} $G$-global $\Gamma$-spaces on the one hand and connective $G$-global spectra on the other hand.

In Chapter~4 we introduce $G$-global algebraic $K$-theory and compare it to both equivariant and global algebraic $K$-theory. Finally, we use almost all of the theory developed in the previous chapters together with our results in \cite{sym-mon-global} to prove Theorems~\ref{introthm:thomason-general}--\ref{introthm:thomason-equivariant}.

\subsection*{Conventions} We assume Grothendieck's Axiom of Universes and fix once and for all universes $\cat{U}\in\cat{V}$; we will refer to elements of $\cat{U}$\nomenclature[aU]{$\cat{U}$}{implicitly chosen Grothendieck universe} as \emph{sets} and to elements of $\cat{V}$ as \emph{classes}. Accordingly, a \emph{small category} will be a $\cat{U}$-small category (i.e.~a category with sets of objects and morphisms), while we will refer to $\cat{V}$-small categories simply as `categories.'

By common abuse of language, we use the term `topological space' to mean \emph{compactly generated space} in the sense of \cite[§2]{mccord}.

\subsection*{Acknowledgements}
The results in this monograph were obtained as part of my PhD thesis at the University of Bonn. I would like to thank my advisor Stefan Schwede for suggesting global and equivariant algebraic $K$-theory as a thesis topic as well as for various instructive discussions about the results of this paper. I am moreover grateful to Peter May, Stefan Schwede, and the anonymous referee for helpful remarks on previous versions of this monograph.

This work began as a project to prove a global version of Thomason's theorem. I am indebted to Markus Hausmann, whose earlier suggestion that there should be a notion of $G$-global homotopy theory proved to be an illuminating perspective completely changing the form of this monograph and leading to the proof of the equivariant version of Thomason's theorem.
\endgroup

\chapter[Unstable $G$-global homotopy theory]{Unstable \texorpdfstring{\for{toc}{$G$}\except{toc}{$\bm G$}}{G}-global homotopy theory}\label{chapter:unstable}
In this chapter we will introduce several models of \emph{unstable $G$-global homotopy theory}, generalizing Schwede's unstable global homotopy theory \cite[Chapter~1]{schwede-book}. These models are already geared towards the study of $G$-global algebraic $K$-theory, and in particular, while we will be ultimately interested in \emph{stable} $G$-global homotopy theory and in the theory of $G$-global infinite loop spaces, the comparisons proven here will be instrumental in establishing results on the latter.

\section{Equivariant homotopy theory for monoids}\label{section:equiv-models}
\index{G-equivariant homotopy theory@$G$-equivariant homotopy theory}
Let $G$ be a discrete group. In unstable $G$-equivariant homotopy theory one is usually interested in $G$-spaces or $G$-simplicial sets up to so-called \emph{(genuine) $G$-weak equivalences},\index{G-weak-equivalence@$G$-weak equivalence|textbf}\index{genuine G-weak equivalence@genuine $G$-weak equivalence|seeonly{$G$-weak equivalence}} i.e.~$G$-equivariant maps that induce weak homotopy equivalences on $H$-fixed points for all subgroups $H\subset G$. More generally, one can consider the \emph{$\mathcal F$-weak equivalences}\index{F-weak equivalence@$\mathcal F$-weak equivalence|textbf} for any collection $\mathcal F$\nomenclature[aF]{$\mathcal F$}{generic collection of subgroups} of subgroups of $G$, i.e.~those maps that induce weak homotopy equivalences on $H$-fixed points for all $H\in\mathcal F$. If $\mathcal F=\mathcal A\ell\ell$\nomenclature[aAll]{$\mathcal A\ell\ell$}{collection of all subgroups} is the collection of all subgroups, this recovers the previous notion; at the other extreme, if $\mathcal F$ consists only of the trivial subgroup, then the $\mathcal F$-weak equivalences are precisely the underlying weak homotopy equivalences. We will at several points encounter \emph{proper equivariant homotopy theory}, where one considers the class $\mathcal Fin$\nomenclature[aFin]{$\mathcal Fin$}{collection of finite subgroups} of \emph{finite} subgroups of $G$; of course, $\mathcal Fin=\mathcal A\ell\ell$ if $G$ is finite, but for infinite $G$ these differ.
\index{G-equivariant homotopy theory@$G$-equivariant homotopy theory!proper|seeonly{proper $G$-equivariant homotopy theory}}
\index{proper G-equivariant homotopy theory@proper $G$-equivariant homotopy theory|textbf}

Our first approach to unstable $G$-global homotopy theory will rely on a generalization of this to actions of \emph{simplicial monoids}, and this section is devoted to generalizing several basic results from unstable equivariant homotopy theory to this context. In particular, we will construct equivariant model structures, prove a version of Elmendorf's Theorem, and discuss functoriality with respect to monoid homomorphisms.

\subsection{Equivariant model structures} There are several approaches to the construction of the usual equivariant model structures for group actions, for example the criteria of Dwyer and Kan \cite{dwyer-kan-equivariant} or Stephan \cite{cellular}. In this subsection, we will use the work of Dwyer and Kan to obtain the following analogue for the category $\cat{$\bm M$-SSet}$ of $M$-simplicial sets and $M$-equivariant maps for any simplicial monoid $M$:

\begin{defi}
Let $\mathcal F$ be a collection of finite subgroups of $M_0$. A map $f\colon X\to Y$ in $\cat{$\bm M$-SSet}$ is called an \emph{$\mathcal F$-weak equivalence}\index{F-weak equivalence@$\mathcal F$-weak equivalence|textbf} or \emph{$\mathcal F$-fibration} if $f^H$ is a weak homotopy equivalence or Kan fibration, respectively, for every $H\in\mathcal F$.
\end{defi}

We will also refer to the above also simply as \emph{$M$-weak equivalences}\index{M-weak equivalence@$M$-weak equivalence|seeonly {$\mathcal F$-weak equivalence}} and \emph{$M$-fibrations} (if $\mathcal F$ is clear from the context) or as \emph{equivariant weak equivalences}\index{equivariant weak equivalence|seeonly{$\mathcal F$-weak equivalence}} and \emph{equivariant fibrations} (if also $M$ is clear from the context).

\begin{prop}\label{prop:equiv-model-structure}
The $\mathcal F$-weak equivalences and $\mathcal F$-fibrations are part of a unique model structure on $\cat{$\bm M$-SSet}$, which we call the \emph{$\mathcal F$-model structure}\index{F-model structure@$\mathcal F$-model structure|textbf}\index{F-model structure@$\mathcal F$-model structure!injective|seeonly{injective $\mathcal F$-model structure}}. It is simplicial, combinatorial, proper, and filtered colimits in it are homotopical. A possible set of generating cofibrations is given by
\begin{equation*}
I=\{ M/H\times\del\Delta^n\hookrightarrow M/H\times\Delta^n : n\ge0,\text{ }H\in\mathcal F\},
\end{equation*}
and a possible set of generating acyclic cofibrations is given by
\begin{equation*}
J=\{ M/H\times\Lambda^n_k\hookrightarrow M/H\times\Delta^n : 0\le k\le n,\text{ }H\in\mathcal F\}.
\end{equation*}
\end{prop}

The finiteness condition on the subgroups $H$ is not necessary for the existence of the model structure, but it guarantees that filtered colimits are fully homotopical (i.e.~preserve arbitrary weak equivalences), which will simplify several arguments later. It is crucial for our argument that we only test weak equivalences and fibrations with respect to discrete sub\emph{groups}.

The proof of Proposition~\ref{prop:equiv-model-structure} will be given below. However, we already note:

\begin{lemma}\label{lemma:equivariant-weak-equivalences-prod-homotopical}
The weak equivalences of the above model structure are closed under finite products and small (i.e.~set-indexed) coproducts.
\begin{proof}
As fixed points commute with products and coproducts, this is immediate from the corresponding statement for ordinary simplicial sets.
\end{proof}
\end{lemma}

In order to conveniently formulate the result of Dwyer and Kan, we introduce the following notion:

\begin{defi}\label{defi:cellular-family}\index{cellular family|textbf}
Let $\mathscr D$ be a category enriched and tensored over $\cat{SSet}$. We call a family $(\Phi_i)_{i\in I}$ (for some set $I$) of enriched functors $\mathscr D\to\cat{SSet}$ \emph{cellular} if the following conditions are satisfied:
\begin{enumerate}
\item Each $\Phi_i$ preserves filtered colimits.
\item Each $\Phi_i$ is corepresentable in the enriched sense.
\item For each $i,j\in I$, $n\ge 0$, and some (hence any) $X_i$ corepresenting $\Phi_i$, the functor $\Phi_j$ sends pushouts along $\incl\otimes X_i\colon\del\Delta^n\otimes X_i\to\Delta^n\otimes X_i$ to homotopy pushouts in $\cat{SSet}$ (with respect to the weak homotopy equivalences).\label{item:cf-gluing}
\end{enumerate}
\end{defi}

By Proposition~\ref{prop:U-sharp-closure} (or an easy direct argument), $\Phi_j$ more generally sends pushouts along $f\otimes X_i\colon K\otimes X_i\to L\otimes X_i$ for any cofibration $f\colon K\to L$ of simplicial sets to homotopy pushouts. Thus, the above conditions immediately imply that the corepresenting elements $\{X_i:i\in I\}$ form a \emph{set of orbits} in the sense of \cite[2.1]{dwyer-kan-equivariant}. Our terminology is instead motivated by \cite[Proposition~2.6]{cellular}.

\begin{thm}[Dwyer \& Kan]\label{thm:dk-equivariant}
Let $\mathscr D$ be a complete and cocomplete category that is in addition enriched, tensored, and cotensored over $\cat{SSet}$. Let $(\Phi_i)_{i\in I}$ be a cellular family, and fix for each $i\in I$ an $X_i\in\mathscr D$ corepresenting $\Phi_i$.

Then there is a unique model structure on $\mathscr D$ in which a map $f$ is a weak equivalence or fibration if and only if $\Phi_i(f)$ is a weak equivalence or fibration, respectively, in $\cat{SSet}$ for each $i\in I$. This model structure is simplicial and moreover cofibrantly generated with generating cofibrations
\begin{equation*}
\{\incl\otimes X_i\colon \del\Delta^n\otimes X_i\to\Delta^n\otimes X_i : n\ge 0,i\in I\}
\end{equation*}
and generating acyclic cofibrations
\begin{equation*}
\{\incl\otimes X_i\colon \Lambda^n_k\otimes X_i\to \Delta^n\otimes X_i : 0\le k\le n,i\in I\}.
\end{equation*}
Finally, filtered colimits in $\mathscr D$ are homotopical.
\begin{proof}
\cite[Theorem~2.2]{dwyer-kan-equivariant} and its proof show that the model structure exists, that it is simplicial, and that it is cofibrantly generated by the above sets of (acyclic) cofibrations. The final statement follows immediately from the corresponding statement for $\cat{SSet}$ as each $\Phi_i$ preserves filtered colimits.
\end{proof}
\end{thm}

Next, we want to establish some additional properties of this model structure:

\begin{prop}\label{prop:dk-equivariant-addendum}
In the above situation, $\mathscr D$ is a proper model category. A commutative square in $\mathscr D$ is a homotopy pushout if and only if its image under $\Phi_i$ is a homotopy pushout in $\cat{SSet}$ for each $i\in I$.
\end{prop}

Here we use the term \emph{homotopy pushout}\index{homotopy pushout|textbf} for the dual of what Bousfield and Friedlander called a \emph{homotopy fibre square} in \cite[Appendix~A.2]{bousfield-friedlander}.

To prove the proposition, it will be convenient to realize the above model structure on $\mathscr D$ as a \emph{transferred model structure}:

\begin{defi}\index{transferred model structure|textbf}
Let $\mathscr D$ be a complete and cocomplete category, let $\mathscr C$ be a model category, and let $F\colon\mathscr C\rightleftarrows\mathscr D :\!U$ be an (ordinary) adjunction. The \emph{model structure transferred along $F\dashv U$} on $\mathscr D$ is the (unique if it exists) model structure where a morphism $f$ is a weak equivalence or fibration if and only if $Uf$ is.
\end{defi}

Transferred model structures will play a role at several points in this monograph and we recall some basic facts about them in Appendix~\ref{appendix:transfer}.

\begin{constr}\nomenclature[aPhi1]{$\Phi$}{right adjoint in Elmendorf's Theorem}
Fix for each $i\in I$ an $X_i\in\mathscr D$ corepresenting $\Phi_i$, and let $\cat{O}_{\Phi_\bullet}\subset\mathscr D$ be the full simplicial subcategory spanned by the $X_i$'s. We define the enriched functor $\Phi\colon\mathscr D\to\FUN(\textbf{O}_{\Phi_\bullet}^\op,\cat{SSet})$ as the composition
\begin{equation*}
\mathscr D\xrightarrow{\textup{enriched Yoneda}}\FUN(\mathscr D^\op,\cat{SSet})\xrightarrow{\textup{restriction}}\FUN(\textbf{O}_{\Phi_\bullet}^\op, \cat{SSet}),
\end{equation*}
where $\FUN$\nomenclature[aFun2]{$\FUN$}{enriched category of enriched functors} denotes the enriched category of simplicially enriched functors.

In other words, $\Phi(X)(Y)=\Maps_{\mathscr D}(Y,X)$ with the obvious functoriality. In particular, if $Y\in\textbf{O}_{\Phi_\bullet}$ corepresents $\Phi_i$, then we have an enriched isomorphism
\begin{equation}\label{eq:Phi-level}
\ev_Y\circ \Phi\cong \Phi_i.
\end{equation}
\end{constr}

\begin{constr}
It is well-known---see e.g.~\cite[Theorem~4.51]{enriched-presheaves} together with \cite[Theorem~3.73${}^\op$]{enriched-presheaves} for a statement in much greater generality---that for any essentially small simplicial category $T$, any cocomplete category $\mathscr D$ enriched and tensored over $\cat{SSet}$, and any simplicially enriched functor $F\colon T\to\mathscr D$, there exists an induced simplicial adjunction $\FUN(T^\op,\cat{SSet})\rightleftarrows\mathscr D$ with right adjoint $R$ given by $R(Y)(t)=\Maps(F(t),Y)$ for all $t\in T$, $Y\in\mathscr D$ with the obvious functoriality in each variable.

The left adjoint $L$ can be computed by the simplicially enriched coend
\begin{equation*}
L(X)=\int^{t\in T} X(t)\otimes F(t)
\end{equation*}
for any enriched presheaf $X$, together with the evident functoriality.
\end{constr}

In particular, applying this to the inclusion of $\textbf{O}_{\Phi_\bullet}$ yields:

\begin{cor}\nomenclature[aLambda]{$\Lambda$}{left adjoint in Elemendorf's Theorem}
The functor $\Phi$ has a simplicial left adjoint $\Lambda$.\qed
\end{cor}

By $(\ref{eq:Phi-level})$, the model structure from Theorem~\ref{thm:dk-equivariant} is transferred along $\Lambda\dashv\Phi$ from the projective model structure, allowing us to use the general results from Appendix~\ref{appendix:transfer}.

\begin{proof}[Proof of Proposition~\ref{prop:dk-equivariant-addendum}]
Right properness of the model structure on $\mathscr D$ is immediate from Lemma~\ref{lemma:transferred-properties}-$(\ref{item:tpr-proper})$.

For left properness and the characterization of homotopy pushouts, we observe that homotopy pushouts in $\FUN(\cat{O}_{\Phi_\bullet}^\op,\cat{SSet})$ can be checked levelwise, so that it is enough by Lemma~\ref{lemma:U-pushout-preserve-reflect} that $\Phi$ sends pushouts along cofibrations to homotopy pushouts. By Proposition~\ref{prop:U-sharp-closure} it suffices to check this for a set of generating cofibrations, which is then an instance of Condition~$(\ref{item:cf-gluing})$ of Definition~\ref{defi:cellular-family}.
\end{proof}

We now want to apply this to construct the $M$-equivariant model structure, for which we use the following well-known observation, cf.~e.g.~\cite[1.2]{dwyer-kan-equivariant} for the case of topological groups acting on spaces or \cite[Example~2.14]{cellular} for discrete groups acting on simplicial sets.

\begin{lemma}\label{lemma:fixed-points-cellular}
Let $M$ be a simplicial monoid and let $\mathcal F$ be a collection of finite subgroups of $M_0$. For any $H\in\mathcal F$, the enriched functor $(\blank)^H\colon\cat{$\bm M$-SSet}\to\cat{SSet}$ is corepresented by $M/H$ via evaluation at the class of $1\in M$. It preserves filtered colimits and pushouts along underlying cofibrations. In particular, the family $\big((\blank)^H\big)_{H\in\mathcal F}$ is cellular.
\begin{proof}
The corepresentability statement is obvious. As limits and colimits in $\cat{$\bm M$-SSet}$ are created in $\cat{SSet}$, filtered colimits commute with all finite limits, hence in particular with fixed points with respect to finite groups. Similarly, one reduces the statement about pushouts to the corresponding statement in $\cat{Set}$, which is easy.

As each of the maps $\del\Delta^n\times M/H\hookrightarrow\Delta^n\times M/H$ is in particular an underlying cofibration, the above immediately implies that the family of fixed point functors is cellular, finishing the proof.
\end{proof}
\end{lemma}

\begin{proof}[Proof of Proposition~\ref{prop:equiv-model-structure}]
By the previous lemma, we may apply Theorem~\ref{thm:dk-equivariant} and Proposition~\ref{prop:dk-equivariant-addendum}, so it only remains to show that this model category is combinatorial. But we know it is cofibrantly generated, and as an ordinary category $\cat{$\bm M$-SSet}$ is just the category of enriched functors of the category $BM$ with one object and endomorphism space $M$ into $\cat{SSet}$,\nomenclature[aBM]{$BM$}{(enriched) category with one object and hom set/mapping space $M$} hence locally presentable.
\end{proof}

\subsection{Elmendorf's Theorem}
The classical Elmendorf Theorem\index{Elmendorf's Theorem!classical|textbf} \cite{elmendorf} explains how $G$-equivariant homotopy theory (with respect to a fixed topological or simplicial group $G$) can be modelled in terms of fixed point data. Dwyer and Kan \cite[Theorem~3.1]{dwyer-kan-equivariant} provided a generalization of this to the above context:

\begin{thm}[Dwyer \& Kan]
If $(\Phi_i)_{i\in I}$ is any cellular family on $\mathscr D$, then the simplicial adjunction
\begin{equation*}
\Lambda\colon\FUN(\cat{O}_{\Phi_\bullet}^\op,\cat{SSet})\rightleftarrows\mathscr D :\!\Phi
\end{equation*}
is a Quillen equivalence for the projective model structure on the source.\qed
\end{thm}

If $M$ is a simplicial monoid and $\mathcal F$ is a collection of finite subgroups of $M_0$, then we write $\cat{O}_{\mathcal F}$ (or simply $\cat{O}_M$ if $\mathcal F$ is clear from the context)\nomenclature[aOFOM]{$\cat{O}_{\mathcal F}$ (also $\cat{O}_{M}$)}{(simplicially enriched) orbit category} for the full subcategory of $\cat{$\bm M$-SSet}$ spanned by the $M/H$ for $H\in\mathcal F$. The above then specializes to:

\begin{cor}\label{cor:elmendorf}\index{Elmendorf's Theorem!for monoids|textbf}
The simplicial adjunction
\begin{equation*}
\Lambda\colon\FUN(\cat{O}_{\mathcal F}^\op,\cat{SSet})\rightleftarrows\cat{$\bm M$-SSet}_{\textup{$\mathcal F$-equivariant}} :\!\Phi
\end{equation*}
is a Quillen equivalence for any simplicial monoid $M$ and any collection $\mathcal F$ of finite subgroups of $M_0$. Here $\Phi(X)(M/H)=\Maps_M(M/H,X)\cong X^H$ with the evident functoriality.\qed
\end{cor}

\subsection{Injective model structures}
For a group $G$, it is an easy observation that the cofibrations of the $\mathcal A\ell\ell$-model structure on the category $\cat{$\bm G$-SSet}$ of $G$-simplicial sets and $G$-equivariant maps are precisely the underlying cofibrations. In the case of a general collection $\mathcal F$, this will of course no longer be true---for example, if $\mathcal F$ is closed under subconjugates, then we can explicitly characterize the cofibrations as those injections such that all simplices not in the image have isotropy in $\mathcal F$, see e.g.~\cite[Proposition~2.16]{cellular}.

However, it is well-known that there is still a model structure with the same weak equivalences and whose cofibrations are the underlying cofibrations, called the \emph{mixed} or \emph{injective $\mathcal F$-model structure}, see e.g.~\cite[Proposition~1.3]{shipley-mixed} for a pointed version. We will now construct an analogue of this in our situation, which will use:

\begin{lemma}\label{lemma:homotopy-pushout-M-SSet}
Pushouts in $\cat{$\bm M$-SSet}$ along underlying cofibrations are homotopy pushouts (for any collection $\mathcal F$ of finite subgroups of $M_0$).
\begin{proof}
By Lemma~\ref{lemma:fixed-points-cellular}, each $(\blank)^H$ sends such a pushout to a pushout again. As taking fixed points moreover obviously preserves underlying cofibrations, this is then a homotopy pushout in $\cat{SSet}$, so the claim follows from Proposition~\ref{prop:dk-equivariant-addendum}.
\end{proof}
\end{lemma}

\begin{prop}\label{prop:equivariant-injective-model-structure}\index{injective F-model structure@injective $\mathcal F$-model structure|textbf}\index{equivariant injective model structure|seeonly {injective $\mathcal F$-model structure}}
Let $M$ be any simplicial monoid and let $\mathcal F$ be a collection of finite subgroups of $M_0$. Then there is a unique model structure on $\cat{$\bm{M}$-SSet}$ whose weak equvialences are the $\mathcal F$-weak equivalences and whose cofibrations are the injective cofibrations (i.e.~levelwise injections). We call this the \emph{injective $\mathcal F$-model structure} (or \emph{equivariant injective model structure} if $\mathcal F$ is clear from the context). It is combinatorial, simplicial, proper, and filtered colimits in it are homotopical.
\begin{proof}
As an ordinary category, $\cat{$\bm{M}$-SSet}$ is just a category of enriched functors into $\cat{SSet}$, and hence the usual injective model structure (which has weak equivalences the underlying non-equivariant weak homotopy equivalences) on it exists and is combinatorial. On the other hand, the $\mathcal F$-equivariant model structure is combinatorial, and its weak equivalences are stable under filtered colimits as well as pushouts along underlying cofibrations by the previous lemma.

Thus, we can apply Corollary~\ref{cor:mix-model-structures} to combine the cofibrations of the injective model structure with the $\mathcal F$-weak equivalences, yielding the desired model structure and proving that it is combinatorial, proper, and that filtered colimits in it are homotopical. It only remains to prove that it is simplicial, which means verifying the Pushout Product Axiom.\index{Pushout Product Axiom!for simplicial model categories} So let $i\colon K\to L$ be a cofibration of simplicial sets and let $f\colon X\to Y$ be an underlying cofibration of  $M$-simplicial sets. Because $\cat{SSet}$ is a simplicial model category, we immediately see that the pushout product map\nomenclature[zo]{$\ppo$}{pushout product}
\begin{equation*}
\begin{tikzcd}
K\times X\arrow[d, "i\times X"']\arrow[r, "K\times f"] & K\times Y\arrow[d]\arrow[rdd, "i\times Y", bend left=15pt]\\
L\times X\arrow[r]\arrow[rrd, "L\times f"', bend right=10pt] & (K\times Y)\amalg_{K\times X}(L\times X)\arrow[rd, "i\ppo f" description, dashed]\\
& & L\times Y
\end{tikzcd}
\end{equation*}
is again an underlying cofibration. It only remains to prove that this is a weak equivalence provided that either $i$ or $f$ is. For this we observe that the equivariant weak equivalences are stable under finite products by Lemma~\ref{lemma:equivariant-weak-equivalences-prod-homotopical}; moreover, a weak homotopy equivalence between simplicial sets with trivial $M$-action is already an equivariant weak equivalence. Hence, if $i$ is an acyclic cofibrations of simplicial sets, then the cofibration $i\times X$ is actually acyclic in the equivariant injective model structure, and so is $i\times Y$. Moreover, $K\times Y\to(K\times Y)\amalg_{K\times X}(L\times X)$ is also an acyclic cofibration as the pushout of an acyclic cofibration. It follows by $2$-out-of-$3$ that also $i\ppo f$ is an equivariant weak equivalence. The argument for the case that $f$ is an acyclic cofibration is analogous, and this finishes the proof.
\end{proof}
\end{prop}

\subsection{Functoriality} We will now explain how the above model structures for different monoids or collections of subgroups relate to each other.

\subsubsection{Change of monoid}\index{functoriality in homomorphisms!for M-SSet@for $\cat{$\bm{M}$-SSet}$|(}
If $\alpha\colon H\to G$ is any group homomorphism, then $\alpha^*$\nomenclature[aalphaaupperstar]{$\alpha^*$}{restriction along $\alpha$} obviously preserves cofibrations, fibrations, and weak equivalences of the $\mathcal A\ell\ell$-model structures. It follows immediately that the simplicial adjunctions $\alpha_!\dashv\alpha^*$\nomenclature[aalphalowershriek]{{$\alpha_{"!}$}}{left adjoint to restriction} and $\alpha^*\dashv\alpha_*$\nomenclature[aalphalowerstar]{{$\alpha_*$}}{right adjoint to restriction} are Quillen adjunctions. For monoids and general $\mathcal F$ one instead has to distinguish between the usual $\mathcal F$-model structure and the injective one:

\begin{lemma}\label{lemma:alpha-shriek-projective}
Let $\alpha\colon M\to N$ be any monoid homomorphism, let $\mathcal F$ be a collection of finite subgroups of $M_0$, and let $\mathcal F'$ be a collection of finite subgroups of $N_0$ such that $\alpha(H)\in\mathcal F'$ for all $H\in\mathcal F$. Then $\alpha^*$ sends $\mathcal F'$-weak equivalences to $\mathcal F$-weak equivalences and it is part of a simplicial Quillen adjunction
\begin{equation*}
\alpha_!\colon\cat{$\bm M$-SSet}_{\textup{$\mathcal F$-equivariant}}\rightleftarrows\cat{$\bm N$-SSet}_{\textup{$\mathcal F'$-equivariant}} :\!\alpha^*.
\end{equation*}
\begin{proof}
If $f$ is any morphism in $\cat{$\bm N$-SSet}$ and $H\subset M_0$ is any subgroup, then $(\alpha^*f)^H=f^{\alpha(H)}$. Thus, the claim follows immediately from the definition of the weak equivalences and fibrations of the equivariant model structures.
\end{proof}
\end{lemma}

\begin{lemma}\label{lemma:alpha-star-injective}
In the situation of Lemma~\ref{lemma:alpha-shriek-projective}, also
\begin{equation*}
\alpha^*\colon\cat{$\bm N$-SSet}_{\textup{$\mathcal F'$-equivariant injective}}\rightleftarrows\cat{$\bm M$-SSet}_{\textup{$\mathcal F$-equivariant injective}} :\!\alpha_*.
\end{equation*}
is a simplicial Quillen adjunction.
\begin{proof}
We have seen in the previous lemma that $\alpha^*$ is homotopical. Moreover, it obviously preserves injective cofibrations.
\end{proof}
\end{lemma}

The questions when $\alpha_!$ is left Quillen for the injective model structures or when $\alpha_*$ is right Quillen for the usual model structures are more complicated. The following propositions will cover the cases of interest to us:

\begin{prop}\label{prop:alpha-shriek-injective}
Let $\alpha\colon H\to G$ be an \emph{injective} homomorphism of discrete groups and let $M$ be any simplicial monoid. Let $\mathcal F$ be any collection of finite subgroups of $M_0\times H$ and let $\mathcal F'$ be a collection of finite subgroups of $M_0\times G$ such that the following holds: for any $K\in\mathcal F'$, $g\in G$ also $(M\times\alpha)^{-1}(gKg^{-1})\in\mathcal F$. Then
\begin{equation*}
\alpha_!\colon\cat{$\bm{(M\times H)}$-SSet}_{\textup{$\mathcal F$-equiv.~inj.}}\rightleftarrows\cat{$\bm{(M\times G)}$-SSet}_{\textup{$\mathcal F'$-equiv.~inj.}} :\!\alpha^*=(M\times\alpha)^*
\end{equation*}
is a simplicial Quillen adjunction; in particular, $\alpha_!$ is homotopical.
\begin{proof}
While we have formulated the result above in the way we later want to apply it, it will be more convenient for the proof to switch the order in which we write the actions, i.e.~to work with $(H\times M)$- and $(G\times M)$-simplicial sets.

We may assume without loss of generality that $H$ is a subgroup of $G$ and $\alpha$ is its inclusion. Then $G\times_H\blank$ (with the $M$-action pulled through via enriched functoriality) is a model for $\alpha_!$; if $X$ is any $H$-simplicial set, then the $n$-simplices of $G\times_HX$ are of the form $[g,x]$ with $g\in G$ and $x\in X_n$ where $[g,x]=[g',x']$ if and only if there exists an $h\in H$ with $g'=gh$ and $x=h.x'$.

Now let $K\subset G\times M_0$ be any subgroup. We set $S\mathrel{:=}\{g\in G: g^{-1}Kg\subset H\times M_0\}$ and observe that this is a right $H$-subset of $G$. We fix representatives $(s_i)_{i\in I}$ of the orbits. If now $X$ is any $(H\times M)$-simplicial set, then we define
\begin{equation*}
\iota\colon\coprod_{i\in I} X^{s_i^{-1}Ks_i} \to G\times_HX
\end{equation*}
as the map that is given on the $i$-th summand by $x\mapsto[s_i,x]$.

The following splitting ought to be well-known:

\begin{claim*}
The map $\iota$ is natural in $X$ (with respect to the evident functoriality on the left hand side) and it defines an isomorphism onto $(G\times_HX)^K$.

\medskip
Before we prove this, let us do a quick reality check: if $M=1$ is the trivial monoid, then $S$ consists precisely of the $g\in G$ with $g^{-1}Kg\subset H$, i.e.~those $g\in G$ for which right multiplication by $G$ defines a map $G/K\to G/H$. Two such elements of $S$ belong to the same $H$-orbit if and only if they correspond to the same map $G/K\to G/H$, i.e.~we have an isomorphism $\Hom_{\cat{$\bm G$-Set}}(G/K,G/H)\to S/H$ via evaluating at the coset $[1]$. For $X=*$ this then implies that $\iota$ can be identified with the map $\Hom_{\cat{$\bm G$-Set}}(G/K,G/H)\to G\times_H*\cong G/H, f\mapsto f[1]$, which is clearly a bijection onto the $K$-fixed points.

\begin{proof}
The naturality part is obvious. Moreover, it is clear from the choice of the $s_i$ as well as the above description of the equivalence relation that $\iota$ is injective, so that it only remains to prove that its image equals $(G\times_H X)^K$.

Indeed, assume $[g,x]$ is a $K$-fixed $n$-simplex. Then in particular $k_1g\in gH$ for any $k=(k_1,k_2)\in K$ by the above description of the equivalence relation, hence $g^{-1}k_1g\in H$ which is equivalent to $g^{-1}kg\in H\times M_0$. Letting $k$ vary, we conclude that $g\in S$, and after changing the representative if necessary we may assume that $g=s_i$ for some $i\in I$. But then
\begin{align*}
[s_i,x]&=k.[s_i,x]=[k_1s_i, k_2.x]=[s_i (s_i^{-1}k_1s_i), k_2.x]\\
&=[s_i, (s_i^{-1}k_1s_i, k_2).x]=[s_i,(s_i^{-1}ks_i).x]
\end{align*}
for any $k\in K$, and hence $x\in X^{s_i^{-1}Ks_i}$ as $H$ acts faithfully on $G$. Thus, $\im\iota$ contains all $K$-fixed points. Conversely, going through the above equation backwards shows that $[s_i,x]$ is $K$-fixed for any $(s_i^{-1}Ks_i)$-fixed $x$, i.e.~also $\im\iota\subset(G\times_HX)^K$.
\end{proof}
\end{claim*}
In particular, for $K=1$ this recovers the fact that non-equivariantly $G\times_HX$ is given as disjoint union of copies of $X$; we immediately conclude that $\alpha_!$ preserves injective cofibrations. On the other hand, if $K\in\mathcal F'$, then we conclude from the claim that for any morphism $f$ in $\cat{$\bm{(H\times M)}$-SSet}$ the map $(G\times_Hf)^K$ is conjugate to $\coprod_{i\in I} f^{s_i^{-1}Ks_i}$ for some $s_i\in G$ with $s_i^{-1}Ks_i\subset H\times M_0$ for all $i\in I$. Then by assumptions on $\mathcal F$ already $s_i^{-1}Ks_i\in\mathcal F$, so that each $f^{s_i^{-1}Ks_i}$ is a weak homotopy equivalence whenever $f$ is a $\mathcal F$-weak equivalence. As coproducts of simplicial sets are fully homotopical, we conclude that $(G\times_Hf)^K$ is a weak homotopy equivalence, and letting $K$ vary this shows $G\times_Hf$ is an $\mathcal F'$-weak equivalence.
\end{proof}
\end{prop}

\begin{prop}\label{prop:alpha-lower-star-homotopical}
In the situation of Proposition~\ref{prop:alpha-shriek-injective}, also
\begin{equation*}
\alpha^*\colon\cat{$\bm{(M\times G)}$-SSet}_{\textup{$\mathcal F'$-equivariant}}\rightleftarrows\cat{$\bm{(M\times H)}$-SSet}_{\textup{$\mathcal F$-equivariant}} :\!\alpha_*.
\end{equation*}
is a simplicial Quillen adjunction. Moreover, if the index $(G:\im\alpha)$ is finite, then $\alpha_*$ is fully homotopical.
\begin{proof}
We may again assume that $\alpha$ is the inclusion of a subgroup, so that $\alpha_*$ can be modelled as usual by $\Maps^H(G,\blank)$.

Let $K\subset M_0\times G$ be any subgroup, and let $K_2$ be its projection to $G$. We pick a system of representatives $(g_i)_{i\in I}$ of $H\backslash G/K_2$, and we let $L_i=(M\times H)\cap (g_iKg^{-1}_i)$. Similarly to the previous proposition, one checks that we have an isomorphism
$\coprod_{i\in I} (M\times H)/L_i \to (M\times G)/K$ given on summand $i$ by $[m,h]\mapsto[m,hg_i]$. Together with the canonical isomorphism $\Maps^H(G,\blank)^K\cong\Maps^{M\times H}((M\times G)/K,\blank)$ induced by the projection, this shows that for any $(M\times H)$-equivariant map $f\colon X\to Y$ the map $\alpha_*(f)^K$ is conjugate to $\prod_{i\in I}f^{L_i}$. If $K\in\mathcal F'$, then the assumptions guarantee that $L_i\in\mathcal F$, so $\alpha_*$ is obviously right Quillen. If in addition $(G:H)<\infty$, then $H\backslash G$ is finite, and hence so is $I$. As finite products in $\cat{SSet}$ are homotopical, so is $\alpha_*$ in this case.
\end{proof}
\end{prop}

\begin{rk}\index{graph subgroup|textbf}
Let $A,B$ be groups. We recall that a \emph{graph subgroup} $C\subset A\times B$ is a subgroup of the form $\{(a,\phi(a)) : a\in A'\}$ for some subgroup $A'\subset A$ and some group homomorphism $\phi\colon A'\to B$; note that this is not symmetric in $A$ and $B$. Both $A'$ and $\phi$ are uniquely determined by $C$, and we write $C\mathrel{=:}\Gamma_{A',\phi}$.\nomenclature[aGammaHphi]{$\Gamma_{H,\phi}$}{graph subgroup $\{(h,\phi(h)):h\in H\}$} A subgroup $C\subset A\times B$ is a graph subgroup if and only if $C\cap (1\times B)=1$.

If $A$ and $B$ are monoids, then we can define its graph subgroups as the graph subgroups of the maximal subgroup $\core(A\times B)$\nomenclature[acore]{$\core$}{maximal subgroup of a monoid} of $A\times B$. If $A'\subset\core(A)$ and $\phi\colon A'\to B$ is a homomorphism, then we will abbreviate $(\blank)^\phi\mathrel{:=}(\blank)^{\Gamma_{A',\phi}}$.\nomenclature[aphi]{$(\blank)^\phi$}{fixed points for the graph subgroup corresponding to $\phi$}
\end{rk}

\begin{ex}
Let $\mathcal E$ be any collection of finite subgroups of $M_0$ closed under taking subconjugates. Then the assumptions of the previous two propositions are in particular satisfied if we take $\mathcal F=\mathcal G_{\mathcal E,H}$\nomenclature[aG1]{$\mathcal G_{A,B}$}{graph subgroups of $A\times B$ for homomorphisms $A\to B$}\nomenclature[aG2]{$\mathcal G_{\mathcal E,B}$}{graph subgroups $\Gamma_{H,\phi}\subset A\times B$ with $H\in\mathcal E$} to be the collection of those graph subgroups $\Gamma_{K,\phi}$ of $M_0\times H$ with $K\in\mathcal E$, and similarly $\mathcal F'=\mathcal G_{\mathcal E,G}$.
\end{ex}

Let us consider a general homomorphism $\alpha\colon H\to G$ now. Then $(M\times\alpha)^*$ is right Quillen with respect to the $\mathcal G_{\mathcal E,H}$- and $\mathcal G_{\mathcal E,G}$-model structures for $\mathcal E$ as above, so Ken Brown's Lemma implies that $\alpha_!$ preserves weak equivalences between cofibrant objects. On the other hand, Proposition~\ref{prop:alpha-shriek-injective} says that $\alpha_!$ is fully homotopical if $\alpha$ is injective. The following proposition interpolates between these two results:

\begin{prop}\label{prop:free-quotient-general}
Let $\mathcal E$ be a collection of finite subgroups of $M_0$ closed under subconjugates, let $\alpha\colon H\to G$ be a homomorphism, and let $f\colon X\to Y$ be a $\mathcal G_{\mathcal E,H}$-weak equivalence in $\cat{$\bm{(M\times H)}$-SSet}$ such that $\ker(\alpha)$ acts freely on both $X$ and $Y$. Then $\alpha_!(f)$ is a $\mathcal G_{\mathcal E,G}$-weak equivalence.
\begin{proof}
By Proposition~\ref{prop:alpha-shriek-injective} we may assume without loss of generality that $\alpha$ is the quotient map $H\to H/\ker(\alpha)$, so that the functor $\alpha_!$ can be modelled by quotiening out the action of the normal subgroup $K\mathrel{:=}\ker(\alpha)$.

The following splitting follows from a simple calculation similar to the above arguments, which we omit. It can also be obtained from the discrete special case of \cite[Lemma~A.1]{hausmann-equivariant} by adding disjoint basepoints:

\begin{claim*}
Let $L\subset M_0$ be any subgroup and let $\phi\colon L\to H/K$ be any homomorphism. Then we have for any $(L\times H)$-simplicial set $Z$ on which $K$ acts freely a natural isomorphism
\begin{equation*}
\coprod_{[\psi\colon L\to H]} Z^{\psi}/(\centralizer_H(\im\psi)\cap K) \xrightarrow\cong (Z/K)^\phi
\end{equation*}
given on each summand by $[z]\mapsto[z]$. Here the coproduct runs over $K$-conjugacy classes of homomorphisms lifting $\phi$, and $\centralizer_H$ denotes the centralizer in $H$.\nomenclature[aCH]{$\centralizer_H(K)$}{centralizer of $K$ in $H$}\qed
\end{claim*}

We can now prove the proposition. Let $L\in\mathcal E$ and let $\phi\colon L\to H$. In order to show that $(f/K)^\phi$ is a weak homotopy equivalence it suffices by the claim that $f^\psi/(\centralizer_H(\im\psi)\cap K)$ be a weak homotopy equivalence for all  $\psi\colon L\to H$ lifting $\phi$. But indeed, as $K$ acts freely on $X$ and $Y$, so does $\centralizer_H(\im\psi)\cap K$; in particular, it also acts freely on $X^\psi$ and $Y^\psi$. The claim follows as $f^\psi$ is a weak homotopy equivalence by assumption and since free quotients preserve weak homotopy equivalences of simplicial sets (for example by the special case $M=1$ of the above discussion).
\end{proof}
\end{prop}
\index{functoriality in homomorphisms!for M-SSet@for $\cat{$\bm{M}$-SSet}$|)}

\subsubsection{Change of subgroups} We now turn to the special case that the monoid $M$ is fixed (i.e.~$\alpha=\id_M$), but the collection $\mathcal F$ is allowed to vary. For this we will use the notion of \emph{quasi-localizations}, which we recall in Appendix~\ref{appendix:quasi-loc}; in particular, we will employ the notation $\mathscr C^\infty_W$ (or simply $\mathscr C^\infty$) introduced there for `the' quasi-localization of a category $\mathscr C$ at a class $W$ of maps.

\begin{prop}\label{prop:change-of-family-sset}
Let $M$ be a simplicial monoid and let $\mathcal F,\mathcal F'$ be collections of finite subgroups of $M_0$ such that $\mathcal F'\subset\mathcal F$. Then the identity descends to a quasi-localization
\begin{equation}\label{eq:to-smaller-family}
\cat{$\bm M$-SSet}_{\textup{$\mathcal F$-weak equivalences}}^\infty\to\cat{$\bm M$-SSet}_{\textup{$\mathcal F'$-weak equivalences}}^\infty
\end{equation}
at the $\mathcal F'$-weak equivalences, and this functor admits both a left adjoint $\lambda$\nomenclature[alambda]{$\lambda$}{left adjoint to forgetful functor along change of subgroups} as well as a right adjoint $\rho$.\nomenclature[arho1]{$\rho$}{right adjoint to forgetful functor along change of subgroups} Both $\lambda$ and $\rho$ are fully faithful.
\begin{proof}
The identity obviously descends to the quasi-localization $(\ref{eq:to-smaller-family})$. It then only remains to construct the desired adjoints, as they will automatically be fully faithful as adjoints of quasi-localizations, see e.g.~\cite[Proposition~7.1.17]{cisinski-book}.

But indeed, Lemma~\ref{lemma:alpha-shriek-projective} specializes to yield a Quillen adjunction
\begin{equation*}
\id\colon\cat{$\bm M$-SSet}_{\textup{$\mathcal F'$-equivariant}}\rightleftarrows\cat{$\bm M$-SSet}_{\textup{$\mathcal F$-equivariant}} :\!\id
\end{equation*}
so that the left derived functor $\textbf{L}\id$ in the sense of Theorem~\ref{thm:derived-adjunction} defines the desired left adjoint. To construct the right adjoint, we observe that while
\begin{equation*}
\id\colon\cat{$\bm M$-SSet}_{\textup{$\mathcal F$-equivariant}}\rightleftarrows\cat{$\bm M$-SSet}_{\textup{$\mathcal F'$-equivariant}} :\!\id
\end{equation*}
is typically \emph{not} a Quillen adjunction with respect to the usual model structures, it becomes one if we use Corollary~\ref{cor:enlarge-generating-cof} to enlarge the cofibrations on the right hand side to contain all generating cofibrations of the $\mathcal F$-model structure (which we are allowed to do by Lemma~\ref{lemma:homotopy-pushout-M-SSet}), or alternatively that it is a Quillen adjunction for the corresponding injective model structures by Lemma~\ref{lemma:alpha-star-injective}.
\end{proof}
\end{prop}

\begin{rk}\label{rk:change-of-family-sset-explicit}
By the above proof, $\lambda$ can be modelled by taking a cofibrant replacement with respect to the $\mathcal F'$-model structure.
\end{rk}

\section[$G$-global homotopy theory via monoid actions]{\texorpdfstring{\except{toc}{$\bm G$}\for{toc}{$G$}}{G}-global homotopy theory via monoid actions}\label{sec:equivariant-models}
\subsection{The universal finite group}
Schwede \cite{schwede-orbi} proved that unstable global homotopy theory with respect to all compact Lie groups can be modelled by spaces with the action of a certain topological monoid $\mathcal L$, that he calls the \emph{universal compact Lie group},\index{universal compact Lie group} and which we will recall in Section~\ref{sec:global-vs-g-global}. For unstable global homotopy theory with respect to \emph{finite} groups, we will instead be interested in a certain discrete analogue $\mathcal M$, which (under the name $M$) also plays a central role in Schwede's approach \cite{schwede-k-theory} to global algebraic $K$-theory.

\begin{defi}
We write $\omega=\{0,1,2,\dots\}$,\nomenclature[aomega]{$\omega$}{set of non-negative integers} and we denote by $\mathcal M$\nomenclature[aM]{$\mathcal M$}{`universal finite group,' monoid of self-injections of $\omega$} the monoid (under composition) of all injections $\omega\to\omega$.
\end{defi}

Analogously to the Lie group situation \cite[Definition~1.6]{schwede-orbi}, when we model a `global space' by an $\mathcal M$-simplicial set $X$, we do not expect the fixed point spaces $X^H$ for all finite $H\subset\mathcal M$ to carry homotopical information, but only those for certain so-called \emph{universal} $H$. To define these, we first need the following terminology, cf.~\cite[Definition~2.16]{schwede-k-theory}:

\begin{defi}\label{defi:set-universe}
Let $H$ be any finite group. A countable $H$-set $\mathcal U$ is called a \emph{complete $H$-set universe}\index{complete H-set universe@complete $H$-set universe|textbf} if the following equivalent conditions hold:
\begin{enumerate}
\item Every finite $H$-set embeds $H$-equivariantly into $\mathcal U$.
\item Every countable $H$-set embeds $H$-equivariantly into $\mathcal U$.
\item There exists an $H$-equivariant isomorphism
\begin{equation*}
\mathcal U\cong\coprod_{i=0}^\infty\coprod_{\substack{{\scriptstyle K\subset H}\\{\scriptstyle\text{subgroup}}}} H/K.
\end{equation*}\label{item:su-concrete-example}
\item Every subgroup $K\subset H$ occurs as stabilizer of infinitely many distinct elements of $\mathcal U$.
\end{enumerate}
\end{defi}

The proof that the above conditions are indeed equivalent is easy and we omit it. For all of these statements except for the second one this also appears without proof as \cite[Proposition~2.17 and Example~2.18]{schwede-k-theory}.

\begin{rk}
As will become more apparent in the alternative models of unstable $G$-global homotopy theory studied in Section~\ref{sec:diag-spaces} as well as later in the stable context (Chapter~\ref{chapter:stable}), the complete $H$-set universes should be thought of as combinatorial analogues of the \emph{complete $H$-universes} (countably infinite dimensional $H$-representations into which all finite $H$-representations embed) considered e.g.~in the approach to equivariant stable homotopy theory based on orthogonal $H$-spectra, also see Remark~\ref{rk:set-universe-vs-universe} and Proposition~\ref{prop:Fin-global-we} near the end of this chapter.
\end{rk}

The following lemmas are again straightforward to prove from the definitions, and they also appear without proof as part of \cite[Proposition~2.17]{schwede-k-theory}.

\begin{lemma}\label{lemma:supersets-of-universe}
Let $\mathcal U\subset \mathcal V$ be $H$-sets, assume $\mathcal U$ is a complete $H$-set universe and $\mathcal V$ is countable. Then also $\mathcal V$ is a complete $H$-set universe.\qed
\end{lemma}

\begin{lemma}\label{lemma:restriction-universe}
Let $\mathcal U$ be a complete $H$-set universe and let $\alpha\colon K\to H$ be an injective group homomorphism. Then $\alpha^*\mathcal U$ (i.e.~$\mathcal U$ with $K$-action given by $k.x=\alpha(k).x$) is a complete $K$-set universe.\qed
\end{lemma}

\begin{defi}
A finite subgroup $H\subset\mathcal M$ is called \emph{universal}\index{universal subgroup!for M@for $\mathcal M$|textbf} if the restriction of the tautological $\mathcal M$-action on $\omega$ to $H$ makes $\omega$ into a complete $H$-set universe.
\end{defi}

Lemma~\ref{lemma:restriction-universe} immediately implies:

\begin{cor}\label{cor:subgroup-universal-subgroup-universal}
Let $K\subset H\subset\mathcal M$ be subgroups and assume that $H$ is universal. Then also $K$ is universal.\qed
\end{cor}

\begin{lemma}\label{lemma:uniqueness-of-universal-homom}
Let $H$ be any finite group. Then there exists an injective monoid homomorphism $i\colon H\to\mathcal M$ with universal image. Moreover, if $j\colon H\to\mathcal M$ is another such homomorphism, then there exists a $\phi\in\core\mathcal M$ such that
\begin{equation}\label{eq:phi-def-condition}
j(h)=\phi i(h)\phi^{-1}
\end{equation}
for all $h\in H$.
\begin{proof}
This is similar to \cite[Proposition~1.5]{schwede-orbi}: first, we observe that there exists a complete $H$-set universe, for example
\begin{equation*}
\mathcal U\mathrel{:=}\coprod_{i=0}^\infty\coprod_{\substack{{\scriptstyle K\subset H}\\{\scriptstyle\text{subgroup}}}} H/K.
\end{equation*}
As this is countable, we can pick a bijection of sets $\omega\cong\mathcal U$, which yields an $H$-action on $\omega$ turning it into a complete $H$-set universe. The $H$-action amounts to a homomorphism $H\to\Sigma_\omega\subset\mathcal M$ which is injective as the action is faithful, providing the desired homomorphism $i$.

If $j$ is another such homomorphism, then both $i^*\omega$ and $j^*\omega$ are complete $H$-set universes, and hence there exists an $H$-equivariant isomorphism $\phi\colon i^*\omega\to j^*\omega$, see part $(\ref{item:su-concrete-example})$ of Definition~\ref{defi:set-universe}. The $H$-equivariance of $\phi$ then precisely means that $\phi(i(h)(x))=j(h)(\phi(x))$ for all $h\in H$ and $x\in\omega$, i.e.~$\phi i(h)\phi^{-1}=j(h)$ as desired.
\end{proof}
\end{lemma}

\subsection{Global homotopy theory}\label{subsec:global}\index{global homotopy theory|(}
Before we introduce models of unstable $G$-global homotopy theory based on the above monoid $\mathcal M$ in the next subsection, let us first consider the ordinary global situation in order to present some of the main ideas without being overly technical.

Heuristically, we would like to think of a `global space' $X$ as having for each \emph{abstract} finite group $H$ a fixed point space $X^H$ and for each abstract group homomorphism $f\colon H\to K$ a suitably functorial restriction map $f^*\colon X^K\to X^H$. One way to make this heuristic rigorous is given by Schwede's orbispace model, see \cite[Theorem~2.12]{schwede-cat}---in fact, it turns out that there is also some additional $2$-functoriality. The above lemma already tells us that we can assign to an $\mathcal M$-simplicial set $X$ for each abstract finite group $H$ an essentially unique fixed point space $X^H$ as follows: we pick an injective group homomorphism $i\colon H\to\mathcal M$ with universal image and set $X^H\mathrel{:=}X^{i(H)}$. This is indeed independent of the group homomorphism $i$ up to isomorphism: namely, if $j$ is another such group homomorphism, the lemma provides us with a $\phi$ such that $j(h)=\phi i(h)\phi^{-1}$ for all $h\in H$, and $\phi.\blank\colon X\to X$ clearly restricts to an isomorphism $X^{i(H)}\to X^{j(H)}$.

However, there is an issue here---namely, the element $\phi$ (or more precisely, its action) is not canonical (in particular, it is not clear how to define $1$-functoriality):

\begin{ex}\label{ex:too-many-endomorphisms}
Let us consider the special case $H=1$, so that there is in particular only one homomorphism $i=j\colon H\to\mathcal M$. Then \emph{any} $\phi\in\core\mathcal M$ satisfies the condition $(\ref{eq:phi-def-condition})$ and hence produces an endomorphism $\phi.\blank$ of $X=X^{\{1\}}$.

Looking at the orbispace model, we should expect all such endomorphisms of $X^{\{1\}}$ to be homotopically trivial. However, for $X=\mathcal M$ any $\phi\not=1$ gives us a non-trivial endomorphism.
\end{ex}

This suggests that $\mathcal M$-simplicial sets with respect to maps inducing weak homotopy equivalences on fixed points for universal subgroups are not yet a model of unstable global homotopy theory. In order to solve the issue raised in the example, we will enhance $\mathcal M$ to a simplicial (or categorical) monoid in particular trivializing the above action. This uses:

\begin{constr}
Let $X$ be any set. We write $EX$ for the (small) category with objects $X$ and precisely one morphism $x\to y$ for each $x,y\in X$, which we denote by $(y,x)$. We extend $E$ to a functor $\cat{Set}\to\cat{Cat}$ in the obvious way.\nomenclature[aE]{$E$}{right adjoint to $\Ob\colon\cat{Cat}\to\cat{Set}$; right adjoint to $\ev_0\colon\cat{SSet}\to\cat{Set}$}

We will moreover also write $EX$ for the simplicial set given in degree $n$ by
\begin{equation*}
(EX)_n=X^{\times (1+n)}\cong\Maps(\{0,\dots,n\},X)
\end{equation*}
with structure maps via restriction and with the evident functoriality in $X$. We remark that the simplicial set $EX$ is indeed canonically isomorphic to the nerve of the category $EX$, justifying the clash of notation. In fact, it will be useful at several points to switch between viewing $EX$ as a category or as a simplicial set.
\end{constr}

\begin{rk}
It is clear that the category $EX$ is a groupoid and that it is contractible for $X\not=\varnothing$. In particular, the simplicial set $EX$ is a Kan complex, again contractible unless $X=\varnothing$.
\end{rk}

The functor $E\colon\cat{Set}\to\cat{Cat}$ is right adjoint to the functor $\Ob\colon\cat{Cat}\to\cat{Set}$ sending a small category to its set of objects; likewise $E\colon\cat{Set}\to\cat{SSet}$ is right adjoint to the functor sending a simplicial set to its set of zero simplices. In particular, $E$ preserves products, so $E\mathcal M$\nomenclature[aEM]{$E\mathcal M$}{simplicial (or categorical) monoid obtained from $\mathcal M$} is canonically a simplicial monoid. As it is contractible, any two translations $u.\blank,v.\blank$ for $u,v\in\mathcal M$ are homotopic on any $E\mathcal M$-object $X$; in fact, there is unique edge $(v,u)$ from $u$ to $v$ in $E\mathcal M$, and acting with this gives an explicit homotopy $u.\blank\Rightarrow v.\blank$.

We conclude that $E\mathcal M$ avoids the issue detailed in Example~\ref{ex:too-many-endomorphisms}. Indeed, Theorems~\ref{thm:script-I-global-model-structure} and~\ref{thm:strict-global-I-model-structure} together with Theorem~\ref{thm:L-vs-I} will show that $\cat{$\bm{E\mathcal M}$-SSet}$ is a model of global homotopy theory in the sense of \cite{schwede-book} with respect to finite groups. Maybe somewhat surprisingly, the main result of this subsection (Theorem~\ref{thm:m-vs-em}) will be that the same homotopy theory can be modelled by $\mathcal M$-simplicial sets with respect to a slightly intricate notion of weak equivalence.

\begin{ex}\label{ex:global-classifying space}
Let $H$ be any finite group and let $A$ be a countable faithful $H$-set. The set $\Inj(A,\omega)$ of injections $A\to\omega$ \nomenclature[aInj]{$\Inj(A,B)$}{set of injective maps $A\to B$} has a natural $\mathcal M$-action via postcomposition and a commuting $H$-action via precomposition, inducing an $E\mathcal M$-action on $E\Inj(A,\omega)/H$.

Note that $H$ acts freely from the right on $\Inj(A,\omega)$ and hence on the contractible simplicial set $E\Inj(A,\omega)$ as injections of sets are monomorphisms. In particular, ignoring the $E\mathcal M$-action, $E\Inj(A,\omega)/H$ is just an Eilenberg-Mac Lane space of type $K(H,1)$. However, the $E\mathcal M$-simplicial set $E\Inj(A,\omega)/H$ contains interesting additional equivariant information compared to an ordinary $K(H,1)$ equipped with trivial $E\mathcal M$-action, and we call it `the' \emph{global classifying space} of $H$.
\end{ex}

The basis for our comparison between $\cat{$\bm{E\mathcal M}$-SSet}$ and $\cat{$\bm{\mathcal M}$-SSet}$ will be the following model categories provided by Proposition~\ref{prop:equiv-model-structure}:

\begin{defi}
A map $f\colon X\to Y$ in $\cat{$\bm{\mathcal M}$-SSet}$ is called a \emph{universal weak equivalence}\index{universal weak equivalence|textbf} or \emph{universal fibration} if $f^H$ is a weak homotopy equivalence or fibration, respectively, for every universal $H\subset\mathcal M$.
\end{defi}

\begin{cor}\label{cor:m-model-structure}
The universal weak equivalences and fibrations are part of a unique model structure on $\cat{$\bm{\mathcal M}$-SSet}$, which we call the \emph{universal model structure}\index{universal model structure|textbf}. It is simplicial, combinatorial, proper, and filtered colimits in it are homotopical. A possible set of generating cofibrations is given by
\begin{equation*}
I=\{ \mathcal M/H\times\del\Delta^n\hookrightarrow \mathcal M/H\times\Delta^n : n\ge0,\text{ $H\subset\mathcal M$ universal}\}
\end{equation*}
and a possible set of generating acyclic cofibrations is given by
\begin{equation*}
J=\{ \mathcal M/H\times\Lambda^n_k\hookrightarrow \mathcal M/H\times\Delta^n : 0\le k\le n,\text{ $H\subset\mathcal M$ universal}\}.\qedNOW
\end{equation*}
\end{cor}

As already mentioned above, this will not yet model global equivariant homotopy theory, so we reserve the names `global model structure' and `global weak equivalences' for a different model structure.

\begin{defi}
A map $f\colon X\to Y$ of $E\mathcal M$-simplicial sets is called a \emph{global weak equivalence}\index{global weak equivalence|seealso{$G$-global weak equivalence}}\index{global weak equivalence!in EM-SSet@in $\cat{$\bm{E\mathcal M}$-SSet}$|textbf} or \emph{global fibration} if it is a universal weak equivalence or universal fibration, respectively, when considered as a map in $\cat{$\bm{\mathcal M}$-SSet}$. 
\end{defi}

\begin{cor}\label{cor:EM-SSet-model-structure}
The global weak equivalences and global fibrations are part of a unique model structure on $\cat{$\bm{E\mathcal M}$-SSet}$, which we call the \emph{global model structure}\index{global model structure|seealso{$G$-global model structure}}\index{global model structure!on EM-SSet@on $\cat{$\bm{E\mathcal M}$-SSet}$|textbf}. It is simplicial, proper, combinatorial, and filtered colimits in it are homotopical. A possible set of generating cofibrations is given by
\begin{equation*}
I=\{(E\mathcal M)/H\times\del\Delta^n\hookrightarrow (E\mathcal M)/H\times\Delta^n : n\ge0,\text{ $H\subset\mathcal M$ universal}\}
\end{equation*}
and a possible set of generating acyclic cofibrations by
\begin{equation*}
J=\{(E\mathcal M)/H\times\Lambda^n_k\hookrightarrow (E\mathcal M)/H\times\Delta^n : 0\le k\le n,\text{ $H\subset\mathcal M$ universal}\}.\qedNOW
\end{equation*}
\end{cor}

\begin{constr}
The forgetful functor $\cat{$\bm{E\mathcal M}$-SSet}\to\cat{$\bm{\mathcal M}$-SSet}$ admits both a simplicial left and a simplicial right adjoint. While they exist for abstract reasons (e.g.~as simplicially enriched Kan extensions), they are also easy to make explicit:

Let $X$ be any $\mathcal M$-simplicial set. We write $E\mathcal M\times_{\mathcal M}X$\nomenclature[aEMM]{$E\mathcal M\times_{\mathcal M}\blank$}{left adjoint to forgetful functor $\cat{$\bm{E\mathcal M}$-SSet}\to\cat{$\bm{\mathcal M}$-SSet}$} for the following $E\mathcal M$-simplicial set: as a simplicial set, this is the quotient of $E\mathcal M\times X$ under the equivalence relation generated in degree $n$ by $(u_0v,\dots,u_nv;x)\sim(u_0,\dots,u_n;v.x)$ for all $u_0,\dots,u_n,v\in\mathcal M$ and $x\in X_n$. As usual, we denote the class of $(u_0,\dots,u_n;x)$ by $[u_0,\dots,u_n;x]$. The $E\mathcal M$-action on $E\mathcal M\times_{\mathcal M}X$ is induced by the obvious $E\mathcal M$-action on the first factor. If $f\colon X\to Y$ is any $\mathcal M$-equivariant map, then $E\mathcal M\times_{\mathcal M}f$ is induced by $E\mathcal M\times f$, and similarly for higher cells $\Delta^n\times X\to Y$. We omit the easy verification that is well-defined.

We have a natural $\mathcal M$-equivariant map $\eta\colon X\to\forget(E\mathcal M\times_{\mathcal M}X)$ given in degree $n$ by sending an $n$-simplex $x$ to the class $[1,\dots,1;x]$. Moreover, if $Y$ is an $E\mathcal M$-simplicial set, we have an $E\mathcal M$-equivariant map $\epsilon\colon E\mathcal M\times_{\mathcal M}(\forget Y)\to Y$ giving in degree $n$ by acting, i.e.~$[u_0,\dots,u_n;y]\mapsto(u_0,\dots,u_n).y$. We leave the easy verification to the reader that these are well-defined, enriched natural, and define unit and counit, respectively, of a simplicial adjunction $E\mathcal M\times_{\mathcal M}\blank\dashv\forget$.

Similarly, the forgetful functor has a simplicial right adjoint $\Maps^{\mathcal M}(E\mathcal M,\blank)$\nomenclature[amapsMEM]{$\Maps^{\mathcal M}(E\mathcal M,\blank)$}{right adjoint to forgetful functor $\cat{$\bm{E\mathcal M}$-SSet}\to\cat{$\bm{\mathcal M}$-SSet}$} (the simplicial set of $\mathcal M$-equivariant maps), with $E\mathcal M$-action induced by the right $E\mathcal M$-action on itself via postcomposition.
\end{constr}

By definition of the model structures we immediately get:

\begin{cor}\label{cor:forget-right-Quillen}
The simplicial adjunction
\begin{equation}\label{eq:em-times-forget}
E\mathcal M\times_{\mathcal M}\blank\colon \cat{$\bm{\mathcal M}$-SSet}\rightleftarrows\cat{$\bm{E\mathcal M}$-SSet}:\!\forget
\end{equation}
is a Quillen adjunction and the right adjoint creates weak equivalences.\qed
\end{cor}

As already suggested by our heuristic above, this is not a Quillen equivalence; more precisely, $\forget^\infty$ is not essentially surjective:

\begin{ex}
If $X$ is any $E\mathcal M$-simplicial set, then all $u\in\mathcal M$ act on $X$ by weak homotopy equivalences, and hence the same will be true for any $\mathcal M$-simplicial set $Y$ weakly equivalent to $\forget X$. On the other hand, $\mathcal M$ considered as a discrete $\mathcal M$-simplicial set with $\mathcal M$-action by postcomposition does not satisfy this, and hence can't lie in the essential image.
\end{ex}

In the example we have only looked at the underlying non-equivariant homotopy type of a given $\mathcal M$-simplicial set. However, in order to have a sufficient criterion for the essential image of the forgetful functor, we should better take equivariant information into account. While the translation maps will usually not be $\mathcal M$-equivariant (related to the fact that $\mathcal M$ is highly non-commutative), we can look at those parts of the action that we can still expect to be preserved:

\begin{defi}\label{defi:semistable}
An $\mathcal M$-simplicial set $X$ is called \emph{semistable}\index{semistable!M-simplicial set@$\mathcal M$-simplicial set|textbf}\index{semistable|seealso{$G$-semistable}} if for each universal subgroup $H\subset\mathcal M$ and each $u\in\mathcal M$ centralizing $H$ the translation map
\begin{equation*}
u.\blank\colon X\to X
\end{equation*}
is an $H$-equivariant weak equivalence.
\end{defi}

We observe that this is equivalent to demanding that for any such $u\in\mathcal M$ and $H\subset\mathcal M$ the restriction of $u.\blank$ to $X^H\to X^H$ is a weak homotopy equivalence (this uses Corollary~\ref{cor:subgroup-universal-subgroup-universal}).

\begin{ex}\label{ex:em-semistable}
Strengthening the previous example, any $E\mathcal M$-simplicial set $X$ is in fact semistable when viewed as an $\mathcal M$-simplicial set: namely, if $u\in\mathcal M$ centralizes an (arbitrary) subgroup $H\subset\mathcal M$, then $(u,1)$ provides an $H$-equivariant homotopy between $u.\blank$ and the identity.
\end{ex}

\begin{rk}
The term `semistable' refers to Schwede's characterization \cite[Theorem~4.1 and Lemma~2.3-(iii)]{schwede-semistable} of semistable symmetric spectra, i.e.~symmetric spectra whose na\"ive homotopy groups agree with their true homotopy groups,\index{semistable!symmetric spectrum} as those spectra for which a certain canonical $\mathcal M$-action on the na\"ive homotopy groups is given by isomorphisms, also cf.~\cite[Corollary~3.32 and Proposition~3.16]{hausmann-equivariant} for a similar characterization in the equivariant case due to Hausmann.

Both Schwede and Hausmann prove that in the respective situation the action on na\"ive homotopy groups is actually trivial, i.e.~all elements of $\mathcal M$ act by the identity. Likewise, it will follow from Theorem \ref{thm:m-vs-em} together with the argument from Example~\ref{ex:em-semistable} above that for a semistable $\mathcal M$-simplicial set the translation $u.\blank$, $u\in\mathcal M$ centralizing some universal subgroup $H\subset\mathcal M$, is in fact the identity in the $H$-equivariant homotopy category.
\end{rk}

Obviously, semistability is invariant under universal weak equivalences, so it is a necessary condition to lie in the essential image of $\forget$ by the above example. As the main result of this subsection, we will show that it is also sufficient, and moreover the above is everything that prevents $\forget^\infty$ from being an equivalence:

\begin{thm}\label{thm:m-vs-em}
The adjunction $(\ref{eq:em-times-forget})$ induces a Bousfield localization
\begin{equation*}
E\mathcal M\times_{\mathcal M}^{\textbf{\textup L}}\blank\colon \cat{$\bm{\mathcal M}$-SSet}^\infty\rightleftarrows\cat{$\bm{E\mathcal M}$-SSet}^\infty:\!\forget^\infty;
\end{equation*}
in particular, $\forget^\infty$ is fully faithful. Moreover, its essential image consists precisely of the semistable $\mathcal M$-simplicial sets.
\end{thm}

Here we have tacitly identified the objects of $\cat{$\bm{\mathcal M}$-SSet}^\infty$ with those of $\cat{$\bm{\mathcal M}$-SSet}$, cf.~Remark~\ref{rk:associated-quasi}.

\begin{rk}
The theorem in particular tells us that the forgetful functor identifies $\cat{$\bm{E\mathcal M}$-SSet}^\infty$ with the full subcategory of $\cat{$\bm{\mathcal M}$-SSet}^\infty$ spanned by the semistable objects. In view of Proposition~\ref{prop:localization-subcategory}, the latter is canonically identified with the quasi-localization of semistable $\mathcal M$-simplicial sets at the universal weak equivalences, so we do not have to be careful to distinguish them. Put differently: semistable $\mathcal M$-simplicial sets with respect to the universal weak equivalences are a model of unstable global homotopy theory.
\end{rk}

While we usually employ model categorical techniques to prove results about quasi-categories in this monograph, our proof of the above theorem proceeds the other way round: namely, we will exhibit the \emph{global orbit category} $\textbf{O}_{E\mathcal M}$\index{global orbit category|textbf}\index{global orbit category|seealso{$G$-global orbit category}} as an explicit simplicial localization\index{simplicial localization} of $\textbf{O}_{\mathcal M}$ in the sense of Definition~\ref{defi:simplicial-localization}, and then deduce the theorem from the universal property of simplicial localizations (or more precisely its model categorical manifestation Theorem~\ref{thm:simpl-localization-model-cat}) together with Elemendorf's Theorem (Corollary~\ref{cor:elmendorf}).

To do so, let us begin by understanding these categories a bit better:

\index{global orbit category|(}
\begin{lemma}\label{lemma:group-hom-associated}
Let $H,K\subset\mathcal M$, and let $u_0,\dots,u_n\in\mathcal M$ such that $[u_0,\dots,u_n]$ is $H$-fixed in $(E\mathcal M/K)_n$. Then there exists for any $h\in H$ a unique $\sigma(h)\in K$ such that $hu_i=u_i\sigma(h)$ for all $i$.  Moreover, $\sigma\colon H\to K$ is a homomorphism.

Conversely, whenever such a map $\sigma$ exists, $[u_0,\dots,u_n]$ is $H$-fixed in $(E\mathcal M)/K$.
\begin{proof}
As $[u_0,\dots,u_n]$ is $H$-fixed, $(u_0,\dots,u_n)\sim (hu_0,\dots,hu_n)$ for any $h\in H$,
so by definition there indeed exists some $\sigma(h)\in K$ such that $hu_i=u_i\sigma(h)$ for all $i$; moreover, $\sigma(h)$ is unique as $K$ acts freely from the right on $\mathcal M$.

To check that $\sigma$ is a group homomorphism, let $h_1,h_2\in H$ arbitrary. Then $h_1h_2u=h_1u\sigma(h_2)=u\sigma(h_1)\sigma(h_2)$, hence $\sigma(h_1h_2)=\sigma(h_1)\sigma(h_2)$ by uniqueness.

The proof of the converse is trivial.
\end{proof}
\end{lemma}

\begin{constr}
By definition, $\textbf{O}_{E\mathcal M}\subset\cat{$\bm{E\mathcal M}$-SSet}$ is the full simplicial subcategory spanned by the $E\mathcal M/H$ for universal $H\subset\mathcal M$. We have seen in Lemma~\ref{lemma:fixed-points-cellular} that the simplicial set $\Maps_{\textbf{O}_{E\mathcal M}}(E\mathcal M/H, E\mathcal M/K)$ is isomorphic to $(E\mathcal M/K)^H$ via evaluation at $[1]\in E\mathcal M/H$. On $0$-simplices, this gives $\Maps_{\textbf{O}_{E\mathcal M}}(E\mathcal M/H, E\mathcal M/K)_0\cong (\mathcal M/K)^H$; an inverse is then  given by sending $u\in(\mathcal M/K)^H$ to $\blank\cdot u\colon [v_0,\dots,v_n]\mapsto[v_0u,\dots,v_nu]$.

More generally, an $n$-simplex of the mapping space $\Maps_{\textbf{O}_{E\mathcal M}}(E\mathcal M/H, E\mathcal M/K)$ can be represented by an $(n+1)$-tuple $(u_0,\dots,u_n)$ such that $[u_0,\dots,u_n]\in E\mathcal M/K$ is $H$-fixed, which by Lemma~\ref{lemma:group-hom-associated} is equivalent to the existence of a group homomorphism $\sigma\colon H\to K$ such that $hu_i=u_i\sigma(h)$ for all $i=0,\dots,n$. Two tuples represent the same morphism iff they become equal in $E\mathcal M/K$, i.e.~iff they only differ by right multiplication with some $k\in K$.

Moreover, one immediately sees by direct inspection, that if the $n$-simplex $f\in\Maps_{\textbf{O}_{E\mathcal M}}(E\mathcal M/H,E\mathcal M/K)_n$ is represented by $(u_0,\dots,u_n)$ and the $n$-simplex $f'\in\Maps_{\textbf{O}_{E\mathcal M}}(E\mathcal M/K,E\mathcal M/L)_n$ is represented by $(u_0',\dots,u_n')$, then their composition $f'f$ is represented by $(u_0u_0',\dots,u_nu_n')$ (note the different order!).

Likewise, $\textbf{O}_{\mathcal M}$ is the full (simplicial) subcategory of $\cat{$\bm{\mathcal M}$-SSet}$ spanned by the $\mathcal M/H$, and we have an isomorphism $\Hom(\mathcal M/H,\mathcal M/K)\cong (\mathcal M/K)^H$ via evaluation at $[1]$; composition is again induced from multiplication in $\mathcal M$.

We now define a functor $i\colon\textbf{O}_{\mathcal M}\to\textbf{O}_{E\mathcal M}$ as follows: an object $\mathcal M/H$ is sent to $E\mathcal M/H$. On morphism spaces, $i$ is given as the composition
\begin{align*}
\Maps_{\textbf{O}_{\mathcal M}}(\mathcal M/H,\mathcal M/K)&\cong (\mathcal M/K)^H\cong\Maps_{\textbf{O}_{E\mathcal M}}(E\mathcal M/H,E\mathcal M/K)_0\\ &\hookrightarrow\Maps_{\textbf{O}_{E\mathcal M}}(E\mathcal M/H,E\mathcal M/K),
\end{align*}
which then just sends the morphism represented by $u$ to the morphism represented by the same element. As an upshot of the above discussion, this is indeed functorial.
\end{constr}

\begin{defi}\index{centralizing morphism|textbf}\index{centralizing morphism|seealso{$G$-centralizing morphism}}
Let $H\subset\mathcal M$ be universal. A map $f\colon \mathcal M/H\to \mathcal M/H$ is called \emph{centralizing} if there exists a $u\in\mathcal M$ centralizing $H$ such that $f$ is given by right multiplication by $u$. Analogously, we define centralizing morphisms in $\textbf{O}_{E\mathcal M}$.
\end{defi}

The following will be the main ingredient to the proof of Theorem~\ref{thm:m-vs-em}:

\begin{prop}\label{prop:em-quasi-localization}\index{global orbit category!as a simplicial localization}
The functor $i\colon\textbf{\textup O}_{\mathcal M}\to \textbf{\textup O}_{E\mathcal M}$ is a simplicial localization at the centralizing morphisms.
\begin{proof}
By construction, $i$ induces an isomorphism onto the underlying category $\und\textbf{O}_{E\mathcal M}=(\textbf{O}_{E\mathcal M})_0$\nomenclature[au]{$\und$}{underlying category (of $0$-simplices) of a simplicially enriched category} of $\textbf{O}_{E\mathcal M}$. Thus, it is enough to prove that $\und\textbf{O}_{E\mathcal M}\hookrightarrow\textbf{O}_{E\mathcal M}$ is a quasi-localization at the centralizing morphisms, for which we will verify the assumptions of Proposition~\ref{prop:enrichment-vs-localization}, i.e.~that $\cat{O}_{E\mathcal M}$ has fibrant mapping spaces, that the centralizing morphisms are homotopy equivalences, and that the functors
\begin{equation}\label{eq:emql-iterated-degeneracy}
s^*\colon\big((\textbf{O}_{E\mathcal M})_0,W\big)\to \big((\textbf{O}_{E\mathcal M})_n,s^*W\big)
\end{equation}
induce equivalences on quasi-localizations, where $W$ denotes the class of centralizing morphisms and $s$ is the unique map $[n]\to[0]$ in $\Delta$.

For the first claim we will show that all mapping spaces in $\cat{O}_{E\mathcal M}$ are actually even nerves of groupoids. Indeed, it suffices to prove that each $E\mathcal M/H$ is, for which it is then in turn enough to observe that $E\mathcal M$ is the nerve of a groupoid by construction and that we can form the quotient by $H$ already in the category of groupoids as the nerve preserves quotients by \emph{free} group actions.

For the second claim we note that if $u$ centralizes $H$, then $\blank\cdot u$ is even homotopic to the identity via the edge $[1,u]$ in $(E\mathcal M/H)^H\cong\Maps_{{\textbf O}_{E\mathcal M}}(E\mathcal M/H,E\mathcal M/H)$.

Finally, for the third claim it is by Corollary~\ref{cor:homotopy-equivalence} enough to prove that $(\ref{eq:emql-iterated-degeneracy})$ is a homotopy equivalence in the sense of Definition~\ref{defi:homotopy-equivalence}. For this we will show that $i\colon[0]\to[n],0\mapsto 0$ induces a homotopy inverse. Indeed, $i^*$ is obviously homotopical and moreover $i^*s^*=(si)^*=\id$ by functoriality. It remains to prove that $s^*i^*$ is homotopic to the identity of $(\textbf{O}_{E\mathcal M})_n$.

We begin by picking for each universal $H\subset\mathcal M$ an $H$-equivariant isomorphism $\omega\amalg\omega\cong\omega$, where $H$ acts on each of the three copies of $\omega$ in the tautological way; such an isomorphism indeed exists as both sides are complete $H$-set universes. Restricting to the two copies of $\omega$ then gives injections $\alpha_H,\beta_H\in\mathcal M$ such that:
\begin{enumerate}
\item $\alpha_H$ and $\beta_H$ centralize $H$\label{item:centralizing}
\item $\omega=\im(\alpha_H)\sqcup\im(\beta_H)$, i.e.~the images of $\alpha_H$ and $\beta_H$ partition $\omega$.\nomenclature[zsqcup]{$\sqcup$}{internal disjoint union of sets}\label{item:disjoint-cover}
\end{enumerate}
We now define $f\colon (\textbf{O}_{E\mathcal M})_n\to(\textbf{O}_{E\mathcal M})_n$ as follows: $f$ is the identity on objects. A morphism $(E\mathcal M)/H\to (E\mathcal M)/K$ represented by $(u_0,\dots,u_n)\in\mathcal M^{n+1}$ is sent to the morphism represented by $(v_0,\dots,v_n)$ where $v_i$ satisfies
\begin{equation}\label{eq:O-M-vs-O-EM-f-defining}
v_i\alpha_K = \alpha_H u_i\qquad\text{and}\qquad v_i\beta_K=\beta_Hu_0.
\end{equation}
We first observe that there is indeed a unique such $v_i$ as $\alpha_K$ and $\beta_K$ are injections whose images form a partition of $\omega$. Moreover, this is an injection as $\alpha_H$ and $\beta_H$ have disjoint image and as both $\alpha_H u_i$ and $\beta_Hu_0$ are injective.

Next, we show that $(v_0,\dots,v_n)$ indeed defines a morphism, i.e.~it represents an $H$-fixed simplex of $(E\mathcal M)/K$. Indeed, as $(u_0,\dots,u_n)$ represents an $H$-fixed simplex, there is a (unique) group homomorphism $\sigma\colon H\to K$ such that $hu_i=u_i\sigma(h)$ for all $h\in H$. But then
\begin{equation*}
hv_i\alpha_K = h\alpha_H u_i=\alpha_H hu_i=\alpha_H u_i\sigma(h)=v_i\alpha_K\sigma(h)=v_i\sigma(h)\alpha_K,
\end{equation*}
where we have used Condition~$(\ref{item:centralizing})$ twice as well as the definition of $v_i$. Analogously one shows $hv_i\beta_K= v_i\sigma(h)\beta_K$; as the images of $\alpha_K$ and $\beta_K$ together cover $\omega$, we conclude that $hv_i=v_i\sigma(h)$ for all $i$ and all $h$, and hence $[v_0,\dots,v_n]\in (E\mathcal M)/K$ is $H$-fixed as desired.

Moreover, this is independent of the choice of representative: if we pick any other representative $(u_0',\dots,u_n')$, then there is some $k\in K$ such that $u_i'=u_ik$ for all $i$, and thus the associated $v_i'$ satisfy
\begin{equation*}
v_i'\alpha_K=\alpha_H u_i'=\alpha_H u_i k=v_i\alpha_K k = v_ik\alpha_K,
\end{equation*}
where we have used the definitions of $v_i$ and $v_i'$ as well as $(\ref{item:centralizing})$. Analogously one shows $v_i'\beta_K= v_ik\beta_K$; as before we conclude that $v_i'=v_ik$, so that $(v_0,\dots,v_n)$ and $(v_0',\dots,v_n')$ represent the same morphism.

With this established, one easily checks that $f$ is a functor $(\textbf{O}_{E\mathcal M})_n\to(\textbf{O}_{E\mathcal M})_n$. As our notion of homotopy equivalence does not require the intermediate functors to be homotopical (although $f$ actually is by a computation analogous to the above), it only remains to prove that $f$ is homotopic to both the identity and $s^*i^*$. For this we observe that we have by the defining equation $(\ref{eq:O-M-vs-O-EM-f-defining})$ natural transformations
\begin{equation*}
\id\Leftarrow f\Rightarrow s^*i^*
\end{equation*}
where the left hand transformation is given on $(E\mathcal M)/H$ by the morphism corresponding to $s^*[\alpha_H]=[\alpha_H,\dots,\alpha_H]$ while the right hand transformation corresponds to $s^*[\beta_H]$, and these are levelwise weak equivalences by Condition~$(\ref{item:centralizing})$ above.
\end{proof}
\end{prop}

\begin{cor}
In the Quillen adjunction
\begin{equation*}
(i^\op)_!\colon \FUN(\textbf{\textup O}_{\mathcal M}^\op,\cat{SSet})\rightleftarrows
\FUN(\textbf{\textup O}_{E\mathcal M}^\op,\cat{SSet}) :\!(i^\op)^*,
\end{equation*}
the right adjoint is homotopical and the induced functor between associated quasi-categories is fully faithful with essential image precisely those simplicial pre\-sheaves on $\textbf{\textup O}_{\mathcal M}$ that invert the centralizing morphisms.
\begin{proof}
It is clear, that $(i^\op)^*$ is homotopical. By the previous proposition, $i\colon\textbf{\textup O}_{\mathcal M}\to\textbf{\textup O}_{E\mathcal M}$ is a simplicial localization at the centralizing morphisms, and hence $i^\op$ is a simplicial localization at their opposites. As both source and target of $i^\op$ are small and locally fibrant, the claim now follows from Theorem~\ref{thm:simpl-localization-model-cat}.\index{global orbit category|)}
\end{proof}
\end{cor}

\begin{proof}[Proof of Theorem~\ref{thm:m-vs-em}]
We already know from Corollary~\ref{cor:forget-right-Quillen} that $E\mathcal M\times_{\mathcal M}(\blank)\dashv\forget$ is a Quillen adjunction with homotopical right adjoint. It therefore suffices to prove that $\forget^\infty$ is fully faithful with essential image the semistable $\mathcal M$-simplicial sets.

\begin{claim*}
The diagram
\begin{equation}\label{diag:forget-vs-restr}
\begin{tikzcd}[column sep=large]
\cat{$\bm{E\mathcal M}$-SSet} \arrow[r, "\forget"]\arrow[d, "\Phi"'] & \cat{$\bm{\mathcal M}$-SSet}\arrow[d, "\Phi"]\\
\FUN(\textbf{\textup O}_{E\mathcal M}^\op, \cat{SSet})\arrow[r, "(i^\op)^*"'] &[4em] \FUN(\textbf{\textup O}_{\mathcal M}^\op, \cat{SSet})
\end{tikzcd}
\end{equation}
of homotopical functors commutes up to natural isomorphism.
\begin{proof}
An explicit choice of such an isomorphism $\tau$ is given as follows: if $X$ is any $E\mathcal M$-simplicial set and $H\subset\mathcal M$ is universal, then $\tau_X(\mathcal M/H)$ is the composition
\begin{align*}
\Phi(X)(i(\mathcal M/H))=\Phi(X)(E\mathcal M/H)&\xrightarrow{\ev_{[1]}} X^H=(\forget X)^H\\ &\xrightarrow{(\ev_{[1]})^{-1}}\Phi(\forget X)(\mathcal M/H).
\end{align*}
To see that this is well-defined, let $H,K\subset\mathcal M$ universal, and let $u\in\mathcal M$ define an $K$-fixed point of $\mathcal M/H$. Then we have commutative diagrams
\begin{equation}\label{diag:Phi-vs-fixed-points-E}
\begin{tikzcd}
\Phi(X)(E\mathcal M/H)=\Maps^{E\mathcal M}(E\mathcal M/H, X)\arrow[d, "(\blank\cdot u)^*"']\arrow[r, "\ev_{[1]}", "\cong"'] & X^H\arrow[d, "u.\blank"]\\
\Phi(X)(E\mathcal M/K)=\Maps^{E\mathcal M}(E\mathcal M/K, X)\arrow[r, "\ev_{[1]}", "\cong"'] & X^K
\end{tikzcd}
\end{equation}
and, for each $Y\in\cat{$\bm{\mathcal M}$-SSet}$,
\begin{equation}\label{diag:Phi-vs-fixed-points}
\begin{tikzcd}
\Phi(Y)(\mathcal M/H)=\Maps^{\mathcal M}(\mathcal M/H, Y)\arrow[d, "(\blank\cdot u)^*"']\arrow[r, "\ev_{[1]}", "\cong"'] & Y^H\arrow[d, "u.\blank"]\\
\Phi(Y)(\mathcal M/K)=\Maps^{\mathcal M}(\mathcal M/K, Y)\arrow[r, "\ev_{[1]}", "\cong"'] & Y^K.
\end{tikzcd}
\end{equation}
Taking $Y=\forget X$ in $(\ref{diag:Phi-vs-fixed-points})$, these two together then show that $\tau_X$ is natural (and hence defines a morphism in $\FUN(\textbf{\textup O}_{\mathcal M}^\op, \cat{SSet})$). It is then obvious that $\tau$ is natural (say, in the unenriched sense), as it is levelwise given by a composition of natural transformations.
\end{proof}
\end{claim*}

We can now deduce the theorem: the vertical functors in $(\ref{diag:forget-vs-restr})$ induce equivalences on quasi-localizations by Corollary~\ref{cor:elmendorf}, and by the previous corollary the bottom horizontal arrow is fully faithful with essential image those pre\-sheaves that invert centralizing isomorphisms. Thus, $\forget^\infty$ is fully faithful with essential image those $\mathcal M$-simplicial sets $X$ such that $\Phi(X)$ inverts centralizing isomorphisms. Taking $Y=X$ and a $u\in\mathcal M$ centralizing $H=K$ in $(\ref{diag:Phi-vs-fixed-points})$, we see that $\Phi(X)$ inverts centralizing morphisms if and only if $X$ is semistable, finishing the proof.
\end{proof}
\index{global homotopy theory|)}

\subsection[$G$-global model structures]{\texorpdfstring{$\bm G$}{G}-global model structures}
Let us fix some (possibly infinite) discrete group $G$. We now want to extend the above discussion to yield a \emph{$G$-global model structure} on the category $\cat{$\bm{E\mathcal M}$-$\bm G$-SSet}$ of simplicial sets with a $G$-action and a commuting $E\mathcal M$-action, which we can equivalently think of as simplicial sets with an action of $E\mathcal M\times G$, or as the category of $G$-objects in $\cat{$\bm{E\mathcal M}$-SSet}$.

For $G=1$ this will recover the previous model; however, as soon as $G$ contains torsion, the weak equivalences will be strictly finer than the underlying global weak equivalences. In particular, we will show later in Theorem~\ref{thm:G-global-vs-proper-sset} that the weak equivalences are fine enough that one can recover proper $G$-equivariant homotopy theory as a Bousfield localization.

\begin{defi}\index{universal graph subgroup|textbf}
A graph subgroup $\Gamma=\Gamma_{H,\phi}$ of $\mathcal M\times G$ is called \emph{universal} if the corresponding subgroup $H\subset\mathcal M$ is universal. A graph subgroup $\Gamma\subset(E\mathcal M\times G)_0$ is universal if it is universal as a subgroup of $\mathcal M\times G$.
\end{defi}

\begin{ex}\nomenclature[aphiG]{$\blank\times_\phi G$}{quotient of $\blank\times G$ by graph subgroup}\label{ex:EM-fixed-points-corep}
Let $H\subset\mathcal M$ be a subgroup and let $\phi\colon H\to G$ be a group homomorphism. We write
\begin{equation*}
\mathcal M\times_\phi G\mathrel{:=} (\mathcal M\times G)/\Gamma_{H,\phi}
\qquad\text{and}\qquad
E\mathcal M\times_\phi G\mathrel{:=} (E\mathcal M\times G)/\Gamma_{H,\phi}
\end{equation*}
and call them \emph{$G$-global classifying spaces}.

It follows immediately from the constructions that $\mathcal M\times_\phi G$ corepresents $(\blank)^\phi$ on $\cat{$\bm{\mathcal M}$-$\bm G$-SSet}$ and that $E\mathcal M\times_\phi G$ corepresents $(\blank)^\phi$ on $\cat{$\bm{E\mathcal M}$-$\bm G$-SSet}$.
\end{ex}

As in the global situation we will now specialize Proposition~\ref{prop:equiv-model-structure}:

\begin{defi}
A map $f\colon X\to Y$ in $\cat{$\bm{\mathcal M}$-$\bm G$-SSet}$ is called a \emph{$G$-universal weak equivalence}\index{G-universal weak equivalence@$G$-universal weak equivalence|textbf}\index{universal weak equivalence|seealso{$G$-universal weak equivalence}} or \emph{$G$-universal fibration} if $f^\phi$ is a weak homotopy equivalence or Kan fibration, respectively, for each universal $H\subset\mathcal M$ and each homomorphism $\phi\colon H\to G$.

A map of $E\mathcal M$-$G$-simplicial sets is called a \emph{$G$-global weak equivalence}\index{G-global weak equivalence@$G$-global weak equivalence!in EM-SSet@in $\cat{$\bm{E\mathcal M}$-$\bm G$-SSet}$|textbf} or \emph{$G$-global fibration} if it is a $G$-universal weak equivalence or $G$-universal fibration, respectively, when viewed as a map in $\cat{$\bm{\mathcal M}$-$\bm G$-SSet}$.
\end{defi}

\begin{cor}
The $G$-universal weak equivalences and fibrations are part of a unique model structure on $\cat{$\bm{\mathcal M}$-$\bm G$-SSet}$, which we call the \emph{$G$-universal model structure}.\index{G-universal model structure@$G$-universal model structure|textbf}\index{universal model structure|seealso{$G$-universal model structure}} It is simplicial, combinatorial, proper, and filtered colimits in it are homotopical. A possible set of generating cofibrations is given by
\begin{equation*}
\{(\mathcal M\times_\phi G)\times(\del\Delta^n\hookrightarrow\Delta^n) : n\ge 0,H\subset\mathcal M\text{ universal},\phi\colon H\to G\text{ homomorphism}\}
\end{equation*}
and a possible set of generating acyclic cofibrations by
\begin{equation*}
\pushQED{\qed}
\{(\mathcal M\times_\phi G)\times(\Lambda^n_k\hookrightarrow\Delta^n) : 0\le k\le n,H\subset\mathcal M\text{ universal},\phi\colon H\to G\}.\qedhere
\popQED
\end{equation*}
\end{cor}

\begin{cor}\label{cor:EM-G-model-structure}
The $G$-global weak equivalences and fibrations are part of a unique model structure on $\cat{$\bm{E\mathcal M}$-$\bm G$-SSet}$, which we call the \emph{$G$-global model structure}.\index{G-global model structure@$G$-global model structure!on EM-G-SSet@on $\cat{$\bm{E\mathcal M}$-$\bm G$-SSet}$|textbf} It is simplicial, combinatorial, proper, and filtered colimits in it are homotopical. A possible set of generating cofibrations is given by
\begin{equation*}
\{(E\mathcal M\times_\phi G)\times(\del\Delta^n\hookrightarrow\Delta^n) : n\ge 0,H\subset\mathcal M\text{ universal},\phi\colon H\to G\text{ homomorphism}\}
\end{equation*}
and a possible set of generating acyclic cofibrations by
\begin{equation*}
\pushQED{\qed}
\{(E\mathcal M\times_\phi G)\times(\Lambda^n_k\hookrightarrow\Delta^n) : 0\le k\le n,H\subset\mathcal M\text{ universal},\phi\colon H\to G\}.\qedhere
\popQED
\end{equation*}
\end{cor}

\begin{ex}\label{ex:G-globally-contractible}
As a concrete instance of Example~\ref{ex:EM-fixed-points-corep}, we can consider the case where $G$ itself is a (universal) subgroup of $\mathcal M$, and $\phi$ is the identity. In this case, $E\mathcal M\times_\id G\cong E\mathcal M$ with $G$ acting via precomposition. We claim that this is $G$-globally weakly contractible in the sense that the unique map to the terminal object is a $G$-global weak equivalence. Indeed, if $H$ is universal and $\psi\colon H\to G$, then $(E\mathcal M)^\psi\cong E(\mathcal M^\psi)$, so it suffices that $\mathcal M^\psi\not=\varnothing$, i.e.~that there exists an $H$-equivariant injection $\psi^*\omega\to\omega$. This is immediate from universality of $H$.

More generally, the same argument shows that $E\Inj(A,\omega)$ is $G$-globally weakly contractible for any discrete group $G$ and any countable $G$-set $A$.
\end{ex}

Analogously to the global situation, we will now show that the forgetful functor $\cat{$\bm{E\mathcal M}$-$\bm G$-SSet}\to\cat{$\bm{\mathcal M}$-$\bm G$-SSet}$ is fully faithful on associated quasi-categories and characterize its essential image.

\begin{defi}\label{defi:G-semistable}
An $\mathcal M$-$G$-simplicial set $X$ is called \emph{$G$-semistable}\index{G-semistable@$G$-semistable!M-G-simplicial set@$\mathcal M$-$G$-simplicial set|textbf} if the $(H\times G)$-equivariant map $u.\blank\colon X\to X$ is a $\mathcal G_{H,G}$-weak equivalence for any universal $H\subset\mathcal M$ and any $u\in\mathcal M$ centralizing $H$, i.e.~for any homomorphism $\phi\colon H\to G$ the induced map $u.\blank\colon X^\phi\to X^\phi$ is a weak homotopy equivalence.
\end{defi}

\subsubsection{Combinatorics of the $G$-global orbit category}\index{G-global orbit category@$G$-global orbit category|(}
As before, our comparison will proceed indirectly via the respective orbit categories. The mapping spaces in $\cat{O}_{E\mathcal M\times G}$ are again given by fixed points of quotients, and this section is devoted to clarifying their structure. In fact we will do all of this in greater generality now as we will need the additional flexibility later. For this let us fix a group $K$ together with a (finitely or infinitely) countable faithful $K$-set $A$ and a $G$-$K$-biset $X$.

We begin with the following trivial observation:

\begin{lemma}\label{lemma:emg-basic-general}
Let $(u_0,\dots,u_n;x),(v_0,\dots,v_n;y)\in\Inj(A,\omega)^{n+1}\times X$.
\begin{enumerate}
\item These represent the same $n$-simplex of $E\Inj(A,\omega)\times_KX$ (where `${\times_K}\mskip-1.25\thinmuskip$' means that we divide out the diagonal right $K$-action) if and only if there exists some $k\in K$ such that $y=x.k$ and $v_i=u_i.k$ for all $i=0,\dots,n$.\label{item:emgbg-def}
\item Assume these indeed represent the same $n$-simplex and that there is some $i$ such that $u_i=v_i$. Then $x=y$ and $u_j=v_j$ for \emph{all} $j=0,\dots,n$.\label{item:emgbg-unique}
\item Assume $x=y$. Then these represent the same $n$-simplex if and only if there is a $k\in\stabilizer_K(x)$ (where $\stabilizer_K$ denotes the stabilizer)\nomenclature[aStabHx]{$\stabilizer_H(x)$}{stabilizer of $x$ with respect to a given action of $H$} such that $v_i=u_i.k$ for $i=0,\dots,n$.\label{item:emgbg-stab}
\end{enumerate}
\begin{proof}
The first statement holds by definition and the third one is obviously a special case of this. Finally, the second statement follows from the first by freeness of the right $K$-action on $\Inj(A,\omega)$.
\end{proof}
\end{lemma}

We can now characterize the $\phi$-fixed points for any universal $H\subset\mathcal M$ and any homomorphism $\phi\colon H\to G$, generalizing Lemma~\ref{lemma:group-hom-associated}:

\begin{lemma}\label{lemma:emg-fixed-point-characterization-general}
Let $(u_0,\dots,u_n;x)\in\Inj(A,\omega)\times X$ such that $[u_0,\dots,u_n;x]\in (E\Inj(A,\omega)\times_KX)^\phi$. Then there exists for any $h\in H$ a unique $\sigma(h)\in K$ with
\begin{equation}\label{emg-gen:sigma-defining}
hu_i=u_i.\sigma(h)\qquad\text{for all $i=0,\dots,n$},
\end{equation}
and this satisfies
\begin{equation}\label{emg-gen:sigma-vs-phi}
x.\sigma(h)=\phi(h).x.
\end{equation}

The converse holds: whenever there exists a set map $\sigma\colon H\to K$ satisfying $(\ref{emg-gen:sigma-defining})$ and $(\ref{emg-gen:sigma-vs-phi})$ for all $h\in H$, then $[u_0,\dots,u_n;x]$ is a $\phi$-fixed point. Moreover, $\sigma$ is automatically a group homomorphism in this case.
\begin{proof}
That $(u_0,\dots,u_n;x)$ is $\phi$-fixed means by definition that for each $h\in H$
\begin{equation*}
(u_0,\dots,u_n;x)\sim(h,\phi(h)).(u_0,\dots,u_n;x)=(hu_0,\dots,hu_n;\phi(h).x)
\end{equation*}
which again means by definition that there exists a $\sigma(h)\in K$ such that $hu_i=u_i.\sigma(h)$ and moreover $\phi(h).x=x.\sigma(h)$. This $\sigma(h)$ is already uniquely characterized by the first property (for $i=0$) as $K$ acts freely from the right on $\Inj(A,\omega)$, proving the first half of the proposition.

Conversely, if such a $\sigma$ exists, then $[u_0,\dots,u_n;x]$ is clearly $\phi$-fixed. Moreover,
\begin{equation*}
u_0.\sigma(hh')=hh'u_0=hu_0.\sigma(h')=u_0.\sigma(h)\sigma(h'),
\end{equation*}
and hence $\sigma(hh')=\sigma(h)\sigma(h')$ by freeness of the right $K$-action.
\end{proof}
\end{lemma}

In the situation of Lemma~\ref{lemma:emg-fixed-point-characterization-general} we write $\sigma_{(u_0,\dots,u_n)}$\nomenclature[asigmau0un]{$\sigma_{(u_0,\dots,u_n)}$}{homomorphism witnessing that $[u_0,\dots,u_n]$ is a fixed point} for the unique (homomorphism) $H\to K$ satisfying $(\ref{emg-gen:sigma-defining})$. One can show that the lemma provides a complete characterization of the homomorphisms arising this way, which we will only need on the level of vertices:

\begin{cor}\label{cor:emg-equiv-group-hom-realization-general}
Let $\sigma\colon H\to K$ be any group homomorphism. Then there exists $u\in\Inj(A,\omega)$ such that $hu=u.\sigma(h)$ for all $h\in H$. Moreover, if $x\in X$ satisfies $(\ref{emg-gen:sigma-vs-phi})$, then $[u;x]$ is a $\phi$-fixed point of $E\Inj(A,\omega)\times_KX$.
\begin{proof}
When equipped with the tautological $H$-action, $\omega$ is a complete $H$-set universe; on the other hand, $\sigma^*A$ is a countable $H$-set by assumption, so that there exists an $H$-equivariant injection $u\colon\sigma^*A\to\omega$. The $H$-equivariance of $u$ directly translates to $hu=u.\sigma(h)$, and Lemma~\ref{lemma:emg-fixed-point-characterization-general} then proves that $[u;x]$ is $\phi$-fixed.
\end{proof}
\end{cor}

\subsubsection{The comparison} We can now describe the orbit category $\textbf{O}_{E\mathcal M\times G}$ (with respect to the universal graph subgroups) in more concrete terms:
\begin{rk}
The objects of $\textbf{O}_{E\mathcal M\times G}$ are the $E\mathcal M$-$G$-simplicial sets of the form $E\mathcal M\times_\phi G$ where $H\subset\mathcal M$ is universal and $\phi\colon H\to G$ is a homomorphism. If $K\subset\mathcal M$ is another universal subgroup and $\psi\colon K\to G$ a homomorphism, then Lemma~\ref{lemma:emg-fixed-point-characterization-general} tells us that any $n$-simplex of $\Maps_{\textbf{O}_{E\mathcal M\times G}}(E\mathcal M\times_\phi G, E\mathcal M\times_\psi G)$ can be represented by a tuple $(u_0,\dots,u_n;g)\in \mathcal M^{1+n}\times G$ such that there exists a (necessarily unique) group homomorphism $\sigma\colon H\to K$ satisfying
\begin{equation}\label{eq:fixed-point-relation}
hu_i=u_i\sigma(h) \qquad\text{and}\qquad \phi(h)g=g\psi(\sigma(h))
\end{equation}
for all $i=0,\dots,n$ and $h\in H$. Another such tuple $(u_0',\dots,u_n';g')$ represents the same morphism if and only if there exists a $k\in K$ such that $u_i'=u_ik$ for all $i=0,\dots,n$ and $g'=g\psi(k)$.

If $L\subset\mathcal M$ is another universal subgroup, $\theta\colon L\to G$ a group homomorphism, and if $(u_0',\dots,u_n';g')$ represents a morphism $E\mathcal M\times_\psi G\to E\mathcal M\times_\theta G$, then the composition $(u_0',\dots,u_n';g')(u_0,\dots,u_n;g)$ is represented by $(u_0u_0',\dots,u_nu_n';gg')$.

Similarly, objects of $\textbf{O}_{\mathcal M\times G}$ are the $\mathcal M$-$G$-sets $\mathcal M\times_\phi G$ with $\phi$ as above, and maps can be represented by pairs $(u;g)$ with $u\in\mathcal M$ and $g\in G$ satisfying analogous conditions to the above. Compositions are once more given by multiplication in $\mathcal M$ and $G$. In particular, we again get a functor $i\colon  \textbf{O}_{\mathcal M\times G}\to\textbf{O}_{E\mathcal M\times G}$, sending $\mathcal M\times_\phi G$ to $E\mathcal M\times_\phi G$ and a morphism $\mathcal M\times_\phi G\to\mathcal M\times_\psi G$ represented by $(u;g)$ to the morphism $E\mathcal M\times_\phi G\to E\mathcal M\times_\psi G$ represented by the same pair $(u;g)$.
\end{rk}

\begin{defi}\index{G-centralizing morphism@$G$-centralizing morphism|textbf}
A morphism $f\colon\mathcal M\times_\phi G\to \mathcal M\times_\phi G$ in $\textbf{O}_{\mathcal M\times G}$ is called \emph{$G$-centralizing} if there exists a $u\in\mathcal M$ centralizing $H$ such that $f$ is represented by $(u;1)$. Analogously, we define $G$-centralizing morphisms in $\textbf{O}_{E\mathcal M\times G}$.
\end{defi}

\begin{prop}\index{G-global orbit category@$G$-global orbit category!as a simplicial localization}
The above functor $i\colon\textbf{\textup O}_{\mathcal M\times G}\to\textbf{\textup O}_{E\mathcal M\times G}$ is a simplicial localization at the $G$-centralizing morphisms.
\begin{proof}
Let us write $W\subset(\textbf{\textup O}_{E\mathcal M\times G})_0$ for the subcategory of $G$-centralizing morphisms. As in the purely global setting (Proposition~\ref{prop:em-quasi-localization}), this consists of homotopy equivalences and $\cat{O}_{E\mathcal M\times G}$ has fibrant mapping spaces, so that it is enough to prove that for each $n\ge 0$ the homotopical functor
\begin{equation*}
s^*\colon ((\textbf{\textup O}_{E\mathcal M\times G})_0, W)\to ((\textbf{\textup O}_{E\mathcal M\times G})_n, s^*W)
\end{equation*}
induced by the unique map $s\colon[n]\to[0]$ is a homotopy equivalence. As before we have a strict left inverse given by restriction along $i\colon [0]\to[n],0\mapsto 0$ and it suffices to construct a zig-zag of levelwise weak equivalences between $s^*i^*$ and the identity.

For this we recall from the proof of Proposition~\ref{prop:em-quasi-localization} that we can choose for each universal subgroup $H\subset\mathcal M$ injections $\alpha_H,\beta_H\in\mathcal M$ centralizing $H$ and such that $\omega=\im(\alpha_H)\sqcup\im(\beta_H)$.

Now let $H,K\subset\mathcal M$ be universal and let $\phi\colon H\to G$ and $\psi\colon K\to G$ be group homomorphisms. Then any morphism $E\mathcal M\times_\phi G\to E\mathcal M\times_\psi G$ in $(\textbf{\textup O}_{E\mathcal M\times G})_n$ can be represented by a tuple $(u_0,\dots,u_n;g)$ such that there exists a group homomorphism $\sigma\colon H\to K$ satisfying the relations $(\ref{eq:fixed-point-relation})$. We recall from the proof of the ordinary global case that there exists for each $i=0,\dots,n$ a unique $v_i$ such that
\begin{equation*}
v_i\alpha_K=\alpha_Hu_i\qquad\text{and}\qquad v_i\beta_K=\beta_Hu_0.
\end{equation*}
We claim that $(v_0,\dots,v_n;g)$ again defines a morphism, i.e.~its class in $E\mathcal M\times_\psi G$ is $\phi$-fixed. Indeed, we have seen in the global case that $hv_i=v_i\sigma(h)$, hence
\begin{align*}
(h,\phi(h)).(v_0,\dots,v_n;g)&=(hv_0,\dots,hv_n;\phi(h)g)\\ &=(v_0\sigma(h),\dots,v_n\sigma(h);g\psi(\sigma(h)))\sim (v_0,\dots,v_n;g)
\end{align*}
as desired. Similarly, one uses the argument from the non-equivariant case to show that the morphism represented by $(v_0,\dots,v_n;g)$ does not depend on the chosen representative $(u_0,\dots,u_n;g)$.

We now define a functor $f\colon(\textbf{\textup O}_{E\mathcal M\times G})_n\to(\textbf{\textup O}_{E\mathcal M\times G})_n$ as follows: $f$ is the identity on objects and on morphisms given by the above construction. Using that the above is independent of choices one easily checks that $f$ is indeed a functor. As before, we have by construction natural transformations $\id\Leftarrow f\Rightarrow s^*i^*$,
where the left hand transformation is given on $E\mathcal M\times_\phi G$ by $s^*[\alpha_H;1]=[\alpha_H,\dots,\alpha_H;1]$ and the right hand one by $s^*[\beta_H;1]$. As $\alpha_H$ and $\beta_H$ centralize $H$ by definition, these are weak equivalences, finishing the proof.
\end{proof}
\end{prop}

By the same arguments as in the ordinary global setting we deduce:

\begin{thm}\label{thm:em-vs-m-equiv}
The adjunction
\begin{equation*}
E\mathcal M\times_{\mathcal M}^{\textbf{\textup L}}\blank\colon\cat{$\bm{\mathcal M}$-$\bm G$-SSet}_{\textup{$G$-universal}}^\infty\rightleftarrows\cat{$\bm{E\mathcal M}$-$\bm G$-SSet}^\infty :\!\forget^\infty
\end{equation*}
is a Bousfield localization;
in particular, $\forget^\infty$ is fully faithful. Moreover, its essential image consists precisely of the $G$-semistable $\mathcal M$-$G$-simplicial sets.\index{G-global orbit category@$G$-global orbit category|)}\qed
\end{thm}

\subsubsection{Additional model structures}
We can also lift the above comparison of quasi-categories to the level of model categories:

\begin{defi}
A map $f\colon X\to Y$ in $\cat{$\bm{\mathcal M}$-$\bm G$-SSet}$ is called a \emph{$G$-global weak equivalence} if $E\mathcal M\times_{\mathcal M}^{\textbf{\textup L}}f$ is an isomorphism in $\Ho(\cat{$\bm{E\mathcal M}$-$\bm G$-SSet})$.\index{G-global weak equivalence@$G$-global weak equivalence!in M-G-SSet@in $\cat{$\bm{\mathcal M}$-$\bm G$-SSet}$|textbf}
\end{defi}

\begin{cor}\label{cor:em-vs-m-equiv-model-cat}\index{G-global model structure@$G$-global model structure!on M-G-SSet@on $\cat{$\bm{\mathcal M}$-$\bm G$-SSet}$|textbf}
The $G$-universal cofibrations and $G$-global weak equivalences are part of a unique model structure on $\cat{$\bm{\mathcal M}$-$\bm G$-SSet}$, which we call the \emph{$G$-global model structure}. An object $X\in\cat{$\bm{\mathcal M}$-$\bm G$-SSet}$ is fibrant in this model structure if and only if it is fibrant in the $G$-universal model structure and moreover $G$-semistable in the sense of Definition~\ref{defi:G-semistable}.

Moreover, this model structure is combinatorial with generating cofibrations
\begin{equation*}
\{(\mathcal M\times_\phi G)\times(\del\Delta^n\hookrightarrow\Delta^n) : n\ge 0, H\subset\mathcal M\text{ universal},\phi\colon H\to G\text{ homomorphism}\},
\end{equation*}
simplicial, left proper, and filtered colimits in it are homotopical.

Finally, the simplicial adjunction
\begin{equation*}
E\mathcal M\times_\mathcal M\blank\colon \cat{$\bm{\mathcal M}$-$\bm G$-SSet}_{\textup{$G$-global}} \rightleftarrows \cat{$\bm{E\mathcal M}$-$\bm G$-SSet} :\!\forget
\end{equation*}
is a Quillen equivalence with homotopical right adjoint.
\begin{proof}
Theorem~\ref{thm:em-vs-m-equiv} allows us to invoke Lurie's localization criterion (Theorem~\ref{thm:lurie-localization-criterion}) which proves all of the above claims except for the part about filtered colimits, which is in turn an instance of Lemma~\ref{lemma:filtered-still-homotopical}.
\end{proof}
\end{cor}

As a special case of Proposition~\ref{prop:equivariant-injective-model-structure}, the $G$-global weak equivalences of $E\mathcal M$-$G$-simplicial sets are part of an injective model structure.\index{injective G-global model structure@injective $G$-global model structure!on EM-G-SSet@on $\cat{$\bm{E\mathcal M}$-$\bm G$-SSet}$} We will now prove the analogue of this for $\cat{$\bm{\mathcal M}$-$\bm G$-SSet}$:

\begin{thm}\label{thm:G-M-injective-semistable-model-structure}\index{G-global model structure@$G$-global model structure!injective|seeonly{injective $G$-global model structure}}\index{injective G-global model structure@injective $G$-global model structure!on M-G-SSet@on $\cat{$\bm{\mathcal M}$-$\bm G$-SSet}$|textbf}
There exists a unique model structure on $\cat{$\bm{\mathcal M}$-$\bm G$-SSet}$ whose cofibrations are the underlying cofibrations and whose weak equivalences are the $G$-global weak equivalences. We call this the \emph{injective $G$-global model structure}. It is combinatorial, simplicial, left proper, and filtered colimits in it are homotopical.

Finally, the simplicial adjunction
\begin{equation*}
\forget\colon\cat{$\bm{E\mathcal M}$-$\bm G$-SSet}_{\textup{inj.~$G$-global}}\rightleftarrows\cat{$\bm{\mathcal M}$-$\bm G$-SSet}_{\textup{inj.~$G$-global}} :\!\Maps^{\mathcal M}(E\mathcal M,\blank)
\end{equation*}
is a Quillen equivalence with homotopical left adjoint.
\begin{proof}
We first claim that the $G$-global weak equivalences are stable under pushout along \emph{injective} cofibrations. For this we let $f\colon A\to B$ be a $G$-global weak equivalence and $i\colon A\to C$ an injective cofibration. Applying the factorization axiom in the $G$-global model structure, we can factor $f=pk$ where $k$ is an acyclic cofibration and $p$ is an acyclic fibration. As the $G$-global and $G$-universal model structures on $\cat{$\bm{\mathcal M}$-$\bm G$-SSet}$ have the same cofibrations, they also have the same acyclic fibrations; in particular, $p$ is a $G$-universal weak equivalence.

We now consider the iterated pushout
\begin{equation*}
\begin{tikzcd}
A \arrow[rr, bend left=25pt, "f"]\arrow[r, "k"']\arrow[d, "i"'] & X\arrow[r, "p\vphantom{k}"']\arrow[d, "j"] & B\arrow[d]\\
C \arrow[r, "\ell"'] & Y \arrow[r, "q\vphantom{k}"'] & D.
\end{tikzcd}
\end{equation*}
Then $\ell$ is an acyclic cofibration in the $G$-global model structure as a pushout of an acyclic cofibration. Moreover, $j$ is an injective cofibration as a pushout of an injective cofibration, so $q$ is a $G$-universal (and hence in particular $G$-global) weak equivalence by left properness of the equivariant injective model structure. The claim follows as $q\ell$ is a pushout of $f$ along $i$.

We therefore conclude from Corollary~\ref{cor:mix-model-structures} that the model structure exists and that it is combinatorial and left proper. Moreover, Lemma~\ref{lemma:filtered-still-homotopical} shows that filtered colimits in it are still homotopical, so it only remains to verify the Pushout Product Axiom for the simplicial tensoring.

For this we can argue precisely as in the proof of Proposition~\ref{prop:equivariant-injective-model-structure} once we show that for any simplicial set $K$ the functor $K\times\blank$ preserves $G$-global weak equivalences, and that for any $\mathcal M$-$G$-simplicial set $X$ the functor $\blank\times X$ sends weak equivalences of simplicial sets to $G$-global weak equivalences.

For the second statement it is actually clear that $\blank\times X$ sends weak equivalences even to $G$-universal weak equivalences. For the first statement it is similarly clear that $K\times\blank$ preserves $G$-universal weak equivalences, but it also preserves acyclic cofibrations in the usual $G$-global model structure as the latter is simplicial. The claim again follows as any $G$-global weak equivalence can be factored as a $G$-global acyclic cofibration followed by a $G$-universal weak equivalence.
\end{proof}
\end{thm}

\begin{rk}
With a bit of (combinatorial) work one can show that $\forget\dashv\Maps^{\mathcal M}(E\mathcal M,\blank)$ is also a Quillen equivalence for the usual $G$-global model structures. On the other hand, it is not clear whether $E\mathcal M\times_{\mathcal M}\blank$ preserves injective cofibrations or $G$-global weak equivalences in general, and in particular whether $E\mathcal M\times_{\mathcal M}\blank\dashv\forget$ is also a Quillen equivalence for the injective $G$-global model structures.
\end{rk}

\subsection{An explicit \texorpdfstring{$\bm G$}{G}-semistable replacement}\label{sec:g-semistable-replacement}\index{G-semistable@$G$-semistable!M-G-simplicial set@$\mathcal M$-$G$-simplicial set!replacement|(}
If one wants to check if a morphism $f\colon X\to Y$ in $\cat{$\bm{\mathcal M}$-$\bm G$-SSet}$ is a $G$-global weak equivalence straight from the definition, one runs into trouble as soon as at least one of $X$ or $Y$ is not cofibrant because computing $E\mathcal M\times_{\mathcal M}^{\cat{L}}f$ then involves cofibrant replacements, and the ones provided by the small object argument are completely intractable.

On the other hand, the \emph{$G$-universal} weak equivalences are much easier to understand, so one could instead try to take $G$-semistable replacements of $X$ and $Y$ and then apply the following general observation about Bousfield localizations:

\begin{lemma}\label{lemma:semistable-between-semistable}
A morphism $f\colon X\to Y$ of $G$-semistable $\mathcal M$-$G$-simplicial sets is a $G$-global weak equivalence if and only if it is a $G$-universal weak equivalence.\qed
\end{lemma}

However, this leaves us with the problem of finding (functorial) $G$-semistable replacements. While the $G$-global model structure asserts that these exist, they are again completely inexplicit. We could also try to construct them by means of $E\mathcal M\times_{\mathcal M}^{\textbf L}\blank$, but then we would of course be back where we started.

In this subsection we will remedy this situation by constructing explicit $G$-semistable replacements, which will become crucial later (see e.g.~Theorem~\ref{thm:tame-M-sset-vs-EM-sset}) both to identify $G$-global weak equivalences as well as to compute the left derived functor ${E\mathcal M}\times_{\mathcal M}^\cat{L}\blank$.

\begin{rk}
Before we come to the construction, let us think about how it should look like. The simplicial set $E\mathcal M$ is (weakly) contractible, and one can conclude from this, see e.g.~\cite[Proposition~4.2.4.4]{htt}, that $\cat{$\bm{E\mathcal M}$-SSet}$ models ordinary non-equivariant homotopy theory when equipped with the \emph{underlying} weak equivalences. More precisely, with respect to these weak equivalences, the homotopical functors
\begin{equation*}
\const\colon\cat{SSet}\to\cat{$\bm{E\mathcal M}$-SSet}\qquad\text{and}\qquad
\forget\colon\cat{$\bm{E\mathcal M}$-SSet}\to\cat{SSet}
\end{equation*}
induce mutually quasi-inverse equivalences on associated quasi-categories. It follows that the composition
\begin{equation*}
\cat{$\bm{\mathcal M}$-SSet}^\infty_{\textup{universal}}\xrightarrow{E\mathcal M\times_{\mathcal M}^{\textbf L}\blank}\cat{$\bm{E\mathcal M}$-SSet}^\infty_{\textup{global}}\xrightarrow{\forget}\cat{SSet}^\infty
\end{equation*}
is equivalent to taking homotopy colimits over $\mathcal M$, also cf.~\cite[Definition~3.2]{I-vs-M-1-cat}.
\end{rk}

\subsubsection{Action categories}
The remark suggests that we might be lucky and succeed in constructing the replacement by means of a suitable equivariant enhancement of one of the standard constructions of homotopy colimits. This will indeed work for the model of what is usually called the \emph{action groupoid} (although it won't be a groupoid in our case), which we now recall:

\begin{constr}\nomenclature[aM]{$(\blank)\hq\mathcal M$}{$G$-global action category/bar construction}
Let $X$ be any $\mathcal M$-set. We write $X\hq\mathcal M$ for the \emph{action category} (also known as the \emph{translation category}),\index{action category|textbf} i.e.~the category with set of objects $X$ and for any $x\in X$ and $u\in\mathcal M$ a morphism $u\colon x\to u.x$; we emphasize that this means that if $u\not=v$ with $u.x=v.x$, then $u$ and $v$ define two distinct morphisms $x\to u.x=v.x$. The composition in $X\hq\mathcal M$ is given by multiplication in $\mathcal M$.
\end{constr}

The $\mathcal M$-action on $X$ immediately gives an $\mathcal M$-action on $\Ob(X\hq\mathcal M)$; however, it is not entirely clear how to extend this to morphisms. For an invertible element $\alpha\in\core\mathcal M$, the condition that $\alpha.f$ for $f\colon x\to y$ should be a morphism $\alpha.x\to \alpha.y$ naturally leads to the guess $\alpha.f=\alpha f\alpha^{-1}$. While general elements of $\mathcal M$ are not invertible, there is still a notion of \emph{conjugation}. This is made precise by the following construction, which is implicit in \cite[proof of Lemma~5.2]{schwede-semistable} (which Schwede attributes to Strickland) and also appeared in an earlier version of \cite{schwede-k-theory}:

\begin{constr}\label{constr:M-conjugation}\index{conjugation on M@conjugation on $\mathcal M$|textbf}
Let $\alpha\in\mathcal M$. We define for any $u\in\mathcal M$ the \emph{conjugation} $c_\alpha(u)$\nomenclature[acalpha]{$c_\alpha$}{conjugation by $\alpha\in\mathcal M$} of $u$ by $\alpha$ via
\begin{equation*}
c_\alpha(u)(x)=\begin{cases}
\alpha u(y) & \text{if }x=\alpha(y)\\
x & \text{if }x\notin\im \alpha.
\end{cases}
\end{equation*}
We remark that this is well-defined (as $\alpha$ is injective), and one easily checks that this is again injective, so that we get a map $c_\alpha\colon\mathcal M\to\mathcal M$.
\end{constr}

Put differently, $c_\alpha(u)$ is the unique element of $\mathcal M$ such that
\begin{equation*}
c_\alpha(u)\alpha=\alpha u\qquad\text{and}\qquad c_\alpha(u)(x)=x\text{ for $x\notin\im \alpha$}.
\end{equation*}
The first condition justifies the name `conjugation,' and it means in particular that $c_\alpha(u)=\alpha u\alpha^{-1}$ for invertible $\alpha$.

For a group $G$, conjugation by a fixed element $g$ defines an endomorphism of $G$, and letting $g$ vary this yields an action of $G$ on itself. The analogous statement holds for the above construction:

\begin{lemma}\label{lemma:conjugation-homomorphism}
For any $\alpha\in\mathcal M$, the map $c_\alpha\colon\mathcal M\to\mathcal M$ is a monoid homomorphism. Moreover, for varying $\alpha$ this defines an action of $\mathcal M$ on itself, i.e.~$c_1=\id_{\mathcal M}$ and $c_\alpha\circ c_\beta=c_{\alpha\beta}$ for all $\alpha,\beta\in\mathcal M$.
\begin{proof}
We will only prove the first statement, the calculations for the other claims being similar. For this let $u,v\in\mathcal M$. Then
\begin{equation*}
c_\alpha(uv)\alpha=\alpha uv= c_\alpha(u)\alpha v=c_\alpha(u)c_\alpha(v)\alpha.
\end{equation*}
On the other hand, if $x\notin\im(\alpha)$, then $c_\alpha(uv)(x)=x$ and $c_\alpha(v)(x)=x\notin\im\alpha$, hence also $(c_\alpha(u)c_\alpha(v))(x)=x$. We conclude that $c_\alpha(uv)=c_\alpha(u)c_\alpha(v)$. Moreover, clearly $c_\alpha(1)=1$, so that $c_\alpha$ is indeed a monoid homomorphism.
\end{proof}
\end{lemma}

\begin{constr}
We define an $\mathcal M$-action on $X\hq\mathcal M$ as follows: the action of $\mathcal M$ on objects is the one on $X$, and on morphisms $\alpha\in\mathcal M$ acts by sending
\begin{equation*}
x\xrightarrow{u}u.x\qquad\text{to}\qquad \alpha.x\xrightarrow{c_\alpha(u)}\alpha.(u.x);
\end{equation*}
note that this is indeed a morphism as $c_\alpha(u).(\alpha.x)=(c_\alpha(u)\alpha).x=(\alpha u).x=\alpha.(u.x)$ by construction of $c_\alpha$. By the previous lemma, this is then an endofunctor of $X\hq\mathcal M$, and for varying $\alpha$ this yields an $\mathcal M$-action.

If $f\colon X\to Y$ is a map of $\mathcal M$-sets, then we write $f\hq\mathcal M$ for the functor that is given on objects by $f$ and that sends a morphism
\begin{equation*}
x\xrightarrow{u}u.x\qquad\text{to}\qquad f(x)\xrightarrow{u}u.f(x)=f(u.x).
\end{equation*}
One easily checks that this is well-defined, functorial in $f$, and that $f\hq\mathcal M$ is $\mathcal M$-equivariant. Postcomposing with the nerve we therefore get a functor $\cat{$\bm{\mathcal M}$-Set}\to\cat{$\bm{\mathcal M}$-SSet}$ that we denote by $(\blank)\hq\mathcal M$ again. If $G$ is any group, we moreover get an induced functor $\cat{$\bm{\mathcal M}$-$\bm G$-Set}\to\cat{$\bm{\mathcal M}$-$\bm G$-SSet}$ by pulling through the $G$-action (again denoted by the same symbol).

If $X$ is any $\mathcal M$-set, then there is a unique functor from $X$ (viewed as a discrete category) to $X\hq\mathcal M$ that is the identity on objects. This then yields a natural transformation $\pi\colon\discr\Rightarrow(\blank)\hq\mathcal M$ from the functor that sends an $\mathcal M$-$G$-set $X$ to the discrete simplicial set $X$ with the induced action.
\end{constr}

\begin{constr}
Let $X$ be any $\mathcal M$-$G$-simplicial set. Applying the above construction levelwise yields a bisimplicial set $X\#\mathcal M$ with $(X\#\mathcal M)_{n,\bullet}=X_n\hq\mathcal M$, and this receives a map from the bisimplicial set $\Discr X$ with $(\Discr X)_{n,\bullet}=\discr X_n$. This yields a functor $\cat{$\bm{\mathcal M}$-$\bm G$-SSet}\to\cat{$\bm{\mathcal M}$-$\bm G$-BiSSet}$\nomenclature[aBiSSet]{$\cat{BiSSet}$}{category of bisimplicial sets} receiving a natural transformation $\Pi\colon\Discr\Rightarrow(\blank)\#\mathcal M$.

Taking diagonals, we then obtain a functor $\cat{$\bm{\mathcal M}$-$\bm G$-SSet}\to\cat{$\bm{\mathcal M}$-$\bm G$-SSet}$, that we again denote by $(\blank)\hq\mathcal M$, together with a natural transformation $\id\Rightarrow(\blank)\hq\mathcal M$, that we again denote by $\pi$. We remark that on $\mathcal M$-$G$-sets (viewed as discrete simplicial sets) this recovers the previous construction.
\end{constr}

\subsubsection{The detection result}
Now we are ready to state the main result of this subsection:

\begin{thm}\label{thm:hq-M-semistable-replacement}
\begin{enumerate}
\item For any $X\in\cat{$\bm{\mathcal M}$-$\bm G$-SSet}$, $X\hq\mathcal M$ is $G$-semistable and the map $\pi_X\colon X\to X\hq\mathcal M$ is a $G$-global weak equivalence.
\item For any $f\colon X\to Y$ in $\cat{$\bm{\mathcal M}$-$\bm G$-SSet}$ the following are equivalent:
\begin{enumerate}
\item $f$ is a $G$-global weak equivalence.
\item $f\hq\mathcal M$ is a $G$-global weak equivalence.
\item $f\hq\mathcal M$ is a $G$-universal weak equivalence.
\end{enumerate}
\end{enumerate}
\end{thm}

The proof of the theorem will occupy the rest of this subsection.

\begin{rk}
We defined the bisimplicial set $X\#\mathcal M$ in terms of the action category construction. We can also look at this bisimplicial set `from the other side,' which recovers the bar construction (just as in the usual construction of non-equivariant homotopy quotients):

If $Y$ is any $\mathcal M$-$G$-set, then $(Y\hq\mathcal M)_m$ consists by definition of the $m$-chains
\begin{equation*}
y\xrightarrow{u_1}u_1.y\xrightarrow{u_2}(u_2u_1).y\xrightarrow{u_3}\cdots\xrightarrow{u_m}(u_m\cdots u_1).y
\end{equation*}
of morphisms in $Y\hq\mathcal M$. Such a chain is obviously uniquely described by the source $y\in Y$ together with the injections $u_m,\dots, u_2,u_1\in\mathcal M$, which yields a bijection $(Y\hq\mathcal M)_m\cong \mathcal M^m\times Y$. This bijection becomes $\mathcal M$-$G$-equivariant, when we let $G$ act via its action on $X$ and $\mathcal M$ via its action on $X$ and the conjugation action on each of the $\mathcal M$-factors.

The assignment $Y\mapsto\mathcal M^m\times Y$ becomes a functor in $Y$ in the obvious way, and with respect to this the above bijection is clearly natural in $\mathcal M$-equivariant maps. Applying this levelwise, we therefore get a natural isomorphism
\begin{equation}\label{eq:x-hash-m-other-levels}
(X\#\mathcal M)_{\bullet,m}\cong\mathcal M^m\times X
\end{equation}
of $\mathcal M$-$G$-simplicial sets. While we will not need this below, we remark that unravelling the definitions, one can work out that under the isomorphism $(\ref{eq:x-hash-m-other-levels})$ the simplicial structure maps of $X\#\mathcal M$ indeed correspond to those of the usual bar construction.
\end{rk}

By construction and the previous remark, we understand the bisimplicial set $X\#\mathcal M$ in both its simplicial directions individually. In non-equivariant simplicial homotopy theory, the \emph{Diagonal Lemma} then often allows to leverage this to prove statements about the diagonal ($X\hq\mathcal M$ in our case). Luckily, this immediately carries over to our situation:

\begin{lemma}\label{lemma:equivariant-diagonal}\index{Diagonal Lemma!equivariant|textbf}\index{Diagonal Lemma}
Let $M$ be a monoid, let $\mathcal F$ be a collection of subgroups of $M$, and let $f\colon X\to Y$ be a map of $M$-bisimplicial sets. Assume that for each $n\ge 0$ the map $f_{n,\bullet}\colon X_{n,\bullet}\to Y_{n,\bullet}$ is an $\mathcal F$-weak equivalence, or that for each $m\ge 0$ the map $f_{\bullet,m}\colon X_{\bullet,m}\to Y_{\bullet,m}$ is. Then $\diag f\colon\diag X\to\diag Y$ is an $\mathcal F$-weak equivalence.
\begin{proof}
By symmetry it suffices to consider the first case. If $H\in\mathcal F$ is any subgroup, then $(f_{n,\bullet})^H$ literally agrees with $(f^H)_{n,\bullet}$ (if we take the usual construction of fixed points), and likewise $(\diag f)^H=\diag(f^H)$. Thus, the claim follows immediately from the usual Diagonal Lemma, see e.g.~\cite[Proposition~IV.1.7]{goerss}.
\end{proof}
\end{lemma}

Let us draw some non-trivial consequences from this:

\begin{cor}\label{cor:X-hq-M-homotopical}
The above functor $(\blank)\hq\mathcal M$ preserves $G$-universal weak equivalences of $\mathcal M$-$G$-simplicial sets.
\begin{proof}
By the previous lemma (applied to the universal graph subgroups of the monoid $\mathcal M\times G$), it suffices that each $(\blank\#\mathcal M)_{\bullet,m}$ does, which follows from the isomorphism~$(\ref{eq:x-hash-m-other-levels})$ together with Lemma~\ref{lemma:equivariant-weak-equivalences-prod-homotopical}.
\end{proof}
\end{cor}

\begin{lemma}\label{lemma:X-hq-M-semistable}
Let $X$ be any $\mathcal M$-$G$-simplicial set. Then $X\hq\mathcal M$ is $G$-semistable.
\begin{proof}
Let $H\subset\mathcal M$ be any subgroup and let $\alpha\in\mathcal M$ centralize $H$. We will show that $\alpha.\blank\colon X\hq\mathcal M\to X\hq\mathcal M$ is even an $(H\times G)$-weak equivalence.

By Lemma~\ref{lemma:equivariant-diagonal} we reduce this to the case that $X$ is an \emph{$\mathcal M$-$G$-set}, in which case $X\hq\mathcal M$ is the nerve of the action category. We then observe that the maps $\alpha\colon x\to \alpha.x$ assemble into a natural transformation $a\colon\id\Rightarrow(\alpha.\blank)$ by virtue of the identity $c_\alpha(u)\alpha=\alpha u$. This natural transformation is obviously $G$-equivariant, but it is also $H$-equivariant: if $h\in H$ and $x\in X$ are arbitrary, then $h.a_x$ is by definition the map $h\alpha h^{-1}\colon h.x\to h.(\alpha.x)$; as $h$ commutes with $\alpha$, this agrees with $a_{h.x}$ as desired. Upon taking nerves, $a$ therefore induces an $(H\times G)$-equvariant homotopy between the identity and $\alpha.\blank$; in particular, $\alpha.\blank$ is an $(H\times G)$-weak equivalence.
\end{proof}
\end{lemma}

\begin{lemma}\label{lemma:hq-M-homotopy-po}
The functor $(\blank)\hq\mathcal M\colon\cat{$\bm{\mathcal M}$-$\bm G$-SSet}\to\cat{$\bm{\mathcal M}$-$\bm G$-SSet}$ is cocontinuous and it preserves injective cofibrations.
\begin{proof}
As colimits in $\cat{$\bm{\mathcal M}$-$\bm G$-SSet}$ are defined levelwise and by definition of the injective cofibrations, it is enough to show that the functor $\cat{$\bm{\mathcal M}$-$\bm G$-SSet}\to\cat{Set}, X\mapsto (X\hq\mathcal M)_n$ is cocontinuous and sends levelwise injections to injections for all $n\ge0$. Plugging in the definition, this is the functor $X\mapsto (X\#\mathcal M)_{n,n}$, which by $(\ref{eq:x-hash-m-other-levels})$ is isomorphic to $X\mapsto X_n\times\mathcal M^n$. The claim follows as the latter clearly has the required properties.
\end{proof}
\end{lemma}

\begin{lemma}\label{lemma:hq-M-tensors}
The functor $(\blank)\hq\mathcal M\colon\cat{$\bm{\mathcal M}$-$\bm G$-SSet}\to\cat{$\bm{\mathcal M}$-$\bm G$-SSet}$ has a natural simplicial enrichment with respect to which it preserves tensors. The natural transformation $\pi$ is simplicially enriched.
\begin{proof}
If $X\in\cat{$\bm{\mathcal M}$-$\bm G$-SSet}$, then the $\mathcal M$-$G$-set $(X\hq\mathcal M)_n$ by construction only depends on the $n$-simplices $X_n$, and similarly for maps in $\cat{$\bm{\mathcal M}$-$\bm G$-SSet}$. It follows formally that for any simplicial set $K$ the maps
\begin{equation}\label{eq:comparison-hq-M-tensor}
K_n\times (X\hq\mathcal M)_n \to ((K_n\times X)\hq\mathcal M)_n = ((K\times X)\hq\mathcal M)_n
\end{equation}
induced by the inclusions $X\to K_n\times X, x \mapsto (k,x)$ for varying $k\in K_n$ assemble into a map $K\times (X\hq\mathcal M)\to (K\times X)\hq\mathcal M$ of $\mathcal M$-$G$-simplicial sets, and that this makes $(\blank)\hq\mathcal M$ lax tensored (and hence enriched) over $\cat{SSet}$. Moreover, as $(\blank)\hq\mathcal M$ preserves coproducts by the previous lemma, the maps $(\ref{eq:comparison-hq-M-tensor})$ are actually isomorphisms, i.e.~with respect to this enrichment $(\blank)\hq\mathcal M$ preserves tensors.

Finally, to see that $\pi$ is simplicially enriched, we note that again by construction the induced map on $n$-simplices $(\pi_X)_n$ only depends on the $n$-simplices $X_n$, i.e.~the right hand square in the diagram
\begin{equation*}
\begin{tikzcd}
K_n\times X_n\arrow[r, equal]\arrow[d, "K_n\times(\pi_X)_n"'] & (K_n\times X)_n\arrow[r, equal]\arrow[d, "(\pi_{K_n\times X})_n"] & (K\times X)_n\arrow[d, "(\pi_{K\times X})_n"]\\
K_n\times (X\hq\mathcal M)_n\arrow[r, "\cong"'] & ((K_n\times X)\hq\mathcal M)_n\arrow[r, equal] & ((K\times X)\hq\mathcal M)_n
\end{tikzcd}
\end{equation*}
commutes for every $n\ge0$. Since also the left hand square commutes by naturality and the universal property of coproducts, the claim follows.
\end{proof}
\end{lemma}

The only remaining input that we need to prove the theorem is an explicit computation of $(\blank)\hq\mathcal M$ on the $\mathcal M$-$G$-sets $\mathcal M\times_\phi G$ for universal $H\subset\mathcal M$ and any group homomorphism $\phi\colon H\to G$. We will actually do this in greater generality, which will become important later in the proof of Theorem~\ref{thm:tame-M-sset-vs-EM-sset}:

\begin{thm}\label{thm:hq-M-computation}
Let $H\subset\mathcal M$ be universal, let $A$ be a countable faithful $H$-set, and let $X$ be any $G$-$H$-biset. If we consider $\Inj(A,\omega)$ as an $\mathcal M$-$H$-biset in the obvious way, then the unique functor
$\Inj(A,\omega)\hq\mathcal M \to E\Inj(A,\omega)$ that is the identity on objects induces a $G$-universal weak equivalence
\begin{equation}\label{eq:hq-m-comparison}
\big(\Inj(A,\omega)\times_HX\big)\hq\mathcal M\cong
\big(\Inj(A,\omega)\hq\mathcal M\big)\times_HX
\xrightarrow{\sim} \big(E\Inj(A,\omega)\big)\times_HX.
\end{equation}
Here the unlabelled isomorphism comes from the cocontinuity of $(\blank)\hq\mathcal M$.
\end{thm}

\begin{rk}
The theorem of course also holds for $A$ uncountable (as then both sides are just empty). However, some lemmas we will appeal to only hold for countable $A$ and this is also the only case we will be interested in. Therefore, we have decided to state the theorem in the above form.
\end{rk}

The proof of the theorem needs some preparations and will be given at the end of this subsection. Before that, let us already use it to deduce the detection result. This will  involve a standard `cell induction'\index{cell induction|textbf} argument; as similar lines of reasoning will appear again later, we formalize this part once and for all:

\begin{lemma}\label{lemma:saturated-objects}
Let $\mathscr C$ be a cocomplete category, let $I$ be a class of morphisms, and let $\mathscr S$ be a class of objects in $\mathscr C$ such that the following conditions are satisfied:
\begin{enumerate}
\item $\mathscr S$ contains the initial object $\varnothing$\label{item:so-empty}.
\item If
\begin{equation}\label{diag:pushout-generating-cof}
\begin{tikzcd}
A\arrow[d]\arrow[r, "i"] & B\arrow[d]\\
C \arrow[r] & D
\end{tikzcd}
\end{equation}
is a pushout with $C\in\mathscr S$ and $i\in I$, then $D\in\mathscr S$.\label{item:so-pushout-closure}
\item $\mathscr S$ is closed under filtered colimits.\label{item:so-closure-filtered}
\end{enumerate}
Then $\mathscr S$ contains all $I$-cell complexes. If in addition also
\begin{enumerate}
\setcounter{enumi}{3}
\item $\mathscr S$ is closed under retracts,\label{item:so-retract-closure}
\end{enumerate}
then $\mathscr S$ contains all retracts of $I$-cell complexes. In particular, if $\mathscr C$ is a cofibrantly generated model category and $I$ a set of generating cofibrations, then $\mathscr S$ contains all cofibrant objects.
\begin{proof}
It is enough to prove the first statement. For this let $X$ be an $I$-cell complex, i.e.~there exists an ordinal $\alpha$ and a functor $X_\bullet\colon\{\beta<\alpha\}\to\mathscr C$ with colimit $X$ such that the following conditions hold:
\begin{enumerate}
\item[(A)] $X_0=\varnothing$
\item[(B)] For each $\beta$ with $\beta+1<\alpha$, the map $X_\beta\to X_{\beta+1}$ can be written as a pushout of some $i\in I$.
\item[(C)] If $\beta<\alpha$ is a limit ordinal, then the maps $X_\gamma\to X_\beta$ for $\gamma<\beta$ exhibit $X_\beta$ as a (filtered) colimit.
\end{enumerate}
If $\alpha$ is a limit ordinal, we extend this to $\{\beta\le\alpha\}$ via $X_\alpha\mathrel{:=}\colim_{\beta<\alpha}X_\beta=X$ together with the obvious structure maps; if $\alpha$ is a successor ordinal instead, we replace $\alpha$ by its predecessor $\beta$ (and $X_\beta$ by the isomorphic object $X$). In both cases we obtain a functor $X_\bullet\colon\{\beta\le\alpha\}\to\mathscr C$ satisfying conditions (A)--(C) above (with `$<$' replaced by `$\le$') and such that $X_\alpha=X$. We will prove by transfinite induction that $X_\beta\in\mathscr S$ for all $\beta\le\alpha$ which will then imply the claim.

By Conditions $(\ref{item:so-empty})$ and (A) we see that $X_0=\varnothing\in\mathscr S$. Now assume $\beta>0$, and we know the claim for all $\gamma<\beta$. If $\beta$ is a successor ordinal, $\beta=\gamma+1$, then Conditions $(\ref{item:so-pushout-closure})$ and (B) together with the induction hypothesis imply that $X_\beta\in\mathscr S$. On the other hand, if $\beta$ is a limit ordinal, then the maps $X_\gamma\to X_\beta$ for $\gamma<\beta$ express $X_\beta$ as a filtered colimit of elements of $\mathscr S$ by Condition (C) together with the induction hypothesis. Thus, Condition~$(\ref{item:so-closure-filtered})$ immediately implies the claim, finishing the proof.
\end{proof}
\end{lemma}

\begin{cor}\label{cor:saturated-trafo}
Let $\mathscr C$ be a cocomplete category and let $\mathscr D$ be a left proper model category such that filtered colimits in it are homotopical. Let $I$ be any collection of morphisms in $\mathscr C$, and let $F,G\colon\mathscr C\to\mathscr D$ be functors together with a natural transformation $\tau\colon F\Rightarrow G$. Assume the following:
\begin{enumerate}
\item $\tau_\varnothing$ is a weak equivalence and for every map $(X\to Y)\in I$ both $\tau_X$ and $\tau_Y$ are weak equivalences. \label{item:st-base}
\item Both $F$ and $G$ send pushouts along maps $i\in I$ to homotopy pushouts in $\mathscr D$.\label{item:st-pushout}
\item $F$ and $G$ preserve filtered colimits up to weak equivalence in the sense that for any small filtered category $J$ and any diagram $X\colon J\to\mathscr C$ the canonical maps $\colim_J F\circ X\to F(\colim_JX)$ and $\colim_J G\circ X\to G(\colim_JX)$ are weak equivalences.\label{item:st-filtered}
\end{enumerate}
Then $\tau$ is a weak equivalence on all retracts of $I$-cell complexes. In particular, if $\mathscr C$ is a cofibrantly generated model category and $I$ is a set of generating cofibrations, then $\tau_X$ is a weak equivalence for all cofibrant $X$.
\begin{proof}
Let $\mathscr S$ be the class of objects $X\in\mathscr C$ such that $\tau_X$ is a weak equivalence. It will suffice to verify the conditions of the previous lemma. Indeed, Condition $(\ref{item:so-empty})$ of the lemma is in immediate consequence of our Condition~$(\ref{item:st-base})$.

In order to verify Condition~$(\ref{item:so-pushout-closure})$ of the lemma, we consider a pushout as in $(\ref{diag:pushout-generating-cof})$. By naturality, this induces a commutative cube
\begin{equation*}
\begin{tikzcd}[row sep=small, column sep=small]
& GA\arrow[dd] \arrow[rr, "Gi"] && GB\arrow[dd]\\
FA\arrow[dd]\arrow[ur] \arrow[rr, "Fi" near end, crossing over] && FB\arrow[ur]\\
& GC\arrow[rr] && GD\\
FC\arrow[ur]\arrow[rr] && FD\arrow[from=uu, crossing over]\arrow[ur]
\end{tikzcd}
\end{equation*}
with all the diagonal maps coming from $\tau$. The assumption $C\in\mathscr S$ implies that the front-to-back map at the lower left corner is a weak equivalence, and our Condition $(\ref{item:st-base})$ tells us that the two top front-to-back maps are weak equivalences. Moreover, both front and back square are homotopy pushouts by Condition~$(\ref{item:st-pushout})$. We conclude that also the lower right front-to-back map is a weak equivalence, just as desired.

To verify Condition~$(\ref{item:so-closure-filtered})$ of the lemma, we observe that for any small filtered category $J$ and any $X_\bullet\colon J\to\mathscr C$ the map $\tau_{\colim_{j\in J} X_j}$
fits into a commutative diagram
\begin{equation*}
\begin{tikzcd}[column sep=huge]
\colim_{j\in J} F(X_j)\arrow[d]\arrow[r, "{\colim_{j\in J}\tau_{X_j}}"] & \colim_{j\in J} G(X_j)\arrow[d]\\
F(\colim_{j\in J} X_j) \arrow[r, "{\tau_{\colim_{j\in J} X_j}}"'] & G\big(\colim_{j\in J} X_j\big)
\end{tikzcd}
\end{equation*}
where the vertical maps are the canonical comparison maps and hence weak equivalences by our Condition~($\ref{item:st-filtered}$). If now all $X_i$ lie in $\mathscr S$, then the top horizontal map is a filtered colimit of weak equivalences and hence a weak equivalence by assumption, proving Condition $(\ref{item:so-closure-filtered})$ of the lemma.

Finally, Condition~($\ref{item:so-retract-closure}$) for $\mathscr S$ is automatic as weak equivalences in any model category are closed under retracts. This finishes the proof.
\end{proof}
\end{cor}

\begin{proof}[Proof of Theorem~\ref{thm:hq-M-semistable-replacement}]
We will only prove the first statement; the second one will then follow formally from this (cf.~Lemma~\ref{lemma:semistable-between-semistable}).

We already know by Lemma~\ref{lemma:X-hq-M-semistable} that $X\hq\mathcal M$ is $G$-semistable for any $X\in\cat{$\bm{\mathcal M}$-$\bm G$-SSet}$. It therefore only remains to prove that $\pi_X\colon X\to X\hq\mathcal M$ is a $G$-global weak equivalence. Since we have seen in Corollary~\ref{cor:X-hq-M-homotopical} that $(\blank)\hq\mathcal M$ preserves $G$-universal weak equivalences, it suffices to prove this for cofibrant $X$, and Lemma~\ref{lemma:hq-M-homotopy-po} together with Corollary~\ref{cor:saturated-trafo} reduces this further to the case that $X$ is the source or target of one of the standard generating cofibrations, i.e.~$X=(\mathcal M\times_\phi G)\times\del\Delta^n$ or $X=(\mathcal M\times_\phi G)\times\Delta^n$ for some $n\ge 0$. By Lemma~\ref{lemma:hq-M-tensors} we are then further reduced to $X=\mathcal M\times_\phi G$.

In this case we observe that the unit $\eta$ of the adjunction $E\mathcal M\times_{\mathcal M}(\blank)\dashv\forget$, evaluated at $\mathcal M\times_\phi G$ can be factored as
\begin{equation*}
\mathcal M\times_\phi G\xrightarrow{\pi} (\mathcal M\times_\phi G)\hq\mathcal M\xrightarrow{(\ref{eq:hq-m-comparison})}E\mathcal M\times_\phi G\cong\forget\big(E\mathcal M\times_{\mathcal M}(\mathcal M\times_\phi G)\big),
\end{equation*}
where the middle arrow comes from applying Theorem~\ref{thm:hq-M-computation} for $A=\omega$ and $X=G$ with its left regular $G$-action and $H$ acting from the right via $\phi$, and the final isomorphism uses that the left adjoint $E\mathcal M\times_{\mathcal M}\blank$ is cocontinuous. The above composition is a $G$-global weak equivalence by Corollary~\ref{cor:em-vs-m-equiv-model-cat}, and the middle arrow is even a $G$-universal weak equivalence by Theorem~\ref{thm:hq-M-computation}. Therefore, $2$-out-of-$3$ implies that $\pi$ is a $G$-global weak equivalence as desired. By the above reduction, this completes the proof.\index{G-semistable@$G$-semistable!M-G-simplicial set@$\mathcal M$-$G$-simplicial set!replacement|)}
\end{proof}

It remains to prove that the map $(\ref{eq:hq-m-comparison})$ is a $G$-universal weak equivalence.

\begin{ex}
As the proof of the theorem will be quite technical, let us begin with something much easier that will nevertheless give an idea of the general argument: we will show that the map in question is a non-equivariant weak equivalence for $A$ finite and $H=1$, which amounts to saying that $\Inj(A,\omega)\hq\mathcal M$ is weakly contractible as a simplicial set. This argument also appears as \cite[Example~3.3]{I-vs-M-1-cat}.

By construction, the category $\Inj(A,\omega)\hq\mathcal M$ has objects the injections $i\colon A\to\omega$ and it has for every $u\in\mathcal M$ and $i\in\Inj(A,\omega)$ a map $u\colon i\to ui$. If now $i,j\in\Inj(A,\omega)$ are arbitrary, then we can pick a bijection $u\in\mathcal M$ with $ui=j$, i.e.~$i$ and $j$ are isomorphic in $\Inj(A,\omega)\hq\mathcal M$. It follows that for our favourite $i\in\Inj(A,\omega)$ the induced functor $B\End(i)\to \Inj(A,\omega)\hq\mathcal M$ is an equivalence of categories. Therefore it suffices that the monoid $\End(i)$ has weakly contractible classifying space. However, $\End(i)$ consists precisely of those $u\in\mathcal M$ that fix $\im(i)$ pointwise. Picking a bijection $\omega\cong \omega\setminus\im(i)$ therefore yields an isomorphism of monoids $\mathcal M\cong\End(i)$. Since $\mathcal M$ has weakly contractible classifying space by \cite[Lemma~5.2]{schwede-semistable} (whose proof Schwede attributes to Strickland), this finishes the proof.
\end{ex}

\subsubsection{Equivariant analysis of action categories}
Fix a universal subgroup $K\subset\mathcal M$ and a homomorphism $\phi\colon K\to G$. By definition, the simplicial set $(\Inj(A,\omega)\times_HX)\hq\mathcal M$ is the nerve of the category of the same name, and as the nerve is a right adjoint, this is compatible with $\phi$-fixed points. In the following, we want to understand these fixed point categories better and in particular describe them as disjoint unions of monoids in analogy with the above example.

However, in Theorem~\ref{thm:hq-M-computation} we allow $A$ to be infinite (and $A=\omega$ is the case we actually used in the proof of Theorem~\ref{thm:hq-M-semistable-replacement}). In this case, there are `too many isomorphism classes' in $(\Inj(A,\omega)\hq\mathcal M)\times_HX$: for example, not all objects in $\mathcal M\hq\mathcal M$ are isomorphic, though they all receive a map from $1\in\mathcal M$. To salvage this situation we introduce the full subcategory $C_K\subset(\Inj(A,\omega)\times_HX)\hq\mathcal M$ spanned by those $[u,x]$ for which $\im(u)^c\subset\omega$ contains a complete $K$-set universe (this is independent of the choice of representative as $\im(u)=\im(u.h)$ for all $h\in H$).

\begin{lemma}\label{lemma:c-K-whe}
The inclusion $C_K\hookrightarrow(\Inj(A,\omega)\hq\mathcal M)\times_HX$ induces a $(K\times G)$-homotopy equivalence on nerves.
\begin{proof}
Let $\alpha\in\mathcal M$ be $K$-equivariant such that $\im(\alpha)^c$ contains a complete $K$-set universe. Then $\alpha.\blank$ is $(K\times G)$-equivariant and it takes all of $(\Inj(A,\omega)\hq\mathcal M)\times_HX$ to $C_K$. We claim that this is a $(K\times G)$-homotopy inverse to the inclusion. Indeed, the proof of Lemma~\ref{lemma:X-hq-M-semistable} shows that the maps $\alpha\colon x\to\alpha.x$ define a natural transformation from the identity to $\alpha.\blank$ as endofunctors of $(\Inj(A,\omega)\hq\mathcal M)\times_HX$, showing that $\alpha.\blank$ is a right $(K\times G)$-homotopy inverse. However, as $C_K$ is a full subcategory, these also define such a transformation for the other composite, proving that $\alpha.\blank$ is also a left $(K\times G)$-homotopy inverse.
\end{proof}
\end{lemma}

The next lemma in particular tells us that $C_K$ avoids the aforementioned issue: 

\begin{lemma}\label{lemma:map-implies-isomorphic}
Let $p,q\in C_K^\phi$ and fix a representative $(u,x)\in\Inj(A,\omega)\times X$ of $p$. Then the following are equivalent:
\begin{enumerate}
\item There exists a map $f\colon p\to q$ in $C_K^\phi$.
\item There exists a representative of $q$ of the form $(v,x)$ such that in addition $\sigma_u=\sigma_v$ (see the discussion after Lemma~\ref{lemma:emg-fixed-point-characterization-general}).
\item There exists an isomorphism $f'\colon p\to q$ in $C_K^\phi$.
\end{enumerate}
\end{lemma}

For the proof we need the following notation:

\begin{defi}
Let $A,B$ be sets, let $A=A_1\sqcup A_2$ be any partition, and let $f_i\colon A_i\to B$ ($i=1,2$) be any maps of sets. Then we write $f_1+f_2$ for the unique map $A\to B$ that agrees on $A_1$ with $f_1$ and on $A_2$ with $f_2$.
\end{defi}

By slight abuse, we will also apply the above when $f_1$ and $f_2$ are maps into subsets of $B$. Obviously, $f_1+f_2$ will be injective if and only if $f_1$ and $f_2$ are injections with disjoint image.

\begin{proof}[Proof of Lemma~\ref{lemma:map-implies-isomorphic}]
Obviously $(3)\Rightarrow(1)$; we will prove that also $(1)\Rightarrow(2)$ and $(2)\Rightarrow(3)$.

Assume $f\colon p\to q$ is any morphism in $C_K^\phi$. Then $q=f.p$, so that $(fu,x)$ is a representative of $q$. We claim that this has the desired property, i.e.~$v\mathrel{:=}fu$ satisfies $\sigma_v=\sigma_u$. But indeed, as $f$ is $\phi$-fixed, it has to be $K$-equivariant by definition of the action. Then
\begin{equation*}
v.\sigma_v(k)=kv=kfu=fku=fu.\sigma_u(k)=v.\sigma_u(k),
\end{equation*}
and hence $\sigma_v(k)=\sigma_u(k)$ as $H$ acts freely on $\Inj(A,\omega)$. This proves $(1)\Rightarrow(2)$.

For the proof of $(2)\Rightarrow(3)$, we observe that $\im(u)$ and $\im(v)$ are both $K$-subsets of $\omega$: indeed, $k.u(a)=(ku)(a)=u(\sigma_u(k).a)$ for all $a\in A$, and similarly for $v$. Thus, $\omega\setminus\im(u)$ and $\omega\setminus\im(v)$ are $K$-sets in their own right; as they are obviously countable and moreover contain a complete $K$-set universe each by definition of $C_K$, they are themselves complete $K$-set universes. It follows that there exists a $K$-equivariant bijection $f'_1\colon\omega\setminus\im(u)\cong\omega\setminus\im(v)$.

On the other hand, $f'_0=vu^{-1}\colon\im(u)\to\im(v)$ is also $K$-equivariant because
\begin{equation*}
f'_0(k.u(a))=f'_0(u(\sigma_u(k).a))=v(\sigma_u(k).a)=v(\sigma_v(k).a)=k.v(x)=k.f'_0(u(a))
\end{equation*}
for all $k\in K$, $a\in A$. As it is moreover obviously bijective, we conclude that $f'\mathrel{:=}f'_0+f'_1$ defines an isomorphism $[u,x]\to[v,x]$ in $C_K^\phi$ as desired.
\end{proof}

We fix for each isomorphism class of $C_K^\phi$ a representative $p\in C_K^\phi$ and for each such $p$ in turn a representative $(u,x)\in\Inj(A,\omega)\times X$. Let us write $\mathscr I$ for the set of all these. The above lemma then implies:

\begin{cor}\label{cor:c-K-phi-decomposition}
The functor $\Phi\colon\coprod_{(u,x)\in\mathscr I}B\End_{C_K^\phi}([u,x])\to C_K^\phi$ induced by the evident embeddings is an equivalence of categories.\qed
\end{cor}

Let us now study these monoids more closely:

\begin{prop}\label{prop:tau-group-completion}
Let $(u,x)\in\Inj(A,\omega)\times X$ such that $[u,x]\in C_K^\phi$.
\begin{enumerate}
\item Let $f\colon[u,x]\to[u,x]$ be any map in $C_K^\phi$. Then there exists a unique $\tau(f)\in H$ such that $fu=u.\tau(f)$. Moreover, $\tau(f)$ centralizes $\im\sigma_u$ and stabilizes $x\in X$. (We caution the reader that $\tau$ depends on the chosen representative.)
\item The assignment $f\mapsto\tau(f)$ defines a monoid homomorphism
\begin{equation*}
\tau\colon\End_{C_K^\phi}([u,x])\to\centralizer_H(\im\sigma_u)\cap\stabilizer_H(x).
\end{equation*}
\item The homomorphism $\tau$ induces a weak equivalence on classifying spaces.
\end{enumerate}
\end{prop}

For the proof of the proposition we will need:

\begin{lemma}\label{lemma:M-K-classifying}
The monoid $\mathcal M^K$ of $K$-equivariant injections $\omega\to\omega$ has weakly contractible classifying space.
\begin{proof}
This is an equivariant version of~\cite[proof of Lemma~5.2]{schwede-semistable}. As $K$ is universal, we can pick a $K$-equivariant bijection $\omega\cong\omega\amalg\omega$ which yields two $K$-equivariant injections $\alpha,\beta\in\mathcal M$ whose images partition $\omega$. Then the conjugation homomorphism $c_\alpha\colon\mathcal M\to\mathcal M$ satisfies
\begin{equation*}
c_\alpha(u)\alpha=\alpha u\qquad\text{and}\qquad c_\alpha(u)\beta=\beta
\end{equation*}
for all $u\in\mathcal M$. The first equality proves that $\alpha$ defines a natural transformation from the identity to $B(c_\alpha)\colon B\mathcal M\to B\mathcal M$ (also cf.~the proof of Lemma~\ref{lemma:X-hq-M-semistable}) while the second one shows that $\beta$ defines a natural transformation from the constant functor to it. Upon taking nerves, we therefore get a zig-zag of homotopies between the identity and a constant map, proving the claim.
\end{proof}
\end{lemma}

\begin{proof}[Proof of Proposition~\ref{prop:tau-group-completion}]
For the first statement we observe that there exists at most one such $\tau(f)$ by freeness of the action. On the other hand, $f$ being a morphism in particular means that $f.[u,x]=[u,x]$. Plugging in the definition of the action and of the equivalence relation we divided out, this means that there exists $\tau(f)\in H$ with $(fu,x)=(u.\tau(f),x.\tau(f))$. Thus, it only remains to show that $\tau(f)$ centralizes $\im(\sigma_u)$. Indeed, if $k\in K$ is arbitrary, then on the one hand
\begin{equation*}
kfu=fku=fu.\sigma_u(k)=u.(\tau(f)\sigma_u(k))
\end{equation*}
(where we have used that $f$ is $K$-equivariant since it is a morphism in $C_K^\phi$), and on the other hand
\begin{equation*}
kfu=ku.\tau(f)=u.(\sigma_u(k)\tau(f)).
\end{equation*}
Thus, $u.(\tau(f)\sigma_u(k))=u.(\sigma_u(k)\tau(f))$, whence indeed $\tau(f)\sigma_u(k)=\sigma_u(k)\tau(f)$ by the freeness of the right $H$-action. This finishes the proof of $(1)$.

For the second statement, we observe that $ff'u=fu.\tau(f')=u.(\tau(f)\tau(f'))$ and hence $\tau(ff')=\tau(f)\tau(f')$ by the above uniqueness statement. Analogously, $1u=u=u.1$ shows $\tau(1)=1$.

For the final statement, we will first prove:

\begin{claim*} The assignment
\begin{align*}
\textup{T}\colon\End_{C_K^\phi}([u,x]) &\to \Inj(\omega\setminus\im u,\omega\setminus\im u)^K\times\big(\centralizer_H(\im\sigma_u)\cap\stabilizer_H(x)\big)\\
f &\mapsto \big(f|_{\omega\setminus\im(u)}\colon\omega\setminus\im(u)\to\omega\setminus\im(u), \tau(f)\big)
\end{align*}
defines an isomorphism of monoids.
\begin{proof}
This is well-defined by the first part and since the injection $f$ restricts to a self-bijection of $\im(u)$ by the above, so that it also has to preserve $\omega\setminus\im(u)$. It is then obvious (using the second part of the proposition) that $\textup{T}$ is a monoid homomorphism.

We now claim that it is actually an isomorphism of monoids. For injectivity, it suffices to observe that if $\tau(f)=\tau(f')$, then $fu=u\tau(f)=u\tau(f')=f'u$.

For surjectivity, we let $f_1\colon\omega\setminus\im(u)\to\omega\setminus\im(u)$ be any $K$-equivariant injection and we let $t\in\centralizer_H(\im\sigma_u)\cap\stabilizer_H(x)$. Then there is a unique map $f_0\colon\im(u)\to\im(u)$ with $f_0u=u.t$ and this is automatically injective (in fact, even bijective). We claim that it is also $K$-equivariant. Indeed, if $k\in K$ is arbitrary, then
\begin{equation*}
k.f_0(u(x))=ku(t.x)=u(\sigma_u(k)t.x)=u(t\sigma_u(k).x)=f_0(u(\sigma_u(k).x))=f_0(k.u(x))
\end{equation*}
where we have used that $t$ commutes with $\sigma_u(k)$. Thus, $f\mathrel{:=}f_0+f_1$ defines a $K$-equivariant injection $\omega\to\omega$ and we claim that this is an endomorphism of $[u,x]$ in $C_K$, hence the desired preimage of $(f_1,t)$. Indeed, $f.[u,x]$ is represented by
\begin{equation*}
(fu,x)=(u.t,x)\sim (u,x.t^{-1})=(u,x),
\end{equation*}
where we have used that $t$ and hence also $t^{-1}$ stabilizes $x$.
\end{proof}
\end{claim*}

To prove that $\tau$ induces a weak homotopy equivalence on classifying spaces, we now observe that the induced map factors as
\begin{align*}
\nerve\big(B\End([u,x])\big)&\xrightarrow{\nerve(B\textup{T})}\nerve\big(B(\Inj(\omega\setminus\im u,\omega\setminus\im u)^K\times L)\big)\\
&\cong\nerve\big(B(\Inj(\omega\setminus\im u,\omega\setminus\im u)^K)\big)\times\nerve(BL)\xrightarrow{\pr}\nerve(BL),
\end{align*}
where $L\mathrel{:=}\centralizer_H(\im\sigma_u)\cap\stabilizer_H(x)$. The first map is an isomorphism by the previous statement, so it suffices that $\Inj(\omega\setminus\im u,\omega\setminus\im u)^K$ has trivial classifying space. But indeed, as in the proof of Lemma~\ref{lemma:map-implies-isomorphic} we see that $\omega\setminus\im u$ is a complete $K$-set universe, so that $\Inj(\omega\setminus\im u,\omega\setminus\im u)^K\cong\mathcal M^K$ as monoids. Thus, the claim follows from Lemma~\ref{lemma:M-K-classifying}.
\end{proof}

\subsubsection{Equivariant analysis of quotient categories}
We recall that the simplicial set $E\Inj(A,\omega)\times X$ is canonically and equivariantly isomorphic to the nerve of the groupoid of the same name. On the other hand, the right $H$-action on $E\Inj(A,\omega)\times X$ is free, and as the nerve preserves quotients by \emph{free} group actions, the simplicial set $E\Inj(A,\omega)\times_H X$ is again canonically identified with the nerve of the corresponding quotient of categories, which we denote by the same name.

In the following we want to devise a description of this category and its fixed points analogous to the above results. For this we first observe that while hom sets in quotient categories are in general hard to describe, the situation is easier here because this particular quotient is preserved by the nerve: namely, any morphism $p\to q$ (for $p,q\in \Inj(A,\omega)\times_KX$) is represented by a triple $(u,v;x)$ with $u,v\in\Inj(A,\omega)$, $x\in X$ such that $[u;x]=q$ and $[v,x]=p$; moreover, a triple $(u',v';x')$ represents the same morphism if and only if there exists an $h\in H$ with $u'=u.h$, $v'=v.h$, and $x'=x.h$. The following lemma gives a more concrete description once we have fixed representatives of $p$ and $q$:

\begin{lemma}\label{lemma:emg-equiv-hom-sets-general}
Let $(u,x),(v,y)\in\Inj(A,\omega)\times X$. Then the assignment
\begin{align*}
\{h\in H: y.h=x\} &\to\Hom_{E\Inj(A,\omega)\times_HX}([u;x],[v;y])\\
h&\mapsto [v.h,u;x]=[v,u.h^{-1},y]
\end{align*}
is well-defined and bijective. In particular, the assignment
\begin{equation}\label{eq:emgehs-endo}
\begin{aligned}
\stabilizer_H(x)&\to\End_{E\Inj(A,\omega)\times_HX}([u;x])\\
h&\mapsto [u.h,u;x]
\end{aligned}
\end{equation}
is bijective; in fact, this even defines an isomorphism of groups.
\begin{proof}
Let us denote the above map by $\alpha$. We first observe that for any $h$ on the left hand side the representative $(v.h,u;x)$ differs from $(v,u.h^{-1},y)$ only by right multiplication by $h$, so that the two given definitions of $\alpha(h)$ indeed agree. In particular, they define an edge from $[u;x]$ to $[v;y]$, proving that $\alpha$ is well-defined.

Lemma~\ref{lemma:emg-basic-general}-$(\ref{item:emgbg-unique})$ already implies that $\alpha$ is injective. For surjectivity, we pick an edge on the right hand side and let $(a,b;c)$ be a representative. By definition $[b;c]=[u;x]$, so after acting suitably from the right by $H$ on $(a,b;c)$, we may assume $b=u$ and $c=x$, i.e.~our chosen representative takes the form $(a,u;x)$. But by assumption this represents an edge to $[v;y]$ and hence $[a;x]=[v;y]$, i.e.~there exists an $h\in H$ such that $a=v.h$ and $x=y.h$, which then provides the desired preimage.

Specializing to $(v;y)=(u;x)$ shows that $(\ref{eq:emgehs-endo})$ is bijective. The calculation
\begin{equation*}
[u.hh',u;x]=[u.hh',u.h';x][u.h',u;x]=[u.h,u;x][u.h',u;x]
\end{equation*}
(where the final equality uses that $x.(h')^{-1}=x$) then shows that it is in fact an isomorphism of groups.
\end{proof}
\end{lemma}

Similarly, we can describe the hom sets in the fixed point categories:

\begin{lemma}\label{lemma:emg-equiv-hom-sets-general-phi}
Let $(u,x),(v,y)\in\Inj(A,\omega)\times X$ such that $[u;x],[v;y]$ are $\phi$-fixed points of $E\Inj(A,\omega)\times_HX$. Then
\begin{align*}
\{h\in H: y.h=x, h\sigma_u(k)h^{-1}=\sigma_v(k)\text{ $\forall k\in K$}\} &\to
\Hom_{(E\Inj(A,\omega)\times_HX)^\phi}([u;x],[v;y])\\
h&\mapsto [v.h,u;x]=[v,u.h^{-1},y]
\end{align*}
is well-defined and bijective. In particular, this yields a bijection
\begin{align*}
\stabilizer_H(x)\cap\centralizer_H(\im\sigma_u)&\to\End_{(E\Inj(A,\omega)\times_HX)^\phi}([u;x])\\
h&\mapsto [u.h,u;x].
\end{align*}
This map is in fact even an isomorphism of groups.
\begin{proof}
By the previous lemma we are reduced to proving that $[v.h,u;x]$ is $\phi$-fixed if and only if $\sigma_v(k)=h\sigma_u(k)h^{-1}$ for all $k\in K$. Indeed,
\begin{align*}
(k,\phi(k)).(v.h,u;x)&=(kv.h,k.u;\phi(k).x)=(v.\sigma_v(k)h,u.\sigma_u(k);x.\sigma_u(k))\\
&\sim\big(v.(\sigma_v(k)h\sigma_u(k)^{-1}),u;x\big)
\end{align*}
where we have applied Lemma~\ref{lemma:emg-fixed-point-characterization-general} to $(u;x)$. By freeness of the right $H$-action this represents the same element as $(v.h,u;x)$ if and only if $\sigma_v(k)h\sigma_u(k)^{-1}=h$, which is obviously equivalent to the above condition.
\end{proof}
\end{lemma}

\begin{prop}\label{prop:EInj-decomposition}
The functor
\begin{equation*}
\Psi\colon\coprod_{(u,x)\in\mathscr I}B\big(\centralizer_H(\im\sigma_u)\cap\stabilizer_H(x)\big)\to \big(E\Inj(A,\omega)\times_HX\big)^\phi
\end{equation*}
given on the $(u,x)$-summand by sending $t\in\centralizer_H(\im\sigma_u)\cap\stabilizer_H(x)$ to the morphism $[u.t,u;x]$ is an equivalence of groupoids.
\begin{proof}
Lemma \ref{lemma:emg-equiv-hom-sets-general-phi} implies that this indeed lands in the $\phi$-fixed points and that each of the maps
\begin{equation*}
\centralizer_H(\im\sigma_u)\cap\stabilizer_H(x)\to\End_{\big(E\Inj(A,\omega)\times_HX\big)^\phi}([u;x])
\end{equation*}
is a group isomorphism; in particular, $\Psi$ is a functor. As both source and target of $\Psi$ are groupoids, it suffices then to show that $\mathscr I$ also is a system of (representatives of) representatives of the isomorphism classes on the right hand side.

To see that $\mathscr I$ hits every isomorphism class at most once, assume $(u,x),(v,y)\in\mathscr I$ represent isomorphic elements in $(E\Inj(A,\omega)\times_HX)^\phi$. Lemma~\ref{lemma:emg-equiv-hom-sets-general-phi} implies that there exists an $h\in H$ with $y.h=x$ and $h^{-1}\sigma_v(k)h=\sigma_u(k)$ for all $k\in K$. Then $(w,x)\mathrel{:=}(v.h,y.h)$ represents the same element as $(v,y)$ in both $C_K^\phi$ as well as $(E\Inj(A,\omega)\times_HX)^\phi$, and one easily checks that $\sigma_w=\sigma_u$. Thus, $[u,x]\cong[w,x]=[v,y]$ in $C_K^\phi$ by Lemma~\ref{lemma:map-implies-isomorphic}, and hence $(u,x)=(v,y)$ by definition of $\mathscr I$.

But $\mathscr I$ also covers the isomorphism classes of $(E\Inj(A,\omega)\times_HX)^\phi$: if $(v;y)$ represents an arbitrary element of it and $\alpha\in\mathcal M$ is again $K$-equivariant with $\omega\setminus\im\alpha$ a complete $K$-set universe, then one easily checks that $[\alpha v,v;y]$ is $\phi$-fixed, so that it witnesses $[v;y]\cong[\alpha v;y]$. In other words, we may assume that $v$ misses a complete $K$-set universe. But in this case it defines an element of $C_K$, which is then by definition isomorphic in $C_K$ to some $[u,x]$ with $(u,x)\in\mathscr I$, and hence also in $(E\Inj(A,\omega)\times_HX)^\phi$ by functoriality.
\end{proof}
\end{prop}

\subsubsection{Quotient vs.~action categories} Putting all of the above results together, we finally get:

\begin{proof}[Proof of Theorem~\ref{thm:hq-M-computation}]
As before let $K\subset\mathcal M$ be a universal subgroup, and let $\phi\colon K\to G$ be any group homomorphism. We have to show that $(\ref{eq:hq-m-comparison})$ induces a weak homotopy equivalence on $\phi$-fixed points.

To this end, we consider the diagram of categories and functors
\begin{equation*}
\begin{tikzcd}[column sep=-15pt]
\smash{\coprod\limits_{(u,x)\in\mathscr I}}B\End_{C_K^\phi}([u,x]) \arrow[rr, "\coprod B\tau"]\arrow[d, "\Phi"'] && \smash{\coprod\limits_{(u,x)\in\mathscr I}}B\big(\centralizer_H(\im\sigma_u)\cap\stabilizer_H(x)\big)\arrow[d, "\Psi"]\\
C_K^\phi\arrow[r, hook] & \big((\Inj(A,\omega)\times_HX)\hq\mathcal M)^\phi\arrow[r] & (E\Inj(A,\omega)\times_HX)^\phi
\end{tikzcd}
\end{equation*}
where the unlabelled arrow is induced by the map in question (viewed as a functor). The top path through this diagram sends $f\colon[u,x]\to[u,x]$ to $[u.\tau(f),u;x]$ while the lower one sends it to $[fu,u;x]$. As $u.\tau(f)=fu$ by definition of $\tau$, these two agree, i.e.~the above diagram commutes.

We now observe that the top map is a weak homotopy equivalence by Proposition~\ref{prop:tau-group-completion}, that the vertical maps are equivalences by Corollary~\ref{cor:c-K-phi-decomposition} and Proposition~\ref{prop:EInj-decomposition}, respectively, and that the lower left inclusion is a weak homotopy equivalence by Lemma~\ref{lemma:c-K-whe}. The claim now follows by $2$-out-of-$3$.
\end{proof}

\subsection{Functoriality}
\index{functoriality in homomorphisms!for EM-G-SSet@for $\cat{$\bm{E\mathcal M}$-$\bm G$-SSet}$|(}
We will now study various change-of-group functors for the above models of $G$-global homotopy theory. We begin with the versions for $E\mathcal M$-actions, where Lemma~\ref{lemma:alpha-shriek-projective} and Lemma~\ref{lemma:alpha-star-injective}, respectively, specialize to:

\begin{cor}
Let $\alpha\colon H\to G$ be any group homomorphism. Then
\begin{equation*}
\alpha_!\colon\cat{$\bm{E\mathcal M}$-$\bm{H}$-SSet}_{\textup{$H$-global}}\rightleftarrows\cat{$\bm{E\mathcal M}$-$\bm{G}$-SSet}_{\textup{$G$-global}} :\!\alpha^*
\end{equation*}
is a simplicial Quillen adjunction with homotopical right adjoint.\qed
\end{cor}

\begin{cor}
Let $\alpha\colon H\to G$ be any group homomorphism. Then
\begin{equation*}
\alpha^*\colon\cat{$\bm{E\mathcal M}$-$\bm{G}$-SSet}_{\textup{$G$-global injective}}\rightleftarrows\cat{$\bm{E\mathcal M}$-$\bm{H}$-SSet}_{\textup{$H$-global injective}} :\!\alpha_*
\end{equation*}
is a simplicial Quillen adjunction.\qed
\end{cor}

On the other hand, Propositions~\ref{prop:alpha-shriek-injective} and \ref{prop:alpha-lower-star-homotopical} imply:

\begin{cor}\label{cor:alpha-shriek-injective-EM}
Let $\alpha\colon H\to G$ be an \emph{injective} group homomorphism. Then
\begin{equation*}
\alpha_!\colon\cat{$\bm{E\mathcal M}$-$\bm H$-SSet}_{\textup{$H$-global injective}}\rightleftarrows\cat{$\bm{E\mathcal M}$-$\bm G$-SSet}_{\textup{$H$-global injective}} :\!\alpha^*
\end{equation*}
is a simplicial Quillen adjunction. In particular, $\alpha_!$ is homotopical.\qed
\end{cor}

\begin{cor}\label{cor:alpha-lower-star-injective-EM}
Let $\alpha\colon H\to G$ be an \emph{injective} group homomorphism. Then
\begin{equation*}
\alpha^*\colon\cat{$\bm{E\mathcal M}$-$\bm G$-SSet}_{\textup{$G$-global}}\rightleftarrows \cat{$\bm{E\mathcal M}$-$\bm H$-SSet}_{\textup{$H$-global}}:\!\alpha_*
\end{equation*}
is a simplicial Quillen adjunction. If $(G:\im\alpha)<\infty$, then $\alpha_*$ is homotopical.\qed
\end{cor}

Finally, Proposition~\ref{prop:free-quotient-general} specializes to:

\begin{cor}\label{cor:free-quotient-EM}
Let $\alpha\colon H\to G$ be any group homomorphism. Then
\begin{equation*}
\alpha_!\colon\cat{$\bm{E\mathcal M}$-$\bm H$-SSet}\to\cat{$\bm{E\mathcal M}$-$\bm G$-SSet}
\end{equation*}
preserves weak equivalences between objects with free $\ker(\alpha)$-action.\index{functoriality in homomorphisms!for EM-G-SSet@for $\cat{$\bm{E\mathcal M}$-$\bm G$-SSet}$|)}\qed
\end{cor}

The case of $\mathcal M$-actions needs slightly more work:
\index{functoriality in homomorphisms!for M-G-SSet@for $\cat{$\bm{\mathcal M}$-$\bm G$-SSet}$|(}
\begin{cor}\label{cor:alpha-shriek-projective-M}
Let $\alpha\colon H\to G$ be any group homomorphism. Then
\begin{equation*}
\alpha_!\colon\cat{$\bm{\mathcal M}$-$\bm{H}$-SSet}_{\textup{$H$-global}}\rightleftarrows\cat{$\bm{\mathcal M}$-$\bm{G}$-SSet}_{\textup{$G$-global}} :\!\alpha^*
\end{equation*}
is a simplicial Quillen adjunction.
\begin{proof}
This follows as before for the $H$-universal and $G$-universal model structure, respectively. Thus, it suffices by Proposition~\ref{prop:cofibrations-fibrant-qa} that $\alpha^*$ sends $G$-semistable objects to $H$-semistable ones, which is obvious from the definition.
\end{proof}
\end{cor}

\begin{prop}\label{prop:alpha-star-EM-homotopical}
Let $\alpha\colon H\to G$ be any group homomorphism. Then
\begin{equation*}
\alpha^*\colon\cat{$\bm{\mathcal M}$-$\bm{G}$-SSet}_{\textup{$G$-global injective}}\rightleftarrows\cat{$\bm{\mathcal M}$-$\bm{H}$-SSet}_{\textup{$H$-global injective}} :\!\alpha_*
\end{equation*}
is a simplicial Quillen adjunction. In particular, $\alpha^*$ is homotopical.
\begin{proof}
It is clear that $\alpha^*$ preserves injective cofibrations, so it only remains to show that it is homotopical. However, while $\alpha^*$ commutes with $E\mathcal M\times_{\mathcal M}\blank$, it is not clear a priori that it is also suitably compatible with $E\mathcal M\times_{\mathcal M}^{\textbf{L}}\blank$ since $\alpha^*$ usually does not preserve cofibrant objects of the universal model structures.

Instead, we consider the commutative diagram
\begin{equation*}
\begin{tikzcd}
\cat{$\bm{\mathcal M}$-$\bm G$-SSet}_{\textup{$G$-global}}\arrow[r, "\blank\hq\mathcal M"]\arrow[d, "\alpha^*"'] & \cat{$\bm{\mathcal M}$-$\bm G$-SSet}_{\textup{$G$-universal}}\arrow[d,"\alpha^*"]\\
\cat{$\bm{\mathcal M}$-$\bm H$-SSet}_{\textup{$H$-global}}\arrow[r, "\blank\hq\mathcal M"'] & \cat{$\bm{\mathcal M}$-$\bm H$-SSet}_{\textup{$H$-universal}}.
\end{tikzcd}
\end{equation*}
Then the horizontal arrows create weak equivalences by Theorem~\ref{thm:hq-M-semistable-replacement}, while the right vertical arrow is clearly homotopical. Thus, also the left hand vertical arrow is homotopical.
\end{proof}
\end{prop}

\begin{cor}
Let $\alpha\colon H\to G$ be an \emph{injective} group homomorphism. Then
\begin{equation*}
\alpha_!\colon\cat{$\bm{\mathcal M}$-$\bm H$-SSet}_{\textup{$H$-global injective}}\rightleftarrows\cat{$\bm{\mathcal M}$-$\bm G$-SSet}_{\textup{$H$-global injective}} :\!\alpha^*
\end{equation*}
is a simplicial Quillen adjunction. In particular, $\alpha_!$ is homotopical.
\begin{proof}
By Proposition~\ref{prop:alpha-shriek-injective}, $\alpha_!$ preserves injective cofibrations and it sends $H$-universal weak equivalences to $G$-universal ones. On the other hand, any $H$-global weak equivalence factors as an $H$-global acyclic cofibration followed by an $H$-universal weak equivalence. Since $\alpha_!$ sends the former to $G$-global weak equivalences by Corollary~\ref{cor:alpha-shriek-projective-M}, $2$-out-of-$3$ implies that $\alpha_!$ is also homotopical.\index{functoriality in homomorphisms!for M-G-SSet@for $\cat{$\bm{\mathcal M}$-$\bm G$-SSet}$|)}
\end{proof}
\end{cor}

\subsection{\texorpdfstring{$\bm{G}$}{G}-global homotopy theory vs.~\texorpdfstring{$\bm G$}{G}-equivariant homotopy theory}\label{sec:g-global-vs-g-em}\index{proper G-equivariant homotopy theory@proper $G$-equivariant homotopy theory!vs G-global homotopy theory@vs.~$G$-global homotopy theory|(}
We continue to let $G$ denote a discrete group, possibly infinite. As promised, we will now explain how to exhibit classical proper $G$-equivariant homotopy theory (i.e.~equivariant homotopy theory with respect to the collection of \emph{finite} subgroups of $G$) as a Bousfield localization of our models of $G$-global homotopy theory. In fact, we will construct for both models a chain of four adjoint functors, in particular yielding two Bousfield localizations each. Before we can do any of this however, we need to understand a particular case of the `change of family' adjunction from Proposition~\ref{prop:change-of-family-sset} better:

\begin{defi}\label{defi:class-E}
We define $\mathcal E$ (for `equivariant')\nomenclature[aE]{$\mathcal E$}{collection of graph subgroups $\Gamma_{H,\phi}\subset\mathcal M\times G$ with $H$ universal and $\phi$ injective} as the collection of those graph subgroups $\Gamma_{H,\phi}$ of $\mathcal M\times G$ such that $H$ is universal and $\phi$ is \emph{injective}.
\end{defi}

Below we will need a characterization of the essential image of the left adjoint
\begin{equation}\label{eq:lambda-injective}
\lambda\colon\cat{$\bm{E\mathcal M}$-$\bm G$-SSet}_{\textup{$\mathcal E$-w.e.}}^\infty\to\cat{$\bm{E\mathcal M}$-$\bm G$-SSet}_{\textup{$G$-global w.e.}}^\infty
\end{equation}
of the localization functor. Let us give some intuition for this: on the left hand side, we only remember the fixed points for injective $\phi\colon H\to G$. If now $X$ is any $E\mathcal M$-$G$-simplicial set, then $\lambda(X)^\phi$ and $X^\phi$ are weakly equivalent by abstract nonsense. On the other hand, if $\psi\colon H\to G$ is not necessarily injective, then $\lambda(X)^\psi$ should be somehow computable from the fixed points for injective homomorphisms. A natural guess is that $\lambda(X)^\psi$ be weakly equivalent to $X^{\bar\psi}$ where $\bar\psi\colon H/\ker(\psi)\to G$ is the induced homomorphism (and we have secretly identified $H/\ker(\psi)$ with a universal subgroup of $\mathcal M$ isomorphic to it), and this indeed turns out to be true. However, we of course do not want $\lambda(X)^\psi$ and $\lambda(X)^{\bar\psi}$ to be merely abstractly weakly equivalent, but instead we want some explicit and suitably coherent way to identify them. The following definition turns this vague heuristic into a rigorous notion:

\begin{defi}\label{defi:ker-oblivious}\index{kernel oblivious|textbf}\index{kernel oblivious|(}
An $E\mathcal M$-$G$-simplicial set $X$ is called \emph{kernel oblivious} if the following holds: for any universal $H,H'\subset\mathcal M$, any surjective group homomorphism $\alpha\colon H\to H'$, any arbitrary group homomorphism $\phi\colon H'\to G$, and any $u\in\mathcal M$ such that $hu=u\alpha(h)$ for all $h\in H$, the map
\begin{equation*}
u.\blank\colon X^\phi\to X^{\phi\alpha}
\end{equation*}
is a weak homotopy equivalence of simplicial sets.
\end{defi}

\begin{ex}\label{ex:triv-oblivious}
Any $E\mathcal M$-$G$-simplicial set with trivial $E\mathcal M$-action is kernel oblivious: in fact, in this case $X^\phi=X^{\phi\alpha}$, and $u.\blank$ is just the identity.
\end{ex}

\begin{thm}\label{thm:im-lambda}
The following are equivalent for an $E\mathcal M$-$G$-simplicial set $X$:
\begin{enumerate}
\item $X$ is kernel oblivious.\label{item:il-oblivious}
\item $X$ lies in the essential image of $(\ref{eq:lambda-injective})$.\label{item:il-ess-im}
\item $X$ is \emph{$G$-globally} equivalent to an $\mathcal E$-cofibrant $X'\in\cat{$\bm{E\mathcal M}$-$\bm G$-SSet}$.\label{item:il-some-cof}
\item In any cofibrant replacement $\pi\colon X'\to X$ in the $\mathcal E$-model structure on $\cat{$\bm{E\mathcal M}$-$\bm G$-SSet}$, $\pi$ is actually a $G$-global weak equivalence.\label{item:il-any-cof}
\end{enumerate}
\end{thm}

The proof will require some preparations.

\begin{lemma}\label{lemma:I-between-ker-oblivious}
Let $f\colon X\to Y$ be an $\mathcal E$-weak equivalence in $\cat{$\bm{E\mathcal M}$-$\bm G$-SSet}$ such that $X$ and $Y$ are kernel oblivious. Then $f$ is a $G$-global weak equivalence.
\begin{proof}
Let $H\subset\mathcal M$ be universal and let $\phi\colon H\to G$ be any group homomorphism. We have to show that $f^\phi$ is a weak homotopy equivalence.

For this we choose a universal subgroup $H'\subset\mathcal M$ together with an isomorphism $H'\cong H/\ker\phi$, which gives rise to a surjective homomorphism $\alpha\colon H\to H'$ with $\ker(\alpha)=\ker(\phi)$. By the universal property of quotients, there then exists a unique $\bar\phi\colon H'\to G$ with $\bar\phi\alpha=\phi$; moreover, $\bar\phi$ is injective.

We now appeal to Corollary~\ref{cor:emg-equiv-group-hom-realization-general} to find a $u\in\mathcal M$ such that $hu=u\alpha(h)$ for all $h\in H$, yielding a commutative diagram
\begin{equation*}
\begin{tikzcd}
X^{\bar\phi}\arrow[d, "f^{\bar\phi}"']\arrow[r, "u.\blank"] & X^{\bar\phi\alpha}=X^\phi\arrow[d, "f^\phi", shift left=13pt]\\
Y^{\bar\phi}\arrow[r, "u.\blank"'] & Y^{\bar\phi\alpha}=Y^\phi.
\end{tikzcd}
\end{equation*}
The horizontal maps are weak equivalences as $X$ and $Y$ are kernel oblivious, and the left hand vertical map is a weak equivalence because $f$ is an $\mathcal E$-weak equivalence. Thus, also the right hand map is a weak equivalence as desired.
\end{proof}
\end{lemma}

\begin{prop}\label{prop:EM-inj-discrete}
Let $K\subset\mathcal M$ be any subgroup and let $\psi\colon K\to G$ be an \emph{injective} homomorphism. Then the projection $E\mathcal M\times_\psi G\to G/\im\psi$ is a $G$-global weak equivalence (where the right hand side carries the trivial $E\mathcal M$-action).
\begin{proof}
The map in question is conjugate to the image of the unique map $p\colon E\mathcal M\to *$ under $\psi_!\colon\cat{$\bm{E\mathcal M}$-$\bm{K}$-SSet}\to\cat{$\bm{E\mathcal M}$-$\bm G$-SSet}$, where $K$ acts on $E\mathcal M$ from the right via $k.(u_0,\dots,u_n)=(u_0k^{-1},\dots,u_nk^{-1})$.

As $\psi_!$ is homotopical (Corollary~\ref{cor:alpha-shriek-injective-EM}), it is then enough to show that $p$ is a $K$-global weak equivalence, which is just an instance of Example~\ref{ex:G-globally-contractible}.
\end{proof}
\end{prop}

\begin{cor}
Let $K\subset\mathcal M$ be any subgroup, let $\psi\colon K\to G$ be injective, and let $L$ be any simplicial set. Then $(E\mathcal M\times_\psi G)\times L$ is kernel oblivious.
\begin{proof}
The kernel oblivious $E\mathcal M$-$G$-simplicial sets are obviously closed under $G$-global weak equivalences. Thus, the claim follows from the previous proposition together with Example~\ref{ex:triv-oblivious}.
\end{proof}
\end{cor}

\begin{proof}[Proof of Theorem~\ref{thm:im-lambda}]
The implications $(\ref{item:il-any-cof})\Rightarrow(\ref{item:il-ess-im})$ and $(\ref{item:il-ess-im})\Rightarrow(\ref{item:il-some-cof})$ follow immediately from Remark~\ref{rk:change-of-family-sset-explicit}. For the remaining implications we will use:

\begin{claim*}
If $X$ is cofibrant in $\cat{$\bm{E\mathcal M}$-$\bm G$-SSet}_{\mathcal E}$, then $X$ is kernel oblivious.
\begin{proof}
Let $\alpha,\phi,u$ as in Definition~\ref{defi:ker-oblivious}. It suffices to verify that the natural transformation $u.\blank\colon (\blank)^\phi\Rightarrow(\blank)^{\phi\alpha}$ satisfies the assumptions of Corollary~\ref{cor:saturated-trafo}. But indeed, Condition~$(\ref{item:st-base})$ is an instance of the previous corollary, while all the remaining conditions are part of Lemma~\ref{lemma:fixed-points-cellular}.
\end{proof}
\end{claim*}

The implication $(\ref{item:il-some-cof})\Rightarrow(\ref{item:il-oblivious})$ now follows immediately from the claim. On the other hand, if $X$ is kernel oblivious and $\pi\colon X'\to X$ is any $\mathcal E$-cofibrant replacement, then also $X'$ is kernel oblivious by the claim, so that $\pi$ is a $G$-global weak equivalence by Lemma~\ref{lemma:I-between-ker-oblivious}. This shows $(\ref{item:il-oblivious})\Rightarrow(\ref{item:il-any-cof})$, finishing the proof.\index{kernel oblivious|)}
\end{proof}

We now consider the functor $\triv_{E\mathcal M}\colon\cat{$\bm G$-SSet}\to\cat{$\bm{E\mathcal M}$-$\bm G$-SSet}$ that equips a $G$-simplicial set with the trivial $E\mathcal M$-action. The equality of functors
\begin{equation}\label{eq:triv-fixed-points}
(\blank)^\phi\circ\triv_{E\mathcal M}=(\blank)^{\im\phi}
\end{equation}
(for any $H\subset\mathcal M$ and any homomorphism $\phi\colon H\to G$) shows that this is homotopical with respect to the proper weak equivalences on the source and the $G$-global or $\mathcal E$-weak equivalences on the target.

\begin{cor}\label{cor:triv-vs-lambda-sset}
The diagram
\begin{equation}\label{diag:triv-lambda-triv-sset}
\begin{tikzcd}[column sep=small]
& \cat{$\bm G$-SSet}_{\textup{proper}}^\infty\arrow[dl, "\triv_{E\mathcal M}^\infty"', bend right=10pt]\arrow[dr, "\triv_{E\mathcal M}^\infty", bend left=10pt]\\
\cat{$\bm{E\mathcal M}$-$\bm G$-SSet}_{\mathcal E}^\infty\arrow[rr, "\lambda"'] &&\cat{$\bm{E\mathcal M}$-$\bm G$-SSet}_{\textup{$G$-global}}^\infty
\end{tikzcd}
\end{equation}
commutes up to canonical equivalence.
\begin{proof}
This is is obviously true if we replace $\lambda$ by its right adjoint; in particular, there is a natural transformation filling the above, induced by the unit $\eta$ of $\lambda\dashv\textup{localization}$. To see that this is in fact an equivalence, it suffices (as $\lambda$ is fully faithful) that the right hand arrow lands in the essential image of $\lambda$. But by Theorem~\ref{thm:im-lambda} this is equivalent to demanding that $\triv_{E\mathcal M}X$ be kernel oblivious for every $G$-simplicial set $X$, which is just an instance of Example~\ref{ex:triv-oblivious}.
\end{proof}
\end{cor}

We can now prove the comparison between $G$-global and proper $G$-equivariant homotopy theory:

\begin{thm}\label{thm:G-global-vs-proper-sset}
The functor $\triv_{E\mathcal M}\colon\cat{$\bm G$-SSet}_{\textup{proper}}\to\cat{$\bm{E\mathcal M}$-$\bm G$-SSet}_{\textup{$G$-global}}$\nomenclature[atriv]{$\triv$}{embedding of proper $G$-equivariant into $G$-global homotopy theory via trivial actions} descends to a fully faithful functor on associated quasi-categories with essential image the kernel oblivious $E\mathcal M$-$G$-simplicial sets.

This induced functor admits both a left adjoint $\textbf{\textup L}(E\mathcal M\backslash\blank)$\nomenclature[aLEM]{$\textbf{\textup L}(E\mathcal M\backslash\blank)$}{left adjoint to $\triv$} as well as a right adjoint $(\blank)^{\textbf{\textup R}E\mathcal M}$.\nomenclature[aREM]{$(\blank)^{\textbf{\textup R}E\mathcal M}$}{right adjoint to $\triv$} Moreover, $(\blank)^{\textbf{\textup R}E\mathcal M}$
is a quasi-localization at the $\mathcal E$-weak equivalences, and it in turn admits another right adjoint $\mathcal R$,\nomenclature[aR]{$\mathcal R$}{right adjoint to right adjoint to usual embedding from proper $G$-equivariant into $G$-global homotopy theory ($\triv$ or $\const$)} which is again fully faithful.
\end{thm}

In particular, Example~\ref{ex:triv-oblivious} accounts for all kernel oblivious $E\mathcal M$-$G$-simplicial sets up to homotopy.

\begin{proof}
We first observe that the adjunctions
\begin{equation}\label{eq:triv-EM-G-global}
E\mathcal M\backslash\blank\colon\cat{$\bm{E\mathcal M}$-$\bm{G}$-SSet}_{\textup{$G$-global}}\rightleftarrows
\cat{$\bm{G}$-SSet}_{\textup{proper}} :\!\triv_{E\mathcal M}
\end{equation}
and
\begin{equation}\label{eq:triv-EM-I}
E\mathcal M\backslash\blank\colon\cat{$\bm{E\mathcal M}$-$\bm{G}$-SSet}_{\mathcal E}\rightleftarrows
\cat{$\bm{G}$-SSet}_{\textup{proper}} :\!\triv_{E\mathcal M}
\end{equation}
are Quillen adjunctions with homotopical right adjoints by the equality $(\ref{eq:triv-fixed-points})$, and in particular $(\ref{eq:triv-EM-G-global})$ induces the desired left adjoint of $\triv_{E\mathcal M}^\infty$.

To construct the right adjoint of $\triv_{E\mathcal M}$ it suffices to observe that while
\begin{equation*}
\triv_{E\mathcal M}\colon\cat{$\bm G$-SSet}\rightleftarrows\cat{$\bm{E\mathcal M}$-$\bm G$-SSet} :\!(\blank)^{E\mathcal M}
\end{equation*}
is not a Quillen adjunction with respect to the above model structures, it becomes one once we use Corollary~\ref{cor:enlarge-generating-cof} to enlarge the generating cofibrations of the $G$-global model structure as to contain all $G/H\times\del\Delta^n\hookrightarrow G/H\times\Delta^n$ for $H\subset G$ finite and $n\ge 0$ (which we are allowed to do by Lemma~\ref{lemma:homotopy-pushout-M-SSet}); in particular, $\triv_{E\mathcal M}$ is left Quillen with respect to the injective $G$-global model structure on the target.

To prove the remaining statements we will need:

\begin{claim*}
The Quillen adjunction $(\ref{eq:triv-EM-I})$ is a Quillen equivalence. In particular,
\begin{equation*}
\triv_{E\mathcal M}^\infty\colon\cat{$\bm G$-SSet}_{\textup{proper}}^\infty\to\cat{$\bm{E\mathcal M}$-$\bm{G}$-SSet}_{\textup{$\mathcal E$-w.e.}}^\infty
\end{equation*}
is an equivalence of quasi-categories.
\begin{proof}
It suffices to prove the first statement. The equality $(\ref{eq:triv-fixed-points})$ shows that the right adjoint preserves and reflects weak equivalences. It is therefore enough that the ordinary unit $\eta\colon X\to\triv_{E\mathcal M}(E\mathcal M\backslash X)$, which is given by the projection map, is an $\mathcal E$-weak equivalence for any $\mathcal E$-cofibrant $E\mathcal M$-$G$-simplicial set $X$. By the above, both $E\mathcal M\backslash\blank$ as well as $\triv_{E\mathcal M}$ are left Quillen (after suitably enlarging the cofibrations in the target), and they moreover commute with tensoring with simplicial sets. By Corollary~\ref{cor:saturated-trafo} it therefore suffices that $\eta$ is a weak equivalence for every $E\mathcal M\times_\psi G$ with $K\subset\mathcal M$ universal and $\psi\colon K\to G$ injective.

An easy calculation shows that the projection $E\mathcal M\times G\to G$ descends to an isomorphism $E\mathcal M\backslash(E\mathcal M\times_\psi G)\to G/\im\psi$, so we want to show that the projection $E\mathcal M\times_\psi G\to G/\im\psi$ is an $\mathcal E$-weak equivalence. But this is actually even a $G$-global weak equivalence by Proposition~\ref{prop:EM-inj-discrete}, finishing the proof of the claim.
\end{proof}
\end{claim*}

We now contemplate the diagram $(\ref{diag:triv-lambda-triv-sset})$ from Corollary~\ref{cor:triv-vs-lambda-sset}. By the above claim together with Proposition~\ref{prop:change-of-family-sset} we then conclude that
\begin{equation*}
\triv_{E\mathcal M}^\infty\colon\cat{$\bm G$-SSet}_{\textup{proper}}^\infty\to\cat{$\bm{E\mathcal M}$-$\bm{G}$-SSet}_{\textup{$G$-global}}^\infty
\end{equation*}
is fully faithful, and by Theorem~\ref{thm:im-lambda} its essential image then consists precisely of the kernel oblivious $E\mathcal M$-$G$-simplicial sets. Moreover, we deduce by uniqueness of adjoints that its right adjoint $(\blank)^{\textbf{R}E\mathcal M}$ is canonically equivalent to the composite
\begin{equation}\label{eq:EM-derived-fixpoints-alternative}
\begin{aligned}
\cat{$\bm{E\mathcal M}$-$\bm{G}$-SSet}_{\textup{$G$-global}}^\infty&\xrightarrow{\textup{localization}}
\cat{$\bm{E\mathcal M}$-$\bm{G}$-SSet}_{\textup{$\mathcal E$-w.e.}}^\infty\\
&\xrightarrow{(\triv_{E\mathcal M}^\infty)^{-1}}\cat{$\bm G$-SSet}_{\textup{proper}}
\end{aligned}
\end{equation}
(where the right hand arrow denotes any chosen quasi-inverse) and hence indeed a quasi-localization at the $\mathcal E$-weak equivalences. Finally, $(\ref{eq:EM-derived-fixpoints-alternative})$ has a fully faithful right adjoint by Proposition~\ref{prop:change-of-family-sset}, given explicitly by
\begin{equation*}
\cat{$\bm G$-SSet}_{\textup{proper}}^\infty\xrightarrow{\triv_{E\mathcal M}^\infty}
\cat{$\bm{E\mathcal M}$-$\bm{G}$-SSet}_{\textup{$\mathcal E$-w.e.}}^\infty\xrightarrow{\rho}\cat{$\bm{E\mathcal M}$-$\bm{G}$-SSet}_{\textup{$G$-global}}^\infty
\end{equation*}
which is then also right adjoint to $(\blank)^{\textbf{R}E\mathcal M}$, finishing the proof.
\end{proof}

\begin{rk}
In summary, we in particular have two Bousfield localizations
\begin{equation*}
\textbf{L}(E\mathcal M\backslash\blank)\dashv\triv_{E\mathcal M}\qquad\text{and}\qquad(\blank)^{\textbf{R}E\mathcal M}\dashv\mathcal R.
\end{equation*}
In the ordinary global setting one is mostly interested in the analogue of the adjunction $\triv_{E\mathcal M}\dashv(\blank)^{\textbf{R}E\mathcal M}$ (cf.~Remark~\ref{rk:R-vs-R}) which is a `wrong way' (i.e.~right) Bousfield localization.
\end{rk}

\begin{warn}
Using the Quillen equivalence $(\ref{eq:triv-EM-I})$, we can give another description of the composition $(\ref{eq:EM-derived-fixpoints-alternative})$ and hence of the right adjoint $(\blank)^{\textbf{R}E\mathcal M}$ of $\triv_{E\mathcal M}^\infty$: namely, this can be computed by taking a cofibrant replacement \emph{with respect to the $\mathcal E$-model structure} and then dividing out the left $E\mathcal M$-action.

In contrast to this, the \emph{left} adjoint of $\triv_{E\mathcal M}^\infty$ is computed by taking a cofibrant replacement \emph{with respect to the $G$-global model structure} and then dividing out the action. These two functors are \emph{not} equivalent even for $G=1$: namely, let $H\subset\mathcal M$ be any non-trivial universal subgroup and consider the projection $p\colon E\mathcal M\to E\mathcal M/H$. As $E\mathcal M$ and $E\mathcal M/H$ are both cofibrant in $\cat{$\bm{E\mathcal M}$-SSet}$, we can calculate the value of $\textbf{L}(E\mathcal M\backslash\blank)$ at $p$ simply by $E\mathcal M\backslash p\colon E\mathcal M\backslash E\mathcal M\to E\mathcal M\backslash E\mathcal M/H$, which is a map between terminal objects, hence in particular an equivalence.

On the other hand, $p$ is not an $\mathcal E$-weak equivalence (i.e.~underlying weak equivalence) as $E\mathcal M/H$ is a $K(H,1)$ while $E\mathcal M$ is contractible. But $(\blank)^{\textbf RE\mathcal M}$ is a quasi-localization at the $\mathcal E$-weak equivalences, and these are saturated as they are part of a model structure. Thus $p^{\textbf{R}E\mathcal M}$ is not an equivalence, and in particular it cannot be conjugate to $\textbf{L}(E\mathcal M\backslash\blank)(p)$.
\end{warn}

The above descriptions of the adjoints are not really suitable for computations. However, for finite $G$ we can give easier constructions of $(\blank)^{\textbf{R}E\mathcal M}$ and $\mathcal R$:

\begin{prop}\label{prop:G-finite-simple}
Assume $G$ is finite and choose an injective homomorphism $i\colon G\to\mathcal M$ with universal image, inducing $(i,\id)\colon G\to E\mathcal M\times G$. Then
\begin{equation}\label{eq:G-global-vs-G-finite}
(i,\id)^*\colon\cat{$\bm{E\mathcal M}$-$\bm G$-SSet}_{\textup{$G$-global}}\rightleftarrows\cat{$\bm G$-SSet} :\!(i,\id)_*
\end{equation}
is a Quillen adjunction with homotopical left adjoint, and there are preferred equivalences
\begin{equation*}
\big((i,\id)^*\big)^\infty\simeq(\blank)^{\textbf{\textup R}E\mathcal M}\qquad \text{and}\qquad\textbf{\textup R}(i,\id)_*\simeq\mathcal R.
\end{equation*}
\begin{proof}
It is obvious from the definition that $(i,\id)^*$ sends $\mathcal E$-weak equivalences (and hence in particular $G$-global weak equivalences) to $G$-weak equivalences. Moreover, it preserves cofibrations as the cofibrations on the right hand side are just the underlying cofibrations. We conclude that $(\ref{eq:G-global-vs-G-finite})$ is a Quillen adjunction and that $(i,\id)^*$ descends to $\cat{$\bm{E\mathcal M}$-$\bm G$-SSet}_{\mathcal E}^\infty\to\cat{$\bm G$-SSet}^\infty$.

On the other hand, by Theorem~\ref{thm:G-global-vs-proper-sset} also $(\blank)^{\cat{R}E\mathcal M}$ descends accordingly and the resulting functor is quasi-inverse to the one induced by $\triv_{E\mathcal M}$.
The equality $(i,\id)^*\circ\triv_{E\mathcal M}=\id_{\cat{$\bm G$-SSet}}$ of homotopical functors then also exhibits the functor induced by $(i,\id)^*$ on $\cat{$\bm{E\mathcal M}$-$\bm G$-SSet}_{\mathcal E}^\infty$ as left quasi-inverse to the one induced by $\triv_{E\mathcal M}$ which provides the first equivalence. The second one is then immediate from the uniqueness of adjoints.
\end{proof}
\end{prop}

Finally, let us consider the analogues for $\mathcal M$-actions:

\begin{cor}
The homotopical functor
\begin{equation*}
\triv_{\mathcal M}\colon \cat{$\bm G$-SSet}_{\textup{proper}}\to\cat{$\bm{\mathcal M}$-$\bm G$-SSet}_{\textup{$G$-global}}
\end{equation*}
descends to a fully faithful functor on associated quasi-categories. This induced functor admits both a left adjoint $\textbf{\textup L}(\mathcal M\backslash\blank)$\nomenclature[aLM]{$\textbf{\textup L}(\mathcal M\backslash\blank)$}{left adjoint to $\triv$} as well as a right adjoint $(\blank)^{\textbf{\textup R}\mathcal M}$.\nomenclature[aRM]{$(\blank)^{\textbf{\textup R}\mathcal M}$}{right adjoint to $\triv$} The latter is a quasi-localization at those $f$ such that $E\mathcal M\times_{\mathcal M}^{\textbf{\textup L}}f$ is an $\mathcal E$-weak equivalence, and it in turn admits another right adjoint $\mathcal R$ which is again fully faithful.

Moreover, the diagram
\begin{equation}\label{diag:triv-M-vs-EM}
\begin{tikzcd}[column sep=small]
& \cat{$\bm G$-SSet}_{\textup{proper}}^\infty\arrow[dl, "\triv_{E\mathcal M}^\infty"', bend right=10pt]\arrow[dr, "\triv_{\mathcal M}^\infty", bend left=10pt]\\
\cat{$\bm{E\mathcal M}$-$\bm G$-SSet}^\infty_{\textup{$G$-global}}\arrow[rr, "\forget^\infty"'] & & \cat{$\bm{\mathcal M}$-$\bm G$-SSet}^\infty_{\textup{$G$-global}}
\end{tikzcd}
\end{equation}
commutes up to canonical equivalence.
\end{cor}

It follows formally that the forgetful functor is also compatible with the remaining functors in the two adjoint chains constructed above, and likewise for its own adjoints $E\mathcal M\times^{\cat{L}}_{\mathcal M}\blank$ and $\cat{R}\Maps^{\mathcal M}(E\mathcal M,\blank)$.

\begin{proof}
We obviously have a Quillen adjunction
\begin{equation*}
\triv_{\mathcal M}\colon\cat{$\bm G$-SSet}_{\textup{proper}}\rightleftarrows\cat{$\bm{\mathcal M}$-$\bm G$-SSet}_{\textup{injective $G$-global}} :\!(\blank)^{\mathcal M},
\end{equation*}
justifying the above description of the right adjoint; in fact, as in the proof of Theorem~\ref{thm:G-global-vs-proper-sset}, it would have been enough to enlarge the generating cofibrations by the maps $G/H\times\del\Delta^n\hookrightarrow G/H\times\Delta^n$ for $H$ finite and $n\ge0$.

For the left adjoint, we now want to prove that also
\begin{equation*}
\mathcal M\backslash\blank\colon\cat{$\bm{\mathcal M}$-$\bm G$-SSet}_{\textup{$G$-global}}\rightleftarrows\cat{$\bm G$-SSet}_{\textup{proper}} :\!\triv_{\mathcal M}
\end{equation*}
is a Quillen adjunction. For this we observe that this is true for the $G$-universal model structure on the left hand side (as $\triv_{\mathcal M}$ is then obviously right Quillen); in particular $\mathcal M\backslash\blank$ preserves cofibrations and $\triv_{\mathcal M}$ sends fibrant $G$-simplicial sets to $G$-universally fibrant $\mathcal M$-$G$-simplicial sets. As this adjunction has an obvious simplicial enrichment, it therefore suffices by Proposition~\ref{prop:cofibrations-fibrant-qa} and the characterization of the fibrant objects provided in Corollary~\ref{cor:em-vs-m-equiv-model-cat} that $\triv_{\mathcal M}$ has image in the $G$-semistable $\mathcal M$-$G$-simplicial sets, which is in fact obvious from the definition.

To prove that $(\ref{diag:triv-M-vs-EM})$ commutes up to canonical equivalence, it suffices to observe that the evident diagram of homotopical functors inducing it actually commutes on the nose. All the remaining statements then follow formally from the commutativity of $(\ref{diag:triv-M-vs-EM})$ as before.\index{proper G-equivariant homotopy theory@proper $G$-equivariant homotopy theory!vs G-global homotopy theory@vs.~$G$-global homotopy theory|)}
\end{proof}

\section{Tameness}\label{sec:tame}
In this section we will be concerned with the notion of \emph{tameness} for $\mathcal M$-actions and $E\mathcal M$-actions, and we will in particular show that the models from the previous sections have tame analogues, that still model unstable $G$-global homotopy theory.

Our reason to study tameness here is twofold: firstly, the categories of tame $\mathcal M$- and $E\mathcal M$-simplicial sets carry interesting symmetric monoidal structures, which we will study in Chapter~\ref{chapter:coherent}, and which play a central role in the construction of global algebraic $K$-theory \cite{schwede-k-theory}; secondly, tame $\mathcal M$- and $E\mathcal M$-simplicial sets are intimately connected to the diagram spaces we will consider in the next section, and the theory developed here will be crucial in establishing those models.

\subsection{A reminder on tame \texorpdfstring{$\bm{\mathcal M}$}{M}-actions}
We begin with the notion of tame $\mathcal M$-actions, which first appeared (for actions on abelian groups) as \cite[Definition~1.4]{schwede-semistable}, and which were then further studied (for actions on sets and simplicial sets) in \cite{I-vs-M-1-cat}, see in particular~\cite[Definition~2.2 and Definition~3.1]{I-vs-M-1-cat}.

\index{support!for M-actions@for $\mathcal M$-actions|(}
\begin{defi}\label{defi:M-set-tame-support}
Let $A\subset\omega$ be any finite set. We write $\mathcal M_A\subset\mathcal M$ for the submonoid of those $u\in\mathcal M$ that fix $A$ pointwise, i.e.~$u(a)=a$ for all $a\in A$.

Let $X$ be any $\mathcal M$-set. An element $x\in X$ is said to be \emph{supported on $A$} if $u.x=x$ for all $u\in\mathcal M_A$; we write $X_{[A]}\subset X$ for the subset of those elements that are supported on $A$, i.e.~$X_{[A]}=X^{\mathcal M_A}$.\nomenclature[aA]{$(\blank)_{[A]}$}{(simplicial) subset of elements supported on $A$}

We call $x$ \emph{finitely supported} if it is supported on some finite set, and we write
\begin{equation*}
X^\tau\mathrel{:=}\bigcup_{A\subset\omega\textup{ finite}} X_{[A]}
\end{equation*}
for the subset of all finitely supported elements. We call $X$ \emph{tame} if $X=X^\tau$.\index{tame!M-action@$\mathcal M$-action|textbf}\nomenclature[atau]{$(\blank)^\tau$}{(simplicial) subset of finitely supported elements; subcategory of tame objects}
\end{defi}

On $\mathcal M$-simplicial sets everything can be extended levelwise:

\begin{defi}\label{defi:tame-M-simplicial-set}
Let $X$ be an $\mathcal M$-simplicial set. An $n$-simplex $x$ is \emph{supported} on the finite set $A\subset\omega$ if it is supported on $A$ as an element of the $\mathcal M$-set $X_n$ of $n$-simplices of $X$. Analogously, $x$ is said to be \emph{finitely supported} if it is so as an element of $X_n$. We define $X^\tau$ via $(X^\tau)_n=(X_n)^\tau$, i.e.~as the family of all finitely supported simplices. We call $X$ \emph{tame} if $X^\tau=X$.\index{tame!M-action@$\mathcal M$-action|textbf}
\end{defi}

\subsubsection{Basic properties} Let us record some basic facts about the above notions. All of these can be found explicitly in \cite{I-vs-M-1-cat} for $\mathcal M$-sets and are easily extended to $\mathcal M$-simplicial sets (for which they also appear implicitly in \emph{op.~cit.}).

\begin{lemma}\label{lemma:support-vs-morphism-M}\label{lemma:first-M-basic}
\begin{enumerate}
\item If $f\colon X\to Y$ is a map of $\mathcal M$-sets and $A\subset\omega$ is finite, then $f$ restricts to $f_{[A]}\colon X_{[A]}\to Y_{[A]}$, hence in particular to $f^\tau\colon X^\tau\to Y^\tau$.
\item If $X$ is any $\mathcal M$-simplicial set and $A\subset\omega$ is finite, then $X_{[A]}$ is a simplicial subset of $X$. In particular, $X^\tau$ is a simplicial subset.
\item If $f\colon X\to Y$ is a map of $\mathcal M$-simplicial sets and $A\subset\omega$ is finite, then $f$ restricts to $f_{[A]}\colon X_{[A]}\to Y_{[A]}$. In particular, it restricts to $f^\tau\colon X^\tau\to Y^\tau$.
\end{enumerate}
\begin{proof}
The first statement is a trivial calculation which we omit; this also appears without proof in \cite[discussion before Lemma~2.6]{I-vs-M-1-cat}. The second statement follows by applying the first one to the structure maps, and the third one follows then by applying the first one levelwise.
\end{proof}
\end{lemma}

\begin{rk}\label{rk:tame-M-colimits}
The above lemma provides us with a functor $(\blank)^\tau\colon\cat{$\bm{\mathcal M}$-SSet}\to\cat{$\bm{\mathcal M}$-SSet}^\tau$ right adjoint to the inclusion of the full subcategory of tame $\mathcal M$-simplicial sets. It follows formally that $\cat{$\bm{\mathcal M}$-SSet}^\tau$ is complete and cocomplete, with colimits created in $\cat{SSet}$, also see~\cite[Lemma~2.6]{I-vs-M-1-cat}.

In fact, it is trivial to check that tame $\mathcal M$-simplicial sets are preserved by \emph{finite} products and passing to $\mathcal M$-subsets, so also finite limits in $\cat{$\bm{\mathcal M}$-SSet}^\tau$ are created in $\cat{SSet}$, also cf.~\cite[Example~4.11]{schwede-k-theory}.
\end{rk}

\begin{defi}\index{support!for M-actions@for $\mathcal M$-actions|textbf}
Let $X$ be any $\mathcal M$-set and let $x\in X$ be finitely supported. Then the \emph{support} $\supp(x)$\nomenclature[asupp]{$\supp$}{support (with respect to $\mathcal M$- or $E\mathcal M$-action)} is the intersection of all finite sets $A\subset\omega$ on which $x$ is supported.
\end{defi}

\begin{lemma}
In the above situation, $x$ is supported on $\supp(x)$.
\begin{proof}
This is immediate from \cite[Proposition~2.3]{I-vs-M-1-cat}.
\end{proof}
\end{lemma}

\begin{ex}\label{ex:M-Inj-support}
Let $A$ be a finite set. Then the $\mathcal M$-set $\Inj(A,\omega)$ is tame and the support of an injection $i\colon A\to\omega$ is simply its image, also see~\cite[Example~2.9]{I-vs-M-1-cat}: namely, it is clear from the definition that any injection $i\colon A\to\omega$ is supported on its image; on the other hand, if $B\not\supset\im(i)$ is any finite set, then we pick an $a\in\im(i)\setminus B$ together with an injection $u\in\mathcal M_B$ such that $a\notin\im(u)$. Then $a\notin\im(ui)$, so $u.i\not=i$, and hence $i$ cannot be supported on $B$.

On the other hand, if $A$ is countably infinite, then a similar computation shows that $\Inj(A,\omega)$ is not tame, and in fact even $\Inj(A,\omega)^\tau=\varnothing$. In particular, $\mathcal M$ itself is not tame.
\end{ex}

\begin{lemma}\label{lemma:support-vs-action-M}
Let $u\in\mathcal M$ and let $X$ be any $\mathcal M$-simplicial set. Then $\supp(u.x)=u(\supp(x))$ for any finitely supported simplex $x$. In particular, $u.\blank\colon X\to X$ restricts to $X_{[A]}\to X_{[u(A)]}$ for any finite $A\subset\omega$, and $X^\tau$ is an $\mathcal M$-simplicial subset of $X$.
\begin{proof}
The first statement is \cite[Proposition~2.5-(ii)]{I-vs-M-1-cat}, which immediately implies the second one. The final statement then in turn follows from the second one, also cf.~\cite[discussion after Proposition~2.5-(ii)]{I-vs-M-1-cat}.
\end{proof}
\end{lemma}

\begin{lemma}\label{lemma:support-agree-M}\label{lemma:last-M-basic}
Let $X\in\cat{$\bm{\mathcal M}$-SSet}$, $u,u'\in\mathcal M$, and let $x\in X_n$. Assume that $x$ is supported on some finite set $A\subset\omega$ such that $u|_A=u'|_A$. Then $u.x=u'.x$.
\begin{proof}
This is \cite[Proposition~2.5-(i)]{I-vs-M-1-cat}, applied to the $\mathcal M$-set $X_n$.\index{support!for M-actions@for $\mathcal M$-actions|)}
\end{proof}
\end{lemma}

\subsubsection{The structure of tame \texorpdfstring{$\mathcal M$-$G$}{M-G}-simplicial sets}
\cite[Theorem~2.11]{I-vs-M-1-cat}, which we recall below, describes how tame $\mathcal M$-sets decompose into some standard pieces. As a consequence of this we will prove:

\begin{thm}\label{thm:structure-tame-M-G-SSet}
Let us define
\begin{equation}\label{eq:definition-I-tame}
\begin{aligned}
I_{\textup{tame}}\mathrel{:=}&\{(\Inj(A,\omega)\times_{\Sigma_A}X)\times\del\Delta^n\hookrightarrow(\Inj(A,\omega)\times_{\Sigma_A}X)\times\Delta^n :\\
&\hphantom{\{}\text{$A\subset\omega$ finite, $X$ a $G$-$\Sigma_A$-biset, $n\ge 0$}\}.
\end{aligned}
\end{equation}
Then the $I_{\textup{tame}}$-cell complexes in $\cat{$\bm{\mathcal M}$-$\bm G$-SSet}$ are precisely the tame $\mathcal M$-$G$-simplicial sets.
\end{thm}

If $X$ is a tame $\mathcal M$-set, let us write $s_n(X)$ for the subset of those $x\in X$ with $\supp(x)=\{1,\dots,n\}$. Lemma~\ref{lemma:support-agree-M} provides for any $u\in\Inj(\{1,\dots,n\},\omega)$ a well-defined map $s_n(X)\to X$ obtained by acting with any extension of $u$ to an $\bar{u}\in\mathcal M$, and Lemma~\ref{lemma:support-vs-action-M} shows that this restricts to a $\Sigma_n$-action on $s_n(X)$.

\begin{thm}\label{thm:structure-tame-m-set}
Let $X$ be any tame $\mathcal M$-set. Then the map
\begin{equation}\label{eq:tame-decomposition}
\coprod_{n=0}^\infty\Inj(\{1,\dots,n\},\omega)\times_{\Sigma_n}s_n(X)\to X
\end{equation}
induced by the above construction is well-defined and an isomorphism of $\mathcal M$-sets.
\begin{proof}
See \cite[Theorem~2.11]{I-vs-M-1-cat}.
\end{proof}
\end{thm}

\begin{cor}\label{cor:structure-tame-m-g-set}
Let $X$ be any tame $\mathcal M$-$G$-set. Then each $s_n(X)$ is a $G$-subset of $X$, and the map $(\ref{eq:tame-decomposition})$ is an isomorphism of $\mathcal M$-$G$-sets.
\begin{proof}
In order to prove that $s_n(X)$ is closed under the action of $G$, we have to show that $\supp(g.x)=\supp(x)$ for all $g\in G$ and $x\in X$. The inclusion `$\subset$' is an instance of Lemma~\ref{lemma:support-vs-morphism-M} because the $G$-action commutes with the $\mathcal M$-action. The inclusion `$\supset$' then follows by applying the same argument to $g^{-1}$ and $g.x$, or by simply observing that \emph{injective} $\mathcal M$-equivariant maps strictly preserve supports.

Again using that the $\mathcal M$-action commutes with the $G$-action, we see that $(\ref{eq:tame-decomposition})$ is $G$-equivariant, hence an isomorphism of $\mathcal M$-$G$-sets by the previous theorem.
\end{proof}
\end{cor}

The only remaining ingredient for the proof of Theorem~\ref{thm:structure-tame-M-G-SSet} is the following:

\begin{lemma}\label{lemma:tame-complement}
Let $X$ be a tame $\mathcal M$-$G$-set and let $Y\subset X$ be an $\mathcal M$-$G$-subset. Then also $X\setminus Y$ is an $\mathcal M$-$G$-subset.
\end{lemma}

We caution the reader that the above is not true in general for non-tame actions---for example, the subset $Y\subset\mathcal M$ of \emph{non-surjective} maps is closed under the left regular $\mathcal M$-action, but its complement is not.

\begin{proof}
The set $X\setminus Y$ is obviously closed under the $G$-action. The fact that it is moreover closed under the $\mathcal M$-action appeared in a previous version of \cite{I-vs-M-1-cat}; let us give the argument for completeness. Assume $x\in X\setminus Y$ and let $u\in\mathcal M$ such that $u.x\in Y$. By tameness, there exists a finite set $A$ on which $x$ is supported. We now pick any invertible $v\in\mathcal M$ such that $v|_A=u|_A$. By Lemma~\ref{lemma:support-agree-M} we then have $v.x=u.x\in Y$ and hence also $x=v^{-1}.(v.x)\in Y$, which is a contradiction.
\end{proof}

\begin{proof}[Proof of Theorem~\ref{thm:structure-tame-M-G-SSet}]
Obviously all sources and targets of maps in $I_{\textup{tame}}$ are tame. As the tame $\mathcal M$-$G$-simplicial sets are closed under all colimits, we see that all $I_{\textup{tame}}$-cell complexes are tame (cf.~Lemma~\ref{lemma:saturated-objects}).

Conversely, let $X$ be any tame $\mathcal M$-$G$-simplicial set; we consider the usual skeleton filtration $\varnothing=X^{(-1)}\subset X^{(0)}\subset X^{(1)}\subset\cdots$ of $X$. It suffices to prove that each $X^{(n-1)}\to X^{(n)}$ is a relative $I_{\textup{tame}}$-cell complex.

For this we contemplate the (a priori non-equivariant) pushout
\begin{equation*}
\begin{tikzcd}
X_n^{\textup{nondeg}}\times\del\Delta^n\arrow[d]\arrow[r, hook] &X_n^{\textup{nondeg}}\times\Delta^n\arrow[d]\\
X^{(n-1)}\arrow[r,hook]&X^{(n)}
\end{tikzcd}
\end{equation*}
where $X_n^{\textup{nondeg}}$ denotes the subset of nondegenerate $n$-simplices. The \emph{degenerate} simplices obviously form an $\mathcal M$-$G$-subset of $X_n$, and hence so do the nondegenerate ones by the previous lemma. With respect to this action, all maps in the above square are obviously $\mathcal M$-$G$-equivariant so that this is a pushout in $\cat{$\bm{\mathcal M}$-$\bm G$-SSet}$. But applying Corollary~\ref{cor:structure-tame-m-g-set} to $X_n^{\textup{nondeg}}$ expresses the top horizontal map as a coproduct of maps in $I_{\textup{tame}}$, finishing the proof.
\end{proof}

\subsection{Tame \texorpdfstring{$\bm{E\mathcal M}$}{EM}-actions}
Next, we will introduce and study analogues of tameness and support for $E\mathcal M$-simplicial sets.

\begin{defi}\label{defi:tame-EM-simplicial-set}
Let $A\subset\omega$ be finite and let $X$ be an $E\mathcal M$-simplicial set. An $n$-simplex $x$ of $X$ is said to be \emph{supported on $A$} if $E(\mathcal M_A)$ acts trivially on $x$, i.e.~the composition
\begin{equation*}
E(\mathcal M_A)\times\Delta^n\hookrightarrow E\mathcal M\times\Delta^n\xrightarrow{E\mathcal M\times x} E\mathcal M\times X\xrightarrow{\textup{act}} X
\end{equation*}
agrees with
\begin{equation*}
E(\mathcal M_A)\times\Delta^n\xrightarrow{\pr} \Delta^n\xrightarrow{x} X.
\end{equation*}
The simplex $x$ is \emph{finitely supported} if it is supported on some finite set $A$.

We write $X_{[A]}$ for the family of simplices supported on a finite set $A\subset\omega$ and $X^\tau\mathrel{:=}\bigcup_{A\subset\omega\textup{ finite}} X_{[A]}$ for the family of all finitely supported simplices. We call $X$ \emph{tame} if $X=X^\tau$, i.e.~if all its simplices are finitely supported.\index{support!for EM-actions@for $E\mathcal M$-actions}
\index{tame!EM-actions@$E\mathcal M$-actions|textbf}
\end{defi}

\begin{ex}\label{ex:EM-Inj-support}
Let $A$ be a finite set. Then $E\Inj(A,\omega)$ is tame, and the support of an $n$-simplex $(i_0,\dots,i_n)$ is given by $B\mathrel{:=}\im(i_0)\cup\cdots\cup\im(i_n)$: namely, to prove that $(i_0,\dots,i_n)$ is supported on $B$ we have to show that
\begin{equation*}
(u_0,\dots,u_m).f^*(i_0,\dots,i_n)=f^*(i_0,\dots,i_n)
\end{equation*}
for any $f\colon [m]\to[n]$ in $\Delta$ and any $u_0,\dots,u_m\in\mathcal M_B$. But plugging in the definition, the left hand side evaluates to $(u_0i_{f(0)},\dots,u_mi_{f(m)})$ while the right hand side is given by $(i_{f(0)},\dots,i_{f(m)})$, so the claim is obvious. Conversely, one argues as in Example~\ref{ex:M-Inj-support} that $(i_0,\dots,i_n)$ is not supported on any finite $C\not\supset B$, or one simply notes that if $(i_0,\dots,i_n)$ is supported on $C$, then so are the individual injections $i_0,\dots,i_n$ as elements of the $\mathcal M$-set $\Inj(A,\omega)$.

Again, one similarly shows that $(E\mathcal M)^\tau=\varnothing$; in particular, $E\mathcal M$ is not tame.
\end{ex}

\begin{lemma}\label{lemma:EM-support-preserve}
Let $f\colon X\to Y$ be a map of $E\mathcal M$-simplicial sets, and let $A\subset\omega$ be any finite set. Then $X_{[A]}$ and $Y_{[A]}$ are simplicial subsets of $X$ and $Y$, respectively, and $f$ restricts to $f_{[A]}\colon X_{[A]}\to Y_{[A]}$.
\begin{proof}
The first statement is clear from the definition.

For the second statement we have to show that $(u_0,\dots,u_m).g^*f(x)=g^*f(x)$ for any $g\colon[m]\to[n]$ in $\Delta$ and any $u_0,\dots,u_m\in\mathcal M_A$. But using that $f$ is simplicial and $E\mathcal M$-equivariant, the left hand side equals $f((u_0,\dots,u_m).g^*x)$, while the right hand side equals $f(g^*x)$, so the claim follows from the definitions.
\end{proof}
\end{lemma}

It is a straight-forward but somewhat lengthy endeavor to also verify the analogues of the other basic properties of tame $\mathcal M$-actions established above in the world of $E\mathcal M$-simplicial sets. We will not do this at this point as they have shorter proofs once we know the following theorem, that is also of independent interest:

\begin{thm}\label{thm:support-EM-vs-M}\index{support!these notions agree|textbf}
Let $X$ be an $E\mathcal M$-simplicial set, let $n\ge 0$, and let $A\subset\omega$ be any finite set. Then $x\in X_n$ is supported on $A$ in the sense of Definition~\ref{defi:tame-EM-simplicial-set} if and only if it is supported on $A$ as a simplex of the underlying $\mathcal M$-simplicial set of $X$ (see Definition~\ref{defi:tame-M-simplicial-set}). In other words, $X_{[A]}=(\forget X)_{[A]}$ as simplicial sets, where $\forget$ denotes the forgetful functor $\cat{$\bm{E\mathcal M}$-SSet}\to\cat{$\bm{\mathcal M}$-SSet}$.

Moreover, the subfunctor $(\blank)_{[A]}\colon\cat{$\bm{E\mathcal M}$-SSet}\to\cat{SSet}$ of the forgetful functor is corepresented in the enriched sense by $E\Inj(A,\omega)$ via evaluation at the $0$-simplex given by the inclusion $A\hookrightarrow\omega$.
\end{thm}

The proof requires some combinatorial preparations:

\begin{prop}
Let $A\subset\omega$ be any finite set, and let $u_0,\dots,u_n\in\mathcal M$. Then there exists a $\chi\in\mathcal M_A$ such that $\im(u_0\chi)\cup\cdots\cup\im(u_n\chi)$ has infinite complement in $\omega$.
\begin{proof}
We will construct strictly increasing chains $B_0\subsetneq B_1\subsetneq\cdots$ and $C_0\subsetneq C_1\subsetneq\cdots$ of finite subsets of $\omega\setminus A$ and $\omega$, respectively, such that for all $j\ge 0$
\begin{equation}\label{eq:induction-hypothesis-B-C}
C_j\cap\bigcup_{i=0}^nu_i(B_j)=\varnothing.
\end{equation}
Let us first show how this yields the proof of the claim: we set $B_\infty\mathrel{:=}\bigcup_{j=0}^\infty B_j$ and $C_\infty\mathrel{:=}\bigcup_{j=0}^\infty C_j$. Then both of these are infinite sets, and moreover each $u_i(B_\infty)$ misses $C_\infty$ by Condition $(\ref{eq:induction-hypothesis-B-C})$ and since the unions are increasing. As $B_\infty$ is infinite, we find an injection $c\colon\omega\setminus A\to\omega$ with image $B_\infty$. We claim that $\chi\mathrel{:=}c+\id_A$ has the desired properties: indeed, this is again an injection as $B_\infty\cap A=\varnothing$, and it is the identity on $A$ by construction. On the other hand,
\begin{equation*}
\bigcup_{i=0}^n \im(u_i\chi)=\underbrace{\bigcup_{i=0}^n u_i(A)}_{{}\mathrel{=:}A'}\cup\underbrace{\bigcup_{i=0}^n u_i(B_\infty)}_{{}\mathrel{=:}B'}
\end{equation*}
and $B'$ has infinite complement in $\omega$ (namely, at least $C_\infty$) whereas $A'$ is even finite; we conclude that also their union has infinite complement as desired.

It therefore only remains to construct the chains $B_0\subsetneq\cdots$ and $C_0\subsetneq\cdots$, for which we will proceed by induction. We begin by setting $B_0=C_0=\varnothing$ which obviously has all the required properties. Now assume we've already constructed the finite sets $B_j$ and $C_j$ satisfying $(\ref{eq:induction-hypothesis-B-C})$.

The set $A\cup B_j\cup\bigcup_{i=0}^n u_i^{-1}(C_j)$ is finite as $A,B_j,C_j$ are finite and each $u_i$ is injective, so we can pick a $b\in\omega$ not contained in it. We now set $B_{j+1}\mathrel{:=}B_j\cup\{b\}$, which is obviously finite and a proper superset of $B_j$ by construction. We moreover observe that $C_j$ misses all $u_i(B_{j+1})$ as it misses $u_i(B_j)$ by the induction hypothesis and moreover $u_i(b)\notin C_j$ for any $i$ by construction.

By the same argument we can pick $c\in\omega\setminus\left(C_j\cup \bigcup_{i=0}^n u_i(B_{j+1})\right)$ and define $C_{j+1}\mathrel{:=}C_j\cup\{c\}$; this is obviously again finite and a proper superset of $C_j$. We claim that Condition $(\ref{eq:induction-hypothesis-B-C})$ is satisfied for $B_{j+1}$ and $C_{j+1}$. Indeed, we have already seen that $C_j$ misses all of $u_i(B_{j+1})$. On the other hand, also $c\notin u_i(B_{j+1})$ for each $i$ by construction, verifying the condition. This finishes the inductive construction and hence the proof of the proposition.
\end{proof}
\end{prop}

\begin{prop}
Let $A\subset\omega$ be finite, and let $(u_0,\dots,u_n), (v_0,\dots,v_n)\in\mathcal M^{1+n}$ such that $u_i|_A=v_i|_A$ for $i=0,\dots,n$. Then $[u_0,\dots,u_n]=[v_0,\dots,v_n]$ in $\mathcal M^{1+n}/\mathcal M_A$.
\begin{proof}
Applying the above to the $2n+2$ injections $u_0,\dots,u_n,v_0,\dots,v_n$, we may assume without loss of generality that
\begin{equation*}
B\mathrel{:=}\omega\setminus\left(\bigcup_{i=0}^n\im u_i\cup\bigcup_{i=0}^n\im v_i\right)
\end{equation*}
is infinite. We can therefore choose an injection $\phi\colon\omega\setminus A\to\omega$ with image in $B$, and we moreover pick a bijection $\omega\setminus A\cong (\omega\setminus A)\amalg(\omega\setminus A)$, yielding two injections $j_1,j_2\colon\omega\setminus A\to\omega\setminus A$ whose images partition $\omega\setminus A$. We now define for $i=0,\dots,n$
\begin{equation*}
w_i(x)\mathrel{:=}\begin{cases}
u_i(x) & \textup{if }x\in A\\
u_i(y) & \textup{if }x = j_1(y)\textup{ for some $y\in\omega\setminus A$}\\
\phi(y) & \textup{if }x = j_2(y)\textup{ for some $y\in\omega\setminus A$.}
\end{cases}
\end{equation*}
This is indeed well-defined as $\omega$ is the disjoint union $A\sqcup \im(j_1)\sqcup\im(j_2)$ and since $j_1$ and $j_2$ are injective. We now claim that $w_i$ is injective (and hence an element of $\mathcal M$): indeed, assume $w_i(x)=w_i(x')$ for $x\not=x'$. Since $\im(u_i)$ is disjoint from $\im\phi$ and since $\phi$ is injective, we conclude that $x,x'\notin\im j_2$. As moreover $w_i|_A=u_i|_A$ and $w_ij_1=u_i|_{\omega\setminus A}$ are injective, we can assume without loss of generality that $x\in A$ and $x'=j_1(y')$ for some $y'\in\omega\setminus A$. But then $u_i(x)=w_i(x)=w_i(x')=u_i(y')$, which contradicts the injectivity of $u_i$ as $y'\notin A$ and hence in particular $y'\not=x$.

We now observe that by construction $w_i(\incl_A+j_1)=u_i$ and $w_i(\incl_A+j_2)=u_i|_A+\phi$. On the other hand, $\incl_A+j_1$ and $\incl_A+j_2$ are obviously injections fixing $A$ pointwise, so that they witness the equalities
\begin{equation*}
[u_0,\dots,u_n]=[w_0,\dots,w_n]=[u_0|_A+\phi,\dots,u_n|_A+\phi]
\end{equation*}
in $\mathcal M^{1+n}/\mathcal M_A$. Analogously, one shows that $[v_0,\dots,v_n]=[v_0|_A+\phi,\dots,v_n|_A+\phi]$, and as $v_i|_A=u_i|_A$ by assumption, this further agrees with $[u_0|_A+\phi,\dots,u_n|_A+\phi]$, finishing the proof.
\end{proof}
\end{prop}

\begin{cor}\label{cor:support-agree-EM}
Let $(u_0,\dots,u_n),(u_0',\dots,u_n')\in (E\mathcal M)_n$, let $X$ be an $E\mathcal M$-simplicial set, and let $x\in X_n$. Assume that $x$ is supported \emph{as an element of the $\mathcal M$-set $X_n$} on some finite set $A\subset\omega$ such that $u_i|_A=u_i'|_A$ for $i=0,\dots,n$.

Then $(u_0,\dots,u_n).x = (u_0',\dots,u_n').x$.
\begin{proof}
We begin with the special case that there exists an $\alpha\in\mathcal M_A$ such that $u_i'=u_i\alpha$ for all $i$. In this case
\begin{equation*}
(u_0',\dots,u_n').x=(u_0\alpha,\dots,u_n\alpha).x=(u_0,\dots,u_n).\alpha.x=(u_0,\dots,u_n).x
\end{equation*}
as desired, where the last step uses that $x$ is fixed by $\alpha\in\mathcal M_A$.

We conclude that $\mathcal M^{1+n}\to X_n$, $(u_0,\dots,u_n)\mapsto (u_0,\dots,u_n).x$ descends to $\mathcal M^{1+n}/\mathcal M_A$; the claim therefore follows from the previous proposition.
\end{proof}
\end{cor}

\begin{proof}[Proof of Theorem~\ref{thm:support-EM-vs-M}]\index{support!these notions agree}
Let $\chi\colon E\Inj(A,\omega)\times\Delta^n\to X$ be any $E\mathcal M$-equivariant map. We claim that the image of the $n$-simplex $(i,\dots,i;\id_{[n]})$, where $i$ denotes the inclusion $A\hookrightarrow\omega$, is supported on $A$. Indeed, $(i,\dots,i;\id_{[n]})$ is supported on $A$ by the argument from Example~\ref{ex:EM-Inj-support}, hence so is its image by Lemma~\ref{lemma:EM-support-preserve}.

On the other hand, let $x\in X_n$ be supported on $A$ \emph{with respect to the underlying $\mathcal M$-action}. We define $\chi_m\colon (E\Inj(A,\omega)\times\Delta^n)_m\to X_m$ as follows: we send a tuple $(u_0,\dots,u_m;f)$, where the $u_i$ are injections $A\to\omega$ and $f\colon[m]\to[n]$ is any map in $\Delta$, to $(\bar{u}_0,\dots,\bar{u}_m).f^*x$ where each $\bar{u}_i$ is an extension of $u_i$ to all of $\omega$, i.e.~to an element of $\mathcal M$. Such extensions can certainly be chosen as $A$ is finite, and we claim that this is in fact independent of the choice of extension: indeed, as $x$ is fixed by $\mathcal M_A$, so is $f^*x$, and hence this follows from the previous corollary. With this established, it is trivial to prove that the $\chi_m$ assemble into a simplicial map $E\Inj(A,\omega)\times\Delta^n\to X$ and that this is $E\mathcal M$-equivariant. Moreover, a possible extension of $i\colon A\hookrightarrow\omega$ is given by the identity of $\omega$ and hence we see that $\chi_m(i,\dots,i;\id_{[n]})=x$.

We conclude that if $x\in X_n$ is supported on $A$ with respect to the underlying $\mathcal M$-action, then it is obtained by evaluating some $E\mathcal M$-equviariant $\chi\colon E\Inj(A,\omega)\times\Delta^n\to X$ at the canonical element (by the second paragraph), and hence it is actually supported on $A$ with respect to the $E\mathcal M$-action (by the first one). As the converse holds for trivial reasons, we conclude that the two notions of `being supported on $A$' indeed agree. It then follows from the above that evaluation at the canonical element defines a surjection
\begin{equation}\label{eq:corepr-supported-on-A}
\Maps(E\Inj(A,\omega),X)\to X_{[A]}.
\end{equation}

As we have seen in Lemma~\ref{lemma:EM-support-preserve}, $X\mapsto X_{[A]}$ defines a subfunctor of the forgetful functor $\cat{$\bm{E\mathcal M}$-SSet}\to\cat{SSet}$, and with respect to this $(\ref{eq:corepr-supported-on-A})$ is obviously natural. To finish the proof of the claimed corepresentability result, it therefore suffices that $(\ref{eq:corepr-supported-on-A})$ is also injective. For this we let $\chi,\chi'\colon E\Inj(A,\omega)\times\Delta^n\to X$ with $\chi(i,\dots,i;\id_{[n]})=\chi'(i,\dots,i;\id_{[n]})$. Then we have for any $(u_0,\dots,u_m)\in\mathcal M^{1+m}$ and $f\colon[m]\to[n]$ in $\Delta$
\begin{align*}
\chi(u_0,\dots,u_m;f)&=(\bar{u}_0,\dots,\bar{u}_m).\chi(i,\dots,i;f)=(\bar{u}_0,\dots,\bar{u}_m).f^*\chi(i,\dots,i;\id_{[n]})\\
&=(\bar{u}_0,\dots,\bar{u}_m).f^*\chi'(i,\dots,i;\id_{[n]})=\chi'(u_0,\dots,u_m;f)
\end{align*}
where again $\bar{u}_i$ is any extension of $u_i$ to an element of $\mathcal M$. This finishes the proof of corepresentability and hence of the theorem.
\end{proof}

The above theorem has the following computational consequence:

\begin{cor}\label{cor:E-Inj-corepr}
For any simplicial set $K$, any finite set $A$, and any $G$-$\Sigma_A$-biset $X$, the map
\begin{equation*}
E\mathcal M\times_{\mathcal M}\big((\Inj(A,\omega)\times_{\Sigma_A}X)\times K\big)\to(E\Inj(A,\omega)\times_{\Sigma_A}X)\times K
\end{equation*}
adjoint to the product of the inclusion of the $0$-simplices with the identity of $K$ is an isomorphism of $E\mathcal M$-$G$-simplicial sets.
\begin{proof}
As $E\mathcal M\times_{\mathcal M}\blank$ is a simplicial left adjoint, we are reduced to proving that the corresponding map $E\mathcal M\times_{\mathcal M}\Inj(A,\omega)\to E\Inj(A,\omega)$ is an isomorphism (it is obviously left-$\mathcal M$-right-$\Sigma_A$-equivariant). For this we observe that by the above theorem $E\Inj(A,\omega)$ corepresents $(\blank)_{[A]}$ by evaluating at the inclusion $\iota\colon A\hookrightarrow\omega$. On the other hand, it is obvious that $\Inj(A,\omega)$ corepresents $(\blank)_{[A]}\colon\cat{$\bm{\mathcal M}$-SSet}\to\cat{SSet}$ by evaluating at the same element, also see \cite[Example~2.9]{I-vs-M-1-cat}. By adjointness, $E\mathcal M\times_{\mathcal M}\Inj(A,\omega)$ therefore corepresents $(\blank)_{[A]}\circ\forget$ via evaluating at $[1;\iota]$. By the previous theorem this agrees with $(\blank)_{[A]}$, and as the above map sends $[1;\iota]$ to $\iota$, this completes the proof.
\end{proof}
\end{cor}

We now very easily prove the $E\mathcal M$-analogue of Lemma~\ref{lemma:support-vs-action-M}:

\begin{lemma}\label{lemma:support-vs-action-EM}
Let $(u_0,\dots,u_n)\in (E\mathcal M)_n$, let $A\subset\omega$ be finite, and let $X$ be any $E\mathcal M$-simplicial set. Then the composition
\begin{equation}\label{eq:EM-support-restricted-action}
\Delta^n\times X_{[A]}\xrightarrow{(u_0,\dots,u_n)\times\incl}E\mathcal M\times X\xrightarrow{\textup{act}} X
\end{equation}
has image in $X_{[u_0(A)\cup\cdots\cup u_n(A)]}$; in particular, $X^\tau$ is an $E\mathcal M$-simplicial subset of $X$.
\begin{proof}
It suffices to prove the first statement. For this we observe that any simplex in the image of $(\ref{eq:EM-support-restricted-action})$ can be written as $\big(f^*(u_0,\dots,u_n)\big).x$ for some $m$-simplex $x$ of $X_{[A]}$ and some $f\colon[m]\to[n]$ in $\Delta$. We now calculate for any $v\in\mathcal M$
\begin{equation}\label{eq:double-action}
v.(f^*(u_0,\dots,u_n).x)=v.(u_{f(0)},\dots,u_{f(m)}).x=(vu_{f(0)},\dots,vu_{f(n)}).x.
\end{equation}
If now $v$ is the identity on $u_0(A)\cup\cdots\cup u_n(A)$, then $vu_{f(i)}$ and $u_{f(i)}$ agree on $A$ for all $i=0,\dots,m$. Hence, if $x$ is supported on $A$, then
\begin{equation*}
(vu_{f(0)},\dots, vu_{f(m)}).x=(u_{f(0)},\dots,u_{f(m)}).x=f^*(u_0,\dots,u_n).x
\end{equation*}
by Corollary~\ref{cor:support-agree-EM}. Together with $(\ref{eq:double-action})$ this precisely yields the claim.
\end{proof}
\end{lemma}

\begin{cor}\label{cor:EM-tau-colim-lim}
The full simplicial subcategory $\cat{$\bm{E\mathcal M}$-SSet}^\tau\subset\cat{$\bm{E\mathcal M}$-SSet}$ is complete and cocomplete, and it is closed under all small colimits and finite limits.
\begin{proof}
As in Remark~\ref{rk:tame-M-colimits}, the previous lemma provides a right adjoint $(\blank)^\tau$ to the inclusion $\cat{$\bm{E\mathcal M}$-SSet}^\tau\hookrightarrow\cat{$\bm{E\mathcal M}$-SSet}$, which shows that $\cat{$\bm{E\mathcal M}$-SSet}^\tau$ is complete and cocomplete with colimits created in $\cat{$\bm{E\mathcal M}$-SSet}$. As the forgetful functor $\cat{$\bm{E\mathcal M}$-SSet}\to\cat{$\bm{\mathcal M}$-SSet}$ creates limits and preserves and reflects tameness by Theorem~\ref{thm:support-EM-vs-M}, the identification of finite limits follows from Remark~\ref{rk:tame-M-colimits}.
\end{proof}
\end{cor}

\begin{rk}
Of course, we could have just turned Theorem~\ref{thm:support-EM-vs-M} into the definition instead---however, this does not buy as anything, as we  then instead would have had to do all of the above work in order to prove Corollary~\ref{cor:support-agree-EM} (or equivalently the corepresentability statement), which will be crucial below. Moreover, defining support by means of the mere $\mathcal M$-action is `evil' from a homotopical viewpoint, as we \emph{a priori} throw away a lot of higher information: for example, if $u\in\mathcal M$ and $x$ is a vertex of an $E\mathcal M$-simplicial set, then we can think of the edge $(u,1).(s^*x)$, where $s\colon[1]\to[0]$ is the unique morphism in $\Delta$, as providing a natural comparison between $x$ and $u.x$. From the homotopical viewpoint we should then always be interested in this edge itself and not only in its endpoints.
\end{rk}

A hindrance to understanding the $G$-global weak equivalences of general $\mathcal M$-simplicial sets was that it is not clear whether $E\mathcal M\times_{\mathcal M}\blank$ is fully homotopical. Using the above theory, we can prove the following comparison result, which in particular tells us that this issue goes away when restricting to tame actions:

\begin{thm}\label{thm:tame-M-sset-vs-EM-sset}
The simplicial adjunction $E\mathcal M\times_{\mathcal M}\blank\dashv\forget$ restricts to
\begin{equation}\label{eq:EM-adjunction-tame}
E\mathcal M\times_{\mathcal M}\blank\colon\cat{$\bm{\mathcal M}$-$\bm G$-SSet}^\tau\rightleftarrows\cat{$\bm{E\mathcal M}$-$\bm G$-SSet}^\tau :\!\forget.
\end{equation}
Both functors in $(\ref{eq:EM-adjunction-tame})$ preserve and reflect $G$-global weak equivalences, and they descend to mutually inverse equivalences on associated quasi-categories.
\begin{proof}
The forgetful functor obviously restricts to the full subcategories of tame objects. To see that also $E\mathcal M\times_{\mathcal M}\blank$ restricts accordingly, we appeal to Theorem~\ref{thm:structure-tame-M-G-SSet}: as the tame $E\mathcal M$-$G$-simplicial sets are closed under all colimits, we are reduced (Lemma~\ref{lemma:saturated-objects}) to showing that $E\mathcal M\times_{\mathcal M}\big((\Inj(A,\omega)\times_{\Sigma_A}X)\times K\big)$ is tame for every $G$-$\Sigma_A$-biset $X$ and any simplicial set $K$, which follows easily from Corollary~\ref{cor:E-Inj-corepr}.

The forgetful functor creates weak equivalences (even without the tameness assumption) as it is homotopical and part of a Quillen equivalence by Corollary~\ref{cor:em-vs-m-equiv-model-cat}.

Let us now prove that the unit $\eta\colon Y\to\forget E\mathcal M\times_{\mathcal M}Y$ is a $G$-global weak equivalence for every $Y\in\cat{$\bm{\mathcal M}$-$\bm G$-SSet}^\tau$. We caution the reader that this is not just a formal consequence of Corollary~\ref{cor:em-vs-m-equiv-model-cat} because we are \emph{not} deriving $E\mathcal M\times_{\mathcal M}\blank$ in any way here. Instead, we will use Theorem~\ref{thm:hq-M-semistable-replacement} together with the full generality of Theorem~\ref{thm:hq-M-computation}:

Both $E\mathcal M\times_{\mathcal M}\blank$ as well as $\forget$ are left adjoints and hence cocontinuous. As both functors preserve tensors, the composition $\forget(E\mathcal M\times_{\mathcal M}\blank)$ sends the maps in $I_{\textup{tame}}$ to underlying cofibrations. Corollary~\ref{cor:saturated-trafo} therefore reduces this to the case that $Y=(\Inj(A,\omega)\times_{\Sigma_A}X)\times K$ for $H,A,X$ as above and $K$ any simplicial set. Again using that both $E\mathcal M\times_{\mathcal M}\blank$ and $\forget$ preserve tensors (and that the above is a simplicial adjunction) we reduce further to the case that $Y=\Inj(A,\omega)\times_{\Sigma_A}X$. After postcomposing with the isomorphism from Corollary~\ref{cor:E-Inj-corepr}, the unit simply becomes the inclusion of the $0$-simplices, and this factors as
\begin{equation*}
\Inj(A,\omega)\times_{\Sigma_A}X\xrightarrow{\pi}\big(\Inj(A,\omega)\times_{\Sigma_A}X\big)\hq\mathcal M\to E\Inj(A,\omega)\times_{\Sigma_A}X,
\end{equation*}
where the right hand map is the $G$-universal weak equivalence from Theorem~\ref{thm:hq-M-computation}. As the left hand map is moreover a $G$-global weak equivalence by Theorem~\ref{thm:hq-M-semistable-replacement}, thus so is the unit.

Now let $f\colon X\to Y$ be any map in $\cat{$\bm{\mathcal M}$-$\bm G$-SSet}^\tau$. In the naturality square
\begin{equation*}
\begin{tikzcd}[column sep=1in]
X\arrow[d, "\eta"']\arrow[r, "f"] & Y\arrow[d,"\eta"]\\
\forget E\mathcal M\times_{\mathcal M} X\arrow[r, "\forget E\mathcal M\times_{\mathcal M}f"'] & \forget E\mathcal M\times_{\mathcal M}Y
\end{tikzcd}
\end{equation*}
both vertical maps are $G$-global weak equivalences by the above. Thus, if $f$ is a $G$-global weak equivalence, then so is the lower horizontal map. But $\forget$ reflects these, hence also $E\mathcal M\times_{\mathcal M}f$ is a $G$-global weak equivalence, proving that $E\mathcal M\times_{\mathcal M}\blank$ is homotopical. Conversely, if $E\mathcal M\times_{\mathcal M}f$ is a $G$-global weak equivalence, then so is $\forget E\mathcal M\times_{\mathcal M}f$ and hence also $f$ by the above square, i.e.~$E\mathcal M\times_{\mathcal M}\blank$ also reflects $G$-global weak equivalences.

Finally, if $X\in\cat{$\bm{E\mathcal M}$-$\bm G$-SSet}^\tau$ is arbitrary, then the triangle identity for adjunctions shows that $\forget\epsilon_X\colon\forget E\mathcal M\times_{\mathcal M}(\forget X)\to\forget X$ is right inverse to $\eta_{\forget X}$, hence a $G$-global weak equivalence. As $\forget$ reflects these, we conclude that also $\epsilon$ is levelwise a $G$-global weak equivalence, proving that the functors in $(\ref{eq:EM-adjunction-tame})$ induce mutually inverse equivalences of quasi-categories.
\end{proof}
\end{thm}

\subsection{The Taming of the Shrew} As promised, we can now finally prove that also tame $E\mathcal M$-$G$-simplicial sets and tame $\mathcal M$-$G$-simplicial sets are models of $G$-global homotopy theory. At this point, we will only consider them as categories with weak equivalences whereas suitable $G$-global model structures are the subject of Subsection~\ref{subsec:tame-model-structures}.

\begin{thm}\label{thm:shrew}
The inclusions
\begin{equation*}
\cat{$\bm{E\mathcal M}$-$\bm G$-SSet}^\tau\hookrightarrow\cat{$\bm{E\mathcal M}$-$\bm G$-SSet}
\qquad\text{and}\qquad
\cat{$\bm{\mathcal M}$-$\bm G$-SSet}^\tau\hookrightarrow\cat{$\bm{\mathcal M}$-$\bm G$-SSet}
\end{equation*}
are homotopy equivalences with respect to the $G$-global weak equivalences.
\end{thm}

For the proof we will need:

\begin{prop}\label{prop:restriction-to-faithful}
Let $H\subset\mathcal M$ be a universal subgroup, let $A\subset\omega$ be a faithful $H$-subset, and let $\phi\colon H\to G$ be any homomorphism. Then:
\begin{enumerate}
\item The restriction $r\colon E\mathcal M\times_\phi G\to E\Inj(A,\omega)\times_\phi G$ is a $G$-global weak equivalence in $\cat{$\bm{E\mathcal M}$-$\bm G$-SSet}$.
\item The restriction $r\colon \mathcal M\times_\phi G\to \Inj(A,\omega)\times_\phi G$ is a $G$-global weak equivalence in $\cat{$\bm{\mathcal M}$-$\bm G$-SSet}$.
\end{enumerate}
\begin{proof}
For the first statement, we note that the restriction $E\mathcal M\to E\Inj(A,\omega)$ is an $H$-global weak equivalence (with respect to $H$ acting on $\mathcal M$ and $A$ via its tautological action on $\omega$) by Example~\ref{ex:G-globally-contractible}. On the other hand, the $H$-action on both sides is free, so Corollary~\ref{cor:free-quotient-EM} implies that its image under $\phi_!$ is a $G$-global weak equivalence. But this is clearly conjugate to $r$, finishing the proof of the first statement.

For the second statement, we note that the inclusion $\mathcal M\times_\phi G\hookrightarrow E\mathcal M\times_\phi G$ is a $G$-global weak equivalence by Theorem~\ref{thm:em-vs-m-equiv}, and so is $\Inj(A,\omega)\times_\phi G\hookrightarrow E\Inj(A,\omega)\times_\phi G$ by the proof of Theorem~\ref{thm:tame-M-sset-vs-EM-sset}. The claim therefore follows from the first statement together with $2$-out-of-$3$.
\end{proof}
\end{prop}

\begin{cor}\label{cor:injective-almost-tame}
Let $H\subset\mathcal M$ be a universal subgroup, and let $A\subset\omega$ be a finite faithful $H$-subset.
\begin{enumerate}
\item Let $X$ be fibrant in the injective $G$-global model structure on $\cat{$\bm{E\mathcal M}$-$\bm G$-SSet}$. Then the $H$-action on $X$ restricts to $X_{[A]}$, and $X_{[A]}\hookrightarrow X$ is a $\mathcal G_{H,G}$-weak equivalence. In particular, $X^\tau\hookrightarrow X$ is a $G$-global weak equivalence.
\item Let $Y$ be fibrant in the injective $G$-global model structure on $\cat{$\bm{\mathcal M}$-$\bm G$-SSet}$. Then the $H$-action on $Y$ restricts to $Y_{[A]}$, and $Y_{[A]}\hookrightarrow Y$ is a $\mathcal G_{H,G}$-weak equivalence. In particular, $Y^\tau\hookrightarrow Y$ is a $G$-universal weak equivalence.
\end{enumerate}
\begin{proof}
We will only prove the claims for $E\mathcal M$-actions, the proof of the other ones being analogous.

Let $\phi\colon H\to G$ be any homomorphism. Using Theorem~\ref{thm:support-EM-vs-M}, we see that the inclusion $X_{[A]}^\phi\hookrightarrow X^\phi$ agrees up to conjugation by isomorphisms with
\begin{equation}\label{eq:inclusion-corep}
r^*\colon\Maps^{E\mathcal M\times G}(E\Inj(A,\omega)\times_\phi G,X)\to\Maps^{E\mathcal M\times G}(E\mathcal M\times_\phi G,X),
\end{equation}
where $r$ is as in the previous proposition. As $r$ is a $G$-global weak equivalence between injectively cofibrant $E\mathcal M$-$G$-simplicial sets, and since the injective $G$-global model structure is simplicial by Proposition~\ref{prop:equivariant-injective-model-structure}, $(\ref{eq:inclusion-corep})$ is a weak homotopy equivalence, proving the first statement.

For the second statement we consider the commutative diagram
\begin{equation*}
\begin{tikzcd}
\colim_{A} X_{[A]}\arrow[r]\arrow[d] & \colim_{A}X\arrow[d]\\
X^\tau\arrow[r,hook] & X,
\end{tikzcd}
\end{equation*}
where the colimits run over the filtered poset of finite faithful $H$-subsets $A\subset\omega$, and all maps are induced by the inclusions. The right hand vertical arrow is an isomorphism as the indexing category is filtered, and so is the left hand vertical map since a simplex supported on some finite set $B$ is also supported on the finite faithful $H$-subset $HB\cup F$, where $F\subset\omega$ is any chosen free $H$-orbit (which exists by universality). The claim therefore follows from the first statement as $\mathcal G_{H,G}$-weak equivalences are closed under filtered colimits.
\end{proof}
\end{cor}

\begin{proof}[Proof of Theorem~\ref{thm:shrew}]
Again, we will only prove the first statement. For this, we factor the inclusion through the full subcategory $\cat{$\bm{E\mathcal M}$-$\bm G$-SSet}^{w\tau}$ of those $E\mathcal M$-$G$-simplicial sets $X$ for which $X^\tau\hookrightarrow X$ is a $G$-global weak equivalence. It suffices to prove that both intermediate inclusions are homotopy equivalences.

Indeed, a homotopy inverse to $\cat{$\bm{E\mathcal M}$-$\bm G$-SSet}^{\tau}\hookrightarrow\cat{$\bm{E\mathcal M}$-$\bm G$-SSet}^{w\tau}$ is obviously given by $(\blank)^\tau$, whereas the previous corollary implies that taking functorial fibrant replacements in the injective $G$-global model structure provides a homotopy inverse to the remaining inclusion.
\end{proof}

\section[$G$-global homotopy theory via diagram spaces]{\texorpdfstring{\except{toc}{$\bm G$}\for{toc}{$G$}}{G}-global homotopy theory via diagram spaces}\label{sec:diag-spaces}
We will now consider models of unstable $G$-global homotopy theory in terms of more general `diagram spaces,' i.e.~functors from suitable indexing categories to $\cat{SSet}$. These models will in particular be useful in Chapter~\ref{chapter:stable} to connect unstable and stable $G$-global homotopy theory.

\begin{defi}
We write $I$\nomenclature[aI1]{$I$}{category of finite sets and injections} for the category of finite sets and injections, and we write $\mathcal I$\nomenclature[aI2]{$\mathcal I$}{simplicially enriched category built from $I$} for the simplicial category obtained by applying $E$ to the hom sets. An \emph{$I$-simplicial set} is a functor $I\to\cat{SSet}$, and we write $\cat{$\bm I$-SSet}$ for the simplicially enriched functor category $\FUN(I,\cat{SSet})$. An $\emph{$\mathcal I$-simplicial set}$ is a simplicially enriched functor $\mathcal I\to\cat{SSet}$, and we write $\cat{$\bm{\mathcal I}$-SSet}\mathrel{:=}\FUN(\mathcal I,\cat{SSet})$.
\end{defi}

In the literature, the category $I$ is also denoted by $\mathbb I$ \cite{lind} and unfortunately also by $\mathcal I$ \cite{sagave-schlichtkrull}. Sagave and Schlichtkrull proved that $\cat{$\bm I$-SSet}$ models ordinary homotopy theory, see~\cite[Theorem~3.3]{sagave-schlichtkrull}.

We can also view $I$ as a discrete analogue of the topological category $L$ used in Schwede's model of unstable global homotopy theory in terms of $L$-spaces, see \cite[Sections~1.1--1.2]{schwede-book}. In an earlier version, Schwede also sketched that $I$-spaces model global homotopy theory (with respect to finite groups), also cf.~\cite[Section~6.1]{hausmann-global}, for which we will give a full proof as Theorem~\ref{thm:L-vs-I}.

\subsection{Model structures} Next, we will introduce $G$-global model structures on $G$-$I$- and $G$-$\mathcal I$-simplicial sets (i.e.~$G$-objects in $\cat{$\bm I$-SSet}$ or $\cat{$\bm{\mathcal I}$-SSet}$, respectively), and prove that they are equivalent to the models from the previous sections. As usual in this context, we begin by constructing a suitable level model structure that we will then later Bousfield localize at the desired weak equivalences:

\begin{defi}
A map $f\colon X\to Y$ in $\cat{$\bm{G}$-$\bm{I}$-SSet}$ is called a \emph{strict level weak equivalence}\index{strict level weak equivalence!in G-I-SSet@in $\cat{$\bm G$-$\bm{I}$-SSet}$|textbf} or \emph{strict level fibration} if $f(A)\colon X(A)\to Y(A)$ is a $\mathcal G_{\Sigma_A,G}$-weak equivalence or fibration, respectively, for every finite set $A$. Put differently, for each finite group $H$ acting faithfully on $A$ and each homomorphism $\phi\colon H\to G$ the induced map $X(A)^\phi\to Y(A)^\phi$ is a weak equivalence or fibration, respectively.

A map of $G$-$\mathcal I$-simplicial sets is called a strict level weak equivalence or fibration if it is so when viewed as a map in $\cat{$\bm{G}$-$\bm{I}$-SSet}$.
\end{defi}

\index{G-global model structure@$G$-global model structure!strict level|seeonly{strict level model structure}}
\index{strict level model structure!on G-I-SSet@on $\cat{$\bm G$-$\bm I$-SSet}$|textbf}
\begin{prop}\label{prop:strict-model-structure}
The strict level weak equivalences and fibrations are part of a unique model structure on $\cat{$\bm{G}$-$\bm{I}$-SSet}$, which we call the the \emph{strict level model structure}. It is proper, simplicial, combinatorial, and filtered colimits in it are homotopical. A possible set of generating cofibrations is given by the maps
\begin{equation*}
(I(A,\blank)\times_\phi G)\times\del\Delta^n\hookrightarrow(I(A,\blank)\times_\phi G)\times\Delta^n
\end{equation*}
for $n\ge 0$ and $H$, $A$, and $\phi$ as above, and a set of generating acyclic cofibrations is likewise given by the maps
\begin{equation*}
(I(A,\blank)\times_\phi G)\times\Lambda^n_k\hookrightarrow(I(A,\blank)\times_\phi G)\times\Delta^n
\end{equation*}
with $0\le k\le n$.

\index{strict level model structure!on G-II-SSet@on $\cat{$\bm G$-$\bm{\mathcal I}$-SSet}$|textbf}
Similarly, the strict level weak equivalences and fibrations are part of a unique model structure on $\cat{$\bm{G}$-$\bm{\mathcal I}$-SSet}$, which has the same properties; we again call it the \emph{strict level model structure}. A possible set of generating cofibrations is given by the maps
\begin{equation*}
(\mathcal I(A,\blank)\times_\phi G)\times\del\Delta^n\hookrightarrow(\mathcal I(A,\blank)\times_\phi G)\times\Delta^n
\end{equation*}
for $n\ge 0$ and $H$, $A$, and $\phi$ as above, and a set of generating acyclic cofibrations is likewise given by the maps
\begin{equation*}
(\mathcal I(A,\blank)\times_\phi G)\times\Lambda^n_k\hookrightarrow(\mathcal I(A,\blank)\times_\phi G)\times\Delta^n
\end{equation*}
with $0\le k\le n$.

Finally, the forgetful functor is part of a simplicial Quillen adjunction\nomenclature[aII]{$\mathcal I\times_I\blank$}{left adjoint to forgetful functor $\cat{$\bm{\mathcal I}$-SSet}\to\cat{$\bm I$-SSet}$}
\begin{equation}\label{eq:I-vs-script-I-restricted}
\mathcal I\times_{I}\blank\colon\cat{$\bm{G}$-$\bm{I}$-SSet}_{\textup{strict level}}\rightleftarrows\cat{$\bm{G}$-$\bm{\mathcal I}$-SSet}_{\textup{strict level}} :\!\forget.
\end{equation}
\end{prop}

\begin{rk}
To be entirely precise, the above generating (acyclic) cofibrations of course do not form a set since there are too many finite groups $H$ and finite faithful $H$-sets $A$. However, this issue is easily solved by restricting to a set of finite groups $H$ and sets of finite $H$-sets such that these cover all isomorphism classes, and we will tacitly do so. Similar caveats apply to several of the other model structures considered below.
\end{rk}

\subsubsection{Generalized projective model structures}
We will obtain these model structures as an instance of a more general construction from \cite{schwede-book} of `generalized projective model structures' for suitable indexing categories. For this we will need the following terminology:

\begin{defi}\index{dimension function|seealso{generalized projective model structure, for categories with dimension function}}
\index{dimension function|textbf}
Let $\mathscr C$ be a complete and cocomplete closed symmetric monoidal category. We say that a $\mathscr C$-enriched category $\mathscr I$ \emph{has a dimension function} if there exists a function $\dim\colon\Ob(\mathscr I)\to\mathbb N$ such that
\begin{enumerate}
\item $\Hom(d,e)$ is initial in $\mathscr C$ whenever $\dim(e)<\dim(d)$.
\item If $\dim(d)=\dim(e)$, then $d\cong e$.
\end{enumerate}
\end{defi}

\begin{ex}
Both $I$ and $\mathcal I$ have dimension functions; a canonical choice is the function sending a finite set $A$ to its cardinality $|A|$. Moreover, if $G$ is any discrete group, then composing with the projection $I\times BG\to I$ or $\mathcal I\times BG\to\mathcal I$ yields a dimension function on $I\times BG$ or $\mathcal I\times BG$, respectively.
\end{ex}

If now $\mathscr I$ is a $\mathscr C$-enriched category, then we consider the enriched functor category $\FUN(\mathscr I,\mathscr C)$ where $\mathscr C$ is viewed as enriched over itself via the closed symmetric monoidal structure. If $\mathscr C$ carries a suitably nice model structure, then we can equip $\FUN(\mathscr I, \mathscr C)$ with the {projective model structure}, in which a map $f$ is a weak equivalence or fibration if and only if $f(A)$ is a weak equivalence or fibration, respectively, in $\mathscr C$ for every $A\in\mathscr I$. The following proposition refines this to the situation, where we instead are given individual model structures on the enriched functor categories $\End(A)\text{--}\mathscr C\mathrel{:=}\FUN(B\End(A),\mathscr C)$ of $\End(A)$-objects in $\mathscr C$ for all $A\in\mathscr I$; here $B\End(A)$ denotes the $\mathscr C$-enriched category with a single object whose endomorphisms are given by $\End(A)$ or, equivalently, the full subcategory of $\mathscr I$ spanned by $A$.

\begin{prop}\label{prop:generalized-projective-dim}\index{generalized projective model structure!for categories with dimension function|textbf}
Let $\mathscr I$ be a $\mathscr C$-enriched category with dimension function $\dim$, and assume we are given for each $A\in\mathscr I$ a model structure on the category $\End(A)\textup{--}\mathscr C$ such that the following `consistency condition' holds: if $\dim(A)\le\dim(B)$, then any pushout of $\Hom(A,B)\otimes_{\End(A)}i$ is a weak equivalence in the given model structure on $\End(B)\textup{--}\mathscr C$ for any acyclic cofibration $i$ in $\End(A)\textup{--}\mathscr C$.

Then there exists a unique model structure on $\FUN(\mathscr I,\mathscr C)$ such that a map $f\colon X\to Y$ is a weak equivalence or fibration if and only if $f(A)\colon X(A)\to Y(A)$ is a weak equivalence or fibration, respectively, in $\End(A)\textup{--}\mathscr C$ for each $A\in\mathscr I$.

Moreover, if each $\End(A)\textup{--}\mathscr C$ is cofibrantly generated with set of generating cofibrations $I_A$ and set of generating acyclic cofibrations $J_A$, then the resulting model structure is cofibrantly generated with set of generating cofibrations
\begin{equation*}
\{\Hom(A,\blank)\otimes_{\End(A)}i\colon A\in \mathscr I,i\in I_A\},
\end{equation*}
and generating acyclic cofibrations
\begin{equation*}
\{\Hom(A,\blank)\otimes_{\End(A)}j\colon A\in \mathscr I,j\in J_A\}.
\end{equation*}
\begin{proof}
See \cite[Proposition~C.23]{schwede-book}.
\end{proof}
\end{prop}

Here, $\Hom(A,\blank)\otimes_{\End(A)}\blank$ is the left adjoint of the functor $\FUN(\mathscr I,\cat{SSet})\to\End(A)\textup{--}\mathscr C$ given by evaluation at $A$. For $\mathscr C=\cat{SSet}$, we have a concrete model of $\Hom(A,B)\otimes_{\End(A)}X$ as the \emph{balanced product} $\Hom(A,B)\times_{\End(A)}X$, i.e.~the quotient of the ordinary product by the equivalence relation generated in each simplicial degree by $(f\sigma,x)\sim (f,X(\sigma)(x))$.

We can also chacterize the cofibrations of the above model structure in analogy with the usual characterization in Reedy model structures. For this we need the following notion from \cite[Construction C.13~and Definition~C.15]{schwede-book}:

\begin{constr}\label{constr:latching-general}
\index{latching object|textbf}\index{latching map|textbf}
For any $m\ge 0$, we define $\mathscr I^{{}\le m}\subset\mathscr I$ as the full subcategory of those objects $A$ satisfying $\dim(A)\le m$. Then restriction along $i_m\colon\mathscr I^{{}\le m}\hookrightarrow\mathscr I$ admits a left adjoint $i_{m!}\colon\Fun(\mathscr I^{{}\le m},\mathscr C)\to\Fun(\mathscr I,\mathscr C)$ via enriched left Kan extension, and we write $\epsilon_m$ for the counit of the adjunction $i_{m!}\dashv i_m^*$. We now define the \emph{$A^{\text{th}}$ latching object} of $X\colon\mathscr I\to\mathscr C$ via $L_A(X)\mathrel{:=}(i_{(\dim(A)-1)!}i_{\dim(A)-1}^*X)(A)$ and the \emph{$A^{\text{th}}$ latching map} $\ell_A\colon L_A(X)\to X(A)$ as $\epsilon_{\dim(A)-1}(A)$.
\end{constr}

\begin{rk}\index{latching category|textbf}\index{latching object|textbf}\index{latching map|textbf}
If $\mathscr I$ is an ordinary category with dimension function, then the \emph{latching category} $\del(\mathscr I\downarrow A)$ is defined for $A\in\mathscr I$ as the full subcategory of the slice $\mathscr I\downarrow A$ on all objects $B\to A$ that are not isomorphisms. In analogy with the usual terminology in Reedy categories, we can then define the \emph{$A^{\text{th}}$ latching object}\nomenclature[aLA]{$L_A$}{$A^{\text{th}}$ latching object} as
\begin{equation*}
L_A(X)\mathrel{:=}\colim\limits_{B\to A\in \del(\mathscr I\downarrow A)} X(B) \in\End(A)\textup{--}\mathscr C,
\end{equation*}
and the maps $X(B)\to X(A)$ induced by the given maps $B\to A$ assemble via the universal property of colimits into a map $\ell_A\colon L_A(X)\to X(A)$ that we again call the \emph{$A^{\text{th}}$ latching map}.\nomenclature[alA]{$\ell_A$}{$A^{\text{th}}$ latching map} Note that this is indeed a special case of the above general definition by the usual pointwise formula for left Kan extension.
\end{rk}

\begin{ex}
While we view $BG\times I$ as a simplicially enriched category, the category $\FUN(BG\times I,\mathscr C)$ of enriched functors is isomorphic to the usual functor category, so the previous remark applies to this setting to express the latching objects as colimits over the latching categories.

However, we can give an even simpler description in this case: namely, the evident inclusion of the poset $\{B\subsetneq A\}$ into $\del(BG\times I\downarrow A)$ is an equivalence for any $A\in I$, in particular cofinal. Thus, we can describe the $A^{\textup{th}}$ latching object for $X\in\cat{$\bm G$-$\bm I$-SSet}$ as $\colim_{B\subsetneq A}X(B)$ with the induced $(\Sigma_A\times G)$-action, and the latching map is again induced by the inclusions $B\hookrightarrow A$.
\end{ex}

\begin{prop}\label{prop:characterization-cofibration-general}\index{generalized projective model structure!for categories with dimension function!cofibrations|textbf}
A map $f\colon X\to Y$ is a cofibration in the model structure of Proposition~\ref{prop:generalized-projective-dim} if and only if the map
\begin{equation*}
X(A)\amalg_{L_A(X)}L_A(Y)\xrightarrow{(f(A),\ell_A)} Y(A)
\end{equation*}
(where the pushout is taken over the maps $\ell_A\colon L_A(X)\to X(A)$ and $L_A(f)$) is a cofibration in the given model structure on $\End(A)\textup{--}\mathscr C$ for all $A\in\mathscr I$.
\begin{proof}
This is part of \cite[Proposition~C.23]{schwede-book}
\end{proof}
\end{prop}

\subsubsection{Strict level model structures} Using the above criterion we now get:

\begin{proof}[Proof of Proposition~\ref{prop:strict-model-structure}]
We will only construct the model structure on $\cat{$\bm{G}$-$\bm{I}$-SSet}$, the construction for \cat{$\bm{G}$-$\bm{\mathcal I}$-SSet} being analogous. For this we want to appeal to Proposition~\ref{prop:generalized-projective-dim} (for $\mathscr C=\cat{SSet}$, $\mathscr I=BG\times I$), so we have to check the consistency condition. To this end we claim that in the above situation, the functor $(BG\times I)(A,B)\times_{G\times\Sigma_A}\blank\cong I(A,B)\times_{\Sigma_A}\blank$ sends acyclic cofibrations of the usual $\mathcal G_{\Sigma_A,G}$-model structure to acyclic cofibrations in the \emph{injective} $\mathcal G_{\Sigma_B,G}$-model structure. But indeed, by cocontinuity it suffices to check this on generating acyclic cofibrations, where this is obvious.

This already shows that the model structure on $\cat{$\bm G$-$\bm I$-SSet}$ exists, and that it is cofibrantly generated with generating (acyclic) cofibrations as claimed above. As $\cat{$\bm G$-$\bm I$-SSet}$ is locally presentable, we conclude that it is in fact combinatorial.

The category $\cat{$\bm G$-$\bm I$-SSet}$ is enriched, tensored, and cotensored over $\cat{SSet}$ in the obvious way, and as the $\mathcal G_{\Sigma_A,G}$-model structures are simplicial and since pullbacks, cotensors, and (acyclic) fibrations are defined levelwise, also $\cat{$\bm G$-$\bm I$-SSet}$ is simplicial. Similarly one proves right properness and the preservation of weak equivalences under filtered colimits.

Moreover, the forgetful functor admits a simplicial left adjoint (via simplicially left Kan extension along $I\to\mathcal I$), which we denote by $\mathcal I\times_I\blank$; explicitly, we can arrange that $\mathcal I\times_I(I(A,\blank)\times K)=\mathcal I(A,\blank)\times K$ for all $A\in I$ and $K\in\cat{$\bm G$-SSet}$ with the evident functoriality, and this in turn describes $\mathcal I\times_I\blank$ up to canonical (simplicial) isomorphism. It is then obvious from the definition that $\forget$ is right Quillen, so that $(\ref{eq:I-vs-script-I-restricted})$ is a Quillen adjunction.

It only remains to establish left properness, for which we observe that any of the above generating cofibrations is a levelwise cofibration, and hence so is any cofibration of the strict level model structure. The claim therefore follows from Lemma~\ref{lemma:homotopy-pushout-M-SSet}.
\end{proof}

Proposition~\ref{prop:characterization-cofibration-general} specializes to:

\begin{cor}\label{cor:characterization-i-cof}
A map $f\colon X\to Y$ in $\cat{$\bm G$-$\bm I$-SSet}$ is a cofibration in the strict level model structure if and only if for each finite set $A$ the map
\begin{equation}\label{eq:latching-map}
X(A)\amalg_{L_A(X)}L_A(Y)\xrightarrow{(f(A),\ell_A)} Y(A)
\end{equation}
is a cofibration in the $\mathcal{G}_{\Sigma_A,G}$-model structure on $\cat{$\bm{(\Sigma_A\times G)}$-SSet}$.

In particular, $X$ is cofibrant in the strict level model structure if and only if for each finite set $A$ the latching map $\colim_{B\subsetneq A} X(B)\to X(A)$ is a cofibration in the above model category.\qedhere\qed
\end{cor}

\begin{ex}
If $G=1$, the $\mathcal G_{\Sigma_A,G}$-cofibrations are precisely the underlying cofibrations of simplicial sets. In this case, the strict level cofibrations on $\cat{$\bm I$-SSet}$ have been considered non-equivariantly under the name \emph{flat cofibrations} \cite[Definition~3.9]{sagave-schlichtkrull}.\index{flat cofibration!in I-SSet@in $\cat{$\bm I$-SSet}$}
\end{ex}

\subsubsection{$G$-global weak equivalences} Before we can introduce the $G$-global weak equivalences for the above models, we need some preparations.

\begin{constr}\label{constr:i-extension}\nomenclature[aI1z]{$\overline{I}$}{extension of $I$ to all sets}\nomenclature[aI2z]{$\overline{\mathcal I}$}{extension of $\mathcal I$ to all sets}
We write $\overline{I}$ for the category of all sets and injections, and we write $\overline{\mathcal I}$ for the simplicial category obtained by applying $E$ to each hom set. Then $I$ and $\mathcal I$ are full (simplicial) subcategories of $\overline{I}$ and $\overline{\mathcal I}$, respectively. We will now explain how to extend any $I$- or $\mathcal I$-simplicial set to $\overline{I}$ or $\overline{\mathcal I}$, respectively:

In the case of $X\colon I\to\cat{SSet}$ we define
\begin{equation*}
\overline{X}(A)=\colim\limits_{B\subset A\textup{ finite}} X(B)
\end{equation*}
for any set $A$. If $i\colon A\to A'$ is any injection, we define the structure map $\overline X(A)\to\overline X(A')$ as the map induced via the universal property of the above colimit by the family
\begin{equation*}
X(B)\xrightarrow{X(i|_B)} X(i(B))\to \colim\limits_{B'\subset A'\textup{ finite}} X(B')=\overline{X}(A')
\end{equation*}
for all finite $B\subset A$, where the unlabelled arrow is the structure map of the colimit for the term indexed by $i(B)\subset A'$. We omit the trivial verification that this functorial.

In the case of a simplicially enriched $X\colon\mathcal I\to\cat{SSet}$ we define the extension analogously on objects and morphisms. If now $(i_0,\dots,i_n)$ is a general $n$-cell of $\overline{\mathcal I}(A,A')$, then we define $\overline X(i_0,\dots,i_n)$ as the composition
\begin{equation*}
\Delta^n\times\colim\limits_{B\subset A\textup{ finite}} X(B)\cong
\colim\limits_{B\subset A\textup{ finite}}\Delta^n\times X(B)\to
\colim\limits_{B'\subset A'\textup{ finite}} X(B')
\end{equation*}
where the isomorphism is the canonical one and the unlabelled arrow is induced by $X(i_0|_B,\dots, i_n|_B)\colon\Delta^n\times X(B)\to X(i_0(B)\cup i_1(B)\cup\cdots\cup i_n(B))$ for all finite $B\subset A$. We omit the easy verification that this defines a simplicially enriched functor $\overline{\mathcal I}\to\cat{SSet}$.

One moreover easily checks that these become simplicially enriched functors by sending $F\colon\Delta^n\times X\to Y$ to the transformation given on a set $A$ by
\begin{equation*}
\Delta^n\times\colim\limits_{B\subset A\textup{ finite}} X(B)\cong
\colim\limits_{B\subset A\textup{ finite}}\Delta^n\times X(B)\xrightarrow{\colim F(B)}
\colim\limits_{B\subset A\textup{ finite}} Y(B),
\end{equation*}
and that with respect to this the structure maps
\begin{equation*}
X(A)\to \colim\limits_{B\subset A\textup{ finite}} X(B)=\overline X(A)
\end{equation*}
for finite $A$ define simplicially enriched natural isomorphisms; accordingly, we will from now on no longer distinguish notationally between the extension and the original object.
\end{constr}

\begin{rk}
It is not hard to check that the above is a model for the (simplicially enriched) left Kan extension; however, we will at several points make use of the above explicit description.
\end{rk}

\begin{rk}
Simply by functoriality, the above construction lifts to provide simplicially enriched extension functors
\begin{equation*}
\cat{$\bm G$-$\bm I$-SSet}\to\cat{$\bm G$-$\bm{\overline{I}}$-SSet}\qquad\text{and}\qquad
\cat{$\bm G$-$\bm{\mathcal I}$-SSet}\to\cat{$\bm G$-$\bm{\overline{\mathcal I}}$-SSet}
\end{equation*}
for any group $G$.
\end{rk}

For later use we record:

\begin{lemma}\label{lemma:evaluation-h-universe}
Let $\mathcal U$ be a complete $H$-set universe, let $A$ be any $H$-set, and let $i\colon\mathcal U\to A$ be an $H$-equivariant injection. Then $X(i)\colon X(\mathcal U)\to X(A)$ is an $(H\times G)$-weak equivalence for any $X\in\cat{$\bm G$-$\bm{\mathcal I}$-SSet}$.
\end{lemma}

The proof will rely on the following easy observation:

\begin{lemma}\label{lemma:extension-filtered-colimit}
If $X$ is a $G$-$I$- or $G$-$\mathcal I$-simplicial set, then its extension preserves filtered colimits, i.e.~if $J$ is a small filtered category and $A_\bullet\colon J\to \overline{I}$ any functor, then the canonical map
\begin{equation*}
\colim_j X(A_j)\to X\big(\colim_j A_j\big)
\end{equation*}
is an isomorphism. (Here it does not matter whether we form the colimit on the right hand side in $\cat{Set}$ or in $\overline{I}$).
\begin{proof}
Unravelling the definition, the left hand side is given by the double colimit $\colim_j\colim_{B\subset A_j\textup{ finite}} X(B)$. Write $\overline{A}\mathrel{:=}\colim_{j\in J}A_j$; we define a map in the other direction as follows: a finite subset $B\subset\overline{A}$ is contained in the image of some structure map $i\colon A_j\to\overline{A}$, and we send $X(B)$ via $X(i^{-1})$ to $X(i^{-1}(B))$ in the $(j, i^{-1}(B))$-term on the left hand side. We omit the easy verification that this is well-defined and inverse to the above map.
\end{proof}

\begin{proof}[Proof of Lemma~\ref{lemma:evaluation-h-universe}]
Let us first assume that $A$ is countable. As $\mathcal U$ is a complete $H$-set universe by assumption, we can therefore find an $H$-equivariant injection $j\colon A\to\mathcal U$. Then $X(j)$ is an $(H\times G)$-homotopy inverse to $X(i)$, as exhibited by the $(H\times G)$-equivariant homotopies $X(\id_A, ij)$ and $X(\id_{\mathcal U}, ji)$, finishing the proof of the special case.

In the general case we now observe that the map in question factors as
\begin{equation*}
X(\mathcal U)\cong\colim\limits_{i(\mathcal U)\subset B\subset A\textup{ countable $H$-set}} X(\mathcal U)\xrightarrow{\sim}
\colim\limits_{i(\mathcal U)\subset B\subset A\textup{ countable $H$-set}} X(B) \cong X(A),
\end{equation*}
where the left hand map is induced by the inclusion of the term indexed by $i(\mathcal U)$ (which is an isomorphism because $\{i(\mathcal U)\subset B\subset A\}$ is a filtered poset, so that it has connected nerve), the second map uses the above special case levelwise, and the final isomorphism comes from the previous lemma. The claim follows immediately.
\end{proof}
\end{lemma}

\begin{constr}\nomenclature[aevomega]{$\ev_\omega$}{functor $\cat{$\bm I$-SSet}\to\cat{$\bm{\mathcal M}$-SSet}$ or $\cat{$\bm{\mathcal I}$-SSet}\to\cat{$\bm{E\mathcal M}$-SSet}$ given by evaluating at $\omega$}
By applying Construction~\ref{constr:i-extension} and then restricting along the inclusion $B\mathcal M\hookrightarrow\overline{I}$ or $B(E\mathcal M)\to\overline{\mathcal I}$, respectively, sending the unique object to $\omega$, we get functors
\begin{equation*}
\ev_\omega\colon \cat{$\bm{G}$-$\bm{I}$-SSet}\to\cat{$\bm{\mathcal M}$-$\bm{G}$-SSet}
\qquad\text{and}\qquad
\ev_\omega\colon \cat{$\bm{G}$-$\bm{\mathcal I}$-SSet}\to\cat{$\bm{E\mathcal M}$-$\bm{G}$-SSet}.
\end{equation*}
\end{constr}

\begin{defi}\index{G-global weak equivalence@$G$-global weak equivalence!in G-II-SSet@in $\cat{$\bm G$-$\bm{\mathcal I}$-SSet}$|textbf}\index{G-global weak equivalence@$G$-global weak equivalence!in G-I-SSet@in $\cat{$\bm G$-$\bm{I}$-SSet}$|textbf}
A map $f\colon X\to Y$ in $\cat{$\bm G$-$\bm I$-SSet}$ or $\cat{$\bm G$-$\bm{\mathcal I}$-SSet}$ is called a \emph{$G$-global weak equivalence} if $f(\omega)=\ev_\omega f$ is a $G$-global weak equivalence in $\cat{$\bm{\mathcal M}$-$\bm G$-SSet}$ or $\cat{$\bm{E\mathcal M}$-$\bm G$-SSet}$, respectively.
\end{defi}

\begin{lemma}\label{lemma:strict-are-global}
Let $f\colon X\to Y$ be a strict level weak equivalence in $\cat{$\bm{G}$-$\bm{\mathcal I}$-SSet}$ or $\cat{$\bm{G}$-$\bm{I}$-SSet}$. Then $f$ is also a $G$-global weak equivalence.
\begin{proof}
It suffices to consider the second case. Let $H\subset\mathcal M$ be universal; we will show that $X(\omega)\to Y(\omega)$ is a $\mathcal G_{H,G}$-weak equivalence.

Pick a free $H$-orbit $F$ inside $\omega$ (which exists by universality). We then observe that we have a commutative diagram
\begin{equation*}
\begin{tikzcd}
\colim\limits_{F\subset A\subset\omega\textup{ finite $H$-set}} X(A)\arrow[d, "\colim f(A)"'] \arrow[r] & \colim\limits_{A\subset\omega\textup{ finite}} X(A)\arrow[d, "\colim f(A)"]\\
\colim\limits_{F\subset A\subset\omega\textup{ finite $H$-set}} Y(A)\arrow[r] & \colim\limits_{A\subset\omega\textup{ finite}} Y(A)
\end{tikzcd}
\end{equation*}
where the horizontal maps are induced from the inclusion of filtered posets $\{F\subset A\subset\omega\textup{ finite $H$-set}\}\hookrightarrow\{A\subset\omega\textup{ finite}\}$. This is is cofinal: any finite subset $A\subset\omega$ is contained in the finite $H$-set $F\cup HA$ which is an element of the left hand side.

Thus, the horizontal maps are isomorphisms, and it therefore suffices that the left hand vertical map is a $\mathcal G_{H,G}$-weak equivalence. But on this side, $H$ simply acts on each term of the colimit by functoriality, and each $f(A)$ is a $\mathcal G_{H,G}$-weak equivalence with respect to this action as $A$ is in particular faithful. The claim follows as the $\mathcal G_{H,G}$-weak equivalences are closed under filtered colimits.
\end{proof}
\end{lemma}

In order to later characterize the fibrant objects in the $G$-global model structure, we introduce:

\begin{defi}
A $G$-$I$-simplicial set (or $G$-$\mathcal I$-simplicial set) $X$ is called \emph{static},\index{static|textbf} if for all finite faithful $H$-sets $A$ and each $H$-equivariant injection $i\colon A\to B$ into another finite $H$-set, the induced map $X(i)\colon X(A)\to X(B)$ is a $\mathcal G_{H,G}$-weak equivalence.
\end{defi}

\begin{lemma}
Let $f\colon X\to Y$ be a map of static $G$-$I$-simplicial sets or $G$-$\mathcal I$-simplicial sets. Then $f$ is a $G$-global equivalence if and only if $f$ is a strict level weak equivalence.
\begin{proof}
Again it suffices to consider the first case.

The implication `$\Leftarrow$' holds without any assumptions by Lemma~\ref{lemma:strict-are-global}. For the remaining implication we first observe:

\begin{claim*}
Both $X(\omega)$ and $Y(\omega)$ are $G$-semistable.
\begin{proof}
It suffices to prove the first statement. Let $H\subset\mathcal M$ be universal, let $u\in\mathcal M$ centralize $H$, and pick a finite faithful $H$-subset $A\subset\omega$. We now consider the commutative diagram
\begin{equation*}
\begin{tikzcd}
X(A)\arrow[d, "X(u|_A)"'] \arrow[r] & \colim\limits_{A\subset B\subset\omega\textup{ finite $H$-set}} X(B) \arrow[r, "\cong"] & X(\omega)\arrow[d, "u.\blank"]\\
X(u(A)) \arrow[r] & \colim\limits_{u(A)\subset B\subset\omega\textup{ finite $H$-set}} X(B) \arrow[r, "\cong"'] & X(\omega)
\end{tikzcd}
\end{equation*}
where the isomorphisms on the right come from cofinality again, and the left hand horizontal maps are structure maps of the respective colimits, hence $\mathcal G_{H,G}$-weak equivalences as all transition maps are. As the left hand vertical map is an isomorphism for trivial reasons, we conclude that also the right hand vertical map is a $\mathcal G_{H,G}$-weak equivalence, i.e.~$X(\omega)$ is semistable as desired.
\end{proof}
\end{claim*}

Let $H$ be a finite group and $A$ a finite faithful $H$-set; we have to show that $f(A)$ is a $\mathcal G_{H,G}$-weak equivalence, for which we may assume without loss of generality that $H$ is a universal subgroup of $\mathcal M$ and $A$ an $H$-subset of $\omega$. But the same argument as above then shows that $f(A)$ agrees up to conjugation by $\mathcal G_{H,G}$-weak equivalences with $f(\omega)$. The latter is a $G$-global weak equivalence between $G$-semistable $\mathcal M$-$G$-simplicial sets by the above claim, hence a $G$-universal weak equivalence. Thus, also $f(A)$ is a $\mathcal G_{H,G}$-weak equivalence by $2$-out-of-$3$ as desired.
\end{proof}
\end{lemma}

\subsubsection{Connection to tame actions}
\index{tame!M-action@$\mathcal M$-action!tame M-simplicial sets vs flat I-simplicial sets@tame $\mathcal M$-simplicial sets vs.~flat $I$-simplicial sets|(}
In order to construct the $G$-global model structures on $\cat{$\bm G$-$\bm I$-SSet}$ and $\cat{$\bm G$-$\bm{\mathcal I}$-SSet}$ and to compare them to our previous models, we will exploit a close connection between them on the point-set level. Namely, Sagave and Schwede showed that the $1$-category of tame $\mathcal M$-simplicial sets is equivalent to the full subcategory of $\cat{$\bm I$-SSet}$ spanned by the flat\index{flat!I-simplicial set@$I$-simplicial set} (i.e.~globally cofibrant) objects. In order to state their precise comparison, we need:

\begin{constr}\nomenclature[zbullet]{$(\blank)_\bullet$}{right adjoint to $\ev_\omega$, \textit{see also} $(\blank)_{[A]}$}
Let $X$ be any $\mathcal M$-simplicial set. We write $X_\bullet$ for the $I$-simplicial set with $(X_\bullet)(A)=X_{[A]}$ for each finite $A\subset\omega$; an injection $j\colon A\to B$ acts by extending it to an injection $\bar{\jmath}\in\mathcal M$ and then using the $\mathcal M$-action (this is well-defined by Lemma~\ref{lemma:support-vs-action-M} together with Lemma~\ref{lemma:support-agree-M}). As any finite set is isomorphic to a subset of $\omega$, there is an essentially unique way to extend this to a functor $I\to\cat{SSet}$, and we fix any such extension.

This becomes a simplicial functor by the enriched functoriality of the individual $X_{[A]}$; in particular, we get an induced functor $(\blank)_\bullet\colon\cat{$\bm{\mathcal M}$-$\bm G$-SSet}\to\cat{$\bm G$-$\bm I$-SSet}$ for any group $G$.

The inclusions $X_{[A]}\hookrightarrow X$ assemble into an enriched natural map $\epsilon\colon X_\bullet(\omega)\to X$ for any $\mathcal M$-$G$-simplicial set $X$. Moreover, if $Y$ is a $G$-$I$-simplicial set, then the structure map $Y(A)\to Y(\omega)$ for any finite $A\subset\omega$ factors through $Y(\omega)_{[A]}$, and for varying $A$ these assemble into an enriched natural map $Y\to Y(\omega)_\bullet$.
\end{constr}

The above construction is a `coordinate free' version of \cite[Construction~5.5]{I-vs-M-1-cat} applied in each simplicial degree (and with $G$-actions pulled through), also cf.~\cite[discussion before Corollary~5.7]{I-vs-M-1-cat}. In particular,~\cite[Proposition~5.6]{I-vs-M-1-cat} implies, also cf.~\cite[Corollary~5.7]{I-vs-M-1-cat}:

\begin{lemma}\label{lemma:M-SSet-vs-I-SSet-sagave-schwede}
The above defines a simplicial adjunction
\begin{equation}\label{eq:evaluation-support-adjunction}
\ev_\omega\colon\cat{$\bm G$-$\bm I$-SSet}\rightleftarrows\cat{$\bm{\mathcal M}$-$\bm G$-SSet} :\!(\blank)_\bullet
\end{equation}
where the left adjoint has image in $\cat{$\bm{\mathcal M}$-$\bm G$-SSet}^\tau$, and for any $\mathcal M$-$G$-simplicial set $X$ the counit $X_\bullet(\omega)\to X$ factors through an isomorphism onto $X^\tau$. In particular, $(\ref{eq:evaluation-support-adjunction})$ restricts to a Bousfield localization $\cat{$\bm G$-$\bm I$-SSet}\rightleftarrows\cat{$\bm{\mathcal M}$-$\bm G$-SSet}^\tau$.

Finally, the right adjoint has essential image the flat $G$-$I$-simplicial sets, i.e.~the unit $\eta\colon Y\to Y(\omega)_\bullet$ is an isomorphism if and only if $Y$ is cofibrant in $\cat{$\bm I$-SSet}$.\index{tame!M-action@$\mathcal M$-action!tame M-simplicial sets vs flat I-simplicial sets@tame $\mathcal M$-simplicial sets vs.~flat $I$-simplicial sets|)}\qed
\end{lemma}

Using Theorem~\ref{thm:support-EM-vs-M}, we will now give an analogous comparison between $\cat{$\bm G$-$\bm{\mathcal I}$-SSet}$ and $\cat{$\bm{E\mathcal M}$-$\bm G$-SSet}$.

\index{tame!EM-action@$E\mathcal M$-action!tame EM-simplicial sets vs flat II-simplicial sets@tame $E\mathcal M$-simplicial sets vs.~flat $\mathcal I$-simplicial sets|(}
\begin{constr}
Let $X$ be any $E\mathcal M$-simplicial set. We write $X_\bullet$ for the $\mathcal I$-simplicial set with $(X_\bullet)(A)=X_{[A]}$ for every finite $A\subset\omega$; an $(n+1)$-tuple of injections $j_0,\dots,j_n\colon A\to B$ acts by extending each of them to an injection $\bar{\jmath}_k\colon\omega\to\omega$ and then using the $E\mathcal M$-action (which is well-defined by Corollary~\ref{cor:support-agree-EM} together with Lemma~\ref{lemma:support-vs-action-EM}). Again we fix an extension to all of $\mathcal I$.

This becomes a simplicial functor by enriched functoriality of the individual $X_{[A]}$; in particular, if $G$ is any group we then again get an induced simplicial functor $(\blank)_\bullet\colon\cat{$\bm{E\mathcal M}$-$\bm G$-SSet}\to\cat{$\bm G$-$\bm{\mathcal I}$-SSet}$.

We define $\epsilon\colon X_\bullet(\omega)\to X$ as the map induced by the inclusions $X_{[A]}\hookrightarrow X$. Moreover, if $Y\in\cat{$\bm G$-$\bm{\mathcal I}$-SSet}$, then $Y(A)\to Y(\omega)$ factors through $Y(\omega)_{[A]}$ by definition of the action, and one easily checks that these assemble into $\eta\colon Y\to Y(\omega)_\bullet$.
\end{constr}

\begin{rk}
By Theorem~\ref{thm:support-EM-vs-M}, the diagram
\begin{equation}\label{eq:forget-vs-blank-bullet}
\begin{tikzcd}
\cat{$\bm{E\mathcal M}$-$\bm G$-SSet}\arrow[r, "(\blank)_\bullet"]\arrow[d, "\forget"'] & \cat{$\bm G$-$\bm{\mathcal I}$-SSet}\arrow[d, "\forget"]\\
\cat{$\bm{\mathcal M}$-$\bm G$-SSet}\arrow[r, "(\blank)_\bullet"'] & \cat{$\bm G$-$\bm{I}$-SSet}
\end{tikzcd}
\end{equation}
commutes strictly, and the same can be arranged for $\ev_\omega$ instead of $(\blank)_\bullet$ by construction. Under these identifications, also the unit and counit are preserved, i.e.~$\forget(\eta_Y)=\eta_{\forget Y}$ and $\forget(\epsilon_X)=\epsilon_{\forget X}$.
\end{rk}

As the notions of tameness in $\cat{$\bm{E\mathcal M}$-$\bm G$-SSet}$ and $\cat{$\bm{\mathcal M}$-$\bm G$-SSet}$ agree by another application of Theorem~\ref{thm:support-EM-vs-M}, we conclude from Lemma~\ref{lemma:M-SSet-vs-I-SSet-sagave-schwede}:

\begin{lemma}\label{lemma:evaluation-support-adjunction-E}
The above yields a simplicially enriched adjunction
\begin{equation}\label{eq:evaluation-support-adjunction-E}
\ev_\omega\colon\cat{$\bm G$-$\bm{\mathcal I}$-SSet}\rightleftarrows\cat{$\bm{E\mathcal M}$-$\bm G$-SSet} :\!(\blank)_\bullet
\end{equation}
where the left adjoint has image in $\cat{$\bm{E\mathcal M}$-$\bm G$-SSet}^\tau$. Moreover,
for any $E\mathcal M$-$G$-simplicial set $X$ the counit $(X_\bullet)(\omega)\to X$ factors through an isomorphism onto $X^\tau$, so that $(\ref{eq:evaluation-support-adjunction-E})$ restricts to a Bousfield localization $\cat{$\bm G$-$\bm{\mathcal I}$-SSet}\rightleftarrows\cat{$\bm{E\mathcal M}$-$\bm G$-SSet}^\tau$.

Finally, the essential image of $(\blank)_\bullet$ consists precisely of those $G$-$\mathcal I$-simplicial sets that are flat,\index{flat!II-simplicial set@$\mathcal I$-simplicial set} i.e.~whose underlying $I$-simplicial sets are globally cofibrant.\index{tame!EM-action@$E\mathcal M$-action!tame EM-simplicial sets vs flat II-simplicial sets@tame $E\mathcal M$-simplicial sets vs.~flat $\mathcal I$-simplicial sets|)}\qed
\end{lemma}

As an upshot of this, we can now very easily prove the following alternative description of the $G$-global weak equivalences of $G$-$I$-simplicial sets, that will become useful at several points later:

\begin{thm}\label{thm:G-global-we-I-characterization}
\index{G-global weak equivalence@$G$-global weak equivalence!in G-I-SSet@in $\cat{$\bm G$-$\bm{I}$-SSet}$}
The following are equivalent for a map $f$ in $\cat{$\bm G$-$\bm I$-SSet}$:
\begin{enumerate}
\item $f$ is a $G$-global weak equivalence in $\cat{$\bm G$-$\bm{I}$-SSet}$.
\item $f(\omega)$ is a $G$-global weak equivalence in $\cat{$\bm{\mathcal M}$-$\bm G$-SSet}$.
\item $\mathcal I\times_If$ is a $G$-global weak equivalence in $\cat{$\bm G$-$\bm{\mathcal I}$-SSet}$.\label{item:I-underived}
\item $E\mathcal M\times_{\mathcal M} f(\omega)$ is a $G$-global weak equivalence in $\cat{$\bm{E\mathcal M}$-$\bm G$-SSet}$.\label{item:EM-underived}
\end{enumerate}
\end{thm}

We emphasize that the above functors are \emph{not} derived in any way.

\begin{proof}
The equivalence $(1)\Leftrightarrow(2)$ holds by definition, and $(2)\Leftrightarrow(4)$ is an instance of Theorem~\ref{thm:tame-M-sset-vs-EM-sset}. It therefore only remains to show that $(3)\Leftrightarrow(4)$, for which we observe that the total mate of $(\ref{eq:forget-vs-blank-bullet})$ provides a natural isomorphism filling
\begin{equation*}
\begin{tikzcd}
\cat{$\bm{E\mathcal M}$-$\bm G$-SSet}^\tau &\arrow[l, "\ev_\omega"'] \cat{$\bm G$-$\bm{\mathcal I}$-SSet}\\
\cat{$\bm{\mathcal M}$-$\bm G$-SSet}^\tau\arrow[u, "E\mathcal M\times_{\mathcal M}\blank"] &\arrow[l, "\ev_\omega"] \cat{$\bm G$-$\bm{I}$-SSet}.\arrow[u, "\mathcal I\times_I\blank"']
\end{tikzcd}
\end{equation*}
The claim then follows immediately from the definitions.
\end{proof}

\subsubsection{$G$-global model structures} Using the above as well as our knowledge about tame $\mathcal M$- and $E\mathcal M$-actions we can now prove:

\begin{thm}\label{thm:script-I-global-model-structure}\index{G-global model structure@$G$-global model structure!on G-II-SSet@on $\cat{$\bm G$-$\bm{\mathcal I}$-SSet}$|textbf}
There is a unique model structure on $\cat{$\bm G$-$\bm{\mathcal I}$-SSet}$ whose weak equivalences are the $G$-global weak equivalences and with cofibrations those of the strict level model structure. We call this model structure the \emph{$G$-global model structure}. It is proper, combinatorial, simplicial, and filtered colimits in it are homotopical. Moreover, the fibrant objects of this model structure are precisely the strictly level fibrant static $G$-$\mathcal I$-simplicial sets.

Finally, the simplicial adjunction $(\ref{eq:evaluation-support-adjunction-E})$ is a Quillen equivalence with respect to the $G$-global injective model structure on $\cat{$\bm{E\mathcal M}$-$\bm G$-SSet}$.
\begin{proof}
By Lemma~\ref{lemma:strict-are-global}, $\ev_\omega$ sends strict level weak equivalences to $G$-global weak equivalences, and it clearly sends generating cofibrations to injective cofibrations. In particular, the simplicial adjunction
\begin{equation*}
\ev_\omega\colon\cat{$\bm G$-$\bm{\mathcal I}$-SSet}_{\textup{strict level}}\rightleftarrows\cat{$\bm{E\mathcal M}$-$\bm G$-SSet}_{\textup{injective $G$-global}} :\!(\blank)_\bullet
\end{equation*}
is a Quillen adjunction. We now want to apply Lurie's localization criterion (Theorem~\ref{thm:lurie-localization-criterion}) to this, for which we have to show that $\cat{R}(\blank)_\bullet$ is fully faithful with essential image the static $G$-$\mathcal I$-simplicial sets.

Indeed, Corollary~\ref{cor:injective-almost-tame} shows that $\cat{R}(\blank)_\bullet$ restricts accordingly, and that the counit $X_\bullet(\omega)\to X$ is a $G$-global weak equivalence for any injectively fibrant $X$. On the other hand, for any $G$-$\mathcal I$-simplicial set $Y$, the map $\eta(\omega)\colon Y(\omega)\to Y(\omega)_\bullet(\omega)$ is a one-sided inverse of $\epsilon_{Y(\omega)}$, hence an isomorphism by Lemma~\ref{lemma:evaluation-support-adjunction-E} as $Y(\omega)$ is tame. Since $(\blank)_\bullet$ clearly preserves $G$-global weak equivalences in $\cat{$\bm{E\mathcal M}$-$\bm G$-SSet}^{w\tau}$ (i.e.~the subcategory of those $E\mathcal M$-$G$-simplicial sets $X$ for which $X^\tau\hookrightarrow X$ is a $G$-global weak equivalence), and as this contains both the injectively fibrant objects by Corollary~\ref{cor:injective-almost-tame} as well as the tame $E\mathcal M$-$G$-simplicial set $Y(\omega)$ for trivial reasons, we conclude that the derived unit $\eta_Y$, represented by the composition $Y\to Y(\omega)_\bullet\to Z_\bullet$ for some injectively fibrant replacement $Y(\omega)\to Z$, is a $G$-global weak equivalence. Thus, if $Y$ is static, then the derived unit is a $G$-global weak equivalence between static $G$-$\mathcal I$-simplicial sets, hence a strict level weak equivalence by Lemma~\ref{lemma:strict-are-global}. This completes the proof of the claim.

Lurie's criterion then shows that the desired model structure exists, and that it is left proper, combinatorial, simplicial, has the fibrant objects described above, and that $(\ref{eq:evaluation-support-adjunction-E})$ becomes a Quillen equivalence for this model structure. Moreover, Lemma~\ref{lemma:filtered-still-homotopical} shows that filtered colimits in this model structure are homotopical.

Finally, we consider a pullback square
\begin{equation*}
\begin{tikzcd}
P\arrow[r, "f"]\arrow[d, "p"'] & X\arrow[d, "q"]\\
Y\arrow[r, "g"'] & Z
\end{tikzcd}
\end{equation*}
in $\cat{$\bm{G}$-$\bm{\mathcal I}$-SSet}$ such that $q$ is a strict level fibration and $g$ is a $G$-global weak equivalence. We will show that also $f$ is a $G$-global weak equivalence, which will in particular imply right properness of the $G$-global model structure.

As finite limits in $\cat{SSet}$ commute with filtered colimits, and as limits commute with each other, we get for any universal subgroup $H\subset\mathcal M$ and any group homomorphism $\phi\colon H\to G$ a pullback square
\begin{equation*}
\begin{tikzcd}[column sep=large]
P(\omega)^\phi\arrow[r, "f(\omega)^\phi"]\arrow[d, "p(\omega)^\phi"'] & X(\omega)^\phi\arrow[d, "q(\omega)^\phi"]\\
Y(\omega)^\phi\arrow[r, "g(\omega)^\phi"'] & Z(\omega)^\phi
\end{tikzcd}
\end{equation*}
in $\cat{SSet}$. The map $g(\omega)^\phi$ is a weak equivalence by definition, and we have to show that also $f(\omega)^\phi$ is. For this it suffices by right properness of $\cat{SSet}$ that $q(\omega)^\phi$ is a Kan fibration. But as before, after picking a free $H$-orbit $F\subset\omega$, it can be identified with the filtered colimit
\begin{equation*}
\colim\limits_{F\subset A\subset\omega\textup{ finite $H$-set}} q(A)^\phi
\end{equation*}
of Kan fibrations, and hence is itself a Kan fibration as desired.
\end{proof}
\end{thm}

\begin{thm}\label{thm:strict-global-I-model-structure}
There is a unique model structure on $\cat{$\bm G$-$\bm I$-SSet}$ whose cofibrations are the strict level cofibrations and whose weak equivalences are the $G$-global weak equivalences. We call this the \emph{$G$-global model structure}.\index{G-global model structure@$G$-global model structure!on G-I-SSet@on $\cat{$\bm G$-$\bm I$-SSet}$|textbf} It is left proper, combinatorial, simplicial, and filtered colimits in it are homotopical. Moreover, its fibrant objects are precisely the static strictly level fibrant ones.

Finally, the simplicial adjunctions
\begin{align}
\ev_\omega\colon\cat{$\bm G$-$\bm I$-SSet}_{\textup{$G$-global}}&\rightleftarrows\cat{$\bm{\mathcal M}$-$\bm G$-SSet}_{\textup{injective $G$-global}} :\!(\blank)_\bullet\nonumber\\
\label{eq:I-vs-script-I}
\mathcal I\times_I\blank\colon\cat{$\bm G$-$\bm I$-SSet}_{\textup{$G$-global}}&\rightleftarrows\cat{$\bm G$-$\bm{\mathcal I}$-SSet}_{\textup{$G$-global}}:\!\forget
\end{align}
are Quillen equivalences, and both functors in $(\ref{eq:I-vs-script-I})$ are fully homotopical.
\begin{proof}
All statements except for those about the adjunction $(\ref{eq:I-vs-script-I})$ are proven just as in the previous theorem.

For the remaining statements, we observe that $\mathcal I\times_I\blank$ preserves cofibrations by Proposition~\ref{prop:strict-model-structure} while it is homotopical by Theorem~\ref{thm:G-global-we-I-characterization}; in particular it is left Quillen. On the other hand, the forgetful functor is clearly homotopical. We then consider the diagram
\begin{equation*}
\begin{tikzcd}
\cat{$\bm G$-$\bm{\mathcal I}$-SSet}\arrow[d, "\ev_\omega"']\arrow[r, "\forget"] &[1em] \cat{$\bm G$-$\bm{I}$-SSet}\arrow[d, "\ev_\omega"]\\
\cat{$\bm{E\mathcal M}$-$\bm G$-SSet} \arrow[r, "\forget"'] & \cat{$\bm{\mathcal M}$-$\bm G$-SSet}
\end{tikzcd}
\end{equation*}
of homotopical functors (with respect to the $G$-global weak equivalences everywhere), which commutes up to canonical isomorphism. By the above, the vertical maps induce equivalences on associated quasi-categories, and so does the lower horizontal map by Corollary~\ref{cor:em-vs-m-equiv-model-cat}. The claim follows by $2$-out-of-$3$.
\end{proof}
\end{thm}

\subsubsection{Further model structures}
For later use we record the existence of \emph{positive $G$-global model structures}, which can be constructed in precisely the same way as above; we leave the details to the reader.

\begin{defi}
A $G$-$I$-simplicial set $X$ is called \emph{positively static}\index{static!positive|seeonly{positively static}}\index{positively static|textbf} if the map $X(i)\colon X(A)\to X(B)$ is a $\mathcal G_{H,G}$-weak equivalence for every finite group $H$ and every injection $i\colon A\to B$ of finite \emph{non-empty} faithful $H$-sets. It is called \emph{positively level fibrant} if $X(A)$ is $\mathcal G_{\Sigma_A,G}$-equivariantly fibrant for every \emph{non-empty} finite set $A$.

A $G$-$\mathcal I$-simplicial set is called positively static or positively level fibrant if it is so as a $G$-$I$-simplicial set.
\end{defi}

\begin{thm}\label{thm:positive-G-global-script-I}\index{positive G-global model structure@positive $G$-global model structure!on G-II-SSet@on $\cat{$\bm G$-$\bm{\mathcal I}$-SSet}$|textbf}
\index{G-global model structure@$G$-global model structure!positive|seeonly{positive $G$-global model structure}}
There is a unique cofibrantly generated model structure on $\cat{$\bm G$-$\bm{\mathcal I}$-SSet}$ with weak equivalences the $G$-global weak equivalences and generating cofibrations the maps
\begin{equation*}
\mathcal I(A,\blank)\times_\phi G\times\del\Delta^n\hookrightarrow\mathcal I(A,\blank)\times_\phi G\times\Delta^n
\end{equation*}
for $n\ge 0$, finite groups $H$, homomorphisms $\phi\colon H\to G$, and \emph{non-empty} finite faithful $H$-sets $A$. We call this the \emph{positive $G$-global model structure}. It is combinatorial, simplicial, proper, and filtered colimits in it are homotopical. Moreover, its fibrant objects are precisely the positively static positively level fibrant ones.

Finally, the identity adjunction $\cat{$\bm G$-$\bm{\mathcal I}$-SSet}_{\textup{positive $G$-global}}\rightleftarrows\cat{$\bm G$-$\bm{\mathcal I}$-SSet}_{\textup{$G$-global}}$ is a Quillen equivalence.\qed
\end{thm}

\begin{thm}\label{thm:positive-global-I-model-structure}
There is a unique cofibrantly generated model structure on $\cat{$\bm G$-$\bm I$-SSet}$ with weak equivalences the $G$-global weak equivalences and generating cofibrations the maps
\begin{equation*}
I(A,\blank)\times_\phi G\times\del\Delta^n\hookrightarrow I(A,\blank)\times_\phi G\times\Delta^n
\end{equation*}
for $n\ge 0$, finite groups $H$, homomorphisms $\phi\colon H\to G$, and \emph{non-empty} finite faithful $H$-sets $A$. We call this the \emph{positive $G$-global model structure}.\index{positive G-global model structure@positive $G$-global model structure!on G-I-SSet@on $\cat{$\bm G$-$\bm I$-SSet}$|textbf} It is left proper, combinatorial, simplicial, and filtered colimits in it are homotopical. Moreover, its fibrant objects are precisely the positively static positively level fibrant ones. Finally, the simplicial adjunctions
\begin{align*}
\mathcal I\times_I\blank\colon\cat{$\bm G$-$\bm{I}$-SSet}_{\textup{positive $G$-global}}&\rightleftarrows\cat{$\bm G$-$\bm{\mathcal I}$-SSet}_{\textup{positive $G$-global}} :\!\forget\\
\id\colon\cat{$\bm G$-$\bm{I}$-SSet}_{\textup{positive $G$-global}}&\rightleftarrows\cat{$\bm G$-$\bm{I}$-SSet}_{\textup{$G$-global}} :\!\id
\end{align*}
are Quillen equivalences.\qed
\end{thm}

\begin{rk}
Again, there are suitable \emph{positive level model structures}\index{positive level model structure} in the background of which the above are Bousfield localizations.
\end{rk}

\begin{lemma}\label{lemma:pos-cof-script-I}
Let $f\colon X\to Y$ be a cofibration in either of the $G$-global positive model structures. Then $f(\varnothing)\colon X(\varnothing)\to Y(\varnothing)$ is an isomorphism.
\begin{proof}
The class of such maps is obviously closed under retracts, pushouts, and transfinite compositions. Thus, it suffices to verify the claim for each generating cofibration $i$. But in this case both source and target are obviously empty in degree $\varnothing$, in particular $i(\varnothing)$ is an isomorphism.
\end{proof}
\end{lemma}

Finally, we come to injective model structures:

\begin{thm}\label{thm:script-I-vs-EM-injective}
\index{injective G-global model structure@injective $G$-global model structure!on G-II-SSet@on $\cat{$\bm G$-$\bm{\mathcal I}$-SSet}$|textbf}
There exists a unique model structure on $\cat{$\bm G$-$\bm{\mathcal I}$-SSet}$ whose cofibrations are the injective cofibrations and whose weak equivalences are the $G$-global weak equivalences. We call this the \emph{injective $G$-global model structure}. It is combinatorial, proper, simplicial, and filtered colimits in it are homotopical.
\end{thm}

For the proof we will need:

\begin{lemma}\label{lemma:G-global-I-injective-pushout}
The $G$-global weak equivalences in $\cat{$\bm G$-$\bm{\mathcal I}$-SSet}$ and $\cat{$\bm G$-$\bm{I}$-SSet}$ are stable under pushout along injective cofibrations, and a commutative square in either of these is a homotopy pushout if and only if its image under $\ev_\omega$ is.
\begin{proof}
As the left adjoint functor $\ev_\omega$ preserves injective cofibrations and creates $G$-global weak equivalences, this is simply an instance of Lemma~\ref{lemma:U-pushout-preserve-reflect} together with the existence of the injective $G$-global model structures on $\cat{$\bm{\mathcal M}$-$\bm G$-SSet}$ and $\cat{$\bm{E\mathcal M}$-$\bm G$-SSet}$ (Theorem~\ref{thm:G-M-injective-semistable-model-structure} and Proposition~\ref{prop:equivariant-injective-model-structure}, respectively).
\end{proof}
\end{lemma}

\begin{proof}[Proof of Theorem~\ref{thm:script-I-vs-EM-injective}]
As observed in the proof of Proposition~\ref{prop:strict-model-structure}, the cofibrations of the $G$-global model structure are in particular injective cofibrations. On the other hand, pushouts along injective cofibrations preserve $G$-global weak equivalences by Lemma~\ref{lemma:G-global-I-injective-pushout}. Corollary~\ref{cor:mix-model-structures} therefore shows that the model structure exists, that it is combinatorial and proper, and that filtered colimits in it are homotopical.

It only remains to verify the Pushout Product Axiom for the simplicial tensoring, which in turn follows immediately for the Pushout Product Axiom for $\cat{SSet}$ and for the injective model structure on $\cat{$\bm{E\mathcal M}$-$\bm G$-SSet}$.
\end{proof}

\begin{thm}\label{thm:I-vs-M-injective}
\index{injective G-global model structure@injective $G$-global model structure!on G-I-SSet@on $\cat{$\bm G$-$\bm{I}$-SSet}$|textbf}
There exists a unique model structure on $\cat{$\bm G$-$\bm{I}$-SSet}$ whose cofibrations are the injective cofibrations and whose weak equivalences are the $G$-global weak equivalences. We call this the \emph{injective $G$-global model structure}. It is combinatorial, left proper, simplicial, and filtered colimits in it are homotopical. Moreover, the forgetful functor is part of a simplicial Quillen equivalence
\begin{equation}\label{eq:forget-lQ-injective-I-SSet}
\forget\colon\cat{$\bm G$-$\bm{\mathcal I}$-SSet}_{\textup{inj.~$G$-global}}\rightleftarrows\cat{$\bm G$-$\bm I$-SSet}_{\textup{inj.~$G$-global}} :\Maps_I(\mathcal I,\blank).
\end{equation}
\begin{proof}
The first part is analogous to the proof of the previous theorem.

For the final statement, we observe that the forgetful functor admits a simplicial right adjoint given by simplicially enriched right Kan extension along $I\to\mathcal I$, which we denote by $\Maps_I(\mathcal I,\blank)$.\nomenclature[amapsII]{$\Maps_I(\mathcal I,\blank)$}{right adjoint to forgetful functor $\cat{$\bm{\mathcal I}$-SSet}\to\cat{$\bm{I}$-SSet}$}
As the forgetful functor preserves weak equivalences and injective cofibrations for trivial reasons, it is left Quillen. Thus, $(\ref{eq:forget-lQ-injective-I-SSet})$ is a Quillen equivalence by Theorem~\ref{thm:strict-global-I-model-structure}.
\end{proof}
\end{thm}

\subsection{Functoriality}\index{functoriality in homomorphisms!for G-II-SSet@for $\cat{$\bm G$-$\bm{\mathcal I}$-SSet}$|(}
We will now discuss change of group functors for the above models.

\begin{lemma}\label{lemma:alpha-shriek-projective-script-I}
Let $\alpha\colon H\to G$ be any group homomorphism. Then
\begin{equation}\label{eq:alpha-shriek-projective-script-I}
\alpha_!\colon\cat{$\bm H$-$\bm{\mathcal I}$-SSet}_{\textup{$H$-global}}\rightleftarrows\cat{$\bm G$-$\bm{\mathcal I}$-SSet}_{\textup{$G$-global}} :\!\alpha^*
\end{equation}
is a simplicial Quillen adjunction with fully homotopical right adjoint, and likewise for the positive model structures on either side.
\begin{proof}
For the first statement, one immediately checks that $\alpha^*$ is right Quillen with respect to the strict level model structures. On the other hand, it obviously sends static $G$-$\mathcal I$-simplicial sets to static $H$-$\mathcal I$-simplicial sets, so that also the simplicial adjunction $(\ref{eq:alpha-shriek-projective-script-I})$ is a Quillen adjunction by Proposition~\ref{prop:cofibrations-fibrant-qa}.

To see that $\alpha^*$ is homotopical, we observe that in the commutative diagram
\begin{equation*}
\begin{tikzcd}
\cat{$\bm G$-$\bm{\mathcal I}$-SSet}_\textup{$G$-global}\arrow[r, "\ev_\omega"]\arrow[d, "\alpha^*"'] & \cat{$\bm{E\mathcal M}$-$\bm G$-SSet}_\textup{$G$-global}\arrow[d, "\alpha^*"]\\
\cat{$\bm H$-$\bm{\mathcal I}$-SSet}_\textup{$H$-global}\arrow[r, "\ev_\omega"'] & \cat{$\bm{E\mathcal M}$-$\bm H$-SSet}_\textup{$H$-global}
\end{tikzcd}
\end{equation*}
the horizontal arrows create weak equivalences by definition while the right hand vertical arrow obviously preserves weak equivalences; the claim follows immediately.

The proof for the positive model structures is analogous.
\end{proof}
\end{lemma}

\begin{cor}
Let $\alpha\colon H\to G$ be any group homomorphism. Then
\begin{equation*}
\alpha^*\colon\cat{$\bm G$-$\bm{\mathcal I}$-SSet}_{\textup{$G$-global injective}}\rightleftarrows\cat{$\bm H$-$\bm{\mathcal I}$-SSet}_{\textup{$H$-global injective}} :\!\alpha_*
\end{equation*}
is a simplicial Quillen adjunction.
\begin{proof}
It is obvious that $\alpha^*$ preserves injective cofibrations, and it is moreover homotopical by the previous lemma, hence left Quillen.
\end{proof}
\end{cor}

\begin{lemma}\label{lemma:alpha-shriek-injective-script-I}
Let $\alpha\colon H\to G$ be an \emph{injective} group homomorphism. Then
\begin{equation*}
\alpha_!\colon\cat{$\bm H$-$\bm{\mathcal I}$-SSet}_{\textup{$H$-global injective}}\rightleftarrows\cat{$\bm G$-$\bm{\mathcal I}$-SSet}_{\textup{$G$-global injective}} :\!\alpha^*
\end{equation*}
is a simplicial Quillen adjunction. In particular, $\alpha_!$ is homotopical.
\begin{proof}
We may assume without loss of generality that $H$ is a subgroup of $G$ and that $\alpha$ is its inclusion, in which case $\alpha_!$ can be modelled by applying $G\times_H\blank$ levelwise. We immediately see that $\alpha_!$ preserves injective cofibrations. To finish the proof it suffices now to show that it is also homotopical, for which we consider
\begin{equation*}
\begin{tikzcd}
\cat{$\bm H$-$\bm{\mathcal I}$-SSet}\arrow[d, "\alpha_!"']\arrow[r, "\ev_\omega"] & \cat{$\bm{E\mathcal M}$-$\bm H$-SSet}\arrow[d, "\alpha_!"]\\
\cat{$\bm G$-$\bm{\mathcal I}$-SSet}\arrow[r, "\ev_\omega"'] & \cat{$\bm{E\mathcal M}$-$\bm G$-SSet}.
\end{tikzcd}
\end{equation*}
We claim that this commutes up to isomorphism. Indeed, as $\alpha_!$ is cocontinuous, there is a canonical $G$-equivariant isomorphism filling this, and one easily checks that this isomorphism is also $E\mathcal M$-equivariant.

But the horizontal arrows in the above diagram preserve and reflect weak equivalences by definition and the right hand arrow is homotopical by Corollary~\ref{cor:alpha-shriek-injective-EM}, so the claim follows immediately.
\end{proof}
\end{lemma}

\begin{lemma}\label{lemma:alpha-lower-star-injective-script-I}
Let $\alpha\colon H\to G$ be an \emph{injective} group homomorphism. Then
\begin{equation*}
\alpha^*\colon\cat{$\bm G$-$\bm{\mathcal I}$-SSet}_{\textup{$G$-global}}\rightleftarrows
\cat{$\bm H$-$\bm{\mathcal I}$-SSet}_{\textup{$H$-global}} :\!\alpha_*
\end{equation*}
is a simplicial Quillen adjunction, and likewise for the positive model structures. If $(G:\im\alpha)<\infty$, then $\alpha_*$ is homotopical.
\begin{proof}
Let us first show that this is a Quillen adjunction. We already know that $\alpha^*$ is homotopical, so it suffices to show that the above is a Quillen adjunction for the strict level model structures, which follows in turn by applying Proposition~\ref{prop:alpha-lower-star-homotopical} levelwise.

Finally, if $(G:\im\alpha)<\infty$, then $\alpha_*$ is non-equivariantly just given by a finite product. Using that filtered colimits commute with finite products in $\cat{SSet}$, one concludes similarly to the argument from the previous lemma that $\alpha_*$ commutes with $\ev_\omega$. The claim follows as $\alpha_*\colon\cat{$\bm{E\mathcal M}$-$\bm H$-SSet}\to\cat{$\bm{E\mathcal M}$-$\bm G$-SSet}$ is homotopical by Corollary~\ref{cor:alpha-lower-star-injective-EM}.

The proof for the positive model structures is again analogous.\index{functoriality in homomorphisms!for G-II-SSet@for $\cat{$\bm G$-$\bm{\mathcal I}$-SSet}$|)}
\end{proof}
\end{lemma}

As an application of the calculus just developed we can now prove:

\begin{prop}\label{prop:injective-is-really-static-script-I}
Let $X$ be fibrant in the injective $G$-global model structure on $\cat{$\bm G$-$\bm{\mathcal I}$-SSet}$ and let $i\colon A\to B$ be a $G$-equivariant injection of (not necessarily finite) $G$-sets. Then $X(i)\colon X(A)\to X(B)$ is a proper $G$-equivariant weak equivalence with respect to the diagonal $G$-action.
\begin{proof}
Fix a finite subgroup $H\subset G$; we have to show that $X(i)^H\colon X(A)^H\to X(B)^H$ is a weak homotopy equivalence. By $2$-out-of-$3$ we may assume without loss of generality that $A=\varnothing$, and filtering $B$ by its finite $H$-subsets we may assume that $B$ itself is finite.

We now observe that $X$ is also fibrant in the $H$-global injective model structure by Lemma~\ref{lemma:alpha-shriek-injective-script-I}. On the other hand, by the Yoneda Lemma $X(i)^H$ agrees up to conjugation by isomorphisms with
\begin{equation*}
\Maps^H(p, X)\colon \Maps^H(*,X)\to\Maps^H(\mathcal I(B,\blank), X)
\end{equation*}
where $p\colon\mathcal I(B,\blank)\to *$ is the unique map, and $H$ acts on $\mathcal I(B,\blank)$ via $B$.

As the injective $H$-global model structure is simplicial and since all its objects are cofibrant, it therefore suffices that $p$ is an $H$-global weak equivalence, which by definition amounts to saying that $E\Inj(B,\omega)\to *$ is an $H$-global weak equivalence of $E\mathcal M$-$H$-simplicial sets. This is however just the content of Example~\ref{ex:G-globally-contractible}.
\end{proof}
\end{prop}

\index{functoriality in homomorphisms!for G-I-SSet@for $\cat{$\bm G$-$\bm I$-SSet}$|(}
Most of the above functoriality properties have analogues for the models based on $I$-simplicial sets. Let us demonstrate this for a selection of these:

\begin{cor}\label{cor:alpha-shriek-left-Quillen-I-SSet}
Let $\alpha\colon H\to G$ be any group homomorphism. Then
\begin{equation*}
\alpha_!\colon\cat{$\bm H$-$\bm{I}$-SSet}_{\textup{$H$-global}}\rightleftarrows\cat{$\bm G$-$\bm{I}$-SSet}_{\textup{$G$-global}} :\!\alpha^*
\end{equation*}
is a simplicial Quillen adjunction with fully homotopical right adjoint, and likewise for the corresponding positive model structures on either side.
\begin{proof}
One proves as in Lemma~\ref{lemma:alpha-shriek-projective-script-I} that these are Quillen adjunctions. To prove that $\alpha^*$ is homotopical, we consider the commutative diagram
\begin{equation*}
\begin{tikzcd}
\cat{$\bm G$-$\bm{I}$-SSet}\arrow[r, "\mathcal I\times_I\blank"]\arrow[d, "\alpha^*"'] & \cat{$\bm G$-$\bm{\mathcal I}$-SSet}\arrow[d, "\alpha^*"]\\
\cat{$\bm H$-$\bm{I}$-SSet}\arrow[r, "\mathcal I\times_I\blank"'] & \cat{$\bm H$-$\bm{\mathcal I}$-SSet}.
\end{tikzcd}
\end{equation*}
The horizontal arrows preserve and reflect weak equivalences by Theorem~\ref{thm:G-global-we-I-characterization} while the right hand vertical arrow is homotopical by Lemma~\ref{lemma:alpha-shriek-projective-script-I}; the claim follows immediately.
\end{proof}
\end{cor}

\begin{cor}
Let $\alpha\colon H\to G$ be any group homomorphism. Then
\begin{equation*}
\alpha^*\colon\cat{$\bm G$-$\bm{I}$-SSet}_{\textup{$G$-global injective}}\rightleftarrows\cat{$\bm H$-$\bm{I}$-SSet}_{\textup{$H$-global injective}} :\!\alpha_*
\end{equation*}
is a simplicial Quillen adjunction.
\begin{proof}
By the previous corollary $\alpha^*$ is homotopical and it obviously preserves injective cofibrations, so it is left Quillen.
\end{proof}
\end{cor}

\begin{cor}
Let $\alpha\colon H\to G$ be an \emph{injective} group homomorphism. Then
\begin{equation*}
\alpha_!\colon\cat{$\bm H$-$\bm{I}$-SSet}_{\textup{$H$-global injective}}\rightleftarrows\cat{$\bm G$-$\bm{I}$-SSet}_{\textup{$G$-global injective}} :\!\alpha^*
\end{equation*}
is a simplicial Quillen adjunction. In particular, $\alpha_!$ is homotopical.
\begin{proof}
One proves analogously to Lemma~\ref{lemma:alpha-shriek-injective-script-I} that $\alpha_!$ preserves injective cofibrations. We now consider the commutative square on the left in
\begin{equation*}
\begin{tikzcd}
\cat{$\bm H$-$\bm{I}$-SSet} &\arrow[l, "\forget"'] \cat{$\bm H$-$\bm{\mathcal I}$-SSet}\\
\cat{$\bm G$-$\bm{I}$-SSet}\arrow[u, "\alpha^*"] &\arrow[l, "\forget"] \cat{$\bm G$-$\bm{\mathcal I}$-SSet}\arrow[u, "\alpha^*"']
\end{tikzcd}
\qquad
\begin{tikzcd}
\cat{$\bm H$-$\bm{I}$-SSet}\arrow[d, "\alpha_!"']\arrow[r, "\mathcal I\times_I\blank"] & \cat{$\bm H$-$\bm{\mathcal I}$-SSet}\arrow[d, "\alpha_!"]\\
\cat{$\bm G$-$\bm{I}$-SSet}\arrow[r, "\mathcal I\times_I\blank"'] & \cat{$\bm G$-$\bm{\mathcal I}$-SSet}.
\end{tikzcd}
\end{equation*}
Passing to total mates yields a canonical isomorphism filling the square on the right. But the horizontal arrows in this preserve and reflect weak equivalences by Theorem~\ref{thm:G-global-we-I-characterization} while the right hand vertical arrow is homotopical by Lemma~\ref{lemma:alpha-shriek-injective-script-I}. The claim follows immediately.
\end{proof}
\end{cor}

\begin{rk}
By direct computation, the (homotopical) restriction functors $\alpha^*$ for any homomorphism $\alpha\colon H\to G$ are compatible with $\ev_\omega\colon\cat{$\bm G$-$\bm{\mathcal I}$-SSet}\to\cat{$\bm{E\mathcal M}$-$\bm G$-SSet}$, $\ev_\omega\colon\cat{$\bm G$-$\bm{I}$-SSet}\to\cat{$\bm{\mathcal M}$-$\bm G$-SSet}$, as well as all the forgetful functors.

It follows by abstract nonsense that $(\alpha^*)^\infty$ is actually compatible with all the equivalences of associated quasi-categories constructed above, and so are its adjoints $\cat{L}\alpha_!$ and $\cat{R}\alpha_*$.\index{functoriality in homomorphisms!for G-I-SSet@for $\cat{$\bm G$-$\bm I$-SSet}$|)}
\end{rk}

Using the characterization of the cofibrations given in Corollary~\ref{cor:characterization-i-cof} and the above functoriality properties, we can now prove:

\begin{thm}\label{thm:forget-left-quillen-global}
The simplicial adjunction
\begin{equation*}
\forget\colon\cat{$\bm G$-$\bm{\mathcal I}$-SSet}_{\textup{$G$-global}}\rightleftarrows\cat{$\bm G$-$\bm{I}$-SSet}_{\textup{$G$-global}} :\!\Maps_I(\mathcal I,\blank)
\end{equation*}
is a Quillen equivalence.
\begin{proof}
As the forgetful functor is homotopical and descends to an equivalence of associated quasi-categories (Theorem~\ref{thm:strict-global-I-model-structure}), it only remains to prove that it sends generating cofibrations to cofibrations.

\begin{claim*}
Let $A,B$ be finite sets and let $n\ge 0$. Then the latching map
\begin{equation*}
\colim_{C\subsetneq B} I(A,C)^{n+1}\to I(A,B)^{n+1}
\end{equation*}
is injective.
\begin{proof}
Let $(f_0,\dots,f_n)$ be a family of injections $A\to C$ and let $(f_0',\dots,f_n')$ be a family of injections $A\to C'$ for proper subsets $C,C'\subsetneq B$, such that both are sent to the same element of $I(A,B)^{n+1}$, i.e.~for each $a\in A$ and $i=0,\dots,n$
\begin{equation*}
C\ni f_i(a)=f_i'(a)\in C'.
\end{equation*}
We conclude that $f_i$ and $f_i'$ both factor through the same injection $f_i''\colon A\to C\cap C'$. But then obviously $(f_0,\dots,f_n)$ represents the same element of the colimit as $(f_0'',\dots,f_n'')$, and so does $(f_0',\dots,f_n')$, finishing the proof of the claim.
\end{proof}
\end{claim*}
Now let $A$ be a finite faithful $H$-set and let $B$ be any finite set. We can then view $I(A,B)^{n+1}$ as a $(\Sigma_B\times H)$-set. We claim that the isotropy group of \emph{any} $(f_0,\dots,f_n)\in I(A,B)^{n+1}$ is contained in $\mathcal G_{\Sigma_B,H}$. Indeed, this just amounts to saying that $H$ acts freely on $I(A,B)$ via its action on $A$, which is trivial to check.

We are now ready to finish the proof of the proposition: let $A$ be any finite faithful $H$-set, and let $B$ be any finite set. By the above claim, the latching map
\begin{equation}\label{eq:latching-base}
\ell_B\colon\colim_{C\subsetneq B} \mathcal I(A,C)\to\mathcal I(A,B)
\end{equation}
is injective, and the argument from the previous paragraph tells us in particular that any simplex not in the image has isotropy a graph subgroup of $\Sigma_B\times H$. This precisely means that $(\ref{eq:latching-base})$ is a cofibration for the $\mathcal G_{\Sigma_B,H}$-model structure on $\cat{$\bm{(\Sigma_B\times H)}$-SSet}$. We conclude from Corollary~\ref{cor:characterization-i-cof} that $\mathcal I(A,\blank)$ (with $H$ acting via $A$) is cofibrant in the $H$-global model structure on $\cat{$\bm{H}$-$\bm{I}$-SSet}$.

If now $\phi\colon H\to G$ is any homomorphism, then $\phi_!\colon\cat{$\bm H$-$\bm I$-SSet}\to\cat{$\bm G$-$\bm I$-SSet}$ is left Quillen for the $H$-global and $G$-global model structure, respectively, by Corollary~\ref{cor:alpha-shriek-left-Quillen-I-SSet}, so $\mathcal I(A,\blank)\times_\phi G\cong\phi_!\mathcal I(A,\blank)$ is cofibrant. As the $G$-global model structure is simplicial, we conclude that $\mathcal I(A,\blank)\times_\phi G\times\del\Delta^n\hookrightarrow\mathcal I(A,\blank)\times_\phi G\times\Delta^n$ is indeed a cofibration in $\cat{$\bm G$-$\bm I$-SSet}$ as desired.
\end{proof}
\end{thm}

\subsection{Another connection to monoid actions} We will now provide another comparison between $\cat{$\bm G$-$\bm{\mathcal I}$-SSet}$ and $\cat{$\bm{E\mathcal M}$-$\bm G$-SSet}$ in terms of a certain `reparametrization functor' appearing in the construction of global algebraic $K$-theory \cite[Constructions~3.3 and~8.2]{schwede-k-theory}.

\begin{constr}\label{constr:omega-bullet}\nomenclature[aomegabullet]{$(\blank)[\omega^\bullet]$}{reparametrization of $E\mathcal M$-action; $\mathcal I$-simplicial set built from an $E\mathcal M$-simplicial set this way}
We write $\overline{I}_\omega\subset\overline{I}$ for the full subcategory of countably infinite sets, and analogously $\overline{\mathcal I}_\omega\subset\overline{\mathcal I}$. Then the inclusion $B\mathcal M\to\overline I$ factors through an equivalence $i\colon B\mathcal M\to \overline{I}_\omega$. We pick once and for all a retraction $r$; this is then automatically quasi-inverse to $i$ and we fix $\tau\colon ir\cong\id$ with $\tau_\omega=\id_\omega$. The functor $r$ uniquely extends to a simplicially enriched functor $\overline{\mathcal I}_\omega\to BE\mathcal M$, which we denote by the same symbol; this is automatically a retraction of the inclusion $i$, and $\tau$ is a simplicially enriched isomorphism $ir\cong\id$.

For any $E\mathcal M$-$G$-simplicial set $X$, we now write $X[\blank]\mathrel{:=} X\circ r\colon\overline{\mathcal I}_\omega\to\cat{$\bm G$-SSet}$. We then define $X[\omega^\bullet]$ to be the following $G$-$\mathcal I$-simplicial set: if $A\not=\varnothing$, then $X[\omega^A]=X(r(\omega^A))$ as above. The $G$-action is as before, and for any injection $i\colon A\to B$ we take the structure map to be $X[i_!]\colon X[\omega^A]\to X[\omega^B]$, i.e.~it is given by applying $X\circ r$ to the `extension by zero map' $i_!\colon\omega^A\to\omega^B$ with
\begin{equation*}
i_!(f)(b)=\begin{cases}
f(a) & \text{if $b=i(a)$}\\
0 & \text{if $b\notin\im i$}.
\end{cases}
\end{equation*}
More generally, we let an $n$-simplex $(i_0,\dots,i_n)$ act by $X[i_{0!},\dots,i_{n!}]$. We remark that this means that as a $G$-simplicial set $X[\omega^A]=X$, and all of the above (higher) structure maps are given by acting with certain (inexplicit and mysterious) elements of $E\mathcal M$.

Finally, define $X[\omega^\varnothing]\mathrel{:=}X_{[\varnothing]}$. The structure maps $X[\omega^\varnothing]\to X=X[\omega^A]$ are given by the inclusions, and we choose all higher cells to be trivial.  As all the remaining structure is given by acting with elements of $E\mathcal M$, this is easily seen to be functorial, yielding a $G$-$\mathcal I$-simplicial set.

This extends to a simplicially enriched functor $\cat{$\bm{E\mathcal M}$-$\bm G$-SSet}\to\cat{$\bm G$-$\bm{\mathcal I}$-SSet}$ by sending an $n$-simplex $f\colon\Delta^n\times X\to Y$ of the mapping space to the transformation given in non-empty degree by $f$ itself and in degree $\varnothing$ by $f_{[\varnothing]}$.
\end{constr}

\begin{rk}
There is an alternative `coordinate-free' perspective on the above construction, that we briefly sketch; for the arguments in this paper we will however only be interested in the above version of the construction.

If $A$ is any non-empty set, then $E\Inj(\omega^A,\omega)$ is left $E\mathcal M$-isomorphic to $E\mathcal M$ by precomposing with the isomorphism $\tau$; for $A=\varnothing$, $E\Inj(\omega^\varnothing,\omega)\cong E\omega$ corepresents the functor $(\blank)_{[\{0\}]}$.

It is then not hard to produce a natural map from $X[\omega^\bullet]$ to the $G$-$\mathcal I$-simplicial set sending $A\in \mathcal I$ to $\Maps^{E\mathcal M}(E\Inj(\omega^A,\omega),X)$ with the functoriality in $A$ is as above. This map is an isomorphism in all positive degrees and hence in particular a $G$-global weak equivalence. This `coordinate free' description is then analogous to \cite[construction after Proposition~3.5]{schwede-orbi}, also cf.~\cite[Section~8]{lind}.
\end{rk}

We can now state our comparison:

\begin{prop}\label{prop:omega-bullet-inverse}
The functors
\begin{equation*}
\ev_\omega\colon \cat{$\bm G$-$\bm{\mathcal I}$-SSet}\rightleftarrows\cat{$\bm{E\mathcal M}$-$\bm G$-SSet} :\!(\blank)[\omega^\bullet]
\end{equation*}
are mutually inverse homotopy equivalences.
\end{prop}

For the proof we will need the following example of a complete $H$-set universe from \cite[Proposition~2.19]{schwede-k-theory}:

\begin{lemma}\label{lemma:exponential-universe}
Let $A$ be a finite $H$-set containing a free $H$-orbit. Then $\omega^A$ with left $H$-action via $(h.f)(a)=f(h^{-1}.a)$ is a complete $H$-set universe.\qed
\end{lemma}

\begin{rk}\label{rk:theta}
We can give an alternative description of $\ev_\omega\circ(\blank)[\omega^\bullet]$ in non-empty degrees as follows: for any $A\not=\varnothing$ the isomorphism $\tau\colon\omega\to\omega^A$ from the construction of $(\blank)[\omega^\bullet]$ induces
\begin{equation}\label{eq:theta-isomorphism}
X(\omega^A)\xrightarrow{X(\tau^{-1})} X(\omega)=X(\omega)[\omega^A]
\end{equation}
and these are by definition compatible with all the relevant (higher) structure maps and moreover natural in $X$.

Using this, we define $\theta_X\colon X\to X(\omega)[\omega^\bullet]$ in degree $A\not=\varnothing$ as the composition
\begin{equation*}
X(A)\xrightarrow{X(e)} X(\omega^A)\xrightarrow{(\ref{eq:theta-isomorphism})} X(\omega)[\omega^A],
\end{equation*}
where $e\colon A\to\omega^A$ sends $a\in A$ to its characteristic function, i.e.~$e(a)(a)=1$ and $e(a)(b)=0$ otherwise.
In degree $\varnothing$, we define $\theta_X(\varnothing)\colon X(\varnothing)\to X(\omega)_{[\varnothing]}$ via the unit of $\ev_\omega\dashv(\blank)_\bullet$. We omit the easy verification that $\theta_X$ is a map of $G$-$\mathcal I$-simplicial sets and natural in $X$.
\end{rk}

\begin{proof}[Proof of Proposition~\ref{prop:omega-bullet-inverse}]
We first show that $(\blank)[\omega^\bullet]$ is homotopical, for which we let $f$ be any $G$-global weak equivalence in $\cat{$\bm{E\mathcal M}$-$\bm G$-SSet}$. If now $H\subset\mathcal M$ is any universal subgroup, and $A$ is any finite $H$-set containing a free $H$-orbit (in particular $A\not=\varnothing$), then Lemma~\ref{lemma:exponential-universe} shows that $\omega^A$ is a complete $H$-set universe. We can therefore pick an $H$-equivariant isomorphism $\omega\cong\omega^A$, which shows that $f[\omega^A]$ agrees up to conjugation by $(H\times G)$-equivariant isomorphism with $f[\omega]=f$; in particular, $f[\omega^A]$ is a $\mathcal G_{H,G}$-weak equivalence. The argument from Lemma~\ref{lemma:strict-are-global} thus shows that $f[\omega^\bullet]$ is a $G$-global weak equivalence as desired.

Next, let us show that $\theta_X\colon X\to X(\omega)[\omega^\bullet]$ is a $G$-global weak equivalence for any $X\in\cat{$\bm G$-$\bm{\mathcal I}$-SSet}$. As $\ev_\omega$ and $(\blank)[\omega^\bullet]$ are homotopical, we may assume without loss of generality that $X$ is static. We then let $A$ be any finite non-empty faithful $H$-set, and observe that $\theta_X(A)$ agrees by the usual cofinality argument up to $(H\times G)$-equivariant isomorphism with the composition
\begin{equation*}
X(A)\xrightarrow{X(e)} X(e(A))\to\colim\limits_{e(A)\subset B\subset\omega^A\text{ finite $H$-set}} X(B),
\end{equation*}
where the unlabelled arrow is the structure map. As the left hand map is an isomorphism and all transition maps of the colimit are $\mathcal G_{H,G}$-weak equivalences since $X$ is static, this is a $\mathcal G_{H,G}$-weak equivalence, so $\theta_X$ is a $G$-global weak equivalence as before. In particular, $\ev_\omega$ is right homotopy inverse to $(\blank)[\omega^\bullet]$.

To see that it is also left homotopy inverse, we consider the following zig-zag for each $E\mathcal M$-$G$-simplicial set $Y$:
\begin{equation*}
Y[\omega^\bullet](\omega)=\colim\limits_{A\subset\omega\textup{ finite}} Y[\omega^A]
\xrightarrow{\alpha} \colim\limits_{A\subset\omega\textup{ finite}}
Y[\omega^A\amalg\omega] \xleftarrow{\beta}
\colim\limits_{A\subset\omega\textup{ finite}} Y[\omega]\cong Y.
\end{equation*}
The transition maps on the two left hand colimits come from the extension by zero maps $\omega^A\to\omega^B$, and the transition maps of the remaining colimit are trivial. The group $G$ acts by its action on $Y$ everywhere. Moreover, $E\mathcal M$ acts on all the colimits analogously to Construction~\ref{constr:i-extension} (observe that this part of the action is trivial for the rightmost colimit), in addition on the middle colimit by its tautological action on the $\omega$-summand of $Y[\omega^A\amalg\omega]$, and finally by its given action on $Y[\omega]=Y$ for the final colimit. The maps $\alpha$ and $\beta$ are given in each degree by the inclusions $\omega^A\hookrightarrow\omega^A\amalg\omega\hookleftarrow\omega$. One immediately checks that they are well-defined and $E\mathcal M$-$G$-equivariant. Moreover, we can make the middle term into a functor in $Y$ in the obvious way and with respect to this the maps $\alpha$ and $\beta$ are clearly natural.

It only remains to prove that $\alpha$ and $\beta$ are $G$-global weak equivalences. We prove this for $\alpha$, the other argument being similar. For this let $H\subset\mathcal M$ be universal. We pick a free $H$-orbit $F$ inside $\omega$; by the same argument as before, $\alpha$ agrees up to conjugation by $(H\times G)$-equivariant isomorphisms with the map
\begin{equation}\label{eq:alpha-H-G-focus}
\colim\limits_{F\subset A\subset\omega\textup{ finite $H$-set}} Y[\omega^A]\to \colim\limits_{F\subset A\subset\omega\textup{ finite $H$-set}} Y[\omega^A\amalg\omega]
\end{equation}
still induced in each degree by the inclusion $i\colon\omega^A\hookrightarrow\omega^A\amalg\omega$. We claim that each of these maps is even an $(H\times G)$-equivariant homotopy equivalence. But indeed, as both $\omega^A$ and $\omega^A\amalg\omega$ are complete $H$-set universes by Lemma~\ref{lemma:exponential-universe}, there exists an $H$-equivariant injection $j\colon\omega^A\amalg\omega\to\omega^A$, and $Y[j]$ is $(H\times G)$-equivariantly homotopy inverse to $Y[i]$ as witnessed by the equivariant homotopies $Y[ij,\id]$ and $Y[ji,\id]$.
\end{proof}

\subsection{Comparison to $\bm G$-equivariant homotopy theory}\label{subsec:G-global-vs-G-equiv-script-I}\index{proper G-equivariant homotopy theory@proper $G$-equivariant homotopy theory!vs G-global homotopy theory@vs.~$G$-global homotopy theory|(}
We can now prove the analogues of the results of Section~\ref{sec:g-global-vs-g-em} for $G$-$\mathcal I$-simplicial sets:

\begin{defi}
A map $f$ in $\cat{$\bm G$-$\bm{\mathcal I}$-SSet}$ is called an \emph{\emph{$\mathcal E$}-weak equivalence} if $f(\omega)$ is an $\mathcal E$-weak equivalence of $E\mathcal M$-$G$-simplicial sets (see Definition~\ref{defi:class-E}).
\end{defi}

\begin{cor}\label{cor:G-script-I-SSet-vs-G-SSet}
The homotopical functor
\begin{equation*}
\const\colon\cat{$\bm G$-SSet}_{\textup{proper}}\to\cat{$\bm G$-$\bm{\mathcal I}$-SSet}_{\textup{$G$-global}}
\end{equation*}
induces a fully faithful functor on associated quasi-categories. This induced functor admits both a left adjoint $\textbf{\textup L}{\colim_{\mathcal I}}$ as well as a right adjoint $\textbf{\textup R}\ev_\varnothing$. The latter is a quasi-localization at the $\mathcal E$-weak equivalences, and it in turn admits another right adjoint $\mathcal R$, which is again fully faithful.

Finally, the diagram
\begin{equation}\label{diag:triv-vs-const}
\begin{tikzcd}[column sep=small]
& \cat{$\bm G$-SSet}_{\textup{proper}}^\infty\arrow[dl, bend right=10pt, "\const^\infty"']\arrow[dr, bend left=10pt, "\triv_{E\mathcal M}^\infty"]\\
\cat{$\bm G$-$\bm{\mathcal I}$-SSet}_{\textup{$G$-global}}^\infty\arrow[rr, "(\ev_\omega)^\infty"'] &&
\cat{$\bm G$-$\bm{E\mathcal M}$-SSet}_{\textup{$G$-global}}^\infty
\end{tikzcd}
\end{equation}
commutes up to canonical equivalence.
\begin{proof}
The adjunction
\begin{equation*}
\const\colon\cat{$\bm G$-SSet}_{\textup{proper}}\rightleftarrows\cat{$\bm G$-$\bm{\mathcal I}$-SSet}_{\textup{$G$-global injective}} :\!\ev_\varnothing
\end{equation*}
is easily seen to be a Quillen adjunction, providing the above description of the right adjoint. Similarly, for the left adjoint we want to prove that the simplicial adjunction
\begin{equation}\label{eq:colim-vs-const-G-I}
\colim_{\mathcal I}\colon\cat{$\bm G$-$\bm{\mathcal I}$-SSet}_{\textup{$G$-global}}\rightleftarrows\cat{$\bm G$-SSet}_{\textup{proper}} :\!\const
\end{equation}
is a Quillen adjunction. For this we first observe that this is a Quillen adjunction when we equip the left hand side with the strict level model structure (as $\const$ is then obviously right Quillen). In particular, $\colim_{\mathcal I}$ preserves $G$-global (i.e.~strict level) cofibrations, and $\const$ sends fibrant $G$-simplicial sets to strictly fibrant $G$-$\mathcal I$-simplicial sets. It then suffices by Proposition~\ref{prop:cofibrations-fibrant-qa} and the characterization of the fibrant objects in $\cat{$\bm G$-$\bm{\mathcal I}$-SSet}_{\textup{$G$-global}}$ provided by Theorem~\ref{thm:script-I-global-model-structure} that $\const$ sends fibrant objects to static ones, which is immediate from the definitions.

In order to construct a canonical equivalence filling $(\ref{diag:triv-vs-const})$, it suffices to observe that the evident diagram of homotopical functors inducing it even commutes up to canonical isomorphism. The remaining statements then follow formally from Theorem~\ref{thm:G-global-vs-proper-sset} and commutativity of $(\ref{diag:triv-vs-const})$ as $(\ev_\omega)^\infty$ is an equivalence.
\end{proof}
\end{cor}

In the world of $G$-$\mathcal I$-simplicial sets there are rather explicit point-set models of the right adjoints $\textbf{R}\ev_\varnothing$ and $\mathcal R$, which we will introduce now. To this end, we define
\begin{equation}\label{eq:proper-set-universe}
\mathcal U_G\mathrel{:=}\coprod_{i=0}^\infty\coprod_{H\subset G,H\textup{ finite}} G/H.
\end{equation}
\nomenclature[aUG]{$\mathcal U_G$}{union of countably infinitely many copies of all transitive $G$-sets with finite isotropy up to isomorphism (replacement for complete set universe in the proper context)}

If $G$ is finite, the above is just a particular construction of a complete $G$-set universe. However, in general this need not be countable and it can even have uncountably many orbits (e.g.~for $G=\bigoplus_{\mathbb R}\mathbb Z/2\mathbb Z$). Despite these words of warning we always have:

\begin{lemma}\label{lemma:proper-set-universe}
Let $H$ be any finite group and let $\iota\colon H\to G$ be any injective group homomorphism. Then the $H$-set $\iota^*\mathcal U_G$ contains a complete $H$-set universe.
\begin{proof}
We may assume without loss of generality that $\iota$ is literally the inclusion of a finite subgroup $H$ of $G$. We then have an $H$-equivariant injection
\begin{equation*}
\coprod_{i=0}^\infty\coprod_{K\subset H} H/K\to\mathcal U_G
\end{equation*}
induced from $H\hookrightarrow G$ (and using that each subgroup $K\subset H$ is in particular a finite subgroup of $G$). As the left hand side is a complete $H$-set universe, this finishes the proof.
\end{proof}
\end{lemma}

By the universal property of enriched presheaves, the simplicial functor $\ev_{\mathcal U_G}$ given by evaluation at $\mathcal U_G$ has a simplicial right adjoint given explicitly by
\nomenclature[aR]{$R$}{point-set level model of $\mathcal R$ for $G$-$\mathcal I$-simplicial sets}
\begin{equation*}
R(X)(A)=\Maps(E\Inj(A,\mathcal U_G),X)
\end{equation*}
(where $\Maps$ denotes the simplicial set of all maps, with $G$ acting by conjugation) with the obvious functoriality in each variable. We can now prove:

\begin{prop}\label{prop:evaluation-proper-vs-evaluation-empty}
The simplicial adjunction
\begin{equation}\label{eq:evaluate-proper-universe}
\und_G\mathrel{:=}\ev_{\mathcal U_G}\colon\cat{$\bm G$-$\bm{\mathcal I}$-SSet}_{\textup{$G$-global}}\rightleftarrows\cat{$\bm G$-SSet}_{\textup{proper}} :\!R
\end{equation}
\nomenclature[aug]{$\und_G$}{underlying proper $G$-equivariant space of a $G$-$\mathcal I$-simplicial set}%
is a Quillen adjunction with homotopical left adjoint, and there are preferred equivalences
\begin{equation*}
\big(\ev_{\mathcal U_G}\big)^\infty\simeq\textbf{\textup{R}}\ev_\varnothing\qquad\text{and}\qquad \textbf{\textup R}R\simeq\mathcal R.
\end{equation*}
\begin{proof}
Because of the canonical isomorphism $\ev_{\mathcal U_G}\circ\const\cong\id$ it suffices that $(\ref{eq:evaluate-proper-universe})$ is a Quillen adjunction and that $\ev_{\mathcal U_G}$ is homotopical in $\mathcal E$-weak equivalences.

If $H\subset G$ is any finite subgroup, then we pick an injective homomorphism $\iota\colon H'\to G$ with image $H$ from a universal subgroup $H'\subset\mathcal M$. With this notation we then have for any $G$-$\mathcal I$-simplicial set $X$ an actual equality
\begin{equation*}
\big(X(\mathcal U_G)\big)^H=X\big(\iota^*(\mathcal U_G)\big)^\iota,
\end{equation*}
and analogously for morphisms. By the previous lemma there exists an $H'$-equi\-variant embedding $\omega\to \iota^*(\mathcal U_G)$ (with respect to the tautological action on the left hand side), and by Lemma~\ref{lemma:evaluation-h-universe} we conclude that the induced natural transformation $(\blank)^\iota\circ\ev_\omega\Rightarrow(\blank)^H\circ\ev_{\mathcal U_G}$ is a weak equivalence. Thus, $\ev_{\mathcal U_G}$ is homotopical in $\mathcal E$-weak equivalences (and hence in particular in $G$-global weak equivalences).

It only remains to prove that $\ev_{\mathcal U_G}$ sends the standard generating cofibrations to proper cofibrations. As the proper model structure is simplicial, we are reduced to showing that the $G$-simplicial set $E\Inj(A,\mathcal U_G)\times_\phi G$ is cofibrant in the proper model structure (i.e. has finite isotropy groups) for any finite group $H$, any finite (faithful) $H$-set $A$ and any group homomorphism $\phi\colon H\to G$. For this we let $(f_0,\dots,f_n;g)$ represent any $n$-simplex. If $g'$ fixes $[f_0,\dots,f_n;g]$, then in particular $g'g=g\phi(h)$ for some $h\in H$, and hence $g'=g\phi(h)g^{-1}$. As the right hand side can only take finitely many values for any fixed $g$, the claim follows.
\end{proof}
\end{prop}

\begin{rk}\label{rk:R-vs-R}
For $G=1$ the above adjunction is the $\mathcal I$-analogue of \cite[Remark~1.2.24 and Proposition~1.2.27]{schwede-book}. As we will show in Section~\ref{sec:global-vs-g-global}, there exists a zig-zag of homotopical functors
\begin{equation*}
\cat{$\bm{\mathcal I}$-SSet}\xrightarrow{\forget}\cat{$\bm I$-SSet}\xrightarrow{|\blank|}\cat{$\bm I$-Top}\xleftarrow{\text{forget}}\cat{$\bm{L}$-Top}
\end{equation*}
where $\cat{$\bm L$-Top}$ denotes Schwede's orthogonal spaces
(with respect to the global weak equivalences for the class of \emph{finite} groups), and all of these induce equivalences on associated quasi-categories. It is then not hard to check directly that under this identification the adjunction $\const^\infty\dashv\cat{R}\ev_\varnothing\dashv\mathcal R$ corresponds to the adjunction $\textbf{L}_1\dashv\textbf{R}\ev_0\dashv\textbf{R}R$ considered in \emph{loc.~cit.} for the trivial group.

In fact, there is also a completely abstract way to prove this: namely, in both cases the leftmost adjoint under consideration preserves the terminal object (in our case because we have explicitly constructed a further left adjoint, in Schwede's case by direct inspection). If we now have \emph{any} equivalence $\Phi$ between $\cat{$\bm L$-Top}^\infty$ and $\cat{$\bm{\mathcal I}$-SSet}^\infty$, then both ways through the diagram
\begin{equation*}
\begin{tikzcd}[column sep=small]
& \cat{SSet}^\infty\arrow[dl, "\const^\infty"', bend right=10pt]\arrow[dr, "\textbf{L}_1", bend left=10pt]\\
\cat{$\bm{\mathcal I}$-SSet}^\infty\arrow[rr, "\Phi"'] & & \cat{$\bm L$-Top}^\infty
\end{tikzcd}
\end{equation*}
are cocontinuous and send the terminal object to the terminal object. By the universal property of spaces \cite[Theorem~5.1.5.6]{htt}, there is therefore a contractible space of natural equivalences filling this, and as before $\Phi$ is then also compatible with the other adjoints.\index{proper G-equivariant homotopy theory@proper $G$-equivariant homotopy theory!vs G-global homotopy theory@vs.~$G$-global homotopy theory|)}
\end{rk}

\subsection{Model structures for tame actions}\label{subsec:tame-model-structures}
Above, we have used our understanding of tame $\mathcal M$- and $E\mathcal M$-actions to introduce the $G$-global model structures on $\cat{$\bm G$-$\bm I$-SSet}$ and $\cat{$\bm G$-$\bm{\mathcal I}$-SSet}$. Conversely, we will now use the results of this section to construct $G$-global model structures on $\cat{$\bm{E\mathcal M}$-$\bm G$-SSet}^\tau$ and $\cat{$\bm{\mathcal M}$-$\bm G$-SSet}^\tau$:

\begin{thm}\index{positive G-global model structure@positive $G$-global model structure!on EM-G-SSettau@on $\cat{$\bm{E\mathcal M}$-$\bm G$-SSet}^\tau$|textbf}\label{thm:pos-G-global-EM-tau}
There exists a unique model structure on $\cat{$\bm{E\mathcal M}$-$\bm G$-SSet}^\tau$ in which a map $f$ is a weak equivalence, fibration, or cofibration if and only if $f_\bullet$ is a weak equivalence, fibration, or cofibration, respectively, in the positive $G$-global model structure on $\cat{$\bm G$-$\bm{\mathcal I}$-SSet}$. We call this the \emph{positive $G$-global model structure}. Its weak equivalences are precisely the $G$-global weak equivalences, and hence they are in particular closed under filtered colimits.

This model structure is proper, simplicial, and combinatorial with generating cofibrations
\begin{equation}\label{eq:EM-tau-generating-cof}
(E\Inj(A,\omega)\times_\phi G)\times\del\Delta^n\hookrightarrow(E\Inj(A,\omega)\times_\phi G)\times\Delta^n
\end{equation}
where $H$ runs through finite groups, $A\not=\varnothing$ is a finite faithful $H$-set, $\phi$ is a homomorphism $H\to G$, and $n\ge 0$. Finally, the simplicial adjunctions
\begin{align}\label{eq:em-tau-vs-em}
\incl\colon\cat{$\bm{E\mathcal M}$-$\bm G$-SSet}^\tau&\rightleftarrows\cat{$\bm{E\mathcal M}$-$\bm G$-SSet}_{\textup{injective $G$-global}} :\!(\blank)^\tau\\
\intertext{and}
\label{eq:em-tau-vs-script-I}
\ev_\omega\colon\cat{$\bm G$-$\bm{\mathcal I}$-SSet}&\rightleftarrows\cat{$\bm{E\mathcal M}$-$\bm G$-SSet}^\tau :\!(\blank)_\bullet
\end{align}
are Quillen equivalences.
\begin{proof}\begingroup\parskip=\the\parskip plus .76pt
\baselineskip=\the\baselineskip plus .1pt
The adjunction $\ev_\omega\dashv(\blank)_\bullet$ exhibits $\cat{$\bm{E\mathcal M}$-$\bm G$-SSet}^\tau$ as accessible Bousfield localization of the locally presentable category $\cat{$\bm G$-$\bm{\mathcal I}$-SSet}$, so it is itself locally presentable, cf.~\cite[Remark~5.5.1.6]{htt}.

Let us now show that the above defines a model structure on $\cat{$\bm{E\mathcal M}$-$\bm G$-SSet}^\tau$. Instead of using Crans' Transfer Criterion (Proposition~\ref{prop:transfer-criterion}), we can actually easily verify the model structure axioms directly here as the adjunction $\ev_\omega\dashv(\blank)_\bullet$ is already well-behaved $1$-categorically, also see \cite[Theorem~5.10 and Corollary~5.11]{I-vs-M-1-cat} where a similar argument is used to non-equivariantly relate $\cat{$\bm I$-SSet}$ and $\cat{$\bm{\mathcal M}$-SSet}^\tau$:

The $2$-out-of-$3$ property for weak equivalences as well as the closure under retracts for all three classes are obvious. Moreover, as $(\blank)_\bullet$ is fully faithful, the lifting axioms are inherited from the lifting axioms for $\cat{$\bm{G}$-$\bm{\mathcal I}$-SSet}$.

It only remains to verify the factorization axioms, for which we let $f\colon X\to Y$ be any map of tame $E\mathcal M$-$G$-simplicial sets. Then we can factor $f_\bullet$ as a cofibration $i\colon X_\bullet\to Z$ followed by an acyclic fibration $p\colon Z\to Y_\bullet$. But cofibrations in the positive $G$-global model structure are in particular cofibrations in $\cat{$\bm I$-SSet}$ (Theorem~\ref{thm:forget-left-quillen-global} together with Lemma~\ref{lemma:alpha-lower-star-injective-script-I}), so $Z$ is flat again and hence lies in the essential image of $(\blank)_\bullet$. We can therefore assume without loss of generality that $Z=Z'_\bullet$ for some $Z'\in\cat{$\bm{E\mathcal M}$-$\bm G$-SSet}^\tau$. By full faithfulness of $(\blank)_\bullet$ we can then write $i=i'_\bullet$, $p=p'_\bullet$, which yields the desired factorization $f=p'i'$. The remaining factorization axiom is proven analogously.

This completes the proof of the existence of the positive $G$-global model structure. By definition, $f\colon X\to Y$ is a weak equivalence if and only if $f_\bullet$ is, which in turn is equivalent by definition to $f_\bullet(\omega)$ being a $G$-global weak equivalence of $E\mathcal M$-$G$-simplicial sets. As $(\blank)_\bullet$ is fully faithful, $f_\bullet(\omega)$ is conjugate to $f$, which shows that the weak equivalences are precisely the $G$-global weak equivalences. In particular, they are closed under filtered colimits (which are created in $\cat{$\bm{E\mathcal M}$-$\bm G$-SSet}$).

The model structure is right proper since it is transferred from a right proper model stucture, see Lemma~\ref{lemma:transferred-properties}-$(\ref{item:tpr-proper})$. Moreover, $\cat{$\bm{E\mathcal M}$-$\bm G$-SSet}^\tau$ is tensored and cotensored over $\cat{SSet}$ with the tensoring given by the tensoring on $\cat{$\bm{E\mathcal M}$-$\bm G$-SSet}$ and the cotensoring given by applying $(\blank)^\tau$ to the usual cotensoring. As $(\ref{eq:em-tau-vs-script-I})$ is a simplicial adjunction, the positive $G$-global model structure on $\cat{$\bm{E\mathcal M}$-$\bm G$-SSet}^\tau$ is then simplicial by Lemma~\ref{lemma:transferred-properties}-$(\ref{item:tpr-simplicial})$.

It is clear from the construction that $(\blank)_\bullet$ is right Quillen so that $(\ref{eq:em-tau-vs-script-I})$ is a Quillen adjunction. Moreover, both adjoints preserve and reflect weak equivalences by definition and the above characterization of the weak equivalences. To show that $(\ref{eq:em-tau-vs-script-I})$ is a Quillen equivalence it is therefore enough that the counit $X_\bullet(\omega)\to X$ be a $G$-global weak equivalence for each tame $E\mathcal M$-$G$-simplicial set $X$, but we already know that it is even an isomorphism.

Now let $I_\cat{$\bm G$-$\bm{\mathcal I}$-SSet}$, $J_\cat{$\bm G$-$\bm{\mathcal I}$-SSet}$ be sets of generating cofibrations and generating acyclic cofibrations, respectively, of the positive $G$-global model structure on $\cat{$\bm G$-$\bm{\mathcal I}$-SSet}$. We claim that the positive $G$-global model structure on $\cat{$\bm{E\mathcal M}$-$\bm G$-SSet}^\tau$ is cofibrantly generated (and hence combinatorial) with generating cofibrations $I_\cat{$\bm G$-$\bm{\mathcal I}$-SSet}(\omega)\mathrel{:=}\{i(\omega) : i\in I_\cat{$\bm G$-$\bm{\mathcal I}$-SSet}\}$ and with generating acyclic cofibrations $J_\cat{$\bm G$-$\bm{\mathcal I}$-SSet}(\omega)$. Indeed, $I_\cat{$\bm G$-$\bm{\mathcal I}$-SSet}(\omega)$ and $J_\cat{$\bm G$-$\bm{\mathcal I}$-SSet}(\omega)$ permit the small object argument as $\cat{$\bm{E\mathcal M}$-$\bm G$-SSet}^\tau$ is locally presentable, and they detect acyclic fibrations and fibrations, respectively, by adjointness. Taking $I_\cat{$\bm G$-$\bm{\mathcal I}$-SSet}$ to be the usual set of generating cofibrations and using the canonical isomorphism $\mathcal I(A,\blank)(\omega)\cong E\Inj(A,\omega)$ shows that $(\ref{eq:EM-tau-generating-cof})$ is a set of generating cofibrations.

Next, let us show that $(\ref{eq:em-tau-vs-em})$ is a simplicial Quillen adjunction. It is obvious that the left adjoint preserves tensors, so that this is indeed a simplicial adjunction. Moreover, $\incl$ sends the above generating cofibrations to injective cofibrations, and it is homotopical by the above characterization of the weak equivalences, hence left Quillen. Finally, the inclusion descends to an equivalence on homotopy categories by Theorem~\ref{thm:shrew}, i.e.~$(\ref{eq:em-tau-vs-em})$ is a Quillen equivalence.

It only remains to show that $\cat{$\bm{E\mathcal M}$-$\bm G$-SSet}^\tau$ is left proper. But indeed, we have seen that the weak equivalences are precisely the $G$-global weak equivalences, and since $(\ref{eq:em-tau-vs-em})$ is a Quillen adjunction, the cofibrations are in particular injective cofibrations. As pushouts in $\cat{$\bm{E\mathcal M}$-$\bm G$-SSet}^\tau$ can be computed inside all of $\cat{$\bm{E\mathcal M}$-$\bm G$-SSet}$, the claim therefore follows from the left properness of the injective $G$-global model structure on $\cat{$\bm{E\mathcal M}$-$\bm G$-SSet}$.
\endgroup
\end{proof}
\end{thm}

\begin{rk}
While we will be almost exclusively interested in the positive $G$-global model structure considered above, the analogous statement for the $G$-global model structure on $\cat{$\bm G$-$\bm{\mathcal I}$-SSet}$ holds and can be proven in the same way. We call the resulting model structure the \emph{$G$-global model structure}\index{G-global model structure@$G$-global model structure!on EM-G-SSettau@on $\cat{$\bm{E\mathcal M}$-$\bm G$-SSet}^\tau$}; its generating cofibrations are again given by evaluating the usual generating cofibrations for the corresponding model structure on $\cat{$\bm G$-$\bm{\mathcal I}$-SSet}$ at $\omega$.
\end{rk}

\begin{rk}\label{rk:support-em-tau-model-structure}
Similarly to Corollary~\ref{cor:injective-almost-tame}, if $X$ is fibrant in the $G$-global positive model structure, then all its simplices have `small support up to weak equivalence.' More precisely, let $H\subset\mathcal M$ be universal and let $A\subset\omega$ be a non-empty finite faithful $H$-set. Then the $H$-action on $X$ restricts to $X_{[A]}$, and the inclusion $X_{[A]}\hookrightarrow X$ is a $\mathcal G_{H,G}$-weak equivalence as $X_\bullet$ is fibrant in the positive $G$-global model structure on $\cat{$\bm G$-$\bm{\mathcal I}$-SSet}$, hence in particular positively static.

Of course, the analogous statement for the $G$-global model structure holds.
\end{rk}

For later use, we record two properties of the above cofibrations:

\begin{lemma}\label{lemma:g-global-pos-cof}
Let $f\colon X\to Y$ be a map in $\cat{$\bm{E\mathcal M}$-$\bm G$-SSet}^\tau$.
\begin{enumerate}
\item If $f$ is a cofibration in the $G$-global positive model structure, then $f$ restricts to an isomorphism $f_{[\varnothing]}\colon X_{[\varnothing]}\to Y_{[\varnothing]}$.
\item If $f$ is a $G$-global cofibration (for example, if $f$ is a positive $G$-global cofibration), then it is also a $\mathcal G_{H,G}$-cofibration for any (finite) subgroup $H\subset\mathcal M$, i.e.~$f$ is injective and $G$ acts freely on $Y$ outside the image of $f$.
\end{enumerate}
\begin{proof}
For the first statement we observe that $f_\bullet$ is a $G$-global positive cofibration by definition, so $f_{[\varnothing]}=f_\bullet(\varnothing)$ is an isomorphism by Lemma~\ref{lemma:pos-cof-script-I}.

For the second statement, it suffices to prove this for the generating cofibrations. As the $\mathcal G_{H,G}$-model structure is simplicial, it suffices further to show that $E\Inj(A,\omega)\times_\phi G$ is cofibrant in the $\mathcal G_{H,G}$-model structure, i.e.~that $G$ acts freely on it. This is immediate from Lemma~\ref{lemma:emg-basic-general}-(\ref{item:emgbg-unique}).
\end{proof}
\end{lemma}

While the above argument for the construction of the positive $G$-global model structure does not apply for the injective $G$-global model structure on $\cat{$\bm G$-$\bm{\mathcal I}$-SSet}$, we still have:

\begin{thm}
There is a unique model structure on $\cat{$\bm{E\mathcal M}$-$\bm G$-SSet}^\tau$ with cofibrations the injective cofibrations and weak equivalences the $G$-global weak equivalences. We call this the \emph{injective $G$-global model structure}.\index{injective G-global model structure@injective $G$-global model structure!on EM-G-SSettau@on $\cat{$\bm{E\mathcal M}$-$\bm G$-SSet}^\tau$|textbf}  It is combinatorial, proper, simplicial, and filtered colimits in it are homotopical. Moreover, the simplicial adjunction
\begin{equation}\label{eq:ev-omega-quillen-tame}
\ev_\omega\colon\cat{$\bm G$-$\bm{\mathcal I}$-SSet}_{\textup{inj.~$G$-global}}\rightleftarrows\cat{$\bm{E\mathcal M}$-$\bm G$-SSet}^\tau_{\textup{inj.~$G$-global}} :\!(\blank)_\bullet
\end{equation}
is a Quillen equivalence.
\begin{proof}
Let $\hat I$ be a set of generating cofibrations for the injective $G$-global model structure on $\cat{$\bm G$-$\bm{\mathcal I}$-SSet}$, and define $I\mathrel{:=}\hat I(\omega)=\{i(\omega):i\in \hat I\}$.

\begin{claim*}
The injective cofibrations are precisely the retracts of relative $I$-cell complexes.
\begin{proof}
It is clear that any element of $I$ is an injective cofibration. As the latter are closed under pushouts, retracts, and transfinite composition, it suffices to show conversely that any injective cofibration is a retract of an $I$-cell complex.

But indeed, if $f$ is an injective cofibration, then so is $f_\bullet$ by direct inspection, so that it can be written as a retract of a relative $\hat I$-cell complex. As $\ev_\omega$ is a left adjoint, we conclude that $f_\bullet(\omega)$ is a retract of a relative $\hat I(\omega)=I$-cell complex. But by full faithfulness of $(\blank)_\bullet$, $f_\bullet(\omega)$ is conjugate to $f$, which completes the proof.
\end{proof}
\end{claim*}

We have seen in the proof of Theorem~\ref{thm:pos-G-global-EM-tau} that pushouts along injective cofibrations in $\cat{$\bm{E\mathcal M}$-$\bm G$-SSet}^\tau$ preserve $G$-global weak equivalences, and that positive $G$-global cofibrations are in particular injective cofibrations. As the injective cofibrations are generated by the set $I$, Corollary~\ref{cor:enlarge-generating-cof} therefore shows that the model structure exists, that is combinatorial and proper, and that filtered colimits in it are homotopical.

To prove that the model structure is also simplicial, it suffices to observe that colimits and tensors can be computed in all of $\cat{$\bm{E\mathcal M}$-$\bm G$-SSet}$ and that the latter is simplicial by Proposition~\ref{prop:equivariant-injective-model-structure}.

Finally, it is clear that $\ev_\omega$ preserves injective cofibrations and weak equivalences, so that $(\ref{eq:ev-omega-quillen-tame})$ is a Quillen adjunction, hence a Quillen equivalence by Theorem~\ref{thm:pos-G-global-EM-tau}.
\end{proof}
\end{thm}

Next, we come to an analogue of Theorem~\ref{thm:pos-G-global-EM-tau} for tame $\mathcal M$-actions, which can be proven in exactly the same way:

\begin{thm}\label{thm:M-G-tau-model-structure}
There exists a model structure on $\cat{$\bm{\mathcal M}$-$\bm G$-SSet}^\tau$ in which a map $f$ is a weak equivalence, fibration, or cofibration, if and only if $f_\bullet$ is a weak equivalence, fibration, or cofibration, respectively, in the positive $G$-global model structure on $\cat{$\bm G$-$\bm{I}$-SSet}$. We call this the \emph{positive $G$-global model structure}.\index{positive G-global model structure@positive $G$-global model structure!on M-G-SSettau@on $\cat{$\bm{\mathcal M}$-$\bm G$-SSet}^\tau$|textbf}  Its weak equivalences are precisely the $G$-global weak equivalences.

This model structure is left proper, simplicial, combinatorial with generating cofibrations
\begin{equation*}
\begin{aligned}
&\{(\Inj(A,\omega)\times_\phi G)\times\del\Delta^n\hookrightarrow(\Inj(A,\omega)\times_\phi G)\times\Delta^n :{}\\
&\qquad\text{$H$ finite group, $A\not=\varnothing$ finite faithful $H$-set, $\phi\colon H\to G$ homomorphism}\},
\end{aligned}
\end{equation*}
and filtered colimits in it are homotopical. Finally, the simplicial adjunctions
\begin{align*}
\incl\colon\cat{$\bm{\mathcal M}$-$\bm G$-SSet}^\tau&\rightleftarrows\cat{$\bm{\mathcal M}$-$\bm G$-SSet}_{\textup{injective $G$-global}} :\!(\blank)^\tau\\
\ev_\omega\colon\cat{$\bm G$-$\bm{I}$-SSet}&\rightleftarrows\cat{$\bm{\mathcal M}$-$\bm G$-SSet}^\tau :\!(\blank)_\bullet
\end{align*}
are Quillen equivalences.\qed
\end{thm}

\begin{cor}
The simplicial adjunction
\begin{equation*}
E\mathcal M\times_{\mathcal M}\blank\colon\cat{$\bm{\mathcal M}$-$\bm G$-SSet}^\tau\rightleftarrows\cat{$\bm{E\mathcal M}$-$\bm G$-SSet}^\tau :\forget
\end{equation*}
is a Quillen equivalence, and both adjoints are homotopical.
\begin{proof}
We have seen in Theorem~\ref{thm:tame-M-sset-vs-EM-sset} that both adjoints are homotopical, and that they descend to equivalences on homotopy categories. It only remains to show that $E\mathcal M\times_{\mathcal M}\blank$ sends the above generating cofibrations to cofibrations, which is immediate from  Corollary~\ref{cor:E-Inj-corepr}.
\end{proof}
\end{cor}

\begin{rk}
Again we get an analogous result for the $G$-global model structure\index{G-global model structure@$G$-global model structure!on M-G-SSettau@on $\cat{$\bm{\mathcal M}$-$\bm G$-SSet}^\tau$}. Moreover, one can construct an injective $G$-global model structure\index{injective G-global model structure@injective $G$-global model structure!on M-G-SSettau@on $\cat{$\bm{\mathcal M}$-$\bm G$-SSet}^\tau$} by an argument similar to the above. We leave the details to the interested reader.
\end{rk}

Finally, let us discuss functoriality for $\cat{$\bm{E\mathcal M}$-$\bm G$-SSet}^\tau$:

\index{functoriality in homomorphisms!for EM-G-SSettau@for $\cat{$\bm{E\mathcal M}$-$\bm G$-SSet}^\tau$|(}
\begin{lemma}\label{lemma:alpha-lower-shriek-tau}
Let $\alpha\colon H\to G$ be any group homomorphism. Then
\begin{equation*}
\alpha_!\colon\cat{$\bm{E\mathcal M}$-$\bm H$-SSet}\rightleftarrows\cat{$\bm{E\mathcal M}$-$\bm G$-SSet} :\!\alpha^*
\end{equation*}
restricts to a Quillen adjunction
\begin{equation}\label{eq:alpha-shriek-tau}
\alpha_!\colon\cat{$\bm{E\mathcal M}$-$\bm H$-SSet}^\tau_{\textup{positive $H$-global}}\rightleftarrows\cat{$\bm{E\mathcal M}$-$\bm G$-SSet}^\tau_{\textup{positive $G$-global}} :\!\alpha^*.
\end{equation}
The right adjoint is fully homotopical, and if $\alpha$ is injective, so is the left adjoint.
\begin{proof}
It is clear that $\alpha^*$ preserves tameness, and so does $\alpha_!$ as the full subcategory $\cat{$\bm{E\mathcal M}$-SSet}^\tau\subset\cat{$\bm{E\mathcal M}$-SSet}$ is closed under colimits.

To see that $(\ref{eq:alpha-shriek-tau})$ is a Quillen adjunction with homotopical right adjoint, it suffices to observe that $\alpha^*$ commutes with $(\blank)_\bullet$ on the nose, so that the claim follows from Lemma~\ref{lemma:alpha-shriek-projective-script-I}. Finally, if $\alpha$ is injective, then $\alpha_!$ sends $H$-global weak equivalences to $G$-global weak equivalences by Corollary~\ref{cor:alpha-shriek-injective-EM}.
\end{proof}
\end{lemma}

Similarly, one deduces from Corollary~\ref{cor:free-quotient-EM}:

\begin{cor}\label{cor:free-quotient-tame}
In the above situation, $\alpha_!$ preserves $H$-global weak equivalences between objects with free $\ker(\alpha)$-action.\qed
\end{cor}

The situation for the right adjoint is a bit more complicated: of course, $\alpha^*$ still has a right adjoint $\alpha_*$, and $\alpha^*\dashv\alpha_*$ is a Quillen adjunction for the injective model structures---however, $\alpha_*\colon\cat{$\bm{E\mathcal M}$-$\bm H$-SSet}^\tau\to\cat{$\bm{E\mathcal M}$-$\bm G$-SSet}^\tau$ will usually not be given as restriction of $\alpha_*\colon\cat{$\bm{E\mathcal M}$-$\bm H$-SSet}\to\cat{$\bm{E\mathcal M}$-$\bm G$-SSet}$ because the inclusion $\cat{$\bm{E\mathcal M}$-SSet}^\tau\hookrightarrow\cat{$\bm{E\mathcal M}$-SSet}$ does not preserve infinite limits.

However, we still have:

\begin{lemma}\label{lemma:alpha-lower-star-tame}
If $\alpha\colon H\to G$ is injective with $(G:\im\alpha)<\infty$, then
the usual adjunction $\alpha^*\colon\cat{$\bm{E\mathcal M}$-$\bm G$-SSet}\rightleftarrows\cat{$\bm{E\mathcal M}$-$\bm H$-SSet}:\alpha_*$ restricts to a Quillen adjunction
\begin{equation*}
\alpha^*\colon\cat{$\bm{E\mathcal M}$-$\bm G$-SSet}^\tau_{\textup{positive $G$-global}}\rightleftarrows\cat{$\bm{E\mathcal M}$-$\bm H$-SSet}^\tau_{\textup{positive $H$-global}}:\alpha_*
\end{equation*}
in which both functors are fully homotopical.
\begin{proof}
We already know that $\alpha^*$ preserves tameness and is fully homotopical. To see that also $\alpha_*$ preserves tameness, we observe that as an $E\mathcal M$-simplicial set, $\alpha_*X$ is just a $(G:\im\alpha)$-fold product of copies of $X$, and that $\cat{$\bm{E\mathcal M}$-SSet}^\tau\subset\cat{$\bm{E\mathcal M}$-SSet}$ is closed under \emph{finite} limits by Corollary~\ref{cor:EM-tau-colim-lim}.

It only remains to show that $\alpha_*$ preserves fibrations as well as weak equivalences. But indeed, as $\alpha^*$ commutes with $\ev_\omega$, $\alpha_*$ commutes with $(\blank)_\bullet$ up to (canonical) isomorphism, so these follow from Lemma~\ref{lemma:alpha-lower-star-injective-script-I}.
\end{proof}
\end{lemma}

By abstract nonsense, (finite) products of fibrations in $\cat{$\bm{E\mathcal M}$-$\bm G$-SSet}^\tau$ are fibrations, and we have seen that also finite products of $G$-global weak equivalences are weak equivalences. If now $S$ is any finite set and $X$ is any $E\mathcal M$-$G$-simplicial set, then $X^{\times S}\mathrel{:=}\prod_{s\in S}X$ carries a natural $\Sigma_S$-action by permuting the factors, and this way $(\blank)^{\times S}$ lifts to $\cat{$\bm{E\mathcal M}$-$\bm G$-SSet}^\tau\to\cat{$\bm{E\mathcal M}$-$\bm{(G\times\Sigma_S)}$-SSet}^\tau$. We will later need the following (non-formal) strengthening of this observation:

\begin{cor}\label{cor:twisted-product}
The above lift of $(\blank)^{\times S}$ sends $G$-global weak equivalences or fibrations to $(G\times\Sigma_S)$-global weak equivalences or fibrations, respectively.
\begin{proof}
The claim is trivial if $S$ is empty. Otherwise, we pick $s_0\in S$, and we write $\Sigma_S^0\subset\Sigma_S$ for the subgroup of permutations fixing $s_0$, $p\colon G\times\Sigma_S^0\to G$ for the projection to the second factor, and $i\colon G\times\Sigma_S^0\hookrightarrow G\times\Sigma_S$ for the inclusion.

\begin{claim*}
The functor $(\blank)^{\times S}$ is isomorphic to $i_*\circ p^*$.
\begin{proof}
Fix for each $s\in S$ a permutation $\sigma_s\in\Sigma_S$ with $\sigma_s(s_0)=s$. Then we consider for $X\in \cat{$\bm{E\mathcal M}$-$\bm G$-SSet}^\tau$ the natural map $\Maps^{G\times\Sigma_S^0}(G\times\Sigma_S, X)\to X^{\times S}$ given on the $s$-th factor by evaluating at $(1,\sigma_s^{-1})$. We omit the easy verification that this is $(G\times\Sigma_S)$-equivariant and an isomorphism.
\end{proof}
\end{claim*}

Thus, the statement follows from Lemmas~\ref{lemma:alpha-lower-shriek-tau} and~\ref{lemma:alpha-lower-star-tame}.\index{functoriality in homomorphisms!for EM-G-SSettau@for $\cat{$\bm{E\mathcal M}$-$\bm G$-SSet}^\tau$|)}%
\end{proof}
\end{cor}

\begin{rk}
Again, there are similar functoriality properties for the remaining models. We leave the details to the interested reader.
\end{rk}

\section{Comparison to global homotopy theory}\label{sec:global-vs-g-global}
In this section we prove as promised that our theory generalizes Schwede's \emph{unstable global homotopy theory} with respect to finite groups. More precisely, we will give a chain of Quillen adjunctions between $\cat{$\bm I$-SSet}$ and Schwede's \emph{orthogonal spaces} that on associated quasi-categories exhibits the former as a (right Bousfield) localization with respect to an explicit class of `$\mathcal Fin$-global weak equivalences.'

\subsection{A reminder on orthogonal spaces}
Orthogonal spaces are based on a certain topological analogue $L$ of the categories $I$ and $\mathcal I$ considered above. Explicitly, the objects of $L$\nomenclature[aL3]{$L$}{topologically enriched category of finite dimensional real inner product spaces and linear isometric embeddings} are the finite dimensional real inner product spaces $V$, and as a set $L(V,W)$ is given by the linear isometric embeddings $V\to W$. For $V=W$, $L(V,W)$ carries the topology of the orthogonal group $O(V)$, and in general $L(V,W)$ is topologized as a Stiefel manifold; since we can completely black box the topology, we omit the details and refer the curious reader to \cite[discussion before Definition~1.1.1]{schwede-book} instead.

\begin{defi}\label{defi:orthogonal-space}\index{orthogonal space|textbf}
An \emph{orthogonal space} is a topologically enriched functor $L\to\cat{Top}$. We write $\cat{$\bm L$-Top}\mathrel{:=}\FUN(L,\cat{Top})$.
\end{defi}

Schwede \cite[Definition~1.1.1]{schwede-book} denotes the above category by `\textit{spc}' and he constructs a global model structure on it that we will recall now.

\subsubsection{The strict level model structure}
\index{global homotopy theory|(}
As before, we begin by constructing a suitable level model structure.

\begin{defi}
A map $f\colon X\to Y$ of orthogonal spaces is called a \emph{strict level weak equivalence}\index{strict level weak equivalence|textbf} or \emph{strict level fibration} if $f(V)$ is a weak equivalence or fibration, respectively, in the equivariant model structure on $\cat{$\bm{O(V)}$-Top}$ with respect to all closed subgroups for every $V\in L$, i.e.~$f(V)^H$ is a weak homotopy equivalence or Serre fibration, respectively, for every closed subgroup $H\subset O(V)$.
\end{defi}

Schwede \cite[Definition~1.1.8]{schwede-book} calls the above `strong level equivalences' and `strong level fibrations,' respectively.

If $H$ is a compact Lie group, then an \emph{orthogonal $H$-representation} is an object $V\in L$ together with a continuous homomorphism $\rho\colon H\to O(V)$; equivalently, we can view this as a pair of a finite dimensional real inner product space $V$ together with a continuous $H$-action by linear isometries. Restricting along $\rho$, $X(V)$ is then naturally an $H$-space for every orthogonal space $X$, and $f(V)\colon X(V)\to Y(V)$ is $H$-equivariant for any map $f\colon X\to Y$ of orthogonal spaces. The above condition can then be rephrased as saying that $f(V)$ should be a weak equivalence or fibration, respectively, in $\cat{$\bm H$-Top}$ for every compact Lie group $H$ and every orthogonal $H$-representation $V$, see~\cite[Lemmas~1.2.7 and~1.2.8]{schwede-book}.

\begin{prop}
There is a unique model structure on $\cat{$\bm L$-Top}$ with weak equivalences the strict level weak equivalences and fibrations the strict level fibrations. It is topological and cofibrantly generated with generating cofibrations
\begin{equation*}
\{L(V,\blank)/H\times\del D^n\hookrightarrow L(V,\blank)/H\times D^n : V\in L, H\subset O(V)\text{ closed},n\ge0\}.
\end{equation*}
\begin{proof}
This is \cite[Proposition~1.2.10]{schwede-book} and the discussion after it.
\end{proof}
\end{prop}

\subsubsection{Global weak equivalences}
The global weak equivalences of orthogonal spaces are slightly intricate to define due to some point-set topological issues. Intuitively speaking, however, they should again be created by `evaluating at $\mathbb R^\infty$,' analogously to the approach for $I$- and $\mathcal I$-simplicial sets:

\begin{constr}
We write $\mathcal L$ (the `universal compact Lie group')\nomenclature[aL4]{$\mathcal L$}{`universal compact Lie group,' topological monoid of linear isometric embeddings $\mathbb R^\infty\to\mathbb R^\infty$}\index{universal compact Lie group|textbf} for the topological monoid of linear isometric embeddings $\mathbb R^\infty\to\mathbb R^\infty$ under composition; here the scalar product on $\mathbb R^\infty=\mathbb R^{(\omega)}$ (the vector space of functions $\omega\to\mathbb R$ vanishing almost everywhere) is so that the canonical basis consisting of the characteristic functions of elements of $\omega$ is orthonormal. The topology on $\mathcal L$ is given as a subspace of the mapping space.

Let $X$ be an orthogonal space. Then \cite[Construction~3.2]{schwede-orbi} describes how $X$ yields a space $X(\mathbb R^\infty)$ with a continuous $\mathcal L$-action. Explicitly, we define
\begin{equation*}
X(\mathbb R^\infty)\mathrel{:=} \colim_{V\subset\mathbb R^\infty\text{ finite dimensional}} X(V),
\end{equation*}
where the structure maps of the colimit system are induced via $X$ from the inclusions. A continuous $\mathcal L$-action is given as follows: if $x$ is contained in the image of $X(V)\to X(\mathbb R^\infty)$, then $u.x$ is the image of $x$ under the composition
\begin{equation*}
X(V)\xrightarrow{X(u|_V\colon V\to u(V))} X(u(V))\to X(\mathbb R^\infty)
\end{equation*}
for any $u\in\mathcal L$, where the right hand map is the structure map of the colimit.

By functoriality of colimits this yields a functor $\ev_{\mathbb R^\infty}\colon\cat{$\bm L$-Top}\to\cat{$\bm{\mathcal L}$-Top}$.
\end{constr}

Unlike for simplicial sets, filtered (and even sequential) colimits of topological spaces do not preserve weak equivalences in general, which already suggests that $\ev_{\mathbb R^\infty}$ is not homotopically meaningful. \cite[Definition~1.1.2]{schwede-book} avoids this problem by defining the weak equivalences in terms of a `homotopy extension lifting property' instead. We will take a different approach following \cite{schwede-orbi} here:

Namely, non-equivariantly we could solve this issue by replacing the above colimit by a homotopy colimit, but this of course does not retain equivariant information: surely, whether $f\colon X\to Y$ is a global weak equivalence should not only depend on $f(\mathbb R^\infty)$ as a map of ordinary topological spaces, but it should take into account the $\mathcal L$-action and in particular suitable actions of all compact Lie groups. Concretely, see~\cite[Definitions~1.3, 1.4, and~1.6]{schwede-orbi}:

\begin{defi}
A compact subgroup $H\subset\mathcal L$ is called \emph{universal}\index{universal subgroup!for L@for $\mathcal L$|textbf} if it admits the structure of a Lie group and the tautological $H$-action on $\mathbb R^\infty$ makes the latter into a complete $H$-universe,\index{complete H-universe@complete $H$-universe} i.e.~any finite dimensional orthogonal $H$-representation embeds $H$-equivariantly linearly isometrically into $\mathbb R^\infty$.

A map $f\colon X\to Y$ of $\mathcal L$-spaces is a \emph{global weak equivalence}\index{global weak equivalence!in L1-Top@in $\cat{$\bm{\mathcal L}$-Top}$|textbf} if $f^H\colon X^H\to Y^H$ is a weak homotopy equivalence for every universal $H\subset\mathcal L$.
\end{defi}

\begin{rk}
Analogously to the situation for $\mathcal M$, any compact Lie group is isomorphic to a universal subgroup of $\mathcal L$, and any two such embeddings are conjugate \cite[Proposition~1.5]{schwede-orbi}.
\end{rk}

One way to calculate sequential homotopy colimits is by replacing the diagram in question by a sequence of closed embeddings. The same strategy works \emph{mutatis mutandis} in our situation:

\begin{defi}\label{defi:closed-orth-space}
An orthogonal space $X$ is called \emph{closed}\index{orthogonal space!closed|textbf} if for every map $\phi\colon V\to W$ in $L$ the induced map $X(\phi)\colon X(V)\to X(W)$ is a closed embedding.
\end{defi}

\begin{ex}\label{ex:cofibrant-closed}
Any orthogonal space that is cofibrant in the strict level model structure is closed, see~\cite[Proposition~1.2.11-(iii)]{schwede-book}.
\end{ex}

\begin{defi}
Let $f\colon X\to Y$ be a map of orthogonal spaces and let
\begin{equation*}
\begin{tikzcd}
\hat X\arrow[d,"\hat f"']\arrow[r, "\sim"] & X\arrow[d, "f"]\\
\hat Y\arrow[r, "\sim"'] & Y
\end{tikzcd}
\end{equation*}
be a commutative diagram such that the horizontal maps are strict level weak equivalences and $\hat X,\hat Y$ are closed. Then $f$ is called a \emph{global weak equivalence}\index{global weak equivalence!in L2-Top@in $\cat{$\bm L$-Top}$|textbf} if $\hat f(\mathbb R^\infty)$ is a global weak equivalence of $\mathcal L$-spaces.
\end{defi}

Note that we can always find such a square by just taking functorial cofibrant replacements in the strict level model structure. Moreover, \cite[Proposition~3.5]{schwede-orbi} together with \cite[Proposition~1.1.9-(i)]{schwede-book} shows that the above is independent of the choice of replacement and equivalent to Schwede's original definition.

We can now finally introduce the \emph{global model structure} on $\cat{$\bm L$-Top}$ \cite[Theorem~1.2.21]{schwede-book}:

\begin{defi}
An orthogonal space $X$ is called \emph{static} if $X(\phi)\colon X(V)\to X(W)$ is an $H$-equivariant weak equivalence for every compact Lie group $H$ and every $H$-equivariant linear isometric embedding $\phi\colon V\to W$ of faithful finite dimensional orthogonal $H$-representations.
\end{defi}

\begin{thm}\label{thm:L-global-model-structure}\index{global model structure!on L-Top@on $\cat{$\bm L$-Top}$|textbf}\index{orthogonal space!global model structure|seeonly{global model structure, on $\cat{$\bm L$-Top}$}}
There is a unique model structure on $\cat{$\bm L$-Top}$ with the same cofibrations as the strict level model structure and with the global weak equivalences as weak equivalences. This model structure is topological, proper, and cofibrantly generated with generating cofibrations
\begin{equation*}
\{L(V,\blank)/H\times\del D^n\hookrightarrow L(V,\blank)/H\times D^n : V\in L, H\subset O(V)\text{ closed}, n\ge 0\}.
\end{equation*}
Moreover, an orthogonal space is fibrant if and only if it is static.\qed
\end{thm}

\begin{rk}
Again, there is also a \emph{positive global model structure}\index{positive global model structure!on L-Top@on $\cat{$\bm L$-Top}$}\index{orthogonal space!positive global model structure|seeonly{positive global model structure, on $\cat{$\bm L$-Top}$}} where one restricts the generating cofibrations by demanding in addition that $V\not=0$, see~\cite[Proposition~1.2.23]{schwede-book}.
\end{rk}

By design, our models of global homotopy theory considered in the previous sections only see equivariant information with respect to finite groups, so we should not hope for them to be equivalent to the above model category. Instead, we consider the following coarser notion of weak equivalence:

\begin{defi}\label{defi:Fin-global-we}
A map of $\mathcal L$-spaces is called a \emph{$\mathcal Fin$-global weak equivalence}\index{Fin-global weak equivalence@$\mathcal Fin$-global weak equivalence|textbf}\index{global weak equivalence!in L2-Top@in $\cat{$\bm L$-Top}$|seealso{$\mathcal Fin$-global weak equivalence}} if it restricts to an $H$-equivariant weak equivalence for each \emph{finite} universal $H\subset\mathcal L$.

A map $f$ of orthogonal spaces is called a \emph{$\mathcal Fin$-global weak equivalence} if $\hat f(\mathbb R^\infty)$ is a $\mathcal Fin$-global weak equivalence in $\cat{$\bm{\mathcal L}$-Top}$ for some (hence any) replacement of $f$ up to strict level weak equivalence by a map $\hat f$ between closed orthogonal spaces.
\end{defi}

\begin{rk}\label{rk:Fin-global-model-structure}
As remarked without proof in \cite[Remark~3.11]{schwede-orbi}, there is a version for the global model structure on $\cat{$\bm L$-Top}$ which only sees representations of groups belonging to a given \emph{global family} $\mathcal F$, i.e.~a collection of compact Lie groups closed under isomorphisms and subquotients. For $\mathcal F=\mathcal Fin$ the family of finite groups this precisely recovers the above $\mathcal Fin$-global weak equivalences.
\end{rk}

\subsection{$\bm I$-spaces} The intermediate step in our comparison will be a global model structure on $\cat{$\bm I$-Top}$. For this the following terminology will be useful, see e.g.~\cite[Definition~A.28 and discussion after it]{schwede-book}:

\begin{defi}
Let $\mathscr C$ be a category enriched and tensored over $\cat{Top}$. Then a map $f\colon A\to B$ in $\mathscr C$ is an \emph{h-cofibration}\index{h-cofibration|textbf} if the natural map $\big(A\times[0,1]\big)\amalg_A B\to B\times[0,1]$ from the mapping cyclinder of $f$ admits a retraction.
\end{defi}

\begin{ex}
For $\mathscr C=\cat{Top}$ with the usual enrichment, the h-cofibrations are the classical (Hurewicz) cofibrations.
\end{ex}

\begin{lemma}
Let $G$ be a finite group and let
\begin{equation}\label{diag:pushout-along-h-cof}
\begin{tikzcd}
A\arrow[d]\arrow[r, "i"] & B\arrow[d]\\
C\arrow[r] & D
\end{tikzcd}
\end{equation}
be a pushout in $\cat{$\bm G$-Top}$ such that $i$ is an h-cofibration. Then $\Sing\colon\cat{$\bm G$-Top}\to\cat{$\bm G$-SSet}$ sends $(\ref{diag:pushout-along-h-cof})$ to a homotopy pushout.
\begin{proof}
Let $H\subset G$ be any subgroup. The functor $(\blank)^H$ preserves pushouts along closed embeddings by \cite[Proposition~B.1-(i)]{schwede-book} and it clearly preserves tensors; it easily follows that it preserves h-cofibrations, also cf.~\cite[Corollary~A.30-(ii)]{schwede-book}. We conclude that $i^H$ is a Hurewicz cofibration and that the square on the left in
\begin{equation*}
\begin{tikzcd}
A^H\arrow[d]\arrow[r, "i^H"] & B^H\arrow[d]\\
C^H\arrow[r] & D^H
\end{tikzcd}
\qquad\qquad
\begin{tikzcd}[column sep=large]
(\Sing A)^H\arrow[d]\arrow[r,"(\Sing i)^H"] & (\Sing B)^H\arrow[d]\\
(\Sing C)^H\arrow[r]&(\Sing D)^H
\end{tikzcd}
\end{equation*}
is a pushout in $\cat{Top}$. It is well-known that $\Sing$ sends pushouts along Hurewicz cofibrations to homotopy pushouts (which also follows from Lemma~\ref{lemma:U-pushout-preserve-reflect} applied to $|\blank|$, using that the counit $\epsilon\colon |\Sing X|\to X$ is a weak homotopy equivalence for any $X\in\cat{Top}$), and as it moreover commutes with fixed points, we conclude that the square on the right is a homotopy pushout in $\cat{SSet}$.

The claim now follows as homotopy pushouts in $\cat{$\bm G$-SSet}$ can be checked on fixed points by Proposition~\ref{prop:dk-equivariant-addendum}.
\end{proof}
\end{lemma}

\begin{cor}\label{cor:h-cofibration-I-po}
Let $i\colon A\to B$ be any h-cofibration of $I$-spaces. Then $\Sing\colon\cat{$\bm I$-Top}\to\cat{$\bm I$-SSet}$ sends pushouts along $i$ to homotopy pushouts in the global model structure.
\begin{proof}
The evaluation functors $\ev_A\colon\cat{$\bm I$-Top}\to\cat{$\bm{\Sigma_A}$-Top}$ for $A\in I$ are cocontinuous and preserve tensors, so h-cofibrations of $I$-spaces are in particular levelwise h-cofibrations by \cite[Corollary~A.30-(ii)]{schwede-book}.

Thus, if we are given a pushout in $\cat{$\bm I$-Top}$ along $i$, then applying the previous lemma levelwise shows that its image under $\Sing$ is a levelwise homotopy pushout in the sense that evaluating at any $A\in I$ yields a homotopy pushout in $\cat{$\bm{\Sigma_A}$-SSet}$. As one easily concludes from the existence of the injective global model structure on $\cat{$\bm I$-SSet}$, such a levelwise homotopy pushout is in particular a homotopy pushout, finishing the proof.
\end{proof}
\end{cor}

\begin{defi}
A map $f\colon X\to Y$ in $\cat{$\bm I$-Top}$ is called a \emph{global weak equivalence}\index{global weak equivalence|in I-Top@in $\cat{$\bm I$-Top}$|textbf} or \emph{global fibration} if $\Sing f$ is a weak equivalence or fibration, respectively, in the global model structure on $\cat{$\bm I$-SSet}$.
\end{defi}

\begin{defi}
An $I$-space $X$ is called \emph{static} if $\Sing X$ is a static $I$-simplicial set, i.e.~the map $X(i)^H\colon X(A)^H\to X(B)^H$ is a weak homotopy equivalence of topological spaces for every finite group $H$ and every injection $i\colon A\to B$ of finite faithful $H$-sets.
\end{defi}

\begin{prop}\label{prop:I-Top-global}
The global weak equivalences and fibrations are part of a unique model structure on $\cat{$\bm I$-Top}$, which we call the \emph{global model structure}.\index{global model structure!on I-Top@on $\cat{$\bm I$-Top}$|textbf} Its fibrant objects are precisely the static $I$-spaces.

The global model structure is topological, left proper, and cofibrantly generated with generating cofibrations
\begin{equation*}
I(A,\blank)/H\times\del D^n\hookrightarrow
I(A,\blank)/H\times D^n,
\end{equation*}
where $A$ and $H$ are as above and $n\ge 0$. Moreover, pushouts along h-cofibrations are homotopy pushouts.

Finally, the adjunction
\begin{equation}\label{eq:real-sing-I}
|\blank|\colon\cat{$\bm I$-SSet}_{\textup{global}}
\rightleftarrows\cat{$\bm I$-Top}_{\textup{global}} :\!\Sing
\end{equation}
is a Quillen equivalence in which both adjoints are homotopical.
\begin{proof}
The adjunction $|\blank|\colon\cat{SSet}\rightleftarrows\cat{Top} :\!\Sing$ is a Quillen equivalence in which both adjoints are homotopical. In particular, the unit $\eta\colon X\to\Sing|X|$ is a weak homotopy equivalence for any space $X$. As both adjoints commute with finite limits, it follows further that for any group $K$ acting on $X$, the unit $\eta_X$ induces weak homotopy equivalences on $K'$-fixed points for all finite $K'\subset K$, so that the unit of $(\ref{eq:real-sing-I})$ is a strict level weak equivalence for all $X\in\cat{$\bm I$-SSet}$.

Let us now construct the model structure, for which we will verify the conditions of Crans' Transfer Criterion (Proposition~\ref{prop:transfer-criterion}): we let $I_\Delta$ be the usual set of generating cofibrations, and we pick any set $J_\Delta$ of generating acyclic cofibrations of $\cat{$\bm I$-SSet}_{\textup{global}}$. As each cofibration is in particular a levelwise injection, we conclude that $|I_\Delta|$ and $|J_\Delta|$ consist of closed embeddings. As $\Sing$ preserves transfinite compositions along closed embeddings (see e.g.~\cite[Proposition~2.4.2]{hovey}), the local presentability of $\cat{$\bm I$-SSet}$ implies that $|I_\Delta|$ and $|J_\Delta|$ permit the small object argument.

It remains to show that any relative $|J_\Delta|$-cell complex is sent to a global weak equivalence under $\Sing$. Again using that $\Sing$ preserves transfinite compositions along closed embeddings, it suffices to show that pushouts of maps in $|J_\Delta|$ are global weak equivalences.

It is clear that the maps in $|I_\Delta|$ are h-cofibrations, hence so is the geometric realization of any cofibration by \cite[Corollary~A.30-(i)]{schwede-book} and in particular any map in $|J_\Delta|$. As the unit is a levelwise global (in fact, even strict level) weak equivalence, all maps in $|J_\Delta|$ are sent to global weak equivalences under $\Sing$, and hence so is any pushout of a map in $|J_\Delta|$ by Corollary~\ref{cor:h-cofibration-I-po}. This completes the proof of the existence of the model structure. Moreover, again using that the unit is a levelwise global weak equivalence and in addition that the right adjoint creates weak equivalences by definition, we immediately conclude that $(\ref{eq:real-sing-I})$ is a Quillen equivalence and that also $|\blank|$ is homotopical.

By Corollary~\ref{cor:h-cofibration-I-po}, $\Sing$ sends pushouts along h-cofibrations (hence in particular along cofibrations of the above model structure) to homotopy pushouts. Thus, Lemma~\ref{lemma:U-pushout-preserve-reflect} implies that $\cat{$\bm I$-Top}$ is left proper and that $\Sing$ creates homotopy pushouts. In particular, pushouts along h-cofibrations are homotopy pushouts.

As a functor category, $\cat{$\bm I$-Top}$ is enriched, tensored, and cotensored over $\cat{Top}$ in the obvious way. Restricting along the adjunction $|\blank|\colon\cat{SSet}\rightleftarrows\cat{Top}:\!\Sing$ therefore makes it into a category enriched, tensored, and cotensored over $\cat{SSet}$. With respect to this, $(\ref{eq:real-sing-I})$ is naturally a simplicial adjunction, so $\cat{$\bm I$-Top}$ is a simplicial model category by Lemma~\ref{lemma:transferred-properties}-$(\ref{item:tpr-simplicial})$. To see that it is topological, we then simply observe that it suffices to verify the Pushout Product Axiom for generating cofibrations and generating acyclic cofibrations, and that the usual ones for $\cat{Top}$ agree up to conjugation by isomorphisms with the images under geometric realization of the standard generating (acyclic) cofibrations of $\cat{SSet}$.
\end{proof}
\end{prop}

\begin{rk}\label{rk:I-Top-closed}
In analogy to Definition~\ref{defi:closed-orth-space}, let us call an $I$-space $X$ \emph{closed} if all structure maps $X(A)\to X(B)$ are closed embeddings. As $\Sing$ preserves sequential colimits along closed embeddings, it easily follows that $\Sing$ commutes with $\ev_\omega$ on the subcategory of \emph{closed} $I$-spaces. In particular, if $f\colon X\to Y$ is a map of closed $I$-spaces that is a \emph{weak equivalence at infinity}\index{weak equivalence at infinity} in the sense that $f(\omega)^H$ is a weak homotopy equivalence for each universal $H\subset\mathcal M$, then $f$ is already a $G$-global weak equivalence.
\end{rk}

\subsection{Proof of the comparison}
\index{orthogonal space!vs I-spaces@vs.~$I$-spaces|(}
In order to compare $I$-spaces to orthogonal spaces, we will relate their indexing categories:

\begin{constr}\label{constr:I-to-L}
We define a functor $\mathbb R^\bullet: I\to L$ as follows: a finite set $A$ is sent to the real vector space $\mathbb R^A$ of maps $A\to\mathbb R$, where the inner product on $\mathbb R^A$ is the unique one such that the characteristic functions of elements of $A$ form an orthonormal basis.

If $f\colon A\to B$ is an injection of finite sets, then $\mathbb R^f\colon\mathbb R^A\to\mathbb R^B$ is defined to be the unique $\mathbb R$-linear map sending the characteristic function of $a\in A$ to the characteristic function of $f(a)\in B$. It is clear that $\mathbb R^f$ is isometric, and that this makes $\mathbb R^\bullet$ into a well-defined functor.
\end{constr}

Restricting along $\mathbb R^\bullet$ yields a forgetful functor $\cat{$\bm L$-Top}\to\cat{$\bm I$-Top}$. By topologically enriched left Kan extension, this admits a topological left adjoint $L\times_I\blank$ satisfying $L\times_I I(A,\blank) = L(\mathbb R^A,\blank)$ for any finite set $A$; the unit is then given on such corepresentables by $\mathbb R^\bullet\colon I(A,\blank)\to L(\mathbb R^A,\mathbb R^\bullet)= \forget L(\mathbb R^A,\blank)$.

\begin{rk}\label{rk:set-universe-vs-universe}
Let $H$ be a finite group and let $U$ be a countable set containing infinitely many free $H$-orbits (for example if $U$ is a complete $H$-set universe). Then the $\mathbb R$-linearization $\mathbb R^{(U)}$ contains infinitely many copies of the regular representation $\mathbb R^H$, so it is a complete $H$-universe in the usual sense. In particular, if $H\subset\mathcal M$ is universal, then the induced $H$-action on $\mathbb R^\infty=\mathbb R^{(\omega)}$ makes the latter into a complete $H$-universe.
\end{rk}

\begin{prop}\label{prop:L-vs-I-unit-corep}
The map $\mathbb R^\bullet\colon I(A,\blank)/H\times X\to {L}(\mathbb R^A,\mathbb R^\bullet)/H\times X$ is a global weak equivalence in $\cat{$\bm I$-Top}$ for any finite group $H$, any finite faithful $H$-set $A$, and any CW-complex $X$.
\begin{proof}
Taking products with $X$ obviously preserves strict level weak equivalences, and it preserves global acyclic cofibrations as $\cat{$\bm I$-Top}$ is topological. Thus, $\blank\times X$ is fully homotopical, and it suffices that $\mathbb R^\bullet\colon I(A,\blank)/H\to L(\mathbb R^A,\mathbb R^\bullet)/H$ is a global weak equivalence. 

For this we consider the commutative diagram
\begin{equation}\label{diag:semistable-replacements}
\begin{tikzcd}[row sep=2.2em]
 & {|\mathcal I(A,\blank)|/H}\\
{I(A,\blank)/H}\arrow[ur, bend left=12.5pt]\arrow[r]\arrow[dr, bend right=12.5pt] & {\big(|\mathcal I(A,\blank)|\times {L}(\mathbb R^A,\mathbb R^\bullet)\big)/H}\arrow[u, "\pr"']\arrow[d, "\pr"]\\
 & {{L}(\mathbb R^A,\mathbb R^\bullet)/H}
\end{tikzcd}
\end{equation}
where the maps from left to right are induced by the inclusion $I(A,\blank)\hookrightarrow|\mathcal I(A,\blank)|$ and by $\mathbb R^\bullet\colon I(A,\blank)\to {L}(\mathbb R^A,\mathbb R^\bullet)$. We begin by showing that the vertical arrows on the right are global weak equivalences.

It is clear that the $I$-space $|\mathcal I(A,\blank)|$ is closed, and so is ${L}(\mathbb R^A,\mathbb R^\bullet)$ by Example~\ref{ex:cofibrant-closed}. From this we can conclude by \cite[Proposition~B.13-(iii)]{schwede-book} that all the $I$-spaces on the right of $(\ref{diag:semistable-replacements})$ are closed, so that it suffices that the vertical maps are weak equivalences at infinity (see Remark~\ref{rk:I-Top-closed}).

Let $K\subset\mathcal M$ be a universal subgroup. Obviously, $|\mathcal I(A,\blank)|(\omega)\cong|E\Inj(A,\omega)|$ is $\mathcal A\ell\ell$-cofibrant in $\cat{$\bm{(K\times H)}$-Top}$, and so is ${L}(\mathbb R^A,\mathbb R^\bullet)(\omega)$ by \cite[Proposition~1.1.19-(ii)]{schwede-book} together with \cite[Proposition~A.5-(ii)]{schwede-orbi}. We claim that both are classifying spaces for $\mathcal G_{K,H}$ in the sense that their $T$-fixed points for $T\subset K\times H$ are contractible if $T\in\mathcal G_{K,H}$, and empty otherwise. Indeed, for ${L}(\mathbb R^A,\mathbb R^\bullet)(\omega)$ this is an instance of \cite[Proposition~1.1.26-(i)]{schwede-book}. On the other hand, if $K'\subset K$, $\phi\colon K'\to H$, then $|\mathcal I(A,\blank)|(\omega)^\phi\cong |E(\Inj(A,\omega)^\phi|$ is weakly contractible by Example~\ref{ex:G-globally-contractible}. Finally, $H$ acts freely on $\Inj(A,\omega)$, so $\mathcal I(A,\omega)^T$ has no vertices when $T\subset K\times H$ is not contained in $\mathcal G_{K,H}$, whence $|\mathcal I(A,\blank)|(\omega)^T$ has to be empty.

Thus, the projections
\begin{equation*}
|\mathcal I(A,\blank)|(\omega)\gets\big(|\mathcal I(A,\blank)|\times {L}(\mathbb R^A,\mathbb R^\bullet)\big)(\omega)\to{L}(\mathbb R^A,\mathbb R^\bullet)(\omega)
\end{equation*}
are weak equivalences of cofibrant objects in $\cat{$\bm{(K\times H)}$-Top}$ with respect to the $\mathcal A\ell\ell$-model structure. Applying Ken Brown's Lemma to the Quillen adjunction
\begin{equation*}
p_!\colon\cat{$\bm{(K\times H)}$-Top}_{\textup{$\mathcal A\ell\ell$}}\rightleftarrows\cat{$\bm{K}$-Top}_{\mathcal A\ell\ell} :\!p^*
\end{equation*}
where $p\colon K\times H\to K$ denotes the projection, therefore shows that after evaluation at $\omega$ the vertical maps in $(\ref{diag:semistable-replacements})$ become $K$-weak equivalences. The claim follows by letting $K$ vary.

To finish the proof, we observe now that $I(A,\blank)/H\to|\mathcal I(A,\blank)|/H$ agrees up to conjugation by isomorphisms with the image of
\begin{equation}\label{eq:semistabilization-simplicial}
I(A,\blank)/H\hookrightarrow\mathcal I(A,\blank)/H
\end{equation}
under geometric realization. The proposition follows as $(\ref{eq:semistabilization-simplicial})$ is a global weak equivalence by Theorem~\ref{thm:strict-global-I-model-structure} and since geometric realization is fully homotopical by Proposition~\ref{prop:I-Top-global}.
\end{proof}
\end{prop}

We define a monoid homomorphism $i\colon\mathcal M\to\mathcal L$ by sending $f\colon\omega\to\omega$ to $\mathbb R^f\colon\mathbb R^\infty\to\mathbb R^\infty$, i.e.~the unique linear isometry sending the $i$-th standard basis vector to the $f(i)$-th one.

\begin{lemma}
Let $X$ be an $\mathcal L$-space. Then $\Sing(i^*X)\in\cat{$\bm{\mathcal M}$-SSet}$ is semistable.
\begin{proof}
Let $H\subset\mathcal M$ be universal, and let $u\in\mathcal M$ centralize $H$. It suffices to show that $i(u).\blank\colon X\to X$ is an $i(H)$-equivariant weak equivalence.

Indeed, Remark~\ref{rk:set-universe-vs-universe} implies that $i(H)$ is a universal subgroup of $\mathcal L$, so $\mathcal L^{i(H)}$ is contractible by \cite[Proposition~1.1.26-(i)]{schwede-book} together with \cite[Proposition~A.10]{schwede-orbi}. In particular, there exists a path $\gamma\colon[0,1]\to \mathcal L^{i(H)}$ connecting $i(u)$ to the identity. The composition
\begin{equation*}
[0,1]\times X\xrightarrow{\gamma\times X} \mathcal L\times X\xrightarrow{\textup{action}} X
\end{equation*}
is then an $i(H)$-equivariant homotopy from $i(u).\blank$ to the identity, so $i(u).\blank$ is in particular an $i(H)$-equivariant (weak) homotopy equivalence.
\end{proof}
\end{lemma}

\begin{prop}\label{prop:Fin-global-we}\index{Fin-global weak equivalence@$\mathcal Fin$-global weak equivalence}
Let $f\colon X\to Y$ be a map of orthogonal spaces. Then the following are equivalent:
\begin{enumerate}
\item $f$ is a $\mathcal Fin$-global weak equivalence (Definition~\ref{defi:Fin-global-we}).
\item $\forget f$ is a global weak equivalence of $I$-spaces.
\end{enumerate}
Moreover, if $X$ and $Y$ are closed, then also the following statements are equivalent to the above:
\begin{enumerate}
\item[(3)] $\Sing\big(i^*f(\mathbb R^\infty)\big)$ is a universal weak equivalence of $\mathcal M$-simplicial sets.
\item[(4)] $\Sing\big(i^*f(\mathbb R^\infty)\big)$ is a global weak equivalence of $\mathcal M$-simplicial sets.
\end{enumerate}
\begin{proof}
Let us first assume that $X$ and $Y$ are closed; we will show that all of the above statements are equivalent.

For the equivalence $(1)\Leftrightarrow(3)$ we observe once more that $i\colon\mathcal M\to\mathcal L$ sends universal subgroups to universal subgroups. As any two abstractly isomorphic universal subgroups of $\mathcal L$ are conjugate \cite[Proposition~1.5]{schwede-orbi}, the claim now follows from the definitions.

The equivalence $(3)\Leftrightarrow(4)$ is immediate from the previous lemma. Finally, for $(2)\Leftrightarrow(4)$ it suffices by Remark~\ref{rk:I-Top-closed} to show that the diagram
\begin{equation*}
\begin{tikzcd}
\cat{$\bm L$-Top}\arrow[r, "\ev_{\mathbb R^\infty}"]\arrow[d, "\forget"'] & \cat{$\bm{\mathcal L}$-Top}\arrow[d, "i^*"]\\
\cat{$\bm I$-Top}\arrow[r, "\ev_\omega"'] & \cat{$\bm{\mathcal M}$-Top}
\end{tikzcd}
\end{equation*}
commutes up to natural isomorphism. This follows easily from the definitions once we observe that the map of posets $\{A\subset\omega\text{ finite}\}\to\{V\subset\mathbb R^\infty\text{ finite dimensional}\}$ sending $A$ to the image of the canonical map $\mathbb R^A\to\mathbb R^\infty$ is cofinal: if $v\in\mathbb R^\infty$, then there is only a finite set $S(v)$ of $i\in\omega$ such that $\pr_i(v)\not=0$ for the projection $\pr_i\colon\mathbb R^\infty\to\mathbb R$ to the $i$-th summand. Thus, if $V\subset\mathbb R^\infty$ is any finite dimensional subspace, and $v_1,\dots,v_n$ is a basis, then $S\mathrel{:=}S(v_1)\cup\cdots\cup S(v_n)$ is finite, and obviously $V$ is contained in the image of $\mathbb R^S\to\mathbb R^\infty$.

Now assume $X$ and $Y$ are not necessarily closed. If $j\colon A\to B$ is a strict level weak equivalence of orthogonal spaces, then $j$ is in particular a $\mathcal Fin$-global weak equivalence; moreover $\Sing(\forget(j))$ is obviously a strict level weak equivalence of $I$-simplicial sets, hence in particular a global weak equivalence. Thus, $(1)\Leftrightarrow(2)$ follows from the above special case together with $2$-out-of-$3$.
\end{proof}
\end{prop}

\begin{thm}\label{thm:L-vs-I}\index{orthogonal space!vs I-spaces@vs.~$I$-spaces|textbf}
The topologically enriched adjunction
\begin{equation}\label{eq:L-vs-I}
L\times_I\blank\colon \cat{$\bm I$-Top}\rightleftarrows\cat{$\bm L$-Top}:\forget
\end{equation}
is a Quillen adjunction. The induced adjunction of associated quasi-categories is a right Bousfield localization with respect to the $\mathcal Fin$-global weak equivalences.
\end{thm}

Non-equivariantly, this comparison was proven by Lind \cite[Theorem~6.2]{lind}.

\begin{proof}
Let us first show that $(\ref{eq:L-vs-I})$ is a Quillen adjunction. As it is a topologically enriched adjunction of topological model categories, it in particular becomes a simplicial adjunction of simplicial model categories when we restrict along the usual adjunction $\cat{SSet}\rightleftarrows\cat{Top}$. It therefore suffices (Proposition~\ref{prop:cofibrations-fibrant-qa}) to observe that $L\times_I\blank$ obviously preserves generating cofibrations, and that the forgetful functor preserves fibrant objects by the characterizations given in Theorem~\ref{thm:L-global-model-structure} and Proposition~\ref{prop:I-Top-global}, respectively.

By the previous proposition, the forgetful functor is homotopical and it precisely inverts the $\mathcal Fin$-global weak equivalences. It therefore only remains to show that the unit $X\to\forget(L\times_I X)$ is a global weak equivalence for any cofibrant $X\in\cat{$\bm I$-Top}$.

This will again be a cell induction argument (although we cannot literally apply Corollary~\ref{cor:saturated-trafo}): if $X$ is the source or target of one of the standard generating cofibrations of $\cat{$\bm I$-Top}$, then the claim is an instance of Proposition~\ref{prop:L-vs-I-unit-corep}. On the other hand, any pushout
\begin{equation*}
\begin{tikzcd}
I(A,\blank)/H\times\del D^{n-1}\arrow[r, hook]\arrow[d] & I(A,\blank)/H\times D^n\arrow[d]\\
X\arrow[r] & Y
\end{tikzcd}
\end{equation*}
along a generating cofibration is a homotopy pushout by left properness, and its image under $\forget\circ (L\times_I\blank)$ is a pushout along an h-cofibration by \cite[Corollary~A.30-(ii)]{schwede-book}, hence again a homotopy pushout by Proposition~\ref{prop:I-Top-global} above. We conclude that $\eta_Y$ is a global weak equivalence if $\eta_X$ is.

Using that transfinite compositions \emph{of closed embeddings} in $\cat{$\bm I$-Top}$ are homotopical (as they are preserved by $\Sing$), we see that $\eta_Z$ is a global weak equivalence for any cell complex $Z$ in the generating cofibrations. The claim follows as any cofibrant object of $\cat{$\bm I$-Top}$ is a retract of such a cell complex.
\end{proof}

\begin{cor}
The functor $\Sing\circ\forget:\cat{$\bm L$-Top}\to\cat{$\bm I$-SSet}$ preserves global fibrations and global weak equivalences, and it induces a quasi-localization $\cat{$\bm L$-Top}^\infty\to \cat{$\bm I$-SSet}^\infty$ at the $\mathcal Fin$-global weak equivalences.\qed
\end{cor}

\begin{rk}
It is not hard to show that $(\ref{eq:L-vs-I})$ is a Quillen adjunction with respect to the $\mathcal Fin$-global model structure on $\cat{$\bm L$-Top}$ mentioned in Remark~\ref{rk:Fin-global-model-structure}, hence a Quillen equivalence by the above theorem.

On the other hand, one can easily adapt the above proof to transfer the global model structure from $\cat{$\bm I$-Top}$ to $\cat{$\bm L$-Top}$. This way one obtains a model structure with the $\mathcal Fin$-global weak equivalences as weak equivalences, but slightly fewer cofibrations than the model structure sketched by Schwede.
\end{rk}
\index{global homotopy theory|)}
\index{orthogonal space!vs I-spaces@vs.~$I$-spaces|)}

\chapter{Coherent commutativity}\label{chapter:coherent}
Several equivalent notions of `commutative monoids up to coherent homotopies' have been studied classically, in particular in relation to infinite loop spaces and algebraic $K$-theory \cite{may-loop-spaces,segal-gamma,may-permutative}. In this chapter, we introduce a variety of models of `$G$-globally coherently commutative monoids,'\index{coherently commutative monoid!G-global@$G$-global|seeonly{$G$-globally coherently commutative monoid}} either based on the notion of \emph{ultra-commutativity} studied by Schwede \cite{schwede-book} in the global context, or on \emph{$\Gamma$-spaces}, as generalized equivariantly by Shimakawa \cite{shimakawa}. As the main result of this chapter (Theorem~\ref{thm:gamma-vs-uc}), we prove that all these approaches are equivalent.

\section{Ultra-commutativity}
This section is concerned with various \emph{box products} on our different models of $G$-global homotopy theory, and in particular we will revisit the box products on $\cat{$\bm{\mathcal M}$-SSet}^\tau$ and $\cat{$\bm I$-SSet}$ of \cite{I-vs-M-1-cat,sagave-schlichtkrull} from a $G$-global perspective. Among other things, we will show that all these box products are in fact fully homotopical, in particular yielding $G$-global versions of \cite[Theorems~1.2 and~4.8]{I-vs-M-1-cat}.

Our main goal is then to lift the equivalences of models of unstable $G$-global homotopy theory to equivalences between the corresponding homotopy theories of commutative monoids. In particular, this will give a $G$-global refinement of the equivalence between the non-equivariant homotopy theories of the categories of commutative monoids $\CMon(\cat{$\bm{\mathcal M}$-SSet}^\tau)$ and $\CMon(\cat{$\bm I$-SSet})$\nomenclature[aCMon]{$\CMon$}{category of commutative monoids} due to Schwede and Sagave, see \cite[proof of Theorem~5.13]{I-vs-M-1-cat}.

\subsection{A reminder on box products}
We briefly recall the box products of tame $\mathcal M$- and of $I$-simplicial sets, which will serve as blueprints for our definition of the box products of tame $E\mathcal M$- and $\mathcal I$-simplicial sets.

\begin{constr}
The coproduct of (finite) sets enhances $I$ to a symmetric monoidal category. More precisely, we have a functor $I\times I\to I$ given on objects by $(A,B)\mapsto A\amalg B$, and on morphisms by sending a pair $f\colon A\to A',g\colon B\to B'$ to $f\amalg g\colon A\amalg B\to A'\amalg B'$. The usual unitality, associativity, and symmetry isomorphisms of the cocartesian symmetric monoidal structure on $\cat{Set}$ then induce natural isomorphisms for the corresponding functors on $I$, and these satisfy the usual coherence conditions.

This symmetric monoidal structure then induces a \emph{Day convolution product} \cite{day-convolution}\index{Day convolution|textbf}\index{Day convolution|seealso{box product}}\index{box product!on I-SSet@on $\cat{$\bm I$-SSet}$|textbf}\nomenclature[zx]{$\boxtimes$}{box product} on $\cat{$\bm{I}$-SSet}$ as follows: by the universal property of enriched presheaves, there is a simplicially enriched functor $\blank\boxtimes\blank\colon\cat{$\bm{I}$-SSet}\times\cat{$\bm{I}$-SSet}\to\cat{$\bm{I}$-SSet}$ that preserves tensors and small colimits in each variable separately and such that $I(A,\blank)\boxtimes I(B,\blank)=I(A\amalg B,\blank)$ with the evident (contravariant) functoriality in $A$ and $B$. The functor $\boxtimes$ is unique up to a unique simplicial isomorphism that is the identity on pairs of corepresentables. We fix such a choice and call it `the' \emph{box product} on $\cat{$\bm{I}$-SSet}$. This is the symmetric monoidal product of a preferred simplicial symmetric monoidal structure on $\cat{$\bm{I}$-SSet}$, where all the structure isomorphisms are induced from the corresponding structure isomorphisms of $\amalg$. More precisely, the symmetry isomorphism $\tau$ is the unique simplicially enriched natural isomorphism such that $\tau\colon I(A,\blank)\boxtimes I(B,\blank)\to I(B,\blank)\boxtimes I(A,\blank)$ agrees for all finite $A,B$ with restriction along the inverse of the flip $A\amalg B\to B\amalg A$, and similarly for associativity and unitality.
\end{constr}

\begin{rk}
Another common perspective on the Day convolution of two functors $X,Y$---and in particular the one taken in our references \cite{I-vs-M-1-cat, schwede-book}---is as `the' object corepresenting \emph{bimorphisms},\index{bimorphism} i.e.~families of maps $X(A)\times Y(B)\to Z(A\amalg B)$ for $A,B\in I$ that are natural in both variables.

However, this approach is equivalent to the above description as the resulting monoidal product again preserves small colimits and tensors in each variable (see e.g.~\cite[Remark~C.12 and discussion after Theorem~C.10]{schwede-book}), has the correct behaviour on corepresentables \cite[Remark~C.11]{schwede-book}, and since also the structure isomorphisms of the symmetric monoidal structure are compatible with tensoring by direct computation and are given on corepresentables as above, while the monoidal unit is again corepresented by the unit of the indexing category.
\end{rk}

Sagave and Schlichtkrull \cite[Theorem~1.2]{sagave-schlichtkrull} showed that strictly commutative monoids for the box product on $\cat{$\bm{I}$-SSet}$ model all coherently commutative monoids in non-equivariant spaces.\index{coherently commutative monoid!commutative monoid in I-SSet@commutative monoid in $\cat{$\bm I$-SSet}$} Moreover, Schwede \cite[Chapter~2]{schwede-book} used a variant of this for orthogonal spaces as the basis for his approach to coherent commutativity in the global context; we will revisit his construction later in Subsection~\ref{subsec:uc-vs-uc}.

Next, we come to the box product of tame $\mathcal M$-sets and $\mathcal M$-simplicial sets, which was introduced and studied in \cite{I-vs-M-1-cat}.

\begin{defi}\index{box product!on M-SSet@on $\cat{$\bm{\mathcal M}$-SSet}$|textbf}
Let $X,Y$ be tame $\mathcal M$-sets. Their \emph{box product} is the $\mathcal M$-subset $X\boxtimes Y\subset X\times Y$ consisting of precisely those $(x,y)\in X\times Y$ such that $\supp(x)\cap\supp(y)=\varnothing$.
\end{defi}

By \cite[Proposition~2.13]{I-vs-M-1-cat} this is indeed an $\mathcal M$-subset of $X\times Y$; moreover, it becomes a subfunctor of the cartesian product, and the usual unitality, associativity, and symmetry isomorphisms restrict to make the box product the tensor product of a preferred symmetric monoidal structure on $\cat{$\bm{\mathcal M}$-Set}^\tau$ \cite[Proposition~2.14]{I-vs-M-1-cat}. Applying the box product levelwise and pulling through the $G$-actions we therefore also get a box product on $\cat{$\bm{\mathcal M}$-$\bm G$-SSet}^\tau$, and the remaining data then make this into a simplicial symmetric monoidal category.

\begin{ex}\label{ex:boxtimes-schwede}
The map $\Inj(A\amalg B,\omega)\to\Inj(A,\omega)\times\Inj(B,\omega)$ induced by the restrictions factors through an isomorphism $\Inj(A\amalg B,\omega)\cong\Inj(A,\omega)\boxtimes\Inj(B,\omega)$, see~\cite[Example~2.15]{I-vs-M-1-cat}.
\end{ex}

\begin{rk}\label{rk:ev-symmetric-monoidal}
In \cite[Proposition 4.7]{I-vs-M-1-cat}, Sagave and Schwede used the above example to construct a preferred strong symmetric monoidal structure on $\ev_\omega\colon\cat{$\bm{I}$-Set}\to\cat{$\bm{\mathcal M}$-Set}^\tau$, which yields a simplicial strong symmetric monoidal structure on $\ev_\omega\colon\cat{$\bm{I}$-SSet}\to\cat{$\bm{\mathcal M}$-SSet}^\tau$. It follows formally that the right adjoint $(\blank)_\bullet$ acquires a simplicial lax symmetric monoidal structure (which happens to be strong in this case) given by the composites
\begin{equation*}
X_\bullet\boxtimes Y_\bullet\xrightarrow{\eta}(X_\bullet\boxtimes Y_\bullet)(\omega)_\bullet\xrightarrow{\nabla^{-1}_\bullet}(X_\bullet(\omega)\boxtimes Y_\bullet(\omega))_\bullet\xrightarrow{\epsilon\boxtimes\epsilon} (X\boxtimes Y)_\bullet
\end{equation*}
\nomenclature[anabla]{$\nabla$ (nabla)}{multiplicativity map of a lax (symmetric) monoidal functor}\nomenclature[znabla]{$\nabla$}{\textit{see} $\nabla$ (nabla)\nomnorefpage}%
(i.e.~the canonical mates of the inverse structure isomorphism of $\ev_\omega$) and the unique map $\iota\colon*\to *_\bullet$.\nomenclature[aiota]{$\iota$}{unit map of a lax (symmetric) monoidal functor} In particular, we get an induced simplicial adjunction
\begin{equation*}
\ev_\omega\colon\CMon(\cat{$\bm I$-SSet})\rightleftarrows\CMon(\cat{$\bm{\mathcal M}$-SSet}^\tau):(\blank)_\bullet,
\end{equation*}
and Sagave and Schwede show in \cite[proof of Theorem~5.13]{I-vs-M-1-cat} that this is an equivalence of homotopy theories with respect to suitable weak equivalences.
\end{rk}

\subsection{The box product of tame \texorpdfstring{$\except{toc}{\bm{E\mathcal M}}\for{toc}{E\mathcal M}$}{EM}-simplicial sets}
\index{support!for EM-actions@for $E\mathcal M$-actions|(}
In order to define the box product of tame $E\mathcal M$-$G$-simplicial sets as introduced in \cite{sym-mon-global}, we need a finer notion of support:

\begin{defi}\nomenclature[asuppk]{$\supp_k$}{support of a simplex of an $E\mathcal M$-simplicial set with respect to $(k+1)$-th $\mathcal M$-action}
Let $X$ be an $E\mathcal M$-simplicial set, and let $0\le k\le n$. Then we say that $x\in X_n$ is \emph{$k$-supported} on a finite set $A\subset\omega$ if it is supported on $A$ as an element of the $\mathcal M$-set $i_k^*X_n$ where $i_k\colon\mathcal M\to \mathcal M^{n+1}$ denotes the inclusion of the $(k+1)$-th factor. We say that $x$ is \emph{$k$-finitely supported} if it is $k$-supported on some finite set $A$, in which case we define its \emph{$k$-support} $\supp_k(x)$ as its support as an element of $i_k^*X_n$.
\end{defi}

The following lemma in particular shows that our notions of support and tameness (Definition~\ref{defi:tame-EM-simplicial-set}) agree with the ones considered in \cite[Definition~2.5]{sym-mon-global}:

\begin{lemma}\label{lemma:supp-vs-supp-k}
Let $X$ be an $E\mathcal M$-simplicial set, $n\ge 0$, and $x\in X_n$. Then $X$ is supported on the finite set $A\subset\omega$ if and only if it is $k$-supported on $A$ for all $0\le k\le n$. In particular, $x$ is finitely supported if and only if it is $k$-finitely supported for all $0\le k\le n$, in which case $\supp(x)=\bigcup_{k=0}^n\supp_k(x)$.
\begin{proof}
It suffices to prove the first statement. For this we observe that if $x$ is $k$-supported on $A$ for all $0\le k\le n$, then $u.x=i_0(u).i_1(u).\dots.i_n(u).x=x$ for all $u\in\mathcal M_A$ by an immediate inductive argument, i.e.~$x$ is supported on $A$ by Theorem~\ref{thm:support-EM-vs-M}. Conversely, if $x$ is supported on $A$ and $u\in\mathcal M_A$, then $i_k(u)\in\mathcal M_A^{n+1}$, so $i_k(u).x=x$ by definition.
\end{proof}
\end{lemma}
\index{support!for EM-actions@for $E\mathcal M$-actions|)}

\begin{defi}\index{box product!on EM-SSettau@on $\cat{$\bm{E\mathcal M}$-SSet}^\tau$|textbf}
Let $X,Y$ be tame $E\mathcal M$-simplicial sets. Their \emph{box product} $X\boxtimes Y$ is the subcomplex of $X\times Y$ whose $n$-simplices are precisely those pairs $(x,y)\in X_n\times Y_n$ such that $\supp_k(x)\cap\supp_k(y)=\varnothing$ for all $0\le k\le n$.
\end{defi}

This is indeed a subcomplex and closed under the $E\mathcal M$-action by \cite[Proposition~2.17]{sym-mon-global}, and we have in addition shown at \emph{loc.~cit.} that the box product defines a subfunctor of the cartesian product on $\cat{$\bm{E\mathcal M}$-SSet}^\tau$, which is then clearly simplicial and preserves tensors. Pulling through the actions, we can extend this to $\cat{$\bm{E\mathcal M}$-$\bm G$-SSet}^\tau$, and by \cite[Proposition~2.18]{sym-mon-global} the unitality, associativity, and symmetry isomorphisms of the cartesian product restrict to yield a preferred symmetric monoidal structure on $\cat{$\bm{E\mathcal M}$-$\bm G$-SSet}^\tau$, which is then again simplicial.

\begin{ex}\label{ex:boxtimes-std}
Let $A,B$ be finite sets. In analogy to Example~\ref{ex:boxtimes-schwede}
we have an isomorphism $E\Inj(A\amalg B,\omega)\cong E\Inj(A,\omega)\boxtimes E\Inj(B,\omega)$ induced by the inclusions $A\hookrightarrow A\amalg B\hookleftarrow B$, see \cite[Lemma~5.2-(4)]{sym-mon-global}.
\end{ex}

\begin{defi}
We write $\cat{$\bm G$-ParSumSSet}\mathrel{:=}\CMon(\cat{$\bm{E\mathcal M}$-$\bm G$-SSet}^\tau)$\nomenclature[aParSumSSet]{$\cat{ParSumSSet}$}{category of parsummable simplicial sets} for the category whose objects are the commutative monoids (with respect to the box product) in $\cat{$\bm{E\mathcal M}$-$\bm G$-SSet}^\tau$, and whose morphisms are the monoid homomorphisms. We call its objects \emph{$G$-parsummable simplicial sets}.\index{parsummable|seeonly{$G$-parsummable}}\index{G-parsummable simplicial set@$G$-parsummable simplicial set|textbf}
\index{G-globally coherently commutative monoid@$G$-globally coherently commutative monoid!G-parsummable simplicial set@$G$-parsummable simplicial set|seeonly{$G$-parsummable simplicial set}}
\end{defi}

\begin{rk}\label{rk:box-unbiased}
Let $X_1,\dots,X_n$ be tame $E\mathcal M$-simplicial sets. One easily shows by induction that the image of the canonical injection $X_1\boxtimes\cdots\boxtimes X_n\to\prod_{i=1}^n X_i$ is independent of the chosen bracketing on the left and consists in degree $m$ of precisely those $(x_1,\dots,x_n)$ such that $\supp_k(x_i)\cap\supp_k(x_j)=\varnothing$ for all $0\le k\le m$ and $1\le i<j\le n$. We will from now on use the notation $X_1\boxtimes\cdots\boxtimes X_n$ for this `unbiased' iterated box product.
\end{rk}

\begin{convention}
In order to avoid excessive notation, we agree that for any category $\mathscr C$, any object $X\in\mathscr C$, and any singleton set $S$, the product $\prod_S X=X^{\times S}$ is taken to be $X$ itself, with projection the identity of $X$. In particular $X^{\boxtimes S}=X$ for any singleton set $S$ and any tame $E\mathcal M$-simplicial set $X$.
\end{convention}

\subsubsection{Homotopical properties}
\index{box product!on EM-SSettau@on $\cat{$\bm{E\mathcal M}$-SSet}^\tau$!homotopical properties|(}
We will now study the behaviour of the $G$-global weak equivalences under the above box product. While the proofs are similar to the categorical situation considered for $G=1$ in \cite{schwede-k-theory}, we nevertheless include them for completeness.

\begin{thm}\label{thm:boxtimes-em-homotopical}
Let $X,Y\in\cat{$\bm{E\mathcal M}$-$\bm G$-SSet}^\tau$. Then the inclusion $X\boxtimes Y\hookrightarrow X\times Y$ is a $G$-global weak equivalence. In particular, $\boxtimes$ preserves $G$-global weak equivalences in each variable.
\end{thm}

For the proof we will use:

\index{support!for EM-actions@for $E\mathcal M$-actions}
\begin{lemma}\label{lemma:multisupport}
Let $X$ be a tame $E\mathcal M$-simplicial set, let $0\le k\le n$, let $x\in X_n$, and let $u_0,\dots,u_n\in\mathcal M$. Then $\supp_k((u_0,\dots,u_n).x)= u_k(\supp_k(x))$.
\begin{proof}
We have
\begin{equation*}
(u_0,\dots,u_n).x=i_0(u_0).\dots.i_{k-1}(u_{k-1}).i_{k+1}(u_{k+1}).\dots.i_n(u_n).i_k(u_k).x.
\end{equation*}
For each $\ell\not=k$, the map $i_\ell(u_\ell).\blank\colon i_k^*(X_n)\to i_k^*(X_n)$ is $\mathcal M$-equivariant, and it is moreover injective by \cite[Proposition~2.7]{I-vs-M-1-cat}. It easily follows that $i_\ell(u_\ell).\blank$ preserves $k$-supports, so that $\supp_k((u_0,\dots,u_n).x)=\supp_k(i_k(u_k).x)$. The claim now follows from Lemma~\ref{lemma:support-vs-action-M}.
\end{proof}
\end{lemma}

\begin{proof}[Proof of Theorem~\ref{thm:boxtimes-em-homotopical}]
It suffices to prove the first statement. We let $H\subset\mathcal M$ be any universal subgroup and we will show that the inclusion $i$ is an $(H\times G)$-equivariant homotopy equivalence.

To this end we choose an $H$-equivariant isomorphism $\omega\cong\omega\amalg\omega$, yielding $\alpha,\beta\in\mathcal M$ centralizing $H$ and such that $\im(\alpha)\cap\im(\beta)=\varnothing$. We define $r\colon X\times Y\to X\boxtimes Y$ via $(x,y)\mapsto (\alpha.x,\beta.y)$ for all $(x,y)\in (X\times Y)_n$. This is well-defined by the previous lemma; it is moreover obviously $G$-equivariant and it is $H$-equivariant as $\alpha$ and $\beta$ centralize $H$.

The composition $ri$ is by definition given by $(\alpha.\blank)\times(\beta.\blank)\colon X\times Y\to X\times Y$. Acting with $(\alpha,1)\in E\mathcal M$ on $X$ yields an equivariant homotopy $\Delta^1\times X\to X$ from the identity to $\alpha.\blank$; analogously, $(\beta,1)\in E\mathcal M$ yields $\id\Rightarrow\beta.\blank$. Altogether we get an $(H\times G)$-equivariant homotopy $\Phi$ from the identity to $ri$ as desired.

To show that also $ir$ is $(H\times G)$-equivariantly homotopic to the identity, we only have to show that $\Phi$ restricts to $\Delta^1\times(X\boxtimes Y)\to X\boxtimes Y$. Unravelling definitions this means that for all $-1\le k\le n$ and all $x,y\in(X\boxtimes Y)_n$ also
\begin{equation}\label{eq:gen-simplex-homotopy-box}
\big((\underbrace{\alpha,\dots,\alpha}_{k+1\text{ times}},1,\dots,1).x,(\underbrace{\beta,\dots,\beta}_{k+1\text{ times}},1,\dots,1).y)\in (X\boxtimes Y)_n.
\end{equation}
But indeed, for any $0\le\ell\le k$, the previous lemma implies that
\begin{align*}
\supp_\ell((\alpha,\dots,\alpha,1,\dots,1).x)&=\alpha(\supp_\ell(x)),\\
\supp_\ell((\beta,\dots,\beta,1,\dots,1).y)&=\beta(\supp_\ell(y)),
\end{align*}
and these are disjoint as $\im\alpha\cap\im\beta=\varnothing$. On the other hand, if $\ell>k$, then by the same argument
\begin{equation*}
\supp_\ell((\alpha,\dots,\alpha,1,\dots,1).x)=\supp_\ell(x),\quad
\supp_\ell((\alpha,\dots,\alpha,1,\dots,1).y)=\supp_\ell(y)
\end{equation*}
which are again disjoint by assumption. Thus, $\Phi$ restricts to the desired homotopy between $ir$ and the identity, finishing the proof of the theorem.
\end{proof}

Applying the theorem inductively, we in particular see that the inclusion $X_1\boxtimes\cdots\boxtimes X_n\hookrightarrow X_1\times\cdots\times X_n$ is a $G$-global weak equivalence for all $X_1,\dots,X_n\in\cat{$\bm{E\mathcal M}$-$\bm G$-SSet}^\tau$. Specializing to $X_1=\cdots=X_n\mathrel{=:}X$, we see that $X^{\boxtimes n}\hookrightarrow X^{\times n}$ is a $G$-global weak equivalence. However, the right hand side has an additional $\Sigma_n$-action given by permuting the factors, and the left hand side is preserved by this. The following theorem then refines the above comparison to take this additional action into account:

\begin{thm}\label{thm:boxpower-symmetric-group}
For any $X\in\cat{$\bm{E\mathcal M}$-$\bm G$-SSet}^\tau$ the inclusion $i\colon X^{\boxtimes n}\hookrightarrow X^{\times n}$ is a $(G\times\Sigma_n)$-global weak equivalence.
\begin{proof}
Let $H\subset\mathcal M$ be any universal subgroup and let $\phi\colon H\to\Sigma_n\times G$ be a homomorphism. We have to show that $i$ induces a weak homotopy equivalence on $\phi$-fixed points. To this end let us write $\phi_1\colon H\to\Sigma_n$ for the projection of $\phi$ to the first factor. Then it is obviously enough to show that $i$ is an $(H\times G)$-equivariant homotopy equivalence where $H$ acts on both sides via $\phi_1$ and the $E\mathcal M$-action, and $G$ acts as before.

We now pick an $H$-equivariant injection $u\colon\{1,\dots,n\}\times\omega\to\omega$ where $H$ acts on $\{1,\dots,n\}$ via $\phi_1$ and on $\omega$ as before; this exists since the source is countable and the target is a complete $H$-set universe. If we write $u_i\mathrel{:=}u(i,\blank)\colon\omega\to\omega$, then the $u_i$'s are elements of $\mathcal M$ with pairwise disjoint image and the $H$-equivariance of $u$ translates into the relation
\begin{equation*}
u_{i}h=hu_{h^{-1}.i}
\end{equation*}
for all $i=1,\dots,n$. Let us now define $r\colon X^{\times n}\to X^{\boxtimes n}$ as the restriction of
\begin{equation*}
\prod_{i=1}^n (u_i.\blank)\colon X^{\times n}\to X^{\times n};
\end{equation*}
note that this indeed lands in $X^{\boxtimes n}\subset X^{\times n}$ as for any $n$-tuple $(x_1,\dots,x_n)$ in the image already $\supp(x_i)\cap\supp(x_j)=\varnothing$ for $i\not=j$ since the images of the $u_i$'s are pairwise disjoint.

It is obvious that $r$ is $G$-equivariant. But it is also $H$-equivariant: if $(x_1,\dots,x_n)$ is any $m$-simplex of $X^{\times n}$, then
\begin{equation}\label{eq:r-G-H-equivariant}
\begin{aligned}
r(h.(x_1,\dots,x_n))&=
r(h.x_{h^{-1}.1},\dots, h.x_{h^{-1}.n})\\
&=\big(u_1.(h.x_{h^{-1}.1}),\dots,u_n.(h.x_{h^{-1}.n})\big)\\
&=(h.u_{h^{-1}.1}.x_{h^{-1}.1},\dots,h.u_{h^{-1}.n}.x_{h^{-1}.n})\\
&=h.(u_1.x_1,\dots,u_n.x_n)=h.\big(r(x_1,\dots,x_n)\big).
\end{aligned}
\end{equation}
Thus it only remains to construct $(H\times G)$-equivariant homotopies $ir\simeq\id$ and $ri\simeq\id$. For the first one we observe that we have for each $i=1,\dots,n$ a homotopy $\id\Rightarrow u_i$ given by the action of $(u_i,1)\in (E\mathcal M)_1$, and these assemble into a homotopy $\Phi\colon\id\Rightarrow ir$. This homotopy is obviously $G$-equivariant, and it is moreover $H$-equivariant by a similar calculation as in $(\ref{eq:r-G-H-equivariant})$, so that it provides the desired $(H\times G)$-equivariant homotopy.

To finish the proof it is then enough to show that $\Phi$ also restricts to a homotopy $\id\Rightarrow ri$, for which one argues precisely as in the proof of Theorem~\ref{thm:boxtimes-em-homotopical}.
\end{proof}
\end{thm}

Together with Lemma~\ref{lemma:alpha-lower-star-tame} we immediately conclude:

\begin{cor}\label{cor:box-power-EM}
Let $f\colon X\to Y$ be a $G$-global weak equivalence of tame $E\mathcal M$-$G$-simplicial sets. Then $f^{\boxtimes n}\colon X^{\boxtimes n}\to Y^{\boxtimes n}$ is a $(G\times\Sigma_n)$-global weak equivalence.\qed
\end{cor}

The following corollary will be crucial later for establishing model structures on categories of commutative monoids:

\begin{cor}\label{cor:sym-homotopical-EM}
Let $f\colon X\to Y$ be a $G$-global weak equivalence of tame $E\mathcal M$-$G$-simplicial sets, and assume that neither $X$ nor $Y$ contain vertices of empty support. Then $\Sym^n(f)\mathrel{:=}f^{\boxtimes n}/\Sigma_n$\nomenclature[aSymn]{$\Sym^n(f)$}{$n$-th symmetric power of $f$, $f^{\otimes n}/\Sigma_n$} is a $G$-global weak equivalence for all $n\ge 0$.
\begin{proof}
The previous corollary asserts that $f^{\boxtimes n}$ is a $(G\times\Sigma_n)$-global weak equivalence. By Corollary~\ref{cor:free-quotient-tame} it is therefore enough to show that $\Sigma_n$ acts freely on both $X^{\boxtimes n}$ and $Y^{\boxtimes n}$.

We prove the first statement, the other one being analogous. For this it suffices to show that $\Sigma_n$ acts freely on $(X^{\boxtimes n})_0$; but indeed, if $(x_1,\dots,x_n)$ is any vertex and $i\not=j$, then $\supp(x_i)\cap\supp(x_j)=\varnothing\not=\supp(x_j)$, so $\supp(x_i)\not=\supp(x_j)$ and hence in particular $x_i\not=x_j$. If now $\sigma\in\Sigma_n$ is non-trivial, then we find an $i$ with $\sigma^{-1}(i)\not=i$, so $(x_1,\dots,x_n)$ and $\sigma.(x_1,\dots,x_n)=(x_{\sigma^{-1}(1)},\dots,x_{\sigma^{-1}(n)})$ differ in the $i$-th component.
\end{proof}
\end{cor}
\index{box product!on EM-SSettau@on $\cat{$\bm{E\mathcal M}$-SSet}^\tau$!homotopical properties|)}

\subsubsection{The box product as an operadic product} We will now give an alternative description of the box product of tame $E\mathcal M$-simplicial sets analogous to \cite[Proposition~A.17]{I-vs-M-1-cat}, following a suggestion by Stefan Schwede.

\begin{constr}
Let $X_1,\dots,X_n$ be $E\mathcal M$-simplicial sets. We consider
\begin{equation}\label{eq:operadic-box-EM}
E\Inj(\bm{n}\times\omega,\omega)\times_{(E\mathcal M)^n}(X_1\times\cdots\times X_n),
\end{equation}
i.e.~the quotient of $E\Inj(\bm{n}\times\omega,\omega)\times(X_1\times\cdots\times X_n)$ by the equivalence relation generated on $m$-simplices by
\begin{align*}
&(f_0,\dots,f_m;(u_0^{(1)},\dots,u_m^{(1)}).x_1,\dots,(u_0^{(n)},\dots,u_m^{(n)}).x_n)\\
&\quad\sim(f_0u_0^\bullet,\dots,f_mu_m^\bullet;x_1,\dots,x_n)
\end{align*}
for all injections $f_0,\dots,f_m\colon\bm{n}\times\omega\to\omega$, all $u_{k}^{(j)}\in\mathcal M$ ($k=0,\dots,m$; $j=1,\dots,n$), and $(x_1,\dots,x_n)\in (X_1\times\cdots\times X_n)_m$; here we write $u_k^\bullet\colon\bm{n}\times\omega\to\bm{n}\times\omega$ for the injection with $u_k^\bullet(j,t)=(j,u_k^{(j)}(t))$.

The simplicial set $(\ref{eq:operadic-box-EM})$ has a natural $E\mathcal M$-action given by postcomposition. Moreover, one easily checks that for all $\alpha_j\colon X_j\to X_j'$ ($j=1,\dots,n$) the map $E\Inj(\bm{n}\times\omega,\omega)\times (\alpha_1\times\cdots\times \alpha_n)$ descends to
\begin{equation*}
E\Inj(\bm{n}\times\omega,\omega)\times_{(E\mathcal M)^n}(X_1\times\cdots\times X_n)\to E\Inj(\bm{n}\times\omega,\omega)\times_{(E\mathcal M)^n}(X_1'\times\cdots\times X_n'),
\end{equation*}
and that this yields a functor $(\cat{$\bm{E\mathcal M}$-SSet})^n\to\cat{$\bm{E\mathcal M}$-SSet}$.

Finally, we have a natural $E\mathcal M$-equivariant map
\begin{equation*}
\Phi\colon E\Inj(\bm{n}\times\omega,\omega)\times_{(E\mathcal M)^n}(X_1\times\cdots\times X_n)\to X_1\times\cdots\times X_n
\end{equation*}
given on $m$-simplices by
\begin{equation*}
[f_0,\dots,f_m;x_1,\dots,x_n]\mapsto ((f_0\iota_1,\dots,f_m\iota_1).x_1,\dots,(f_0\iota_n,\dots,f_m\iota_n).x_n)
\end{equation*}
where $\iota_j\colon\omega\to\bm{n}\times\omega$ is defined via $\iota_j(t)=(j,t)$.
\end{constr}

\begin{thm}\label{thm:operadic-box-EM}\index{box product!on EM-SSettau@on $\cat{$\bm{E\mathcal M}$-SSet}^\tau$!as an operadic product|textbf}
For any $X_1,\dots,X_n\in\cat{$\bm{E\mathcal M}$-SSet}^\tau$ the map $\Phi$ restricts to an isomorphism $E\Inj(\bm{n}\times\omega,\omega)\times_{(E\mathcal M)^n}(X_1\times\cdots\times X_n)\to X_1\boxtimes\cdots\boxtimes X_n$.
\end{thm}

For the proof we will need:

\begin{lemma}\label{lemma:operadic-box-EM-restriction}
In the above situation, let $(x_1,\dots,x_n)\in(X_1\times\cdots\times X_n)_m$ and $f_0,\dots,f_m;g_0,\dots,g_m\in\Inj(\bm{n}\times\omega,\omega)$ such that $f_k$ and $g_k$ agree for every $k$ on $\bigcup_{j=1}^n\{j\}\times\supp_k(x_j)$. Then $[f_0,\dots,f_m;x_1,\dots,x_n]=[g_0,\dots,g_m;x_1,\dots,x_n]$ in $E\Inj(\bm{n}\times\omega,\omega)\times_{(E\mathcal M)^n}(X_1\times\cdots\times X_n)$.
\begin{proof}
By induction we may assume that there exists a $k'$ such that $f_k=g_k$ for all $k\not=k'$, and by symmetry we may assume that $k'=0$. We now consider
\begin{align*}
\alpha\colon\Inj(\bm{n}\times\omega,\omega)&\to \big(E\Inj(\bm{n}\times\omega,\omega)\times_{(E\mathcal M)^n}(X_1\times\cdots\times X_n)\big)_m\\
h&\mapsto [h, f_1,\dots,f_m;x_1,\dots,x_n]=[h,g_1,\dots,g_m;x_1,\dots,x_n].
\end{align*}
We claim that this factors over $\Inj(\bm{n}\times\omega,\omega)/(\mathcal M_{\supp_0(x_1)}\times\cdots\times\mathcal M_{\supp_0(x_n)})$: indeed, given $u^{(j)}\in\mathcal M_{\supp_0(x_j)}$ for $1\le j\le n$,
\begin{align*}
[h, f_1,\dots,f_m;x_1,\dots,x_n]&=[h,f_1,\dots,f_m,(u^{(1)},1,\dots,1).x_1,\dots,(u^{(n)},1,\dots,1).x_n]\\
&=[h u^\bullet,f_1,\dots,f_m;x_1,\dots,x_n]
\end{align*}
by definition, i.e.~$\alpha$ is compatible with generating relations and hence factors over the quotient. But $[f_0]=[g_0]$ in $\Inj(\bm{n}\times\omega,\omega)/(\mathcal M_{\supp_0(x_1)}\times\dots\times\mathcal M_{\supp_0(x_n)})$ by \cite[Lemma~A.5]{I-vs-M-1-cat}, which completes the proof the lemma.
\end{proof}
\end{lemma}

\begin{proof}[Proof of Theorem~\ref{thm:operadic-box-EM}]
Lemma~\ref{lemma:multisupport} readily implies that the image of $\Phi$ is contained in $X_1\boxtimes\cdots\boxtimes X_n$. Let us now prove the other inclusion: if $(x_1,\dots,x_n)$ is any $m$-simplex of $X_1\boxtimes\cdots\boxtimes X_n$, then we pick for each $k$ an injection $f_k\colon\bm{n}\times\omega\to\omega$ with $f_k(i,t)=t$ for all $t\in\supp_k(x_i)$; this is possible as the sets $\supp_k(x_i)$ are pairwise disjoint and finite for fixed $k$ and varying $i$. It is then easy to check that $\Phi[f_0,\dots,f_m;x_1,\dots,x_n]=(x_1,\dots,x_n)$.

For the proof of injectivity, we pick any $f_0,\dots,f_m;g_0,\dots,g_m\in\Inj(\bm{n}\times\omega,\omega)$ and $(x_1,\dots,x_n),\allowbreak(y_1,\dots,y_n)\in X_1\times\cdots\times X_n$ such that $\Phi[f_0,\dots,f_m;x_1,\dots,x_n]=\Phi[g_0,\dots,g_m;y_1,\dots,y_n]$.

\begin{claim*}
If $f_k=g_k$ for all $k=0,\dots,m$, then also $x_j=y_j$ for all $j=1,\dots,n$.
\begin{proof}
For any $j=1,\dots,n$, we have
\begin{equation*}
(f_0\iota_j,\dots,f_m\iota_j).x_j=(g_0\iota_j,\dots,g_m\iota_j).y_j=(f_0\iota_j,\dots,f_m\iota_j).y_j.
\end{equation*}
The claim follows as $(f_0\iota_j,\dots,f_m\iota_j).\blank\colon X_1\times\cdots\times X_n$ is injective as composition of the injections $i_0(f_0\iota_j).\blank,\dots,i_m(f_m\iota_j).\blank$ (see \cite[Proposition~2.7]{I-vs-M-1-cat}).
\end{proof}
\end{claim*}

In the general case, we observe that for each $1\le j\le n$ and $0\le k\le m$ in particular
\begin{equation*}
\supp_k\big((f_0\iota_j,\dots,f_m\iota_j).x_j\big)=\supp_k\big((g_0\iota_j,\dots,g_m\iota_j).y_j\big),
\end{equation*}
hence $(f_k\iota_j)(\supp_kx_j)=(g_k\iota_j)(\supp_ky_j)$ by Lemma~\ref{lemma:multisupport}. By injectivity of $f_k\iota_j$ and $g_k\iota_j$ it follows that there exists a (unique) bijection $\sigma^{(j)}_{k}\colon\supp_k(x_j)\to\supp_k(y_j)$ with $g_k\iota_j\sigma^{(j)}_k=f_k\iota_j|_{\supp_k(x_j)}$.

We now pick an extension of $\sigma^{(j)}_k$ to an injection $s^{(j)}_k\in\mathcal M$. Letting $j$ vary, we conclude from the above definition that $g_ks^\bullet_k$ and $f_k$ agree on $\bigcup_{j=1}^n\{j\}\times\supp_k(x_j)$. Letting $k$ vary, Lemma~\ref{lemma:operadic-box-EM-restriction} therefore shows
\begin{align*}
{}[f_0,\dots,f_m;x_1,\dots,x_n]&=[g_0s^\bullet_0,\dots,g_ms^\bullet_m;x_1,\dots,x_n]\\
&=[g_0,\dots,g_m;(s^{(1)}_0,\dots,s^{(1)}_m).x_1,\dots,(s^{(n)}_0,\dots,s^{(n)}_m).x_n].
\end{align*}
But then
\begin{align*}
\Phi[g_0,\dots,g_m;y_1,\dots,y_n]&=
\Phi[f_0,\dots,f_m;x_1,\dots,x_n]\\
&=\Phi[g_0,\dots,g_m;(s^{(1)}_0,\dots,s^{(1)}_m).x_1,\dots,(s^{(n)}_0,\dots,s^{(n)}_m).x_n],
\end{align*}
so the above claim implies that $(s^{(j)}_0,\dots,s^{(j)}_m).x_j=y_j$, finishing the proof.
\end{proof}

\begin{cor}\label{cor:box-EM-cocontinuous}
The box product of tame $E\mathcal M$-simplicial sets is cocontinuous in each variable.
\begin{proof}
By Theorem~\ref{thm:operadic-box-EM} it suffices to check the corresponding statement for $E\Inj(\bm{2}\times\omega,\omega)\times_{(E\mathcal M)^2}(\blank\times\blank)$, where this is obvious.
\end{proof}
\end{cor}

\subsection{Ultra-commutative monoids}
\index{model structure on commutative monoids|(}
Purely by abstract nonsense, there is a unique way to extend the symmetric monoidal structure on $I$ to a simplicially enriched symmetric monoidal structure on $\mathcal I$. In particular, we again get a simplicial symmetric monoidal structure on $\cat{$\bm{\mathcal I}$-SSet}$ with monoidal product characterized by $\mathcal I(A,\blank)\boxtimes\mathcal I(B,\blank)=\mathcal I(A\amalg B,\blank)$ and structure isomorphisms induced by those of $\mathcal I$.\index{box product!on II-SSet@on $\cat{$\bm{\mathcal I}$-SSet}$}

\begin{prop}\label{prop:ev-omega-strong-sym-mon}
There exists a unique enhancement of $\ev_\omega\colon\cat{$\bm{\mathcal I}$-SSet}\to\cat{$\bm{E\mathcal M}$-SSet}^\tau$ to a simplicial strong symmetric monoidal functor such that for all finite sets $A,B$ the inverse structure isomorphism
\begin{equation}\label{eq:boxtimes-ev-representable}
\big(\mathcal I(A,\blank)\boxtimes\mathcal I(B,\blank)\big)(\omega)=\mathcal I(A\amalg B,\blank)(\omega)\to \mathcal I(A,\blank)(\omega)\boxtimes\mathcal I(B,\blank)(\omega)
\end{equation}
is induced by the map $\mathcal I(A\amalg B,\blank)\to\mathcal I(A,\blank)\times\mathcal I(B,\blank)$ given by restriction.
\begin{proof}
As seen in Example~\ref{ex:boxtimes-std}, the restrictions induce an isomorphism $E\Inj(A\amalg B,\omega)\cong E\Inj(A,\omega)\boxtimes E\Inj(B,\omega)$. Using the canonical isomorphisms $\mathcal I(C,\blank)(\omega)\cong E\Inj(C,\omega)$ for $C\in\mathcal I$, we therefore conclude that $(\ref{eq:boxtimes-ev-representable})$ is well-defined and an isomorphism.

As $\ev_\omega$ is cocontinuous and preserves tensors, the universal property of enriched presheaves now easily implies that the inverses of $(\ref{eq:boxtimes-ev-representable})$ extend to a unique simplicial natural isomorphism $\nabla\colon\ev_\omega\boxtimes\ev_\omega\Rightarrow\ev_\omega(\blank\boxtimes\blank)$. On the other hand, there is a unique isomorphism $\iota\colon *\to\ev_\omega\mathcal I(\varnothing,\blank)$ as both sides are terminal.

Finally, to show that $\nabla$ and $\iota$ endow $\ev_\omega$ with the structure of a simplicial strong symmetric monoidal functor, we simply observe that it again suffices to check the coherence conditions on corepresentables, where this is a trivial calculation.
\end{proof}
\end{prop}

We can define the box product of $G$-$\mathcal I$-simplicial sets by pulling through the $G$-action again. The above then induces a simplicial strong symmetric monoidal structure on $\ev_\omega\colon\cat{$\bm G$-$\bm{\mathcal I}$-SSet}\to\cat{$\bm{E\mathcal M}$-$\bm G$-SSet}^\tau$.

\begin{cor}\label{cor:boxtimes-homotopical-script-I}
The box product of $G$-$\mathcal I$-simplicial sets is homotopical in each variable.
\begin{proof}
Let $f\colon X\to X'$ and $g\colon Y\to Y'$ be $G$-global weak equivalences. We have to show that $f\boxtimes g$ is a $G$-global weak equivalence, i.e.~that $(f\boxtimes g)(\omega)$ is a $G$-global weak equivalence of $E\mathcal M$-$G$-simplicial sets. But by the previous proposition this is conjugate to $f(\omega)\boxtimes g(\omega)$, and by assumption both $f(\omega)$ and $g(\omega)$ are $G$-global weak equivalences. The claim therefore follows from Theorem~\ref{thm:boxtimes-em-homotopical}.
\end{proof}
\end{cor}

\begin{cor}\label{cor:sym-homotopical-script-I}
Let $f\colon X\to Y$ be a $G$-global weak equivalence of $G$-$\mathcal I$-simplicial sets such that $X$ and $Y$ are cofibrant in the \emph{positive} $G$-global model structure. Then $\Sym^nf\mathrel{:=}f^{\boxtimes n}/\Sigma_n$ is a $G$-global weak equivalence for all $n\ge 0$.
\begin{proof}
As in the previous corollary, we see that $f^{\boxtimes n}(\omega)$ agrees with $f(\omega)^{\boxtimes n}$ up to conjugation by $(G\times\Sigma_n)$-equivariant isomorphisms (the $\Sigma_n$-equivariance uses that we have a strong \emph{symmetric} monoidal structure). As $\ev_\omega$ is a left adjoint, it commutes with quotients, so we conclude that also $(\Sym^nf)(\omega)$ agrees with $\Sym^n(f(\omega))$ up to conjugation.

By Corollary~\ref{cor:sym-homotopical-EM} we are therefore reduced to showing that $Z(\omega)$ has no elements of empty support for any $Z$ cofibrant in the positive $G$-global model structure. This follows from Lemma~\ref{lemma:g-global-pos-cof} as $\ev_\omega\colon\cat{$\bm G$-$\bm{\mathcal I}$-SSet}\to\cat{$\bm{E\mathcal M}$-$\bm G$-SSet}^\tau$ is left Quillen for the positive $G$-global model structures.
\end{proof}
\end{cor}

\begin{defi}\label{defi:ultra-commutative-monoid}\index{ultra-commutative monoid!in II-SSet@in $\cat{$\bm{\mathcal I}$-SSet}$|seeonly{$G$-ultra-commutative monoid}}
\index{G-ultra-commutative monoid@$G$-ultra-commutative monoid|textbf}
\index{G-globally coherently commutative monoid@$G$-globally coherently commutative monoid!G-ultracommutative monoid@$G$-ultra-commutative monoid|seeonly{$G$-ultra-commutative monoid}}
\nomenclature[aUCom]{$\cat{UCom}$}{category of ultra-commutative monoids}
We write $\cat{$\bm G$-UCom}\mathrel{:=}\CMon(\cat{$\bm G$-$\bm{\mathcal I}$-SSet})$ and call its objects \emph{$G$-ultra-commutative monoids}.
\end{defi}

\subsubsection{A reminder on model structures for commutative monoids}\label{subsubsection:cmon-model-structures}
We want to construct a suitable $G$-global model structure on $\cat{$\bm G$-UCom}$, for which we need to recall the general machinery of \cite{white-cmon} and \cite{sym-powers}.

\begin{defi}
A model category $\mathscr C$ equipped with a closed symmetric monoidal structure is called a \emph{symmetric monoidal model category}\index{symmetric monoidal model category|textbf} if the following two conditions are satisfied:
\begin{enumerate}
\item (\textit{Unit Axiom})\index{Unit Axiom|textbf} If $X\in\mathscr C$ is cofibrant, then $I\otimes X\to \textbf{1}\otimes X$ is a weak equivalence for some (hence any) cofibrant replacement $I\xrightarrow\sim\textbf{1}$ of the unit.
\item (\textit{Pushout Product Axiom})\index{Pushout Product Axiom!for symmetric monoidal model categories|textbf} If $i\colon X\to X'$ and $j\colon Y\to Y'$ are cofibrations, then so is the pushout product map $i\ppo j\colon(X\otimes Y')\amalg_{X\otimes Y}(X'\otimes Y)\to X'\otimes Y'$. Moreover, if at least one of $i$ and $j$ is acyclic, then so is $i\ppo j$.
\end{enumerate}
\end{defi}

It is a well-known fact that for cofibrantly generated $\mathscr C$ it suffices to check the Pushout Product Axiom on some chosen generating (acyclic) cofibrations, see e.g~\cite[Lemma~3.5-(1)]{schwede-shipley-monoidal}.

\begin{defi}
A symmetric monoidal model category $\mathscr{C}$ satisfies the \emph{Monoid Axiom}\index{Monoid Axiom|textbf} if every transfinite composition of pushouts of maps of the form $X\otimes j$ with $X\in\mathscr C$ and $j$ an acyclic cofibration is a weak equivalence again.
\end{defi}

Again, if $\mathscr C$ is cofibrantly generated, then it is enough to restrict to the case that $j$ belongs to a chosen set $J$ of generating acyclic cofibrations, i.e.~the above condition can be reformulated as saying that any relative $(\mathscr C\otimes J)$-cell complex is a weak equivalence, see~\cite[Lemma~3.5-(2)]{schwede-shipley-monoidal}.

It is a classical result of Schwede and Shipley \cite[Theorem~4.1-(3)]{schwede-shipley-monoidal} that for a combinatorial symmetric monoidal model category $\mathscr C$ satisfying the Monoid Axiom the transferred model structure on the category $\Mon(\mathscr C)$ of not necessarily commutative monoids exists, i.e.~$\Mon(\mathscr C)$ can be equipped with a model structure in which the fibrations and weak equivalences are created in $\mathscr C$. We will need a version of this for commutative monoids due to White \cite{white-cmon}. This relies on the following additional notion:

\begin{constr}\label{constr:iterated-ppos}
Let $\mathscr C$ be a cocomplete symmetric monoidal category. As usual in the context of model structures on commutative monoids, we will suppress the associativity isomorphisms to keep the notation tractable. The reader uncomfortable with this may rest assured that we will only ever apply this to box products, which have a preferred unbiased $n$-ary tensor product anyhow.

For $n\ge 1$ we let $C_n$ denote the $n$-cube, i.e.~the poset of subsets of $\{1,\dots,n\}$. If $f\colon X\to Y$ is any morphism in $\mathscr C$, then we have a functor $K^n_n(f)\colon C_n\to\mathscr C$ sending $I\subset\{1,\dots,n\}$ to $Z_1\otimes\cdots\otimes Z_n$ where $Z_i=Y_i$ if $i\in I$, and $Z_i=X_i$ otherwise; the structure maps are given by the obvious tensor products of $f$ and the respective identities.

We write $Q^n_{n-1}(f)$ for the colimit of the subdiagram obtained by removing the terminal vertex of $C_n$. This has a natural $\Sigma_n$-action obtained by the $\Sigma_n$-action on the punctured $n$-cube and the symmetry isomorphisms of the tensor product on $\mathscr C$. The natural map $f^{\ppo n}\colon Q^n_{n-1}(f)\to Y^{\otimes n}$ induced by the remaining structure maps of $K^n_n$ is then equivariant with respect to the $\Sigma_n$-action on the target via permuting the tensor factors.
\end{constr}

\begin{rk}
Note that $f^{\ppo n}$ would be usually used to denote the iterated pushout product of $f$ with itself, which is indeed canonically and $\Sigma_n$-equivariantly conjugate to the above whenever $\otimes$ is cocontinuous in each variable, see e.g.~\cite[Section~3]{sym-powers} (in fact, it suffices that pushouts are preserved). However, we will only care about the above description except in the proof of Corollary~\ref{cor:iterated-ppo-injective}.
\end{rk}

\begin{defi}
A symmetric monoidal model category $\mathscr C$ satisfies the \emph{Strong Commutative Monoid Axiom}\index{Strong Commutative Monoid Axiom|textbf} if $i^{\ppo n}/\Sigma_n$ is a cofibration for all $n\ge 0$ and any cofibration $i$ in $\mathscr C$, and an acyclic cofibration if $i$ is.
\end{defi}

Although this is much harder to prove than the analogous statements for the previous axioms we still have:

\begin{lemma}
Let $\mathscr C$ be a cofibrantly generated symmetric monoidal model category, let $I$ be a set of generating cofibrations, and let $J$ be a set of generating acyclic cofibrations. Then $\mathscr C$ satisfies the Strong Commutative Monoid Axiom provided that $i^{\ppo n}/\Sigma_n$ is a cofibration for each $i\in I$ and that $j^{\ppo n}/\Sigma_n$ is an acyclic cofibration for all $j\in J$.
\begin{proof}
See \cite[Lemma~A.1]{white-cmon} or \cite[Corollary~9]{sym-powers}.
\end{proof}
\end{lemma}

In many practical situations the following lemma due to Gorchinskiy and Guletski\u\i\ further simplifies the verification of the Strong Commutative Monoid Axiom:

\begin{lemma}\label{lemma:sym-instead-of-ppo-powers}
Let $\mathscr C$ be a cofibrantly generated model category, let $I$ be a set of generating cofibrations, and let $J$ be a set of generating acyclic cofibrations such that all maps in $J$ \emph{have cofibrant source}. Then $\mathscr C$ satisfies the Strong Commutative Monoid Axiom provided that $i^{\ppo n}/\Sigma_n$ is a cofibration for all $i\in I$ and that $\Sym^nj=j^{\otimes n}/\Sigma_n$ is a weak equivalence for all $j\in J$.
\begin{proof}
By the argument from the proof of the previous lemma, i.e.~by \cite[Corollary~9]{sym-powers}, $i^{\ppo n}/\Sigma_n$ is actually a cofibration for \emph{all} cofibrations $i$. It remains to show that the cofibration $j^{\ppo n}/\Sigma_n$ is acyclic for each $j\in J$. But as $j$ is a map of cofibrant objects, this follows from \cite[Corollaries~10 and~23]{sym-powers}.
\end{proof}
\end{lemma}

The forgetful functor $\CMon(\mathscr C)\to\mathscr C$ has a left adjoint $\textbf P$\nomenclature[aP]{$\cat P$}{free commutative monoid} given on objects by
\begin{equation*}
X\mapsto\coprod_{n\ge 0}X^{\otimes n}/\Sigma_n=\coprod_{n\ge 0}\Sym^nX
\end{equation*}
and analogously on morphisms; the monoid structure is given by concatenation.

\begin{thm}[White]\label{thm:cmon-existence}\index{model structure on commutative monoids|textbf}
Let $\mathscr C$ be a combinatorial symmetric monoidal model category satisfying both the Monoid Axiom and the Strong Commutative Monoid Axiom. Then there exists a unique model structure on $\CMon(\mathscr C)$ in which a map is a weak equivalence or fibration if and only if it so in $\mathscr C$. This model structure is again combinatorial, and when $I$ and $J$ are generating cofibrations and generating acyclic cofibrations of $\mathscr C$, then $\free I$ and $\free J$ are sets of generating cofibrations and generating acyclic cofibrations, respectively, for $\CMon(\mathscr C)$.
\begin{proof}
See~\cite[Theorem~3.2]{white-cmon} and the discussion after it.
\end{proof}
\end{thm}

\begin{thm}[White]\label{thm:cmon-left-proper}
In the situation of the previous theorem, $\CMon(\mathscr C)$ is left proper provided that the following additional conditions are satisfied:
\begin{enumerate}
\item $\mathscr C$ is left proper and filtered colimits in it are homotopical. Moreover, there exists a set of generating cofibrations \emph{with cofibrant sources}.
\item If $X\in\mathscr C$ and $i$ is any cofibration, then pushouts in $\mathscr C$ along $X\otimes i$ are homotopy pushouts. Moreover, if $X$ is cofibrant, then $X\otimes\blank$ is homotopical.
\end{enumerate}
\begin{proof}
This is a special case of \cite[Theorem~4.17]{white-cmon}, also see \cite[discussion after Definition~4.15]{white-cmon} and \cite[discussion after Definition~2.4]{batanin-berger}.
\end{proof}
\end{thm}
\index{model structure on commutative monoids|)}
\index{model structure on commutative monoids|seealso{positive $G$-global model structure}}

\subsubsection{Construction of the model structure} We can now prove:

\begin{thm}\label{thm:G-UCom}
There is a unique model structure on $\cat{$\bm G$-UCom}$ in which a map is a weak equivalence or fibration if and only if it is so in the positive $G$-global model structure on $\cat{$\bm G$-$\bm{\mathcal I}$-SSet}$.

\index{G-ultra-commutative monoid@$G$-ultra-commutative monoid!model structure|seeonly{positive $G$-global model structure, on $\cat{$\bm G$-UCom}$}}
\index{positive G-global model structure@positive $G$-global model structure!on G-UCom@on $\cat{$\bm G$-UCom}$|textbf}\index{ultra-commutative monoid!in II-SSet@in $\cat{$\bm{\mathcal I}$-SSet}$!model structure|seeonly{positive $G$-global model structure, on $\cat{$\bm G$-UCom}$}}
We call this the \emph{positive $G$-global model structure}. It is proper, simplicial, and combinatorial with generating cofibrations (say, as maps in $\cat{$\bm G$-$\bm{\mathcal I}$-SSet}$)
\begin{equation*}
\coprod_{n\ge 0}\big(\mathcal I(\bm{n}\times A,\blank)\times G^n\times(\del\Delta^m\hookrightarrow\Delta^m)^{\times n}\big)/(\Sigma_n\wr H),
\end{equation*}
where $m\ge 0$, $H$ runs through all finite groups, and $A$ through finite faithful non-empty $H$-sets. Moreover, filtered colimits in this model category are homotopical.

Finally, the adjunction
\begin{equation*}
\textbf{\textup P}\colon\cat{$\bm G$-$\bm{\mathcal I}$-SSet}_{\textup{positive $G$-global}}\rightleftarrows\cat{$\bm G$-UCom}_{\textup{positive $G$-global}} :\!\forget
\end{equation*}
is a simplicial Quillen adjunction.
\begin{proof}
Let us first establish the model structure, for which it suffices to verify the assumptions of Theorem~\ref{thm:cmon-existence}.

The Unit Axiom is immediate because $\boxtimes$ is homotopical by Corollary~\ref{cor:boxtimes-homotopical-script-I}.

In order to verify the Pushout Product Axiom for cofibrations, we may restrict to the standard generating cofibrations. We therefore let $H_1,H_2$ be finite groups, we let $A_i$ ($i=1,2$) be a finite faithful non-empty $H_i$-set, we let $\phi_i\colon H_i\to G$ ($i=1,2$) be a group homomorphism, and we let $n_1,n_2\ge0$. We consider $\mathcal I(A_1,\blank)\times G$ with $H_1$ acting from the right on the first factor via its action on $A_1$ and on the second factor via $\phi_1$; similarly, we equip $\mathcal I(A_2,\blank)\times G$ with a right $H_2$-action. Then
\begin{equation*}
\big(\mathcal I(A_1,\blank)\times G\times(\del\Delta^{n_1}\hookrightarrow\Delta^{n_1})\big)\ppo
\big(\mathcal I(A_2,\blank)\times G\times(\del\Delta^{n_2}\hookrightarrow\Delta^{n_2})\big)
\end{equation*}
agrees up to conjugation by the evident isomorphisms with
\begin{equation*}
\mathcal I(A_1\amalg A_2,\blank)\times\big((G\times\del\Delta^{n_1}\hookrightarrow G\times\Delta^{n_1})\ppo(G\times\del\Delta^{n_2}\hookrightarrow G\times\Delta^{n_2})\big).
\end{equation*}
These isomorphisms are equivariant in $H_1$, $H_2$, and $G$ if we let $G$ act in the obvious way, $H_1$ via its action on $A_1$ and its action on the first $G$-factor, and $H_2$ via its action on $A_2$ and its action on the second $G$-factor. All the $H_1$-actions commute with all the $H_2$-actions, so that they assemble into an $(H_1\times H_2)$-action. As $\boxtimes$ preserves colimits in each variable, we conclude that
\begin{equation*}
\big((\mathcal I(A_1,\blank)\times_{\phi_1}G)\times(\del\Delta^{n_1}\hookrightarrow\Delta^{n_1})\big)\ppo\big((\mathcal I(A_2,\blank)\times_{\phi_2}G)\times(\del\Delta^{n_2}\hookrightarrow\Delta^{n_2})\big)
\end{equation*}
agrees up to conjugation by isomorphisms with
\begin{equation*}
\big(\mathcal I(A_1\amalg A_2,\blank)\times\big((G\times\del\Delta^{n_1}\hookrightarrow G\times\Delta^{n_1})\ppo(G\times\del\Delta^{n_2}\hookrightarrow G\times\Delta^{n_2})\big)\big)/(H_1\times H_2).
\end{equation*}
The $(H_1\times H_2)$-set $A_1\amalg A_2$ is obviously faithful and non-empty, while the map $(G\times\del\Delta^{n_1}\hookrightarrow G\times\Delta^{n_1})\ppo(G\times\del\Delta^{n_2}\hookrightarrow G\times\Delta^{n_2})$ is injective by the Pushout Product Axiom for $\cat{SSet}$; moreover, $G$ clearly acts freely on the target. We claim that this already implies that the above is a positive $G$-global cofibration. Indeed, for any faithful $H$-set $A\not=\varnothing$ the adjunction
\begin{equation*}
\mathcal I(A,\blank)\times_H\blank\colon\cat{$\bm{(G\times H^\op)}$-SSet}\rightleftarrows\cat{$\bm G$-$\bm{\mathcal I}$-SSet}:\ev_A
\end{equation*}
is a Quillen adjunction with respect to the $\mathcal G_{H^\op,G}$-equivariant model structure on the source and the positive level model structure on the target, hence in particular with respect to the $G$-global positive model structure. The claim follows as the cofibrations on the left hand side are precisely the injections $i$ such that $G$ acts freely outside the image of $i$.

This proves the Pushout Product Axiom for cofibrations. For the part about acyclic cofibrations we will use:

\begin{claim*}
Let $i\colon A\to B$ be an injective cofibration in $\cat{$\bm G$-$\bm{\mathcal I}$-SSet}$ and let $X$ be arbitrary. Then any pushout along $X\boxtimes i$ is a homotopy pushout.
\begin{proof}
Lemma~\ref{lemma:G-global-I-injective-pushout} shows that $\ev_\omega$ creates homotopy pushouts, so it is enough to show that any pushout along $(X\boxtimes i)(\omega)$ in $\cat{$\bm{E\mathcal M}$-$\bm G$-SSet}^\tau$ is a homotopy pushout. But $\ev_\omega$ is strong symmetric monoidal, so this is conjugate to $X\boxtimes i(\omega)$ which is evidently an injective cofibration. The claim follows immediately.
\end{proof}
\end{claim*}

If now $i_1\colon X_1\to Y_1$ is any cofibration and $i_2\colon X_2\to Y_2$ is any acyclic cofibration, then we consider the commutative diagram
\begin{equation*}
\begin{tikzcd}
X_1\boxtimes X_2\arrow[dr, phantom, "\ulcorner", very near end]\arrow[d, "X_1\boxtimes i_2"']\arrow[r, "i_1\boxtimes X_2"] &[1em] Y_1\boxtimes X_2\arrow[d]\arrow[ddr, "Y_1\boxtimes i_2", bend left=10pt]\\
X_1\boxtimes Y_2\arrow[r]\arrow[rrd, "i_1\boxtimes Y_2"', bend right=10pt] & P\arrow[dr, "i_1\ppo i_2" description]\\
& & Y_1\boxtimes Y_2.
\end{tikzcd}
\end{equation*}
By the above, the top square is a homotopy pushout, and the left hand arrow as well as the rightmost arrow are weak equivalences as $\boxtimes$ is homotopical (Corollary~\ref{cor:boxtimes-homotopical-script-I}). We conclude that $Y_1\boxtimes X_2\to P$ is a weak equivalence, and hence so is $i_1\ppo i_2$ by $2$-out-of-$3$. Since we already know that it is a cofibration, it is therefore an acyclic cofibration, which completes the verification of the Pushout Product Axiom.

Let us now verify the Monoid Axiom. If $j$ is any acyclic cofibration, then it is in particular an injective cofibration. If now $X$ is any object, then pushouts along $X\boxtimes j$ are therefore homotopy pushouts by the above claim. However, $X\boxtimes j$ is also a weak equivalence because $\boxtimes$ is homotopical, so any pushout of it is a weak equivalence. The claim follows as filtered colimits in $\cat{$\bm G$-$\bm{\mathcal I}$-SSet}$ are homotopical.

Next, we consider the Strong Commutative Monoid Axiom. We let $H$ be a finite group, $A$ a finite faithful non-empty $H$-set, and $m,n\ge 0$. Then
\begin{equation}\label{eq:ppo-power-script-I}
\big((\mathcal I(A,\blank)\times G)\times(\del\Delta^m\hookrightarrow\Delta^m)\big)^{\ppo n}
\end{equation}
agrees up to conjugation by the evident isomorphisms with
\begin{equation}\label{eq:ppo-power-script-I-simplified}
\mathcal I(\bm n\times A,\blank)\times (G\times\del\Delta^m\hookrightarrow G\times\Delta^m)^{\ppo n}.
\end{equation}
Now let $\phi\colon H\to G$ be any group homomorphism. There are $n$ commuting $H$-actions on $(\ref{eq:ppo-power-script-I})$, the $i$-th one of which is given by acting on the $i$-th $\ppo$-factor in the obvious way. We similarly have $n$ commuting $H$-actions on $\bm{n}\times A$ with the $i$-th action given by acting in the prescribed way on the $i$-th copy of $A$ and trivially on all other copies, and we have $n$ further commuting $H$-actions on $(G\times\del\Delta^m\hookrightarrow\Delta^m)^{\ppo n}$ analogously to the above. By taking the respective diagonals we get $n$ commuting $H$-actions on $(\ref{eq:ppo-power-script-I-simplified})$, and the above isomorphism is equivariant with respect to all these actions. Finally, $\Sigma_n$ acts on $\bm{n}\times A$ via its tautological action on $\bm{n}$ and on the $\ppo$-powers via permuting the factors. Again taking the diagonal for $(\ref{eq:ppo-power-script-I-simplified})$ we get $\Sigma_n$-actions compatible with the above identification.

We now observe that the $\Sigma_n$-actions and the $H$-actions assemble into $(\Sigma_n\wr H)$-actions on both $(\ref{eq:ppo-power-script-I})$ and $(\ref{eq:ppo-power-script-I-simplified})$, so we altogether conclude that
\begin{equation*}
\big((\mathcal I(A,\blank)\times_\phi G)\times(\del\Delta^m\hookrightarrow\Delta^m)\big)^{\ppo n}/\Sigma_n
\end{equation*}
agrees up to conjugation by isomorphisms with
\begin{equation}\label{eq:ppo-power-script-I-quotient}
\big(\mathcal I(\bm n\times A,\blank)\times(G\times\del\Delta^m\hookrightarrow G\times\Delta^m)^{\ppo n}\big)/(\Sigma_n\wr H)
\end{equation}
But the $(\Sigma_n\wr H)$-action on $\bm{n}\times A$ is faithful: if $(\sigma;h_1,\dots,h_n)$ acts trivially on $\bm{n}\times A$, then each $h_i$ has to act trivially on $A$, so $h_i=1$ for all $i$. But then if $a\in A$ is arbitrary (here we used that $A\not=\varnothing$), then $(\sigma;h_1,\dots,h_n)=\sigma$ sends $(i,a)$ to $(\sigma(i),a)$ for each $i$, hence also $\sigma=1$. Thus, $(\ref{eq:ppo-power-script-I-quotient})$ is a positive $G$-global cofibration for $n>0$ by the same argument as in the verification of the Pushout Product Axiom; for $n=0$, it is obviously even an isomorphism.

For the Strong Commutative Monoid Axiom for acyclic cofibrations, we observe that the positive $G$-global model structure on $\cat{$\bm G$-$\bm{\mathcal I}$-SSet}$ is combinatorial and that the standard generating cofibrations have cofibrant sources. By \cite[Corollary~2.7]{barwick-tractable} we may conclude that there exists a set of generating acyclic cofibrations $J$ with cofibrant sources. By Lemma~
\ref{lemma:sym-instead-of-ppo-powers} it therefore suffices that for each $j\in J$ and $n\ge 0$ the map $\Sym^nj$ is a weak equivalence. But as the source (and hence also the target) of $j$ was assumed to be cofibrant in the positive $G$-global model structure, this is an instance of Corollary~\ref{cor:sym-homotopical-script-I}.

This completes the verification of the assumptions of Theorem~\ref{thm:cmon-existence}; we conclude that the positive $G$-global model structure on $\cat{$\bm G$-UCom}$ exists and that it is combinatorial with generating cofibrations $\textbf{P}I$ and generating acyclic cofibrations $\textbf{P}J$, where $I$ and $J$ are sets of generating cofibrations and generating acyclic cofibrations for $\cat{$\bm G$-$\bm{\mathcal I}$-SSet}$. If we take $I$ to be the standard generating cofibrations, then a calculation analogous to the verification of the Strong Commutative Monoid Axiom identifies $\textbf{P}I$ with the proposed generating cofibrations.

This model structure is right proper by Lemma~\ref{lemma:transferred-properties}-$(\ref{item:tpr-proper})$ because it is transferred from a right proper model structure, and it is also left proper by an immediate application of Theorem~\ref{thm:cmon-left-proper}, all of whose assumptions have been verified above. Moreover, filtered colimits in it are homotopical, as the forgetful functor to $\cat{$\bm G$-$\bm{\mathcal I}$-SSet}$ creates both weak equivalences as well as filtered colimits.

Moreover, $\CMon(\cat{$\bm G$-$\bm{\mathcal I}$-SSet})$ has an evident simplicial enrichment, and it is cotensored over $\cat{SSet}$ with cotensors formed in $\cat{$\bm G$-$\bm{\mathcal I}$-SSet}$. For fixed $K$, $(\blank)^K$ preserves limits and sufficiently highly filtered colimits (since these are created by the forgetful functor), so it admits a left adjoint by the Special Adjoint Functor Theorem, i.e.~$\CMon(\cat{$\bm G$-$\bm{\mathcal I}$-SSet})$ is also tensored over $\cat{SSet}$. Since the forgetful functor has an obvious enrichment with respect to which it preserves cotensors, $\textbf P\dashv \forget$ becomes a simplicial adjunction. Thus, we conclude from \hbox{Lemma~\ref{lemma:transferred-properties}-$(\ref{item:tpr-simplicial})$} that the positive $G$-global model structure on $\cat{$\bm G$-UCom}$ is simplicial.

Finally, the forgetful functor is right Quillen by design, so that $\textbf{P}\dashv\forget$ is a simplicial Quillen adjunction.
\end{proof}
\end{thm}

\subsection{The model structure on \texorpdfstring{$\bm G$}{G}-parsummable simplicial sets} We now want to lift our comparison between $\cat{$\bm{E\mathcal M}$-$\bm G$-SSet}^\tau$ and $\cat{$\bm G$-$\bm{\mathcal I}$-SSet}$ to commutative monoids, for which we introduce:

\begin{thm}\label{thm:G-ParSumSSet-model}\index{G-parsummable simplicial set@$G$-parsummable simplicial set!model structure|seeonly{positive $G$-global model structure, on $\cat{$\bm G$-ParSumSSet}$}}
\index{positive G-global model structure@positive $G$-global model structure!on G-ParSumSSet@on $\cat{$\bm G$-ParSumSSet}$|textbf}
There exist a unique model structure on $\cat{$\bm G$-ParSumSSet}$ in which a map is a weak equivalence or fibration if and only if it so in the positive $G$-global model structure on $\cat{$\bm{E\mathcal M}$-$\bm G$-SSet}^\tau$.

We call this the \emph{positive $G$-global model structure}. It is proper, simplicial, and combinatorial with generating cofibrations (say, as maps in $\cat{$\bm{E\mathcal M}$-$\bm G$-SSet}^\tau$)
\begin{equation*}
\coprod_{n\ge 0}\big(E\Inj(\bm{n}\times A,\omega)\times G^n\times(\del\Delta^m\hookrightarrow\Delta^m)^{\times n}\big)/(\Sigma_n\wr H),
\end{equation*}
where $m\ge 0$, $H$ runs through all finite groups, and $A$ is a finite faithful non-empty $H$-set; moreover, filtered colimits in this model category are homotopical. Finally,
\begin{equation*}
\textbf{\textup P}\colon(\cat{$\bm{E\mathcal M}$-$\bm G$-SSet}^\tau)_{\textup{positive $G$-global}}\rightleftarrows\cat{$\bm G$-ParSumSSet}_{\textup{positive $G$-global}} :\!\forget
\end{equation*}
is a simplicial Quillen adjunction.
\begin{proof}
We again verify the assumptions of Theorems~\ref{thm:cmon-existence} and~\ref{thm:cmon-left-proper}.

It suffices to check the Pushout Product Axiom on generating (acyclic) cofibrations. As these are given by applying the strong symmetric monoidal left adjoint $\ev_\omega$ to the generating (acyclic) cofibrations of $\cat{$\bm G$-$\bm{\mathcal I}$-SSet}$, this is a formal consequence of the Pushout Product Axiom for $\cat{$\bm G$-$\bm{\mathcal I}$-SSet}$. Analogously, the Strong Commutative Monoid Axiom follows from the one for $\cat{$\bm G$-$\bm{\mathcal I}$-SSet}$.

The Unit Axiom is again automatic as $\boxtimes$ is homotopical (Theorem~\ref{thm:boxtimes-em-homotopical}). Next, we observe that $X\boxtimes\blank$ preserves injective cofibrations, which immediately implies that pushouts along $X\boxtimes i$ are homotopy pushouts for any positive $G$-global cofibration $i$. Thus, the Monoid Axiom holds by the same argument as before.

We conclude as before that the (transferred) model structure exists, is proper, and that filtered colimits in it are homotopical. Moreover, it is combinatorial with generating cofibrations $\textbf PI$ for any set $I$ of generating cofibrations of $\cat{$\bm{E\mathcal M}$-$\bm G$-SSet}^\tau$. Again using that $\ev_\omega$ is cocontinuous and strong symmetric monoidal, the calculation from the previous theorem yields the above description of $\textbf PI$.

Finally, one argues as in the previous theorem to establish that the model structure is simplicial and that $\cat{P}\dashv\forget$ is a simplicial Quillen adjunction.
\end{proof}
\end{thm}

In order to compare this model category to the one from the previous section, we first observe that there is again by abstract nonsense a preferred way to make the right adjoint $(\blank)_\bullet$ into a simplicial lax symmetric monoidal functor, inducing an enriched adjunction $\ev_\omega\colon\cat{$\bm G$-UCom}\rightleftarrows\cat{$\bm G$-ParSumSSet}:(\blank)_\bullet$. By Lemma~\ref{lemma:evaluation-support-adjunction-E} the right adjoint is fully faithful with essential image precisely those $G$-ultracommutative monoids whose underlying $G$-$\mathcal I$-simplicial sets are flat.

\begin{warn}
As opposed to the situation for $I$-simplicial sets, it is not clear whether the above makes $(\blank)_\bullet$ into a \emph{strong} symmetric monoidal functor, i.e.~whether the canonical maps $X_\bullet\boxtimes Y_\bullet\to (X\boxtimes Y)_\bullet$ are isomorphisms. As $(\blank)_\bullet$ is fully faithful with essential image precisely the flat $\mathcal I$-simplicial sets, this is equivalent to the question whether the box product \emph{of $\mathcal I$-simplicial sets} preserves flatness.
\end{warn}

\begin{cor}\label{cor:UCom-vs-ParSumSSet}
The simplicial adjunction
\begin{equation*}
\ev_\omega\colon\cat{$\bm G$-UCom}\rightleftarrows\cat{$\bm G$-ParSumSSet} :\!(\blank)_\bullet
\end{equation*}
is a Quillen equivalence.
\begin{proof}
While this follows from the corresponding result for the underlying adjunction $\cat{$\bm G$-$\bm{\mathcal I}$-SSet}\rightleftarrows\cat{$\bm{E\mathcal M}$-$\bm G$-SSet}^\tau$ by \cite[Theorem~4.19]{white-cmon}, there is an easy direct argument available: by construction, $(\blank)_\bullet$ preserves weak equivalences and fibrations, and $\ev_\omega$ creates weak equivalences, so it suffices that the counit $X_\bullet(\omega)\to X$ be a weak equivalence for every (fibrant) $X$. But this is given by the counit of the underlying adjunction, so it is in fact an isomorphism for all $X$.
\end{proof}
\end{cor}

For later use we record:

\begin{lemma}\label{lemma:ppo-injective-cofibration}
Let $i\colon A\to B$ and $j\colon C\to D$ be injective cofibrations of tame $E\mathcal M$-simplicial sets. Then $i\ppo j$ is an injective cofibration.
\begin{proof}
It is clear that $B\boxtimes j\colon B\boxtimes C\to B\boxtimes D$ and $i\boxtimes D$ are injective.

Now assume $(a,d)\in (A\boxtimes D)_n$ and $(b,c)\in (B\boxtimes C)_n$ have the same image in $(B\boxtimes D)_n$. Then by definition $b=i(a)$ and $d=j(c)$; so to show that $(b,c)$ and $(a,d)$ represent the same element in the pushout $(A\boxtimes D)\amalg_{A\boxtimes C}(B\boxtimes C)$ it suffices that $(a,c)\in (A\boxtimes C)_n$, i.e.~$\supp_k(a)\cap\supp_k(c)=\varnothing$ for all $0\le k\le n$. But indeed, as $i$ is injective and $E\mathcal M$-equivariant, $\supp_k(a)=\supp_k(i(a))=\supp_k(b)$, which is disjoint from $\supp_k(c)$ as $(b,c)\in(B\boxtimes C)_n$.
\end{proof}
\end{lemma}

\begin{cor}\label{cor:iterated-ppo-injective}
Let $i\colon A\to B$ be an injective cofibration of tame $E\mathcal M$-simplicial sets, and let $n\ge 0$. Then $i^{\ppo n}$ and $i^{\ppo n}/\Sigma_n$ are injective cofibrations.
\begin{proof}
Identifying $i^{\ppo n}$ with the iterated pushout product, the first claim follows by applying the previous lemma inductively. The second statement follows as quotients by group actions in $\cat{Set}$ preserve injections.
\end{proof}
\end{cor}

\subsection{\texorpdfstring{$\bm I$}{I} vs.~\texorpdfstring{$\bm{\mathcal I}$}{I} and \texorpdfstring{$\bm{\mathcal M}$}{M} vs.~\texorpdfstring{$\bm{E\mathcal M}$}{EM}}
\index{G-globally coherently commutative monoid@$G$-globally coherently commutative monoid!commutative monoid in G-I-SSet@commutative monoid in $\cat{$\bm G$-$\bm I$-SSet}$|(}
We will now establish analogues of the above model structures for commutative monoids in $\cat{$\bm G$-$\bm I$-SSet}$ and $\cat{$\bm{\mathcal M}$-$\bm G$-SSet}^\tau$, and show that these are equivalent in the evident way to the models considered so far. This will in particular allow us in the next subsection to compare our notion of ultra-commutative monoids to the one introduced by Schwede.

We begin with the following observation that is proven analogously to Proposition~\ref{prop:ev-omega-strong-sym-mon}:

\begin{lemma}
There is a unique way to make $\mathcal I\times_I\blank\colon\cat{$\bm I$-SSet}\to\cat{$\bm{\mathcal I}$-SSet}$ into a simplicial strong symmetric monoidal functor such that the structure isomorphism
\begin{equation*}
\big(\mathcal I\times_I I(A,\blank)\big)\boxtimes\big(\mathcal I\times_I I(B,\blank)\big)\to \mathcal I\times_I\big(I(A,\blank)\boxtimes I(B,\blank)\big)
\end{equation*}
is the identity for all $A,B\in I$.\qed
\end{lemma}

We can now deduce the following $G$-global version of \cite[Theorem~1.1]{I-vs-M-1-cat}:

\begin{thm}\label{thm:box-i-homotopical}
The box product on $\cat{$\bm G$-$\bm I$-SSet}$ is homotopical in $G$-global weak equivalences.
\begin{proof}
By Theorem~\ref{thm:G-global-we-I-characterization} the functor $\mathcal I\times_I\blank$ detects $G$-global weak equivalences (without any need to derive). The claim therefore follows from the previous lemma together with Corollary~\ref{cor:boxtimes-homotopical-script-I}.
\end{proof}
\end{thm}

\begin{thm}\label{thm:CMon-G-I}\index{G-globally coherently commutative monoid@$G$-globally coherently commutative monoid!commutative monoid in G-I-SSet@commutative monoid in $\cat{$\bm G$-$\bm I$-SSet}$!model structure|seeonly{positive $G$-global model structure, on $\CMon(\cat{$\bm G$-$\bm I$-SSet})$}}
\index{positive G-global model structure@positive $G$-global model structure!on CMon(G-I-SSet)@on $\CMon(\cat{$\bm G$-$\bm I$-SSet})$|textbf}
There is a unique model structure on $\CMon(\cat{$\bm G$-$\bm I$-SSet})$ in which a map is a weak equivalence or fibration if and only if it is so in the positive $G$-global model structure on $\cat{$\bm G$-$\bm I$-SSet}$. We call this the \emph{positive $G$-global model structure} again. It is left proper, simplicial, combinatorial, and filtered colimits in it are homotopical. Finally, the simplicial adjunction
\begin{equation}\label{eq:script-I-vs-I-CMon}
\mathcal I\times_I\blank\colon\CMon(\cat{$\bm G$-$\bm I$-SSet})\rightleftarrows\CMon(\cat{$\bm G$-$\bm{\mathcal I}$-SSet})=\cat{$\bm G$-UCom} :\!\forget
\end{equation}
is a Quillen equivalence in which both adjoints are fully homotopical.
\end{thm}

Here we again used that $\forget$ acquires a simplicial lax symmetric monoidal structure from the simplicial strong symmetric monoidal structure on $\mathcal I\times_I\blank$.

\begin{proof}
We will once again verify the assumptions of Theorem~\ref{thm:cmon-existence}: the Pushout Product Axiom for cofibrations is verified analogously to Theorem~\ref{thm:G-UCom}, and for acyclic cofibrations it then follows from the one for $\cat{$\bm G$-$\bm{\mathcal I}$-SSet}$ as $\mathcal I\times_I\blank$ is left Quillen, symmetric monoidal, and reflects weak equivalences. Analogously, one verifies the Strong Commutative Monoid Axiom and the Monoid Axiom. Finally, the Unit Axiom is again automatic by Theorem~\ref{thm:box-i-homotopical}.

We therefore conclude as before that the model structure exists, is simplicial, combinatorial, and that filtered colimits in it are homotopical.

As the forgetful functors create fibrations and weak equivalences, we conclude from Theorem~\ref{thm:positive-global-I-model-structure} that $(\ref{eq:script-I-vs-I-CMon})$ is a Quillen adjunction. As both adjoints in $\mathcal I\times_I\blank\colon\cat{$\bm G$-$\bm I$-SSet}\rightleftarrows\cat{$\bm G$-$\bm{\mathcal I}$-SSet}:\!\forget$ are homotopical (see Theorem~\ref{thm:G-global-we-I-characterization} for the non-trivial case), we conclude that the ordinary unit and counit already represent the derived unit and counit and that they are weak equivalences. From this one easily deduces that both adjoints in $(\ref{eq:script-I-vs-I-CMon})$ are fully homotopical and that unit and counit are weak equivalences, so $(\ref{eq:script-I-vs-I-CMon})$ is in particular a Quillen equivalence.

Finally, the left Quillen functor $\mathcal I\times_I\blank\colon\CMon(\cat{$\bm G$-$\bm I$-SSet})\to\cat{$\bm G$-UCom}$ creates weak equivalences, so left properness of $\CMon(\cat{$\bm G$-$\bm I$-SSet})$ follows from left properness of $\cat{$\bm G$-UCom}$ as before (see Lemma~\ref{lemma:U-pushout-preserve-reflect}).
\end{proof}
\index{G-globally coherently commutative monoid@$G$-globally coherently commutative monoid!commutative monoid in G-I-SSet@commutative monoid in $\cat{$\bm G$-$\bm I$-SSet}$|)}

\index{G-globally coherently commutative monoid@$G$-globally coherently commutative monoid!commutative monoid in M-G-SSettau@commutative monoid in $\cat{$\bm{\mathcal M}$-$\bm G$-SSet}^\tau$|(}
Finally, let us turn to the box product on $\cat{$\bm{\mathcal M}$-SSet}^\tau$:

\begin{lemma}
For any $X,Y\in\cat{$\bm{\mathcal M}$-SSet}^\tau$ the composition
\begin{equation*}
E\mathcal M\times_{\mathcal M}(X\boxtimes Y)\hookrightarrow
E\mathcal M\times_{\mathcal M}(X\times Y)\xrightarrow{\pr_1,\pr_2} (E\mathcal M\times_{\mathcal M}X)\times (E\mathcal M\times_{\mathcal M}Y)
\end{equation*}
factors through an isomorphism $E\mathcal M\times_{\mathcal M}(X\boxtimes Y)\to (E\mathcal M\times_{\mathcal M}X)\boxtimes(E\mathcal M\times_{\mathcal M}Y)$. Together with the unique maps $E\mathcal M\times_{\mathcal M}*\to*$, the inverses of these isomorphisms make $E\mathcal M\times_{\mathcal M}\blank$ into a simplicial strong symmetric monoidal functor.
\begin{proof}
Let $X,Y$ be tame $E\mathcal M$-simplicial sets. Then Lemma~\ref{lemma:supp-vs-supp-k} implies that the box product of their underlying $\mathcal M$-simplicial sets is a subcomplex of their box product as $E\mathcal M$-simplicial sets, so the inclusions define a simplicial natural transformation $\forget(\blank\boxtimes\blank)\Rightarrow\forget(\blank)\boxtimes\forget(\blank)$. As the (simplicial) symmetric monoidal structure isomorphisms of both $\cat{$\bm{\mathcal M}$-SSet}^\tau$ and $\cat{$\bm{E\mathcal M}$-SSet}^\tau$ are defined in terms of the one for the cartesian monoidal structure, it is clear that this makes $\forget$ into a simplicial lax symmetric monoidal functor. One easily checks from the definition that the resulting simplicial oplax structure on $E\mathcal M\times_{\mathcal M}\blank$ consists of precisely the above maps $E\mathcal M\times_{\mathcal M}(X\boxtimes Y)\to (E\mathcal M\times_{\mathcal M}X)\boxtimes(E\mathcal M\times_{\mathcal M}Y)$, so it only remains to show that these are isomorphisms.

But indeed, the box product of $\mathcal M$-simplicial sets is cocontinuous in each variable by \cite[Corollary~2.17]{I-vs-M-1-cat} and it obviously preserves tensors in each variable. The same holds for the box product of $E\mathcal M$-simplicial sets by Corollary~\ref{cor:box-EM-cocontinuous}. As $E\mathcal M\times_{\mathcal M}\blank$ is a simplicial left adjoint, we therefore conclude that the set of all pairs $(X,Y)$ such that the above comparison map is an isomorphism is closed under tensoring and small colimits in each variable. By Theorem~\ref{thm:structure-tame-M-G-SSet} we are therefore reduced to verifying the claim for $X=\Inj(A,\omega)$ and $Y=\Inj(B,\omega)$ for some finite sets $A,B$. This is then a straightforward calculation using Corollary~\ref{cor:E-Inj-corepr} together with Examples~\ref{ex:boxtimes-schwede} and~\ref{ex:boxtimes-std}.
\end{proof}
\end{lemma}

Together with Theorem~\ref{thm:tame-M-sset-vs-EM-sset} we conclude as before:

\begin{cor}
The box product on $\cat{$\bm{\mathcal M}$-$\bm G$-SSet}^\tau$ is homotopical with respect to the $G$-global weak equivalences.\qed
\end{cor}

\begin{thm}\index{G-globally coherently commutative monoid@$G$-globally coherently commutative monoid!commutative monoid in M-G-SSettau@commutative monoid in $\cat{$\bm{\mathcal M}$-$\bm G$-SSet}^\tau$!model structure|seeonly{positive $G$-global model structure, on $\CMon(\cat{$\bm{\mathcal M}$-$\bm G$-SSet}^\tau)$}}
\index{positive G-global model structure@positive $G$-global model structure!on CMon(M-G-SSettau)@on $\CMon(\cat{$\bm{\mathcal M}$-$\bm G$-SSet}^\tau)$|textbf}
There is a unique model structure on $\CMon(\cat{$\bm{\mathcal M}$-$\bm G$-SSet}^\tau)$ in which a map is a weak equivalence or fibration if and only if it is so in the positive $G$-global model structure on $\cat{$\bm{\mathcal M}$-$\bm G$-SSet}^\tau$. We call this the \emph{positive $G$-global model structure}. It is left proper, simplicial, combinatorial, and filtered colimits in it are homotopical. Finally, the simplicial adjunctions
\begin{equation}\label{eq:CMon-M-G-vs-ParSumSSet}
E\mathcal M\times_{\mathcal M}\blank\colon\CMon(\cat{$\bm{\mathcal M}$-$\bm G$-SSet}^\tau)\rightleftarrows\CMon(\cat{$\bm{E\mathcal M}$-$\bm G$-SSet}^\tau) :\!\forget
\end{equation}
and
\begin{equation}\label{eq:CMon-M-G-vs-G-I}
\ev_\omega\colon\CMon(\cat{$\bm G$-$\bm{\mathcal I}$-SSet})\rightleftarrows\CMon(\cat{$\bm{\mathcal M}$-$\bm G$-SSet}^\tau) :\!(\blank)_\bullet
\end{equation}
are Quillen equivalences.
\begin{proof}
Let us verify the assumptions of Theorem~\ref{thm:cmon-existence} again.

For the Pushout Product Axiom, it suffices to check this on generating cofibrations and generating acyclic cofibrations. As the positive $G$-global model structure on $\cat{$\bm{\mathcal M}$-$\bm{G}$-SSet}^\tau$ is transferred from $\cat{$\bm{G}$-$\bm{I}$-SSet}$ and since $\ev_\omega$ is symmetric monoidal (Remark~\ref{rk:ev-symmetric-monoidal}), it is therefore implied by the Pushout Product Axiom for the latter. Analogously, one deduces the Strong Commutative Monoid Axiom.

The Unit Axiom is immediate from the previous corollary, and the Monoid Axiom follows from the one for $\cat{$\bm{E\mathcal M}$-$\bm G$-SSet}^\tau$ (Theorem~\ref{thm:G-ParSumSSet-model}) as $E\mathcal M\times_{\mathcal M}\blank$ is left Quillen, monoidal, and creates weak equivalences.

We therefore conclude as before that the model structure exists and that it is combinatorial and simplicial (with respect to the obvious enrichment).

The forgetful functor $\CMon(\cat{$\bm{E\mathcal M}$-$\bm G$-SSet}^\tau)\mskip0mu minus .2mu\to\mskip0mu minus .2mu\CMon(\cat{$\bm{\mathcal M}$-$\bm G$-SSet}^\tau)$ preserves weak equivalences and fibrations as they are created in the underlying categories, so it is in particular right Quillen. Similarly, $E\mathcal M\times_{\mathcal M}\blank$ is fully homotopical. The ordinary unit and counit of $E\mathcal M\times_{\mathcal M}\blank\colon\cat{$\bm{\mathcal M}$-$\bm G$-SSet}^\tau\rightleftarrows\cat{$\bm{E\mathcal M}$-$\bm G$-SSet}^\tau:\forget$ are $G$-global weak equivalences by Theorem~\ref{thm:tame-M-sset-vs-EM-sset}, so we conclude as in the proof of Theorem~\ref{thm:CMon-G-I} that $(\ref{eq:CMon-M-G-vs-ParSumSSet})$ is a Quillen equivalence.

By definition, $(\blank)_\bullet\colon\CMon(\cat{$\bm{\mathcal M}$-$\bm G$-SSet}^\tau)\mskip0mu minus .3mu\to\mskip0mu minus .3mu\CMon(\cat{$\bm G$-$\bm{I}$-SSet})$ preserves both weak equivalences and fibrations, so $(\ref{eq:CMon-M-G-vs-G-I})$ is a Quillen adjunction. As also $\ev_\omega$ is homotopical, Theorem~\ref{thm:M-G-tau-model-structure} implies by the same arguments as before that $(\ref{eq:CMon-M-G-vs-G-I})$ is a Quillen equivalence.

Finally, left properness of the model structure on $\CMon(\cat{$\bm{\mathcal M}$-$\bm G$-SSet}^\tau)$ again follows from left properness of the one on $\CMon(\cat{$\bm{E\mathcal M}$-$\bm G$-SSet}^\tau)$ as $E\mathcal M\times_{\mathcal M}\blank$ is left Quillen and creates weak equivalences.
\end{proof}
\end{thm}
\index{G-globally coherently commutative monoid@$G$-globally coherently commutative monoid!commutative monoid in M-G-SSettau@commutative monoid in $\cat{$\bm{\mathcal M}$-$\bm G$-SSet}^\tau$|)}

\subsection{Ultra-commutative monoids in $\cat{$\bm I$-SSet}$ vs.~in $\cat{$\bm L$-Top}$}\label{subsec:uc-vs-uc}
The orthogonal spaces of Definition~\ref{defi:orthogonal-space} admit a Day convolution product similarly to our models of $G$-global homotopy theory, i.e.~there is an essentially unique topologically enriched functor
\begin{equation*}
\blank\boxtimes\blank\colon\cat{$\bm L$-Top}\times\cat{$\bm L$-Top}\to\cat{$\bm L$-Top}
\end{equation*}
that preserves tensors and small colimits in each variable and satisfies $L(V,\blank)\boxtimes L(W,\blank)=L(V\oplus W,\blank)$ with the obvious functoriality in each variable. There is then again a unique way to make this into the tensor product of a symmetric monoidal structure on $\cat{$\bm L$-Top}$ such that the structure isomorphisms on corepresentables are induced by the structure isomorphisms of the cartesian symmetric monoidal structure on $\cat{Vect}_{\mathbb R}$. Schwede introduced the term \emph{ultra-commutative monoid}\index{ultra-commutative monoid!in L-Top@in $\cat{$\bm L$-Top}$|textbf} for a commutative monoid in $\cat{$\bm L$-Top}$, and he proved as \cite[Theorem~2.1.15-(i)]{schwede-book}:

\begin{thm}\index{ultra-commutative monoid!in L-Top@in $\cat{$\bm L$-Top}$!model structure|seeonly{positive global model structure, on $\CMon(\cat{$\bm L$-Top})$}}\index{positive global model structure!on CMonL-Top@on $\CMon(\cat{$\bm L$-Top})$|textbf}
There is a unique model structure on $\CMon(\cat{$\bm L$-Top})$ in which a map is a weak equivalence or fibration if and only if it is so in the positive global model structure on $\cat{$\bm L$-Top}$. This model structure is proper, topological, and cofibrantly generated.\qed
\end{thm}
\nobreak
In this subsection, we will lift our comparison between $\cat{$\bm I$-SSet}$ and $\cat{$\bm L$-Top}$ to the level of commutative monoids.
\goodbreak

To this end, we first observe that also $\cat{$\bm I$-Top}$ admits a box product induced by the symmetric monoidal structure on $I$ given by disjoint union. As before one then makes $|\blank|\colon\cat{$\bm I$-SSet}\to\cat{$\bm I$-Top}$ and $L\times_I\blank\colon\cat{$\bm I$-Top}\to\cat{$\bm L$-Top}$ into strong symmetric monoidal functors. We conclude that their right adjoints are lax symmetric monoidal, so that we get induced adjunctions of categories of commutative monoids. We can now state our comparison result:

\begin{thm}\label{thm:uc-vs-uc}\index{ultra-commutative monoid!these models are equivalent|textbf}
The adjunction
\begin{equation}\label{eq:uc-vs-uc-adj}
L\times_I|\blank|\colon\CMon(\cat{$\bm I$-SSet})\rightleftarrows\CMon(\cat{$\bm L$-Top}):\!\Sing\circ\forget
\end{equation}
is a Quillen adjunction, and the induced adjunction of associated quasi-categories is a right Bousfield localization with respect to the $\mathcal Fin$-global weak equivalences.
\end{thm}

In order to deduce Theorem~\ref{thm:uc-vs-uc} from the comparison between $\cat{$\bm I$-SSet}$ and $\cat{$\bm L$-Top}$ provided in Section~\ref{sec:global-vs-g-global}, we will need the following cofibrancy property:

\begin{lemma}
Let $M$ be cofibrant in $\CMon(\cat{$\bm G$-$\bm I$-SSet})$. Then $M$ is cofibrant in the \emph{$G$-global model structure} on $\cat{$\bm G$-$\bm I$-SSet}$.
\begin{proof}
We will show that any cofibration $X\to Y$ in $\CMon(\cat{$\bm G$-$\bm I$-SSet})$ with $X$ cofibrant in the $G$-global model structure on $\cat{$\bm G$-$\bm I$-SSet}$ is also a $G$-global cofibration in $\cat{$\bm G$-$\bm I$-SSet}$. Applying this to $0\to X$ then yields the lemma.

While we cannot literally apply \cite[Corollary~3.6]{white-cmon} to the positive $G$-global model structure (as the unit is not cofibrant in the positive $G$-global model structure) nor the $G$-global model structure (because the acyclicity part of the Strong Commutative Monoid Axiom fails), the same argument works in our situation, as also observed for example in \cite[proof of~Theorem~2.1.15-(ii)]{schwede-book} for $\cat{$\bm L$-Top}$:

Namely, we first note that the $G$-global model structure still satisfies the Pushout Product Axiom. On the other hand, it is clear by direct inspection that the forgetful functor sends generating cofibrations of $\CMon(\cat{$\bm G$-$\bm I$-SSet})$ to $G$-global cofibrations in $\cat{$\bm G$-$\bm I$-SSet}$. If now $f\colon X\to Y$ is a pushout of a generating cofibration, then \cite[Proposition~B.2]{white-cmon} (say, applied to the model structure on $\cat{$\bm G$-$\bm I$-SSet}$ in which all maps are cofibrations, but only the isomorphisms are weak equivalences) shows that $f$ can be written as a transfinite composition of pushouts of maps of the form $X\boxtimes i_n^{\ppo n}/\Sigma_n$ for generating cofibrations $i_n$. By the Strong Commutative Monoid Axiom, each $i_n^{\ppo n}/\Sigma_n$ is a positive $G$-global cofibration, and if $X$ is $G$-globally cofibrant, then $X\boxtimes i_n^{\ppo n}/\Sigma_n$ is a $G$-global cofibration by the Pushout Product Axiom for the $G$-global model structure.

Inductively we see that if $X$ is $G$-globally cofibrant and $X\to Y$ is any relative cell complex in the generating cofibrations, then the underlying map in $\cat{$\bm G$-$\bm I$-SSet}$ is a $G$-global cofibration. Finally, if $f\colon X\to Z$ is a general cofibration with cofibrant source, then Quillen's Retract Argument shows that $f$ is a retract of some relative cell complex $g\colon X\to Y$ (note that the sources agree!). As the underlying map of $g$ is a $G$-global cofibration, so is the underlying map of $f$ as desired.
\end{proof}
\end{lemma}

\begin{proof}[Proof of Theorem~\ref{thm:uc-vs-uc}]
Analogously to the arguments in Section~\ref{sec:global-vs-g-global} one shows that $L\times_I|\blank|\colon\cat{$\bm I$-SSet}\rightleftarrows\cat{$\bm L$-Top}:\Sing\circ\forget$ is a Quillen adjunction for the \emph{positive} global model structures. Thus, also $(\ref{eq:uc-vs-uc-adj})$ is a Quillen adjunction.

As weak equivalences are created in the underlying categories, we conclude from Propositions~\ref{prop:I-Top-global} and~\ref{prop:Fin-global-we} that the right adjoint is homotopical and inverts precisely the $\mathcal Fin$-global weak equivalences. It therefore only remains to show that the unit is a weak equivalence on any cofibrant object $X\in\CMon(\cat{$\bm I$-SSet})$. But by the previous lemma, $X$ is cofibrant in the global model structure on $\cat{$\bm I$-SSet}$, so the claim follows from Proposition~\ref{prop:I-Top-global} and Theorem~\ref{thm:L-vs-I}.
\end{proof}

\begin{rk}\index{global E-infinity-algebra@global $E_\infty$-algebra}
Recently, Barrero \cite{barrero} introduced \emph{global $E_\infty$-operads} in $\cat{$\bm L$-Top}$ and proved that the homotopy theory of the corresponding algebras (with respect to the box product) is equivalent to the global homotopy theory of ultra-commutative monoids (with respect to all compact Lie groups, and therefore in particular with respect to finite groups). By the results of this chapter, this operadic approach is then equivalent to all of our models of `globally coherently commutative monoids.'
\end{rk}

\section[$G$-global $\Gamma$-spaces]{\for{toc}{$G$}\except{toc}{\texorpdfstring{$\bm G$}{G}}-global \for{toc}{$\Gamma$}\except{toc}{\texorpdfstring{$\bm\Gamma$}{Gamma}}-spaces}
\index{G-globally coherently commutative monoid@$G$-globally coherently commutative monoid!G-global Gamma-space@$G$-global $\Gamma$-space|seeonly{$G$-global $\Gamma$-space}}
In this section we study a $G$-global version of Segal's theory of \emph{$\Gamma$-spaces}. While we focus on them as models of `$G$-globally coherently commutative monoids' in this section, we will discuss their relation to $G$-global spectra in Section~\ref{sec:delooping-group-completion}.

\subsection[A reminder on $\Gamma$-spaces]{A reminder on $\bm\Gamma$-spaces}
Let us briefly recall the classical non-equivariant and equivariant theory of $\Gamma$-spaces due to Segal \cite{segal-gamma} and Shimakawa \cite{shimakawa}, respectively.

\begin{defi}\nomenclature[aGamma]{$\Gamma$}{category of finite pointed sets and based maps}
We write $\Gamma$ for the category of finite pointed sets and base point preserving maps. For any $n\ge 0$, let $n^+\mathrel{:=}\{0,\dots,n\}$ with base point $0$.

A \emph{$\Gamma$-space}\index{coherently commutative monoid!Gamma-space@$\Gamma$-space|seeonly{$\Gamma$-space}}\index{Gamma-space@$\Gamma$-space|textbf} is a functor $X\colon\Gamma\to\cat{SSet}$ such that $X(0^+)$ is terminal. We write $\cat{$\bm\Gamma$-SSet}_*$ for the full subcategory of $\Fun(\Gamma,\cat{SSet})$ spanned by the $\Gamma$-spaces.
\end{defi}

Segal \cite[Definition~1.2]{segal-gamma} originally considered contravariant functors from a category equivalent to the opposite of the above category $\Gamma$ (into topological spaces), and he reserved the term `$\Gamma$-space' for functors which are special in the sense of Definition~\ref{defi:gamma-special} below.

\begin{rk}\index{Gamma-space@$\Gamma$-space!as enriched functor|textbf}
The category $\cat{SSet}_*$ of pointed simplicial sets has a zero object, so it admits a unique $\cat{Set}_*$-enrichment; explicitly the base point of $\Hom(X,Y)$ is taken to be the constant map. Similarly, there is a unique $\cat{Set}_*$-enrichment of $\Gamma$.

As observed for example in \cite[Lemma~1.17]{may-merling-osorno} (for topological spaces), any $\cat{Set}_*$-enriched functor $X\colon\Gamma\to\cat{SSet}_*$ satisfies $X(0^+)=*$, so its underlying ordinary functor $\Gamma\to\cat{SSet}$ is a $\Gamma$-space; conversely, every $\Gamma$-space $X$ factors uniquely through the category of based simplicial sets by taking the image of $X(0^+)$ as the base point of $X(S_+)$ for any $S_+\in\Gamma$, and the induced functor $\Gamma\to\cat{SSet}_*$ is $\cat{Set}_*$-enriched. Altogether we see that we can equivalently think of $\cat{$\bm\Gamma$-SSet}_*$ as the category of $\cat{Set}_*$-enriched functors $\Gamma\to\cat{SSet}_*$, with the isomormophism given by the forgetting base points and enrichment.
\end{rk}

\begin{defi}\label{defi:gamma-special}\index{Gamma-space@$\Gamma$-space!special|textbf}\index{special|seeonly{$\Gamma$-space, special}}
Let $S$ be a finite set and let $s\in S$. Then we write $p_s\colon S_+\to 1^+$ for the map in $\Gamma$ with $p_s(s)=1$ and $p_s(t)=0$ otherwise.

Now let $X$ be a $\Gamma$-space. The \emph{Segal map}\index{Segal map|textbf}
\begin{equation*}
\rho\colon X(S_+)\to\prod_{s\in S} X(1^+)
\end{equation*}
\nomenclature[arho2]{$\rho$}{(generalized) Segal map}%
\nomenclature[aps]{$p_s$}{map $S_+\to1^+$ in $\Gamma$ projecting away from $s\in S$; used to define the Segal map}%
is the map given on the factor corresponding to $s\in S$ by $X(p_s)$. We call $X$ \emph{special} if the Segal map is a weak homotopy equivalence for every finite set $S$.
\end{defi}

Intuitively, if $X$ is a special $\Gamma$-space, then we want to think of $X(1^+)$ as its `underlying space' with the remaining structure encoding `(higher) additions.' Explicitly, consider for each $n\ge0$ the map $\mu\colon n^+\to 1^+$ sending every non-base point to $1$. Then we have a zig-zag
\begin{equation*}
\begin{tikzcd}[cramped]
X(1^+)^n & \arrow[l, "\sim", "\rho"'] X(n^+) \arrow[r, "X(\mu)"] & X(1^+)
\end{tikzcd}
\end{equation*}
which we think of as $n$-fold addition. In particular, for $n=2$, passing to $\pi_0$ yields a map $\pi_0(X(1^+))^2\cong \pi_0(X(1^+)^2)\to \pi_0(X(1^+))$ and one can show that this equips $\pi_0(X(1^+))$ with the structure of a commutative monoid.

\begin{ex}
If $A$ is any simplicial abelian monoid, then we can define a special $\Gamma$-space $HA$ via $(HA)(S_+)= A^{\times S}$; writing elements of $HA(S_+)$ as formal sums $\sum_{s\in S}a_ss$, the structure map for $f\colon S_+\to T_+$ is given by $\sum_{s\in S}a_ss\mapsto\sum_{f(s)\not=*} a_sf(s)=\sum_{t\in T}\big(\sum_{s\in f^{-1}(t)}a_s)t$.
\end{ex}

\begin{ex}\label{ex:shimada-shimakawa}\nomenclature[aGammaC]{$\Gamma(\mathscr C)$}{$\Gamma$-category built from the small symmetric monoidal category $\mathscr C$}
Building on a construction of Segal \cite[discussion after Corollary~2.2]{segal-gamma}, Shimada and Shimakawa \cite[Definition~2.1]{shimada-shimakawa} (and earlier May \cite[Construction~10]{may-unique} in an important special case) showed how to functorially associate to a small symmetric monoidal category $\mathscr C$ a `categorically special $\Gamma$-category,'\index{special!categorically special} i.e.~a functor $\Gamma(\mathscr C)\colon\Gamma\to\cat{Cat}$ with $\Gamma(\mathscr C)(0^+)=*$ and such that the obvious analogues of the Segal maps $\Gamma(\mathscr C)(S_+)\to\prod_{s\in S}\Gamma(\mathscr C)$ are equivalences of categories \cite[Lemma~2.2]{shimada-shimakawa}. In particular, applying the nerve levelwise yields a special $\Gamma$-space in the above sense.

Somewhat more concretely, an object of $\Gamma(\mathscr C)(S_+)$ consists of an object $X_A\in\mathscr C$ for each $A\subset S$ together with isomorphisms $\alpha_{A,B}\colon X_{A\cup B}\to X_A\otimes X_B$ for all disjoint $A,B\subset S$, such that $X_\varnothing=\textbf1$ is the tensor unit, and the isomorphisms $\alpha$ are suitably compatible with the associativity, unitality, and symmetry isomorphisms of $\mathscr C$. A morphism $(X_\bullet,\alpha_{\bullet,\bullet})\to (Y_\bullet,\beta_{\bullet,\bullet})$ consists of a map $f_A\colon X_A\to Y_A$ for any $A\subset S$ suitably compatible with $\alpha$, $\beta$ and such that $f_\varnothing$ is the identity of $\bm1$.

The structure maps of $\Gamma(\mathscr C)$ are given by restriction along preimages, and $\Gamma$ becomes a functor in strictly unital strong symmetric monoidal functors (i.e.~strong symmetric functors such that the unit isomorphism is the identity) essentially by pushforward, see~\cite[discussion after Definition~2.4]{shimada-shimakawa} for a precise definition.

On the level of underlying categories, $\Gamma(\mathscr C)$ recovers the original symmetric monoidal category $\mathscr C$; more precisely, there is a natural isomorphism $\iota\colon\mathscr C\to\Gamma(\mathscr C)(1^+)$, sending $X\in\mathscr C$ to the unique object $(X_\bullet,\alpha_{\bullet,\bullet})$ with $X_{\{1\}}=X$ and sending $f\colon X\to Y$ to the unique morphism $f_\bullet$ with $f_{\{1\}}=f$.
\end{ex}

\begin{rk}
Let us call a morphism $f\colon X\to Y$ of $\Gamma$-spaces a \emph{level weak equivalence}\index{level weak equivalence!in Gamma-SSet@in $\cat{$\bm\Gamma$-SSet}_*$|textbf} if each $f(S_+)\colon X(S_+)\to Y(S_+)$ is a weak homotopy equivalence. Bousfield and Friedlander \cite[Theorem~3.5]{bousfield-friedlander} showed that the level weak equivalences are part of a `Reedy type' model structure on $\cat{$\bm\Gamma$-SSet}_*$; however, for our purposes only the usual projective model structure on $\cat{$\bm\Gamma$-SSet}_*$ will be relevant.
\end{rk}

For the rest of this subsection, let $G$ be a finite group. Shimakawa \cite{shimakawa} provided a $G$-equivariant generalization of Segal's theory.

\index{coherently commutative monoid!G-equivariant@$G$-equivariant|seeonly{$\Gamma$-$G$-space}}
\begin{defi}\index{Gamma-G-space@$\Gamma$-$G$-space|textbf}
A \emph{$\Gamma$-$G$-space} is a functor $X\colon\Gamma\to\cat{$\bm G$-SSet}$ such that $X(0^+)$ is terminal. We write $\cat{$\bm\Gamma$-$\bm G$-SSet}_*$ for the evident category of $\Gamma$-$G$-spaces.
\end{defi}

\begin{rk}\index{Gamma-G-space@$\Gamma$-$G$-space!evaluation at finite G-sets@evaluation at finite $G$-sets}
Shimakawa originally considered so-called `$\Gamma_G$-spaces,'\index{Gamma-G-space2@$\Gamma_G$-space} which are $\cat{$\bm G$-Set}$-enriched functors from the category $\Gamma_G$ of finite based $G$-sets and not necessarily $G$-equivariant maps to the category $\cat{SSet}_G$ of $G$-simplicial sets and not necessarily equivariant maps; here the enrichment in $G$-sets is given on both sides via the conjugation action on the hom sets. However, he showed in \cite[Theorem~1]{shimakawa-simplify} that restricting along $\Gamma\hookrightarrow\Gamma_G$ provides an equivalence to the category of ordinary functors $\Gamma\to\cat{$\bm G$-SSet}$.

We can also make the quasi-inverse of this equivalence explicit: if $X$ is any functor $\Gamma\to\cat{$\bm G$-SSet}$ (for example a $\Gamma$-$G$-space), and $S$ is a finite $G$-set, then $X(S_+)$ carries two commuting $G$-actions: the exterior action via the action of $G$ on $X$ and the interior action induced by functoriality from the action on $S$, and we equip $X(S_+)$ with the diagonal of these two actions. This way, we can evaluate $X$ more generally at finite pointed $G$-sets. A not necessarily $G$-equivariant map $f\colon S_+\to T_+$ of finite based $G$-sets then induces via the original functoriality in $\Gamma$ a map $X(f)\colon X(S_+)\to X(T_+)$ (not necessarily $G$-equivariant), which provides the desired extension to $\Gamma_G$.

We will only be interested in the case where $f$ is actually $G$-equivariant, in which case $X(f)\colon X(S_+)\to X(T_+)$ is obviously also $G$-equivariant.
\end{rk}

The crucial insight in the theory of $\Gamma$-$G$-spaces is that, while it is enough to specify $\Gamma$-$G$-spaces on trivial $G$-sets, the correct notion of specialness should still take general $G$-sets into account:

\index{special!G-equivariant@$G$-equivariant|seeonly{$\Gamma$-$G$-space, special}}
\begin{defi}\label{defi:equivariant-special}\index{Gamma-G-space@$\Gamma$-$G$-space!special|textbf}
A $\Gamma$-$G$-space $X$ is called \emph{special} if the Segal map $X(S_+)\to X(1^+)^{\times S}=\prod_{s\in S}X(1^+)$ is a $G$-equivariant weak equivalence for every finite $G$-set $S$. Here we equip the left hand side with the diagonal $G$-action as before, and the right hand side carries the diagonal of the $G$-actions on $X(1^+)$ and the permutation action on the factors via the $G$-action on $S$.
\end{defi}

Note that the above is indeed $G$-equivariant as the usual Segal map is $\Sigma_S$-equivariant by virtue of the relation $p_{\sigma^{-1}(s)}=p_s\sigma$ in $\Gamma$ for all $\sigma\in\Sigma_S,s\in S$.

\begin{lemma}\label{lemma:special-via-graph-subgroups}
A $\Gamma$-$G$-space $X$ is special if and only if the Segal map $X(S_+)\to X(1^+)^{\times S}$ is a $\mathcal G_{G,\Sigma_S}$-weak equivalence for every finite set $S$.
\begin{proof}
This can be rephrased as saying that $X(S_+)\to X(1^+)^{\times S}$ be an $H$-equivariant weak equivalence for every finite subgroup $H\subset G$ and every finite $H$-set $S$. In particular, if $X$ satisfies the above condition, then it is special.

Conversely, assume $X$ is special, let $H\subset G$ be any subgroup, and let $S$ be a finite $H$-set. Then the $H$-equivariant inclusion $i\colon S_+\to G_+\smashp_HS_+, s\mapsto [1,s]$ admits an $H$-equivariant retraction $r$ with
\begin{equation*}
r[g,s]=\begin{cases}
    g.s & \text{if $g\in H$}\\
    * & \text{otherwise},
\end{cases}
\end{equation*}
and these fit into a commutative diagram
\begin{equation*}
\begin{tikzcd}
X(S_+)\arrow[r, "X(i)"]\arrow[d, "\rho"'] & X(G_+\smashp_HS_+)\arrow[r, "X(r)"]\arrow[d, "\rho"] & X(S_+)\arrow[d, "\rho"]\\
X(1^+)^{\times S} \arrow[r]  & X(1^+)^{G\times_HS}\arrow[r] & X(1^+)^{\times S}
\end{tikzcd}
\end{equation*}
of $H$-equivariant maps where lower horizontal maps are given by the inclusion and projection, respectively. The middle vertical map is even a $G$-equivariant weak equivalence as $X$ is special, hence the left hand vertical map is an $H$-equivariant weak equivalence (being a retract of an $H$-weak equivalence) as claimed.
\end{proof}
\end{lemma}

\begin{rk}
The above strong version of specialness is necessary in order for special $\Gamma$-$G$-spaces to yield the correct notion of \emph{equivariant infinite loop spaces}, i.e.~spaces that admit deloopings against all representation spheres; a particularly clear explanation of this can be found in \cite[Section~3.5]{blumberg-infinite-loop-spaces}. Moreover, the argument we give in \ref{subsubsec:wirthmueller} (for $G$-global $\Gamma$-spaces) shows how the above yields a \emph{Wirthmüller isomorphism}\index{Wirthmüller isomorphism!in Gamma-G-SSet@in $\cat{$\bm\Gamma$-$\bm G$-SSet}_*$}\index{Gamma-G-space@$\Gamma$-$G$-space!Wirthmüller isomorphism} for $\Gamma$-$G$-spaces, encoding an additional equivariant, or `twisted,' form of semiadditivity.
\end{rk}

\begin{ex}\label{ex:shimada-shimakawa-equivariant}
If $\mathscr C$ is a small symmetric monoidal category with $G$-action through strictly unital strong symmetric monoidal functors, then we can equip the $\Gamma$-category $\Gamma(\mathscr C)$ from Example~\ref{ex:shimada-shimakawa} with the induced $G$-action. However, the nerve of this is usually \emph{not} special in the above sense.

It was a crucial insight of Shimakawa \cite[discussion before Theorem~$\textup{A}'$]{shimakawa}, recently extensively used by Merling \cite{merling}, that we can solve this issue by replacing the $\Gamma$-$G$-category $\Gamma(\mathscr C)$ with $\Fun(EG,\blank)\circ\Gamma(\mathscr C)$, where $G$ acts on $EG$ from the right in the obvious way. Namely, if $H$ is any subgroup, then $\Fun(EG,\blank)^H$ just computes the \emph{categorical homotopy fixed points}, so it sends $G$-equivariant functors that are underlying equivalences of categories to ordinary equivalences of categories, see e.g.~\cite[Corollary~3.7]{merling}. If now $S$ is any finite $G$-set, then the Segal map $\Gamma(\mathscr C)(S_+)\to\Gamma(\mathscr C)(1^+)^{\times S}$ is $G$-equivariant, and it is of course still an underlying equivalence of categories. The image of this under $\nerve\Fun(EG,\blank)$ will therefore be a $G$-equivariant weak equivalence. However, this is clearly conjugate to the Segal map of $\nerve\Fun(EG,\Gamma(\mathscr C))$.
\end{ex}

\begin{rk}\label{rk:Gamma-G-monoid}
Already when $X$ is special in the na\"ive sense---i.e.~when the Segal maps are $G$-weak equivalences for \emph{trivial} $G$-sets or, equivalently, when all fixed points $X^H$ for $H\subset G$ are special in the non-equivariant sense---the zig-zag
\begin{equation*}
\begin{tikzcd}[cramped]
X(1^+)\times X(1^+) & \arrow[l, "\sim", "\rho"'] X(2^+)\arrow[r, "X(\mu)"] & X(1^+)
\end{tikzcd}
\end{equation*}
equips $\pi_0^H(X(1^+))=\pi_0(X(1^+)^H)$ with an abelian monoid structure for all $H\subset G$.
\end{rk}

The special $\Gamma$-$G$-spaces are not invariant under the na\"ive notion of levelwise weak equivalences. Instead, one should require that $f(S_+)$ be a $G$-equivariant weak equivalence for every finite $G$-set $S$ (and not only the trivial $G$-sets); by the same argument as in Lemma~\ref{lemma:special-via-graph-subgroups} above, this is equivalent to the following condition, also see~\cite[Remark~4.11]{equivariant-gamma}:

\begin{defi}
A map $f\colon X\to Y$ of $\Gamma$-$G$-spaces is called a \emph{$G$-equivariant level weak equivalence}\index{level weak equivalence!G-equivariant@$G$-equivariant|seeonly{$G$-equivariant level weak equivalence}}\index{G-equivariant level weak equivalence@$G$-equivariant level weak equivalence!in Gamma-G-SSet@in $\cat{$\bm\Gamma$-$\bm G$-SSet}_*$|textbf} if $f(S_+)$ is a $\mathcal G_{G,\Sigma_S}$-equivariant weak equivalence for every finite set $S$.
\end{defi}

Note that a map of $\Gamma$-$G$-spaces with trivial $G$-action can be a non-equivariant level weak equivalence and still fail to be a $G$-equivariant level weak equivalence as the following example for $G=\Sigma_2$ shows:

\begin{ex}\label{ex:non-equiv-but-not-equiv-level-we}
We define a $\Gamma$-space $X$ as follows: $X(S_+)$ is the nerve of the category with objects $\Gamma(2^+,S_+)$ and a unique morphism $\alpha\to\beta$ (possibly the identity) for any two $\alpha,\beta\colon 2^+\to S_+$ with $\im\alpha=\im\beta$.
One easily checks that $\Gamma(2^+,f)$ extends uniquely to $X(S_+)\to X(T_+)$ for any $f\colon S_+\to T_+$ in $\Gamma$, which then makes $X$ into a $\Gamma$-space.

Let now $Y$ be the quotient obtained from $\Gamma(2^+,\blank)$ by identifying for each $S_+\in\Gamma$ any two $\alpha,\beta\colon 2^+\to S_+$ with the same image. Then the quotient map extends uniquely to $f\colon X\to Y$, and this is a non-equivariant level weak equivalence as it is given by contracting disjoint copies of $E\Sigma_2$.

However, with respect to the usual $\Sigma_2$-action on $2^+$, $f(2^+)^{\Sigma_2}\colon X(2^+)^{\Sigma_2}\to Y(2^+)^{\Sigma_2}$ is a map from a discrete simplicial set with $2$ vertices to one with $3$ vertices, so $f$ is not a $\Sigma_2$-equivariant level weak equivalence.
\end{ex}

\begin{rk}
As remarked by Ostermayr \cite[Theorem~4.7]{equivariant-gamma} without proof, the $G$-equivariant level weak equivalences are part of a `generalized projective model structure' on $\cat{$\bm\Gamma$-$\bm G$-SSet}_*$. The fibrations in this model structure are precisely those maps $f$ such that $f(S_+)$ is a fibration in the $\mathcal G_{G,\Sigma_S}$-model structure on $\cat{$\bm{(G\times\Sigma_S)}$-SSet}$ for all finite sets $S$. In addition, Ostermayr also constructs an equivariant version of the Bousfield-Friedlander model structure in \cite[Theorem~4.12]{equivariant-gamma}.
\end{rk}

\subsection{$\bm G$-global level model structures}\index{G-global model structure@$G$-global model structure!level|seeonly{$G$-global level model structure}}
\index{G-global model structure@$G$-global model structure!positive level|seeonly{positive $G$-global level model structure}}
In this subsection we will introduce various $G$-global analogues of $\Gamma$-spaces together with suitable level model structures. As the arguments for the existence of these model structures are rather similar,  we will formalize them.

For this we begin with the following `undirected' version of \cite[Proposition C.23]{schwede-book} (the relevant parts of which we recalled as Proposition \ref{prop:generalized-projective-dim} above):

\begin{prop}\label{prop:generalized-projective}\index{generalized projective model structure!undirected|textbf}
Let $A$ be an essentially small $\cat{Set}_*$-enriched category and let $\mathscr C$ be a locally presentable category with $0$-object (which then has a unique $\cat{Set}_*$-enrichment). Assume we are given for each $a\in A$ a combinatorial model structure on the ordinary functor category $\Aut_{A}(a)\text{-}\mathscr C\mathrel{:=}\Fun(B\Aut_{A}(a),\mathscr C)$ with sets of generating cofibrations $I_a$ and generating acyclic cofibrations $J_a$, satisfying the following `consistency condition': for all $b\in A$, any relative $\{\Hom_A(a,b)\otimes_{\Aut_A(a)}j : a\in A, j\in J_A\}$-cell complex is a weak equivalence in $\Aut_A(b)\text{-}\mathscr C$.

Then there exists a unique model structure on the category $A\text{-}\mathscr C$ of $\cat{Set}_*$-enriched functors $A\to\mathscr C$ in which a map $f\colon X\to Y$ is a weak equivalence or fibration if and only if $f(a)\colon X(a)\to Y(a)$ is a weak equivalence or fibration, respectively, in $\Aut_A(a)\text{-}\mathscr C$ for each $a\in A$. We call this the \emph{generalized projective model structure}. It is combinatorial with generating cofibrations
\begin{equation*}
I_A\mathrel{:=}\{\Hom(a,\blank)\otimes_{\Aut(a)}i \colon a\in A, i\in I_a\}
\end{equation*}
and generating acyclic cofibrations.
\begin{equation*}
J_A\mathrel{:=}\{\Hom(a,\blank)\otimes_{\Aut(a)}j \colon a\in A, j\in J_a\}.
\end{equation*}
Moreover:
\begin{enumerate}
\item If each of the model categories $\Aut(a)\text{-}\mathscr C$ is simplicial, then so is $A\text{-}\mathscr C$.
\item If each of the model categories $\Aut(a)\text{-}\mathscr C$ is right proper, then so is $A\text{-}\mathscr C$.
\item If filtered colimits are homotopical in each of the $\Aut(a)\text{-}\mathscr C$, then also filtered colimits in $A\text{-}\mathscr C$ are homotopical.
\end{enumerate}
\end{prop}

Here we again use for a pointed $G$-$H$-biset $X$ the notation $X\otimes_H Y$ for the \emph{balanced tensor product},\nomenclature[aH]{$\blank\otimes_H\blank$}{balanced tensor product (quotient of tensoring by diagonal $H$-action)} i.e.~the quotient of the $\cat{Set}_*$-tensoring $X\otimes Y$ by the diagonal $H$-action, with the induced $G$-action. In particular, in the world of pointed simplicial sets $X\otimes_H\blank$ can be identified with the usual balanced smash product $X\smashp_H\blank$.

\begin{proof}[Proof of Proposition~\ref{prop:generalized-projective}]
For $a\in A$, let $E_a\colon A\text{-}\mathscr C\to\Aut_A(a)\text{-}\mathscr C$ denote the evaluation functor. This admits a left adjoint $G_a$, which by the usual coend formula for enriched left Kan extensions (or the Yoneda Lemma) can be calculated as
\begin{equation*}
G_a(X)(b)=\Hom_A(a,b)\otimes_{\Aut_A(a)}X
\end{equation*}
\nomenclature[aGa]{$G_a$}{semifree diagram/spectrum, left adjoint to evaluation at $a$ with induced $\Aut(a)$-action}%
with the obvious functoriality in each variable. Thus,
\begin{equation*}
I_A=\bigcup_{a\in A} G_a(I_a)\qquad\text{and}\qquad J_A=\bigcup_{a\in A}G_a(J_a).
\end{equation*}
To construct the model structure in question and to show that it is cofibrantly generated (hence combinatorial) with generating cofibrations $I_A$ and generating acyclic cofibrations $J_A$, it therefore suffices by Crans' criterion (Proposition~\ref{prop:transfer-criterion}) applied to $\prod_{a\in\tilde A}\Aut_A(a)\text{-}\mathscr{C}$ (where $\tilde A$ is any small essentially wide subcategory of $A$) to show, also cf.~\cite[Theorem~A.1 and Remark~A.2]{cellular}:
\begin{enumerate}
\item The sets $I_A$ and $J_A$ permit the small object argument, \emph{and}
\item Relative $J_A$-cell complexes are weak equivalences (i.e.~sent to weak equivalences in $\Aut_A(a)\text{-}\mathscr C$ under $E_a$ for all $a\in A$).
\end{enumerate}
The first condition is automatically satisfied as $\mathscr C$ and hence also $A\text{-}\mathscr C$ is locally presentable. On the other hand, the second condition is an immediate consequence of the consistency condition.

The additional properties $(1)$--$(3)$ follow easily from the fact that all the relevant constructions are defined levelwise, also see Lemma~\ref{lemma:transferred-properties} for the first two properties.
\end{proof}

Instead of verifying the consistency condition by hand, we will employ the following criterion, which at the same time takes care of left properness:

\begin{prop}\label{prop:left-proper-projective}\index{generalized projective model structure!undirected|textbf}
Let $A$ be an essentially small $\cat{Set}_*$-enriched category and let $\mathscr C$ be a pointed locally presentable category. Assume we are given for each $a\in A$ a cofibrantly generated (hence combinatorial) model structure on $\Aut_{A}(a)\text{-}\mathscr C$ with sets of generating cofibrations $I_a$ and generating acyclic cofibrations $J_a$, such that $\Aut_{A}(a)\text{-}\mathscr C$ is left proper and such that filtered colimits in it are homotopical. Assume moreover:
\begin{enumerate}
\item For all $a,b\in A$ and $j\in J_a$, the map $\Hom(a,b)\otimes_{\Aut(a)}j$ is a weak equivalence in $\Aut(b)\text{-}\mathscr C$.
\item For all $a,b\in A$ and $i\in I_a$, every pushout along $\Hom(a,b)\otimes_{\Aut(a)}i$ is a homotopy pushout in $\Aut(b)\text{-}\mathscr C$.
\end{enumerate}
Then the generalized projective model structure exists and is combinatorial with generating cofibrations $I_A$ and generating acyclic cofibrations $J_A$ as above. Moreover, it is left proper and a square
\begin{equation}\label{diag:hopo-A}
\begin{tikzcd}
W\arrow[r]\arrow[d] & X\arrow[d]\\
Y\arrow[r] & Z
\end{tikzcd}
\end{equation}
in $A\text{-}\mathscr C$ is a homotopy pushout if and only if it is a levelwise homotopy pushout, i.e.~for every $b\in A$ the induced square
\begin{equation*}
\begin{tikzcd}
W(b)\arrow[r]\arrow[d] & X(b)\arrow[d]\\
Y(b)\arrow[r] & Z(b)
\end{tikzcd}
\end{equation*}
is a homotopy pushout in $\Aut(b)\text{-}\mathscr C$. Finally, filtered colimits are homotopical in $A\text{-}\mathscr C$, and the conclusions $(1)$ and $(2)$ of Proposition~\ref{prop:generalized-projective} still hold.
\begin{proof}
We begin with the following observation:

\begin{claim*}
If $i$ is any $I_A$-cofibration, then pushouts along $i(b)$ are homotopy pushouts in $\Aut(b)\text{-}\mathscr C$ for any $b\in A$. In particular, pushouts along $i$ in $A\text{-}\mathscr C$ are levelwise homotopy pushouts.
\begin{proof}
We fix $b\in A$ and consider the class $\mathscr H_b$ of all morphisms $i$ in $A\text{-}\mathscr C$ such that pushouts along $i(b)$ are homotopy pushouts. Using Proposition~\ref{prop:U-sharp-closure} and that $\ev_b$ is cocontinuous, we see that $\mathscr H_b$ is closed under pushouts, transfinite compositions, and retracts. On the other hand,  $I_A\subset\mathscr H_b$ by assumption, so that $\mathscr H_b$ contains all $I_A$-cofibrations as desired.
\end{proof}
\end{claim*}

As each $G_a$ is cocontinuous, $G_a(i)$ is an $I_A$-cofibration for any $a\in A$ and any cofibration $i$ in $\Aut(a)\text{-}\mathscr C$, so we in particular see that the second assumption holds more generally for \emph{all} cofibrations $i$ of $\Aut(a)\text{-}\mathscr C$.

Let us now verify the assumptions of Proposition~\ref{prop:generalized-projective}, i.e.~that for any $b\in B$ any transfinite composition of pushouts of maps of the form $\Hom_A(a,b)\otimes_{\Aut_A(a)}j$ ($a\in A,j\in J_A$) is a weak equivalence. Indeed, any $\Hom_A(a,b)\otimes_{\Aut_A(a)}j$ is a weak equivalence by the first assumption, and the above strengthening of the second assumption then implies that also pushouts of it are weak equivalences. The claim follows as filtered colimits were assumed to be homotopical.

We are therefore allowed to apply Proposition~\ref{prop:generalized-projective}, and it only remains to show that $A\text{-}\mathscr C$ is left proper and that the homotopy pushouts are precisely the levelwise homotopy pushouts. But the the functor $(E_a)_{a\in\tilde A}\colon A\text{-}\mathscr C\to\prod_{a\in\tilde A}\Aut_A(a)\text{-}\mathscr C$ (for $\tilde A$ as above) clearly preserves and reflects weak equivalences, and it sends pushouts along cofibrations to homotopy pushouts by the above claim, so this is simply an instance of Lemma~\ref{lemma:U-pushout-preserve-reflect}.
\end{proof}
\end{prop}

\subsubsection{Pointed $G$-global homotopy theory}
\index{pointed G-global homotopy theory@pointed $G$-global homotopy theory|(}
The above model categories are inherently pointed, so in order to apply them to $G$-global homotopy theory, we need to develop suitable pointed versions of the models of Chapter~\ref{chapter:unstable}. These are in fact formal consequences of the corresponding unbased results, and in order to avoid a long list of similar looking statements, we will be somewhat terse here.

Let $\mathscr C$ be a category with terminal object $*$. Then we write $\mathscr C_*\mathrel{:=}*\downarrow\mathscr C$ for the category of objects under $*$ and call it the category of \emph{pointed objects in $\mathscr C$}. If $\mathscr C$ has binary coproducts, then the forgetful functor $\mathscr C_*\to\mathscr C$ has a left adjoint $(\blank)_+$ sending an object $X$ to $X_+=X\amalg *$\nomenclature[zp]{$(\blank)_+$}{left adjoint to forgetful functor from pointed to unpointed category} (with structure map the inclusion of the second summand), and similar on morphisms.

\begin{thm}\index{G-global model structure@$G$-global model structure!pointed versions}
There is a unique model structure on $\cat{$\bm{E\mathcal M}$-$\bm G$-SSet}_*$ in which a map is a weak equivalence, fibration, or cofibration if and only if it is so in the $G$-global model structure on $\cat{$\bm{E\mathcal M}$-$\bm G$-SSet}$. We call this the \emph{$G$-global model structure} again. It is proper, simplicial, and combinatorial with generating cofibrations
\begin{equation*}
\{(E\mathcal M\times_\phi G\times\del\Delta^n)_+\hookrightarrow(E\mathcal M\times_\phi G\times\Delta^n)_+ : H\subset\mathcal M\text{ universal},\phi\colon H\to G,n\ge0\}
\end{equation*}
and generating acyclic cofibrations
\begin{equation*}
\{(E\mathcal M\times_\phi G\times\Lambda_k^n)_+\mskip0mu minus 2mu\hookrightarrow\mskip0mu minus 2mu(E\mathcal M\times_\phi G\times\Delta^n)_+\! : H\mskip0mu minus 2mu\subset\mskip0mu minus 2mu \mathcal M\text{ universal},\phi\colon H\mskip0mu minus 2mu\to\mskip0mu minus 2mu G,0\le k\le n\}.
\end{equation*}
Moreover, filtered colimits in it are homotopical, and pushouts along injective cofibrations are homotopy pushouts.
\begin{proof}
The model structure exists by \cite[Theorem~7.6.5-(1)]{hirschhorn} applied to the model structure from Corollary~\ref{cor:EM-G-model-structure}. Moreover, the model structure is cofibrantly generated by the above sets according to \cite[Theorem~2.7]{hirschhorn-slice}, hence combinatorial, and it is proper by \cite[Theorem~2.8-(3)]{hirschhorn-slice}. Finally, filtered colimits in $\cat{$\bm{E\mathcal M}$-$\bm G$-SSet}_*$ are homotopical since both weak equivalences as well as connected colimits are created in $\cat{$\bm{E\mathcal M}$-$\bm G$-SSet}$. Likewise, the characterization of homotopy pushouts follows from Lemma~\ref{lemma:U-pushout-preserve-reflect} applied to the forgetful functor to $\cat{$\bm{E\mathcal M}$-$\bm G$-SSet}$.

Finally, with respect to the evident simplicial enrichment $\cat{$\bm{E\mathcal M}$-$\bm G$-SSet}_*$ is tensored and cotensored over $\cat{SSet}$ with the cotensoring created in $\cat{$\bm{E\mathcal M}$-$\bm G$-SSet}$. It follows that the above model structure is simplicial.
\end{proof}
\end{thm}

One gets analogous statements for all the other model structures established in Chapter~\ref{chapter:unstable}; instead of making them explicit, we will freely refer to the corresponding unbased statement whenever we actually need the based statement.

\index{functoriality in homomorphisms!pointed versions|(}
If $F\colon\mathscr C\rightleftarrows\mathscr D:G$ is a Quillen adjunction such that $F$ preserves the terminal object, then we have an induced adjunction $\mathscr C_*\rightleftarrows\mathscr D_*$, which is obviously a Quillen adjunction again. In particular, if $\alpha\colon H\to G$ is any homomorphism, then we can apply this to the various change of group adjunctions $\alpha^*\dashv\alpha_*$ discussed in Chapter~\ref{chapter:unstable}. Moreover, as weak equivalences are created in the corresponding models of unstable $G$-global homotopy theory, all results established above on whether $\alpha^*$ or $\alpha_*$ is fully homotopical, immediately transfer from the unpointed to the pointed setting. We will therefore freely refer to the results from Chapter~\ref{chapter:unstable} when arguing about the corresonding based adjunctions.

While the functors $\alpha_!:H\textup{-}\mathscr C\to G\textup{-}\mathscr C$ usually do not preserve the base point, the functor $\alpha^*:G\textup{-}\mathscr C_*\to H\textup{-}\mathscr C_*$ still has a left adjoint in each case; as above one deduces that $\alpha^*$ is still right Quillen, so these are again Quillen adjunctions. The only results that do not immediately transfer via the above strategy are the criteria on when $\alpha_!$ is suitably homotopical. Let us therefore explicitly prove those statements that we will need later:

\begin{cor}\label{cor:alpha-shriek-injective-EM-based}
Let $\alpha\colon H\to G$ be an \emph{injective} homomorphism. Then $\alpha_!\colon\cat{$\bm{E\mathcal M}$-$\bm H$-SSet}_*\to\cat{$\bm{E\mathcal M}$-$\bm G$-SSet}_*$ and $\alpha_!\colon\cat{$\bm{E\mathcal M}$-$\bm H$-SSet}_*^\tau\to\cat{$\bm{E\mathcal M}$-$\bm G$-SSet}_*^\tau$ are homotopical.
\begin{proof}
It suffices to prove the first statement, for which we observe that $\alpha_!$ is left Quillen for the injective model structures by the above arguments applied to Corollary~\ref{cor:alpha-shriek-injective-EM}. As all pointed $E\mathcal M$-$H$-simplicial sets are injectively cofibrant, the claim follows from Ken Brown's Lemma.
\end{proof}
\end{cor}

\begin{cor}\label{cor:free-quotient-EM-star}
Let $\alpha\colon H\to G$ be any homomorphism. Then
\begin{equation*}
\alpha_!\colon\cat{$\bm{E\mathcal M}$-$\bm H$-SSet}_*\to\cat{$\bm{E\mathcal M}$-$\bm G$-SSet}_*\text{ and }\alpha_!\colon\cat{$\bm{E\mathcal M}$-$\bm H$-SSet}_*^\tau\to\cat{$\bm{E\mathcal M}$-$\bm G$-SSet}_*^\tau
\end{equation*}
preserve weak equivalences between objects with free $\ker(\alpha)$-action outside the base point.
\begin{proof}
It again suffices to prove the first statement, and the previous corollary reduces to the case that $\alpha$ is surjective, i.e.~$\alpha_!$ is given by dividing out the action of $K\mathrel{:=}\ker(\alpha)$. If we let $L\subset\mathcal M$ be any subgroup, then \cite[Lemma~A.1]{hausmann-equivariant} shows that there is for any $\phi\colon L\to G$ a natural isomorphism
\begin{equation*}
\bigvee_{[\psi\colon L\to H]} Z^{\psi}/(\centralizer_H(\im\psi)\cap K)\cong (Z/K)^\phi
\end{equation*}
analogous to Proposition~\ref{prop:free-quotient-general} whenever $K$ acts freely on $Z$ outside the basepoint. Analogously to the unbased case, $\centralizer_H(\im\psi)\cap K$ acts freely on $Z^\psi$ outside the base point, so $Z^\psi$ is cofibrant in the (non-equivariant) projective model structure on $\cat{$\bm{(\centralizer_H(\im\psi)\cap K})$-SSet}_*$. As the quotient functor $\cat{$\bm{(\centralizer_H(\im\psi)\cap K)}$-SSet}_*\to\cat{SSet}_*$ is left Quillen, the claim follows.
\end{proof}
\end{cor}

\begin{rk}\label{rk:free-quotient-based}
In fact, the same argument as above more generally yields a pointed version of Proposition~\ref{prop:free-quotient-general}.\index{functoriality in homomorphisms!pointed versions|)}
\index{pointed G-global homotopy theory@pointed $G$-global homotopy theory|)}
\end{rk}

\subsubsection{Construction of $G$-global level model structures} With this at hand we can now introduce the desired level model structures:

\begin{defi}
A map $f$ in $\cat{$\bm\Gamma$-$\bm{E\mathcal M}$-$\bm G$-SSet}_*$ is called a \emph{$G$-global level weak equivalence}\index{G-global level weak equivalence@$G$-global level weak equivalence!in Gamma-EM-G-SSet@in $\cat{$\bm\Gamma$-$\bm{E\mathcal M}$-$\bm G$-SSet}_*$|textbf}\index{level weak equivalence!G-global@$G$-global|seeonly{$G$-global level weak equivalence}} or \emph{$G$-global level fibration} if $f(S_+)$ is a $(G\times\Sigma_S)$-global weak equivalence or fibration, respectively, for every finite set $S$.
\end{defi}

\begin{thm}\label{thm:Gamma-EM-SSet}\index{G-global level model structure@$G$-global level model structure!on Gamma-EM-G-SSet@on $\cat{$\bm\Gamma$-$\bm{E\mathcal M}$-$\bm G$-SSet}_*$|textbf}\index{G-global Gamma-space@$G$-global $\Gamma$-space}
The $G$-global level weak equivalences and fibrations are part of a unique model structure on $\cat{$\bm\Gamma$-$\bm{E\mathcal M}$-$\bm G$-SSet}_*$, which we call the \emph{$G$-global level model structure}. It is proper, simplicial, and combinatorial with generating cofibrations
\begin{equation*}
\big\{\big(\Gamma(S_+,\blank)\smashp(E\mathcal M\times G)_+\big)/H\smashp(\del\Delta^n\hookrightarrow\Delta^n)_+ :  H\in\mathcal G_{\mathcal U,G}, \text{$S$ finite $H$-set}, n\ge0\big\}
\end{equation*}
(where $\mathcal U$ denotes the collection of universal subgroups) and generating acyclic cofibrations
\begin{equation*}
\big\{\big(\Gamma(S_+,\blank)\smashp(E\mathcal M\times G)_+\big)/H\smashp(\Lambda^n_k\hookrightarrow\Delta^n)_+ :  H\in\mathcal G_{\mathcal U,G}, \text{$S$ finite $H$-set}, 0\le k\le n\big\}
\end{equation*}
Moreover, filtered colimits in it are homotopical, and a square is a homotopy pushout with respect to it if and only if it so levelwise. In particular, pushouts along injective cofibrations (i.e.~levelwise injections) are homotopy pushouts.
\begin{proof}
One immediately proves by inspection that
\begin{equation}\label{eq:induction-m-gamma-sset}
\Gamma(A_+,B_+)\smashp_{\Sigma_A}\blank\colon\cat{$\bm{E\mathcal M}$-$\bm{(G\times\Sigma_A)}$-SSet}_*\to\cat{$\bm{E\mathcal M}$-$\bm{(G\times\Sigma_B)}$-SSet}_*
\end{equation}
sends the generating (acyclic) cofibrations of the usual $(G\times\Sigma_A)$-global model structure on the source to (acyclic) cofibrations in the \emph{injective} $(G\times\Sigma_B)$-global model structure on the target (where we for simplicity confuse $\Aut_\Gamma(S_+)$ with $\Sigma_S$ for any finite set $S$). We may therefore apply Proposition~\ref{prop:left-proper-projective}, which proves all of the above statements except for the description of the generating cofibrations.

Instead, the aforementioned proposition shows that a set of generating cofibrations is given by the maps
\begin{equation*}
\Gamma(S_+,\blank)\smashp_{\Sigma_S} (E\mathcal M\times G\times\Sigma_S\times\del\Delta^n)_+/H\to
\Gamma(S_+,\blank)\smashp_{\Sigma_S} (E\mathcal M\times G\times\Sigma_S\times\Delta^n)_+/H
\end{equation*}
induced by the inclusions $\del\Delta^n\hookrightarrow\Delta^n$, with $H\in\mathcal G_{\mathcal U,G\times\Sigma_S}$, $n\ge 0$ as above, and $S$ an ordinary finite set. This is obviously conjugate to
\begin{equation*}
\Gamma(S_+,\blank)\smashp_{\Sigma_S} (E\mathcal M\times G\times\Sigma_S)_+\smashp(\del\Delta^n\hookrightarrow\Delta^n)_+;
\end{equation*}
on the other hand, $\Gamma(S_+,\blank)\smashp_{\Sigma_S} (E\mathcal M\times G\times\Sigma_S)_+\cong\Gamma(S_+,\blank)\smashp (E\mathcal M\times G)_+$ via $[f,m,g,\sigma]\mapsto [f\sigma,m,g]$, which is right $\Sigma_S$- and both left and right $(E\mathcal M\times G)$-equivariant. As quotients and smash products preserve colimits, we conclude $\Gamma(S_+,\blank)\smashp_{\Sigma_S}(E\mathcal M\times G\times\Sigma_S)/\Gamma_{H,(\phi,\rho)}\cong(\Gamma(S_+,\blank)\smashp(E\mathcal M\times G)_+)/\Gamma_{H,(\phi,\rho)}$
for any $H\subset\mathcal M$ and any homomorphisms $\phi\colon H\to G$ and $\rho\colon H\to\Sigma_S$. Finally, if we view $S$ as an $\Gamma_{H,\phi}$-set via $\rho\circ\pr_H$, then the above equals $(\Gamma(S_+,\blank)\smashp (EM\times G)_+)/\Gamma_{H,\phi}$. This completes the verification of the given set of generating cofibrations; the argument for the generating acyclic cofibrations is analogous.
\end{proof}
\end{thm}

\begin{warn}\label{warn:strong-consistency-fails}
If $|B|\ge3$, then $(\ref{eq:induction-m-gamma-sset})$ is \emph{not} left Quillen with respect to the ordinary $(G\times\Sigma_B)$-global model structure on the target, even for $G=1$. Namely, any cofibrant object in the latter model structure has free $\Sigma_B$-action outside the basepoint while
$\Gamma(A_+,B_+)\smashp_{\Sigma_A} (E\mathcal M\times\Sigma_A)_+$
has non-trivial $\Sigma_B$-fixed points: for example, if $b\in B$ is arbitrary and $\beta\colon A_+\to B_+$ denotes the map sending every non-basepoint to $b$, then $[\beta,\id_\omega,\id_A]$ is fixed by the non-trivial subgroup of $\Sigma_B^b\subset\Sigma_B$ of the bijections fixing $b$.

In particular, if $f$ is a $G$-global cofibration, then $f(B_+)$ need not be a $(G\times\Sigma_B)$-global cofibration for $|B|\ge3$. However, one can at least show that it is a $G$-global cofibration, cf.~Proposition~\ref{prop:tame-gamma-cof-levelwise} below.
\end{warn}

Next, we will give a tame analogue of this:

\begin{defi}
A map $f$ in $\cat{$\bm\Gamma$-$\bm{E\mathcal M}$-$\bm G$-SSet}_*^\tau$ is called a \emph{$G$-global level weak equivalence}\index{G-global level weak equivalence@$G$-global level weak equivalence!in Gamma-EM-G-SSettau@in $\cat{$\bm\Gamma$-$\bm{E\mathcal M}$-$\bm G$-SSet}_*^\tau$|textbf} or \emph{$G$-global positive level fibration} if $f(S_+)$ is a weak equivalence or fibration, respectively, in the {positive} $(G\times\Sigma_S)$-global model structure on $\cat{$\bm{(G\times\Sigma_S)}$-$\bm{E\mathcal M}$-SSet}^\tau$ for every finite set $S$.
\end{defi}

\begin{thm}\label{thm:Gamma-EM-tau}\index{G-global level model structure@$G$-global level model structure!positive|seeonly{positive $G$-global level model structure}}\index{G-global Gamma-space@$G$-global $\Gamma$-space}
\index{positive G-global level model structure@positive $G$-global level model structure!on Gamma-EM-G-SSettau@on $\cat{$\bm\Gamma$-$\bm{E\mathcal M}$-$\bm G$-SSet}^\tau_*$|textbf}
The $G$-global level weak equivalences and $G$-global positive level fibrations are part of a unique model structure on $\cat{$\bm\Gamma$-$\bm{E\mathcal M}$-$\bm G$-SSet}_*^\tau$, which we call the \emph{positive $G$-global level model structure}. It is combinatorial with generating cofibrations
\begin{equation*}
\big(\Gamma(S_+,\blank)\smashp(E\Inj(A,\blank)\times G)_+\big)/H\smashp(\del\Delta^n\hookrightarrow\Delta^n)_+
\end{equation*}
where $H$ runs through finite groups, $S$ through finite $H$-sets, and $A$ through non-empty finite faithful $H$-sets.
Moreover, it is simplicial, proper, filtered colimits in it are homotopical, and a square is a homotopy pushout if and only if it so levelwise. In particular, pushouts along injective cofibrations are homotopy pushouts.
\end{thm}

Here we use the \emph{positive} model structure merely to faciliate the comparison to $G$-parsummable simplicial sets given in the next section.

To prove the theorem we will employ:

\begin{lemma}\label{lemma:smash-g-global-tame}
Let $K,L$ be discrete groups, let $i\colon Y\to Z$ be an acyclic cofibration in the $(G\times K)$-global model structure on $\cat{$\bm{E\mathcal M}$-$\bm{(G\times K)}$-SSet}_*^\tau$, and let $X$ be a left-$L$-right-$K$-simplicial set. Then $X\smashp_Ki$ is a $(G\times L)$-global weak equivalence.
\begin{proof}
Let $H\subset\mathcal M$ be any universal subgroup. By Lemma~\ref{lemma:g-global-pos-cof}, $i$ is a $\mathcal G_{H,G\times K}$-cofibration, and it is obviously acyclic. Thus, it is enough to show that $X\smashp_K\blank$ sends acyclic cofibrations in the $\mathcal G_{H,G\times K}$-model structure to acyclic cofibrations in the injective $\mathcal G_{H,G\times L}$-model structure. But as before it suffices to prove this for the generating acyclic cofibrations, where this is trivial.
\end{proof}
\end{lemma}

\begin{proof}[Proof of Theorem~\ref{thm:Gamma-EM-tau}]
By the same arguments as in Theorem~\ref{thm:Gamma-EM-SSet} it suffices that for every $A_+,B_+\in\Gamma$ the functor
\begin{equation*}
\Gamma(A_+,B_+)\smashp_{\Sigma_A}\blank\colon\cat{$\bm{E\mathcal M}$-$\bm{(G\times\Sigma_A)}$-SSet}_*^\tau\to\cat{$\bm{E\mathcal M}$-$\bm{(G\times\Sigma_B)}$-SSet}_*^\tau
\end{equation*}
sends positive $(G\times\Sigma_A)$-global (acyclic) cofibrations to (acyclic) cofibrations in the injective $(G\times\Sigma_B)$-global model structure: indeed, it is clear that $\Gamma(A_+,B_+)\smashp_{\Sigma_A}i$ is an injective cofibration, and if $i$ is acyclic, then the previous lemma implies that $\Gamma(A_+,B_+)\smashp_{\Sigma_A}i$ is a $(G\times\Sigma_B)$-global weak equivalence as desired.
\end{proof}

As before, if $i$ is a cofibration in the above model structure, then $i(S_+)$ will in general not even be a $(G\times\Sigma_S)$-global cofibration for a finite set $S$ with at least three elements. However we have:

\begin{prop}\label{prop:tame-gamma-cof-levelwise}
Let $f\colon X\to Y$ be a cofibration in $\cat{$\bm\Gamma$-$\bm{E\mathcal M}$-$\bm G$-SSet}_*^\tau$, and let $S$ be any finite set. Then $f(S_+)$ is a positive $G$-global cofibration.
\begin{proof}
It suffices to prove this for generating cofibrations, for which we are then reduced by the same arguments as before to showing that
\begin{equation*}
\big(\Gamma(T_+,S_+)\smashp (G\times E\Inj(A,\omega))_+\big)/H\cong\big(E\Inj(A,\omega)_+\smashp (G_+\smashp \Gamma(T_+,S_+))\big)/H
\end{equation*}
is cofibrant in the positive $G$-global model structure for any finite group $H$, any non-empty finite faithful $H$-set $A$, any finite $H$-set $T$, and any homomorphism $\phi\colon H\to G$, where $H$ acts from the right via its given actions on $A$ and $T$ and via its right action on $G$ via $\phi$.

We now claim that the functor
\begin{equation}\label{eq:cof-producing-adjunction}
\begin{aligned}
\cat{$\bm{(H^\op\times G)}$-SSet}_*&\to\cat{$\bm{E\mathcal M}$-$\bm{G}$-SSet}_*^\tau\\
X&\mapsto (E\Inj(A,\omega)_+\smashp X)/H
\end{aligned}
\end{equation}
is left Quillen with respect to the $\mathcal G_{H^\op,G}$-equivariant model structure on the source and the positive $G$-global model structure on the target. As the cofibrant objects on the source are precisely those simplicial sets with free $G$-action outside the basepoint, this will then immediately imply the proposition.

To prove the claim, we observe that $(\ref{eq:cof-producing-adjunction})$ factors up to isomorphism as the composition of the left Quillen functor $\ev_\omega\colon\cat{$\bm G$-$\bm{\mathcal I}$-SSet}_*\to\cat{$\bm{E\mathcal M}$-$\bm G$-SSet}_*^\tau$ and the functor
\begin{equation*}
\cat{$\bm{(H^\op\times G)}$-SSet}_*\to\cat{$\bm{G}$-$\bm{\mathcal I}$-SSet}_*, X\mapsto\mathcal (\mathcal I(A,\blank)_+\smashp X)/H
\end{equation*}
which is left adjoint to $X\mapsto X(A)^\phi$. The latter is by definition right Quillen with respect to the $G$-global positive level model structure on the source, hence also with respect to the $G$-global positive model structure. This completes the proof of the claim and hence of the proposition.
\end{proof}
\end{prop}

Together with Lemma~\ref{lemma:g-global-pos-cof} we immediately conclude:

\begin{cor}\label{cor:cofibrant-pos-g-global}
Let $i\colon X\to Y$ be a cofibration in the positive $G$-global model structure on $\cat{$\bm\Gamma$-$\bm{E\mathcal M}$-$\bm G$-SSet}_*^\tau$ and let $S$ be any finite set. Then:
\begin{enumerate}
\item $Y(S_+)$ has no simplices of empty support outside the image of $i(S_+)$.\label{item:cpgg-support}
\item $G$ acts freely on $Y(S_+)$ outside the image of $i(S_+)$.\label{item:cpgg-free}\qed
\end{enumerate}
\end{cor}

Finally, let us consider the models based on $\mathcal I$-simplicial sets:

\begin{thm}\label{thm:gamma-script-I}\index{G-global level model structure@$G$-global level model structure!on Gamma-G-II-SSettau@on $\cat{$\bm\Gamma$-$\bm G$-$\bm{\mathcal I}$-SSet}_*$|textbf}\index{G-global Gamma-space@$G$-global $\Gamma$-space}
There exists a unique model structure on $\cat{$\bm\Gamma$-$\bm G$-$\bm{\mathcal I}$-SSet}_*$ in which a map $f$ is a weak equivalence or fibration if and only if $f(S_+)$ is a $(G\times\Sigma_S)$-global weak equivalence or fibration, respectively, for any finite set $S$.

We call this the \emph{$G$-global level model structure} and its weak equivalences the \emph{$G$-global level weak equivalences}.\index{G-global level weak equivalence@$G$-global level weak equivalence!in Gamma-G-I-SSet@in $\cat{$\bm\Gamma$-$\bm G$-$\bm{\mathcal I}$-SSet}_*$|textbf} It is combinatorial with generating cofibrations
\begin{align*}
\big\{\big(\Gamma(S_+,\blank)\smashp(\mathcal I(A,\blank)\times G)_+\big)/H\smashp(\del\Delta^n\hookrightarrow\Delta^n)_+ : {}&\text{$H$ finite group}, \text{$S$ finite $H$-set},\\ &\text{$A$ finite faithful $H$-set}\big\}
\end{align*}
and moreover simplicial, proper, and filtered colimits in it are homotopical. Finally, a square is a homotopy pushout if and only if it so levelwise. In particular, pushouts along injective cofibrations are homotopy pushouts.
\begin{proof}
As before it suffices that $\Gamma(A_+,B_+)\smashp_{\Sigma_A}\blank\colon\cat{$\bm{(G\times\Sigma_A)}$-$\bm{\mathcal I}$-SSet}_*\to\cat{$\bm{(G\times\Sigma_B)}$-$\bm{\mathcal I}$-SSet}_*$ sends $(G\times\Sigma_A)$-global (acyclic) cofibrations to (acyclic) cofibrations in the injective $(G\times\Sigma_B)$-global model structures for all $A_+,B_+\in\Gamma$.

If $i$ is a $(G\times\Sigma_A)$-global cofibration, then it is clear that $\Gamma(A_+,B_+)\smashp_{\Sigma_A}i$ is an injective cofibration. Moreover, if $i$ is acyclic, then $i(\omega)$ is an acyclic cofibration in the $(G\times\Sigma_A)$-global model structure on $\cat{$\bm{E\mathcal M}$-$\bm{(G\times\Sigma_A)}$-SSet}_*^\tau$. As $(\Gamma(A_+,B_+)\smashp_{\Sigma_A}i)(\omega)$ is conjugate to $\Gamma(A_+,B_+)\smashp_{\Sigma_A}(i(\omega))$, the acyclicity part therefore follows from Lemma~\ref{lemma:smash-g-global-tame}.
\end{proof}
\end{thm}

\begin{rk}
The analogue of the above theorem for the positive $G$-global model structure holds and can be proven in the same way. Finally, we note that similar arguments yield various $G$-global level model structures on $\cat{$\bm\Gamma$-$\bm{\mathcal M}$-$\bm{G}$-SSet}_*$, $\cat{$\bm\Gamma$-$\bm{\mathcal M}$-$\bm{G}$-SSet}_*^\tau$, and $\cat{$\bm\Gamma$-$\bm{G}$-$\bm{I}$-SSet}_*$; as they will play no role in the arguments to come, we leave the details to the interested reader.
\end{rk}

\subsubsection{Comparison of level model structures}
Let us now lift some of the equivalences between the models of unstable $G$-global homotopy theory to $\Gamma$-spaces:

\begin{thm}\label{thm:gamma-ev-omega}
The functors
\begin{equation}\label{eq:Gamma-EM-to-Gamma-script-I}
\ev_\omega\colon\cat{$\bm\Gamma$-$\bm G$-$\bm{\mathcal I}$-SSet}_*\rightleftarrows\cat{$\bm\Gamma$-$\bm{E\mathcal M}$-$\bm G$-SSet}_*:\!(\blank)[\omega^\bullet]
\end{equation}
are homotopical, and they descend to mutually inverse equivalences of associated quasi-categories.
\begin{proof}
We first observe that each of the functors
\begin{equation*}
(\blank)[\omega^\bullet]\colon\cat{$\bm{E\mathcal M}$-$\bm{(G\times\Sigma_S)}$-SSet}_*\to\cat{$\bm{(G\times\Sigma_S)}$-$\bm{\mathcal I}$-SSet}_*
\end{equation*}
is homotopical by Proposition~\ref{prop:omega-bullet-inverse}, so that the right hand functor in $(\ref{eq:Gamma-EM-to-Gamma-script-I})$ is homotopical. On the other hand, each of the functors
\begin{equation*}
\ev_\omega\colon\cat{$\bm{(G\times\Sigma_S)}$-$\bm{\mathcal I}$-SSet}_*\to\cat{$\bm{E\mathcal M}$-$\bm{(G\times\Sigma_S)}$-SSet}_*
\end{equation*}
is homotopical by definition, hence so is the left hand functor in $(\ref{eq:Gamma-EM-to-Gamma-script-I})$.

Now the proof of Proposition~\ref{prop:omega-bullet-inverse} constructs a natural transformation from the endofunctor $(\blank)[\omega^\bullet]\circ\ev_\omega$ of $\cat{$\bm{\mathcal I}$-SSet}$ to the identity and shows that for any group $H$ the induced transformation of endofunctors of $\cat{$\bm H$-$\bm{\mathcal I}$-SSet}$ is a levelwise $H$-global weak equivalence; this then provides a levelwise $G$-global level weak equivalence exhibiting $\ev_\omega$ as right inverse to $(\blank)[\omega^\bullet]$ on the level of $\Gamma$-spaces.

Analogously, we have constructed a natural zig-zag between $\ev_\omega\circ(\blank)[\omega^\bullet]$ and the identity of $\cat{$\bm{E\mathcal M}$-SSet}$, and we proved that for any $H$ the induced zig-zag of endofunctors of $\cat{$\bm{E\mathcal M}$-$\bm H$-SSet}$ is a levelwise weak equivalence. We conclude as before that $\ev_\omega$ is left homotopy inverse to $(\blank)[\omega^\bullet]$ on the level of $\Gamma$-spaces.
\end{proof}
\end{thm}

Similarly one deduces from Theorem~\ref{thm:pos-G-global-EM-tau}:

\begin{thm}
The functors
\begin{equation*}
\ev_\omega\colon\cat{$\bm\Gamma$-$\bm G$-$\bm{\mathcal I}$-SSet}_*\rightleftarrows\cat{$\bm\Gamma$-$\bm{E\mathcal M}$-$\bm G$-SSet}_*^\tau:\!(\blank)_\bullet
\end{equation*}
are homotopical, and they descend to mutually inverse equivalences of associated quasi-categories.\qed
\end{thm}

By $2$-out-of-$3$ we conclude from the above two theorems:

\begin{cor}
The inclusion $\cat{$\bm\Gamma$-$\bm{E\mathcal M}$-$\bm G$-SSet}_*^\tau\hookrightarrow\cat{$\bm\Gamma$-$\bm{E\mathcal M}$-$\bm G$-SSet}_*$ induces an equivalence of associated quasi-categories.\qed
\end{cor}

\subsubsection{Connection to equivariant $\Gamma$-spaces}
We now want to clarify the relation between $G$-global $\Gamma$-spaces and the $G$-equivariant $\Gamma$-spaces discussed before. To this end, we introduce the following model structure, which for finite $G$ recovers \cite[Theorem~4.7]{equivariant-gamma}:

\index{G-equivariant model structure@$G$-equivariant model structure!level|seeonly{$G$-equivariant level model structure}}
\begin{prop}\label{prop:Gamma-G-level}\index{G-equivariant level model structure@$G$-equivariant level model structure!on Gamma-G-SSet@on $\cat{$\bm\Gamma$-$\bm G$-SSet}_*$|textbf}
There is a unique model structure on $\cat{$\bm\Gamma$-$\bm G$-SSet}_*$ in which a map $f$ is a weak equivalence or fibration if and only if $f(S_+)$ is a weak equivalence or fibration, respectively, in the $\mathcal G_{\mathcal Fin,\Sigma_S}$-model structure on $\cat{$\bm{(G\times\Sigma_S)}$-SSet}$. We call this the \emph{$G$-equivariant level model structure} and its weak equivalences the \emph{$G$-equivariant level weak equivalences}\index{G-equivariant level weak equivalence@$G$-equivariant level weak equivalence!in Gamma-G-SSet@in $\cat{$\bm\Gamma$-$\bm G$-SSet}_*$}. It is simplicial, proper, and combinatorial with generating cofibrations
\begin{equation*}
\{(\Gamma(S_+,\blank)\smashp G_+)/K\smashp (\del\Delta^n\hookrightarrow\Delta^n)_+ : \text{$K\subset G$ finite, $S$ finite $K$-set}\}.
\end{equation*}
Finally, the $G$-equivariant level weak equivalences are stable under filtered colimits.
\begin{proof}
As before it suffices that
\begin{equation*}
\Gamma(A_+,B_+)\smashp_{\Sigma_A}\blank\colon\cat{$\bm{(G\times\Sigma_A)}$-SSet}_*\to\cat{$\bm{(G\times\Sigma_B)}$-SSet}_*
\end{equation*}
sends the standard generating (acyclic) cofibrations of the $\mathcal G_{\mathcal Fin,\Sigma_A}$-model structure to (acyclic) cofibrations in the injective $\mathcal G_{\mathcal Fin,\Sigma_B}$-model structure, which is clear by direct inspection.
\end{proof}
\end{prop}

The above condition is equivalent to demanding that $f(S_+)$ be an $H$-equivariant weak equivalence or fibration for any finite subgroup $H\subset G$ and any finite $H$-set $S$. We caution the reader that unlike for finite $G$ this is stronger than demanding that $f(S_+)$ be a $\mathcal Fin$-weak equivalence or fibration for every finite $G$-set $S$:

\begin{ex}
Let $f\colon X\to Y$ be a level weak equivalence of $\Gamma$-spaces that is not a $\mathbb Z/n$-level weak equivalence (with respect to the trivial $\mathbb Z/n$-action) for some $n\ge2$, see e.g.~Example~\ref{ex:non-equiv-but-not-equiv-level-we}. Then $f$ is not a $\mathbb Q/\mathbb Z$-level weak equivalence either as $\mathbb Z/n$ embeds into $\mathbb Q/\mathbb Z$. However, $f(S_+)$ is a $\mathbb Q/\mathbb Z$-weak equivalence for every finite $\mathbb Q/\mathbb Z$-set $S$, as $\mathbb Q/\mathbb Z$ admits no non-trivial actions on finite sets.
\end{ex}

We will now discuss some homotopical properties of the passage from $G$-global $\Gamma$-spaces to $H$-equivariant $\Gamma$-spaces along a group homomorphism $\phi\colon H\to G$, which will become crucial later in the comparison between $G$-global $\Gamma$-spaces and $G$-global spectra. While similar results can be achieved for the other models introduced above, we will restrict to the approach via $\mathcal I$-simplicial sets.

\begin{defi}
Let $H$ be any group (not necessarily finite), let $\mathcal U_H$ as in $(\ref{eq:proper-set-universe})$, and let $\phi\colon H\to G$ be any homomorphism. Then the \emph{$\phi$-underlying equivariant $\Gamma$-space}\nomenclature[auphi]{$\und_\phi$ (also $\und_H$)}{underlying $H$-equivariant $\Gamma$-space of a $G$-global $\Gamma$-space for $\phi\colon H\to G$ or $H\subset G$, respectively} $\und_\phi(X)$ of $X\in\cat{$\bm\Gamma$-$\bm{G}$-$\bm{\mathcal I}$-SSet}_*$ is $(\phi^*X)(\mathcal U_H)$. If $f\colon X\to Y$ is any map in $\cat{$\bm\Gamma$-$\bm{G}$-$\bm{\mathcal I}$-SSet}_*$, then we define $\und_\phi(f)\mathrel{:=}(\phi^*f)(\mathcal U_H)$.

If $H\subset G$ and $\phi$ is the inclusion, then we abbreviate $\und_H\mathrel{:=}\und_\phi$.
\end{defi}

\begin{lemma}\label{lemma:gamma-I-vs-equivariant}\index{G-global level weak equivalence@$G$-global level weak equivalence!in Gamma-G-I-SSet@in $\cat{$\bm\Gamma$-$\bm G$-$\bm{\mathcal I}$-SSet}_*$}
The following are equivalent for any map $f$ in $\cat{$\bm\Gamma$-$\bm{G}$-$\bm{\mathcal I}$-SSet}_*$:
\begin{enumerate}
\item $f$ is a $G$-global level weak equivalence.
\item $\und_\phi(f)$ is an $H$-equivariant level weak equivalence for all groups $H$ and all homomorphisms $\phi\colon H\to G$.
\item $\und_\phi(f)$ is an $H$-equivariant level weak equivalence for all \emph{finite} groups $H$ and all homomorphisms $\phi\colon H\to G$
\end{enumerate}
\begin{proof}
Clearly $(2)\Rightarrow(3)$; we will show that also $(1)\Rightarrow(2)$ and $(3)\Rightarrow(1)$.

For $(1)\Rightarrow(2)$, let $H$ be any group and $\phi\colon H\to G$ any homomorphism. By definition, $f(S_+)$ is a $(G\times\Sigma_S)$-global weak equivalence for any finite set $S$; in particular, if $K\subset H$ is any finite subgroup, then $(\phi^*f)(S_+)$ is a $K$-global weak equivalence for any $K$-action on $S$. As $\mathcal U_H$ (viewed as a $K$-set) contains a complete $K$-set universe (Lemma~\ref{lemma:proper-set-universe}), Lemma~\ref{lemma:evaluation-h-universe} implies that $(\und_\phi f)(S_+)=\ev_{\mathcal U_H}(\phi^*f(S_+))$ is a $K$-equivariant weak equivalence. Letting $K$ and $S$ vary, this precisely means that $\und_\phi f$ is an $H$-equivariant level weak equivalence.

For $(3)\Rightarrow(1)$, let $H\subset\mathcal M$ be any universal subgroup together with homomorphisms $\phi\colon H\to G$, $\rho\colon H\to\Sigma_S$; we have to show that $f(S_+)(\omega)^{(\phi,\rho)}$ is a weak homotopy equivalence. But as both $\mathcal U_H$ as well as $\omega$ are complete $H$-set universes, $f(S_+)(\omega)$ agrees up to conjugation by $(\Sigma_S\times H)$-equivariant isomorphisms with $f(S_+)(\mathcal U_H)=(\und_\phi f)(S_+)$. Thus, $f(S_+)(\omega)^{(\phi,\rho)}$ is conjugate to the weak homotopy equivalence $(\und_\phi f)(\rho^* S_+)^H$, hence a weak homotopy equivalence itself.
\end{proof}
\end{lemma}

We now want to give a simpler construction of $\und_H$ under suitable fibrancy conditions. For this we will need:

\begin{cor}\label{cor:Gamma-I-injective}\index{G-global level model structure@$G$-global level model structure!injective|seeonly{injective $G$-global level model structure}}\index{injective G-global level model structure@injective $G$-global level model structure!on Gamma-G-II-SSet@on $\cat{$\bm\Gamma$-$\bm G$-$\bm{\mathcal I}$-SSet}_*$|textbf}\index{G-global model structure@$G$-global model structure!injective level|seeonly{injective $G$-global level model structrue}}
There is a unique model structure on $\cat{$\bm\Gamma$-$\bm G$-$\bm{\mathcal I}$-SSet}_*$ with the injective cofibrations as cofibrations and the $G$-global level weak equivalences as weak equivalences.

We call this the \emph{injective $G$-global level model structure}. It is combinatorial, simplicial, proper, and filtered colimits in it are homotopical.
\begin{proof}
As pushouts along injective cofibrations preserve weak equivalences, Corollary~\ref{cor:mix-model-structures} shows that the model structure exists, that it is proper, and that filtered colimits in it are homotopical. Finally, the Pushout Product Axiom for the simplicial tensoring follows by applying Theorem~\ref{thm:script-I-vs-EM-injective} levelwise (with $G$ replaced by $G\times\Sigma_S$ for varying finite sets $S$).
\end{proof}
\end{cor}

\begin{warn}
In general, fibrations in $\cat{$\bm\Gamma$-$\bm G$-$\bm{\mathcal I}$-SSet}_*$ need not be levelwise injective fibrations, i.e.~$\ev_{T_+}\colon\cat{$\bm\Gamma$-$\bm G$-$\bm{\mathcal I}$-SSet}_*\to\cat{$\bm{(\Sigma_T\times G)}$-$\bm{\mathcal I}$-SSet}_*$ need not be right Quillen with respect to the injective $G$-global level and injective $(\Sigma_T\times G)$-global model structures. Equivalently, the left adjoint $G_{T_+}$ need not be homotopical as the following example for $G=1$ shows: we let $T$ be any finite set with $|T|\ge 2$. If $G_{T_+}$ were homotopical, then the same would be true for
\begin{equation}\label{eq:sigma-t-warn}
\Gamma(T_+,1^+)\smashp_{\Sigma_T}\blank\cong\ev_{1^+}\circ G_{T_+}\colon\cat{$\bm{\Sigma_T}$-$\bm{\mathcal I}$-SSet}_*\to\cat{$\bm{\mathcal I}$-SSet}_*.
\end{equation}
However, the map $f\colon T_+\to1^+$ with $f(t)=1$ for all $t\in T$ is $\Sigma_T$-fixed, and so is the map sending everything to the base point. We conclude that $S^0$ with trivial $\Sigma_T$-action is a $\Sigma_T$-equivariant retract of $\Gamma(T_+,1^+)$, so that $(\blank)/\Sigma_T\cong S^0\smashp_{\Sigma_T}\blank$ is a retract of $(\ref{eq:sigma-t-warn})$; in particular, it would have to be homotopical.

But we have seen in the proof of Proposition~\ref{prop:injective-is-really-static-script-I} that $\mathcal I(T,\blank)_+\to *_+$ (where $\Sigma_T$ acts on the left hand side via its right action on $T$) is a $\Sigma_T$-global weak equivalence, and we claim that $(\mathcal I(T,\blank)/\Sigma_T)_+\to *_+$ is not a global weak equivalence. Indeed, it suffices that $(E\Inj(T,\omega)/\Sigma_T)\to *$ is not a global weak equivalence of $E\mathcal M$-simplicial sets, but this isn't even an underlying weak equivalence as $E\Inj(T,\omega)/\Sigma_T$ is a $K(\Sigma_T,1)$ and hence in particular not weakly contractible.
\end{warn}

However, we at least have:

\begin{prop}\label{prop:Gamma-evaluation-injective}
Let $H\subset G$ and let $S$ be a finite $H$-set. Then
\begin{equation}\label{eq:evaluation-G-set-injective}
\ev_{S_+}\colon(\cat{$\bm\Gamma$-$\bm G$-$\bm{\mathcal I}$-SSet}_*)_{\textup{injective $G$-global level}}\to(\cat{$\bm H$-$\bm{\mathcal I}$-SSet}_*)_{\textup{injective $H$-global}}
\end{equation}
is right Quillen. Here $X(S_+)$ for $X\in\cat{$\bm G$-$\bm\Gamma$-$\bm{\mathcal I}$-SSet}_*$ is as usual equipped with the diagonal of the $H$-action on $S$ and the restriction of the $G$-action on $X$.
\begin{proof}
We first observe that in the adjunction
\begin{equation*}
G_+\smashp_H\blank\colon\cat{$\bm\Gamma$-$\bm H$-$\bm{\mathcal I}$-SSet}_*\rightleftarrows\cat{$\bm\Gamma$-$\bm G$-$\bm{\mathcal I}$-SSet}_* :\!\res^G_H
\end{equation*}
(where $\res^G_H$ denotes restriction from $G$ to $H$)\nomenclature[aresGH]{$\res^G_H$}{restriction of $G$-action to $H$ for $H\subset G$}
the left adjoint obviously preserves injective cofibrations and that it sends $H$-global level weak equivalences to $G$-global level weak equivalences by Lemma~\ref{lemma:alpha-shriek-injective-script-I} (applied to $H\times\Sigma_S\hookrightarrow G\times\Sigma_S$ for varying finite set $S$). In particular, $\res^G_H$ is right Quillen with respect to the injective level model structures, and we may therefore assume without loss of generality that $H=G$.

Let $\rho\colon G\to\Sigma_S$ be the homomorphism classifying the $G$-action on $S$. Then $(\ref{eq:evaluation-G-set-injective})$ factors as
\begin{equation*}
\cat{$\bm\Gamma$-$\bm G$-$\bm{\mathcal I}$-SSet}_*\xrightarrow{\ev_{S_+}}\cat{$\bm{(G\times\Sigma_S)}$-$\bm{\mathcal I}$-SSet}_*\xrightarrow{(\id,\rho)^*}\cat{$\bm G$-$\bm{\mathcal I}$-SSet}_*,
\end{equation*}
so a left adjoint is given by $G_{S_+}\circ(\id,\rho)_!$. This obviously preserves injective cofibrations, so it suffices to prove that it is homotopical, i.e.~if $f$ is a weak equivalence of pointed $G$-$\mathcal I$-simplicial sets, then
\begin{equation}\label{eq:egsi-left-adjoint-level}
(G_{S_+}(\id,\rho)_!f)(T_+)=\Gamma(S_+,T_+)\smashp_{\Sigma_S}(G\times\Sigma_S)_+\smashp_G f
\end{equation}
(where $G$ acts on both $G$ and $\Sigma_S$ from the right) is a $(G\times\Sigma_T)$-global weak equivalence for every finite set $T$.

Let us write $\Gamma(S_+,T_+)^{\textup{conj}}$ for $\Gamma(S_+,T_+)$ with $G$ action via $S$ and $\Sigma_T$-action via $T$. Then an analogous calculation to the one from the proof of Theorem~\ref{thm:Gamma-EM-SSet} shows that $(\ref{eq:egsi-left-adjoint-level})$ is conjugate to $\Gamma(S_+,T_+)^{\textup{conj}}\smashp f$. One easily checks that smashing with any pointed $(G\times\Sigma_T)$-simplicial set sends $G$-global weak equivalences to $(G\times\Sigma_T)$-global weak equivalences, which then completes the proof.
\end{proof}
\end{prop}

\begin{cor}\label{cor:G-injective-constant}
Let $X$ be fibrant in the injective $G$-global level model structure, let $H\subset G$, and let $i\colon A\to B$ be an $H$-equivariant injection of $H$-sets. Then $X(\blank)(A)\to X(\blank)(B)$ is an $H$-equivariant level weak equivalence of $\Gamma$-$H$-spaces.
\begin{proof}
Let $H'\subset H$ be finite and let $S$ be a finite $H'$-set. We have to show that $X(S_+)(i)\colon X(S_+)(A)\to X(S_+)(B)$ is an $H'$-equivariant weak equivalence. But by the previous proposition $X(S_+)$ is an injectively fibrant $H'$-$\mathcal I$-simplicial set, so the claim follows from Proposition~\ref{prop:injective-is-really-static-script-I}.
\end{proof}
\end{cor}

\begin{cor}\label{cor:injective-Gamma-varnothing-UG}
Let $X$ be fibrant in the injective $G$-global level model structure on $\cat{$\bm\Gamma$-$\bm G$-$\bm{\mathcal I}$-SSet}_*$ and let $H\subset G$. Then $\varnothing\to\mathcal U_H$ induces an $H$-equivariant level weak equivalence $\res^G_HX(\blank)(\varnothing)\to X(\blank)(\mathcal U_H)=\und_H(X)$.\qed
\end{cor}

Finally, let us lift the comparison between $G$-$\mathcal I$-simplicial sets and $G$-simplicial sets established in Subsection~\ref{subsec:G-global-vs-G-equiv-script-I} to the level of $\Gamma$-spaces:

\begin{prop}\label{prop:G-global-Gamma-vs-Gamma-G}\index{proper G-equivariant homotopy theory@proper $G$-equivariant homotopy theory!vs G-global homotopy theory@vs.~$G$-global homotopy theory!for Gamma-spaces@for $\Gamma$-spaces}
The simplicial adjunctions
\begin{align}\label{eq:ev-UG-Gamma}
\und_G=\ev_{\mathcal U_G}\colon\cat{$\bm\Gamma$-$\bm G$-$\bm{\mathcal I}$-SSet}_*&\rightleftarrows\cat{$\bm\Gamma$-$\bm G$-SSet}_* :\!R\\
\label{eq:ev-varnothing-Gamma}
\const\colon\cat{$\bm\Gamma$-$\bm G$-SSet}_*&\rightleftarrows(\cat{$\bm\Gamma$-$\bm G$-$\bm{\mathcal I}$-SSet}_*)_{\textup{injective}} :\!\ev_\varnothing
\end{align}
are Quillen adjunctions, $\ev_{\mathcal U_G}$ is homotopical, and $\ev_{\mathcal U_G}^\infty\simeq\cat{R}\ev_\varnothing$. Moreover, the induced adjunctions $\ev_{\mathcal U_G}^\infty\dashv\cat{R}R$ and  $\cat{L}\const\dashv\cat{R}\ev_\varnothing$ are a left and right Bousfield localization, respectively, with respect to the $\und_G$-weak equivalences (i.e.~those maps sent to $G$-equivariant level weak equivalences under $\und_G$).
\end{prop}

For the proof we will need:

\begin{lemma}\label{lemma:equivariant-Gamma-generalized-cof}
Let $H$ be a group, $S$ a finite $H$-set, and $i\colon X\to Y$ a cofibration in the $\mathcal G_{\mathcal Fin,H}$-model structure on $\cat{$\bm{(G\times H)}$-SSet}_*$. Then $(\Gamma(S_+,\blank)\smashp i)/H$ is a cofibration in the $G$-equivariant level model structure on $\cat{$\bm\Gamma$-$\bm G$-SSet}_*$.
\begin{proof}
As $(\Gamma(S_+,\blank)\smashp\blank)/H$ is cocontinuous and preserves tensors, it suffices to prove that $(\Gamma(S_+,\blank)\smashp (G\times H)_+/\Gamma_{K,\phi})/H$ is cofibrant for any finite group $K\subset G$ and any homomorphism $\phi\colon K\to H$. But this is isomorphic to $(\Gamma(\phi^*S_+,\blank)\smashp G_+)/K$ where $K$ acts on $S$ and on $G$ in the obvious way. The claim follows immediately from the description of the generating cofibrations given in Proposition~\ref{prop:Gamma-G-level}.
\end{proof}
\end{lemma}

\begin{proof}[Proof of Proposition~\ref{prop:G-global-Gamma-vs-Gamma-G}]
Lemma~\ref{lemma:gamma-I-vs-equivariant} shows that $\ev_{\mathcal U_G}$ is homotopical. To prove that it is left Quillen, it is then enough by the previous lemma and the explicit description of the generating cofibrations that $E\Inj(A,\mathcal U_G)\times G$ is $\mathcal G_{\mathcal Fin,H}$-cofibrant for every finite faithful $H$-set $A$ and $H$ acting from the right on $G$ via any homomorphism $\phi$. But indeed, $H$ acts freely on the first factor, so every isotropy group belongs to $\mathcal G_{G,H}$; on the other hand, $G$ acts freely on the second factor, hence every isotropy group belongs to $\mathcal G_{H,G}$, so it is in particular finite, as desired.

To prove that $(\ref{eq:ev-varnothing-Gamma})$ is a Quillen adjunction, too, let $f\colon X\to Y$ be a fibration or acyclic fibration in $\cat{$\bm\Gamma$-$\bm G$-$\bm{\mathcal I}$-SSet}_*$; we have to show that $f(S_+)(\varnothing)$ is a fibration or acyclic fibration, respectively, in $\cat{$\bm H$-SSet}_*$ for each $H\subset G$ and each finite $H$-set $S$. But $f(S_+)$ is a fibration or acyclic fibration in the injective $H$-global model structure by Proposition~\ref{prop:Gamma-evaluation-injective}, so the claim follows from the proof of Corollary~\ref{cor:G-script-I-SSet-vs-G-SSet}.

Next, we observe that the inclusion $\varnothing\hookrightarrow\mathcal U_G$ induces a natural transformation $\iota\colon\ev_{\varnothing}\Rightarrow\ev_{\mathcal U_G}$, which is a $G$-equivariant level weak equivalence on injectively fibrant $G$-global $\Gamma$-spaces by Corollary~\ref{cor:injective-Gamma-varnothing-UG}. As $\ev_{\mathcal U_G}$ is homotopical, this yields the desired equivalence $(\ev_{\mathcal U_G})^\infty\simeq\cat{R}\ev_\varnothing$.

It only remains to show that the derived unit $X\to \cat{R}\ev_\varnothing(\const X)$ is a $G$-equivariant level weak equivalence for every (cofibrant) equivariant $\Gamma$-space $X$. For this we fix an injectively fibrant replacement $\kappa\colon\const X\to Y$ and consider the commutative diagram
\begin{equation*}
\begin{tikzcd}
X\arrow[r,"\eta", "\cong"'] & \ev_{\varnothing}\const X \arrow[r, "\iota","\cong"']\arrow[d,"\ev_\varnothing\kappa"']& \ev_{\mathcal U_G}\const X\arrow[d, "\ev_{\mathcal U_G}\kappa", "\sim"']\\
& \ev_\varnothing Y\arrow[r, "\sim","\iota"'] & \ev_{\mathcal U_G}Y
\end{tikzcd}
\end{equation*}
The composition $X\to\ev_\varnothing Y$ represents the derived unit. But the top horizontal arrows are isomorphisms by direct inspection, the lower horizontal arrow is a $G$-equivariant weak equivalence by the above, and the right hand vertical map is a $G$-equivariant weak equivalence as $\ev_{\mathcal U_G}$ is homotopical. The claim follows by $2$-out-of-$3$.
\end{proof}

\begin{warn}
The functor $\const\colon\cat{$\bm\Gamma$-$\bm G$-SSet}_*\to\cat{$\bm\Gamma$-$\bm G$-$\bm{\mathcal I}$-SSet}_*$ is not homotopical with respect to the $G$-global level weak equivalences on the target. However, it is indeed homotopical with respect to the $\und_G$-weak equivalences.
\end{warn}

\subsubsection{Functoriality}
\index{functoriality in homomorphisms!for Gamma-EM-G-SSettau@for $\cat{$\bm\Gamma$-$\bm{E\mathcal M}$-$\bm G$-SSet}^\tau_*$|(}
The results on the change-of-group adjunctions for the models of (pointed) unstable $G$-global homotopy theory discussed in Chapter~\ref{chapter:unstable} easily transfer to statements about the corresponding adjunctions on the level of $\Gamma$-spaces. Instead of making all of these explicit, we will only collect the results here that we will need below:

\begin{cor}\label{cor:gamma-tau-alpha-lower-shriek}
Let $\alpha\colon H\to G$ be any group homomorphism. Then the adjunction
\begin{equation*}
\alpha_!\colon\cat{$\bm\Gamma$-$\bm{E\mathcal M}$-$\bm H$-SSet}_*^\tau\rightleftarrows\cat{$\bm\Gamma$-$\bm{E\mathcal M}$-$\bm G$-SSet}_*^\tau :\!\alpha^*
\end{equation*}
is a Quillen adjunction with fully homotopical right adjoint. Moreover, if $\alpha$ is injective, then also the left adjoint is homotopical.
\begin{proof}
As weak equivalences and fibrations are defined levelwise, this follows by applying Lemma~\ref{lemma:alpha-lower-shriek-tau} to $\alpha\times\Sigma_S$ for varying finite set $S$.
\end{proof}
\end{cor}

Similarly one deduces from Lemma~\ref{lemma:alpha-lower-star-tame}:

\begin{cor}\label{cor:gamma-tau-alpha-lower-star}
Let $\alpha\colon H\to G$ be an injective homomorphism such that $(G:\im\alpha)<\infty$. Then the adjunction
\begin{equation*}
\alpha^*\colon\cat{$\bm\Gamma$-$\bm{E\mathcal M}$-$\bm G$-SSet}_*^\tau\rightleftarrows\cat{$\bm\Gamma$-$\bm{E\mathcal M}$-$\bm H$-SSet}_*^\tau :\!\alpha_*
\end{equation*}
is a Quillen adjunction in which both functors are homotopical.\index{functoriality in homomorphisms!for Gamma-EM-G-SSettau@for $\cat{$\bm\Gamma$-$\bm{E\mathcal M}$-$\bm G$-SSet}^\tau_*$|)}\qed
\end{cor}

\subsection{Special $\bm G$-global $\bm\Gamma$-spaces}
\index{special!G-global@$G$-global|seeonly{$G$-global $\Gamma$-space, special}}
Let $X$ be a $G$-global $\Gamma$-space; just as in the classical equivariant or non-equivariant setting, we want to think of the maps $X(S_+)\to X(1^+)$ induced by the unique maps $S\to\{1\}$ as `generalized multiplications.' However, in order for this intuition to apply, we again have to put suitable conditions on $X$ first:

\begin{defi}\index{G-global Gamma-space@$G$-global $\Gamma$-space!special|textbf}
We call $X\in\cat{$\bm\Gamma$-$\bm{E\mathcal M}$-$\bm G$-SSet}_*$ \emph{special} if the following holds: for all finite sets $S$ the Segal map
\begin{equation*}
X(S_+)\xrightarrow{\rho_S\mathrel{:=}X(p_s)_{s\in S}}\prod_{s\in S} X(1^+)
\end{equation*}
is a $(G\times\Sigma_S)$-global weak equivalence, where $\Sigma_S$ acts on both sides via its tautological action on $S$, and analogously for $X\in\cat{$\bm\Gamma$-$\bm{E\mathcal M}$-$\bm G$-SSet}_*^\tau$ or $X\in\cat{$\bm\Gamma$-$\bm G$-$\bm{\mathcal I}$-SSet}_*$.

We denote the corresponding full subcategories by the superscript `special.'
\end{defi}

One easily proves by direct inspection similarly to Lemma~\ref{lemma:gamma-I-vs-equivariant}:

\begin{lemma}\label{lemma:special-g-global-u-phi}
A $G$-global $\Gamma$-space $X\in\cat{$\bm\Gamma$-$\bm G$-$\bm{\mathcal I}$-SSet}_*$ is special if and only if $\und_\phi X$ is special in the sense of Definition~\ref{defi:equivariant-special} for all finite groups $H$ and all homomorphisms $\phi\colon H\to G$.\qed
\end{lemma}

\begin{ex}\label{ex:shimada-shimakawa-g-global}
In Example~\ref{ex:shimada-shimakawa-equivariant} we have seen how one can assign a special $\Gamma$-$G$-space to a small symmetric monoidal category $\mathscr C$ with $G$-action through strictly unital strong symmetric monoidal functors. By a slight variation, we can actually build a special \emph{$G$-global} $\Gamma$-space from the same data. For this, let us consider the functor $\Fun(E\mathcal M,\blank)\colon\cat{$\bm G$-Cat}\to\cat{$\bm{E\mathcal M}$-$\bm G$-Cat}$, where $E\mathcal M$ acts on itself from the right via precomposition. We claim that $\nerve\Fun(E\mathcal M,\Gamma(\mathscr C))$ is a special $G$-global $\Gamma$-space for any $\mathscr C$ as above.

Indeed, if $H\subset\mathcal M$ is any subgroup, then $H$ acts freely from the right on $\mathcal M$, so there exists a right $H$-equivariant map $r\colon\mathcal M\to H$. It is then trivial to check that $EH\hookrightarrow E\mathcal M$ is an $H$-equivariant equivalence, i.e.~an equivalence in the $2$-category of right $H$-categories, $H$-equivariant functors, and $H$-equivariant natural transformations: a quasi-inverse is given by $Er$, and the unique maps $h\to ri(h)$ in $EH$ assemble into an $H$-equivariant isomorphism $\id\cong (Er)(Ei)$. Analogously, there is a (unique) $H$-equivariant isomorphism $\id\cong (Ei)(Er)$.

It follows that the functor $\Fun(E\mathcal M,\blank)^H$ is just naturally equivalent to the categorical homotopy fixed point functor $\Fun(EH,\blank)^H$; in particular, it sends underlying equivalences of $H$-categories to equivalences of categories. Restricting along all homomorphisms $H\to G\times\Sigma_S$, we therefore conclude that $\nerve\Fun(E\mathcal M,\blank)$ sends underlying equivalences of $(G\times\Sigma_S)$-categories to $(G\times\Sigma_S)$-global weak equivalences in $\cat{$\bm{E\mathcal M}$-$\bm{(G\times\Sigma_S)}$-SSet}$, and it follows as in Example~\ref{ex:shimada-shimakawa-equivariant} that $\nerve\Fun(E\mathcal M,\Gamma(\mathscr C))$ is $G$-globally special in the above sense.
\end{ex}

The comparisons from the previous subsection carry over to the context of special $G$-global $\Gamma$-spaces:

\begin{cor}\label{cor:comparison-special}
All the functors in the diagram
\begin{equation*}
\begin{tikzcd}[column sep=-20pt]
\cat{$\bm\Gamma$-$\bm{E\mathcal M}$-$\bm G$-SSet}^{\tau,\textup{special}}_*\arrow[rr, hook] && \cat{$\bm\Gamma$-$\bm{E\mathcal M}$-$\bm G$-SSet}^{\textup{special}}_*\arrow[ld, "{(\blank)[\omega^\bullet]}", bend left=10pt]\\
& \cat{$\bm\Gamma$-$\bm G$-$\bm{\mathcal I}$-SSet}_*^{\textup{special}}\arrow[ul, "\ev_\omega", bend left=10pt]
\end{tikzcd}
\end{equation*}
are homotopical, and they induce equivalencess of associated quasi-categories. Moreover, the resulting diagram commutes up to preferred equivalence.
\begin{proof}
As specialness is obviously invariant under $G$-global level weak equivalences, Proposition~\ref{prop:localization-subcategory} implies that for all of the above models the inclusion of the full subcategory of special $G$-global $\Gamma$-spaces descends to a fully faithful functor of associated quasi-categories.

Next we observe that all of the functors are defined levelwise. The corresponding functors in each level preserve and reflect $(G\times\Sigma_S)$-global weak equivalences (as they are homotopical and induce equivalences by the results of Chapter~\ref{chapter:unstable}), and they preserve finite products by direct inspection. We therefore conclude that they preserve and reflect specialness.

Thus, the comparisons of $G$-global $\Gamma$-spaces from the previous section imply that all of the above functors descend to equivalences of quasi-localizations. It only remains to show the commutativity up to equivalence, which now follows easily from Theorem~\ref{thm:gamma-ev-omega}.
\end{proof}
\end{cor}

\subsubsection{A model categorical manifestation} For all of our above model categories of $G$-global $\Gamma$-spaces, one can obtain a model of \emph{special} $G$-global $\Gamma$-spaces via Bousfield localization. We will make this explicit for $\cat{$\bm\Gamma$-$\bm{E\mathcal M}$-$\bm G$-SSet}_*^\tau$:

\begin{defi}\index{G-global Gamma-space@$G$-global $\Gamma$-space!special weak equivalence|seeonly{$G$-global special weak equivalence}}
\index{G-global special weak equivalence@$G$-global special weak equivalence!in Gamma-EM-G-SSettau@in $\cat{$\bm\Gamma$-$\bm{E\mathcal M}$-$\bm G$-SSet}_*^\tau$|textbf}
A morphism $f\colon X\to Y$ in $\cat{$\bm\Gamma$-$\bm{E\mathcal M}$-$\bm G$-SSet}_*^\tau$ is called a \emph{$G$-global special weak equivalence} if the induced map $f^*\colon[Y,T]\to[X,T]$ is bijective for all special $G$-global $\Gamma$-spaces $T$. Here $[\,{,}\,]$ denotes the hom sets in the homotopy category of $\cat{$\bm\Gamma$-$\bm{E\mathcal M}$-$\bm G$-SSet}_*^\tau$ with respect to the $G$-global level weak equivalences.
\end{defi}

\begin{thm}\label{thm:special-model-structure}\index{G-global Gamma-space@$G$-global $\Gamma$-space!special model structure|seeonly{$G$-global special model structure}}
\index{G-global special model structure@$G$-global special model structure!on Gamma-EM-G-SSettau@on $\cat{$\bm\Gamma$-$\bm{E\mathcal M}$-$\bm G$-SSet}_*^\tau$|textbf}\index{G-global model structure@$G$-global model structure!special|seeonly{$G$-global special model structure}}
There exists a unique model structure on $\cat{$\bm\Gamma$-$\bm{E\mathcal M}$-$\bm G$-SSet}_*^\tau$ with the same cofibrations as the positive $G$-global level model structure and whose weak equivalences are the $G$-global special weak equivalences. Its fibrant objects are precisely the positively level fibrant \emph{special} $G$-global $\Gamma$-spaces, and we call this the \emph{$G$-global special model structure}.

It is left proper, simplicial, and filtered colimits in it are homotopical. Moreover, it is combinatorial and there exists a set $J$ of generating acyclic cofibrations consisting only of maps between cofibrant objects.
\begin{proof}
For the existence of the model structure, and the proof that it is simplicial, left proper, and combinatorial, it suffices by Theorem~\ref{thm:local-model-structure} that there exists a set $T$ of maps between cofibrant objects such that the special $G$-global $\Gamma$-spaces are precisely the \emph{$T$-local} objects $Z$, i.e.~those such that the induced map $f^*\colon\Maps(Y,Z')\to\Maps(X,Z')$ is a weak homotopy equivalence of simplicial sets for all $f\colon X\to Y$ in $T$ and some (hence any) choice of fibrant replacement $Z'$ of $Z$ in the positive $G$-global level model structure.\index{local object}

For this we let $H$ be a universal group of $\mathcal M$, $S$ a finite set, and $\phi\colon H\to G$, $\psi\colon H\to \Sigma_S$ any group homorphisms. We now fix a free $H$-orbit $F\subset\omega$. Then $\big(\Gamma(S_+,\blank)\smashp (E\Inj(F,\omega)\times G)_+\big)/H$ (where $H$ acts on $G$ from the right via $\phi$, on $S$ via $\psi$, and on $F$ in the tautological way) corepresents $X\mapsto X(S_+)_{[F]}^{(\phi,\psi)}$ in the simplicially enriched sense. Here $H$ acts via $\phi$, $\psi$ and the restriction of the $E\mathcal M$-action to $H$ (which preserves simplices supported on the $H$-set $F$). Similarly, $\left(S_+\smashp\Gamma(1^+,\blank)\smashp(E\Inj(F,\omega)\times G)_+\right)/H$ corepresents $X\mapsto \big(\prod_{s\in S}X(1_+)\big)^{(\phi,\psi)}_{[F]}$. By the Yoneda Lemma we therefore get a map
\begin{equation*}
\lambda_{H,S,\phi,\psi}\hskip-1pt\colon\hskip0pt minus 2pt\big(S_+\hskip0pt minus 2pt\smashp\hskip0pt minus 2pt \Gamma(1^+,\blank)\hskip0pt minus 2pt\smashp\hskip0pt minus 2pt(E\Inj(F,\omega)\hskip0pt minus 2pt\times\hskip0pt minus 2pt G)_+\big)/H\hskip-1pt minus 2pt\to\hskip-1pt minus 2pt\big(\Gamma(S_+,\blank)\hskip0pt minus 2pt\smashp\hskip0pt minus 2pt (E\Inj(F,\omega)\hskip0pt minus 2pt\times\hskip0pt minus 2pt G)_+\big)/H
\end{equation*}
such that for any $Z$ the restriction $\Maps(\lambda_{H,S,\phi,\psi},Z)$ is conjugate to
\begin{equation*}
(\rho_S)_{[F]}^{(\phi,\psi)}\colon Z(S_+)_{[F]}^{(\phi,\psi)}\to\Big(\prod\nolimits_{S} Z(1^+)\Big)^{(\phi,\psi)}_{[F]};
\end{equation*}
explicitly, $\lambda_{H,S,\phi,\psi}$ is induced by the map
\begin{equation}\label{eq:segal-corepr}
S_+\smashp\Gamma(1^+,T_+)\to\Gamma(S_+,T_+),[s,f]\mapsto f\circ p_s
\end{equation}
Now assume $Z$ fibrant in the $G$-global positive level model structure. Then $Z(S_+)$ is fibrant in the positive $(G\times\Sigma_S)$-global model structure, so the inclusion induces a weak equivalence $Z(S_+)_{[F]}^{(\phi,\psi)}\hookrightarrow Z(S_+)^{(\phi,\psi)}$ by Remark~\ref{rk:support-em-tau-model-structure}. On the other hand, $\prod_{s\in S} Z(1^+)$ is fibrant in the $(G\times\Sigma_S)$-global model structure by Corollary~\ref{cor:twisted-product}; we conclude by the same argument as before that the inclusion induces a weak equivalence $\left(\prod_{s\in S} Z(1^+)\right)_{[F]}^{(\phi,\psi)}\to\left(\prod_{s\in S} Z(1^+)\right)^{(\phi,\psi)}$. Altogether we see that $\rho_S^{(\phi,\psi)}$ is a weak equivalence if and only if $\Maps(\lambda_{H,S,\phi,\psi},Z)$ is a weak equivalence. In particular, $\rho_S$ is a $(G\times\Sigma_S)$-global weak equivalence if and only if $\Maps(\lambda_{H,S,\phi,\psi},Z)$ is a weak equivalence for all $\phi$ and $\psi$ as above.

Now we define $T$ to be set of all $\lambda_{H,S,\phi,\psi}$ for all universal $H\subset\mathcal M$, all $S=\{1,\dots,n\}$ for $n\ge 0$, and all homomorphisms $\phi\colon H\to G$, $\psi\colon H\to\Sigma_S$. The above then shows that a positively level fibrant $Z$ is special if and only if $\Maps(f,Z)$ is a weak equivalence for all $f\in T$. On the other hand, the targets of the maps in $S$ are obviously cofibrant, while for the sources we observe that the functor corepresented by them sends acyclic $G$-global cofibrations to acyclic Kan fibrations by Corollary~\ref{cor:twisted-product}, hence in particular to surjections on $0$-simplices. Thus, also the sources are cofibrant, finishing the proof of the existence of the model structure.

Lemma~\ref{lemma:filtered-still-homotopical} immediately implies that filtered colimits in the resulting model structure are homotopical. Finally, if $I$ is the usual set of generating cofibrations of the positive $G$-global level model structure, then $I$ is also a set of generating cofibrations for the new model structure. As $I$ obviously consists of maps between cofibrant objects, \cite[Corollary~2.7]{barwick-tractable} implies that also the set $J$ of generating \emph{acyclic} cofibrations can be chosen to consist of maps between cofibrant objects.
\end{proof}
\end{thm}

\begin{lemma}\label{lemma:we-between-special}
Let $f\colon X\to Y$ be a map in $\cat{$\bm\Gamma$-$\bm{E\mathcal M}$-$\bm G$-SSet}_*^{\tau,\textup{special}}$. Then the following are equivalent:
\begin{enumerate}
\item $f$ is a $G$-global special weak equivalence.
\item $f$ is a $G$-global level weak equivalence.
\item $f(1^+)$ is a $G$-global weak equivalence.
\end{enumerate}
\begin{proof}
It is clear that $(2)\Rightarrow(1)$ and $(2)\Rightarrow(3)$. The implication $(1)\Rightarrow(2)$ is a general fact about Bousfield localizations, using that the special $G$-global $\Gamma$-spaces are precisely the local objects.
Finally, if $S$ is any finite set, then $f(S_+)$ agrees with $\prod_{s\in S}f$ up to conjugation by $(G\times\Sigma_S)$-global weak equivalences. The implication $(3)\Rightarrow(2)$ thus follows from Corollary~\ref{cor:twisted-product}.
\end{proof}
\end{lemma}

For later use we record:

\begin{prop}\label{prop:alpha-shriek-gamma-special-homotopical}\index{functoriality in homomorphisms!for Gamma-EM-G-SSettau@for $\cat{$\bm\Gamma$-$\bm{E\mathcal M}$-$\bm G$-SSet}^\tau_*$}
Let $\alpha\colon H\to G$ be any group homomorphism. Then the simplicial adjunction
\begin{equation*}
\alpha_!\colon(\cat{$\bm\Gamma$-$\bm{E\mathcal M}$-$\bm H$-SSet}_*^\tau)_{\textup{special}}\rightleftarrows(\cat{$\bm\Gamma$-$\bm{E\mathcal M}$-$\bm G$-SSet}^\tau_*)_{\textup{special}} :\!\alpha^*
\end{equation*}
is a Quillen adjunction. If $\alpha$ is injective, then $\alpha_!$ is fully homotopical.
\begin{proof}
For the corresponding level model structures, this was shown as Corollary~\ref{cor:gamma-tau-alpha-lower-shriek}. To prove the first statement, it is therefore enough to show that $\alpha^*$ preserves specialness, which is immediate from Lemma~\ref{lemma:alpha-lower-shriek-tau}.

The second statement then follows from Corollary~\ref{cor:gamma-tau-alpha-lower-shriek} as any $G$-global special weak equivalence factors as a special $G$-global acyclic cofibration followed by a $G$-global level weak equivalence.
\end{proof}
\end{prop}

\begin{lemma}\label{lemma:pushout-injective-infty}
Pushouts along injective cofibrations in $\cat{$\bm\Gamma$-$\bm{E\mathcal M}$-$\bm G$-SSet}_*^\tau$ are homotopy pushouts in the $G$-global special model structure.
\begin{proof}
As the injective cofibrations are closed under pushout, it suffices that pushouts of $G$-global special weak equivalences along injective cofibrations are again $G$-global special weak equivalences, see e.g.~\cite[Proposition~1.6]{batanin-berger}.

Indeed, the $G$-global \emph{level} weak equivalences are stable under pushouts along injective cofibrations by Theorem~\ref{thm:Gamma-EM-tau}, and the general case then follows as in the proof of Theorem~\ref{thm:G-M-injective-semistable-model-structure}.
\end{proof}
\end{lemma}

\subsubsection{$G$-global vs.~$G$-equivariant specialness}\index{Gamma-G-space@$\Gamma$-$G$-space!special}
Let us call a $\Gamma$-$G$-space for a general discrete group $G$ \emph{special}, if it is special as a $\Gamma$-$H$-space for all finite subgroups $H\subset G$. We now want to use the comparison between $G$-global and $G$-equivariant $\Gamma$-spaces from Proposition~\ref{prop:G-global-Gamma-vs-Gamma-G} to prove:

\begin{thm}\label{thm:special-G-global-Gamma-vs-Gamma-G}\index{proper G-equivariant homotopy theory@proper $G$-equivariant homotopy theory!vs G-global homotopy theory@vs.~$G$-global homotopy theory!for Gamma-spaces@for $\Gamma$-spaces|textbf}
The homotopical functor $\ev_{\mathcal U_G}\hskip-1pt minus 1pt\colon\hskip0pt minus 2pt\cat{$\bm\Gamma$\kern-.5pt-$\bm G$-$\bm{\mathcal I}$\kern-.5pt-SSet}_*\hskip0pt minus 2pt\to\hskip0pt minus 2pt\cat{$\bm\Gamma$\kern-.5pt-$\bm G$-SSet}_*$ induces a Bousfield localization $(\cat{$\bm\Gamma$-$\bm G$-$\bm{\mathcal I}$-SSet}_*^{\textup{special}})^\infty\to(\cat{$\bm\Gamma$-$\bm G$-SSet}_*^{\textup{special}})^\infty$.
\end{thm}

The proof of this is slightly more involved because being $G$-globally special is a much stronger condition than being $G$-equivariantly special; in particular, it turns out that $\const\colon \cat{$\bm\Gamma$-$\bm G$-SSet}_*\to\cat{$\bm\Gamma$-$\bm G$-$\bm{\mathcal I}$-SSet}_*$ does \emph{not} preserve specialness. Instead, we will have to consider the \emph{right} adjoint $\cat{R}R$ of $\ev_{\mathcal U_G}^\infty$, but proving that this has the desired properties requires some preparations.

We begin with the following analogue of Theorem~\ref{thm:special-model-structure}:

\begin{prop}\label{prop:equivariant-special-model-structure}
There is a unique model structure on $\cat{$\bm\Gamma$-$\bm G$-SSet}_*$ with the same cofibrations as the $G$-equivariant level model structure and with fibrant objects the $G$-equivariantly level fibrant \emph{special} $\Gamma$-$G$-spaces. This model structure is combinatorial, simplicial, left proper, and filtered colimits in it are homotopical.
\begin{proof}
This is proven in precisely the same way as Theorem~\ref{thm:special-model-structure} by localizing with respect to the maps
\begin{equation}\label{eq:generating-cof-bfl-Gamma-G}
(S_+\smashp\Gamma(1_+,\blank)\smashp G_+)/H\to (\Gamma(S_+,\blank)\smashp G_+)/H
\end{equation}
induced by $(\ref{eq:segal-corepr})$ for every finite $H\subset G$ and every finite $H$-set $S$ (up to isomorphism); observe that the sources of these maps are indeed cofibrant as they corepresent $X\mapsto\Maps(S_+,X(1^+))^H$, which is isomorphic to $X\mapsto \prod_{i=1}^r X(1^+)^{K_i}$ when $S\cong\coprod_{i=1}^n H/K_i$.
\end{proof}
\end{prop}

We call the weak equivalences of the above model structure the \emph{$G$-equivariant special weak equivalences}.\index{G-equivariant special weak equivalence@$G$-equivariant special weak equivalence|textbf}\index{G-equivariant weak equivalence@$G$-equivariant weak equivalence!special|seeonly{$G$-equivariant special weak equivalence}} They can again be detected by mapping into special $\Gamma$-$G$-spaces in the homotopy category (with respect to the $G$-equivariant level weak equivalences); in particular, each of the maps $(\ref{eq:generating-cof-bfl-Gamma-G})$ is a $G$-equivariant special weak equivalence.

\begin{prop}\index{G-global special model structure@$G$-global special model structure!on Gamma-G-II-SSet@on $\cat{$\bm\Gamma$-$\bm G$-$\bm{\mathcal I}$-SSet}$}
There is a unique model structure on $\cat{$\bm\Gamma$-$\bm G$-$\bm{\mathcal I}$-SSet}_*$ with the same cofibrations as the $G$-global level model structure and with fibrant objects the $G$-globally level fibrant \emph{special} $G$-global $\Gamma$-spaces. This model structure is combinatorial, simplicial, left proper, and filtered colimits in it are homotopical.
\begin{proof}
This follows similarly by localizing with respect to the maps
\begin{equation*}
\lambda_{H,S,\phi,\psi}\colon \big(S_+\smashp\Gamma(1^+,\blank)\smashp (\mathcal I(H,\blank)\times G)_+\big)/H\to\big(\Gamma(S_+,\blank)\smashp (\mathcal I(H,\blank)\times G)_+\big)/H
\end{equation*}
analogous to the above.
\end{proof}
\end{prop}

\begin{lemma}
Let $K$ be any group and let $X\in\cat{$\bm{(G\times K})$-SSet}_*$ be cofibrant in the $\mathcal G_{\mathcal Fin,K}$-equivariant model structure. Then the map $(S_+\smashp\Gamma(1_+,\blank)\smashp X)/K\to (\Gamma(S_+,\blank)\smashp X)/K$ induced by $(\ref{eq:segal-corepr})$ is a $G$-equivariant special weak equivalence between cofibrant objects for every finite $K$-set $S$.
\begin{proof}
We first note that the target is indeed cofibrant by Lemma~\ref{lemma:equivariant-Gamma-generalized-cof}, and so is the source by the following observation:

\begin{claim*}
The functor $(S_+\smashp\Gamma(1_+,\blank)\smashp \blank)/K\colon\cat{$\bm{(G\times K)}$-SSet}_*\to\cat{$\bm\Gamma$-$\bm G$-SSet}_*$ sends $\mathcal G_{\mathcal Fin,K}$-cofibrations to cofibrations in the $G$-global level model structure.
\begin{proof}
We will show that $(S_+\smashp\Gamma(1_+,\blank)\smashp (G\times K)_+/\Gamma_{H,\phi})/K$ is cofibrant for any finite $H\subset G$ and any homomorphism $\phi\colon H\to K$; the claim will then follow as in Lemma~\ref{lemma:equivariant-Gamma-generalized-cof}. But indeed, this is isomorphic to $((\phi^*S)_+\smashp\Gamma(1_+,\blank)\smashp G)/H$, which is cofibrant by the proof of Proposition~\ref{prop:equivariant-special-model-structure} above.
\end{proof}
\end{claim*}

We will verify the assumptions of Corollary~\ref{cor:saturated-trafo} for the natural transformation $\tau\colon(S_+\smashp\Gamma(1_+,\blank)\smashp \blank)/K\Rightarrow (\Gamma(S_+,\blank)\smashp \blank)/K$ induced by $(\ref{eq:segal-corepr})$.

If $X=(G\times K)/\Gamma_{H,\phi}$ for some finite subgroup $H\subset G$ and some homomorphism $\phi\colon H\to K$, then $\tau_X$ is clearly conjugate to the map $(\phi^*S_+\smashp\Gamma(1_+,\blank)\smashp \blank)/H\to (\Gamma(\phi^*S_+,\blank)\smashp \blank)/H$ induced by $(\ref{eq:segal-corepr})$, hence a $G$-equivariant special weak equivalence by the above. Moreover, if $Y$ is any simplicial set, then $\tau_{X\times Y}$ is conjugate to $\tau_X\times Y$, so $\tau$ is a $G$-equivariant special weak equivalence for $(G\times K)/\Gamma_{H,\phi}\times Y$ as the $G$-global special model structure is simplicial. In particular, $\tau$ is a $G$-global special weak equivalences on the sources and targets of the generating cofibrations. The claim now follows as $(S_+\smashp\Gamma(1_+,\blank)\smashp \blank)/K$ and $(\Gamma(S_+,\blank)\smashp \blank)/K$ each preserve cofibrations as well as small colimits.
\end{proof}
\end{lemma}

\begin{cor}
The simplicial adjunction
\begin{equation}\label{eq:special-ev-R-model-cat}
\ev_{\mathcal U_G}\colon(\cat{$\bm\Gamma$-$\bm G$-$\bm{\mathcal I}$-SSet}_*)_{\textup{$G$-global special}}\rightleftarrows(\cat{$\bm\Gamma$-$\bm G$-SSet}_*)_{\textup{$G$-equiv.~special}} :\!R
\end{equation}
is a Quillen adjunction with fully homotopical left adjoint.
\begin{proof}
Let us first show that this is a Quillen adjunction. Since we already know this for the respective level model structures (Proposition~\ref{prop:G-global-Gamma-vs-Gamma-G}), it suffices to show that $R$ sends fibrant objects to special $G$-global $\Gamma$-spaces.

For this we let $H$ be any finite group and $\phi\colon H\to G$ any homomorphism, yielding a right $H$-action on $G$. Then we have seen in the proof of Proposition~\ref{prop:G-global-Gamma-vs-Gamma-G} that $(E\Inj(H,\mathcal U_G)\times G)_+$ is cofibrant in the $\mathcal G_{\mathcal Fin,H}$-model structure, so the map $(S_+\smashp\Gamma(1^+,\blank)\smashp(E\Inj(H,\mathcal U_G)\times G)_+)/H\to (\Gamma(S_+,\blank)\smashp (E\Inj(H,\mathcal U_G)\times G)_+)/H$
induced by $(\ref{eq:segal-corepr})$ is a $G$-equivariant special weak equivalence of cofibrant objects by the previous lemma. The claim now follows from an easy adjointness argument.

We conclude in particular that $\ev_{\mathcal U_G}$ sends $G$-global special acyclic cofibrations to $G$-equivariant special weak equivalences. As any $G$-global special weak equivalence factors as an acyclic cofibration followed by a $G$-global \emph{level} weak equivalence, we conclude together with Lemma~\ref{lemma:gamma-I-vs-equivariant} that $\ev_{\mathcal U_G}$ is homotopical.
\end{proof}
\end{cor}

\begin{proof}[Proof of Theorem~\ref{thm:special-G-global-Gamma-vs-Gamma-G}]
It is easy to check that $\ev_{\mathcal U_G}$ sends special $G$-global $\Gamma$-spaces to special $\Gamma$-$G$-spaces. On the other hand, the previous corollary implies that the right derived functor
\begin{equation*}
\cat{R}R\colon(\cat{$\bm\Gamma$-$\bm G$-SSet}_*)_{\textup{$G$-equivariant level}}^\infty\to(\cat{$\bm\Gamma$-$\bm G$-$\bm{\mathcal I}$-SSet}_*)_{\textup{$G$-global level}}^\infty
\end{equation*}
sends special $\Gamma$-$G$-special spaces to special $G$-global $\Gamma$-spaces, so that the Bousfield localization $\ev_{\mathcal U_G}^\infty\dashv\cat{R}R$ from Proposition~\ref{prop:G-global-Gamma-vs-Gamma-G} restricts to a Bousfield localization between the full subcategory of $\cat{$\bm\Gamma$-$\bm G$-$\bm{\mathcal I}$-SSet}_*^\infty$ spanned by the special $G$-global $\Gamma$-spaces and the full subcategory of $\cat{$\bm\Gamma$-$\bm G$-SSet}_*^\infty$ spanned by the special $\Gamma$-$G$-spaces. The claim now follows from Proposition~\ref{prop:localization-subcategory}.
\end{proof}

\subsubsection{Semiadditivity} We want to think of special $G$-global $\Gamma$-spaces as $G$-globally homoto\-py coherent versions of commutative monoids, so at the very least their homotopy category should be semiadditive. In this subsection we will prove the following `underived' version of this:

\begin{thm}\label{thm:gamma-special-add}\index{G-global Gamma-space@$G$-global $\Gamma$-space!semiadditivity|textbf}
The canonical map $\iota\colon X\vee Y\to X\times Y$ is a $G$-global special weak equivalence for all $X,Y\in\cat{$\bm\Gamma$-$\bm{E\mathcal M}$-$\bm G$-SSet}^\tau_*$.
\end{thm}

The proof we will give below is a spiced-up version of the usual argument that finite coproducts and products in a category of commutative monoids agree. This requires some preparations.

\begin{constr}\label{constr:rho-mu}
Let $A,B,C\in\Gamma$, and let $f\colon B\to C$ be any morphism. Then we have an induced map $A\smashp f\colon A\smashp B\to A\smashp C$. If $C=1^+$, then we will by slight abuse of notation also denote the composition $A\smashp B\to A\smashp C\cong A$ with the canonical isomorphism $A\smashp C\cong A,[a,1]\mapsto a$ by $A\smashp f$, and similarly for $B=1^+$.

Now let $T\in\cat{$\bm\Gamma$-$\bm{E\mathcal M}$-$\bm G$-SSet}_*^\tau$. Applying the above to the map $\mu\colon 2^+\to 1^+$, $\mu(1)=\mu(2)=1$ induces $T(\mu\smashp\blank)\colon T(2^+\smashp\blank)\to T$, and this is clearly natural in $T$. On the other hand, the assignment $T\mapsto T(2^+\smashp\blank)$ obviously preserves $G$-global level weak equivalences, so it descends to a functor on the homotopy category $\Ho(\cat{$\bm\Gamma$-$\bm{E\mathcal M}$-$\bm G$-SSet}_*^\tau)_{\text{strict}}$ with respect to these. It follows formally that the map $\mu\colon T(2^+\smashp\blank)\to T$ is also natural with respect to maps in the homotopy category.

Similarly the Segal maps assemble into a natural morphism
\begin{equation*}
\rho\mathrel{:=}(T(p_1\smashp\blank),T(p_2\smashp\blank))\colon T(2^+\smashp\blank)\to T\times T;
\end{equation*}
\index{Segal map}%
as the product preserves $G$-global level weak equivalences, this again descends to a natural transformation on the strict homotopy category.
\end{constr}

The following lemma follows easily from the definitions and we omit its proof.

\begin{lemma}\label{lemma:generalized-segal}
Let $T\in\cat{$\bm\Gamma$-$\bm{E\mathcal M}$-$\bm G$-SSet}_*^\tau$ be special. Then $\rho\colon T(2^+\smashp\blank)\to T\times T$ is a $G$-global level weak equivalence.\qed
\end{lemma}

If $T$ is special, then we write $m_T$ for the map $T\times T\to T$ in the strict homotopy category corresponding to the zig-zag
\begin{equation*}
\begin{tikzcd}
T\times T & \arrow[l, "\sim", "\rho"'] T(2^+\smashp\blank) \arrow[r, "T(\mu\smashp\blank)"] &[1em] T.
\end{tikzcd}
\end{equation*}

\begin{proof}[Proof of Theorem~\ref{thm:gamma-special-add}]
Let us fix a special $T\in\cat{$\bm\Gamma$-$\bm{E\mathcal M}$-$\bm G$-SSet}^\tau_*$; we have to show that $\iota^*\colon[X\times Y,T]\to[X\vee Y,T]$ is bijective, where $[\,{,}\,]$ denotes hom sets in the strict homotopy category.

As $\vee$ is homotopical in $G$-global level weak equivalences, $X\vee Y$ is also a coproduct in the strict homotopy category, so that $(i_X^*,i_Y^*)\colon[X\vee Y,T]\to[X,T]\times[Y,T]$ is bijective. It therefore suffices to show that the composition $\alpha\colon[X\times Y,T]\to[X,T]\times[Y,T]$ is bijective; plugging in the definition, we see that this is given by $(j_X^*,j_Y^*)$, with $j_X\colon X\to X\times Y$, $j_Y\colon Y\to X\times Y$ the usual inclusions.

To prove the claim, we define $\beta\colon [X,T]\times[Y,T]\to[X\times Y,T]$ as follows: if $f\in [X,T],g\in[Y,T]$ are arbitrary, then $\beta(f,g)\mathrel{:=}m_T\circ(f\times g)$; note that this indeed makes sense because $\times$ descends to the homotopy category as above.

We claim that $\beta$ is a two-sided inverse of $\alpha$. Indeed, if $f\in[X,T],g\in[Y,T]$ are arbitrary, then $\alpha\beta(f,g)=\alpha(m_T\circ(f\times g))=(m_T\circ (f\times g)\circ j_X, m_T\circ (f\times g)\circ j_Y)$. We will prove that $m_T\circ (f\times g)\circ j_X=f$, the argument for the second component being similar. Indeed, we have $(f\times g)j_X=j_1 f$ as both sides agree after postcomposing with each of the two projections $T\times T\to T$ (here we use that $T\times T$ is a product in the homotopy category). On the other hand, let us consider the diagram
\begin{equation}\label{diag:gamma-add-unitality}
\begin{tikzcd}
&& T\\
T\arrow[urr, bend left=15pt, "="]\arrow[rr, "T(i_1\smashp\blank)"]\arrow[drr, bend right=15pt, "j_1"'] && T(2^+\smashp\blank)\arrow[u, "T(\mu\smashp\blank)"']\arrow[d,"\rho"]\\
&& T\times T
\end{tikzcd}
\end{equation}
(say, on the point-set level), where $i_1\colon 1^+\to 2^+$ is defined by $i_1(1)=1$. Then the top triangle commutes because $\mu i_1=\id_{1^+}$, and we claim that also the lower triangle commutes. Indeed, after postcomposing with the projection to the first factor both paths through the diagram are the identity (as $p_1i_1=\id_{1^+}$). On the other hand, $\pr_2j_1$ is constant at the basepoint, as is $\pr_2\circ\rho\circ T(i_1\smashp\blank)= T(p_2\smashp\blank)\circ T(i_1\smashp\blank)$ since it factors through $T(0^+)=*$. From commutativity of $(\ref{diag:gamma-add-unitality})$ we can now conclude $m_Tj_1=T(\mu\smashp\blank)T(i_1\smashp\blank)=\id$, hence $m_T\circ(f\times g)\circ j_X=f$ as desired.

Finally, let $F\in [X\times Y,T]$ be arbitrary and consider the diagram
\begin{gather}
\begin{aligned}\label{diag:additivity-gamma}
\begin{tikzcd}[column sep=large, ampersand replacement=\&]
\&\& (X\times Y)\times (X\times Y)\arrow[r, "F\times F"] \& T\times T\\
X\times Y\arrow[drr, bend right=12pt, "="']\arrow[urr, bend left=12pt, "j_X\times j_Y"]\arrow[rr, "X(i_1\smashp\blank)\times Y(i_2\smashp\blank)"] \&\& X(2^+\smashp\blank)\times Y(2^+\smashp\blank)\arrow[u,"\rho"']\arrow[d, "X(\mu\smashp\blank)\times Y(\mu\smashp\blank)"]\arrow[r, "F(2^+\smashp\blank)"] \& T(2^+\smashp\blank)\arrow[u, "\rho"']\arrow[d, "T(\mu\smashp\blank)"]\\
\&\& X\times Y\arrow[r, "F"'] \& T
\end{tikzcd}
\end{aligned}\raisetag{-5pt}
\end{gather}
in the strict homotopy category. The right hand portion commutes by the naturality established in Construction~\ref{constr:rho-mu}, and as above one shows that the triangles on the left already commute on the point-set level.

Using the commutativity of $(\ref{diag:additivity-gamma})$ we then compute
\begin{align*}
\beta\alpha(F)&=\beta(Fj_X,Fj_Y)\\
&=m_T\circ\big((Fj_X)\times (Fj_Y)\big)=m_T\circ (F\times F)\circ(j_X\times j_Y)\\
&=m_T\circ\rho\circ F(2^+\smashp\blank)\circ \big(X(i_1\smashp\blank)\times Y(i_2\smashp\blank)\big)\\
&=T(\mu\smashp\blank)\circ F(2^+\smashp\blank)\circ \big(X(i_1\smashp\blank)\times Y(i_2\smashp\blank)\big)=F,
\end{align*}
which completes the proof of the theorem.
\end{proof}

\begin{cor}\label{cor:product-special-homotopical}
Finite coproducts and finite products in $\cat{$\bm\Gamma$-$\bm{E\mathcal M}$-$\bm G$-SSet}_*^\tau$ preserve $G$-global special weak equivalences.
\begin{proof}
It suffices to treat the case of binary coproducts and products, for which we let $f\colon X\to X''$ and $g\colon Y\to Y''$ be any $G$-global special weak equivalences. Then we can factor $f$ as an acyclic cofibration $i\colon X\to X'$ followed by an (automatically acyclic) fibration $p\colon X'\to X''$, and similarly $g=qj$ with an acyclic fibration $q$ and an acyclic cofibration $j$. In the commutative diagram
\begin{equation*}
\begin{tikzcd}[ampersand replacement=\&]
X\vee Y\arrow[r, "i\vee j"]\arrow[d, "\iota"'] \& X'\vee Y'\arrow[d, "\iota"]\arrow[r, "p\vee q"] \& X''\vee Y''\arrow[d,"\iota"]\\
X\times Y\arrow[r, "i\times j"'] \& X'\times Y'\arrow[r, "p\times q"'] \& X''\times Y''
\end{tikzcd}
\end{equation*}
the vertical arrows are $G$-global special weak equivalences by the previous theorem, and so are the top left and bottom right arrows as acyclic cofibrations in any model category are stable under (all) coproducts while acyclic fibrations are stable under products. The claim follows by $2$-out-of-$3$.
\end{proof}
\end{cor}

\begin{warn}
We can define specialness for elements of the ordinary functor category $\cat{$\bm\Gamma$-$\bm{E\mathcal M}$-$\bm G$-SSet}^\tau\mathrel{:=}\Fun(\Gamma,\cat{$\bm{E\mathcal M}$-$\bm G$-SSet}^\tau)$ in the same way as above, which leads to a notion of $G$-global special weak equivalences on $\cat{$\bm\Gamma$-$\bm{E\mathcal M}$-$\bm G$-SSet}^\tau$. These are however \emph{not} stable under finite products: for example the map
\begin{equation*}
f\colon (\Gamma(1^+,\blank)\amalg\Gamma(1^+,\blank))\times E\Inj(*,\omega)\times G\to\Gamma(2^+,\blank)\times E\Inj(*,\omega)\times G
\end{equation*}
induced by restricting along $p_1,p_2$ is a $G$-global special weak equivalence because for any level fibrant $T$ the induced map $[f,T]$ is conjugate to the map $\pi_0 T(2^+)(\omega)\to\pi_0 T(1^+)(\omega)\times\pi_0 T(1^+)(\omega)$ induced by the Segal map. However, $\Gamma(1^+,\blank)\times f$ is conjugate to a map $(\Gamma(2^+,\blank)\amalg\Gamma(2^+,\blank))\times E\Inj(*,\omega)\times G\to\Gamma(3^+,\blank)\times E\Inj(*,\omega)\times G$ by the universal property of coproducts in $\Gamma$, hence not a $G$-global special weak equivalence by a similar calculation as before.

Thus, the fact that also non-special $G$-global $\Gamma$-spaces are trivial in degree $0^+$ (or at least weakly contractible) is crucial for the theorem, which is why we have been a bit more verbose in verifying that certain diagrams commute than usual.
\end{warn}

\subsubsection{The Wirthmüller isomorphism}\label{subsubsec:wirthmueller}
\index{Wirthmüller isomorphism!in Gamma-EM-G-SSettau@in $\cat{$\bm\Gamma$-$\bm{E\mathcal M}$-$\bm G$-SSet}_*^\tau$|(}
\index{G-global Gamma-space@$G$-global $\Gamma$-space!Wirthmüller isomorphism|seeonly{Wirthmüller isomorphism, for $G$-global $\Gamma$-spaces}}
Let $\alpha\colon H\to G$ be an injective homomorphism of finite groups. If $X$ is a genuine $H$-equivariant spectrum, the \emph{Wirthmüller isomorphism} is a specific $G$-weak equivalence $\gamma\colon\alpha_!X\to\alpha_*X$,
see e.g.~\cite[Proposition~3.7]{hausmann-equivariant} and also \cite[Theorem~2.1.10]{proper-equivariant} or \cite[Theorem~II.6.2]{lms} for generalizations to the proper or compact Lie case, respectively.\index{Wirthmüller isomorphism!in G-Spectra@in $\cat{$\bm G$-Spectra}$!G-equivariantly@$G$-equivariantly}

In a precise sense, the Wirthmüller isomorphism marks the distinction between \emph{genuine} stable equivariant homotopy theory (encoding deloopings against all representation spheres) and na\"ive stable equivariant homotopy theory (only admitting deloopings against spheres with trivial actions). As non-equivariantly $\alpha_!$ is given by a $(G:\im\alpha)$-fold wedge, while $\alpha_*$ is a $(G:\im\alpha)$-fold product, we can also view this as some sort of `twisted semiadditivity.' Below, we will construct an analogue of the Wirthmüller map $\gamma$ in our $\Gamma$-space context and prove:

\begin{thm}\label{thm:gamma-special-wirthmueller}\index{Wirthmüller isomorphism!in Gamma-EM-G-SSettau@in $\cat{$\bm\Gamma$-$\bm{E\mathcal M}$-$\bm G$-SSet}_*^\tau$|textbf}
Let $\alpha\colon H\to G$ be an injective homomorphism (of not necessarily finite groups) such that $(G:\im\alpha)<\infty$, and let  $X\in\cat{$\bm\Gamma$-$\bm{E\mathcal M}$-$\bm H$-SSet}^\tau_*$. Then the Wirthmüller map $\gamma\colon \alpha_!X\to\alpha_*X$ is a $G$-global special weak equivalence.
\end{thm}

The above theorem (in the guise of Corollary~\ref{cor:boxtimes-special-power}) will also be instrumental in the proof of the equivalence between $G$-ultra-commutative monoids and special $G$-global $\Gamma$-spaces that we will give in the next section.

\begin{constr}\label{constr:wirthmueller-set}\index{Wirthmüller map|textbf}\index{Wirthmüller map|seealso{Wirthmüller isomorphism}}
Let us first recall the Wirthmüller map, which actually already exists in the based context:

Without loss of generality, we may assume for the rest of this discussion that $H$ is a finite index subgroup of $G$ and that $\alpha$ is its inclusion. Then $\alpha_!$ can be modelled by applying $G_+\smashp_H\blank\colon\cat{$\bm H$-Set}_*\to\cat{$\bm G$-Set}_*$ levelwise and pulling through the $E\mathcal M$-action; the counit of this $\cat{Set}_*$-level adjunction is then given by $i\colon Y\to G_+\smashp_HY,y\mapsto[1,y]$. Similarly, $\alpha_*$ is given by applying $\Maps^H(G,\blank)$ levelwise.

If now $Y$ is any pointed $H$-set, then we define the \emph{Wirthmüller map} as the $G$-equivariant map $\gamma\colon G_+\smashp_HY \to\Maps^H(G,Y)$\nomenclature[agamma]{$\gamma$}{Wirthmüller map} adjoint to the $H$-equivariant map
\begin{equation*}
G_+\smashp_H Y\to Y, [g,y]\mapsto\begin{cases}
g.y & \text{if }g\in H\\
* & \text{otherwise.}
\end{cases}
\end{equation*}
Explicitly,
\begin{equation*}
\gamma[g,y](g')=\begin{cases}
g'g.y & \text{if }g'g\in H\\
* & \text{otherwise}.
\end{cases}
\end{equation*}
It is easy to check that the Wirthmüller map is natural; in particular, we can apply it levelwise to get a natural map $\gamma\colon\alpha_!X\to\alpha_*X$ for any $X\in\cat{$\bm\Gamma$-$\bm{E\mathcal M}$-$\bm{H}$-SSet}_*^\tau$.
\end{constr}

Just like we can think of Theorem~\ref{thm:gamma-special-wirthmueller} as twisted semiadditivity, the basic idea of our proof will be similar to the proof of semiadditivity (Theorem~\ref{thm:gamma-special-add}); however, the actual combinatorics are a bit more complicated and one has to be slightly careful in keeping track of all the actions involved.

\begin{constr}
If $Z$ is a pointed $G$-set and $Y=\res_H^G Z$, then we can partially `untwist' the action on $\alpha_!Y=G_+\smashp_H\res^G_HZ$ by the usual $G$-equivariant \emph{shearing isomorphism}\index{shearing isomorphism|textbf}
\begin{equation*}
\shear\colon (G/H)_+\smashp Z\to G_+\smashp_H\res^G_HZ, [[g],z]\mapsto[g,g^{-1}.z],
\end{equation*}
and dually we have a natural $G$-equivariant \emph{coshearing isomorphism}\index{coshearing isomorphism|textbf}
\begin{equation*}
\coshear\colon Z^{G/H}\to\Maps^H(G,\res^G_HZ)
\end{equation*}
defined by $\coshear(y_\bullet)(g)=g.z_{[g^{-1}]}$. Again, applying this levelwise we can extend these to natural maps of $G$-global $\Gamma$-spaces.
\end{constr}

\begin{constr}
Let $X$ be a $G$-global $\Gamma$-space. We define the \emph{twisted Segal map}\index{Segal map!twisted|textbf} $\varrho\colon X(G/H_+\smashp\blank)\to\Maps^H(G,\res_H^G X)$\nomenclature[arho3]{$\varrho$}{twisted Segal map} as the composition
\begin{equation*}
X(G/H_+\smashp\blank)\xrightarrow{\;\rho\;} X^{G/H}\xrightarrow{\coshear} \Maps^H(G,\res^G_HX).
\end{equation*}
Here the map $\rho$ is obtained by applying in each degree the `generalized Segal map' $X(G/H_+\smashp S_+)\to X(S_+)^{\times(G/H)}$ induced on factor $[g]$ by $p_{[g]}\smashp S_+\colon G/H_+\smashp S_+\to S_+$.
\end{constr}

As in the previous section, $\varrho$ is natural on the point-set level. Moreover, $\Maps^H(G,\res^G_H\blank)$ is homotopical in the $G$-global level weak equivalences by Corollary~\ref{cor:gamma-tau-alpha-lower-star}, and so is $X(G/H_+\smashp\blank)$ for trivial reasons. As before we conclude that $\varrho$ is also natural in maps in the strict homotopy category.

\begin{lemma}
Let $T$ be a special $G$-global $\Gamma$-space. Then $\varrho\colon T(G/H_+\smashp\blank)\to \Maps^H(G,\res^G_HT)$ is a $G$-global level weak equivalence.
\begin{proof}
It suffices to show that the generalized Segal map $\rho\colon T(G/H_+\smashp\blank)\to T^{G/H}$ is a $G$-global level weak equivalence (where $G$ acts via its action on $T$ and on $G/H$). For this we observe that for any finite set $S$ the composition
\begin{align*}
T((G/H\times S)_+)\cong T(G/H_+\smashp S_+)&\xrightarrow{\rho} T(S_+)^{\times(G/H)}\\
&\xrightarrow{\rho_S^{\times(G/H)}} (T(1^+)^{\times S})^{\times(G/H)}\cong T(1^+)^{\times(G/H\times S)}
\end{align*}
agrees with the Segal map $\rho_{G/H\times S}$ for $G/H\times S$, so it is a $(G\times \Sigma_{G/H\times S})$-global weak equivalence. Applying Corollary~\ref{cor:gamma-tau-alpha-lower-shriek} to the homomorphism $G\times\Sigma_S\to G\times\Sigma_{G/H\times S}$ induced by the identity of $G$ and the homomorphism $G\times\Sigma_S\to \Sigma_{G/H\times S}$ classifying the obvious $(G\times\Sigma_S)$-action on $G/H\times S$ then shows that this is in particular a $(G\times\Sigma_S)$-global weak equivalence.

Similarly, $\rho_S$ is a $(G\times\Sigma_S)$-global weak equivalence, so $\rho_S^{\times(G/H)}$ is a $(G\times\Sigma_S\times\Sigma_{G/H})$-global weak equivalence by Corollary~\ref{cor:twisted-product}, hence in particular a $(G\times\Sigma_S)$-global weak equivalence. The claim now follows by $2$-out-of-$3$.
\end{proof}
\end{lemma}

We can now use this to define a twisted version of the `multiplication map' $m_T$ considered before:

\begin{constr}
We write $\nu\colon G/H_+\to1^+$ for the map with $\nu[g]=1$ for all $g\in G$, which induces for any $X\in\cat{$\bm\Gamma$-$\bm{E\mathcal M}$-$\bm G$-SSet}_*^\tau$ a natural map $X(\nu\smashp\blank)\colon X(G/H_+\smashp\blank)\to X$. As before this descends to a natural map on the strict homotopy category.

If $X=T$ is special, then we define $n_T\colon \Maps^H(G,\res^G_HT)\to T$ as the map in the strict homotopy category represented by the zig-zag
\begin{equation*}
\begin{tikzcd}
{\Maps^H(G,\res^G_HT)} & \arrow[l, "\sim", "\varrho"'] T(G/H_+\smashp \blank) \arrow[r, "{T(\nu\smashp\blank)}"] &[1em] T.
\end{tikzcd}
\end{equation*}
\end{constr}

\begin{constr}
We now define a `diagonal map' $\delta\colon \Maps^H(G,X)\to \Maps^H(G,X)(G/H_+\smashp\blank)$ for every $H$-global $\Gamma$-space $X$ as follows: if $S$ is a finite set, $n\ge 0$, and $f\colon G\to X(S_+)_n$ is $H$-equivariant, then $\delta(f)\colon G\to X(G/H_+\smashp S_+)_n$ is defined via $\delta(f)(g)=X(i_{[g^{-1}]}\smashp S_+)(f(g))$, where $i_{[g]}\colon 1^+\to S_+$ sends $1$ to $[g]$.
\end{constr}

\begin{lemma}
The above defines a map in $\cat{$\bm\Gamma$-$\bm{E\mathcal M}$-$\bm G$-SSet}_*^\tau$.
\begin{proof}
We first show that $\delta(f)$ is $H$-equivariant, for which it is important to observe that the $H$-action from the definition of $\Maps^H(G,X)(G/H_+\smashp S_+)$ is via the $H$-action on $X$ only; \emph{the $H$-action on $G/H$ does not come into play yet.} Thus,
\begin{align*}
h.\big(\delta(f)(g)\big)&=h.\big(X(i_{[g^{-1}]})(f(g))\big)= X(i_{[g^{-1}]})\big(h.(f(g))\big)\\
&=X(i_{[g^{-1}]})\big(f(hg)\big)=X(i_{[(hg)^{-1}]})\big(f(hg)\big)=\delta(f)(hg)
\end{align*}
where the second equality uses that $X(i_{[g^{-1}]})$ is equivariant with respect to the $H$-action \emph{coming from $X$ only}, whereas the penultimate equation uses that $[(hg)^{-1}]=[g^{-1}h^{-1}]=[g^{-1}]$ in $G/H$.

Thus, $\delta$ indeed lands in $\Maps^H(G,X)(G/H_+\smashp\blank)$. It is then easy to check that $\delta$ is compatible with the simplicial structure maps, the structure maps of $\Gamma$, and the $E\mathcal M$-action, so it only remains to prove $G$-equivariance. For this we observe that the $G$-action on $\Maps^H(G,X)(G/H_+\smashp S_+)$ is via the diagonal of the right $G$-action on $G$ and the left $G$-action on $G/H$. Thus, the chain of equalities
\begin{align*}
\delta(g'.f)(g)&= X(i_{[g^{-1}]})\big((g'.f)(g)\big)=X(i_{[g^{-1}]})\big(f(gg'))\\
&=X\big((g'.\blank)\smashp\blank\big)X(i_{[(gg')^{-1}]})\big(f(gg')\big)=X\big((g'.\blank)\smashp\blank\big)\big(\delta(f)(gg')\big)
\end{align*}
precisely shows that $\delta$ is $G$-equivariant.
\end{proof}
\end{lemma}

\begin{constr}
We define for any pointed $H$-set $Y$ an $H$-equivariant map $j\colon Y\to\res^G_H\Maps^H(G,Y)$ via
\begin{equation*}
j(y)(g)=\begin{cases}
g.y & \text{if }g\in H\\
* & \text{otherwise}.
\end{cases}
\end{equation*}
This is clearly natural, and in particular we can apply it levelwise to get a natural map $X\to\res^G_H\Maps^H(G,X)$ for any $H$-global $\Gamma$-space $X$.
\end{constr}

\begin{prop}
Let $F\colon\Maps^H(G,X)\to T$ be any morphism in the strict homotopy category of $\cat{$\bm\Gamma$-$\bm{E\mathcal M}$-$\bm G$-SSet}_*^\tau$. Then the diagram
\begingroup\hfuzz=20pt
\begin{equation*}\hskip-3pt
\begin{tikzcd}[column sep=large]
&[-2.75em] \Maps^H(G,\res^G_H\Maps^H(G,X))\arrow[rr, "{\Maps^H(G,\res^G_HF)}\;"] &[-.67em]& \Maps^H(G,\res^G_HT)\\
\Maps^H(G,X)\arrow[ru, bend left=12pt, "{\Maps^H(G,j)}"]\arrow[r, "\delta"]\arrow[rd, "="', bend right=12pt] &\Maps^H(G,X)(G/H_+\smashp\blank)\arrow[rr, "F(G/H_+\smashp\blank)"]\arrow[u, "\varrho"']\arrow[d, "{\Maps^H(G,X)(\nu\smashp\blank)}"] && T(G/H_+\smashp\blank)\arrow[u, "\varrho"']\arrow[d, "T(\nu\smashp\blank)"]\\
& \Maps^H(G,X)\arrow[rr, "F"'] && T
\end{tikzcd}
\end{equation*}
\endgroup
in the strict homotopy category commutes.
\begin{proof}
The right hand portion commutes by the above naturality considerations, and we will now prove that the two triangles on the left already commute on the point-set level.

Let us consider the lower triangle first. If $S$ is a finite set, $n\ge 0$, and $f\colon G\to X(S_+)_n$ is $H$-equivariant, then
\begin{align*}
\Maps^H(G,X)(\nu\smashp S_+)(\delta(f))(g)&=X(\nu\smashp S_+)(\delta(f)(g))\\
&=X(\nu\smashp S_+)X(i_{[g^{-1}]}\smashp S_+)\big(f(g))=f(g)
\end{align*}
for all $g\in G$, i.e.~$\Maps^H(G,X)(\nu\smashp S_+)\circ\delta=\id$.

Similarly, we compute for the upper left triangle
\begin{align*}
\varrho(\delta(f))(g)&=\coshear\big(\rho(\delta(f))\big)(g)=
g.\big(\rho(\delta(f))_{[g^{-1}]})\\
&=g.\big(\Maps^H(G,X)(p_{[g^{-1}]}))(\delta(f))\big)\in\res^G_H\Maps^H(G,X(S_+)_n),
\end{align*}
hence
\begin{align*}
\varrho(\delta(f))(g_1)(g_2)&=\big(g_1.\big(\Maps^H(G,X)(p_{[g^{-1}_1]}))(\delta(f))\big)\big)(g_2)\\
&=\big(\Maps^H(G,X)(p_{[g_1^{-1}]})(\delta(f))\big)(g_2g_1)\\
&=X(p_{[g_1^{-1}]})\big(\delta(f)(g_2g_1)\big)=X(p_{[g_1^{-1}]})X(i_{[g_1^{-1}g_2^{-1}]})\big(f(g_2g_1)\big).
\end{align*}
If $g_2\in H$, then $[g_1^{-1}g_2^{-1}]=[g_1^{-1}]$, hence $X(p_{[g_1^{-1}]})X(i_{[g_1^{-1}g_2^{-1}]})=\id$. On the other hand, if $g_2\notin H$, then $[g_1^{-1}g_2^{-1}]\not=[g_1^{-1}]$ and $X(p_{[g_1^{-1}]})X(i_{[g_1^{-1}g_2^{-1}]})$ factors through the base point. Thus,
\begin{equation}\label{eq:varrho-delta-f}
\varrho(\delta(f))(g_1)(g_2)=\begin{cases}
f(g_2g_1) & \text{if }g_2\in H\\
* & \text{otherwise}
\end{cases}=
\begin{cases}
g_2.\big(f(g_1)\big) & \text{if }g_2\in H\\
* & \text{otherwise}
\end{cases}
\end{equation}
for all $g_1,g_2\in G$, where the second equality uses $H$-equivariance of $f$.

On the other hand, $\big(\Maps^H(G,j)(f)\big)(g)=j(f(g))\in\res^G_H\Maps^H(G,X(S_+)_n)$ for all $g\in G$, hence
\begin{equation*}
\big(\Maps^H(G,j)(f)\big)(g_1)(g_2)=j(f(g_1))(g_2)=\begin{cases}
g_2.\big(f(g_1)\big) & \text{if }g_2\in H\\
* & \text{otherwise}
\end{cases}
\end{equation*}
for all $g_1,g_2\in G$, which agrees with $(\ref{eq:varrho-delta-f})$, finishing the proof.
\end{proof}
\end{prop}

\begin{proof}[Proof of Theorem~\ref{thm:gamma-special-wirthmueller}]
By Corollary~\ref{cor:gamma-tau-alpha-lower-shriek}, $G_+\smashp_H\blank$ and $\res^G_H$ are homotopical in the respective level weak equivalences, and they descend to an adjunction between strict homotopy categories. We conclude that we have for every $X\in\cat{$\bm\Gamma$-$\bm{E\mathcal M}$-$\bm H$-SSet}_*^\tau$ and $T\in\cat{$\bm\Gamma$-$\bm{E\mathcal M}$-$\bm G$-SSet}_*^\tau$ a bijection
\begin{equation*}
[G_+\smashp_HX,T]_G\to [X,T]_H, F\mapsto\res^G_H(F)\circ i;
\end{equation*}
here we write $[\,{,}\,]_H$ for the hom-sets in the strict homotopy category of $H$-global $\Gamma$-spaces, and $[\,{,}\,]_G$ for the corresponding hom-sets of $G$-global $\Gamma$-spaces. Using that $j=\res^G_H(\gamma)\circ i$, we are then reduced to showing that $\theta\colon [\Maps^H(G,X),T]_G\to [X,T]_H, F\mapsto\res^G_H(F)\circ j$ is bijective whenever $T$ is special.

For this, we define an explicit inverse $\zeta\colon[X,T]_H\to [\Maps^H(G,X),T]_G$ via $\zeta(f)=n_T\circ\Maps^H(G,f)$; here we have used again that $\Maps^H(G,\blank)$ descends to strict homotopy categories (Corollary~\ref{cor:gamma-tau-alpha-lower-star}).

If now $f\colon X\to\res^G_H T$ is any map in the strict homotopy category of $H$-global $\Gamma$-spaces, then $\theta\zeta(f)=\theta(n_T\circ\Maps^H(G,f))=\res^G_H(n_T)\circ\res^G_H(\Maps^H(G,f))\circ j= \res^G_H(n_T) \circ j\circ f$ by naturality of $j$. A straight-forward calculation then shows that the diagram
\begin{equation*}
\begin{tikzcd}[column sep=large]
& \res^G_H\Maps^H(G,\res^G_HT)\\
\res^G_HT\arrow[ur, "j", bend left=10pt]\arrow[r, "\res^G_Hk"]\arrow[dr, "\res^G_H(T(i_{[1]}\smashp\blank))"', bend right=10pt] & \res^G_H(T^{\times(G/H)})\arrow[u,"\res^G_H\coshear"']\\
& \res^G_H T(G/H\smashp\blank)\arrow[u, "\rho"']
\end{tikzcd}
\end{equation*}
in $\cat{$\bm\Gamma$-$\bm{E\mathcal M}$-$\bm H$-SSet}_*^\tau$ commutes, where $k\colon T\to T^{\times(G/H)}$ is the inclusion of the factor corresponding to $[1]\in G/H$. From this we immediately conclude that $\res^G_H(n_T)\circ j=\res^G_H\big(T(\nu\smashp\blank)\circ T(i_{[1]}\smashp\blank))=\id$, hence $\theta\zeta(f)=\res^G_H(n_T) \circ j\circ f=f$ as desired.

Finally, $\zeta\theta(F)=\zeta(\res^G_H(F)\circ j)=n_T\circ\Maps^H(G,\res^G_HF)\circ\Maps^H(G,j)$, which agrees with $F$ by the previous proposition.\index{Wirthmüller isomorphism!in Gamma-EM-G-SSettau@in $\cat{$\bm\Gamma$-$\bm{E\mathcal M}$-$\bm G$-SSet}_*^\tau$|)}
\end{proof}

\section{Comparison of the approaches}
A priori, the ultra-commutative or parsummable models are of a very different nature than the models based on $\Gamma$-spaces. However, in this section we will construct specific homotopical functors $\digamma\hskip0pt minus .5pt\colon\hskip0pt minus 1pt\cat{$\bm G$-ParSumSSet}\hskip0pt minus 1pt\to\hskip0pt minus 1pt\cat{$\bm\Gamma$\kern-.5pt-$\bm{E\mathcal M}$-$\bm G$-SSet}^{\tau,\textup{special}}_*$ and $\digamma\colon\cat{$\bm G$-UCom}\to\cat{$\bm\Gamma$-$\bm G$-$\bm{\mathcal I}$-SSet}^{\tau,\textup{special}}_*$, and prove:

\begin{thm}\label{thm:gamma-vs-uc}
The diagram
\begin{equation*}
\begin{tikzcd}
\cat{$\bm G$-UCom}\arrow[d, "\digamma"']\arrow[r, "\ev_\omega"] & \cat{$\bm G$-ParSumSSet}\arrow[d, "\digamma"]\\
\cat{$\bm\Gamma$-$\bm G$-$\bm{\mathcal I}$-SSet}_*^{\textup{special}}\arrow[r, "\ev_\omega"'] & \cat{$\bm\Gamma$-$\bm{E\mathcal M}$-$\bm G$-SSet}^{\tau,\textup{special}}_*
\end{tikzcd}
\end{equation*}
of homotopical functors commutes up to canonical isomorphism. Moreover, all these functors induce equivalences of associated quasi-categories.
\end{thm}

Together with Theorem~\ref{thm:uc-vs-uc} this will in particular show that through the eyes of finite groups, Schwede's ultra-commutative monoids (i.e.~commutative monoids for the box product on $\cat{$\bm L$-Top}$) are equivalent to a suitable notion of `special global $\Gamma$-spaces,' connecting them to classical approaches to equivariant coherent commutativity, see Corollary~\ref{cor:ucom-vs-global-gamma}.

On the other hand, together with Theorem~\ref{thm:CMon-G-I} we can also view the above result as a $G$-global strengthening of the equivalence between commutative monoids for the box product on $\cat{$\bm I$-SSet}$ and coherently commutative monoids in $\cat{SSet}$ due to Sagave and Schlichtkrull \cite[Theorem~1.2]{sagave-schlichtkrull}.

Our construction of the functors $\digamma$ (the archaic Greek letter digamma) is an analogue of \cite[Construction~4.3]{schwede-k-theory} for so-called \emph{parsummable categories}, which we will also recall later in Subsection~\ref{subsec:global-k}. While Schwede uses the letter $\gamma$, this is already taken in our context by the Wirthmüller isomorphism (which will actually play a crucial role in the proof of the theorem).

\begin{constr}\label{constr:digamma}
\nomenclature[adigamma]{$\digamma$ (digamma)}{$G$-global $\Gamma$-space associated to a $G$-parsummable simplicial set/$G$-ultra-commutative monoid}
\nomenclature[aDigamma]{$\Digamma$ (digamma)}{lift of $\digamma$ to $\Gamma$-$G$-parsummable simplicial sets}
\nomenclature[af]{$\digamma$}{\textit{see} $\digamma$ (digamma)\nomnorefpage}
\nomenclature[aF]{$\Digamma$}{\textit{see} $\Digamma$ (digamma)\nomnorefpage}
Let us write $\cat{$\bm\Gamma$-$\bm G$-ParSumSSet}_*$ for the category of $\cat{Set}_*$-enriched functors $\Gamma\to\cat{$\bm G$-ParSumSSet}$, which we can identify as before with ordinary functors $X$ such that $X(0^+)$ is terminal. Then the evaluation functor $\ev\colon\cat{$\bm\Gamma$-$\bm G$-ParSumSSet}_*\to \cat{$\bm G$-ParSumSSet},X\mapsto X(1^+)$ has a left adjoint $\Digamma$ (capital digamma) given by $\cat{Set}_*$-enriched left Kan extension.

Explicitly, $(\Digamma X)(S_+)= X^{\boxtimes S}$ (as $\boxtimes$ is the coproduct on $\cat{$\bm G$-ParSumSSet}$) with the evident functoriality in $X$. The functoriality in $S_+$ is as follows: if $f\colon S_+\to T_+$ is any map in $\Gamma$, then $(\Digamma X)(f)$ is the map $X^{\boxtimes S}\to X^{\boxtimes T}$ given in each simplicial degree by $(x_s)_{s\in S}\mapsto (y_t)_{t\in T}$ with $y_t=\sum_{s\in f^{-1}(t)}x_s$. By direct inspection, the Segal maps $X^{\boxtimes S}=(\Digamma X)(S_+)\to X(1^+)^{\times S}=X^{\times S}$ are precisely the inclusions, so the underlying $G$-global $\Gamma$-space of $\Digamma X$ is special by Theorem~\ref{thm:boxpower-symmetric-group}.

We now write $\digamma$ for the composition
\begin{equation}\label{eq:digamma-tau-def}
\cat{$\bm G$-ParSumSSet}\xrightarrow{\Digamma}
\cat{$\bm\Gamma$-$\bm G$-ParSumSSet}_*\xrightarrow{\forget}\cat{$\bm\Gamma$-$\bm{E\mathcal M}$-$\bm G$-SSet}_*^\tau.
\end{equation}
The construction of $\digamma\colon\cat{$\bm G$-UCom}\to\cat{$\bm\Gamma$-$\bm G$-$\bm{\mathcal I}$-SSet}^\tau_*$ is analogous.
\end{constr}

The main part of the proof of Theorem~\ref{thm:gamma-vs-uc} will be establishing that the composition $(\ref{eq:digamma-tau-def})$ induces an equivalence of homotopy theories. To this end we will introduce a suitable model structure on $\cat{$\bm\Gamma$-$\bm G$-ParSumSSet}_*$ and then show that both $\Digamma$ and the forgetful functor are already equivalences of homotopy theories.

\subsection{Symmetric products}
On $\cat{$\bm\Gamma$-$\bm{E\mathcal M}$-$\bm G$-SSet}_*^\tau$, we can define a \emph{box product}\index{box product!on Gamma-EM-SSet-tau@on $\cat{$\bm\Gamma$-$\bm{E\mathcal M}$-SSet}_*^\tau$} by performing the box product of tame $E\mathcal M$-simplicial sets levelwise; this is indeed well-defined as $*\boxtimes*=*$. We can then identify $\cat{$\bm\Gamma$-$\bm G$-ParSumSSet}_*$ with the category of commutative monoids for $\boxtimes$ on $\cat{$\bm\Gamma$-$\bm{E\mathcal M}$-$\bm G$-SSet}_*^\tau$, which suggests constructing the required model structure via the general machinery recalled in \ref{subsubsection:cmon-model-structures}. While we cannot directly apply this (as $\boxtimes$ no longer preserves initial objects in each variable \emph{separately}), our arguments will still be very close to the usual approach, and in particular we will need homotopical information about $G$-global $\Gamma$-spaces of the form $X^{\boxtimes n}/\Sigma_n$.

\begin{constr}\index{symmetric product|textbf}
Let $X$ be a $G$-global $\Gamma$-space and let $n\ge 1$. The \emph{$n$-th symmetric product} $\SP^nX$\nomenclature[aSPn]{$\SP^n$}{$n$-th symmetric product (in pointed contexts), agrees with $n$-th symmetric power for $n<\infty$} is defined as $X^{\boxtimes n}/\Sigma_n$; we will confuse $\SP^1X$ with $X$.

There is an evident way to make $\SP^n$ into an endofunctor of $\cat{$\bm\Gamma$-$\bm{E\mathcal M}$-$\bm G$-SSet}_*^\tau$. The map $X^{\boxtimes n}\to X^{\boxtimes(n+1)}$ given by inserting the basepoint in the last factor descends to a natural map $\SP^nX\to\SP^{n+1}X$. We define the functor $\SP^\infty$ as the colimit $\SP^1\Rightarrow\SP^2\Rightarrow\cdots$ along these natural maps.
\end{constr}

\begin{rk}
We have previously employed the notation $\Sym^nX$ for $X^{\otimes n}/\Sigma_n$ for any symmetric monoidal model category $\mathscr C$. Our reason for introducing new notation is twofold: firstly, writing $\Sym^\infty$ for $\SP^\infty$ would be ambiguous, as it is often used for $\coprod_{n\ge 0}\Sym^nX$; secondly, the change of notation forces us to remember that we cannot apply any of the previous results on $\Sym^n$ directly as $\boxtimes$ is no longer cocontinuous in each variable.
\end{rk}

The following theorem will be the key ingredient in establishing the model structure on $\cat{$\bm\Gamma$-$\bm G$-ParSumSSet}_*$ and comparing it to $\cat{$\bm\Gamma$-$\bm{E\mathcal M}$-$\bm G$-SSet}_*^\tau$:

\begin{thm}\label{thm:symmetric-products}
Let $X\in\cat{$\bm\Gamma$-$\bm{E\mathcal M}$-$\bm G$-SSet}_*^\tau$ and assume that $X(S_+)$ has no $\mathcal M$-fixed points apart from the base point for any finite set $S$. Then all the maps in
\begin{equation*}
X=\SP^1X\to\SP^2X\to\cdots\to\SP^\infty X
\end{equation*}
are $G$-global special weak equivalences.
\end{thm}

The proof of the theorem will be given below after some preparations. We begin with the following consequences of the results of the previous sections:

\begin{cor}\label{cor:boxtimes-special-power}
Let $n\ge 0$ and $X\in\cat{$\bm\Gamma$-$\bm{E\mathcal M}$-$\bm G$-SSet}^\tau_*$. Then the maps
\begin{equation*}
X^{\vee n}\to X^{\boxtimes n}\hookrightarrow X^{\times n}
\end{equation*}
(where the first map is induced by the $n$ natural inclusions $X\to X^{\boxtimes n}$) are $(G\times\Sigma_n)$-global special weak equivalences.
\begin{proof}
All the above maps are even isomorphisms for $n=0$, so we may assume that $n\ge1$. Theorem~\ref{thm:boxpower-symmetric-group} implies that the right hand map is even a $(G\times\Sigma_n)$-global \emph{level} weak equivalence, so it suffices to prove the claim for the composition $X^{\vee n}\to X^{\times n}$. For this we write $\Sigma_n^1\subset\Sigma_n$ for the subgroup of those permutations that fix $1$, and we let $p\colon G\times\Sigma_n^1\to G$ denote the projection. As in the proof of Corollary~\ref{cor:twisted-product} we have a $(G\times\Sigma_n)$-equivariant isomorphism
\begin{equation*}
\Maps^{G\times\Sigma_n^1}(G\times\Sigma_n,p^*X)\to X^{\times n}
\end{equation*}
given on the $i$-th factor by evaluating at $(1,\sigma_i^{-1})$, where $\sigma_i$ is any fixed permutation with $\sigma_i(1)=i$. Dually, we have a $(G\times\Sigma_n)$-equivariant isomorphism
\begin{equation*}
X^{\vee n}\to (G\times\Sigma_n)_+\smashp_{G\times\Sigma_n^1} p^*X
\end{equation*}
given on the $i$-th summand and in each simplicial degree by $x\mapsto[(1,\sigma_i),x]$. One then easily checks that the natural map $X^{\vee n}\to X^{\times n}$ factors as
\begin{equation*}
X^{\vee n}\xrightarrow{\cong}(G\times\Sigma_n)_+\smashp_{G\times\Sigma_n^1} p^*X\xrightarrow{\gamma}\Maps^{G\times\Sigma_n^1}(G\times\Sigma_n,p^*X)\xrightarrow{\cong} X^{\times n},
\end{equation*}
so that the claim follows from Theorem~\ref{thm:gamma-special-wirthmueller}.
\end{proof}
\end{cor}

As we have already seen in several instances above, quotients by free group actions are often fully homotopical, and in particular one can prove that quotiening out a free $H$-action sends $(G\times H)$-global special weak equivalences to $G$-global special weak equivalences. Unfortunately, this does not yet imply Theorem~\ref{thm:symmetric-products}; namely, while $\Sigma_n$ acts freely on $X^{\boxtimes n}$ for any $E\mathcal M$-simplicial set $X$ without $\mathcal M$-fixed points, in the situation of our theorem the $\Sigma_n$-action is usually not free (outside the base point) for $n\ge 3$ as an $(x_1,\dots,x_n)$ with $1<k\le n$ base point entries has non-trivial isotropy. However, this is the only thing that can go wrong:

\begin{lemma}\label{lemma:non-empty-support-free-permutation}
Let $X$ be a pointed $E\mathcal M$-simplicial set without $\mathcal M$-fixed points apart from the base point. Then the canonical $\Sigma_n$-action on $X^{\boxtimes n}$ is free outside those simplices with at least one base point component for any $n\ge0$.
\begin{proof}
Let $(x_1,\dots,x_n)$ be an $m$-simplex such that no $x_i$ is the base point. We claim that $x_i\not=x_j$ for all $i\not=j$, which will immediately imply that $(x_1,\dots,x_n)$ has trivial isotropy.

For the proof of the claim we let $1\le i<j\le n$ be arbitrary. Then $\supp(x_i)=\bigcup_{k=0}^m \supp_k(x_i)$ is non-empty by assumption, so there exists a $0\le k\le m$ with $\supp_k(x_i)\not=\varnothing$. On the other hand, $\supp_k(x_j)\cap\supp_k(x_i)=\varnothing\not=\supp_k(x_i)$, hence $\supp_k(x_j)\not=\supp_k(x_i)$, and hence in particular $x_j\not=x_i$ as desired.
\end{proof}
\end{lemma}

Using this, we can now salvage the above argument by exploiting the filtration of $X^{\boxtimes n}$ by the number of base point components. For this we will need the following relative version of the above `free quotient' heuristic:

\begin{lemma}\label{lemma:free-quotient-relative}
Let $f\colon X\to Y$ be a $(G\times H)$-global special weak equivalence in $\cat{$\bm\Gamma$-$\bm{E\mathcal M}$-$\bm{(G\times H)}$-SSet}^\tau_*$. Assume that $f$ is an injective cofibration and that $H$ acts freely on $Y(S_+)$ outside the image of $f(S_+)$ for all finite sets $S$. Then $f/H\colon X/H\to Y/H$ is a $G$-global special weak equivalence.
\begin{proof}
We factor $f$ as a $(G\times H)$-global special acyclic cofibration $i$ followed by a fibration $q$ (automatically acyclic). Proposition~\ref{prop:alpha-shriek-gamma-special-homotopical} implies that $i/H$ is a $G$-global special weak equivalence, so to finish the proof it suffices that also $q/H$ is.

We claim that $q/H$ is even a $G$-global level weak equivalence. For this we pick a finite set $S$ and a universal subgroup $K\subset\mathcal M$; we have to show that $q(S_+)/H$ is a $\mathcal G_{K,G\times\Sigma_S}$-weak equivalence, for which we make the following crucial observation:

\begin{claim*}
$q(S_+)$ is an $\mathcal F$-weak equivalence, where $\mathcal F\subset\mathcal G_{K\times H,G\times\Sigma_S}$ is the collection of those subgroups $\Gamma_{L,\phi}$ ($L\subset K\times H$, $\phi\colon L\to G\times\Sigma_S$) such that $\phi(k,h)=\phi(k,h')$ for all $(k,h),(k,h')\in L$.
\begin{proof}
We first observe that $q$ is an acyclic fibration in the positive special $(G\times H)$-global model structure, hence in particular a $(G\times H)$-global level weak equivalence. Thus, $q(S_+)$ is a $(G\times H\times\Sigma_S)$-global weak equivalence, hence in particular a $\mathcal G_{K,H\times G\times\Sigma_S}$-weak equivalence.

Now let $T\mathrel{:=}\Gamma_{L,\phi}\in\mathcal F\setminus\mathcal G_{K,H\times G\times\Sigma_S}$. We claim that $q(S_+)^T$ is even an isomorphism, for which it is enough to show that both $f(S_+)^T$ and $i(S_+)^T$ are isomorphisms.

Indeed, $f(S_+)^T$ is levelwise injective by assumption, so we only have to show that it is also surjective. For this we observe that $T$ contains an element of the form $(1,h,g,\sigma)$ with $(h,g,\sigma)\not=(1,1,1)$ as it is not contained in $\mathcal G_{K,H\times G\times\Sigma_S}$, but on the other hand $(g,\sigma)=\phi(1,h)=\phi(1,1)=(1,1)$ as $T=\Gamma_{L,\phi}$ is contained in $\mathcal F$. Thus, $(1,h,g,\sigma)=(1,h,1,1)$ with $h\not=1$, and hence $T\cap H\not=1$. But $H$ and hence also $T\cap H$ acts freely outside the image of $f(S_+)$ by assumption, so any $T$-fixed simplex of $Y(S_+)$ already lies in the image of $f(S_+)$ (and hence in the image of $f(S_+)^T$ by injectivity) as desired.

On the other hand, Corollary~\ref{cor:cofibrant-pos-g-global}-$(\ref{item:cpgg-free})$ shows that also $i(S_+)$ is an injective cofibration with free $H$-action outside the image, so the same argument shows that $i(S_+)^T$ is an isomorphism, which completes the proof of the claim.
\end{proof}
\end{claim*}
We have to show that $\phi^*(q(S_+)/H)$ is a $K$-equivariant weak equivalence for each $\phi\colon K\to G\times\Sigma_S$ (where $K$ acts via the diagonal of the $K$-action via $\phi$ and the one given by restriction of the $\mathcal M$-action). However, this can be rewritten as ${\big((\phi\times H)^*q(S_+)\big)/H}$, and using that $q(S_+)$ is an $\mathcal F$-weak equivalence by the above claim, one easily checks that $(\phi\times H)^*q(S_+)$ is a $(K\times H)$-equivariant weak equivalence. To finish the proof, we now simply observe that
\begin{equation*}
(\blank)/H\colon\cat{$\bm{(K\times H)}$-SSet}\to\cat{$\bm K$-SSet}
\end{equation*}
is left Quillen for the $\mathcal A\ell\ell$-model structures, hence fully homotopical by Ken Brown's Lemma as all $(K\times H)$-simplicial sets are $\mathcal A\ell\ell$-cofibrant.
\end{proof}
\end{lemma}

The filtration of $X^{\boxtimes n}$ and $\SP^nX$ according to the number of base point components is an instance of a more general construction for tensor powers which we will now recall:

\begin{constr}
Let $\mathscr C$ be a cocomplete symmetric monoidal category and let $f\colon X\to Y$ be a morphism in $\mathscr C$. We recall for each $n\ge 1$ the $n$-cube $C_n$ and the functor $K^n(f)\colon C_n\to\mathscr C$ from Construction~\ref{constr:iterated-ppos}.

For $0\le k\le n$ we let $K^n_k(f)$ denote the subdiagram spanned by all those sets $I$ with $|I|\le k$, and we define $Q^n_k(f)\mathrel{:=}\colim K^n_k(f)$. The inclusions of diagram shapes then induce
\begin{equation*}
X^{\otimes n}\cong Q^n_0(f)\xrightarrow{i_1} Q^n_1(f)\xrightarrow{i_2}\cdots\to Q^n_{n-1}(f)\xrightarrow{i_n} Q^n_n(f)\cong Y^{\otimes n}.
\end{equation*}
\nomenclature[aQkn]{$Q^n_k$}{filtration of $n^{\text{th}}$ tensor power}%
where the outer isomorphisms are induced by the structure maps corresponding to the unique terminal objects of $K^n_0$ and $K^n_n$, respectively.

The composition of these is precisely $f^{\otimes n}$, while the composite map $Q^n_{n-1}(f)\to Y^{\otimes n}$ was previously denoted $f^{\ppo n}$.

For any $0\le k\le n$ there is a $\Sigma_n$-action on $Q^n_k(f)$ induced by the $\Sigma_n$-action on $K_n$ and the symmetry isomorphisms of $\otimes$. All of the above maps are $\Sigma_n$-equivariant, and for $X^{\otimes n}, Y^{\otimes n}$ and $Q^n_{n-1}(f)$ this recovers the actions considered before.
\end{constr}

\begin{thm}[Gorchinskiy \& Guletski\u\i]
Assume that $\mathscr C$ is locally presentable and that $\otimes$ is cocontinuous in each variable. Moreover, let $f\colon X\to Y$ be any morphism in $\mathscr C$ and let $1\le k\le n$. We write $\alpha\colon\Sigma_{n-k}\times\Sigma_{k}\to\Sigma_n$ for the usual block sum embedding. Then we have a $\Sigma_{n}$-equivariant pushout square
\begin{equation*}
\begin{tikzcd}[column sep=1.1in]
\alpha_!(X^{\otimes(n-k)}\otimes Q^k_{k-1}(f))\arrow[d]\arrow[r, "\alpha_!(X^{\otimes(n-k)}\otimes f^{\ppo k})"] & \alpha_!(X^{\otimes(n-k)}\otimes Y^{\otimes k})\arrow[d]\\
Q^{n}_{k-1}\arrow[r, "i_k"'] & Q^{n}_k(f)
\end{tikzcd}
\end{equation*}
and a pushout square
\begin{equation*}
\begin{tikzcd}[column sep=1.4in]
\Sym^{n-k}X\otimes(Q^{k}_{k-1}(f)/\Sigma_k)\arrow[d]\arrow[r, "\Sym^{n-k}X\otimes(f^{\ppo k}/\Sigma_k)"] & \Sym^{n-k}X\otimes\Sym^{k}Y \arrow[d]\\
Q^{n}_{k-1}(f)/\Sigma_n\arrow[r, "i_k/\Sigma_n"'] & Q^n_k(f)/\Sigma_n.
\end{tikzcd}
\end{equation*}
\begin{proof}
The assumptions guarantee that $\mathscr C$ is a symmetric monoidal model category in which all maps are cofibrations (and fibrations). The claim is therefore established in \cite[proof of Theorem~22]{sym-powers}. (In fact, going through their proof one only needs the existence of pushouts and that the tensor product preserves these, but we want to avoid repeating their argument.)
\end{proof}
\end{thm}

While the theorem as stated above does not directly apply to the levelwise box product on $\cat{$\bm\Gamma$-$\bm{E\mathcal M}$-$\bm G$-SSet}_*^\tau$, we can apply it to the levelwise box product on the ordinary functor category $\cat{$\bm\Gamma$-$\bm{E\mathcal M}$-$\bm G$-SSet}^\tau\mathrel{:=}\Fun(\Gamma,\cat{$\bm{E\mathcal M}$-$\bm G$-SSet}^\tau)$. As the full subcategory $\cat{$\bm\Gamma$-$\bm{E\mathcal M}$-$\bm G$-SSet}_*^\tau$ is closed under all connected colimits, we then see a posteriori that we also have the corresponding pushouts there.

With this established, we can adapt the proof strategy of \cite[Corollary~23]{sym-powers} to deduce Theorem~\ref{thm:symmetric-products} from the Wirthmüller isomorphism:

\begin{proof}[Proof of Theorem~\ref{thm:symmetric-products}]
We will prove that the composite $X\to\SP^nX$ (induced by any of the $n$ canonical embeddings $X\to X^{\boxtimes n}$) is a $G$-global special weak equivalence for any $n\ge 1$. By $2$-out-of-$3$ we can then conclude that $\SP^nX\to\SP^{n+1}X$ is a $G$-global special weak equivalence, and so is $X\to\SP^\infty X$ as transfinite composition of $G$-global special weak equivalences.

It remains to prove the claim. For this we let $f\colon *\to X$ denote the inclusion of the basepoint, and we will prove more generally by induction on $n$:
\begin{enumerate}
\item For all $2\le k\le n$ and any group $H$ the map $i_k\colon Q^n_{k-1}(f)\to Q^n_k(f)$ is a $(G\times H\times\Sigma_n)$-global special weak equivalence when we let $H$ act trivially.
\item For all $2\le k\le n$ the induced map $i_k/\Sigma_n$ is a $G$-global weak equivalence.
\end{enumerate}

For $n=1$ there is nothing to prove, so assume $n\ge 2$. If $2\le k<n$, then we let $\alpha\colon\Sigma_{n-k}\times\Sigma_k\to\Sigma_n$ be the evident embedding again. By the induction hypothesis, $f^{\ppo k}$ is in particular a $(G\times H\times\Sigma_{n-k}\times\Sigma_k)$-global special weak equivalence for any group $H$. Moreover, it is an injective cofibration by Corollary~\ref{cor:iterated-ppo-injective} (where we have again used that we can form the corresponding colimit in $\cat{$\bm\Gamma$-$\bm{E\mathcal M}$-$\bm G$-SSet}^\tau$), hence also $\alpha_!(f^{\ppo k})$ is an injective cofibration. On the other hand, Proposition~\ref{prop:alpha-shriek-gamma-special-homotopical} shows that it is a $(G\times H\times\Sigma_n)$-global special weak equivalence. Applying Lemma~\ref{lemma:pushout-injective-infty} to the pushout
\begin{equation*}
\begin{tikzcd}[column sep=large]
\alpha_!Q^k_{k-1}(f)\arrow[r, "\alpha_!(f^{\ppo k})"]\arrow[d]&\alpha_!(X^{\boxtimes k})\arrow[d]\\
Q^{n}_{k-1}(f)\arrow[r, "i_k"'] & Q^n_k(f)
\end{tikzcd}
\end{equation*}
from the above discussion (where we used that $*^{\boxtimes(n-k)}$ is terminal) therefore shows that $i_k$ is a $(G\times H)$-global special weak equivalence as desired.

Similarly, we deduce part (2)~for $2\le k<n$ from the induction hypothesis together with the pushout
\begin{equation*}
\begin{tikzcd}[column sep=large]
Q^k_{k-1}(f)/\Sigma_k\arrow[r, "f^{\ppo k}/\Sigma_k"]\arrow[d] & \SP^kX\arrow[d]\\
Q^{n}_{k-1}(f)/\Sigma_n\arrow[r, "i_k/\Sigma_n"'] & Q^{n}_k(f)/\Sigma_n,
\end{tikzcd}
\end{equation*}
so it only remains to treat the case $k=n$ of both statements.

For the first statement we observe that the composition
\begin{equation*}
X^{\vee n}\cong Q^n_1(f)\to Q^n_2(f)\to\cdots Q^n_{n-1}(f)\to Q^n_n(f)\cong X^{\boxtimes n},
\end{equation*}
where the unlabelled isomorphism on the left is induced on the $i$-th summand by the iterated unitality isomorphism $X\cong K_n(\{i\})$, is precisely the canonical map $X^{\vee n}\to X^{\boxtimes n}$, hence a $(G\times H\times\Sigma_n)$-global special weak equivalence by Corollary~\ref{cor:boxtimes-special-power} applied to $X$ with trivial $H$-action. On the other hand, we have seen above that all the maps in $Q_n^1(f)\to Q^n_2(f)\to\cdots\to Q^n_{n-1}(f)$ are $(G\times H\times\Sigma_n)$-global special weak equivalences, so also $Q^n_{n-1}(f)\to Q^n_n(f)$ is a $(G\times H\times\Sigma_n)$-global special weak equivalence by $2$-out-of-$3$.

For the second statement, we observe that $i_n$ is a $(G\times\Sigma_n)$-global weak equivalence by the first statement, and an injective cofibration as seen above. On the other hand, the image of
$i_n(S_+)\colon Q^n_{n-1}(f)(S_+)\to X^{\boxtimes n}(S_+)$ contains all $n$-tuples of simplices with at least one basepoint component by construction, so the $\Sigma_n$-action is free outside the image of $i_n(S_+)$ by Lemma~\ref{lemma:non-empty-support-free-permutation}. We may therefore conclude from Lemma~\ref{lemma:free-quotient-relative} that also $i_n/\Sigma_n$ is a $G$-global weak equivalence, which completes the proof of the theorem.
\end{proof}

\subsection{The intermediate model}
In this section we will establish the model structure on $\cat{$\bm\Gamma$-$\bm G$-ParSumSSet}_*$ and compare it to $\cat{$\bm\Gamma$-$\bm{E\mathcal M}$-$\bm G$-SSet}_*^\tau$.

\begin{constr}\index{symmetric product!SPinfinity as a Gamma-G-parsummable simplicial set@$\SP^\infty$ as a $\Gamma$-$G$-parsummable simplicial set|textbf}
If $X\in\cat{$\bm\Gamma$-$\bm{E\mathcal M}$-$\bm G$-SSet}^\tau_*$, then $\SP^\infty X$ naturally becomes an element of $\cat{$\bm\Gamma$-$\bm G$-ParSumSSet}_*$ by declaring for each finite set $S$ that the unit of $(\SP^\infty X)(S_+)$ should be the basepoint and that the composition should be induced by the canonical isomorphisms $X^{\boxtimes m}\boxtimes X^{\boxtimes n}\cong X^{\boxtimes(m+n)}$. We omit the easy verification that this lifts $\SP^\infty$ to $\cat{$\bm\Gamma$-$\bm{E\mathcal M}$-$\bm G$-SSet}^\tau_*\to\cat{$\bm\Gamma$-$\bm G$-ParSumSSet}_*$.

The structure map of $\SP^1X$ yields a natural map $\eta\colon X\to\forget\SP^\infty X$, and we have for each $Y\in\cat{$\bm\Gamma$-$\bm G$-ParSumSSet}_*$ a natural map $\epsilon\colon\SP^\infty(\forget Y)\to Y$ induced by the iterated multiplication maps $Y^{\boxtimes n}\to Y$. We omit the easy verification that $\eta$ and $\epsilon$ satisfy the triangle identities, making them into unit and counit, respectively, of an adjunction $\SP^\infty\dashv\forget$.
\end{constr}

\begin{thm}\label{thm:gamma-uc-intermediate}
There is a unique model structure on $\cat{$\bm\Gamma$-$\bm G$-ParSumSSet}_*$ in which a map is a weak equivalence or fibration if and only if it so in the special $G$-global model structure on $\cat{$\bm\Gamma$-$\bm{E\mathcal M}$-$\bm G$-SSet}^\tau_*$. Moreover, the adjunction
\begin{equation*}
\SP^\infty\colon\cat{$\bm\Gamma$-$\bm{E\mathcal M}$-$\bm G$-SSet}^\tau_*\rightleftarrows\cat{$\bm\Gamma$-$\bm G$-ParSumSSet}_* :\!\forget
\end{equation*}
is a Quillen equivalence.
\end{thm}

One might be tempted to prove the theorem by first constructing a suitable level model structure and then Bousfield localizing. However, in this approach it is not clear \emph{a priori} how the resulting weak equivalences in $\cat{$\bm\Gamma$-$\bm G$-ParSumSSet}_*$ relate to the $G$-global special weak equivalences of underlying $G$-global $\Gamma$-spaces. 

Instead, we are going to transfer the special model structure directly, for which Theorem~\ref{thm:symmetric-products} established above will be the key ingredient. However, there are some point-set level issues we will have to deal with first.

\begin{lemma}\label{lemma:boxtimes-special-homotopical}
Let $X$ be any $G$-global $\Gamma$-space. Then $X\boxtimes\blank$ preserves $G$-global special weak equivalences.
\begin{proof}
The functor $X\times\blank$ preserves $G$-global special weak equivalences by Corollary~\ref{cor:product-special-homotopical}. If now $Y$ is any $G$-global $\Gamma$-space, then the inclusion $X\boxtimes Y\hookrightarrow X\times Y$ is a $G$-global level weak equivalence by Theorem~\ref{thm:boxtimes-em-homotopical}. Moreover, this is clearly natural in $Y$, so the claim follows by $2$-out-of-$3$.
\end{proof}
\end{lemma}

\begin{prop}\label{prop:pushout-products-gamma-parsum}
Let $f\colon X\to Y$ be an injective cofibration and a $G$-global special weak equivalence in $\cat{$\bm\Gamma$-$\bm{E\mathcal M}$-$\bm G$-SSet}_*^\tau$ such that the $\mathcal M$-actions on $X(S_+)$ and $Y(S_+)$ have no fixed points apart from the base point for any finite set $S$.

Then $f^{\ppo n}/\Sigma_n$ is a $G$-global special weak equivalence.
\begin{proof}
We first observe that $\Sym^nf=\SP^nf$ is a $G$-global special weak equivalence by Theorem~\ref{thm:symmetric-products} together with $2$-out-of-$3$.

With this established, we can argue by induction on $n$ similarly to the proof of Theorem~\ref{thm:symmetric-products}. For $n\le1$ the claim is trivial; for $n\ge 2$ we will prove more generally that all the maps in
\begin{equation}\label{eq:qns-ppo-sigma-special}
Q^n_0(f)/\Sigma_n\to Q^n_1(f)/\Sigma_n\to\cdots\to Q^n_{n-1}(f)\to Q_n^n(f)/\Sigma_n
\end{equation}
are $G$-global special weak equivalences; the claim will then follow because the right hand map is conjugate to $f^{\ppo n}/\Sigma_n$.

For the proof of the claim we observe that $Q^n_{k-1}(f)/\Sigma_n\to Q^n_k(f)/\Sigma_n$ for $k<n$ is a pushout of $\SP^{n-k}X\boxtimes(f^{\ppo k}/\Sigma_k)$ as seen in the previous subsection. The latter is a $G$-global special weak equivalence by the induction hypothesis together with the previous lemma, and it is an injective cofibration by Corollary~\ref{cor:iterated-ppo-injective} together with Lemma~\ref{lemma:ppo-injective-cofibration}. Thus, also $Q^n_{k-1}(f)/\Sigma_n\to Q^n_k(f)/\Sigma_n$ is a $G$-global weak equivalence by Lemma~\ref{lemma:pushout-injective-infty}.

It only remains to consider the case $k=n$, for which we observe that the composition $(\ref{eq:qns-ppo-sigma-special})$ is conjugate to $\SP^n f$, hence a $G$-global special weak equivalence by the above observation. The claim therefore follows by $2$-out-of-$3$.
\end{proof}
\end{prop}

\begin{prop}\label{prop:pushout-filtration-sp}
Let $f\colon A\to B$ be a map in $\cat{$\bm\Gamma$-$\bm{E\mathcal M}$-$\bm G$-SSet}^\tau_*$, and let
\begin{equation}\label{diag:pushout-sp-infty}
\begin{tikzcd}[column sep=large]
\SP^\infty A\arrow[r, "\SP^\infty f"]\arrow[d] & \SP^\infty B\arrow[d]\\
X\arrow[r, "g"'] & Y
\end{tikzcd}
\end{equation}
be a pushout in $\cat{$\bm\Gamma$-$\bm G$-ParSumSSet}_*$. Then the underlying map $\forget g$ of $G$-global $\Gamma$-spaces can be written as a transfinite composition
\begin{equation*}
X=Y_0\xrightarrow{g_1} Y_1 \xrightarrow{g_2} Y_2\to\cdots\to Y_\infty=Y
\end{equation*}
where each $g_n$ fits into a pushout square
\begin{equation}\label{diag:pushout-filtration-sym}
\begin{tikzcd}[column sep=huge]
X \boxtimes Q^n_{n-1}(f)/\Sigma_n\arrow[r, "X\boxtimes f^{\ppo n}/\Sigma_n"]\arrow[d] & X\boxtimes\SP^n(B)\arrow[d]\\
Y_{n-1}\arrow[r, "g_n"'] & Y_n
\end{tikzcd}
\end{equation}
in $\cat{$\bm\Gamma$-$\bm{E\mathcal M}$-$\bm G$-SSet}^\tau_*$.
\end{prop}

The definition of the vertical maps in $(\ref{diag:pushout-filtration-sym})$ is slightly involved, but fortunately we will never need their explicit form.

\begin{proof}
We recall the left adjoint $\textbf{P}\colon\cat{$\bm\Gamma$-$\bm{E\mathcal M}$-$\bm G$-SSet}^\tau\to\cat{$\bm\Gamma$-$\bm G$-ParSumSSet}$ of the forgetful functor; explicitly, $\textbf{P}X=\coprod_{n\ge 0}\Sym^nX$ with the evident functoriality in $X$. The unit is given by $X\cong\Sym^1X\hookrightarrow\textbf{P} X$. There is a natural map $p\colon\textbf{P}X\to\SP^\infty X$ induced by the structure maps $\Sym^nX=\SP^n X\to\SP^\infty X$, and this is compatible with the adjunction units in the sense that $p\eta=\eta$.

\begin{claim*}
The naturality square
\begin{equation*}
\begin{tikzcd}[column sep=large]
\textbf{P} A\arrow[r, "\textbf{P} f"]\arrow[d, "p"'] & \textbf{P} B\arrow[d, "p"]\\
\SP^\infty A\arrow[r, "\SP^\infty f"'] & \SP^\infty B
\end{tikzcd}
\end{equation*}
is a pushout in $\cat{$\bm\Gamma$-$\bm G$-ParSumSSet}$.
\begin{proof}
We first observe that the pushout of $\SP^\infty A\gets\textbf{P} A\to\textbf{P} B$ formed in $\cat{$\bm\Gamma$-$\bm G$-ParSumSSet}$ belongs to $\cat{$\bm\Gamma$-$\bm G$-ParSumSSet}_*$ as $(\textbf{P} f)(0_+)=\textbf{P} (f(0_+))$ is an isomorphism because $f(0_+)$ is. It therefore suffices to check the universal property with respect to every $T\in\cat{$\bm\Gamma$-$\bm G$-ParSumSSet}_*$, which is an easy diagram chase using the defining properties of $\textbf{P}$ and $\SP^\infty$ as left adjoints and the compatibility of $p$ with the unit maps.
\end{proof}
\end{claim*}

As the inclusion of $\cat{$\bm\Gamma$-$\bm G$-ParSumSSet}_*$ preserves connected colimits (hence in particular pushouts), $(\ref{diag:pushout-sp-infty})$ is also a pushout in $\cat{$\bm\Gamma$-$\bm G$-ParSumSSet}$. Pasting with the pushout from the claim we therefore get a pushout
\begin{equation*}
\begin{tikzcd}[column sep=large]
\textbf{P} A\arrow[r, "\textbf{P} f"]\arrow[d] & \textbf{P} B\arrow[d]\\
X\arrow[r, "g"'] & Y
\end{tikzcd}
\end{equation*}
in $\cat{$\bm\Gamma$-$\bm G$-ParSumSSet}$. As $\cat{$\bm\Gamma$-$\bm{E\mathcal M}$-$\bm G$-SSet}^\tau$ becomes a symmetric monoidal model category with respect to the levelwise box product when we declare all maps to be both cofibrations and fibrations, \cite[Proposition~B.2]{white-cmon} shows that $\forget g$ can be written as a transfinite composition $X=Y_0\xrightarrow{g_1} Y_1\to\cdots$ in $\cat{$\bm\Gamma$-$\bm{E\mathcal M}$-$\bm G$-SSet}^\tau$ with each $g_n$ fitting into a pushout
\begin{equation}\label{diag:pushout-sym-outer}
\begin{tikzcd}[column sep=huge]
X \boxtimes Q^n_{n-1}(f)/\Sigma_n\arrow[r, "X\boxtimes f^{\ppo n}/\Sigma_n"]\arrow[d] & X\boxtimes\Sym^n(B)\arrow[d]\\
Y_{n-1}\arrow[r, "g_n"'] & Y_n
\end{tikzcd}
\end{equation}
in $\cat{$\bm\Gamma$-$\bm{E\mathcal M}$-$\bm G$-SSet}^\tau$. Note that it does not matter whether we form $Q^n_{n-1}(f)/\Sigma_n$ and $f^{\ppo n}/\Sigma_n$ in $\cat{$\bm\Gamma$-$\bm{E\mathcal M}$-$\bm G$-SSet}_*^\tau$ or in $\cat{$\bm\Gamma$-$\bm{E\mathcal M}$-$\bm G$-SSet}^\tau$ as the former is closed under connected colimits. To finish the proof it therefore suffices that $Y_n$ belongs to $\cat{$\bm\Gamma$-$\bm{E\mathcal M}$-$\bm G$-SSet}^\tau_*$ and that $(\ref{diag:pushout-sym-outer})$ is a pushout in $\cat{$\bm\Gamma$-$\bm{E\mathcal M}$-$\bm G$-SSet}^\tau_*$ for all $n\ge1$. This is easily proven by induction using the closure under pushouts.
\end{proof}

\begin{proof}[Proof of Theorem~\ref{thm:gamma-uc-intermediate}]
Let us first prove the existence of the model structure, for which we will use Crans' Transfer Criterion (Proposition~\ref{prop:transfer-criterion}).

The category $\cat{$\bm\Gamma$-$\bm G$-ParSumSSet}_*$ is locally presentable as $\cat{$\bm G$-ParSumSSet}$ is; in particular, any set of maps permits the small object argument, and we therefore only have to show that every relative $\SP^\infty(J)$-cell complex is a weak equivalence for a suitable set $J$ of generating acyclic cofibrations of the positive special model structure on $\cat{$\bm\Gamma$-$\bm{E\mathcal M}$-$\bm G$-SSet}^\tau_*$.

For this we pick any such set $J$ such that all maps in it have cofibrant sources; this is possible by Theorem~\ref{thm:special-model-structure}. As filtered colimits and weak equivalences in $\cat{$\bm\Gamma$-$\bm G$-ParSumSSet}_*$ are created in $\cat{$\bm\Gamma$-$\bm{E\mathcal M}$-$\bm G$-SSet}^\tau_*$ and since filtered colimits in the latter are homotopical, it suffices that in any pushout
\begin{equation*}
\begin{tikzcd}[column sep=large]
\SP^\infty A\arrow[r, "\SP^\infty j"]\arrow[d] & \SP^\infty B\arrow[d]\\
X\arrow[r, "g"'] & Y
\end{tikzcd}
\end{equation*}
with $j\in J$ also $g$ is a weak equivalence in $\cat{$\bm\Gamma$-$\bm G$-ParSumSSet}_*$.

By Lemma~\ref{lemma:boxtimes-special-homotopical} together with Proposition~\ref{prop:pushout-products-gamma-parsum} we see that $X\boxtimes j^{\ppo n}/\Sigma_n$ is a $G$-global special weak equivalence for every $n\ge 1$. As it is moreover an injective cofibration by the same argument as above, also every pushout of such a map is a $G$-global special weak equivalence by Lemma~\ref{lemma:pushout-injective-infty}. We therefore conclude from the previous proposition that the underlying map of $g$ can be written as a transfinite composition of $G$-global special weak equivalences, so $g$ itself is a weak equivalence as filtered colimits are homotopical in $\cat{$\bm\Gamma$-$\bm{E\mathcal M}$-$\bm G$-SSet}_*^\tau$. This completes the verification of Crans' criterion and hence of the existence of the model structure.

The forgetful functor preserves and reflects weak equivalences and fibrations by definition; to see that $\SP^\infty\dashv\forget$ is a Quillen equivalence it therefore only remains to show that the unit $X\to\forget\SP^\infty X$ is a $G$-global special weak equivalence for every cofibrant $X\in\cat{$\bm\Gamma$-$\bm{E\mathcal M}$-$\bm G$-SSet}^\tau_*$. But for any such $X$, the $\mathcal M$-action on $X(S_+)$ has no $\mathcal M$-fixed points apart from the basepoint for any finite set $S$ by Corollary~\ref{cor:cofibrant-pos-g-global}-$(\ref{item:cpgg-support})$; the claim is therefore an instance of Theorem~\ref{thm:symmetric-products}.
\end{proof}

\subsection{Proof of the comparison theorem} In this section we will complete the proof of Theorem~\ref{thm:gamma-vs-uc} by comparing $\cat{$\bm\Gamma$-$\bm{G}$-ParSumSSet}_*$ to $\cat{$\bm G$-ParSumSSet}$.

\begin{prop}\label{prop:digamma-quillen-adjunction}
The adjunction
\begin{equation}\label{eq:ev-right-Quillen-equivalence}
\Digamma\colon\cat{$\bm G$-ParSumSSet}\rightleftarrows\cat{$\bm\Gamma$-$\bm G$-ParSumSSet}_* :\!\ev
\end{equation}
is a Quillen equivalence with fully homotopical left adjoint.
\begin{proof}
The (acyclic) fibrations in $\cat{$\bm\Gamma$-$\bm G$-ParSumSSet}_*$ are created in the $G$-global special model structure on $\cat{$\bm\Gamma$-$\bm{E\mathcal M}$-$\bm G$-SSet}_*^\tau$, hence they are in particular (acyclic) fibrations in the $G$-global level model structure. On the other hand, the (acyclic) fibrations of $\cat{$\bm G$-ParSumSSet}$ are created in $\cat{$\bm{E\mathcal M}$-$\bm G$-SSet}^\tau_*$, so it immediately follows that $\ev$ is right Quillen, i.e.~$(\ref{eq:ev-right-Quillen-equivalence})$ is a Quillen adjunction.

Next, we observe that $\ev$ is homotopical in $G$-global \emph{level} weak equivalences, and hence in particular homotopical in $G$-global special weak equivalences between \emph{special} $\Gamma$-$G$-parsummable simplicial sets by Lemma~\ref{lemma:we-between-special}. As any fibrant object of the above model structure on $\cat{$\bm\Gamma$-$\bm G$-ParSumSSet}_*$ is special, we conclude that $\textbf{R}\ev$ can be computed by taking a special replacement.

On the other hand, Corollary~\ref{cor:box-power-EM} shows that $\Digamma$ is homotopical, and we have already noted in Construction~\ref{constr:digamma} that it takes values in special $\Gamma$-$G$-parsummable simplicial sets. Together with the above we conclude that the ordinary unit $X\to\ev(\Digamma X)$ already represents the derived unit for all $G$-parsummable simplicial sets; this shows that the derived unit is an isomorphism. Finally, Lemma~\ref{lemma:we-between-special} shows that $\textbf{R}\ev$ is conservative; it follows that $\Ho(\Digamma)\dashv\textbf{R}\ev$ is an adjoint equivalence, finishing the proof of the proposition.
\end{proof}
\end{prop}

\begin{proof}[Proof of Theorem~\ref{thm:gamma-vs-uc}]
\index{G-ultra-commutative monoid@$G$-ultra-commutative monoid!vs G-global Gamma-spaces@vs.~$G$-global $\Gamma$-spaces}
\index{G-global Gamma-space@$G$-global $\Gamma$-space!vs G-ultra-commutative monoids@vs.~$G$-ultra-commutative monoids}
The functor $\ev_\omega\hskip 0pt minus 1pt\colon\hskip0pt minus .8pt\cat{$\bm{\mathcal I}$\kern-.6pt-\kern-.1ptSSet}\hskip0pt minus 1.2pt\to\hskip0pt minus 1.2pt\cat{$\bm{E\hskip-.5pt\mathcal M}$\kern-.1pt-\kern-.1ptSSet\kern-.33pt}^\tau$ is strong symmetric monoidal by Proposition~\ref{prop:ev-omega-strong-sym-mon}, so $\ev_\omega\colon\cat{UCom}\to\cat{ParSumSSet}$ preserves finite coproducts. It follows easily that the canonical mate $\Digamma\circ\ev_\omega\Rightarrow\ev_\omega\circ\Digamma$ of the identity transformation
\begin{equation*}
\begin{tikzcd}
\cat{$\bm\Gamma$-$\bm G$-UCom}_*\arrow[r, "\ev_\omega"]\arrow[d, "\ev"'] & \cat{$\bm\Gamma$-$\bm G$-ParSumSSet}_*\arrow[d, "\ev"]\\
\cat{$\bm G$-UCom}\twocell[ur]\arrow[r, "\ev_\omega"'] & \cat{$\bm G$-ParSumSSet}
\end{tikzcd}
\end{equation*}
is an isomorphism. On the other hand, $\ev_\omega$ commutes with the forgetful functors on the nose, so we altogether get an isomorphism filling
\begin{equation*}
\begin{tikzcd}
\cat{$\bm G$-UCom}\arrow[d, "\digamma"']\arrow[r, "\ev_\omega"] & \cat{$\bm G$-ParSumSSet}\arrow[d, "\digamma"]\\
\cat{$\bm\Gamma$-$\bm G$-$\bm{\mathcal I}$-SSet}_*\arrow[r, "\ev_\omega"'] & \cat{$\bm\Gamma$-$\bm{E\mathcal M}$-$\bm G$-SSet}^\tau_*,
\end{tikzcd}
\end{equation*}
hence (as $\ev_\omega$ preserves and reflects specialness) also
\begin{equation*}
\begin{tikzcd}
\cat{$\bm G$-UCom}\arrow[d, "\digamma"']\arrow[r, "\ev_\omega"] & \cat{$\bm G$-ParSumSSet}\arrow[d, "\digamma"]\\
\cat{$\bm\Gamma$-$\bm G$-$\bm{\mathcal I}$-SSet}_*^{\textup{special}}\arrow[r, "\ev_\omega"'] & \cat{$\bm\Gamma$-$\bm{E\mathcal M}$-$\bm G$-SSet}^{\tau,\textup{special}}_*.
\end{tikzcd}
\end{equation*}
commutes up to natural isomorphism.

The horizontal maps in this induce equivalences on associated quasi-categories by Corollary~\ref{cor:UCom-vs-ParSumSSet} and Theorem~\ref{thm:gamma-ev-omega}, respectively. Moreover, the right hand vertical arrow induces an equivalence by Theorem~\ref{thm:gamma-uc-intermediate} together with the previous proposition. By $2$-out-of-$3$ we conclude that also the left hand vertical arrow descends to an equivalence of associated quasi-categories, which completes the proof of the theorem.
\end{proof}

Together with Theorem~\ref{thm:uc-vs-uc} and Corollary~\ref{cor:comparison-special} we immediately conclude:

\begin{cor}\index{ultra-commutative monoid!in L-Top@in $\cat{$\bm L$-Top}$!vs global Gamma-spaces@vs.~global $\Gamma$-spaces|textbf}\label{cor:ucom-vs-global-gamma}
There are preferred equivalences between
\begin{itemize}
\item the quasi-category $\cat{UCom}^\infty$ of ultra-commutative monoids in the sense of Definition~\ref{defi:ultra-commutative-monoid},
\item the quasi-category of ultra-commutative monoids in the sense of \textup{\cite{schwede-book}} with respect to $\mathcal Fin$-global weak equivalences (see Subsection~\ref{subsec:uc-vs-uc}), and
\item the quasi-category $(\cat{$\bm\Gamma$-$\bm{E\mathcal M}$-SSet}_*^{\textup{special}})^\infty$ of special $\Gamma$-$E\mathcal M$-spaces with respect to the global level weak equivalences.\qedhere\qed
\end{itemize}
\end{cor}

On the other hand, Theorem~\ref{thm:gamma-vs-uc} also provides information about the classical equvariant approach as well as its proper equivariant generalization:

\begin{thm}\label{thm:equivariant-ultracommutativity}\index{G-ultra-commutative monoid@$G$-ultra-commutative monoid!vs Gamma-G-spaces@vs.~$\Gamma$-$G$-spaces|textbf}
\index{Gamma-G-space@$\Gamma$-$G$-space!vs G-ultra-commutative monoids@vs.~$G$-ultra-commutative monoids|textbf}
There exists a preferred equivalence between
\begin{itemize}
\item the quasi-category $(\cat{$\bm\Gamma$-$\bm G$-SSet}_*^{\textup{special}})^\infty$ of special $\Gamma$-$G$-spaces (with respect to the $G$-equivariant level weak equivalences), and
\item the quasi-category $\cat{$\bm G$-UCom}^\infty$ of $G$-ultra-commutative monoids with respect to those maps $f$ such that $\und_Gf$ is a proper $G$-equivariant weak equivalence.
\end{itemize}
\end{thm}

As alluded to in the introduction of this section, the special case $G=1$ was known by \cite[Theorem~1.2]{sagave-schlichtkrull} and the usual comparison between $E_\infty$-algebras and $\Gamma$-spaces (see e.g.~\cite[Corollary~7.3]{boavida-moerdijk} for a proof in modern language). It seems that the above generalization has not appeared in the literature before, even for finite $G$.

\begin{proof}
By Theorem~\ref{thm:gamma-vs-uc}, $\digamma$ restricts to $\cat{$\bm G$-UCom}\to\cat{$\bm\Gamma$-$\bm G$-$\bm{\mathcal I}$-SSet}_*^{\textup{special}}$, and this restriction descends to an equivalence of quasi-localizations with respect to the $G$-global weak equivalences and $G$-global level weak equivalences, respectively. In particular, the composition with $\ev_{\mathcal U_G}$ exhibits the special $\Gamma$-$G$-spaces as a (Bousfield) localization of $\cat{$\bm G$-UCom}^\infty$ by Theorem~\ref{thm:special-G-global-Gamma-vs-Gamma-G}.

It only remains to show that this composition precisely inverts the $\und_G$-weak equivalences. But indeed, a map $g\colon X\to Y$ of special $\Gamma$-$G$-spaces is a $G$-equivariant level weak equivalence if and only if $g(1^+)$ is a proper $G$-equivariant weak equivalence; the claim follows as $(\ev_{\mathcal U_G}\circ\digamma)(f)(1^+)$ is actually equal to $\und_Gf=\ev_{\mathcal U_G}(f)$ for any map $f$ of $G$-ultra-commutative monoids.
\end{proof}

\begin{rk}
For a finite group $G$, there is also an operadic approach to `$G$-equivariant coherent commutativity' based on the notion of \emph{genuine $G$-$E_\infty$-algebras}. May, Merling, and Osorno \cite[10.2]{may-merling-osorno} have shown that the corresponding homotopy theory is equivalent to the homotopy theory of special $\Gamma$-$G$-spaces, whence (by the above results) also to all the other models studied here, and in particular to $G$-parsummable simplicial sets.

In fact, as a consequence of the alternative description of the box product given in Theorem~\ref{thm:operadic-box-EM}, $G$-parsummable simplicial sets turn out to be closely related to the operadic approach. We return to this connection in \cite{gnpg}, where we exploit it to compare the models studied above to a suitable (new) theory of \emph{$G$-global $E_\infty$-algebras}. In particular, this yields a direct comparison between $G$-parsummable simplicial sets and genuine $G$-$E_\infty$-algebras that does not pass through special $\Gamma$-$G$-spaces.
\end{rk}

\chapter[Stable $G$-global homotopy theory]{Stable \for{toc}{$G$}\except{toc}{\texorpdfstring{$\bm G$}{G}}-global homotopy theory}\label{chapter:stable}
In this chapter we introduce a model of \emph{stable $G$-global homotopy theory} based on looking at the usual symmetric spectra with $G$-action through a finer notion of weak equivalence than the ordinary $G$-equivariant stable weak equivalences.

We then discuss several connections to the models considered in the previous two chapters, and in particular we will prove a $G$-global version of Segal's classical \emph{Delooping Theorem}, relating $G$-global spectra to $G$-globally coherently commutative monoids.

\section[$G$-global homotopy theory of $G$-spectra]{\for{toc}{$G$}\except{toc}{\texorpdfstring{$\bm G$}{G}}-global homotopy theory of \for{toc}{$G$}\except{toc}{\texorpdfstring{$\bm G$}{G}}-spectra}
\subsection{Recollections on equivariant stable homotopy theory}
We begin by recalling \emph{symmetric spectra} \cite{hss} which will serve as the basis of our models of $G$-global homotopy theory below.

\begin{constr}
We write $\symm$\nomenclature[aSigma]{$\symm$}{indexing category for symmetric spectra} for the following $\cat{SSet}_*$-enriched category: the objects of $\symm$ are the finite sets, and if $A$ and $B$ are finite sets, then
\begin{equation*}
\Maps_{\symm}(A,B)=\bigvee_{i\colon A\to B\text{ injective}}S^{B\setminus i(A)}.
\end{equation*}
If $C$ is yet another finite set, then the composition $\Maps_{\symm}(A,B)\smashp\Maps_{\symm}(B,C)\to\Maps_{\symm}(A,C)$ is given on the summand corresponding to $i\colon A\to B$, $j\colon B\to C$ by
\begin{equation*}
S^{B\setminus i(A)}\smashp S^{C\setminus j(B)}\to S^{j(B)\setminus ji(A)}\smashp S^{C\setminus j(B)}\cong S^{C\setminus ji(A)}\hookrightarrow\Maps_{\symm}(A,C)
\end{equation*}
where the unlabelled arrow on the left is induced by $j$, the isomorphism is the canonical one, and the final map is the inclusion of the summand corresponding to the injection $ji$.
\end{constr}

\begin{defi}\index{symmetric spectrum|textbf}\index{spectrum|seeonly{symmetric spectrum}}
A \emph{symmetric spectrum} (or, by slight abuse of language, `spectrum' for short) is an $\cat{SSet}_*$-enriched functor $\symm\to\cat{SSet}_*$. We write $\cat{Spectra}$ for the $\cat{SSet}_*$-enriched category $\FUN(\symm,\cat{SSet}_*)$ of enriched functors.
\end{defi}

\begin{defi}
Let $G$ be any discrete group, possibly infinite. A \emph{$G$-spectrum} (or, more precisely, a \emph{$G$-symmetric spectrum}) is a $G$-object in $\cat{Spectra}$. We write $\cat{$\bm G$-Spectra}$ for the $\cat{SSet}_*$-enriched category of $G$-spectra.
\end{defi}

For finite $G$, Hausmann \cite{hausmann-equivariant} studied $G$-spectra as a model of $G$-equivariant stable homotopy theory in the sense of \cite{lms}. In the rest of this subsection we will recall some of his results as the $G$-global theory sometimes parallels the equivariant one. We restrict ourselves to the foundations here and will recall other results later as needed.

\begin{rk}
Strictly speaking, \cite{hss} and \cite{hausmann-equivariant} define a symmetric spectrum as a sequence of based $\Sigma_n$-simplicial sets $X_n$, $n\ge 0$, together with suitably equivariant and associative structure maps $S^m\smashp X_n\to X_{m+n}$; the equivalence to the above definition is noted for example in \cite[2.4]{hausmann-equivariant}. In what follows we will tacitly translate the results of \cite{hausmann-equivariant} to the above language where necessary.
\end{rk}

\subsubsection{Equivariant level model structures} Assume for the rest of this subsection that $G$ is finite. Before we can introduce the stable homotopy theory of $G$-spectra, we first need to consider suitable level model structures again.

\begin{prop}
There exists a unique model structure on $\cat{$\bm G$-Spectra}$ in which a map $f$ is weak equivalence or fibration if and only if $f(A)$ is a $\mathcal G_{G,\Sigma_A}$-weak equivalence or fibration, respectively, for every finite set $A$. We call this the \emph{$G$-equivariant projective level model structure} \index{G-equivariant projective level model structure@$G$-equivariant projective level model structure!on G-Spectra@on $\cat{$\bm G$-Spectra}$|textbf}\index{G-equivariant level model structure@$G$-equivariant level model structure!projective|seeonly{$G$-equivariant projective level model structure}} and its weak equivalences the \emph{$G$-equivariant level weak equivalences}.\index{G-equivariant level weak equivalence@$G$-equivariant level weak equivalence!in G-Spectra@in $\cat{$\bm G$-Spectra}$|textbf} It is proper, combinatorial, and simplicial.
\begin{proof}
\cite[Corollary~2.26]{hausmann-equivariant} and the discussion following it show that the model structure exists and that it is proper as well as cofibrantly generated (hence combinatorial). Finally, as all the relevant constructions are levelwise, one immediately checks that the cotensoring is a right Quillen bifunctor, i.e.~the above is a simplicial model category.
\end{proof}
\end{prop}

To describe the generating cofibrations we introduce the following notation:

\begin{constr}\label{constr:G_A}\index{semifree spectrum|textbf}
Let $A$ be any finite set. Then the evaluation functor $\ev_A\colon\cat{Spectra}\to\cat{$\bm{\Sigma_A}$-SSet}_*$ admits a simplicial left adjoint $G_A$ via $\cat{SSet}_*$-enriched left Kan extension. If $X$ is any pointed $\Sigma_A$-simplicial set, then we call the spectrum $G_AX$ a \emph{semifree spectrum}; explicitly, $(G_AX)(B)=\bm\Sigma(A,B)\smashp_{\Sigma_A}X$.
\end{constr}

An easy adjointness argument then shows that the maps $G_A((G\times\Sigma_A)/H_+\smashp(\del\Delta^n\hookrightarrow\Delta^n)_+)$ for $H\in\mathcal G_{G,\Sigma_A}$ form a set of generating cofibrations, and similarly for the generating acyclic cofibrations.

We will be mostly interested in a variant of the above model structure that has more cofibrations.

\begin{defi}\index{flat cofibration!in G-Spectra@in $\cat{$\bm G$-Spectra}$|textbf}
A map $f\colon X\to Y$ of spectra is called a \emph{flat cofibration} if it has the left lifting property against those $p\colon S\to T$ such that $p(A)\colon S(A)\to T(A)$ is an acyclic fibration in the $\Sigma_A$-equivariant model structure for all finite sets $A$.

If $G$ is any group, then a map $f\colon X\to Y$ of $G$-spectra is called a \emph{flat cofibration} if its underlying map of non-equivariant spectra is a flat cofibration in the above sense. A $G$-spectrum $X$ is called \emph{flat}\index{flat!G-spectrum@$G$-spectrum|textbf}\index{flat|seealso{flat cofibration}} if $0\to X$ is a flat cofibration.
\end{defi}

\begin{ex}\label{ex:coprod-S-flat}
Any projective cofibration of $G$-spectra is a flat cofibration; in particular, the sphere spectrum $\mathbb S$ is flat as a $G$-spectrum (with respect to the trivial action) for any $G$. More generally, if $T$ is any set, then $\mathbb S^{\vee T}=\bigvee_T\mathbb S$ is flat since flat cofibrations are characterized by a left lifting property.
\end{ex}

\begin{ex}\label{ex:prod-S-flat}
Let $T$ be a \emph{finite} set. Then it is a non-trivial result that also $\mathbb S^{\times T}=\prod_T\mathbb S$ is flat, see~\cite[Proposition~B.1]{equivariant-gamma}.
\end{ex}

\begin{prop}\label{prop:equivariant-flat-level-model-structure}\index{G-equivariant flat level model structure@$G$-equivariant flat level model structure!on G-Spectra@on $\cat{$\bm G$-Spectra}$}\index{G-equivariant model structure@$G$-equivariant model structure!flat level|seeonly{$G$-equivariant flat level model structure}}\index{G-equivariant level model structure@$G$-equivariant level model structure!flat|seeonly{$G$-equivariant flat level model structure}}
There exists a unique model structure on $\cat{$\bm G$-Spectra}$ in which a map $f$ is weak equivalence or fibration if and only if $f(A)$ is a weak equivalence or fibration, respectively, in the \emph{injective} $\mathcal G_{G,\Sigma_A}$-equivariant model structure for every finite set $A$. The cofibrations of this model category are precisely the flat cofibrations. We call this the \emph{$G$-equivariant flat level model structure}. It is proper, combinatorial, and simplicial.
\begin{proof}
\cite[Corollary~2.25]{hausmann-equivariant} constructs a model structure with the desired weak equivalences and fibrations, and shows that it is proper and cofibrantly generated (hence combinatorial). The same argument as above shows that this model structure is simplicial. It remains to show that the resulting cofibrations are precisely the flat cofibrations. This is mentioned in \cite[Remark~2.20]{hausmann-equivariant}, but we detail the argument here as we will need it again later.

We first observe that the model structure can be constructed as an instance of the generalized projective model structures of Proposition~\ref{prop:generalized-projective-dim}; the consistency condition is verified as \cite[Proposition~2.24]{hausmann-equivariant}. In particular, the cofibrations are those maps $f$ such that each latching map (cf.~Construction~\ref{constr:latching-general}) is a cofibration in the injective $\mathcal G_{G,\Sigma_A}$-equivariant model structure, i.e.~a cofibration of underlying simplicial sets.

As the latching maps are independent of the the group $G$ (they are defined in terms of left Kan extensions, which are compatible with passing to functor categories), this shows that a map in $\cat{$\bm G$-Spectra}$ is a cofibration in the above model structure if and only if it is so in $\cat{$\bm 1$-Spectra}$. But the cofibrations in the latter are characterized by the same left lifting property as the flat cofibrations, which completes the proof.
\end{proof}
\end{prop}

\index{G-equivariant homotopy theory@$G$-equivariant homotopy theory!stable|seeonly{$G$-equivariant stable homotopy theory}}
\index{G-equivariant stable homotopy theory@$G$-equivariant stable homotopy theory|(}
\subsubsection{Equivariant stable homotopy theory}
The main disadvantage of symmetric spectra is that the correct notion of `stable weak equivalence' turns out to be coarser than those maps inducing isomorphisms on the na\"ively defined homotopy groups---even worse, the correct weak equivalences are only obtained indirectly by Bousfield localizing at the \emph{$\Omega$-spectra}.\index{Omega-spectrum@$\Omega$-spectrum} The equivariant situation is analogous:

\begin{defi}
A $G$-spectrum $X$ is called a \emph{$G$-$\Omega$-spectrum}\index{Omega-spectrum@$\Omega$-spectrum!G-Omega-spectrum@$G$-$\Omega$-spectrum|seeonly{$G$-$\Omega$-spectrum}}\index{G-Omega-spectrum@$G$-$\Omega$-spectrum|textbf} if the following holds: for all $H\subset G$ and all finite $H$-sets $A,B$, the adjoint structure map $X(A)\to\cat{R}\Omega^BX(A\amalg B)$ is an $H$-equivariant weak equivalence. Here $H$ acts on both sides via its actions on $X$ and on $A$ and $B$ (hence on $A\amalg B$).
\end{defi}

We emphasize that while we a priori only specify $G$-spectra on finite sets without any action (`na\"ive $G$-spectra'), in the above definition we plug in general finite $H$-sets for $H\subset G$ (using the symmetric group actions), similarly to the notion of specialness for $\Gamma$-$G$-spaces (Definition~\ref{defi:equivariant-special}).

\begin{thm}\label{thm:equivariant-stable}\index{G-equivariant stable projective model structure@$G$-equivariant stable projective model structure|seeonly{$G$-equivariant projective model structure, on $\cat{$\bm G$-Spectra}$}}\index{G-equivariant model structure@$G$-equivariant model structure!projective|seeonly{$G$-equivariant projective model structure, on $\cat{$\bm G$-Spectra}$}}
There is a unique model structure on $\cat{$\bm G$-Spectra}$ whose cofibrations are the $G$-equivariant projective cofibrations and whose fibrant objects are those $G$-$\Omega$-spectra that are fibrant in the $G$-equivariant projective level model structure. This model structure is simplicial, proper, and combinatorial. We call it the \emph{$G$-equivariant stable projective model structure}, or \emph{$G$-equivariant projective model structure}\index{G-equivariant projective model structure@$G$-equivariant projective model structure!on G-Spectra@on $\cat{$\bm G$-Spectra}$|textbf}\index{G-equivariant stable projective model structure@$G$-equivariant stable projective model structure|seeonly{$G$-equivariant projective model structure, on $\cat{$\bm G$-Spectra}$}}\index{G-equivariant stable model structure@$G$-equivariant stable model structure!projective|seeonly{$G$-equivariant projective model structure, on $\cat{$\bm G$-Spectra}$}} for short.

Similarly, there is a unique model structure on $\cat{$\bm G$-Spectra}$ whose cofibrations are the flat cofibrations and whose fibrant objects are those $G$-$\Omega$-spectra that are fibrant in the $G$-equivariant flat level model structure. This model structure is again simplicial, proper, and combinatorial. We call it the \emph{$G$-equivariant stable flat model structure},\index{G-equivariant stable flat model structure@$G$-equivariant stable flat model structure|seeonly{$G$-equivariant flat model structure, on $\cat{$\bm G$-Spectra}$}} or \emph{$G$-equivariant flat model structure}\index{G-equivariant flat model structure@$G$-equivariant flat model structure!on G-Spectra@on $\cat{$\bm G$-Spectra}$|textbf}\index{G-equivariant model structure@$G$-equivariant model structure!flat|seeonly{$G$-equivariant flat model structure}}\index{G-equivariant stable model structure@$G$-equivariant stable model structure!flat|seeonly{$G$-equivariant flat model structure, on $\cat{$\bm G$-Spectra}$}} for short.
\end{thm}

As before, we can by abstract nonsense characterize the weak equivalences of these model structures as the maps $f\colon X\to Y$ such that $[f,T]\colon[Y,T]\to[X,T]$ is an isomorphism for any $G$-$\Omega$-spectrum $T$, where $[\,{,}\,]$ denotes the hom sets in the homotopy category with respect to the $G$-equivariant level weak equivalences. In particular, the two model structures have the same weak equivalences, which we call the \emph{$G$-equivariant (stable) weak equivalences}.\index{G-equivariant weak equivalence@$G$-equivariant weak equivalence!unstably|seeonly{$G$-weak equivalence}}\index{G-equivariant weak equivalence@$G$-equivariant weak equivalence!stably|textbf}\index{G-equivariant stable weak equivalence@$G$-equivariant stable weak equivalence|seeonly{$G$-equivariant weak equivalence, stably}}

\begin{proof}
Proper combinatorial model structures with the corresponding cofibrations and fibrant objects are constructed as \cite[Theorem~4.8]{hausmann-equivariant} and \cite[Theorem~4.7]{hausmann-equivariant}, respectively, and this completely determines the model structures by Proposition~\ref{prop:cofibrations-fibrant-determine}.

It only remains to show that these model structures are simplicial, for which it suffices (by the corresponding statements for the level model structures) that the pushout product $i\ppo j$ of a cofibration $i$ of simplicial sets with a $G$-equivariant flat acyclic cofibration is a $G$-equivariant weak equivalence again, which appears as \cite[Lemma~4.5]{hausmann-equivariant}.
\end{proof}

\subsubsection{Na\"ive homotopy groups}
As in the non-equivariant setting, the above weak equivalences are quite hard to grasp. On the other hand, it is often easier to show that a map $f\colon X\to Y$ of ordinary symmetric spectra induces an isomorphism of the na\"ive homotopy groups, and any such map is indeed a (non-equivariant) stable weak equivalence by \cite[Theorem~3.1.11]{hss}. This motivates looking for a generalization of the stable homotopy groups to the equivariant setting:

\begin{constr}
Let $\mathcal U$ be a complete $H$-set universe and let $Y$ be any $H$-spectrum (e.g.~the underlying $H$-spectrum of a $G$-spectrum $X$ with $H\subset G$). We define
\begin{equation}\label{eq:non-neg-htpy-group}
\pi_k^{\mathcal U}(Y)=\colim\limits_{A\in s(\mathcal U)} [S^{A\amalg\{1,\dots,k\}}, X(A)]_*^H
\end{equation}
for $k\ge 0$ and
\begin{equation*}
\pi_k^{\mathcal U}(Y)=\colim\limits_{A\in s(\mathcal U)} [S^{A}, X(A\amalg\{1,\dots,-k\})]_*^H
\end{equation*}
for $k<0$, where $s(\mathcal U)$\nomenclature[asU]{$s(\mathcal U)$}{(filtered) poset of finite $H$-subsets of the complete $H$-set universe $\mathcal U$} is the filtered poset of finite $H$-subsets of $\mathcal U$ and $[\,{,}\,]_*^H$ denotes the set of maps in the based $H$-equivariant homotopy category. The transition map of $(\ref{eq:non-neg-htpy-group})$ for an inclusion $A\subset B$ is given by
\begin{align*}
[S^{A\amalg\{1,\dots,k\}}, X(A)]_*^H&\xrightarrow{S^{B\setminus A}\smashp\blank} [S^{B\amalg\{1,\dots,k\}}, S^{B\setminus A}\smashp X(A)]_*^H\\
&\xrightarrow{[S^{B\amalg\{1,\dots,k\}},\sigma]^H_*} [S^{B\amalg\{1,\dots,k\}}, X(B)]^H_*,
\end{align*}
where we have secretly identified $S^{B\setminus A}\smashp S^{A\amalg\{1,\dots,k\}}\cong S^{B\amalg\{1,\dots,k\}}$ via the obvious isomorphism. The definition of the transition maps for $k<0$ is similar.
\end{constr}

\begin{notation}
For the purposes of this monograph, we fix for every finite group $H$ once and for all the complete $H$-set universe $\mathcal U_H$ as in $(\ref{eq:proper-set-universe})$, and we abbreviate $\pi_*^H\mathrel{:=}\pi_*^{\mathcal U_H}$.\nomenclature[apiH]{$\pi_*^H$}{$H$-equivariant stable homotopy groups of a $G$-spectrum ($H\subset G$)}
\end{notation}

\begin{defi}
A map $f\colon X\to Y$ of $G$-spectra is called a \emph{$\underline\pi_*$-isomorphism}\index{pi-isomorphism@$\underline\pi_*$-isomorphism!G-equivariantly@$G$-equivariantly|textbf} if $\pi_*^Hf$ is an isomorphism of $\mathbb Z$-graded abelian groups for every subgroup $H\subset G$.
\end{defi}

Analogously to the non-equivariant situation we have:

\begin{thm}\label{thm:equivariant-pi-star}
Any $\underline\pi_*$-isomorphism of $G$-spectra is a $G$-equivariant stable weak equivalence.
\begin{proof}
See \cite[Theorem~3.36]{hausmann-equivariant}.
\end{proof}
\end{thm}

\begin{rk}
As we already know, any two complete $H$-set universes $\mathcal U,\mathcal U'$ are isomorphic, and any choice of such isomorphism yields a natural isomorphism
$\pi_*^{\mathcal U}\cong\pi_*^{\mathcal U'}$; in particular, the notion of $\underline\pi_*$-isomorphism is independent of any choices. However, this isomorphism is \emph{not} canonical, and in a precise sense this non-canonicity captures the failure of $\underline\pi_*$ to compute the `true' homotopy groups, see~\cite[3.3--3.6]{hausmann-equivariant}.
\end{rk}

\subsubsection{Stability} The stable model structures from Theorem~\ref{thm:equivariant-stable} are indeed stable in the model categorical sense, i.e.~the suspension/loop adjunction on the homotopy category is an equivalence. We will need the following point-set level strengthening of this later:

\begin{prop}\label{prop:equivariant-loop-suspension}
The adjunction
\begin{equation*}
\Sigma=S^1\smashp\blank\colon\cat{$\bm G$-Spectra}_{\textup{$G$-equiv.~proj.}}\rightleftarrows\cat{$\bm G$-Spectra}_{\textup{$G$-equiv.~proj.}} : \Maps(S^1,\blank)=\Omega
\end{equation*}
\nomenclature[aSigma]{$\Sigma$}{suspension, $S^1\smashp\blank$}
\nomenclature[aOmega]{$\Omega$}{loop space, based mapping space $\Maps(S^1,\blank)$}
is a Quillen equivalence. Moreover:
\begin{enumerate}
\item $\Sigma$ preserves $G$-equivariant weak equivalences.
\item $\Omega$ preserves $G$-equivariant weak equivalences between $G$-spectra that are fibrant in the $G$-equivariant projective \emph{level} model structure.
\end{enumerate}
\begin{proof}
The above adjunction is a Quillen adjunction as the $G$-equivariant projective model structure is simplicial. Moreover, $S^1\smashp\blank$ obviously preserves $G$-equivariant level weak equivalences; as any $G$-equivariant weak equivalence factors as an acyclic cofibration followed by an $G$-equivariant level weak equivalence, this shows that $\Sigma$ is in fact homotopical.

To see that it is a Quillen equivalence, we consider the category $\cat{$\bm G$-Spectra}_{\cat{Top}}$ of symmetric $G$-spectra in topological spaces, i.e.~$G$-objects in $\FUN(\bm\Sigma,\cat{Top}_*)$. This again comes with a notion of $G$-equivariant weak equivalences \cite[Definition~2.35]{hausmann-equivariant} such that the adjunction $|\blank|\dashv\Sing$ preserves and reflects weak equivalences \cite[Proposition~2.38]{hausmann-equivariant}. As unit and counit of this are even $G$-equivariant level weak equivalences, we conclude that $|\blank|\dashv\Sing$ induces an equivalence between the corresponding homotopy categories. Now consider the diagram
\begin{equation*}
\begin{tikzcd}
\cat{$\bm G$-Spectra}\arrow[d, "|\blank|"']\arrow[r, "S^1\smashp\blank"] & \cat{$\bm G$-Spectra}\arrow[d, "|\blank|"]\\
\cat{$\bm G$-Spectra}_{\cat{Top}}\arrow[r, "|S^1|\smashp\blank"'] & \cat{$\bm G$-Spectra}_{\cat{Top}}
\end{tikzcd}
\end{equation*}
of homotopical functors, commuting up to natural isomorphism. By the above, the vertical arrows induce equivalences of homotopy categories, and so does the lower arrow by \cite[Proposition~4.9]{hausmann-equivariant}. By $2$-out-of-$3$, also the top horizontal arrow descends to an equivalence, i.e.~$\Sigma\dashv\Omega$ is a Quillen equivalence.

It remains to show that $\Omega$ sends any $G$-equivariant weak equivalence $f\colon X\to Y$ of projectively \emph{level} fibrant $G$-spectra to a $G$-equivariant weak equivalence. For this we consider the naturality square
\begin{equation*}
\begin{tikzcd}
X\arrow[r, "\eta"]\arrow[d, "f"'] & \Sing|X|\arrow[d, "\Sing|f|"]\\
Y\arrow[r, "\eta"'] & \Sing|Y|;
\end{tikzcd}
\end{equation*}
as the horizontal arrows are $G$-equivariant level weak equivalences, $\Sing|f|$ is a $G$-equivariant weak equivalence. On the other hand, $\Sing Z$ is obviously level fibrant for any $Z\in\cat{$\bm G$-Spectra}_{\cat{Top}}$; as $\Omega$ preserves $G$-equivariant \emph{level} weak equivalences between level fibrant $G$-spectra, it therefore suffices to show that $\Omega\Sing|f|$, which is conjugate to $\Sing\Omega|f|$, is a $G$-equivariant weak equivalence. But indeed, $|f|$ is a $G$-equivariant weak equivalence, so $\Omega|f|$ is a $G$-equivariant weak equivalence by \cite[Proposition~4.2-(3)]{hausmann-equivariant} (for $V=\mathbb R$ the $1$-dimensional trivial $G$-representation), and hence so is $\Sing\Omega|f|$ as desired.
\end{proof}
\end{prop}
\index{G-equivariant stable homotopy theory@$G$-equivariant stable homotopy theory|)}

\subsection{$\bm G$-global model structures}
In \cite[Theorem~2.17]{hausmann-global}, Hausmann introduced a \emph{global model structure} on the category $\cat{Spectra}$ of ordinary symmetric spectra analogous to Schwede's global model structure on orthogonal spectra \cite[Theorem~4.3.17]{schwede-book}. We will now generalize this construction to yield \emph{$G$-global model structures} on $\cat{$\bm G$-Spectra}$ for any discrete $G$ (possibly infinite); while we cast them differently, our arguments for this are mostly analogous to Hausmann's.

\subsubsection{Level model structures}
Before we can construct the desired model structures on $\cat{$\bm G$-Spectra}$ modelling stable $G$-global homotopy theory, we again need to introduce several models with a stricter notion of weak equivalence:

\begin{defi}
A map $f$ of $G$-spectra is called a \emph{$G$-global level weak equivalence}\index{G-global level weak equivalence@$G$-global level weak equivalence!in G-Spectra@in $\cat{$\bm G$-Spectra}$|textbf} or \emph{$G$-global projective level fibration} if $f(A)$ is a $\mathcal G_{\Sigma_A,G}$-weak equivalence or $\mathcal G_{\Sigma_A,G}$-fibration, respectively, for every finite set $A$, i.e.~for every finite group $H$ acting faithfully on $A$ and every homomorphism $\phi\colon H\to G$ the map $f(A)^\phi$ is a weak homotopy equivalence or Kan fibration, respectively.
\end{defi}

\begin{prop}\label{prop:projective-level}
The $G$-global level weak equivalences and projective level fibrations are part of a unique model structure on $\cat{$\bm G$-Spectra}$, which we call the \emph{$G$-global projective level model structure}\index{G-global projective level model structure@$G$-global projective level model structure!on G-Spectra@on $\cat{$\bm G$-Spectra}$|textbf}\index{G-global level model structure@$G$-global level model structure!projective|seeonly{$G$-global projective level model structure}}\index{G-global model structure@$G$-global model structure!projective level|seeonly{$G$-global projective level model structure}}. It is right proper and moreover combinatorial with generating cofibrations
\begin{equation*}
\big\{G_A\big(((G\times\Sigma_A)/H\times\del\Delta^n)_+\big)\hookrightarrow G_A\big(((G\times\Sigma_A)/H\times\Delta^n)_+\big) :
n\ge 0, H\in\mathcal G_{\Sigma_A,G}\big\}
\end{equation*}
and generating acyclic cofibrations
\begin{equation*}
\big\{G_A\big(((G\times\Sigma_A)/H\times\Lambda^n_k)_+\big)\hookrightarrow G_A\big(((G\times\Sigma_A)/H\times\Delta^n)_+\big) :
0\le k\le n, H\in\mathcal G_{\Sigma_A,G}\big\}.
\end{equation*}
Finally, filtered colimits in this model category are homotopical.
\end{prop}

We will see in Lemma~\ref{lemma:levelwise-pushout} below that this model structure is also left proper.

\begin{warn}
If $G$ is finite, then being a $G$-global projective level fibration or weak equivalence is in some sense orthogonal to being a $G$-equivariant projective level fibration or weak equivalence, respectively: in level $A$, the former is a condition on $H$-fixed points for $H\in\mathcal G_{\Sigma_A,G}$, while the latter one is a condition on $H$-fixed points for $H\in\mathcal G_{G,\Sigma_A}$.
\end{warn}

\begin{proof}[Proof of Proposition~\ref{prop:projective-level}]
We equip $\cat{$\bm{(G\times\Sigma_A)}$-SSet}_*$ with the $\mathcal G_{\Sigma_A,G}$-model structure for every finite set $A$. In order to construct the desired model structure on $\cat{$\bm G$-Spectra}$, to show that it is cofibrantly generated (hence combinatorial) with the above sets of generating cofibrations and generating acyclic cofibrations, and to verify the above characterization of the cofibrations, it then suffices to verify that these model categories satisfy the `consistency condition' of Proposition~\ref{prop:generalized-projective-dim}, i.e.~that for all finite sets $A,B$ with $|A|\le |B|$ and any acyclic cofibration $i$ in the $\mathcal G_{\Sigma_A,G}$-model structure on $\cat{$\bm{(G\times\Sigma_A)}$-SSet}_*$ any pushout of $\bm\Sigma(B,A)\smashp_{\bm\Sigma(A,A)}i$ is a $\mathcal G_{\Sigma_B,G}$-weak equivalence. For this it is again enough to show that this is an acyclic cofibration in the \emph{injective} $\mathcal G_{\Sigma_B,G}$-equivariant model structure for every generating acyclic cofibration, which is obvious.

It remains to show that this model structure is right proper, simplicial (i.e.~that the cotensoring is a right Quillen bifunctor), and that filtered colimits in it are homotopical. However, these follow easily from the corresponding statements for $\cat{$\bm{(G\times\Sigma_A)}$-SSet}$ with varying $A$ as all the relevant constructions are defined levelwise.
\end{proof}

The $G$-global projective level model structure has few cofibrant objects; while this will be useful in several arguments, an unfortunate side effect is that many examples we care about in practice are not projectively cofibrant:

\begin{rk}\label{rk:projective-free}
We claim that the sphere spectrum $\mathbb S$ is not $G$-globally projectively cofibrant unless $G=1$. Indeed, if $A,B$ are any finite sets, then $\Sigma_A$ acts freely on $\bm\Sigma(A,B)$ outside the base point because it freely permutes the wedge summands. We therefore conclude from the explicit description of $G_A$ that each of the standard generating cofibrations is of the form $X\smashp(\del\Delta^n\hookrightarrow\Delta^n)_+$ for some $G$-spectrum $X$ on which $G$ acts levelwise freely outside the base point.

By cell induction one then easily concludes that $G$ acts levelwise freely outside the basepoint on \emph{any} cofibrant object of the $G$-global projective level model structure. In particular, for $G\not=1$ the only cofibrant spectrum with trivial $G$-action is the zero spectrum.
\end{rk}

To salvage this issue, we will introduce another model structure based on the flat cofibrations:

\begin{prop}\label{prop:flat-level}
There is a unique model structure on $\cat{$\bm G$-Spectra}$ whose cofibrations are the flat cofibrations and whose weak equivalences are the $G$-global level weak equivalences. We call this the \emph{$G$-global flat level model structure}.\index{G-global flat level model structure@$G$-global flat level model structure!on G-Spectra@on $\cat{$\bm G$-Spectra}$|textbf}\index{G-global level model structure@$G$-global level model structure!flat|seeonly{$G$-global flat level model structure}}\index{G-global model structure@$G$-global model structure!flat level|seeonly{$G$-global flat level model structure}} It is right proper, simplicial, and  combinatorial with generating cofibrations
\begin{equation*}
\big\{G_A\big(((G\times\Sigma_A)/H\times\del\Delta^n)_+\big)\hookrightarrow G_A\big(((G\times\Sigma_A)/H\times\Delta^n)_+\big) :
n\ge 0, H\subset\Sigma_A\times G\big\}
\end{equation*}
Moreover, filtered colimits in it are homotopical. Finally, a map $f\colon X\to Y$ is an (acyclic) fibration if and only if $f(A)\colon X(A)\to Y(A)$ is an (acyclic) fibration in the injective $\mathcal G_{\Sigma_A,G}$-equivariant model structure on $\cat{$\bm{(G\times\Sigma_A)}$-SSet}$ for all $A$.
\begin{proof}
To construct the model structure, it is enough to show that for all finite sets $A,B$ with $|A|\le |B|$ the functor
\begin{equation}\label{eq:extension-flat}
\symm(A,B)\smashp_{\Sigma_A}\blank\colon\cat{$\bm{(G\times\Sigma_A)}$-SSet}_*\to\cat{$\bm{(G\times\Sigma_B)}$-SSet}_*
\end{equation}
is left Quillen with respect to the injective $\mathcal G_{\Sigma_A,G}$-equivariant model structure on the source and the injective $\mathcal G_{\Sigma_B,G}$-model structure on the target.

For this, we fix a finite set $C$ together with a bijection $B\cong A\amalg C$, which induces an injective group homomorphism $\alpha\colon\Sigma_A\times\Sigma_C\to\Sigma_B$. Then \cite[proof of Proposition~2.24]{hausmann-equivariant} shows that $(\ref{eq:extension-flat})$ is isomorphic to the composition
\begin{equation*}
\cat{$\bm{(G\times\Sigma_A)}$-SSet}_*\xrightarrow{\blank\smashp S^C}\cat{$\bm{(G\times\Sigma_A\times\Sigma_C)}$-SSet}_*\xrightarrow{(G\times\alpha)_!}\cat{$\bm{(G\times\Sigma_B)}$-SSet}_*
\end{equation*}
and we claim that both of these are left Quillen when we equip the middle term with the injective $\mathcal G_{\Sigma_A\times\Sigma_C,G}$-equivariant model structure. Indeed, it is obvious that the first functor preserves injective cofibrations as well as weak equivalences, and for the second functor it suffices to prove the unbased version, which is an easy application of Proposition~\ref{prop:alpha-shriek-injective}.

As in the proof of Proposition~\ref{prop:projective-level} we now get a model structure with the desired weak equivalences, (acyclic) fibrations, and generating cofibrations. Moreover, we conclude as in Proposition~\ref{prop:equivariant-flat-level-model-structure} that the cofibrations are precisely the flat cofibrations. Finally, to see that this model structure is right proper, simplicial, and that filtered colimits in it are homotopical, one argues as in the projective situation above.
\end{proof}
\end{prop}

Again, we will prove later that this model structure is also left proper.

\begin{rk}\label{rk:strong-level-equivalence}
It is clear that if $p$ is an acyclic fibration in the $G$-global flat model structure, then $p(A)$ is in particular a $(G\times\Sigma_A)$-weak equivalence for all finite sets $A$. Such maps will become useful at several points below and we call them \emph{strong level weak equivalences}.\index{strong level weak equivalence|textbf} The factorization axiom for the above model structure then in particular shows that any map of $G$-spectra factors as a $G$-global flat cofibration followed by a strong level weak equivalence.
\end{rk}

\begin{warn}
While the $G$-global flat level model structure has the same cofibrations as the $G$-equivariant flat level model structure (if $G$ is finite), its weak equivalences are once again very different. As a drastic example, if $G=1$, then a morphism $f$ of ordinary spectra is a $G$-equivariant level weak equivalence if and only if it is levelwise an underlying weak equivalence while it is a $G$-global level weak equivalence if and only if each $f(A)$ is a $\Sigma_A$-weak equivalence. For general $G$, the two notions of weak equivalence are incomparable.
\end{warn}

In particular, the fibrations of the $G$-global flat level model structure are very different from the ones of the $G$-equivariant flat level model structure---the cofreeness conditions are once again `orthogonal' to each other. However we have:

\begin{lemma}\label{lemma:flat-implies-equiv-proj}
Let $H$ be a finite group, and let $\phi\colon H\to G$ be any group homomorphism. Then the simplicial adjunction
\begin{equation*}
\phi_!\colon\cat{$\bm H$-Spectra}_{\textup{$H$-equivariant projective level}}\rightleftarrows\cat{$\bm G$-Spectra}_{\textup{$G$-global flat level}} :\!\phi^*
\end{equation*}
is a Quillen adjunction. In particular, if $G$ is finite, then any fibration or acyclic fibration in the $G$-global flat level model structure is also a fibration or acyclic fibration, respectively, in the $G$-equivariant \emph{projective} level model structure.
\begin{proof}
It suffices to prove the first statement. As this is canonically a simplicial adjunction, it only remains to prove that $\phi^*$ is right Quillen.

A map $p$ of $G$-spectra is a fibration or acyclic fibration in the $G$-global flat level model structure if and only if for each finite set $A$ the map $p(A)$ is a fibration or acyclic fibration, respectively, in the injective $\mathcal G_{\Sigma_A,G}$-model structure; in particular, $p(A)$ is a fibration or acyclic fibration in the $\mathcal A\ell\ell$-model structure. But
\begin{equation*}
\phi^*\colon\cat{$\bm{(\Sigma_A\times G)}$-SSet}\to\cat{$\bm{(\Sigma_A\times H)}$-SSet}
\end{equation*}
is right Quillen with respect to the $\mathcal A\ell\ell$-model structures, so $(\phi^*p)(A)=\phi^*(p(A))$ is a fibration or acyclic fibration in the $\mathcal A\ell\ell$-model structure on $\cat{$\bm{(\Sigma_A\times H)}$-SSet}$, hence in particular in the $\mathcal G_{H,\Sigma_A}$-model structure as desired.
\end{proof}
\end{lemma}

\begin{lemma}\label{lemma:levelwise-pushout}
Both of the above model structures on $\cat{$\bm G$-Spectra}$ are left proper. Moreover, pushouts along injective cofibrations (i.e.~levelwise injections) preserve $G$-global level weak equivalences.
\begin{proof}
This follows immediately from left properness of the injective $\mathcal G_{\Sigma_A,G}$-equivariant model structure for any (finite) set $A$.
\end{proof}
\end{lemma}

\subsubsection{Stable model structures}
We now turn to the appropriate \emph{stable} model structures on $\cat{$\bm G$-Spectra}$ that will then model stable $G$-global homotopy theory.

\begin{defi}
A map $f\colon X\to Y$ is called a \emph{$G$-global stable weak equivalence}\index{G-global stable weak equivalence@$G$-global stable weak equivalence|seeonly{$G$-global weak equivalence, in $\cat{$\bm G$-Spectra}$}}\index{G-global weak equivalence@$G$-global weak equivalence!in G-Spectra@in $\cat{$\bm G$-Spectra}$|textbf} if for all finite groups $H$ and all group homomorphisms $\phi\colon H\to G$ the induced map $\phi^*f\colon\phi^*X\to\phi^*Y$ is an $H$-equivariant (stable) weak equivalence.
\end{defi}

For brevity, we will again drop the word `stable' and just use the term `$G$-global weak equivalence.'

\begin{ex}
For $G=1$, the above agrees with the \emph{global equivalences} of \cite[Definition~2.9]{hausmann-global}.
\end{ex}

\begin{rk}
If $G$ is finite, we can in particular take $\phi$ to be the identity homomorphism in the above definition, i.e.~the $G$-global weak equivalences are a refinement of the $G$-equivariant weak equivalences.

For infinite $G$, this is not possible, but we can still look at the inclusions of all finite subgroups. This shows that a $G$-global weak equivalence is in particular an $H$-equivariant weak equivalence after restricting to any finite subgroup $H\subset G$. We call such maps \emph{proper $G$-equivariant weak equivalences}; while we will not pursue this here, it is plausible that the proper $G$-equivariant weak equivalences are part of a model structure on $\cat{$\bm G$-Spectra}$ that is Quillen equivalent to the model structure on orthogonal $G$-spectra from \cite[Theorem~1.2.22]{proper-equivariant} modelling proper $G$-equivariant homotopy theory.
\end{rk}

\begin{rk}
By the $2$-out-of-$6$ property for the $H$-equivariant weak equivalences, also the $G$-global weak equivalences satisfy $2$-out-of-$6$ and hence in particular $2$-out-of-$3$.
\end{rk}

In general, the $G$-global weak equivalences are hard to grasp as already the $H$-equivariant weak equivalences are quite complicated. However, as in the equivariant setting, there is a notion of \emph{$\underline\pi_*$-isomorphism} that is easier to understand and coarse enough to be useful in many practical situations:

\begin{defi}
Let $X$ be any $G$-spectrum and let $\phi\colon H\to G$ be any group homomorphism from a finite group $H$ to $G$. Then we define the \emph{$\phi$-equivariant (stable) homotopy groups of $X$}\nomenclature[apiphi]{$\pi_*^\phi$}{$\phi$-equivariant stable homotopy groups of a $G$-global spectrum ($\phi\colon H\to G$)} as the $\mathbb Z$-graded abelian group
\begin{equation*}
\pi_*^\phi(X)\mathrel{:=}\pi_*^H(\phi^*X).
\end{equation*}
If $f\colon X\to Y$ is a map of $G$-spectra, then we define $\pi_*^\phi(f)\mathrel{:=}\pi_*^H(\phi^*f)$. The map $f$ is called a \emph{$\underline\pi_*$-isomorphism}\index{pi-isomorphism@$\underline\pi_*$-isomorphism!G-globally@$G$-globally|textbf} if $\pi_*^\phi(f)$ is an isomorphism for all homomorphisms $\phi\colon H\to G$ from finite groups $H$ to $G$.
\end{defi}

\begin{lemma}\label{lemma:level-implies-pi-star}
Let $f\colon X\to Y$ be a $G$-global level weak equivalence. Then $f$ is a $\underline{\pi}_*$-isomorphism.
\begin{proof}
Let $H$ be any finite group and let $\phi\colon H\to G$ be a group homomorphism. We will only show that $\pi_0^\phi(f)$ is an isomorphism; the proof in other dimensions is analogous, but requires more notation.

For this we define $s'(\mathcal U_H)$ as the poset of all finite \emph{faithful} $H$-subsets of $\mathcal U_H$. Then $s'(\mathcal U_H)$ is again filtered and the inclusion $s'(\mathcal U_H)\hookrightarrow s(\mathcal U_H)$ is cofinal. Thus, it induces isomorphisms fitting into a commutative diagram
\begin{equation*}
\begin{tikzcd}[column sep=1.2in]
\colim\limits_{A\in s'(\mathcal U_H)}[S^A, (\phi^*X)(A)]^H_*\arrow[d, "\cong"']\arrow[r, "{\colim_{A} [S^A, (\phi^*f)(A)]_*^H}"] & \colim\limits_{A\in s'(\mathcal U_H)}[S^A, (\phi^*Y)(A)]^H_*\arrow[d, "\cong"]\\
\colim\limits_{A\in s(\mathcal U_H)}[S^A, (\phi^*X)(A)]^H_*\arrow[r, "\pi_0^\phi(f)"'] & \colim\limits_{A\in s(\mathcal U_H)}[S^A, (\phi^*Y)(A)]^H_*.
\end{tikzcd}
\end{equation*}
On the other hand, if $A$ is a finite faithful $H$-set, then $(\phi^*f)(A)$ is an $H$-equivariant weak equivalence. Thus, also the top horizontal arrow is an isomorphism, and hence so is $\pi_0^\phi(f)$ as desired.
\end{proof}
\end{lemma}

On the other hand, Theorem~\ref{thm:equivariant-pi-star} implies:

\begin{cor}\label{cor:pi-star-implies-stable}
Let $f$ be a $\underline\pi_*$-isomorphism (for example if $f$ is a $G$-global level weak equivalence). Then $f$ is a $G$-global weak equivalence.\qed
\end{cor}

In the equivariant setting we have seen a definition of the weak equivalences in terms of \emph{$G$-$\Omega$-spectra}, and a similar characterization exists in the global context, see \cite[Sections~2.2--2.3]{hausmann-global}. Here is a $G$-global generalization of this:

\begin{defi}
A $G$-spectrum $X$ is called a \emph{$G$-global $\Omega$-spectrum}\index{Omega-spectrum@$\Omega$-spectrum!G-global Omega-spectrum@$G$-global $\Omega$-spectrum|seeonly{$G$-global $\Omega$-spectrum}}\index{G-global Omega-spectrum@$G$-global $\Omega$-spectrum|textbf} if the following holds: for any finite group $H$, any group homomorphism $\phi\colon H\to G$, any finite \emph{faithful} $H$-set $A$, and any finite $H$-set $B$ the adjoint structure map
\begin{equation*}
(\phi^*X)(A)\to {\textbf{\textup R}\Omega^B}(\phi^*X)(A\amalg B)
\end{equation*}
is an $H$-weak equivalence.
\end{defi}

An equivalent way of formulating the above is the following: for any $H,A,B$ as above the map $X(A)\to{\textbf{\textup R}\Omega^B}X(A\amalg B)$ is a $\mathcal G_{H,G}$-weak equivalence. If $X$ is $G$-globally projectively level fibrant (and hence in particular if it is fibrant in the $G$-global flat level model structure), then $(\phi^*X)(A\amalg B)$ is a fibrant $H$-simplicial set and $X(A\amalg B)$ is fibrant in the $\mathcal G_{H,G}$-model structure on $(G\times H)$-simplicial sets (both of these use that $A\amalg B$ is faithful). Thus, we can in this case just work with non-derived loop space in either of the above equivalent definitions.

\begin{prop}\label{prop:G-global-intrinsic}\index{G-global weak equivalence@$G$-global weak equivalence!in G-Spectra@in $\cat{$\bm G$-Spectra}$}
A map $f\colon X\to Y$ is a $G$-global weak equivalence if and only if for every $G$-global $\Omega$-spectrum $T$ the induced map $[f, T]\colon[Y,T]\to[X,T]$ is bijective, where $[\,{,}\,]$ denotes the hom sets in the homotopy category with respect to the $G$-global level weak equivalences.
\end{prop}

The proof will be given later once we have constructed the relevant model structures.

\begin{lemma}\label{lemma:G-Omega-sufficient}
Let $H$ be a finite group and $\phi\colon H\to G$ any homomorphism. Then the simplicial adjunction
\begin{equation*}
\phi^*\colon\cat{$\bm G$-Spectra}_{\textup{$G$-global flat level}}\rightleftarrows
\cat{$\bm H$-Spectra}_{\textup{$H$-equivariant flat}} :\!\phi_*
\end{equation*}
is a Quillen adjunction with homotopical left adjoint. Moreover, the essential image of $\textbf{\textup R}\phi_*$ is contained in the class of $G$-global $\Omega$-spectra.
\begin{proof}
It is obvious that $\phi^*$ preserves flat cofibrations, and it is actually homotopical by Corollary~\ref{cor:pi-star-implies-stable}, hence in particular left Quillen. To finish the proof it suffices now that for any fibrant spectrum $X$ in the $H$-equivariant flat model structure the $G$-spectrum $\phi_*X$ is a $G$-global $\Omega$-spectrum, which by the Quillen adjunction just established amounts to saying that the adjoint structure map
\begin{equation}\label{eq:GOs-adjoint-structure-phi-star-X}
(\phi_*X)(A)\to\Omega^B(\phi_*X)(A\amalg B)
\end{equation}
is a $\mathcal G_{K,G}$-weak equivalence for every finite group $K$, finite faithful $K$-set $A$, and finite $K$-set $B$. We will show that it is in fact even a $(K\times G)$-weak equivalence.

Indeed, by \cite[proof of Proposition~2.13]{hausmann-global} the adjoint structure map
\begin{equation}\label{eq:GOs-adjoint-structure-X}
X(A)\to\Omega^B X(A\amalg B)
\end{equation}
is actually a $(K\times H)$-weak equivalence. Moreover, both sides are easily seen to be fibrant in the $\mathcal A\ell\ell$-model structure on $\cat{$\bm{(K\times H)}$-SSet}$. As $\phi^*\colon\cat{$\bm{(K\times G)}$-SSet}\to\cat{$\bm{(K\times H)}$-SSet}$ is left Quillen for the $\mathcal A\ell\ell$-model structures on both sides, its right adjoint $\phi_*$ is right Quillen, so it sends $(\ref{eq:GOs-adjoint-structure-X})$ to a $(K\times G)$-weak equivalence by Ken Brown's Lemma. This image agrees with $(\ref{eq:GOs-adjoint-structure-phi-star-X})$ up to conjugation by isomorphisms (as $\phi_*$ preserves cotensors), which completes the proof of the lemma.
\end{proof}
\end{lemma}

\begin{lemma}\label{lemma:G-global-between-G-Omega}
Let $f\colon X\to Y$ be a map of $G$-global $\Omega$-spectra. Then $f$ is a $G$-global level weak equivalence if and only if $f$ is a $G$-global weak equivalence.
\begin{proof}
The implication `$\Rightarrow$' holds by Corollary~\ref{cor:pi-star-implies-stable} even without any assumptions on $X$ and $Y$.

For the implication `$\Leftarrow$', we may assume without loss of generality by $2$-out-of-$3$ that $X$ and $Y$ are fibrant in the $G$-global flat level model structure. Now let $H$ be any finite group, let $\phi\colon H\to G$ be a group homomorphism, and let $A$ be a finite faithful $H$-set. We consider the diagram
\begin{equation*}
\begin{tikzcd}
\phi^*X\arrow[r, "\lambda_X"]\arrow[d, "\phi^*f"'] & \Omega^A\sh^A\phi^*X\arrow[d, "\Omega^A\sh^A\phi^*f"]\\
\phi^*Y\arrow[r, "\lambda_Y"'] & \Omega^A\sh^A\phi^*Y
\end{tikzcd}
\end{equation*}
with $\lambda_X$ given in degree $B$ by the adjoint structure map $\phi^*X(B)\to\Omega^A(\phi^*X)(A\amalg B)$, and analogously for $Y$. Thus, $\lambda_X$ and $\lambda_Y$ are $H$-global level weak equivalences by assumption on $X$ and $Y$, and hence in particular $H$-equivariant weak equivalences by Corollary~\ref{cor:pi-star-implies-stable}. Moreover, $\phi^*f$ is an $H$-equivariant weak equivalence by definition, so that $\Omega^A\sh^A\phi^*f$ is an $H$-equivariant weak equivalence by $2$-out-of-$3$. As both its source and target are $H$-$\Omega$-spectra by definition, $\Omega^A\sh^A\phi^*f$ is an $H$-equivariant level weak equivalence \cite[Remark~2.36]{hausmann-equivariant}. In particular, if $B$ is any finite $H$-set, then $(\Omega^A\sh^A\phi^*f)(B)$ is an $H$-equivariant weak equivalence. But this fits into a commutative diagram
\begin{equation*}
\begin{tikzcd}[column sep=large]
(\phi^*X)(B)\arrow[r, "\lambda_X(B)"]\arrow[d, "(\phi^*f)(B)"'] & (\Omega^A\sh^A\phi^*X)(B)\arrow[d, "(\Omega^A\sh^A\phi^*f)(B)"]\\
(\phi^*Y)(B)\arrow[r, "\lambda_Y(B)"'] & (\Omega^A\sh^A\phi^*X)(B),
\end{tikzcd}
\end{equation*}
and if $B$ is actually \emph{faithful}, then also the horizontal arrows are $H$-weak equivalences by the above. It follows by $2$-out-of-$3$ that in this case also the left hand vertical arrow is an $H$-equivariant weak equivalence as desired.
\end{proof}
\end{lemma}

\begin{constr}
Let $H$ be a finite group, let $A$ be a finite faithful $H$-set, and let $\phi\colon H\to G$ be any homomorphism. Then enriched adjointness yields a natural isomorphism $(\blank)^\phi\circ\ev_A\cong \Maps^G\big(G_A((G\times\Sigma_A)/H_+),\blank\big)$ (where $H$ acts via $A$ and $\phi$). If $B$ is another finite $H$-set, then $G_{A\amalg B}\big((S^B\smashp (G\times\Sigma_{A\amalg B})_+)/H\big)$ similarly corepresents $T\mapsto(\Omega^BT(A\amalg B))^\phi$. Thus the Yoneda Lemma yields a map
\begin{equation*}
\lambda_{A,\phi}^B\colon G_{A\amalg B}\big((S^B\smashp (G\times\Sigma_{A\amalg B})_+)/H\big)\to G_A\big((G\times\Sigma_A)_+/H\big)
\end{equation*}
such that the restriction $\Maps^G(\lambda_{A,\phi}^B,T)$ is conjugate to the $\phi$-fixed points of the adjoint structure map $T(A)\to\Omega^B T(A\amalg B)$ for any $G$-spectrum $T$.
\end{constr}

We can now finally prove:

\begin{thm}\label{thm:G-global-flat}
There is a unique model structure on $\cat{$\bm G$-Spectra}$ whose cofibrations are the flat cofibrations and whose weak equivalences are the $G$-global weak equivalences. We call this the \emph{$G$-global flat model structure}.\index{G-global flat model structure@$G$-global flat model structure!on G-Spectra@on $\cat{$\bm G$-Spectra}$|textbf}\index{G-global model structure@$G$-global model structure!flat|seeonly{$G$-global flat model structure}} It is left proper, combinatorial, simplicial, and filtered colimits in it are homotopical. Moreover, a $G$-spectrum is fibrant in it if and only if it is fibrant in the $G$-global flat level model structure and moreover a $G$-global $\Omega$-spectrum.
\begin{proof}
Clearly, the sources and targets of the maps $\lambda_{A,\phi}^B$ from the previous construction are cofibrant in the projective $G$-global level model structure, hence in particular in the flat $G$-global level model structure. Thus, picking representatives of isomorphism classes of such $H,A,B,\phi$, we may apply Theorem~\ref{thm:local-model-structure} to obtain a left proper, simplicial, and combinatorial model structure on $\cat{$\bm G$-Spectra}$ with the same cofibrations as the flat $G$-global level model structure and whose fibrant objects are precisely those $X$ that are fibrant in the $G$-global flat level model structure and for which $\Maps^G(\lambda_{A,\phi}^B,X)$ is a weak homotopy equivalence for all $A,B,\phi$. By construction, the latter is equivalent to $X$ being a $G$-global $\Omega$-spectrum.

Let us call the weak equivalences of this model structure \emph{$G$-global $\Omega$-weak equivalences} for the time being. Lemma~\ref{lemma:filtered-still-homotopical} shows that they are stable under filtered colimits, and it only remains to show that they agree with the $G$-global weak equivalences, which amounts to saying that
\begin{equation}\label{eq:all-your-phi}
(\phi^*)_{\phi\colon H\to G}\colon\cat{$\bm G$-Spectra}_{\text{$G$-global $\Omega$-w.e.}}\to\prod_{\phi\colon H\to G}\cat{$\bm H$-Spectra}_{\textup{$H$-equiv.~w.e.}}
\end{equation}
preserves and reflects weak equivalences, where the product runs over a system of representatives of finite groups $H$ and all homomorphisms $\phi\colon H\to G$.

Indeed, let us equip the left hand side with the above model structure and the right hand side with the product of the $H$-equivariant flat model structures; we claim that $(\ref{eq:all-your-phi})$ is left Quillen, for which it suffices to show that $(\phi^*)_\phi$ preserves cofibrations and that its right adjoint (which sends a family $(X_\phi)_\phi$ of $H$-spectra to the product $\prod_{\phi}\phi_*X_\phi$) preserves fibrant objects. But as fibrant objects are stable under products, this is immediate from Lemma~\ref{lemma:G-Omega-sufficient}. As $(\phi^*)_\phi$ moreover sends $G$-global \emph{level} weak equivalences to weak equivalences by the same lemma, we conclude that it is homotopical in $G$-global $\Omega$-weak equivalences.

It therefore only remains to show that $(\ref{eq:all-your-phi})$ also reflects weak equivalences. But since we already know that it is homotopical, it suffices to show this for maps between fibrant objects, in which case the claim follows from Lemma~\ref{lemma:G-global-between-G-Omega}.
\end{proof}
\end{thm}

\begin{proof}[Proof of Proposition~\ref{prop:G-global-intrinsic}]
This is immediate from the description of the weak equivalences of the above model structure provided by Theorem~\ref{thm:local-model-structure}.
\end{proof}

\begin{thm}\label{thm:G-global-projective}
There is a unique model structure on $\cat{$\bm G$-Spectra}$ whose cofibrations are the $G$-global projective cofibrations and whose weak equivalences are the $G$-global weak equivalences. We call this the \emph{$G$-global projective model structure}.\index{G-global projective model structure@$G$-global projective model structure!on G-Spectra@on $\cat{$\bm G$-Spectra}$|textbf}\index{G-global model structure@$G$-global model structure!projective|seeonly{$G$-global projective model structure}} It is left proper, combinatorial, simplicial, and filtered colimits in it are homotopical. Moreover, a $G$-spectrum is fibrant in it if and only if it is fibrant in the $G$-global projective level model structure and moreover a $G$-global $\Omega$-spectrum.
\begin{proof}
Analogously to Theorem~\ref{thm:G-global-flat}, one constructs a model structure with the desired cofibrations and fibrant objects, and shows that it is left proper, combinatorial, simplicial, and that filtered colimits in it are homotopical. The weak equivalences of this model structure can then again be detected by mapping into $G$-global $\Omega$-spectra, so Proposition~\ref{prop:G-global-intrinsic} shows that they agree with the $G$-global weak equivalences.
\end{proof}
\end{thm}

\begin{rk}
For $G=1$ the above two model structures agree with each other (as $\mathcal G_{\Sigma_A,1}$ is already the collection of all subgroups of $\Sigma_A\times1$ for any finite set $A$), and they recover Hausmann's \emph{global stable model structure} \cite[Theorem~2.17]{hausmann-global} on symmetric spectra.

Schwede also introduced a global model structure on the category of \emph{orthogonal} spectra that contains equivariant information for all compact Lie groups, see~\cite[Theorem~4.3.17]{schwede-book}. Analogously to our comparison between $L$-spaces and $I$-simplicial sets (Theorem~\ref{thm:L-vs-I}), one can express Hausmann's global model structure as an explicit Bousfield localization of Schwede's model structure with respect to certain `$\mathcal Fin$-global weak equivalences,' see~\cite[Theorem~5.3]{hausmann-global}.
\end{rk}

\begin{lemma}\label{lemma:stable-we-coproducts}
The $G$-global weak equivalences in $\cat{$\bm G$-Spectra}$ are stable under all small coproducts.
\begin{proof}
This is clear for the $G$-global \emph{level} weak equivalences; the claim follows as the acyclic cofibrations in any model structure are stable under coproducts.
\end{proof}
\end{lemma}

\begin{rk}
If we fix a finite group $G$, then we can equip $\cat{$\bm G$-Spectra}$ with either the $G$-global flat model structure or the $G$-equivariant flat model structure. Both of these have the same cofibrations, but the $G$-global weak equivalences are strictly finer than the $G$-equivariant ones even for $G=1$.

As a model structure is characterized by its cofibrations together with its fibrant objects (Proposition~\ref{prop:cofibrations-fibrant-determine}), it follows that being fibrant in the $G$-global flat model structure is \emph{strictly} weaker than being fibrant in the $G$-equivariant flat model structure. However, we will show later in Proposition~\ref{prop:flat-fibrant-G-Omega} that a fibrant object in the $G$-global flat model structure is at least still a $G$-$\Omega$-spectrum.
\end{rk}

\begin{lemma}
Pushouts of $G$-global weak equivalences along injective cofibrations are again $G$-global weak equivalences.
\begin{proof}
This follows from Lemma~\ref{lemma:levelwise-pushout} as in the proof of Theorem~\ref{thm:G-M-injective-semistable-model-structure}.
\end{proof}
\end{lemma}

\begin{cor}
There is a unique model structure on $\cat{$\bm G$-Spectra}$ whose weak equivalences are the $G$-global weak equivalences and whose cofibrations are the injective cofibrations. We call this the \emph{injective $G$-global model structure}.\index{injective G-global model structure@injective $G$-global model structure!on G-Spectra@on $\cat{$\bm G$-Spectra}$|textbf} It is combinatorial, simplicial, left proper, and filtered colimits in it are homotopical.
\begin{proof}
Applying Corollary~\ref{cor:mix-model-structures} to the above, it only remains to show that this model structure is simplicial. To this end, it suffices as before that for each $K\in\cat{SSet}_*$ the functor $K\smashp\blank$ is homotopical and that for each $G$-spectrum $X$ so is $\blank\smashp X\colon\cat{SSet}_*\to\cat{$\bm G$-Spectra}$. But indeed, the second statement is obvious, and for the first we observe that $K\smashp\blank$ preserves $G$-global projective acyclic cofibrations (by Theorem~\ref{thm:G-global-projective}) as well as $G$-global level weak equivalences (by an easy calculation). The claim follows immediately.
\end{proof}
\end{cor}

\begin{prop}
Let
\begin{equation}\label{diag:bQ-pullback}
\begin{tikzcd}
P\arrow[r, "g"]\arrow[d] & X\arrow[d, "p"]\\
Y\arrow[r, "f"'] & Z
\end{tikzcd}
\end{equation}
be a pullback diagram in $\cat{$\bm G$-Spectra}$ such that $p$ is a $G$-global projective level fibration and $f$ is a $G$-global weak equivalence. Then $g$ is a $G$-global weak equivalence, too. In particular, the $G$-global projective model structure, the $G$-global flat model structure, and the $G$-global injective model structure are right proper (hence proper).
\begin{proof}
It suffices to prove the first statement. For this we employ the factorization axiom of the $G$-global flat \emph{level} model structure to write $f$ as a composition
\begin{equation*}
Y\xrightarrow{i}H\xrightarrow{s}Z,
\end{equation*}
where $i$ is a $G$-global level weak equivalence and $s$ is a fibration in the $G$-global flat level model structure. We observe that $s$ is moreover a $G$-global weak equivalence by Corollary~\ref{cor:pi-star-implies-stable} together with $2$-out-of-$3$. The diagram $(\ref{diag:bQ-pullback})$ now factors as
\begin{equation*}
\begin{tikzcd}
P\arrow[r, "j"]\arrow[d] & I\arrow[r, "t"]\arrow[d, "q"'] & X\arrow[d, "p"]\\
Y\arrow[r, "i"'] & H\arrow[r, "s"'] & Z
\end{tikzcd}
\end{equation*}
where both squares are pullbacks. In particular, $q$ is again a $G$-global projective level fibration and hence $j$ is a $G$-global level weak equivalence by right properness of the $G$-global projective level model structure. By another application of Corollary~\ref{cor:pi-star-implies-stable} and $2$-out-of-$3$, it is therefore enough to show that $t$ is a $G$-global weak equivalence.

For this let $\phi\colon H\to G$ be any homomorphism from a finite group $H$ to $G$; then
\begin{equation*}
\begin{tikzcd}
\phi^*I\arrow[r, "\phi^*t"]\arrow[d, "\phi^*q"'] & \phi^*X\arrow[d, "\phi^*p"]\\
\phi^*H\arrow[r, "\phi^*s"'] & \phi^*Z
\end{tikzcd}
\end{equation*}
is a pullback, and $\phi^*s$ is an $H$-equivariant projective level fibration by Lemma~\ref{lemma:flat-implies-equiv-proj} as well as an $H$-equivariant weak equivalence by the above. It therefore follows from \cite[discussion after Proposition~4.2]{hausmann-equivariant} that also $\phi^*t$ is an $H$-equivariant weak equivalence. Letting $\phi$ vary we see therefore see that $t$ is a $G$-global weak equivalence, finishing the proof.
\end{proof}
\end{prop}

Finally we observe:

\begin{prop}\label{prop:stable-model-structures-are-stable}
The $G$-global projective model structure, the $G$-global flat model structure, and the $G$-global injective model structure all are stable.
\begin{proof}
Each of these are simplicial model categories, so this amounts to saying that the suspension/loop adjunction $\Sigma\mathrel{:=}S^1\smashp\blank\dashv\Maps(S^1,\blank)\mathrel{:=}\Omega$ (which is always a Quillen adjunction) is a Quillen equivalence for each of these. As the left adjoint is homotopical and all of these have the same weak equivalences, it suffices to show:
\begin{enumerate}
\item For every $G$-spectrum $X$ and some (hence any) fibrant replacement $\Sigma X\to Y$ in the $G$-global flat model structure, the composition $X\to \Omega\Sigma X\to \Omega Y$ is a $G$-global weak equivalence.
\item For every $G$-spectrum $Y$ that is fibrant in the $G$-global flat model structure, the counit $\Sigma\Omega X\to X$ is a $G$-global weak equivalence.
\end{enumerate}
We will prove the first statement, the argument for the second one being analogous. If $f\colon\Sigma X\to Y$ is any fibrant replacement in the $G$-global flat model structure, then $\phi^*(f)$ is an $H$-equivariant weak equivalence and $\phi^*(Y)$ is fibrant in the $H$-equivariant projective \emph{level} model structure for every homomorphism $\phi\colon H\to G$. By Proposition~\ref{prop:equivariant-loop-suspension}, $\Sigma$ is also homotopical with respect to $H$-equivariant weak equivalences, and $\Omega$ preserves $H$-equivariant weak equivalences between objects that are fibrant in the $H$-equivariant projective \emph{level} model structure; we therefore conclude that $\phi^*(f\eta)$ is a model for the derived unit of the suspension loop adjunction on $\cat{$\bm H$-Spectra}_{\textup{$H$-equiv.~proj.}}$. The claim follows by another application of Proposition~\ref{prop:equivariant-loop-suspension}.
\end{proof}
\end{prop}

\subsection{Functoriality} \index{functoriality in homomorphisms!for G-Spectra@for $\cat{$\bm G$-Spectra}$|(}
We will now discuss how the above models relate to each other when the group $G$ varies.

\begin{lemma}\label{lemma:restriction-right-Quillen-projective}
For any homomorphism $\alpha\colon H\to G$, the simplicial adjunction
\begin{equation*}
\alpha_!\colon \cat{$\bm H$-Spectra}_{\textup{$H$-global projective}}\rightleftarrows
\cat{$\bm G$-Spectra}_{\textup{$G$-global projective}} :\!\alpha^*
\end{equation*}
is a Quillen adjunction, and likewise for the corresponding level model structures.
\begin{proof}
It is clear that $\alpha^*$ preserves $G$-global level weak equivalences and fibrations, proving the second statement. For the first statement it then suffices that $\alpha^*$ preserves fibrant objects, which by the above together with the characterization of fibrant objects amounts to saying that $\alpha^*$ sends $G$-globally projectively fibrant $G$-global $\Omega$-spectra to $H$-global $\Omega$-spectra. This is obvious from the definitions.
\end{proof}
\end{lemma}

\begin{lemma}\label{lemma:restriction-left-Quillen-flat}
For any homomorphism $\alpha\colon H\to G$, the simplicial adjunction
\begin{equation*}
\alpha^*\colon\cat{$\bm G$-Spectra}_{\textup{$G$-global flat}}\rightleftarrows\cat{$\bm H$-Spectra}_{\textup{$H$-global flat}}\ :\!\alpha_*.
\end{equation*}
is a Quillen adjunction.
\begin{proof}
By definition, $\alpha^*$ preserves weak equivalences and flat cofibrations.
\end{proof}
\end{lemma}

Unfortunately, $\alpha^*$ is not in general right Quillen with respect to the corresponding \emph{flat} model structures:

\begin{ex}
We let $X$ be any globally fibrant spectrum whose underlying non-equivariant stable homotopy type agrees with the Eilenberg-Mac Lane spectrum $H\mathbb Z/2$. As globally fibrant spectra are in particular non-equivariant $\Omega$-spectra (the trivial group acts faithfully on any set) and levelwise Kan complexes, we conclude that $X(*)$ represents the first singular cohomology group with coefficients in $\mathbb Z/2$, i.e.~$\pi_0\Maps(\blank, X(*))\cong H^1(\blank,\mathbb Z/2)$.

Now let $\alpha\colon\mathbb Z/2\to1$ be the unique group homomorphism. We claim that $\alpha^*X$ is not even fibrant in the $\mathbb Z/2$-global flat \emph{level} model structure: namely, this would in particular mean that
\begin{equation*}
\Maps^{\mathbb Z/2}(E({\mathbb Z/2}),(\alpha^*X)(*))\simeq\Maps^{\mathbb Z/2}(*, (\alpha^*X)(*))
\end{equation*}
as $E(\mathbb Z/2)\to *$ is a $\mathcal G_{\Sigma_*,\mathbb Z/2}$-equivariant (i.e.~underlying) weak equivalence. Thus,
\begin{align*}
\Maps(E(\mathbb Z/2)/(\mathbb Z/2),X(*))&\cong
\Maps^{\mathbb Z/2}(E(\mathbb Z/2),(\alpha^*X)(*))\\
&\simeq\Maps^{\mathbb Z/2}(*,(\alpha^*X)(*))=\Maps(*,X(*)).
\end{align*}
However, $E(\mathbb Z/2)$ is a contractible space with free $\mathbb Z/2$-action, so that its quotient by $\mathbb Z/2$ is a $K(\mathbb Z/2,1)$ and hence equivalent to $\mathbb R\text{P}^\infty$. Thus taking $\pi_0$ of the above we could conclude that $H^1(\mathbb R\text{P}^\infty,\mathbb Z/2)\cong H^1(*,\mathbb Z/2)$. However, the right hand side is trivial while the left hand side is cyclic of order $2$, yielding the desired contradiction.
\end{ex}

However, we have:

\begin{prop}\label{prop:restriction-injective-right-Quillen}
Let $\alpha\colon H\to G$ be an \emph{injective} group homomorphism. Then the simplicial adjunction
\begin{equation*}
\alpha_!\colon\cat{$\bm H$-Spectra}_{\textup{flat $H$-global}}\rightleftarrows\cat{$\bm G$-Spectra}_{\textup{flat $G$-global}} :\!\alpha^*
\end{equation*}
is a Quillen adjunction and likewise for the corresponding level model structures.
\begin{proof}
For the second statement it suffices to observe that
\begin{equation*}
\alpha_!\colon\cat{$\bm{(\Sigma_A\times H)}$-SSet}_{\textup{inj.~$\mathcal G_{\Sigma_A,H}$-equiv.}}\rightleftarrows\cat{$\bm{(\Sigma_A\times G)}$-SSet}_{\textup{inj.~$\mathcal G_{\Sigma_A,G}$-equiv.}} :\!\alpha^*,
\end{equation*}
is a Quillen adjunction for every finite set $A$ by Proposition~\ref{prop:alpha-shriek-injective}. With this established, it suffices for the first statement to prove that $\alpha^*$ sends fibrant objects to $H$-global $\Omega$-spectra. This holds in fact for \emph{any} group homomorphism $\alpha$ by Lemma~\ref{lemma:restriction-right-Quillen-projective}.
\end{proof}
\end{prop}

The same example as above shows that $\alpha^*$ is not right Quillen with respect to the corresponding injective model structures for general $\alpha$. However, we still have:

\begin{prop}\label{prop:alpha-shriek-homotopical}
Let $\alpha\colon H\to G$ be an \emph{injective} group homomorphism. Then the simplicial adjunction
\begin{equation*}
\alpha_!\colon\cat{$\bm H$-Spectra}_{\textup{injective $H$-global}}\rightleftarrows\cat{$\bm G$-Spectra}_{\textup{injective $G$-global}} :\!\alpha^*
\end{equation*}
is a Quillen adjunction.
\begin{proof}
It is obvious that $\alpha_!$ preserves levelwise injections. To prove that it is homotopical, we use the factorization axiom in the flat $H$-global model structure to factor any given $H$-global weak equivalence $f$ as an acyclic $H$-global flat cofibration $i$ followed by a strong level weak equivalence $p$. By Proposition~\ref{prop:restriction-injective-right-Quillen}, $\alpha_!(i)$ is a $G$-global weak equivalence, and so is $\alpha_!(p)$ for obvious reasons.
\end{proof}
\end{prop}

The functor $\alpha_!$ preserves injective cofibrations in full generality, so the failure of $\alpha^*$ to be right Quillen with respect to the injective model structures can be reinterpreted as a failure of $\alpha_!$ to be homotopical. However, similarly to the unstable situation we have:

\begin{prop}\label{prop:spectra-free-quotient}
Let $\alpha\colon H\to G$ be any homomorphism, and let $f\colon X\to Y$ be an $H$-global weak equivalence such that $\ker(\alpha)$ acts levelwise freely on $X$ and $Y$ outside the base point. Then $\alpha_!f$ is a $G$-global weak equivalence.
\begin{proof}
Lemma~\ref{lemma:restriction-right-Quillen-projective} in particular implies that $\alpha_!$ preserves weak equivalences between $H$-globally projectively cofibrant objects. Choosing functorial factorizations in the $H$-global projective model structure, it therefore suffices to show: if $p\colon X'\to X$ is any $H$-global \emph{level} weak equivalence such that $X'$ is cofibrant in the $H$-global projective model structure, then $\alpha_!(p)$ is a $G$-global weak equivalence.

Indeed, if $A$ is any finite set, then $p(A)$ is a $\mathcal G_{\Sigma_A,H}$-weak equivalence, and $\ker(\alpha)$ acts freely outside the basepoint on $X(A)$ by assumption. On the other hand, as seen in Remark~\ref{rk:projective-free}, all of $H$ acts freely on $X'(A)$, so the claim is an instance of Remark~\ref{rk:free-quotient-based}.
\end{proof}
\end{prop}
\index{functoriality in homomorphisms!for G-Spectra@for $\cat{$\bm G$-Spectra}$|)}

\subsection{Additivity and the Wirthmüller isomorphism}
The stability of the $G$-global model structures in particular implies that the $G$-global stable homotopy category is additive, so that finite coproducts and products agree in it. The following lemma provides an underived version of this:

\begin{lemma}\label{lemma:underived-additivity}
Let $T$ be a finite set and let $(X_t)_{t\in T}$ be any family of $G$-spectra. Then the canonical map $\bigvee_{t\in T} X_t\to\prod_{t\in T} X_t$ is a $G$-global weak equivalence (and in fact a $\underline\pi_*$-isomorphism).
\begin{proof}
Let $\phi\colon H\to G$ be a group homomorphism from a finite group $H$ to $G$. Then restricting the canonical map along $\phi$ agrees with the canonical map $\bigvee_{t\in T} \phi^*(X_t)\to\prod_{t\in T} \phi^*(X_t)$
(on the nose if we use the usual construction of colimits, up to conjugation by isomorphism in general). The claim follows as the latter is an equivariant $\underline\pi_*$-isomorphism by \cite[Proposition~3.6-(3)]{hausmann-equivariant} and hence in particular an $H$-equivariant weak equivalence.
\end{proof}
\end{lemma}

As in Corollary~\ref{cor:product-special-homotopical} we immediately deduce:

\begin{cor}\label{cor:coprod-prod-homotopical}
Finite products preserve $G$-global weak equivalences.\qed
\end{cor}

\begin{rk}\label{rk:products-homotopical-equivariant}
In the same way, one deduces from the equivariant comparison cited above that finite products preserve $H$-equivariant weak equivalences for any finite group $H$. Moreover, they preserve $H$-equivariant $\underline\pi_*$-isomorphisms for trivial reasons.
\end{rk}

If $X$ is any $G$-spectrum, then we can in particular apply Lemma~\ref{lemma:underived-additivity} to the family constant at $X$, proving that
\begin{equation}\label{eq:canonical-coprod-prod}
X^{\vee T}\to X^{\times T}
\end{equation}
is a $G$-global weak equivalence. As for $G$-global $\Gamma$-spaces, we want to strengthen this to a comparison taking the $\Sigma_T$-symmetries into account:

\begin{prop}\label{prop:equivariant-additivity}
The canonical map $(\ref{eq:canonical-coprod-prod})$ is a $(G\times\Sigma_T)$-global weak equivalence (in fact, even a $\underline\pi_*$-isomorphism) with respect to the $\Sigma_T$-actions permuting the summands or factors, respectively.
\end{prop}

While this is not the easiest way to prove the proposition, we will deduce it from a suitable version of the Wirthmüller isomorphism again. More precisely, recall from Construction~\ref{constr:wirthmueller-set} the natural Wirthmüller map $\gamma\colon\alpha_!Y\to\alpha_*Y$ for any injective homomorphism $\alpha\colon H\to G$ and any pointed $H$-set $Y$. Applying this levelwise in the simplicial and spectral directions, we then obtain a Wirthmüller map $\gamma\colon\alpha_!X\to\alpha_*X$ for any $H$-spectrum $X$. Below we will prove:

\begin{thm}\label{thm:wirthmueller}\index{Wirthmüller isomorphism!in G-Spectra@in $\cat{$\bm G$-Spectra}$!G-globally@$G$-globally|textbf}
Let $X$ be any $H$-spectrum, and let $\alpha\colon H\to G$ be an injective homomorphism with $(G:\im\alpha)<\infty$. Then the Wirthmüller map $\gamma\colon\alpha_!X\to\alpha_*X$ is a $G$-global weak equivalence (and in fact even a $\underline\pi_*$-isomorphism).
\end{thm}

The proof of the theorem needs some preparations, but before we come to this, let us observe that it immediately implies the proposition:

\begin{proof}[Proof of Proposition~\ref{prop:equivariant-additivity}]
This follows from Theorem~\ref{thm:wirthmueller} by the same argument as in the proof of Corollary~\ref{cor:boxtimes-special-power}.
\end{proof}

We will establish the Wirthmüller isomorphism by reduction to the equivariant statement; a similar argument appears in \cite[proof of Theorem~2.1.10]{proper-equivariant}. This will need the following classical point-set level manifestation of the Mackey double coset formula whose proof we leave to the interested reader, also see
\cite[discussion after Definition~1.5]{hausmann-equivariant} where this is stated without proof for finite $G$.

\begin{constr}
Let $K,H\subset G$ be any two subgroups and let $A$ be a pointed $H$-set. We define a $K$-equivariant map
\begin{equation*}
\chi\colon\bigvee_{[g]\in K\backslash G/H}K_+\smashp_{K\cap gHg^{-1}} c_g^*(A|_{g^{-1}Kg\cap H})\to (G_+\smashp_HA)|_K,
\end{equation*}
where the wedge runs over a fixed choice of double coset representatives, $(\blank)|_K$ denotes the underlying $K$-set, and $c_g$ is the conjugation homomorphism $K\cap gHg^{-1}\to g^{-1}Kg\cap H,k\mapsto g^{-1}kg$, as follows: on the summand corresponding to $g\in G$, the map is given by $[k,a]\mapsto [kg,a]$.

Analogously, we define
\begin{equation*}
\theta\colon \Maps^H(G,A)\big|_K\to
\prod_{[g]\in K\backslash G/H}\Maps^{K\cap gHg^{-1}}\big(K, c_g^*(A|_{g^{-1}Kg\cap H})\big)
\end{equation*}
as the map given on the factor corresponding to $g\in G$ by restricting along $K\to G, k\mapsto g^{-1}k$.

We omit the easy verification that these are indeed well-defined, $K$-equivariant, and moreover isomorphisms. Again, we get spectral versions (denoted by the same symbols) by applying these levelwise.
\end{constr}

\begin{constr}
Let $\phi\colon L\to G$ be a \emph{surjective} group homomorphism and let $H\subset G$ be any subgroup. Then we define for any pointed $H$-set $A$ a map
\begin{equation*}
\zeta\colon L_+\smashp_{\phi^{-1}H}(\phi|_{\phi^{-1}(H)}^*A)\to \phi^*(G_+\smashp_H A)
\end{equation*}
via $\zeta[\ell,a]=[\phi(\ell),a]$. Moreover, we write
\begin{equation*}
\nu\colon\phi^*\big(\Maps^H(G,A)\big)\to\Maps^{\phi^{-1}(H)}\big(L,(\phi|_{\phi^{-1}(H)})^*A\big)
\end{equation*}
for the restriction along $\phi$. We again omit the easy verification, that these are well-defined, $L$-equivariant, and natural isomorphisms. As before, we get natural isomorphisms of $L$-spectra, denoted by the same symbols.
\end{constr}

\begin{proof}[Proof of Theorem~\ref{thm:wirthmueller}]
We may again assume without loss of generality that $\alpha$ is the inclusion of a subgroup. Let $\phi\colon L\to G$ be any homomorphism from a finite group $L$; we have to show that $\phi^*(\gamma)$ is an $L$-equivariant weak equivalence.

For this we factor $\phi$ as the composition of the surjection $\bar\phi\colon L\to\im(\phi)\mathrel{=:}K$ followed by the inclusion $i\colon K\hookrightarrow G$. One then immediately checks that the diagram
\begin{equation*}
\begin{tikzcd}
\mathop{\smash{\bigvee\limits_{[g]}}}K_+\smashp_{K\cap gHg^{-1}} c_g^*(X|_{g^{-1}Kg\cap H}) \arrow[r, "\chi"]\arrow[dd, "\textup{canonical}"'] &[1em] (G_+\smashp_HX)|_K\arrow[d, "\gamma"]\\
& \Maps^H(G,X)|_K\arrow[d, "\theta"]\\
\prod\limits_{[g]}K_+\smashp_{K\cap gHg^{-1}} c_g^*(X|_{g^{-1}Kg\cap H})\arrow[r, "\prod\gamma\,"'] &\prod\limits_{[g]}\Maps^{K\cap gHg^{-1}}(K, c_g^*(X|_{g^{-1}Kg\cap H}))
\end{tikzcd}
\end{equation*}
in $\cat{$\bm K$-Spectra}$ commutes, where the products and coproducts run over representatives of $K\backslash G/H$ again. As seen above, $\chi$ and $\theta$ are $K$-equivariant isomorphisms; moreover, $K\backslash G/H$ is finite by assumption on $H$, and hence the left hand vertical arrow in the above diagram is a $K$-global weak equivalence (and in fact a $\underline\pi_*$-isomorphism) by Lemma~\ref{lemma:underived-additivity}. It is therefore enough to show that the lower horizontal map becomes an $L$-equivariant weak equivalences after restricting along $\bar\phi$, for which it is enough (Remark~\ref{rk:products-homotopical-equivariant}) to show this for each individual factor.

Up to renaming, we are therefore reduced to the case that $\phi\colon L\to G$ is surjective (and $G=K$ is finite). But in this case, the diagram
\begin{equation*}
\begin{tikzcd}
L_+\smashp_{\phi^{-1}(H)}(\phi|_{\phi^{-1}(H)})^*X\arrow[r,"\gamma"]\arrow[d, "\zeta"'] & \Maps^{\phi^{-1}(H)}\big(L, (\phi|_{\phi^{-1}(H)})^*X\big)\\
\phi^*(G_+\smashp_HX)\arrow[r, "\phi^*(\gamma)"'] & \phi^*\big(\Maps^H(G,X)\big)\arrow[u, "\nu"']
\end{tikzcd}
\end{equation*}
commutes by direct inspection. The vertical arrows were seen to be $L$-equivariant isomorphisms above, and the top horizontal arrow is an $L$-equivariant weak equivalence (and in fact a $\underline\pi_*$-isomorphism) by the usual Wirthmüller isomorphism, see e.g.~\cite[Proposition~3.7]{hausmann-equivariant} for a proof attributed to Schwede.\index{Wirthmüller isomorphism!in G-Spectra@in $\cat{$\bm G$-Spectra}$!G-equivariantly@$G$-equivariantly} The claim therefore follows by $2$-out-of-$3$.
\end{proof}

\subsection{The smash product}\label{subsec:smash-product}
The ordinary smash product of symmetric spectra\index{symmetric spectrum!smash product|seeonly{smash product of symmetric spectra}}\index{smash product of symmtric spectra!G-globally@$G$-globally|(} gives rise to a smash product on $\cat{$\bm G$-Spectra}$ by functoriality, and analogously for its right adjoint, the function spectrum construction $F$.\nomenclature[aF]{$F$}{function spectrum, right adjoint to smash product of symmetric spectra}\index{function spectrum!G-globally@$G$-globally} Explicitly, $F(X,Y)(A)=\Maps(X,\sh^AY)$, where here---and in the discussion below---we agree on the convention that $\Maps$ denotes the simplicial set of \emph{not necessarily $G$-equivariant maps}, equipped with the conjugation action.

We will now discuss some model categorical properties of these two functors:

\begin{prop}\label{prop:flatness}
\begin{enumerate}
\item If $X$ is any flat $G$-spectrum, then $X\smashp\blank$ preserves $G$-global weak equivalences.
\item If $X$ is any spectrum, then $X\smashp\blank$ preserves $G$-global weak equivalences between flat $G$-spectra.
\end{enumerate}
\begin{proof}
For the first statement, we let $f\colon Y\to Y'$ be any $G$-global weak equivalence. If $\phi\colon H\to G$ is any group homomorphism, then we have to show that $\phi^*(X\smashp f)$ is an $H$-equivariant weak equivalence. But this literally agrees with $\phi^*(X)\smashp\phi^*(f)$, and as $\phi^*(X)$ is flat and $\phi^*(f)$ is an $H$-weak equivalence, the claim follows from the usual equivariant Flatness Theorem \cite[Proposition~6.2-(i)]{hausmann-equivariant}.

The second statement follows similarly from \cite[Proposition~6.2-(ii)]{hausmann-equivariant}.
\end{proof}
\end{prop}

\begin{prop}\label{prop:ppo-smash-flat-G-global}
\begin{enumerate}
\item The smash product is a left Quillen bifunctor with respect to the $G$-global flat model structures everywhere.
\item The function spectrum is a right Quillen bifunctor with respect to the $G$-global flat model structures everywhere.
\end{enumerate}
\begin{proof}
It suffices to prove the first statement, for which we want to verify the Pushout Product Axiom. However, as in the previous proof one readily reduces this to the corresponding equivariant statement, which appears for example as \cite[Proposition~6.1]{hausmann-equivariant}.
\end{proof}
\end{prop}

Next, we want to establish the analogue of the above result with respect to the $G$-global projective model structures. In fact, we will prove a stronger `mixed' version of this, that (in the guises of Corollary~\ref{cor:function-spectrum-extra-homotopical} and Theorem~\ref{thm:F-I-vs-F}) will become crucial later in the proof of the $G$-global Delooping Theorem (Theorem~\ref{thm:group-completion}).

\begin{prop}\label{prop:smash-mixed}
\begin{enumerate}
\item The smash product is a left Quillen bifunctor
\begin{equation*}
\blank\smashp\blank\colon \cat{$\bm G$-Spectra}_{\textup{$G$-global flat}}\times\cat{$\bm G$-Spectra}_{\textup{$G$-global proj.}}\to\cat{$\bm G$-Spectra}_{\textup{$G$-global proj.}}.
\end{equation*}
\item The function spectrum is a right Quillen bifunctor
\begin{equation*}
F\colon \cat{$\bm G$-Spectra}_{\textup{$G$-global flat}}^\op\times\cat{$\bm G$-Spectra}_{\textup{$G$-global proj.}}\to\cat{$\bm G$-Spectra}_{\textup{$G$-global proj.}}.
\end{equation*}
\end{enumerate}
\end{prop}

The proof will rely on the following observation:

\begin{lemma}
Let $A$ be a finite faithful $G$-set (in particular $G$ is finite). Then
\begin{equation*}
\sh^A\colon\cat{$\bm G$-Spectra}_{\textup{$G$-global projective}}\to\cat{$\bm G$-Spectra}_{\textup{$G$-global flat}}
\end{equation*}
preserves acyclic fibrations.
\begin{proof}
Let $p\colon X\to Y$ be an acyclic fibration in the $G$-global projective model structure. We have to show that $(\sh^Ap)(B)^H$ is an acyclic Kan fibration for every finite set $B$ and any subgroup $H\subset G\times\Sigma_B$.

But by definition the $(G\times\Sigma_B)$-action on $(\sh^Ap)(B)=p(A\amalg B)$ is given by restricting the $G\times\Sigma_{A\amalg B}$-action on the right hand side along the homomorphism
\begin{equation*}
\begin{tikzcd}[cramped, sep=small, row sep=0pt]
\phi\colon G\times\Sigma_B\arrow[r]& G\times\Sigma_A\times\Sigma_B\arrow[r, hook] & G\times\Sigma_{A\amalg B}\\
\phantom{\phi\colon{}}(g,\sigma)\arrow[r, mapsto] & (g,\rho(g),\sigma)
\end{tikzcd}
\end{equation*}
(where $\rho\colon G\to\Sigma_A$ classifies the action on $A$), i.e.~$(\sh^Ap)(B)^H=p(A\amalg B)^{\phi(H)}$. It therefore suffices that $\phi(H)\in\mathcal G_{\Sigma_{A\amalg B},G}$. But indeed, if $\phi(h,\sigma)=(g,1)$, then $\sigma=1$ and $\rho(h)=1$ by definition of $\phi$. As $A$ is faithful, $\rho$ is injective, so that $h=1$ and hence also $g=1$. This finishes the proof.
\end{proof}
\end{lemma}

\begin{proof}[Proof of Proposition~\ref{prop:smash-mixed}]
Let $i\colon X\to X'$ be a flat cofibration and let $p\colon Y\to Y'$ be an acyclic fibration of the $G$-global projective model structure. We will first show that that the map
\begin{equation*}
(p_*,i^*)\colon F(X',Y)\to F(X',Y')\times_{F(X,Y')} F(X,Y)
\end{equation*}
is an acyclic $G$-global projective fibration, i.e.~if $\phi\colon H\to G$ is any group homomorphism and $A$ a finite faithful $H$-set, then $\big(\phi^*(p_*,i^*)(A)\big)^H$ is an acyclic Kan fibration. Indeed, this map agrees up to conjugation by isomorphisms with
\begin{align*}
&((\sh^A\phi^*(p))_*, \phi^*(i)^*)\colon\Maps^H(\phi^*X',\sh^A\phi^*Y)\\
&\qquad\to \Maps^H(\phi^*X',\sh^A\phi^*Y')\times_{\Maps^H(\phi^*X,\sh^A\phi^*Y')}\Maps^H(\phi^*X,\sh^A\phi^*Y).
\end{align*}
By the previous lemma together with Lemma~\ref{lemma:restriction-right-Quillen-projective}, $\sh^A\phi^*(p)$ is an acyclic fibration in the $H$-global flat model structure, and moreover $\phi^*(i)$ is a flat cofibration by definition. The claim follows because the $H$-global flat model structure is simplicial.

By adjointness we may therefore conclude that the pushout product
\begin{equation*}
j\ppo k\colon (A\smashp B')\amalg_{(A\smashp B)}(A'\smashp B)\to A'\smashp B'
\end{equation*}
of any $G$-global flat cofibration $j\colon A\to A'$ with any $G$-global projective cofibration $k\colon B\to B'$ is a $G$-global projective cofibration. On the other hand, if at least one of $j$ or $k$ is acyclic, then Proposition~\ref{prop:ppo-smash-flat-G-global} shows that $j\ppo k$ is a $G$-global acyclic flat cofibration (as any $G$-global projective cofibration is also a flat cofibration), hence in particular a $G$-global weak equivalence. This proves the first part of the proposition, and the second one then follows by the usual adjointness argument.
\end{proof}

If $X$ is an $H$-spectrum and $Y$ is a $G$-spectrum, then $F(X,Y)$ carries an $(H\times G)$-action, yielding a functor $F\colon\cat{$\bm H$-Spectra}^\op\times\cat{$\bm G$-Spectra}\to\cat{$\bm{(H\times G)}$-Spectra}$. As an application of the above proposition we can now prove:

\begin{cor}\label{cor:function-spectrum-extra-homotopical}
\begin{enumerate}
\item If $X$ is any flat $H$-spectrum, then
\begin{equation}\label{eq:F-extra-action}
F(X,\blank)\colon\cat{$\bm G$-Spectra}_{\textup{$G$-gl.~proj.}}\to\cat{$\bm{(H\times G)}$-Spectra}_{\textup{$(H\times G)$-gl.~proj.}}
\end{equation}
is right Quillen.
\item If $Y$ is fibrant in the $G$-global projective model structure, then
\begin{equation*}
F(\blank,Y)\colon\cat{$\bm H$-Spectra}_{\textup{$H$-global flat}}^\op\to\cat{$\bm{(H\times G)}$-Spectra}_{\textup{$(H\times G)$-global projective}}
\end{equation*}
is right Quillen.
\end{enumerate}
\begin{proof}
We will prove the first statement, the proof of the second being similar.

For this we let $\pr_1\colon H\times G\to H$ and $\pr_2\colon H\times G\to G$ denote the projections. Then $(\ref{eq:F-extra-action})$ agrees with the composition
\begin{align*}
\cat{$\bm G$-Spectra}_{\textup{$G$-global proj.}}&\xrightarrow{\pr_2^*}
\cat{$\bm{(H\times G)}$-Spectra}_{\textup{$(H\times G)$-global proj.}}\\
&\xrightarrow{F(\pr_1^*X,\blank)}\cat{$\bm{(H\times G)}$-Spectra}_{\textup{$(H\times G)$-global proj.}}.
\end{align*}
By Lemma~\ref{lemma:restriction-right-Quillen-projective} the first functor is right Quillen, and so is the second one by Proposition~\ref{prop:smash-mixed} as $\pr_1^*X$ is flat.
\end{proof}
\end{cor}
\index{smash product of symmtric spectra!G-globally@$G$-globally|)}

\begin{rk}
By another adjointness argument, Proposition~\ref{prop:smash-mixed} implies that $F$ is also a right Quillen bifunctor with respect to the $G$-global projective model structures on the source and the $G$-global flat model structure in the target. While this may sound somewhat odd, we emphasize again that the $G$-global projective model structure has `few' cofibrant objects for $G\not=1$. In particular, $\mathbb S$ is not cofibrant, so this does not apply to $F(\mathbb S,\blank)\cong\id$.
\end{rk}

\section[Connections to unstable $G$-global homotopy theory]{Connections to unstable \for{toc}{$G$}\except{toc}{\texorpdfstring{$\bm G$}{G}}-global homotopy theory}
We can now connect the above stable models to the unstable ones discussed in Chapter~\ref{chapter:unstable}, and in particular to $G$-$I$-simplicial sets and $G$-$\mathcal I$-simplicial sets.

\index{symmetric spectrum!suspension spectrum|seeonly{suspension spectrum}}
\subsection{Suspension spectra}\index{suspension spectrum!G-globally@$G$-globally|(}
We begin by recalling the \emph{suspension/loop adjunction} between $\cat{$\bm I$-SSet}_*$ and $\cat{Spectra}$, see e.g.~\cite[discussion before Proposition~3.19]{sagave-schlichtkrull}:

\begin{constr}
If $X$ is any $I$-simplicial set, then we define a symmetric spectrum $\Sigma^\bullet X$ via $(\Sigma^\bullet X)(A)=S^A\smashp X(A)$ together with the evident structure maps; $\Sigma^\bullet$\nomenclature[aSigmab]{$\Sigma^\bullet$}{suspension spectrum of an $I$- or $\mathcal I$-simplicial set} becomes a simplicially enriched functor in the obvious way. On the other hand, if $Y$ is a symmetric spectrum, then we define $\Omega^\bullet Y$ via $(\Omega^\bullet Y)(A)=\Omega^AY(A)$\nomenclature[aOmegab]{$\Omega^\bullet$}{loop $I$-simplicial set of a global spectrum} together with the evident structure maps. Again, this becomes a functor in the obvious way, and we omit the easy verification that $\Sigma^\bullet$ is a simplicial left adjoint to $\Omega^\bullet$. By functoriality, this gives rise to a simplicial adjunction $\cat{$\bm G$-$\bm I$-SSet}_*\rightleftarrows\cat{$\bm G$-Spectra}$ that we denote by $\Sigma^\bullet\dashv\Omega^\bullet$ again.
\end{constr}

\begin{prop}\label{prop:suspension-loop-projective}
The simplicial adjunction
\begin{equation}\label{eq:suspension-loop-projective}
\Sigma^\bullet\colon(\cat{$\bm G$-$\bm I$-SSet}_*)_{\textup{$G$-global}}\rightleftarrows\cat{$\bm G$-Spectra}_{\textup{$G$-global projective}} :\Omega^\bullet
\end{equation}
is a Quillen adjunction, and $\Sigma^\bullet$ is fully homotopical.
\begin{proof}
One immediately checks from the definitions that
\begin{equation*}
S^A\smashp\blank\colon\cat{$\bm{(G\times\Sigma_A)}$-SSet}_*\rightleftarrows\cat{$\bm{(G\times\Sigma_A)}$-SSet}_* :\!\Omega^A
\end{equation*}
(where $\Sigma_A$ acts on $S^A$ via its tautological $A$-action) is a Quillen adjunction with respect to the $\mathcal G_{\Sigma_A,G}$-equivariant model structure on both sides, proving that $(\ref{eq:suspension-loop-projective})$ is a Quillen adjunction with respect to the strict level model structure on the source and the $G$-global projective level model structure on the target, hence also with respect to the $G$-global projective model structure on the target. To prove that also $(\ref{eq:suspension-loop-projective})$ itself is a Quillen adjunction it is then enough that $\Omega^\bullet$ sends fibrant objects to static $G$-$I$-simplicial sets, which is immediate from the definitions.

To show that $\Sigma^\bullet$ is homotopical, it is now enough to observe that any weak equivalence in $\cat{$\bm G$-$\bm I$-SSet}_*$ factors as an acyclic cofibration followed by a strict level weak equivalence. The former are preserved according to the above, and clearly strict level weak equivalences are even sent to $G$-global level weak equivalences.
\end{proof}
\end{prop}

\begin{cor}
The simplicial adjunction $(\ref{eq:suspension-loop-projective})$ is also a Quillen adjunction with respect to the $G$-global \emph{injective} model structures.
\begin{proof}
It is obvious that $\Sigma^\bullet$ preserves injective cofibrations, and by the previous proposition it also preserves $G$-global weak equivalences.
\end{proof}
\end{cor}

\begin{rk}\label{rk:omega-bullet-I-homotopy-groups}
The above adjunction is compatible with (non-negatively indexed) homotopy groups; we will make one special case of this precise that will become relevant later:

\goodbreak
Assume $G$ is finite, let $X$ be any $G$-spectrum, and let $\mathcal U=\mathcal U_G$ be our fixed complete $G$-set universe. Then we have natural maps
\begin{align*}
\pi_0^G\big((\Omega^\bullet X)(\mathcal U)\big)&\cong
\pi_0^G\big(\colim_{A\in s(\mathcal U)}\Omega^AX(A)\big)\\
&\cong\colim_{A\in s(\mathcal U)}\pi_0^G\Omega^AX(A)=\colim_{A\in s(\mathcal U)}\pi_0\Maps^G(S^A,X(A))\\
&\to\colim_{A\in s(\mathcal U)}[S^A,X(A)]^G_*=\pi_0^G(X)
\end{align*}
(where the first isomorphism comes from cofinality, the second one uses that $\pi_0^G$ preserves filtered colimits, and the final map is induced by the localization functor), yielding a natural transformation $\pi_0^G\circ\Omega^\bullet\Rightarrow \pi_0^G$. If $X$ is fibrant in the $G$-global projective level model structure, then $X(A)$ is fibrant in the $G$-equivariant model structure for every finite faithful $G$-set $A$, so that the final map is an isomorphism as a colimit of eventual isomorphisms, and hence so is the above composition.
\end{rk}

\begin{ex}\label{ex:suspension-IAphi}
Let $H$ be a finite group, let $A$ be a finite faithful $H$-set, and let $\phi\colon H\to G$ be a homomorphism. Generalizing \cite[Example~6.3]{hausmann-global}, $\Sigma_+^\bullet I(A,\blank)\times_\phi G$ corepresents the \emph{true} zeroth $\phi$-equivariant homotopy group in the following sense: if $X$ is any $G$-global $\Omega$-spectrum, then a similar adjointness argument to the above shows that evaluating in degree $A$ and then restricting along
\begin{equation*}
\blank\smashp[\id_A,1]\colon S^A\to S^A\smashp (I(A,\blank)\times_\phi G)_+= \Sigma_+^\bullet(I(A,\blank)\times_\phi G)(A)
\end{equation*}
induces a bijection $[\Sigma_+^\bullet(I(A,\blank)\times_\phi G),X]\to\pi_0^\phi(X)$, where $[\,{,}\,]$ denotes the hom set in the $G$-global stable homotopy category. In general, $[\Sigma^\bullet_+(I(A,\blank)\times_\phi G),Y]$ computes $\pi_0^\phi$ of a replacement of $Y$ by a $G$-global $\Omega$-spectrum.
\end{ex}

We now turn to a variant for $G$-$\mathcal I$-simplicial sets:

\begin{cor}
The simplicial adjunction
\begin{equation*}\label{suspension-loop-script-I}
\Sigma^\bullet\mathrel{:=}\Sigma^\bullet\circ\forget\colon(\cat{$\bm G$-$\bm{\mathcal I}$-SSet}_*)_{\textup{$G$-global}}\rightleftarrows\cat{$\bm G$-Spectra}_{\textup{$G$-global proj.}} :\!\Maps_I(\mathcal I,\blank)\circ\Omega^\bullet
\end{equation*}
is a Quillen adjunction, and likewise for the injective $G$-global model structures.
\begin{proof}
This is immediate from the above as the forgetful functor $\cat{$\bm G$-$\bm{\mathcal I}$-SSet}\to\cat{$\bm G$-$\bm{I}$-SSet}$ is left Quillen for the $G$-global model structures as well as for the injective $G$-global ones (Theorems~\ref{thm:forget-left-quillen-global} and~\ref{thm:I-vs-M-injective}, respectively).
\end{proof}
\end{cor}

\begin{rk}
Composing with the adjunction
\begin{equation*}
(\blank)_+\colon\cat{$\bm G$-$\bm{I}$-SSet}\rightleftarrows\cat{$\bm G$-$\bm{I}$-SSet}_* :\!\forget
\end{equation*}
(or its $\mathcal I$-version) we also get unpointed versions of all of the above adjunctions. We denote the right adjoints as before and write $\Sigma^\bullet_+$ for the left adjoints.
\end{rk}

\begin{rk}\label{rk:omega-bullet-script-I-homotopy-groups}
By Theorem~\ref{thm:forget-left-quillen-global}, the counit $\forget\Maps_I(\mathcal I,Y)\to Y$ is a $G$-global weak equivalence (hence strict level weak equivalence) for any $G$-$I$-simplicial set $Y$ that is fibrant in the $G$-global model structure. By Proposition~\ref{prop:suspension-loop-projective}, this in particular applies to $Y=\Omega^\bullet X$ for any $G$-spectrum $X$ fibrant in the $G$-global projective model structure. Composing with the transformation from Remark~\ref{rk:omega-bullet-I-homotopy-groups} we therefore get a natural transformation
\begin{equation*}
\pi_0^G\big(\Maps_I(\mathcal I,\Omega^\bullet(\blank))(\mathcal U)\big)\Rightarrow\pi_0^G,
\end{equation*}
which is an isomorphism for all fibrant objects of $\cat{$\bm{G}$-Spectra}_{\textup{$G$-global projective}}$.
\end{rk}
\index{suspension spectrum!G-globally@$G$-globally|)}

\index{symmetric spectrum!tensoring of G-Spectra over G-I-SSet@tensoring of $\cat{$\bm G$-Spectra}$ over $\cat{$\bm G$-$\bm I$-SSet}$|seeonly{tensoring of $\cat{$\bm G$-Spectra}$ over $\cat{$\bm G$-$\bm I$-SSet}$}}
\index{symmetric spectrum!tensoring of G-Spectra over G-II-SSet@tensoring of $\cat{$\bm G$-Spectra}$ over $\cat{$\bm G$-$\bm{\mathcal I}$-SSet}$|seeonly{tensoring of $\cat{$\bm G$-Spectra}$ over $\cat{$\bm G$-$\bm{\mathcal I}$-SSet}$}}
\subsection{Enrichment and tensoring over \texorpdfstring{$\bm G$-$\bm{\mathcal I}$}{G-I}-simplicial sets}
We will now exhibit another connection between ($G$-)spectra and ($G$-)$I$-simplicial sets. For this we first recall the usual tensoring of $\cat{$\bm G$-Spectra}$ over $\cat{$\bm G$-$\bm I$-SSet}_*$, see e.g.~\cite[Example~2.31]{schwede-symmetric-book} where this is denoted `$\smashp$':

\begin{constr}
Let $X$ be a spectrum and let $Y$ be a pointed $I$-simplicial set. We write $X\otimes Y$ for the spectrum with $(X\otimes Y)(A)=X(A)\smashp Y(A)$ and structure maps
\begin{equation*}
S^{B\setminus i(A)}\smashp\big(X(A)\smashp Y(A)\big)\cong
\big(S^{B\setminus i(A)}\smashp X(A)\big)\smashp Y(A)
\xrightarrow{\sigma\smashp Y(i)}
X(B)\smashp Y(B)
\end{equation*}
\nomenclature[zo2]{$\otimes$}{tensoring of $\cat{Spectra}$ over $\cat{$\bm I$-SSet}$ via levelwise smash product}\index{tensoring of G-Spectra over G-I-SSet@tensoring of $\cat{$\bm G$-Spectra}$ over $\cat{$\bm G$-$\bm I$-SSet}$|textbf}%
for any injection $i\colon A\to B$ of finite sets, where the unlabelled isomorphism is the associativity constraint. The tensor product becomes a simplicially enriched functor in both variables by applying the enriched functoriality of the smash product of pointed simplicial sets levelwise. If $G$ is any discrete group, then taking the diagonal of the two actions promotes the tensor product to a simplicially enriched functor $\cat{$\bm G$-Spectra}\times\cat{$\bm G$-$\bm I$-SSet}_*\to\cat{$\bm G$-Spectra}$ that we denote by the same symbol. Finally, precomposing with the forgetful functor we get a simplicially enriched functor
\begin{equation}\label{eq:otimes-script-I}
\cat{$\bm G$-Spectra}\times\cat{$\bm G$-$\bm{\mathcal I}$-SSet}_*\to\cat{$\bm G$-Spectra}
\end{equation}
that we again denote by $\blank\otimes\blank$.\index{tensoring of G-Spectra over G-II-SSet@tensoring of $\cat{$\bm G$-Spectra}$ over $\cat{$\bm G$-$\bm{\mathcal I}$-SSet}$|textbf}
\end{constr}

\begin{rk}
It is not hard to verify that these indeed become actions of $\cat{$\bm G$-$\bm I$-SSet}_*$ and $\cat{$\bm G$-$\bm{\mathcal I}$-SSet}_*$ (with the symmetric monoidal structure given by the levelwise smash product), respectively, on $\cat{$\bm G$-Spectra}$ in a preferred way; we leave the details to the curious reader.
\end{rk}

\subsubsection{Homotopical properties of the tensor product}\index{tensoring of G-Spectra over G-II-SSet@tensoring of $\cat{$\bm G$-Spectra}$ over $\cat{$\bm G$-$\bm{\mathcal I}$-SSet}$!homotopical behaviour|(}
It will be natural to first study these equivariantly, so let us fix a finite group $H$. We begin with a comparison of homotopy groups:

\begin{lemma}\label{lemma:tensor-vs-smash-homotopy-groups}
Let $X\in\cat{$\bm H$-Spectra}$ and $Y\in\cat{$\bm H$-$\bm{I}$-SSet}_*$. Then there exists a natural isomorphism
\begin{equation*}
\underline\pi_*^{H}(X\otimes Y)\cong\underline\pi_*^{H}(X\smashp Y(\mathcal U_H)).
\end{equation*}
\begin{proof}
This is a standard argument, see e.g.~\cite[Proposition~5.14]{schwede-equivariant}. We will only prove the claim for $\pi_0^{H}$; the general case is done analogously, but requires more notation. For this we consider
\begin{align*}
\pi_0^{H}(X\otimes Y)&=\colim_{A\in s(\mathcal U_H)}[S^A, (X\otimes Y)(A)]^H_*\\
&=\colim_{A\in s(\mathcal U_H)}[S^A, X(A)\smashp Y(A)]^H_*\\
&\cong\colim_{(A,B)\in s(\mathcal U_H)^2}[S^A, X(A)\smashp Y(B)]^H_*\\
&\cong\colim_{A\in s(\mathcal U_H)}\colim_{B\in s(\mathcal U_H)}[S^A, X(A)\smashp Y(B)]^H_*\\
&\cong\colim_{A\in s(\mathcal U_H)}[S^A, X(A)\smashp\colim_{B\in s(\mathcal U_H)} Y(B)]^H_*\\
&\cong\colim_{A\in s(\mathcal U_H)}[S^A, X(A)\smashp Y(\mathcal U_H)]^H_*\\
&=\pi_0^{H}(X\smashp Y(\mathcal U_H)).
\end{align*}
Here the first isomorphism uses that the diagonal $s(\mathcal U_H)\to s(\mathcal U_H)^2$ is cofinal, the second one is the Fubini Theorem for colimits, and the third one uses that $S^A$ is compact (in the derived sense) and that $\smashp$ is cocontinuous in each variable. As all of the above isomorphisms are clearly natural, this finishes the proof.
\end{proof}
\end{lemma}

For the rest of this section we focus on $(\ref{eq:otimes-script-I})$ as it is more tractable:

\begin{prop}\label{prop:tensor-vs-smash}
Let $X\in\cat{$\bm H$-Spectra}$, $Y\in\cat{$\bm H$-$\bm{\mathcal I}$-SSet}_*$. Then the maps
\begin{equation*}
X\otimes Y\to X\otimes Y(\mathcal U_H\amalg\blank)\gets X\otimes\const Y(\mathcal U_H)=X\smashp Y(\mathcal U_H)
\end{equation*}
induced by the inclusions $A\hookrightarrow \mathcal U_H\amalg A\hookleftarrow\mathcal U_H$ are equivariant $\underline\pi_*$-isomorphisms.
\begin{proof}
Replacing $H$ by a subgroup $K\subset H$ if necessary and using that $\mathcal U_H$ is also a complete $K$-set universe, it suffices to show this for $\pi_*^H$. By the previous lemma it is then enough to prove that the induced maps
\begin{equation*}
\pi_*^H(X\smashp Y(\mathcal U_H))\to\pi_*^H(X\smashp Y(\mathcal U_H\amalg\blank)(\mathcal U_H))\gets\pi_*^H(X\smashp(\const Y(\mathcal U_H))(\mathcal U_H))
\end{equation*}
are isomorphisms. For this we simply observe that up to the natural isomorphisms $Y(\mathcal U_H\amalg\blank)(\mathcal U_H)\cong Y(\mathcal U_H\amalg\mathcal U_H)$ (see Lemma~\ref{lemma:extension-filtered-colimit}) and $(\const Y(\mathcal U_H))(\mathcal U_H)\cong Y(\mathcal U_H)$, these are induced from the two inclusions $\mathcal U_H\rightrightarrows\mathcal U_H\amalg\mathcal U_H$, so that the claim follows from Lemma~\ref{lemma:evaluation-h-universe}.
\end{proof}
\end{prop}

\begin{constr}
Let $X$ be an $H$-spectrum and let $Y$ be any $H$-$\mathcal I$-simplicial set. Then we define $\psi\colon X\smashp \Sigma^\bullet Y\to X\otimes Y$ as the map associated to the bimorphism given in degree $A,B$ by 
\begin{align*}
X(A)\smashp (\Sigma^\bullet Y)(B)&=X(A)\smashp S^B\smashp Y(B)
\xrightarrow{\textup{twist}}
S^B\smashp X(A)\smashp Y(B)\\
&\xrightarrow{\sigma_{B,A}\smashp Y(B\hookrightarrow A\amalg B)}
X(B\amalg A)\smashp Y(A\amalg B)\\
&\xrightarrow{X(\textup{twist})\smashp\id}
X(A\amalg B)\smashp Y(A\amalg B)
=(X\otimes Y)(A\amalg B).
\end{align*}
\end{constr}

\begin{prop}\label{prop:smash-vs-tensor-equivariant}
The above map $\psi$ is natural in both variables. Moreover, if $X$ is flat, then $\psi$ is a $\underline\pi_*$-isomorphism, hence in particular an $H$-equivariant weak equivalence.
\begin{proof}
The naturality is obvious. For the second statement, we observe that
\begin{equation*}
\begin{tikzcd}
X\smashp\Sigma^\bullet Y\arrow[r, "\psi"]\arrow[d] & X\otimes Y\arrow[d]\\
X\smashp\Sigma^\bullet Y(\mathcal U_H\amalg\blank)\arrow[r, "\psi"] & X\otimes Y(\mathcal U_H\amalg\blank)\\
X\smashp\Sigma^\bullet\const Y(\mathcal U_H)\arrow[u]\arrow[r, "\psi"'] & X\otimes\const Y(\mathcal U_H)\arrow[u]
\end{tikzcd}
\end{equation*}
commutes by naturality, where the vertical arrows are induced from the zig-zag in Proposition~\ref{prop:tensor-vs-smash} (using that $\Sigma^\bullet=\mathbb S\otimes\blank$). In particular, the proposition tells us that the right hand vertical maps are $\underline\pi_*$-isomorphisms and as $X\smashp\blank$ preserves $\underline\pi_*$-isomorphisms by \cite[Proposition~6.2-(i)]{hausmann-equivariant}, so are the left hand vertical maps.

But the lower horizontal map in the above diagram literally agrees with the canonical comparison map $X\smashp\Sigma^\infty Y(\mathcal U_H)\to X\smashp Y(\mathcal U_H)$, so it is even an isomorphism, see e.g.~\cite[Proposition~3.5]{schwede-symmetric-book}. The claim follows by $2$-out-of-$3$.
\end{proof}
\end{prop}

We can now use this to establish $G$-global properties of the tensor product:

\begin{prop}\label{prop:E-smash-vs-E-tensor-base-case}
Let $X$ be a flat $G$-spectrum, and let $Y$ be any $G$-$\mathcal I$-simplicial set. Then the above map $\psi\colon X\smashp\Sigma^\bullet Y\to X\otimes Y$ is a $G$-global weak equivalence (in fact, even a $\underline\pi_*$-isomorphism).

Moreover, if $Y$ is cofibrant in the $G$-global model structure, then $\alpha_!(\psi)$ is a $G'$-global weak equivalence for any homomorphism $\alpha\colon G\to G'$.
\begin{proof}
For the first statement we let $\phi\colon H\to G$ be any homomorphism. Then $\phi^*(\psi)$ literally agrees with $\psi\colon (\phi^*X)\smashp\Sigma^\bullet(\phi^*Y)\to (\phi^*X)\otimes(\phi^*Y)$. As flatness is a property of the underlying non-equivariant spectrum, $\phi^*X$ is again flat, so the claim follows from the previous proposition.

For the second statement we observe that $X\smashp\Sigma^\bullet Y$ is cofibrant in the $G$-global projective model structure by Proposition~\ref{prop:suspension-loop-projective} together with Proposition~\ref{prop:smash-mixed}, so $G$ acts levelwise freely outside the base point on it by Remark~\ref{rk:projective-free}. Moreover, $G$ also acts levelwise freely outside the base point on $Y$, hence also on $X\otimes Y$. The claim therefore follows from the first statement together with Proposition~\ref{prop:spectra-free-quotient}.\index{tensoring of G-Spectra over G-II-SSet@tensoring of $\cat{$\bm G$-Spectra}$ over $\cat{$\bm G$-$\bm{\mathcal I}$-SSet}$!homotopical behaviour|)}
\end{proof}
\end{prop}

\subsubsection{$G$-global mapping spaces} We now introduce a certain `$G$-global mapping space' adjoint to the above tensor product:
\index{symmetric spectrum!G-global mapping space@$G$-global mapping space|seeonly{$G$-global mappping space}}

\begin{constr}\label{constr:tensor-right-adjoint}\index{G-global mapping space@$G$-global mapping space|textbf}
We define a simplicially enriched functor
\begin{equation*}
F_{\mathcal I}\colon\cat{$\bm G$-Spectra}^\op\times\cat{$\bm G$-Spectra}\to\cat{$\bm G$-$\bm{\mathcal I}$-SSet}_*
\end{equation*}
as follows: if $X$ and $Y$ are $G$-spectra and $A$ is a finite set, then
\begin{equation*}
F_{\mathcal I}(X,Y)(A)=\Maps(X\otimes\mathcal I(A,\blank)_+, Y)
\end{equation*}
\nomenclature[aFII]{$F_{\mathcal I}$}{right adjoint to tensoring of $\cat{Spectra}$ over $\cat{$\bm{\mathcal I}$-SSet}$}%
(we remind the reader that $\Maps$ consists of all not necessarily $G$-equivariant maps and that it is a $G$-simplicial set by conjugation, with based point coming from the constant map). The covariant functoriality in $A$ is given by contravariant functoriality of $\mathcal I(A,\blank)$; the contravariant functoriality in $X$ and the covariant functoriality in $Y$ are the obvious ones.

We define a map $\epsilon_{X,Y}\colon X\otimes F_{\mathcal I}(X,Y)\to Y$ as follows: in degree $A$, $\epsilon_{X,Y}$ is given by the composition
\begin{align*}
X(A)\hskip-1pt\smashp\hskip-1pt\Maps(X\hskip-1pt\otimes\mathcal I(A,\blank)_+,Y)\hskip-1pt&\to\hskip-1pt X(A)\smashp\Maps(X(A)\smashp\mathcal I(A,A)_+, Y(A))\\
&\to\hskip-1pt X\hskip-.5pt(A)\hskip-1pt\smashp \hskip-1pt\mathcal I(A,A)_+\hskip-1pt\smashp\hskip-1pt\Maps(X(A)\hskip-1pt\smashp\hskip-1pt\mathcal I(A,A)_+,\hskip-1pt Y\hskip-1pt(A))\\
&\to\hskip-1pt Y(A)
\end{align*}
where the first map is induced by evaluating in degree $A$, the second one sends $x\in X_n$ to $x\smashp (\id_A,\dots,\id_A)$, and the final one is the counit of the usual smash-hom adjunction in $\cat{SSet}_*$. We omit the easy verification that $\epsilon_{X,Y}$ is a map of $G$-spectra and that it is natural (in the enriched sense) in both variables.

We further define $\eta_{X,Y}\colon Y\to F_{\mathcal I}(X,X\otimes Y)$ as follows: in degree $A$, $\eta_{X,Y}$ is the map $Y(A)\to\Maps(X\otimes\mathcal I(A,\blank)_+,X\otimes Y)$ whose postcomposition with evaluation at $B$ is the map $Y(A)\to\Maps(X(B)\smashp\mathcal I(A,B)_+,X(B)\smashp Y(B))$ adjoint to the map
\begin{equation*}
X(B)\smashp \mathcal I(A,B)_+\smashp Y(A)\to X(B)\smashp Y(B)
\end{equation*}
induced by the enriched functoriality of $Y$ (here we have suppressed the associativity isomorphism for simplicity). We again omit the straight-forward verification that this is well-defined and a simplicially enriched natural transformation. We moreover omit that for any fixed $G$-spectrum $X$ the natural transformation $\epsilon_{X,\blank}$ and $\eta_{X,\blank}$ satisfy the triangle equations, so that they form unit and counit of a simplicially enriched adjunction
\begin{equation*}
X\otimes\blank\colon\cat{$\bm G$-$\bm{\mathcal I}$-SSet}_*\rightleftarrows\cat{$\bm G$-Spectra} :\!F_{\mathcal I}(X,\blank).
\end{equation*}
\end{constr}

Let $X,Y$ be $G$-spectra. Then there is also another way to cook up a `mapping $G$-$\mathcal I$-space' between $X$ and $Y$, namely $\Maps_I(\mathcal I,\Omega^\bullet F(X,Y))$. Our goal for the rest of this section is to compare these two; more generally, we will need a version where we allow $X$ to be an $H$-spectrum and we want a $(G\times H)$-global comparison.

\begin{constr}\label{constr:psi-hat}
If $T$ is any spectrum, then we have by adjointness and the simplicially enriched Yoneda Lemma a sequence of isomorphisms
\begin{align*}
\Maps_I(\mathcal I,\Omega^\bullet T)(A)&\cong\Maps_{\cat{$\bm{\mathcal I}$-SSet}}(\mathcal I(A,\blank), \Maps_I(\mathcal I, \Omega^\bullet T))\\
&\cong \Maps_\cat{$\bm{I}$-SSet}(\mathcal I(A,\blank),\Omega^\bullet T)\cong \Maps_{\cat{Spectra}}(\Sigma^\bullet_+\mathcal I(A,\blank), T)
\end{align*}
natural in $T$ and the finite set $A$. Applying this to $T=F(X,Y)$ for $X$ and $Y$ spectra and appealing to the adjunction between $\smashp$ and $F$ we get an isomorphism
\begin{equation}\label{eq:maps-I-bullet-simplify}
\Maps_I\big(\mathcal I,\Omega^\bullet F(X,Y)\big)(A) \cong\Maps(X\smashp\Sigma^\bullet_+\mathcal I(A,\blank), Y)
\end{equation}
natural in $X$, $Y$, and $A$. In particular, if the group $H$ acts on $X$ and the group $G$ acts on $Y$, then these maps for varying $A$ assemble into an isomorphism in $\cat{$\bm{(G\times H)}$-$\bm{\mathcal I}$-SSet}_*$.

We now define $\widehat\psi\colon F_{\mathcal I}(X,Y)\to\Maps_I(\mathcal I,\Omega^\bullet F(X,Y))$ as the map given in degree $A$ by the composition
\begin{equation*}
\Maps(X\otimes\mathcal I(A,\blank)_+, Y)\xrightarrow{\psi^*}
\Maps(X\smashp\Sigma^\bullet_+\mathcal I(A,\blank),Y)\cong \Maps_I\big(\mathcal I,\Omega^\bullet F(X,Y)\big)(A)
\end{equation*}
where the unlabelled isomorphism is $(\ref{eq:maps-I-bullet-simplify})$. We omit the easy verification that this is well-defined, and natural in both variables.
\end{constr}

Now we are ready to state and prove the main result of this section:

\begin{thm}\label{thm:F-I-vs-F}\index{G-global mapping space@$G$-global mapping space!comparison to function spectrum|textbf}
\index{function spectrum!G-globally@$G$-globally!comparison to G-global mapping space@comparison to $G$-global mapping space|textbf}
If $X$ is a flat $H$-spectrum and $Y$ is fibrant in the $G$-global injective model structure, then
\begin{equation}\label{eq:F-vs-F-I}
\widehat\psi\colon F_{\mathcal I}(X,Y)\to \Maps_I\big(\mathcal I,\Omega^\bullet F(X,Y)\big)
\end{equation}
is a $(G\times H)$-global weak equivalence.
\begin{proof}
We show that it is a $(G\times H)$-global level weak equivalence. For this we let $K$ be any finite group acting faithfully on $A$ and $\phi\colon K\to G\times H$ any group homomorphism. We have to show that $\phi^*(\widehat\psi)(A)^K$ is a weak homotopy equivalence, for which it is enough by construction of $\widehat\psi$ that $\phi^*(\psi^*)(A)^K$ is. But if we write $\phi_1\colon K\to G$ and $\phi_2\colon K\to H$ for the components of $\phi$, then this agrees with
\begin{equation*}
\Maps^K(\psi,\phi_1^*Y)\colon\Maps^K(\phi_2^*X\otimes\mathcal I(A,\blank)_+, \phi_1^*Y)\to\Maps^K(\phi_2^*X\smashp\Sigma^\bullet_+\mathcal I(A,\blank), \phi_1^*Y).
\end{equation*}
Using that $\phi_{1!}$ is a simplicial left adjoint to $\phi_1^*$ we see that this agrees up to conjugation by isomorphisms with
\begin{equation*}
\big(\phi_{1!}(\psi)\big)^*\colon\Maps^G(\phi_{1!}(\phi_2^*X\otimes\mathcal I(A,\blank)_+), Y)\to
\Maps^G(\phi_{1!}(\phi_2^*X\smashp\Sigma^\bullet_+\mathcal I(A,\blank)), Y).
\end{equation*}
But $\phi_{1!}(\psi)$ is a $G$-global weak equivalence by Proposition~\ref{prop:E-smash-vs-E-tensor-base-case}, and $Y$ is fibrant in the injective $G$-global model structure by assumption. The claim follows since the latter is a simplicial model structure.
\end{proof}
\end{thm}

\section[$G$-global spectra vs.~$G$-equivariant spectra]{\texorpdfstring{\except{toc}{$\bm G$}\for{toc}{$G$}}{G}-global spectra vs.~\texorpdfstring{\except{toc}{$\bm G$}\for{toc}{$G$}}{G}-equivariant spectra}
\index{G-equivariant stable homotopy theory@$G$-equivariant stable homotopy theory!vs G-global stable homotopy theory@vs.~$G$-global stable homotopy theory|(}
In this section we will finally prove as promised:

\begin{prop}\label{prop:flat-fibrant-G-Omega}
Let $G$ be a finite group. Then the simplicial adjunctions
\begin{align}\label{eq:flat-global-vs-equivariant}
\id\colon\cat{$\bm G$-Spectra}_{\textup{$G$-equivariant projective}}&\rightleftarrows\cat{$\bm G$-Spectra}_{\textup{$G$-global flat}} :\!\id\\
\label{eq:proj-global-vs-equivariant}
\Sigma^\bullet_+I(A,\blank)\smashp\blank\colon\cat{$\bm G$-Spectra}_{\textup{$G$-equiv.~proj.}}&\rightleftarrows\cat{$\bm G$-Spectra}_{\textup{$G$-global proj.}} :\!\Omega^A\sh^A
\end{align}
are Quillen adjunctions, where $A$ is any finite faithful $G$-set. In particular, if $X$ is fibrant in the $G$-global flat model structure, then $X$ is a $G$-$\Omega$-spectrum.\index{G-global Omega-spectrum@$G$-global $\Omega$-spectrum}

In addition, the adjunctions induced by $(\ref{eq:flat-global-vs-equivariant})$ and $(\ref{eq:proj-global-vs-equivariant})$ on associated quasi-categories are equivalent in a preferred way and moreover right Bousfield localizations at the $G$-equivariant weak equivalences.
\end{prop}

A key ingredient is the following consequence of the results of Subsection~\ref{subsec:smash-product}:

\begin{prop}\label{prop:flat-projective-shift}
Let $A$ be a finite faithful $G$-set. Then
\begin{equation}\label{eq:smash-shift}
\Sigma^\bullet I(A,\blank)\smashp\blank \colon \cat{$\bm G$-Spectra}_{\textup{$G$-global flat}}\rightleftarrows\cat{$\bm G$-Spectra}_{\textup{$G$-global proj.}} :\!\Omega^A\sh^A
\end{equation}
is a simplicial Quillen equivalence. Moreover, the adjunction induced on associated quasi-categories is canonically equivalent to the identity adjunction.
\begin{proof}
Forgetting about the action on $A$, we have for any spectra $X$, $Y$ a sequence of $\cat{SSet}$-enriched natural isomorphisms
\begin{align*}
\Maps(\Sigma^\bullet_+ I(A,\blank)\smashp X, Y)&\hskip-.5pt\cong\hskip-.5pt \Maps(\Sigma^\bullet_+ I(A,\blank), F(X,Y))\hskip-.5pt\cong\hskip-.5pt
\Maps(I(A,\blank), \Omega^\bullet F(X,Y))\\
&\hskip-.5pt\cong\hskip-.5pt\Omega^A\Maps(X,\sh^AY)\cong\Maps(X,\Omega^A\sh^AY)
\end{align*}
where we have used in this order: the smash-function spectrum adjunction; the adjunction $\Sigma^\bullet_+\dashv\Omega^\bullet$; the $\cat{SSet}$-enriched Yoneda Lemma together with the definition of $\Omega^\bullet$ and $F$; the fact that $\Omega^A$ is defined in terms of the cotensoring of $\cat{Spectra}$ over $\cat{SSet}_*$. If $G$ acts on $X$, $Y$, and $A$, then these isomorphisms are $G$-equivariant with respect to the conjugate action by naturality, and taking $G$-fixed points of this thus witnesses that $(\ref{eq:smash-shift})$ is indeed a simplicial adjunction.

As $A$ is faithful, $I(A,\blank)\cong I(A,\blank)\times_\id G$ is cofibrant in the $G$-global model structure on $\cat{$\bm G$-$\bm I$-SSet}$, and hence $\Sigma^\bullet_+I(A,\blank)$ is $G$-globally projectively cofibrant by Proposition~\ref{prop:suspension-loop-projective}; thus, $(\ref{eq:smash-shift})$ is a Quillen adjunction by Proposition~\ref{prop:smash-mixed}.

To finish the proof it suffices to provide an equivalence $\textbf{R}\,\Omega^A\sh^A\simeq\id$, for which it is enough to give a natural levelwise weak equivalence between
\begin{equation*}
\Omega^A\sh^A\colon(\cat{$\bm G$-Spectra}_{\textup{$G$-global projective}})^f\to\cat{$\bm G$-Spectra}_{\textup{$G$-global flat}}
\end{equation*}
and the respective inclusion. A preferred such choice is provided by the maps $\lambda$ from the proof of Lemma~\ref{lemma:G-global-between-G-Omega}.
\end{proof}
\end{prop}

\begin{proof}[Proof of Proposition~\ref{prop:flat-fibrant-G-Omega}]
Let us show that the simplicial adjunction $(\ref{eq:flat-global-vs-equivariant})$ is a Quillen adjunction. Composing with the simplicial Quillen adjunction from Proposition~\ref{prop:flat-projective-shift} will then show that the same holds true for $(\ref{eq:proj-global-vs-equivariant})$.

We already know from Lemma~\ref{lemma:flat-implies-equiv-proj} that this is a Quillen adjunction for the respective level model structures. It therefore suffices to prove that any $X$ fibrant in the $G$-global flat model structure is a $G$-$\Omega$-spectrum.

By the proof of the previous proposition, $\lambda\colon X\to\Omega^A\sh^AX$ is a $G$-global weak equivalence, and both source and target are fibrant in the $G$-global flat model structure, so Brown's Factorization Lemma \cite[I.3]{brown-factorization} asserts that $\lambda$ factors as $\lambda=ps$ where $p$ is an acyclic fibration and $s$ is a section of an acyclic fibration. As acyclic fibrations in the $G$-global flat model structure are in particular strong level weak equivalences, we conclude that also $\lambda$ is a strong level weak equivalence. As $G$-$\Omega$-spectra are obviously closed under these, it then suffices that $\Omega^A\sh^AX$ is a $G$-$\Omega$-spectrum, which is obvious from the explicit characterization of the fibrant objects provided by Theorem~\ref{thm:G-global-flat}. This completes the proof that $(\ref{eq:flat-global-vs-equivariant})$ is a Quillen adjunction.

Proposition~\ref{prop:flat-projective-shift} already implies that the induced adjunctions are canonically equivalent, so it only remains to show that $(\ref{eq:flat-global-vs-equivariant})$ induces a right Bousfield localization at the $G$-equivariant weak equivalences. This is however obvious because
\begin{equation*}
\id^\infty=\textbf{R}\id\colon\cat{$\bm G$-Spectra}_{\textup{$G$-global flat}}^\infty\to\cat{$\bm G$-Spectra}_{\textup{$G$-equivariant projective}}^\infty
\end{equation*}
is evidently a quasi-localization at these.
\end{proof}

On the other hand one immediately concludes from the definitions:

\begin{cor}\label{cor:forget-has-right-adjoint}
Let $G$ be a finite group. Then the simplicial adjunction
\begin{equation*}
\id\colon\cat{$\bm G$-Spectra}_\textup{$G$-global flat}\rightleftarrows\cat{$\bm G$-Spectra}_\textup{$G$-equivariant flat} :\!\id.
\end{equation*}
is a Quillen adjunction. The induced adjunction on associated quasi-categories is a (left) Bousfield localization at the $G$-equivariant weak equivalences.\qed
\end{cor}

As an upshot of the above we conclude:

\begin{cor}\label{cor:und-adjoints}
Let $G$ be a discrete group (not necessarily finite) and let $\phi\colon H\to G$ be a homomorphism from a finite group to $G$. Then the simplicial adjunctions
\begin{align}\label{eq:und-has-right-adjoint}
\und_\phi\mathrel{:=}\phi^*\colon\cat{$\bm G$-Spectra}_{\textup{$G$-global flat}}&\rightleftarrows\cat{$\bm H$-Spectra}_{\textup{$H$-equivariant flat}} :\!\phi_*\\\intertext{and}
\phi_!(\Sigma^\bullet_+I(A,\blank)\smashp\blank)\colon\cat{$\bm H$-Spectra}_{\textup{$H$-equiv.~proj.}}&\rightleftarrows\cat{$\bm G$-Spectra}_{\textup{$G$-gl.~proj.}} :\!\Omega^A\sh^A\phi^*\nonumber
\end{align}
\nomenclature[auphi]{$\und_\phi$ (also $\und_H$)}{underlying $H$-equivariant spectrum of a $G$-global spectrum for $\phi\colon H\to G$ or $H\subset G$}%
are Quillen adjunctions for any finite faithful $H$-set $A$ (e.g.~$A=H$ with the left regular action). In particular, $\und_\phi^\infty$ admits both a left and a right adjoint.
\end{cor}

For $G=1$, this (or rather its analogue in the world of orthogonal spectra and with respect to all compact Lie groups) is sketched in \cite[Remark~4.5.25]{schwede-book} and an alternative proof of the existence of adjoints on the level of the homotopy categories is spelled out as Theorem~4.5.24 of \emph{op.~cit.}

\begin{proof}
It is immediate from Lemma~\ref{lemma:restriction-left-Quillen-flat} together with Corollary~\ref{cor:forget-has-right-adjoint} that $\phi^*$ in $(\ref{eq:und-has-right-adjoint})$ is left Quillen, proving the first statement.

For the second statement it suffices to factor this as
\begin{equation*}
\begin{tikzcd}[cramped]
\cat{$\bm H$-Spectra}_{\textup{$H$-equiv.~proj.}}\arrow[r, "{\Sigma^\bullet_+I(A,\blank)\smashp\blank}\,", shift left=2.5pt] &[2em] \arrow[l, "\Omega^A\sh^A", shift left=2.5pt] \cat{$\bm H$-Spectra}_{\textup{$H$-gl.~proj.}}\arrow[r, "\phi_!", shift left=2.5pt] &[-7pt] \arrow[l, "\phi^*", shift left=2.5pt] \cat{$\bm G$-Spectra}_{\textup{$G$-gl.~proj.}}
\end{tikzcd}
\end{equation*}
and then appeal to Proposition~\ref{prop:flat-fibrant-G-Omega} and Lemma~\ref{lemma:restriction-right-Quillen-projective}. Finally, Proposition~\ref{prop:flat-projective-shift} identifies the right adjoint of the induced adjunction with $\und_\phi^\infty$, which shows that $\textbf{L}\phi_!(\Sigma^\bullet_+I(A,\blank)\smashp\blank)$ is the desired left adjoint.
\end{proof}
\index{G-equivariant stable homotopy theory@$G$-equivariant stable homotopy theory!vs G-global stable homotopy theory@vs.~$G$-global stable homotopy theory|)}

\section{Delooping}\label{sec:delooping-group-completion}
Segal's classical Delooping Theorem \cite[Proposition~3.4]{segal-gamma} exhibits the homotopy theory of connective spectra as an explicit Bousfield localization of the homotopy theory of $\Gamma$-spaces, also see \cite[Theorem~5.8]{bousfield-friedlander} for a model categorical formulation. A $G$-equivariant version of this (for any finite group $G$) appeared as \cite[Theorem~6.5]{equivariant-gamma}. In this section we will strengthen this to a $G$-global comparison for any discrete group $G$.

\subsection{Non-equivariant and equivariant deloopings}
We begin by briefly recalling the usual equivariant and non-equivariant theory.

\begin{constr}
The restriction $\FUN(\cat{SSet}_*,\cat{SSet}_*)\to\cat{$\bm\Gamma$-SSet}_*$ along the evident embedding $\Gamma\hookrightarrow\cat{SSet}_*$ admits a fully faithful left adjoint via $\cat{SSet}_*$-enriched Kan extension. We call this the \emph{prolongation of $X$}.\index{Gamma-space@$\Gamma$-space!prolongation to SSet@prolongation to $\cat{SSet}_*$|textbf} As a Kan extension along a fully faithful functor, the prolongation of $X$ agrees with $X$ on $\Gamma$ up to canonical isomorphism, so we will not distinguish it notationally from the original $\Gamma$-space.

Explicitly, if $X\in\cat{$\bm\Gamma$-SSet}_*$ and $K$ is any pointed simplicial set, then we define $X(K)$ as the $\cat{SSet}_*$-enriched coend (or, equivalently, $\cat{SSet}$-enriched coend)
\begin{equation*}
\int^{S_+\in\Gamma} F(S_+)\times K^{\times S};
\end{equation*}
here the contravariant functoriality of $S_+\mapsto K^{\times S}$ is induced by the canonical identification of $K^{\times S}$ with the simplicial set $\Maps(S_+,K)$ of base-point preserving maps.

By the enriched functoriality of coends, we obtain an $\cat{SSet}_*$-enriched functor $\cat{SSet}_*\to\cat{SSet}_*$, and any map $X\to Y$ of $\Gamma$-spaces induces $X(K)\to Y(K)$ for any $K\in\cat{SSet}_*$.

Now let $G$ be a finite group. If $X$ is a $\Gamma$-$G$-space and $K$ is any pointed $G$-simplicial set, then we make $X(K)$ into a pointed $G$-simplicial set via the diagonal of the two actions.\index{Gamma-G-space@$\Gamma$-$G$-space!prolongation to G-SSet@prolongation to $\cat{$\bm G$-SSet}_*$|textbf}
\end{constr}

It is a classical observation that for finite $K$ there exists a natural isomorphism between $X(K)$ and the diagonal of the bisimplicial set $(m,n)\mapsto X(K_m)_n$ (with the evident functoriality in $m$ and $n$), see e.g.~\cite[p.~331]{schwede-gamma}. Thus, the equivariant Diagonal Lemma (Lemma~\ref{lemma:equivariant-diagonal}) immediately implies, also see~\cite[Lemma~4.8]{equivariant-gamma}:

\begin{cor}\label{cor:prolongation-homotopical}
If $X\to Y$ is a $G$-equivariant level weak equivalence, then the induced map $X(K)\to Y(K)$ is a $G$-equivariant weak equivalence for any finite pointed $G$-simplicial set $K$.\qed
\end{cor}

In fact, we can get rid of the finiteness assumption by filtering $K$ appropriately, but we will only need the above version.

\begin{constr}
We recall that a (symmetric) spectrum $Y$ is an $\cat{SSet}_*$-enriched functor $\bm\Sigma\to\cat{SSet}_*$. As the prolongation of any $\Gamma$-space $X$ is $\cat{SSet}_*$-enriched, we can therefore define $X(Y)\mathrel{:=}X\circ Y\colon\bm\Sigma\to\cat{SSet}_*$. This becomes a functor in $X$ and $Y$ in the obvious way.

Taking the convention that coends and limits in functor categories are computed pointwise, the above literally agrees with the $\cat{SSet}_*$-enriched coend
\begin{equation*}
\int^{S_+\in\Gamma} X(S_+)\smashp Y^{\times S}.
\end{equation*}
We can also describe $X(Y)$ explicitly as follows: if $B$ is a finite set, then $X(Y)(B)= X(Y(B))$ with the evident $\Sigma_B$-action; for any further finite set $A$, the structure map $S^{A}\smashp X(Y)(B)= S^A\smashp X(Y(B))\to X(Y(A\amalg B))= X(Y)(A\amalg B)$ is given by the composite
\begin{equation*}
S^A\smashp X(Y(B))\xrightarrow{\textup{asm}} X(S^A\smashp Y(B))\xrightarrow{X(\sigma)} X(Y(A\amalg B))
\end{equation*}
where the left hand arrow is the \emph{assembly map}\index{assembly map|textbf}\nomenclature[aasm]{asm}{assembly map (for an $\cat{SSet}_*$-enriched functor, e.g.~prolongation of a $\Gamma$-space)} induced by the $\cat{SSet}_*$-enrichment.
\end{constr}

Again, we can evaluate any $\Gamma$-$G$-space $X$ at any $G$-spectrum by applying the above construction and then pulling through the actions.

\begin{defi}
The \emph{associated spectrum} of a $\Gamma$-$G$-space $X$ is the $G$-spectrum $X(\mathbb S)$. We write $\mathcal E_G\colon\cat{$\bm\Gamma$-$\bm G$-SSet}_*\to\cat{$\bm G$-Spectra}$\index{associated spectrum|seealso {delooping}}\index{delooping!G-equivariant@$G$-equivariant|seeonly{$G$-equivariant delooping}}\index{G-equivariant delooping@$G$-equivariant delooping|textbf}\nomenclature[aE3G]{$\mathcal E_G$}{$G$-equivariant delooping of a $\Gamma$-$G$-space (via evaluating at the sphere spectrum)} for the functor $X\mapsto X(\mathbb S)$.
\end{defi}

Note that $X(\mathbb S)(\varnothing)=X(S^0)\cong X(1^+)$ is the `underlying $G$-space' of the $\Gamma$-$G$-space $X$. For suitable $X$, we can think of the remaining data as an \emph{equivariant delooping} of this $G$-space:

\begin{defi}
A $\Gamma$-$G$-space is called \emph{very special}\index{special!very special!G-equivariantly@$G$-equivariantly|seeonly{$\Gamma$-$G$-space, very special}}\index{very special!G-equivariantly@$G$-equivariantly|seeonly{$\Gamma$-$G$-space, very special}}\index{Gamma-G-space@$\Gamma$-$G$-space!very special|textbf} if it is special and the monoid $\pi_0^H(X(1^+))$ from Remark~\ref{rk:Gamma-G-monoid} is a group for all $H\subset G$.
\end{defi}

The following theorem was originally proven in slightly different form by Shimakawa \cite[Theorem~B]{shimakawa} (using the bar construction instead of the above prolongation), generalizing a non-equivariant result due to Segal \cite[Proposition~1.4]{segal-gamma}.

\begin{thm}\label{thm:equivariant-delooping-Omega}
If $X$ is very special, then $X(\mathbb S)$ is a $G$-$\Omega$-spectrum.
\begin{proof}
See~e.g.~\cite[Proposition~5.7]{equivariant-gamma}.
\end{proof}
\end{thm}

In \cite[Theorem~6.2]{equivariant-gamma}, Ostermayr constructs a \emph{stable model structure} on $\cat{$\bm\Gamma$-$\bm G$-SSet}_*$ with the same cofibrations as the $G$-equivariant level model structure (see Proposition~\ref{prop:Gamma-G-level}) and whose fibrant objects are precisely the very special level fibrant $\Gamma$-$G$-spaces.\index{Gamma-G-space@$\Gamma$-$G$-space!stable model structure|seeonly{$G$-equivariant stable model structure, on $\cat{$\bm\Gamma$-$\bm G$-SSet}_*$}}
\index{G-equivariant stable model structure@$G$-equivariant stable model structure!on Gamma-G-SSet@on $\cat{$\bm\Gamma$-$\bm G$-SSet}_*$}

\begin{prop}
The functor $\mathcal E_G$ is part of a simplicial Quillen adjunction
\begin{equation*}
\mathcal E_G\colon(\cat{$\bm\Gamma$-$\bm G$-SSet}_*)_{\textup{stable}}\rightleftarrows\cat{$\bm G$-Spectra}_{\textup{$G$-equivariant}} :\!\Phi_G.
\end{equation*}
\begin{proof}
\cite[Proposition~5.2]{equivariant-gamma} shows this for the level model structures; the claim follows as the right adjoint sends fibrant objects to very special $\Gamma$-$G$-spaces by \cite[proof of Theorem~5.9]{equivariant-gamma}.
\end{proof}
\end{prop}

Together with Corollary~\ref{cor:prolongation-homotopical} we conclude that $\mathcal E_G$ is in fact homotopical with respect to the above model structures.

\begin{prop}
The $G$-spectrum $X(\mathbb S)$ is \emph{connective}\index{connective!G-equivariantly@$G$-equivariantly|seeonly{$G$-equivariantly connective}}\index{G-equivariantly connective@$G$-equivariantly connective|textbf} for any $\Gamma$-$G$-space $X$, i.e.~$\underline\pi_*X'$ vanishes in negative degrees for some (hence any) $G$-equivariant weak equivalence $X(\mathbb S)\to X'$ to a $G$-$\Omega$-spectrum $X'$.
\begin{proof}
\cite[Corollary~5.5]{equivariant-gamma} shows that the negative na\"ive homotopy groups of $Y(\mathbb S)$ vanish for any $Y\in\cat{$\bm\Gamma$-$\bm G$-SSet}_*$.

To prove the proposition, we now simply choose the replacement $X(\mathbb S)\to X'$ as the image under $\mathcal E_G$ of a fibrant replacement $X\to Y$ in the stable model structure on $\cat{$\bm\Gamma$-$\bm G$-SSet}_*$; here we used that $X'\to Y$ is indeed a weak equivalence as $\mathcal E_G$ is homotopical and that $X'=Y(\mathbb S)$ is a $G$-$\Omega$-spectrum by Theorem~\ref{thm:equivariant-delooping-Omega}.
\end{proof}
\end{prop}

\begin{rk}
One can in fact show that the na\"ive homotopy groups $\underline\pi_*(X(\mathbb S))$ agree with the homotopy groups of $X'$ (i.e.~$X(\mathbb S)$ is \emph{semistable})\index{semistable!symmetric spectrum}, but we will not need this below.
\end{rk}

The following comparison is the key result on the homotopy theory of very special $\Gamma$-$G$-spaces:

\begin{thm}\label{thm:equivariant-delooping}
The adjunction $(\mathcal E_G)^\infty\dashv\cat{R}\Phi_G$ restricts to a Bousfield localization
\begin{equation*}
(\mathcal E_G)^\infty\colon(\cat{$\bm\Gamma$-$\bm G$-SSet}_*)_{\textup{level w.e.}}^\infty\rightleftarrows(\cat{$\bm G$-Spectra}^{\ge0})^\infty :\cat{R}\Phi_G
\end{equation*}
where $\cat{$\bm G$-Spectra}^{\ge0}$ denotes the subcategory of connective $G$-spectra. Moreover, the essential image of $\cat{R}\Phi_G$ consists precisely of the very special $\Gamma$-$G$-spaces.
\begin{proof}
This is immediate from the model categorical statement proven in \cite[Proposition~5.2 and Theorem~6.5]{equivariant-gamma}, also see \cite[Theorem~5.8]{bousfield-friedlander} for the special case $G=1$.
\end{proof}
\end{thm}

\subsection{The $\bm G$-global Delooping Theorem} Before we can properly state a $G$-global version of the Delooping Theorem, we first need to understand what it means to be `connective' or `very special' in the $G$-global context.

\begin{defi}\index{connective!G-globally@$G$-globally|seeonly{$G$-globally connective}}\index{G-globally connective@$G$-globally connective|textbf}
A $G$-spectrum $X$ is \emph{$G$-globally connective} if $\phi^*X$ is $H$-equivariantly connective for all finite groups $H$ and all homomorphisms $\phi\colon H\to G$.
\end{defi}

\begin{defi}
We call $X\in\cat{$\bm\Gamma$-$\bm G$-$\bm{\mathcal I}$-SSet}_*$ \emph{very special}\index{very special!G-globally@$G$-globally|seeonly{$G$-global $\Gamma$-space, very special}}\index{G-global Gamma-space@$G$-global $\Gamma$-space!very special|textbf}\index{special!very special!G-globally@$G$-globally|seeonly{$G$-global $\Gamma$-space, very special}} if $\und_\phi X$ is very special for every finite group $H$ and every homomorphism $\phi\colon H\to G$.
\end{defi}

Put differently (see Lemma~\ref{lemma:special-g-global-u-phi}), $X$ is very special if and only if it is special and the abelian monoid structure on $\pi_0^H((\phi^*X)(\mathcal U_H)(1^+))$ induced by
\begin{equation*}
\begin{tikzcd}
X(1^+)\times X(1^+) & \arrow[l, "\sim", "\rho"'] X(2^+)\arrow[r, "X(\mu)"] & X(1^+)
\end{tikzcd}
\end{equation*}
is a group structure for every finite group $H$ and every homomorphism $\phi\colon H\to G$.

Now we introduce the $G$-global delooping functor, cf.~\cite[Construction~3.3]{schwede-k-theory}:

\begin{constr}\index{delooping!G-global@$G$-global|seeonly{$G$-global delooping}}\index{G-global delooping@$G$-global delooping|textbf}
We define $\mathcal E^\otimes\colon\cat{$\bm\Gamma$-$\bm G$-$\bm{\mathcal I}$-SSet}_*\to\cat{$\bm G$-Spectra}$\nomenclature[aE2o]{$\mathcal E^\otimes$}{$G$-global delooping of a $G$-global $\Gamma$-space (in terms of tensoring over $\cat{$\bm G$-$\bm{\mathcal I}$-SSet}$)} via the $\cat{SSet}_*$-enriched coend
\begin{equation*}
\mathcal E^\otimes(X)=\int^{T_+\in\Gamma}\mathbb S^{\times T}\otimes X(T_+),
\end{equation*}
together with the evident $\cat{SSet}_*$-enriched functoriality.
\end{constr}

\begin{rk}\label{rk:E-otimes-levelwise-description}
We again take the convention that the above coend is constructed by forming the levelwise coend. In this case,
\begin{equation*}
\mathcal E^\otimes(X)(A)=\int^{T_+\in\Gamma}(S^A)^{\times T} \smashp X(A)(T_+)=X(A)(S^A),
\end{equation*}
and the structure map $\mathcal E^\otimes(X)(B)\to \mathcal E^\otimes(X)(A\amalg B)$ is given by the diagonal composite
\begin{equation*}
\hskip-8pt\begin{tikzcd}
S^A\smashp X(B)(S^B)\arrow[d, "S^A\smashp X(i)(S^B)"']\arrow[r, "\textup{asm}"] & X(B)(S^A\smashp S^B) \arrow[r, "\cong"] & X(B)(S^{A\amalg B})\arrow[d, "X(i)(S^{A\amalg B})"]\\
S^A\smashp X(A\amalg B)(S^B)\arrow[r, "\textup{asm}"'] & X(A\amalg B)(S^A\smashp S^B)\arrow[r, "\cong"'] & X(A\amalg B)(S^{A\amalg B})
\end{tikzcd}
\end{equation*}
where $i\colon B\hookrightarrow A\amalg B$ is the inclusion.
\end{rk}

\subsubsection{Comparison to equivariant deloopings} We now want to relate this to Shimakawa's equivariant version of Segal's machinery, for which we begin with the following trivial observation:

\begin{rk}\label{rk:E-otimes-E-G-const}
If $X$ is any $\Gamma$-$G$-space, then $\mathcal E^\otimes(\const X)=\mathcal E_G(X)$. More precisely, the diagram
\begin{equation*}
\begin{tikzcd}
\cat{$\bm\Gamma$-$\bm G$-SSet}_*\arrow[d, "\const"']\arrow[r, "\mathcal E_G"] & \cat{$\bm G$-Spectra}\\
\cat{$\bm\Gamma$-$\bm G$-$\bm{\mathcal I}$-SSet}_*\arrow[ur, bend right=15pt, "\mathcal E^\otimes"']
\end{tikzcd}
\end{equation*}
of simplicially enriched functors is strictly commutative.
\end{rk}

\begin{cor}
Let $X\in \cat{$\bm\Gamma$-$\bm G$-$\bm{\mathcal I}$-SSet}_*$ be fibrant in the injective $G$-global model structure. Then the degree zero inclusion induces a $G$-equivariant level weak equivalence $\mathcal E_G(X(\varnothing))=\mathcal E^\otimes(\const X(\varnothing))\to\mathcal E^\otimes(X)$.
\begin{proof}
We have to show that the inclusion induces a $G$-equivariant weak equivalence $X(\varnothing)(S^A)\to X(A)(S^A)$ for every finite $G$-set $A$. But $X(\varnothing)\to X(A)$ is a $G$-equivariant level equivalence of $\Gamma$-$G$-spaces by Corollary~\ref{cor:G-injective-constant}, so the claim follows from Corollary~\ref{cor:prolongation-homotopical}.
\end{proof}
\end{cor}

\begin{cor}\label{cor:E-otimes-G-Omega}
Let $G$ be finite and assume $X\in\cat{$\bm\Gamma$-$\bm G$-$\bm{\mathcal I}$-SSet}_*$ is fibrant in the injective $G$-global model structure and moreover $G$-globally very special. Then $\mathcal E^\otimes(X)$ is a $G$-$\Omega$-spectrum.
\begin{proof}
By the previous corollary, $\mathcal E^\otimes(X)$ is $G$-equivariantly level equivalent to $\mathcal E_G(X(\varnothing))$. The claim now follows from Theorem~\ref{thm:equivariant-delooping-Omega} as $X(\varnothing)$ is $G$-equivariantly level equivalent to $\und_GX$ by Corollary~\ref{cor:injective-Gamma-varnothing-UG}, hence very special by definition.
\end{proof}
\end{cor}

For a general comparison we introduce:

\begin{constr}
Let $H$ be a finite group and let $\phi\colon H\to G$ be a homomorphism. We define $\mathcal E_\phi\colon\cat{$\bm\Gamma$-$\bm G$-$\bm{\mathcal I}$-SSet}_*\to\cat{$\bm H$-Spectra}$ as follows: if $A$ is any finite set and $X$ is any $\Gamma$-$G$-$\mathcal I$-space, then $\mathcal E_\phi(X)(A)= (\phi^*X)(\mathcal U\amalg A)(S^A)$, where $\mathcal U$ is a fixed complete $H$-set universe, e.g.~the one from the definition of the underlying $H$-simplicial set of a $G$-$\mathcal I$-simplicial set. The structure maps of $\mathcal E_\phi(X)$ are defined analogously to Remark~\ref{rk:E-otimes-levelwise-description}, and the functoriality in $X$ is the obvious one.

We now consider the zig-zag of natural transformations
\begin{equation}\label{eq:und-E-otimes-vs-E-H}
{\und_\phi}\circ\mathcal E^\otimes=\phi^*\circ\mathcal E^\otimes\Rightarrow\mathcal E_\phi\Leftarrow \mathcal E_H\circ\ev_{\mathcal U}\circ\phi^*\cong\mathcal E_H\circ\und_\phi
\end{equation}
induced in degree $A$ by the inclusions $A\hookrightarrow\mathcal U\amalg A\hookleftarrow \mathcal U$.
\end{constr}

The following in particular generalizes (the simplicial analogue of) \cite[Theorem~3.13]{schwede-k-theory}:

\begin{lemma}\label{lemma:und-E-vs-E-H}
Both maps in $(\ref{eq:und-E-otimes-vs-E-H})$ are $\underline\pi_*$-isomorphisms (and hence in particular $H$-equivariant weak equivalences).
\begin{proof}
Analogously to Lemma~\ref{lemma:tensor-vs-smash-homotopy-groups} one gets natural isomorphisms
\begin{equation*}
\underline\pi_*\und_\phi\mathcal E^\otimes(X)\cong
\underline\pi_*\mathcal E_H\big((\phi^*X)(\mathcal U)\big)\qquad\text{and}\qquad
\underline\pi_*(\mathcal E_\phi X)\cong \underline\pi_*\mathcal E_H\big((\phi^*X)(\mathcal U\amalg\mathcal U)),
\end{equation*}
for every $X\in\cat{$\bm\Gamma$-$\bm G$-$\bm{\mathcal I}$-SSet}_*$, and under this identification the actions of the maps $(\ref{eq:und-E-otimes-vs-E-H})$ on homotopy groups are induced by the two inclusions $\mathcal U\hookrightarrow\mathcal U\amalg\mathcal U\hookleftarrow\mathcal U$. Thus, it is enough to show that the induced maps $(\phi^*X)(\mathcal U)\rightrightarrows(\phi^*X)(\mathcal U\amalg\mathcal U)$ are $H$-equivariant weak equivalences of $\Gamma$-$H$-spaces, i.e.~if $T$ is any finite $H$-set, then the induced maps $(\phi^*X)(\mathcal U)(T_+)\rightrightarrows(\phi^*X)(\mathcal U\amalg\mathcal U)(T_+)$ are $H$-equivariant weak equivalences. But this is simply an instance of Lemma~\ref{lemma:evaluation-h-universe} applied to the $H$-$\mathcal I$-simplicial set $(\phi^*X)(\blank)(T_+)$.
\end{proof}
\end{lemma}

\begin{cor}\label{cor:E-otimes-homotopical}\label{cor:E-otimes-connective}
The functor $\mathcal E^\otimes$ sends $G$-global level weak equivalences to $G$-global weak equivalences. Moreover, it takes values in $G$-globally connective $G$-spectra.\qed
\end{cor}

Now we can finally state the main result of this section:

\begin{thm}\label{thm:group-completion}
Let $G$ be any group. Then
\begin{equation*}
(\mathcal E^\otimes)^\infty\colon(\cat{$\bm\Gamma$-$\bm G$-$\bm{\mathcal I}$-SSet}_*)_{\textup{$G$-global level}}^\infty\to(\cat{$\bm G$-Spectra}^{\ge0}_{\textup{$G$-global}})^\infty
\end{equation*}
(where the superscript `{\small$\!{}\ge0\mskip-.67\thinmuskip$}' denotes the subcategory of $G$-globally connective spectra) has a fully faithful right adjoint $\cat{R}\Phi^\otimes$,\nomenclature[aPhio]{$\Phi^\otimes$}{right adjoint to $\mathcal E^\otimes$} yielding a Bousfield localization. Moreover:
\begin{enumerate}
\item The essential image of $\textbf{\textup R}\Phi^\otimes$ consists precisely of the very special $G$-global $\Gamma$-spaces.
\item $(\mathcal E^\otimes)^\infty$ inverts a map $f$ if and only if $\und_\phi f$ is a stable weak equivalence for all finite groups $H$ and all homomorphisms $\phi\colon H\to G$.
\item For any finite group $H$ and any homomorphism $\phi\colon H\to G$ there are preferred equivalences filling
\begin{equation*}
\begin{tikzcd}[cramped]
\cat{$\bm\Gamma$-$\bm G$-$\bm{\mathcal I}$-SSet}_*^\infty\arrow[d, "(\mathcal E^\otimes)^\infty"']\arrow[r, "\und_\phi^\infty"] &[-.8em] \cat{$\bm\Gamma$-$\bm H$-SSet}_*^\infty\arrow[d, "\mathcal E_H^\infty"]\\
\cat{$\bm G$-Spectra}_{\textup{$G$-gl.}}^\infty\arrow[r, "\und_\phi^\infty\;"'] & \cat{$\bm H$-Spectra}_{\textup{$H$-equiv.}}^\infty
\end{tikzcd}
\hskip1.1em\hskip-1.925pt
\begin{tikzcd}[cramped]
\cat{$\bm\Gamma$-$\bm G$-$\bm{\mathcal I}$-SSet}_*^\infty\arrow[r, "\und_\phi^\infty"] &[-.8em] \cat{$\bm\Gamma$-$\bm H$-SSet}_*^\infty\\
\cat{$\bm G$-Spectra}_{\textup{$G$-gl.}}^\infty\arrow[r, "\und_\phi^\infty\;"']\arrow[u, "\textbf{\textup R}\Phi^\otimes"] & \cat{$\bm H$-Spectra}_{\textup{$H$-equiv.}}^\infty\arrow[u, "\textbf{\textup R} \Phi_H"']
\end{tikzcd}
\end{equation*}
and these are moreover canonical mates of each other.
\end{enumerate}
\end{thm}

The proof of the theorem is somewhat involved and occupies the rest of this section, so let us briefly outline the general strategy. In \ref{subsubsec:deloop-ra}, we will construct the right adjoint $\cat{R}\Phi^\otimes$ of $(\mathcal E^\otimes)^\infty$ and establish its basic properties, in particular showing that it takes values in very special $G$-global $\Gamma$-spaces. The key argument to the proof is then contained in \ref{subsubsec:deloop-ff}, where we show that the counit of the above adjunction is an equivalence, i.e.~that this is indeed a Bousfield localization. From these two results together with the comparison to the equivariant deloopings given above, we will then be able to deduce the theorem in \ref{subsubsec:deloop-final} purely by abstract nonsense.

\subsubsection{$G$-global $\Gamma$-spaces from $G$-global spectra}\label{subsubsec:deloop-ra}
As the notation suggests, we will obtain $\cat{R}\Phi^\otimes$ as a right derived functor, and we begin by introducing the corresponding point-set level functor:

\begin{constr}
We define $\Phi^\otimes\colon\cat{$\bm G$-Spectra}\to\cat{$\bm\Gamma$-$\bm G$-$\bm{\mathcal I}$-SSet}_*$ as follows: if $X$ is a $G$-spectrum and $T$ is a finite set, then $(\Phi^\otimes X)(T_+)=F_{\mathcal I}(\mathbb S^{\times T},X)$. The functoriality in $T_+$ is induced by the natural identification $\mathbb S^{\times T}\cong\Maps(T_+,\mathbb S)$, and the $\cat{SSet}_*$-enriched functoriality in $X$ is the obvious one.

We define for every $G$-spectrum $X$ the map $\epsilon_X\colon\mathcal E^\otimes\Phi^\otimes X\to X$ as the one induced by the maps $\epsilon_{\mathbb S^{\times T},X}\colon\mathbb S^{\times T}\otimes F_{\mathcal I}(\mathbb S^{\times T},X)\to X$ for varying $T_+\in\Gamma$, where $\epsilon$ comes from Construction~\ref{constr:tensor-right-adjoint}. We omit the easy verification that this is well-defined and a $\cat{SSet}_*$-enriched natural transformation.

For a $G$-global $\Gamma$-space $Y$, we write $\eta_Y\colon Y\to\Phi^\otimes\mathcal E^\otimes Y$ for the map given in degree $T_+$ by the composition
\begin{equation*}
Y(T_+)\xrightarrow{\eta} F_{\mathcal I}\big(\mathbb S^{\times T},\mathbb S^{\times T}\otimes Y(T_+)\big)\to
F_{\mathcal I}\left(\mathbb S^{\times T},{\textstyle\int^{U_+\in\Gamma}\mathbb S^{\times U}\otimes Y(U_+)}\right)
\end{equation*}
where the right hand arrow is induced by the structure map of the coend and the left hand arrow is again from Construction~\ref{constr:tensor-right-adjoint}. We omit the easy verification that also $\eta$ is a well-defined simplicial transformation, and that $\epsilon$ and $\eta$ satisfy the triangle identities, yielding an enriched adjunction
\begin{equation}\label{eq:E-otimes-R-otimes}
\mathcal E^\otimes\colon \cat{$\bm\Gamma$-$\bm G$-$\bm{\mathcal I}$-SSet}_*\rightleftarrows\cat{$\bm G$-Spectra} :\!\Phi^\otimes.
\end{equation}
\end{constr}

\begin{rk}\label{rk:counit-degree-empty}
One immediately checks from the above definitions that the counit $\epsilon_X\colon \mathcal E^\otimes\Phi^\otimes X\to X$ is an \emph{isomorphism} in degree $\varnothing$ for any $G$-spectrum $X$.
\end{rk}

\begin{cor}\label{cor:E-otimes-left-Quillen}
The simplicial adjunction $(\ref{eq:E-otimes-R-otimes})$ is a Quillen adjunction with respect to the $G$-global injective model structures.
\begin{proof}
Using the description from Remark~\ref{rk:E-otimes-levelwise-description} it is obvious that $\mathcal E^\otimes$ preserves injective cofibrations. Moreover, it is homotopical by Corollary~\ref{cor:E-otimes-homotopical}.
\end{proof}
\end{cor}

We can also give an alternative, more tractable description of the right adjoint:

\begin{constr}
Define $\Phi^\Sigma\colon\cat{$\bm G$-Spectra}\to\cat{$\bm\Gamma$-$\bm G$-Spectra}_*$ (where the right hand side again denotes $\cat{Set}_*$-enriched functors) via $\Phi^\Sigma(X)(T_+)=F(\mathbb S^{\times T}, X)$; here the functoriality in $T_+$ is as before, and the $\cat{SSet}_*$-enriched functoriality in $X$ is the obvious one.

Using this, we can now define $\Phi^\smashp$\nomenclature[aPhis]{$\Phi^\smashp$}{right adjoint to $\mathcal E^\smashp$} as the composite
\begin{equation*}
\cat{$\bm G$-Spectra}\xrightarrow{\Phi^\Sigma}\cat{$\bm\Gamma$-$\bm G$-Spectra}_*\xrightarrow{\Omega^\bullet}\cat{$\bm\Gamma$-$\bm G$-$\bm I$-SSet}_*\xrightarrow{\Maps_I(\mathcal I,\blank)}\cat{$\bm\Gamma$-$\bm G$-$\bm{\mathcal I}$-SSet}_*.
\end{equation*}
\end{constr}

\begin{cor}
The maps $\widehat\psi$ from Construction~\ref{constr:psi-hat} define a natural transformation $\widehat\psi\colon\Phi^\otimes\Rightarrow\Phi^\smashp$. If $X$ is fibrant in the $G$-global injective model structure, then $\widehat\psi_X\colon \Phi^\otimes X\to\Phi^\smashp X$ is a $G$-global level weak equivalence.
\begin{proof}
As $\mathbb S^{\times T}$ is a flat $\Sigma_T$-spectrum for any finite set $T$ (Example~\ref{ex:prod-S-flat}), this follows by applying Theorem~\ref{thm:F-I-vs-F} levelwise.
\end{proof}
\end{cor}

\begin{rk}
Also $\Phi^\smashp$ admits a simplicial left adjoint, which can be computed via the coend $\mathcal E^\smashp(X)\mathrel{:=}\int^{T_+\in\Gamma}\mathbb S^{\times T}\smashp\Sigma^\bullet X$, and the resulting adjunction $\mathcal E^\smashp\dashv\Phi^\smashp$\nomenclature[aE2s]{$\mathcal E^\smashp$}{$G$-global delooping of a $G$-global $\Gamma$-space (in terms of smash product of $G$-global spectra)}\index{G-global delooping@$G$-global delooping} is left Quillen with respect to the $G$-global level model structure on the source and the $G$-global projective model structure on the target (this uses Proposition~\ref{prop:smash-mixed} again, also cf.~the argument in the proof of the proposition below).
As we will not need this below, we leave the details to the interested reader.
\end{rk}

\begin{prop}\label{prop:R-smash-very-special}
Let $X\in\cat{$\bm G$-Spectra}$ be fibrant in the $G$-global projective model structure. Then $\Phi^\smashp(X)$ is very special.
\begin{proof}
Let us first prove that $\Phi^\smashp(X)$ is special, for which we let $T$ be any finite set. Up to isomorphism, the Segal map $\Phi^\smashp(X)(T_+)\to\prod_T\Phi^\smashp(X)(1^+)$ is then given by
\begin{equation}\label{eq:segal-map}
\Maps_I\big(\mathcal I,\Omega^\bullet F(\mathbb S^{\times T},X)\big)\xrightarrow{\Maps_I(\mathcal I,\Omega^\bullet c^*)} \Maps_I\big(\mathcal I,\Omega^\bullet F(\mathbb S^{\vee T},X)\big)
\end{equation}
where $c$ is the canonical map $\mathbb S^{\vee T}\to\mathbb S^{\times T}$. The latter is a $\Sigma_T$-global weak equivalence (Proposition~\ref{prop:equivariant-additivity}) between flat $\Sigma_T$-spectra (Examples~\ref{ex:coprod-S-flat} and~\ref{ex:prod-S-flat}, respectively); as $X$ is fibrant in the $G$-global projective model structure, $c^*$ is therefore a $(G\times\Sigma_T)$-global weak equivalence between $(G\times\Sigma_T)$-globally projectively fibrant $(G\times\Sigma_T)$-spectra by Corollary~\ref{cor:function-spectrum-extra-homotopical} and Ken Brown's Lemma. As both $\Omega^\bullet$ and $\Maps_I(\mathcal I,\blank)$ are right Quillen with respect to the corresponding projective model structures, we conclude that also $(\ref{eq:segal-map})$ is a $(G\times\Sigma_T)$-global weak equivalence, and hence so is the Segal map.

Now let $\phi\colon H\to G$ be any homomorphism from a finite group to $G$ and let $\mathcal U=\smash{\mathcal U_H}$ be our fixed complete $H$-set universe. By the above, $\und_\phi(\Phi^\smashp X)=(\phi^*\Phi^\smashp X)(\mathcal U)$ is a special $\Gamma$-$H$-space, and to finish the proof we have to show that the induced monoid structure on $\pi_0^H(\phi^*\Phi^\smashp X(\mathcal U)(1^+))=\pi_0^H\big(\Maps_I(\mathcal I,\Omega^\bullet\phi^*F(\mathbb S,X))(\mathcal U)\big)$ is actually a group structure. As in the classical setting, this is an application of the Eckmann-Hilton argument: first, we consider the commutative diagram
\begin{gather}
\begin{aligned}\label{diag:gamma-vs-spectra-homotopy-groups}
\begin{tikzcd}[ampersand replacement=\&, row sep=1.7em]
\pi_0^H\big(\Maps_I(\mathcal I,\Omega^\bullet\phi^*F(\mathbb S,X))(\mathcal U)\big)\arrow[r] \& \pi_0^H(\phi^*F(\mathbb S,X))\\
\pi_0^H\big(\Maps_I(\mathcal I,\Omega^\bullet\phi^*F(\mathbb S\times\mathbb S,X))(\mathcal U)\big)\arrow[d, "{(p_{1*},p_{2*})}"', "\cong"]\arrow[u, "\mu_*"]\arrow[r] \& \pi_0^H(\phi^*F(\mathbb S\times\mathbb S,X))\arrow[u, "\mu_*"']\arrow[d, "{(p_{1*},p_{2*})}", "\cong"']\\
\pi_0^H\smash{\big(\Maps_I(\mathcal I,\Omega^\bullet\phi^*F(\mathbb S,X))(\mathcal U)\big)^{\times 2}}\arrow[r] \& \pi_0^H\smash{(\phi^*F(\mathbb S,X))^{\times 2}}
\end{tikzcd}
\end{aligned}\raisetag{-5.5pt}
\end{gather}
where the horizontal arrows are those from Remark~\ref{rk:omega-bullet-script-I-homotopy-groups} applied to the $H$-spectra $\phi^*F(\mathbb S,X)$ and $\phi^*F(\mathbb S\times\mathbb S,X)$. As both of these are fibrant in the $H$-global projective model structure by Proposition~\ref{prop:smash-mixed} together with Lemma~\ref{lemma:restriction-right-Quillen-projective}, the remark tells us (by invoking Theorem~\ref{thm:forget-left-quillen-global}) that the horizontal arrows are bijective.

By commutativity, we then conclude that the right hand vertical arrows equip $\pi_0^H(\phi^*F(\mathbb S,X))$ with the structure of an abelian monoid, and that the top horizontal arrow in $(\ref{diag:gamma-vs-spectra-homotopy-groups})$ is an isomorphism with respect to this monoid structure. But on the other hand, the vertical arrows on the right are group homomorphisms with respect to the usual group structure on $\pi_0^H$, so the Eckmann-Hilton argument implies that this monoid structure on $\pi_0^H(\phi^*F(\mathbb S, X))$ agrees with the standard group structure on $\pi_0^H$. We conclude that also the monoid structure on $\pi_0^H\big(\Maps_I(\mathcal I,\Omega^\bullet\phi^*F(\mathbb S,X))(\mathcal U)\big)=\pi_0^H(\phi^*\Phi^\smashp(X)(\mathcal U)(1^+))$ is a group structure, which completes the proof.
\end{proof}
\end{prop}

\begin{cor}\label{cor:R-otimes-very-special}
Let $X\in\cat{$\bm G$-Spectra}$ be fibrant in the $G$-global injective model structure. Then $\Phi^\otimes(X)$ is very special.\qed
\end{cor}

\subsubsection{Full faithfulness of $\cat{R}\Phi^\otimes$}\label{subsubsec:deloop-ff}
The following will be the main step in the proof of Theorem~\ref{thm:group-completion}:

\begin{prop}\label{prop:group-completion-counit}
Let $X\in\cat{$\bm G$-Spectra}^{\ge0}$ be fibrant in the $G$-global injective model structure. Then the counit $\epsilon\colon\mathcal E^\otimes\Phi^\otimes X\to X$ is a $G$-global weak equivalence.
\begin{proof}
Let $H$ be a finite group, let $\phi\colon H\to G$ be a homomorphism, and let $i\colon\phi^*X\to X'$ be a fibrant replacement in the injective $H$-global model structure.

\begin{claim*}
$\Phi^\otimes(i)\colon \Phi^\otimes(\phi^*X)\to\Phi^\otimes(X')$ is an $H$-global level weak equivalence.
\begin{proof}
By Corollary~\ref{cor:R-otimes-very-special}, the $H$-global $\Gamma$-space $\Phi^\otimes(X')$ is special, and so is the $G$-global $\Gamma$-space $\Phi^\otimes(X)$. It follows that also $\Phi^\otimes(\phi^*X)=\phi^*\Phi^\otimes(X)$ is special.

It is therefore enough to show that $\Phi^\otimes(i)(1^+)$ is an $H$-global weak equivalence. But this agrees up to conjugation by isomorphisms with $\Maps_I(\mathcal I,\Omega^\bullet(i))$, cf.~Construction~\ref{constr:psi-hat}. As both $\phi^*X$ and $X'$ are in particular fibrant in the $H$-global \emph{projective} model structure, the claim follows from Ken Brown's Lemma.
\end{proof}
\end{claim*}

To show that $\phi^*(\epsilon_X)=\epsilon_{\phi^*X}$ is an $H$-equivariant weak equivalence, we consider the naturality square
\begin{equation*}
\begin{tikzcd}
\mathcal E^\otimes\Phi^\otimes(\phi^*X)\arrow[r, "\epsilon_{\phi^*X}"]\arrow[d, "\mathcal E^\otimes\Phi^\otimes(i)"'] & \phi^*X\arrow[d, "i"]\\
\mathcal E^\otimes\Phi^\otimes(X')\arrow[r, "\epsilon_{X'}"'] & X'.
\end{tikzcd}
\end{equation*}
The right hand vertical arrow is an $H$-global weak equivalence, and so is the left hand vertical arrow by the above claim together with Corollary~\ref{cor:E-otimes-homotopical}. Thus, it is enough to show that $\epsilon_{X'}$ is an $H$-equivariant weak equivalence.

For this we observe that $\Phi^\otimes(X')$ is fibrant in the injective $H$-global model structure on $\cat{$\bm\Gamma$-$\bm{H}$-$\bm{\mathcal I}$-SSet}_*$ by Corollary~\ref{cor:E-otimes-left-Quillen} and moreover very special by Corollary~\ref{cor:R-otimes-very-special}. Thus, Corollary~\ref{cor:E-otimes-G-Omega} implies that $\mathcal E^\otimes\Phi^\otimes(X')$ is an $H$-$\Omega$-spectrum. In addition, Corollary~\ref{cor:E-otimes-connective} implies that it is also $H$-globally connective.

On the other hand, also $X'$ is an $H$-$\Omega$-spectrum by Proposition~\ref{prop:flat-fibrant-G-Omega}, and it is connective by assumption. As $\epsilon_{X'}(\varnothing)$ is an isomorphism, we immediately see that $\epsilon_{X'}$ is an $H$-equivariant $\underline\pi_*$-isomorphism, finishing the proof.
\end{proof}
\end{prop}

\subsubsection{Proof of the Delooping Theorem}\label{subsubsec:deloop-final}
Instead of the above approach, we could have also tried to prove the proposition by comparing the counit of our delooping adjunction to the counit of the usual equivariant one. However, we have carefully avoided such a comparison because the computations involved can become quite cumbersome. Instead, we will get all these compatibilities for free now thanks to general facts about Bousfield localizations. To formulate these, it will be convenient to forget about the above model structures and move to the world of quasi-categories completely (so that the above derived adjunctions just become ordinary adjunctions):

\begin{lemma}\label{lemma:abstract-mate}
Let
\begin{equation}\label{eq:big-adjunction}
F\colon\mathscr C\rightleftarrows\mathscr D :\!U
\end{equation}
be an adjunction of quasi-categories, let $I$ be a set (or more generally any class), and let $(S_i\colon\mathscr C\to\mathscr C_i)_{i\in I}$ and $(T_i\colon\mathscr D\to\mathscr D_i)_{i\in I}$ be two jointly conservative families. Assume moreover we are given for each $i\in I$ a Bousfield localization
\begin{equation*}
F_i\colon\mathscr C_i\rightleftarrows\mathscr D_i :\!U_i
\end{equation*}
such that $\essim S_iU\subset\essim U_i$, together with a natural equivalence $\phi_i$ filling
\begin{equation}\label{diag:comparison-base}
\begin{tikzcd}
\mathscr C\arrow[r, "S_i"]\arrow[d, "F"'] & \mathscr C_i\arrow[d, "F_i"]\twocell[dl, "\scriptscriptstyle\phi_i"{xshift=-6pt, yshift=4pt}, "\scriptscriptstyle\simeq"'{xshift=4pt, yshift=-4pt}]\\
\mathscr D\arrow[r, "T_i"'] & \mathscr D_i.
\end{tikzcd}
\end{equation}
Then the following are equivalent:
\begin{enumerate}
\item $U$ is fully faithful, i.e.~also $(\ref{eq:big-adjunction})$ is a Bousfield localization.
\item For all $i\in I$, the canonical mate $S_iU\Rightarrow U_iT_i$ of $(\ref{diag:comparison-base})$ is an equivalence.
\end{enumerate}
Moreover, in this case $X\in\mathscr D$ lies in the essential image of $U$ if and only if $S_iX$ lies in the essential image of $U_i$ for all $i\in I$, and a map $f$ is inverted by $F$ if and only if $S_if$ is inverted by $F_i$ for all $i\in I$.
\begin{proof}
$(1)\Rightarrow(2)$: For any $i\in I$, the canonical mate of $\phi_i$ is defined as the pasting
\begin{equation}\label{diag:comparison-mate}
\begin{tikzcd}
\mathscr D\arrow[r, "U"]\arrow[dr, "="'{name="eps-anchor"}, bend right=15pt] & \twocell[to="eps-anchor", "\scriptscriptstyle\epsilon"{xshift=-2pt,yshift=4pt}]
\mathscr C\arrow[r, "S_i"]\arrow[d, "F"'] & \mathscr C_i\arrow[d, "F_i"]\twocell[dl, "\scriptscriptstyle\phi_i"{xshift=-6pt, yshift=4pt}, "\scriptscriptstyle\simeq"'{xshift=4pt, yshift=-4pt}]\arrow[dr, bend left=15pt, "="{name="eta-anchor"}]&\\
& \mathscr D\arrow[r, "T_i"'] & \mathscr D_i\twocell[from="eta-anchor","\scriptscriptstyle\eta"{xshift=-2pt,yshift=4pt}] \arrow[r, "U_i"'] & \mathscr C_i.
\end{tikzcd}
\end{equation}
Since $U$ is fully faithful, $\epsilon$ is an equivalence, and so is $\phi_i$ by assumption. On the other hand, $S_iUX\in\essim U_i$ for any $X\in\mathscr D$ by assumption, so that also $\eta_i$ evaluated at $S_iUX$ is an equivalence because $U_i$ is fully faithful. We conclude that the above pasting is indeed an equivalence as desired.

$(2)\Rightarrow(1)$: We have to show that the counit $\epsilon_X$ of $F\dashv U$ is an equivalence for every $X\in\mathscr D$. As the $T_i$'s are jointly conservative, this is equivalent to $T_i\epsilon_X$ being an equivalence for all $i\in I$, and since $U_i$ is fully faithful, this is in turn equivalent to $U_iT_i\epsilon_X$ being an equivalence. But $(2)$ precisely tells us that `the' composition
\begin{equation*}
S_iUX\xrightarrow{(\eta_i)_{S_iUX}}U_iF_iS_iU\xrightarrow{U_i(\phi_i)_{UX}} U_iT_iFUX\xrightarrow{U_iT_i\epsilon_X} U_iT_iX
\end{equation*}
is an equivalence, and so are the first two arrows in it by the assumption on $\essim S_iU$ and $\phi_i$, respectively. The claim follows by $2$-out-of-$3$.

It only remains to prove the characterizations of the maps inverted by $F$ and the essential image of $U$ under these equivalent assumptions. But indeed, if $f$ is any morphism in $\mathscr C$, then $Ff$ is an equivalence if and only if each $T_iFf$ is, which by the equivalence $(\ref{diag:comparison-base})$ is equivalent to $F_iS_if$ being an equivalence.

Finally, if $X$ lies in the essential image of $U$, then $S_iX\in\essim U_i$ for all $i\in I$ by assumption. For the converse, we observe that $\eta_X\colon X\to UFX$ is inverted by $F$ for all $X\in\mathscr C$, so $F_iS_i\eta_X$ is an equivalence for all $i\in I$ by the above. But $S_iUFX\simeq U_iT_iFX\in\essim U_i$ while $S_iX\in\essim U_i$ by assumption; thus, both source and target of $S_i\eta_X$ lie in $\essim U_i$, so $S_i\eta_X$ is conjugate to the equivalence $U_iF_iS_i\eta_X$, hence itself an equivalence. Letting $i$ vary and using joint conservativity again, we therefore conclude that $\eta_X$ is an equivalence, hence $X\in\essim U$. This completes the proof of the lemma.
\end{proof}
\end{lemma}

We can now use the lemma to finally prove the $G$-global Delooping Theorem. It turns out to be convenient to prove this in parallel with the following result:

\begin{thm}\label{thm:G-global-truncation}\index{G-globally connective@$G$-globally connective!truncation|textbf}
The inclusion $\iota\colon(\cat{$\bm G$-Spectra}^{\ge0}_{\textup{$G$-gl.}})^\infty\hookrightarrow\cat{$\bm G$-Spectra}_{\textup{$G$-gl.}}^\infty$ of the $G$-globally connective $G$-spectra admits a right adjoint $\tau$. Moreover, for any homomorphism $\phi$ the canonical mate ${\und_\phi^\infty}\tau\Rightarrow\tau\und_\phi^\infty$ of the natural equivalence
\begin{equation*}
\begin{tikzcd}[column sep=small]
(\cat{$\bm G$-Spectra}_{\textup{$G$-global}}^{\ge0})^\infty\arrow[r, hook]\arrow[d, "\und_\phi^\infty"'] & \cat{$\bm G$-Spectra}_{\textup{$G$-global}}^\infty\arrow[d, "\und_\phi^\infty"]\\
(\cat{$\bm H$-Spectra}_{\textup{$H$-equivariant}}^{\ge0})^\infty\arrow[r, hook]\twocell[ur] & \cat{$\bm H$-Spectra}^\infty_{\textup{$H$-equivariant}}
\end{tikzcd}
\end{equation*}
induced by the identity is an equivalence.
\end{thm}

\begin{proof}[Proof of Theorems~\ref{thm:group-completion} and~\ref{thm:G-global-truncation}]
We will break up the argument into several steps.

\textit{Step 1.} We prove a version of Theorem~\ref{thm:group-completion} for connective $G$-global/$H$-equivariant spectra. For this we let $I$ be the class of all homomorphisms $\phi\colon H\to G$ from finite groups to $G$. We want to apply Lemma~\ref{lemma:abstract-mate} to the adjunction $(\ref{eq:g-global-group-completion-adjunction})$ together with the family of the Bousfield localizations
\begin{equation}\label{eq:g-global-group-completion-adjunction}
\mathcal E_H^\infty\colon(\cat{$\bm\Gamma$-$\bm H$-SSet}_*)_{\textup{$H$-equivariant}}^\infty\rightleftarrows(\cat{$\bm H$-Spectra}^{\ge0}_{\textup{$H$-equivariant}})^\infty :\!\textbf{\textup R}\Phi_H
\end{equation}
from Theorem~\ref{thm:equivariant-delooping} for all $\phi\colon H\to G$ and the jointly conservative families
\begin{align*}
\big(\mathord{\und_\phi^\infty}\colon(\cat{$\bm\Gamma$-$\bm G$-$\bm{\mathcal I}$-SSet}_*)_{\textup{$G$-global level}}^\infty&\to(\cat{$\bm\Gamma$-$\bm H$-SSet}_*)_{\textup{$H$-equivariant level}}^\infty\big)_{\phi\in I}\\
\big(\mathord{\und_\phi^\infty}\colon(\cat{$\bm G$-Spectra}_{\textup{$G$-global}}^{\ge0})^\infty&\to(\cat{$\bm H$-Spectra}^{\ge0}_{\textup{$H$-equivariant}})^\infty\big)_{\phi\in I}.
\end{align*}
The equivalences $(\ref{diag:comparison-base})$ are provided by Lemma~\ref{lemma:und-E-vs-E-H}, whereas Corollary~\ref{cor:R-otimes-very-special} together with Theorem~\ref{thm:equivariant-delooping} shows that $\smash{\essim\und_\phi^\infty\cat{R}\Phi^{\otimes}\subset \essim\cat{R}\Phi_H}$. Finally, Proposition~\ref{prop:group-completion-counit} verifies Condition~$(1)$ of the lemma. Together with another application of Theorem~\ref{thm:equivariant-delooping}, we may therefore conclude that the essential image of $\textbf{R}\Phi^\otimes$ consists precisely of the very special $G$-global $\Gamma$-spaces, that a map is inverted by $(\mathcal E^\otimes)^\infty$ if and only if it sent under $\und_\phi$ to a stable weak equivalence for all $\phi$, and that the canonical mates
\begin{equation}\label{eq:connective-mate}
\begin{tikzcd}
(\cat{$\bm\Gamma$-$\bm G$-$\bm{\mathcal I}$-SSet}_*)_{\textup{$G$-global level}}^\infty\arrow[d, "\und_\phi^\infty"']\twocell[dr] &[1em]\arrow[l, "\textbf{R}\Phi^\otimes"'] (\cat{$\bm G$-Spectra}_{\textup{$G$-global}}^{\ge0})^\infty\arrow[d, "\und_\phi^\infty"]\\
(\cat{$\bm\Gamma$-$\bm H$-SSet}_*)_{\textup{$H$-equivariant level}}^\infty&\arrow[l, "\textbf{R}\Phi_H"] (\cat{$\bm H$-Spectra}_{\textup{$H$-equivariant}}^{\ge0})^\infty
\end{tikzcd}
\end{equation}
of the given equivalences are again equivalences. Note that this almost completes the proof of Theorem~\ref{thm:group-completion}, except that we want the equivalences $(\ref{eq:connective-mate})$ without the connectivity requirement.

\textit{Step 2.} We will now prove Theorem~\ref{thm:G-global-truncation}.
In order to avoid confusion, let us momentarily write $(e^\otimes)^\infty$ for $(\mathcal E^\otimes)^\infty$ viewed as a functor to $(\cat{$\bm G$-Spectra}^{\ge0}_{\textup{$G$-global}})^\infty$, i.e.~$(\mathcal E^\otimes)^\infty=\iota\circ (e^\otimes)^\infty$. It is then a purely formal calculation (using the results of the previous step) that $\tau\mathrel{:=}(e^\otimes)^\infty\circ\textbf{R}\Phi^\otimes$ is the desired right adjoint; the counit can be taken to be the counit of $(\mathcal E^\otimes)^\infty\dashv\textbf{R}\Phi^\otimes$ (which is indeed a natural transformation from $(\mathcal E^\otimes)^\infty\circ\textbf{R}\Phi^\otimes=\iota\circ(e^\otimes)^\infty\circ\textbf{R}\Phi^\otimes=\iota\circ\tau$ to the identity). The rest of the claim then follows by an easy application of the opposite of Lemma~\ref{lemma:abstract-mate}.

\textit{Step 3.} We can now finish the proof of Theorem~\ref{thm:group-completion}: we have to show that the canonical mate of the equivalence
\begin{equation}\label{diag:mate-combined}
\begin{tikzcd}
(\cat{$\bm\Gamma$-$\bm G$-$\bm{\mathcal I}$-SSet}_*)_{\textup{$G$-global level}}^\infty\arrow[r, "(\mathcal E^\otimes)^\infty"]\arrow[d, "\und_\phi^\infty"'] & \cat{$\bm G$-Spectra}_{\textup{$G$-global}}^\infty\arrow[d, "\und_\phi^\infty"]\\
(\cat{$\bm\Gamma$-$\bm H$-SSet}_*)_{\textup{$H$-equivariant level}}^\infty\arrow[r, "\mathcal E_H^\infty"']\twocell[ur]& \cat{$\bm H$-Spectra}_{\textup{$H$-equivariant}}^\infty
\end{tikzcd}
\end{equation}
considered above is itself an equivalence. However, $(\ref{diag:mate-combined})$ agrees up to canonical higher homotopy with the pasting
\begin{equation*}
\begin{tikzcd}[cramped]
(\cat{$\bm\Gamma$-$\bm G$-$\bm{\mathcal I}$-SSet}_*)_{\textup{$G$-gl.~level}}^\infty\arrow[r, "(\mathcal E^\otimes)^\infty"]\arrow[d, "\und_\phi^\infty"'] & (\cat{$\bm G$-Spectra}^{\ge0}_{\textup{$G$-global}})^\infty\arrow[d, "\und_\phi^\infty"]\arrow[r, hook] &[-1.5em] \cat{$\bm G$-Spectra}_{\textup{$G$-global}}^\infty\arrow[d, "\und_\phi^\infty"]\\
(\cat{$\bm\Gamma$-$\bm H$-SSet}_*)_{\textup{$H$-equiv.~level}}^\infty\arrow[r, "\mathcal E_H^\infty"']\twocell[ur] & \cat{$\bm H$-Spectra}_{\textup{$H$-equiv.}}^\infty\arrow[r, hook]\twocell[ur]& (\cat{$\bm H$-Spectra}^{\ge0}_{\textup{$H$-equiv.}})^\infty
\end{tikzcd}
\end{equation*}
where the left hand square is filled with the restriction of the previous equivalence while the right hand equivalence is induced by the identity. Thus, the canonical mate of $(\ref{diag:mate-combined})$ is equivalent to the pasting of the individual mates of these two squares, and both of these were seen to be equivalences above.
\end{proof}

\begin{cor}\label{cor:E-otimes-special}
The homotopical functor $\mathcal E^\otimes$ induces a (Bousfield) localization
\begin{equation}\label{eq:delooping-restricted}
(\cat{$\bm\Gamma$-$\bm G$-$\bm{\mathcal I}$-SSet}_*^{\textup{special}})_{\textup{$G$-global level}}^\infty\to(\cat{$\bm G$-Spectra}_{\textup{$G$-global}}^{\ge0})^\infty.
\end{equation}
\begin{proof}
By the above theorem, the fully faithful right adjoint $\cat{R}\Phi^\otimes$ of
\begin{equation*}
(\mathcal E^\otimes)^\infty\colon(\cat{$\bm\Gamma$-$\bm G$-$\bm{\mathcal I}$-SSet}_*)^\infty_{\textup{$G$-global level}}\to(\cat{$\bm G$-Spectra}_{\textup{$G$-global}}^{\ge0})^\infty
\end{equation*}
factors through the (very) special $G$-global $\Gamma$-spaces, so it restricts to a fully faithful right adjoint to $(\ref{eq:delooping-restricted})$ by Proposition~\ref{prop:localization-subcategory}.
\end{proof}
\end{cor}

\chapter[$G$-global algebraic $K$-theory]{\texorpdfstring{$\except{toc}{\bm G}\for{toc}{G}$}{G}-global algebraic \texorpdfstring{$\except{toc}{\bm K}\for{toc}{K}$}{K}-theory}
In this final chapter, we construct \emph{$G$-global algebraic $K$-theory}, generalizing Schwede's global algebraic $K$-theory \cite{schwede-k-theory} and refining classical $G$-equivariant algebraic $K$-theory. We then employ the theory of $G$-globally coherently commutative monoids developed in the previous chapters to prove $G$-global refinements of the Barratt-Priddy-Quillen Theorem and of Thomason's theorem \cite[Theorem~5.1]{thomason} that symmetric monoidal categories model all connective stable homotopy types.

\section{Definition and basic properties}
We will give two equivalent constructions of $G$-global algebraic $K$-theory: one as a straightforward generalization of Schwede's global algebraic $K$-theory, and the other one by adapting Shimakawa's construction of $G$-equivariant algebraic $K$-theory.

\subsection{$\bm G$-global $\bm K$-theory of $\bm G$-parsummable categories}\label{subsec:global-k} Let us begin by recalling the basic setup \cite[Definitions~2.11, 2.12, and~2.32]{schwede-k-theory} of Schwede's approach to global algebraic $K$-theory:

\begin{defi}
A small $E\mathcal M$-category $\mathcal C$ (i.e.~a category with a strict action of the categorical monoid $E\mathcal M$) is \emph{tame} if its underlying $\mathcal M$-set $\Ob(\mathcal C)$ of objects is tame in the sense of Definition~\ref{defi:M-set-tame-support}. The \emph{support} $\supp(c)$ of an object $c\in\mathcal C$ is its support as an element of the $\mathcal M$-set $\Ob(\mathcal C)$.\index{support!for EM-actions@for $E\mathcal M$-actions|textbf}
\end{defi}

If $\mathcal C$ is an $E\mathcal M$-category, then we denote the action of $u\in\mathcal M$ by $u_*\colon\mathcal C\to\mathcal C$ and the action of the morphism $(v,u)$ by $[v,u]\colon u_*\Rightarrow v_*$.

\begin{defi}
The \emph{box product}\index{box product!on EM-Cat-tau@on $\cat{$\bm{E\mathcal M}$-Cat}^\tau$|textbf} $\mathcal C\boxtimes\mathcal D$ of two tame $E\mathcal M$-categories is the full subcategory of $\mathcal C\times\mathcal D$ spanned by all pairs $(c,d)$ such that $\supp(c)\cap\supp(d)=\varnothing$.
\end{defi}

We write $\cat{$\bm{E\mathcal M}$-Cat}^\tau$ for the category of tame $E\mathcal M$-categories and strictly $E\mathcal M$-equivariant functors. The box product then becomes a subfunctor of the cartesian product on $\cat{$\bm{E\mathcal M}$-Cat}^\tau$, and \cite[Proposition~2.35]{schwede-k-theory} shows that the structure isomorphisms of the cartesian symmetric monoidal structure restrict to make $\boxtimes$ into the tensor product of a preferred symmetric monoidal structure on $\cat{$\bm{E\mathcal M}$-Cat}^\tau$.

We can now define \cite[Definition~4.1]{schwede-k-theory}:

\begin{defi}
We write $\cat{ParSumCat}\mathrel{:=}\CMon(\cat{$\bm{E\mathcal M}$-Cat}^\tau)$ and call its objects \emph{parsummable categories}.\index{G-parsummable category@$G$-parsummable category|textbf}\nomenclature[aParSumCat]{$\cat{ParSumCat}$}{category of parsummable categories}
\end{defi}

In the next subsection, we will see how any small symmetric monoidal category gives rise to a parsummable category. For now, let us recall Schwede's construction \cite[Construction~10.1]{schwede-k-theory} of the parsummable category $\mathcal P(R)$ associated to a (not necessarily commutative) ring $R$:

\begin{ex}\label{ex:parsummable-module-category}\nomenclature[aPR]{$\mathcal P(R)$}{parsummable category associated to a ring $R$}
We write $R^{(\omega)}$ for the left $R$-module of functions $\omega\to R$ vanishing almost everywhere; this is free with basis given by the characteristic functions $e_0,e_1,\dots,$ with $e_i(j)=1$ for $i=j$ and $e_i(j)=0$ otherwise.

An object of $\mathcal P(R)$ is a finitely generated $R$-submodule $M\subset R^{(\omega)}$ such that the inclusion admits an $R$-linear retraction (which is not part of the data); note that any such $M$ is in particular projective. The morphisms $M\to N$ in $\mathcal P(R)$ are the (abstract) $R$-linear isomorphisms, and the composition in $\mathcal P(R)$ is the evident one. In particular, we have a fully faithful embedding $\mathcal P(R)\hookrightarrow\mathscr P(R)$ into the usual category of finitely generated projective left $R$-modules and $R$-linear isomorphisms; in fact, it is not hard to show that this functor is even an equivalence of categories, see \cite[proof of Theorem~10.3-(i)]{schwede-k-theory} or cf.~the argument in Example~\ref{ex:module-saturated} (specialized to $H=G=1$) below.\nomenclature[aPR]{$\mathscr P(R)$}{symmetric monoidal category of finitely generated projective $R$-modules and $R$-linear isomorphisms}

Any $u\in\mathcal M$ defines an $R$-linear embedding $u_!\colon R^{(\omega)}\to R^{(\omega)}$ via $u_!(e_i)=e_{u(i)}$. We define the $u$-action on objects via $u_*M=u_!(M)\subset R^{(\omega)}$, and we set $[u,1]_{M}=u_!|_M\colon M\to u_!(M)$. These data in fact completely determine the $E\mathcal M$-action \cite[Proposition~2.6]{schwede-k-theory}, and we omit the easy verification that this is well-defined and tame.

We define the sum of two disjointly supported $M,N\in\mathcal P(R)$ as their internal sum as submodules of $R^{(\omega)}$. One easily shows that this sum is actually direct, which allows us to define the sum of two morphisms $f\colon M\to M'$, $g\colon N\to N'$ with disjointly supported sources and disjointly supported targets via $(f+g)(m+n)=f(m)+g(n)$ for all $m\in M$, $n\in N$. Finally, the additive unit is given by the trivial submodule $0\subset R^{(\omega)}$.
\end{ex}

The construction of the global algebraic $K$-theory of a parsummable category can be broken up into two steps: first, we construct a `global $\Gamma$-category,' which we will then later deloop to a global spectrum.

\begin{constr}\index{global Gamma-category@global $\Gamma$-category|seeonly{$G$-global $\Gamma$-category}}\index{global Gamma-category@global $\Gamma$-category}
Let us write $\cat{$\bm\Gamma$-$\bm{E\mathcal M}$-Cat}_*$ for the category of functors $X\colon\Gamma\to\cat{$\bm{E\mathcal M}$-Cat}$ with $X(0^+)=*$.

Analogously to Construction~\ref{constr:digamma}, we can now define $\digamma\colon\cat{ParSumCat}\to\cat{$\bm\Gamma$-$\bm{E\mathcal M}$-Cat}_*$ with $\digamma(\mathcal C)(n^+)=\mathcal C^{\boxtimes n}$ and the evident functoriality in $\mathcal C$; the structure maps of $\digamma(\mathcal C)$ are again given by the universal property of $\boxtimes$ as coproduct in $\cat{ParSumCat}$, also see~\cite[Construction~4.3]{schwede-k-theory}.
\end{constr}

We obtain a functor $\cat{$\bm\Gamma$-$\bm{E\mathcal M}$-Cat}_*\to\cat{$\bm\Gamma$-$\bm{E\mathcal M}$-SSet}_*$ by applying the nerve levelwise. From there, one can then use the following construction from \cite[Construction~3.3]{schwede-k-theory} to pass to symmetric spectra:

\begin{constr}
Let $X\in\cat{$\bm\Gamma$-$\bm{E\mathcal M}$-SSet}_*$. The \emph{associated (symmetric) spectrum}\index{associated spectrum|textbf} $X\langle\mathbb S\rangle$ is given by
\begin{equation*}
X\langle\mathbb S\rangle(A)= X[\omega^A](S^A)
\end{equation*}
(see Construction~\ref{constr:omega-bullet}). If $i\colon A\to B$ is an injection, then the structure map $S^{B\setminus i(A)}\smashp X\langle\mathbb S\rangle(A)\to X\langle\mathbb S\rangle(B)$ is the composition
\begin{equation*}
S^{B\setminus i(A)}\smashp X[\omega^A](S^A)\hskip0pt minus 2pt\xrightarrow{\!\textup{asm}}\hskip0pt minus 2ptX[\omega^A](S^{B\setminus i(A)}\smashp S^A)\hskip0pt minus 2pt\cong\hskip0pt minus 2ptX[\omega^A](S^B)\hskip0pt minus 2pt\xrightarrow{X[i_!](S^B)}\hskip0pt minus 2pt X[\omega^B](S^B).
\end{equation*}
This becomes a functor in $X$ in the obvious way.
\end{constr}

Put differently (see Remark~\ref{rk:E-otimes-levelwise-description}), $(\blank)\langle\mathbb S\rangle$ agrees with the composition
\begin{equation*}
\cat{$\bm\Gamma$-$\bm{E\mathcal M}$-SSet}_*\xrightarrow{(\blank)[\omega^\bullet]}\cat{$\bm\Gamma$-$\bm{\mathcal I}$-SSet}_*\xrightarrow{\mathcal E^\otimes} \cat{Spectra}.
\end{equation*}

\begin{rk}
To be entirely precise, Schwede instead works with the category $\cat{Spectra}_{\cat{Top}}$ of symmetric spectra in \emph{topological spaces}, and he first passes to $\cat{$\bm\Gamma$-$\bm{|E\mathcal M|}$-Top}_*$ via geometric realization. We omit the routine verification that after postcomposing with $|\blank|\colon\cat{Spectra}\to\cat{Spectra}_{\cat{Top}}$ the above agrees with his construction up to isomorphism; the only non-trivial ingredient is the comparison between prolongations of $\Gamma$-spaces in the simplicial and in the topological world, which appears for example as \cite[Proposition~B.29]{schwede-book}.
\end{rk}

We are now ready to define, see~\cite[Definition~4.14]{schwede-k-theory}:

\begin{defi}
We write $\cat{K}_{\textup{gl}}$\nomenclature[aKgl]{$\cat{K}_{\text{gl}}$}{global algebraic $K$-theory}\index{global algebraic K-theory@global algebraic $K$-theory|seealso{$G$-global algebraic $K$-theory}}\index{algebraic K-theory@algebraic $K$-theory!global|seeonly{global algebraic $K$-theory}}\index{global algebraic K-theory@global algebraic $K$-theory!of parsummable categories|textbf} for the composition
\begin{equation*}
\cat{ParSumCat}\xrightarrow{\digamma}\cat{$\bm\Gamma$-$\bm{E\mathcal M}$-Cat}_*\xrightarrow{\nerve}\cat{$\bm\Gamma$-$\bm{E\mathcal M}$-SSet}_*\xrightarrow{(\blank)\langle\mathbb S\rangle}\cat{Spectra}.
\end{equation*}
For any parsummable category $\mathcal C$ we call $\cat{K}_{\textup{gl}}(\mathcal C)$ the \emph{global algebraic $K$-theory} of $\mathcal C$.
\end{defi}

\begin{ex}\label{ex:global-K-rings}\index{global algebraic K-theory@global algebraic $K$-theory!of rings|textbf}
Let $R$ be any ring. Applying the above to the parsummable category $\mathcal P(R)$ from Example~\ref{ex:parsummable-module-category} yields Schwede's definition \cite[Definition~10.2]{schwede-k-theory} of the global algebraic $K$-theory $\cat{K}_{\textup{gl}}(R)$ of $R$. Its underlying non-equivariant spectrum recovers the usual algebraic $K$-theory of $R$, see~\cite[Theorem~10.3-(ii)]{schwede-k-theory}; we will more generally identify its underlying $G$-spectrum for any finite $G$ in Remark~\ref{rk:global-K-rings-vs-G-equivariant} below.
\end{ex}

The correct notion of equivariant algebraic $K$-theory of small symmetric mon\-oidal categories with $G$-action is \emph{not} given by simply applying the non-equivariant Segal-May-Shimada-Shimakawa construction and pulling through the $G$-action; instead, this required Shimakawa's insight explained in Example~\ref{ex:shimada-shimakawa-equivariant} above. Interestingly, it turns out that this is not necessary when generalizing from global to $G$-global algebraic $K$-theory:

\begin{defi}\index{algebraic K-theory@algebraic $K$-theory!G-global@$G$-global|seeonly{$G$-global algebraic $K$-theory}}\index{G-global algebraic K-theory@$G$-global algebraic $K$-theory!of G-parsummable categories@of $G$-parsummable categories|textbf}\index{G-parsummable category@$G$-parsummable category!$G$-global algebraic K-theory@$G$-global algebraic $K$-theory|seeonly{$G$-global algebraic $K$-theory, of $G$-parsummable categories}}
Let $\mathcal C$ be a $G$-parsummable category (i.e.~a $G$-object in $\cat{ParSumCat}$). Its \emph{$G$-global algebraic $K$-theory} $\cat{K}_{\textup{$G$-gl}}(\mathcal C)$\nomenclature[aKGgl]{$\cat{K}_{\textup{$G$-gl}}$}{$G$-global algebraic $K$-theory} is the $G$-global spectrum obtained by equipping $\cat{K}_{\textup{gl}}(\mathcal C)$ with the $G$-action induced by functoriality. We write $\cat{K}_{\textup{$G$-gl}}\colon\cat{$\bm G$-ParSumCat}\to\cat{$\bm G$-Spectra}$ for the resulting functor.
\end{defi}

We will explain in the next subsection how this in particular yields the $G$-global algebraic $K$-theory of a small symmetric monoidal category with $G$-action, and prove later in Theorem~\ref{thm:G-global-K-vs-G-equiv} that this indeed refines usual $G$-equivariant algebraic $K$-theory. Moreover, we can now use this to introduce the $G$-global algebraic $K$-theory of $G$-rings:

\begin{ex}\label{ex:G-global-K-rings}\index{G-global algebraic K-theory@$G$-global algebraic $K$-theory!of G-rings@of $G$-rings|textbf}
Schwede's construction of $\mathcal P$ recalled above is actually strictly functorial in ring homomorphisms $\alpha\colon R\to S$ as follows:

We write $\alpha_{\lozenge}$ for the unique $S$-linear map $S\otimes_R R^{(\omega)}\to S^{(\omega)}$ sending $s\otimes e_i$ to $s\cdot e_i$ for all $i\in\omega, s\in S$; this map is clearly an isomorphism. Then $\mathcal P(\alpha)$ is given on objects by $\mathcal P(\alpha)(M)=\alpha_{\lozenge}(S\otimes_R M)$ where we identify $S\otimes_R M$ with its image in $S\otimes_R R^{(\omega)}$, which secretly uses the \emph{existence} of a retraction to $M\hookrightarrow R^{(\omega)}$. Similarly, if $f\colon M\to N$ is an $R$-linear isomorphism, then $\mathcal P(\alpha)(f)$ is obtained from $S\otimes_R f$ by conjugating with the restrictions of $\alpha_{\lozenge}$. We omit the straightforward verification that this is well-defined and functorial.

In particular, if $R$ is a ring equipped with a $G$-action through ring homomorphisms, then $\mathcal P(R)$ inherits a $G$-action. Explicitly, we have for each $g\in G$ a map $g.\blank\colon R^{(\omega)}\to R^{(\omega)}$ sending $\sum_{i\in\omega} r_i\cdot e_i$ to $\sum_{i\in\omega} (g.r_i)\cdot e_i$, and this is \emph{semilinear} in the sense that $g.(r\cdot x)=(g.r)\cdot (g.x)$ for all $x\in R^{(\omega)}$. The $G$-action on objects is then defined via $g.M=(g.\blank)(M)$, and if $f\colon M\to N$ is a morphism, then $g.f=(g.\blank)\circ f\circ (g^{-1}.\blank)\colon g.M\to g.N$.

Thus, the above construction yields for any such $G$-ring $R$ a $G$-global algebraic $K$-theory spectrum $\cat{K}_{\textup{$G$-gl}}(R)$, generalizing Schwede's global algebraic $K$-theory. In particular, if $E/F$ is any Galois extension (of fields, or more generally of rings), then we obtain a $\Gal(E/F)$-global algebraic $K$-theory spectrum $\cat{K}_{\textup{$\Gal(E/F)$-gl}}(E)$ (viewing $\Gal(E/F)$ as a discrete group only, even if $E/F$ is infinite).
\end{ex}

We close this discussion by establishing a basic invariance property of $G$-global $K$-theory:

\begin{defi}\index{G-global weak equivalence@$G$-global weak equivalence!in G-ParSumCat@in $\cat{$\bm G$-ParSumCat}$|textbf}\index{G-parsummable category@$G$-parsummable category!G-global weak equivalence@$G$-global weak equivalence|seeonly{$G$-global weak equivalence, in $\cat{$\bm G$-ParSumCat}$}}
A map $f\colon\mathcal C\to\mathcal D$ of $G$-parsummable categories is called a \emph{$G$-global weak equivalence} if the induced functor $\mathcal C^\phi\to\mathcal D^\phi$ is a weak homotopy equivalence (i.e.~weak equivalence on nerves) for each universal subgroup $H\subset\mathcal M$ and each homomorphism $\phi\colon H\to G$.
\end{defi}

For $G=1$, Schwede \cite[Definition~2.26]{schwede-k-theory} considered these under the name `global equivalence.' We can now prove the following generalization of \cite[Theorem~4.16]{schwede-k-theory}:

\begin{prop}\label{prop:K-G-gl-homotopical}
The functor $\cat{K}_{\textup{$G$-gl}}$ preserves $G$-global weak equivalences.
\end{prop}

For the proof it will be convenient to slightly reformulate the above construction of $\cat{K}_{\textup{$G$-gl}}$. This will use:

\begin{lemma}
The nerve of a tame $E\mathcal M$-category is tame again, and the induced functor $\nerve\colon\cat{$\bm{E\mathcal M}$-$\bm G$-Cat}^\tau\to\cat{$\bm{E\mathcal M}$-$\bm G$-SSet}^\tau$ admits a preferred strong symmetric monoidal structure.
\begin{proof}
See~\cite[Example~2.7 and Proposition~2.20]{sym-mon-global}.
\end{proof}
\end{lemma}

In particular, we get a canonical lift of $\nerve$ to a functor $\cat{$\bm G$-ParSumCat}\to\cat{$\bm G$-ParSumSSet}$ that we again denote by $\nerve$.

\begin{cor}\label{cor:compatibility-digammas}
The diagram
\begin{equation*}
\begin{tikzcd}
\cat{$\bm G$-ParSumCat}\arrow[d, "\nerve"']\arrow[rr, "\digamma"] && \cat{$\bm\Gamma$-$\bm{E\mathcal M}$-$\bm G$-Cat}_*\arrow[d, "\nerve"]\\
\cat{$\bm G$-ParSumSSet}\arrow[r, "\digamma"'] & \cat{$\bm\Gamma$-$\bm{E\mathcal M}$-$\bm G$-SSet}_*^\tau\arrow[r, hook] & \cat{$\bm\Gamma$-$\bm{E\mathcal M}$-$\bm G$-SSet}_*
\end{tikzcd}
\end{equation*}
commutes up to canonical natural isomorphism.
\begin{proof}
The top arrow obviously factors through $\cat{$\bm\Gamma$-$\bm{E\mathcal M}$-$\bm G$-Cat}_*^\tau$, so it suffices to construct a natural isomorphism filling
\begin{equation*}
\begin{tikzcd}
\cat{$\bm G$-ParSumCat}\arrow[d, "\nerve"']\arrow[r, "\digamma"] & \cat{$\bm\Gamma$-$\bm{E\mathcal M}$-$\bm G$-Cat}_*^\tau\arrow[d, "\nerve"]\\
\cat{$\bm G$-ParSumSSet}\arrow[r, "\digamma"'] & \cat{$\bm\Gamma$-$\bm{E\mathcal M}$-$\bm G$-SSet}_*^\tau.
\end{tikzcd}
\end{equation*}
But as in the proof of Theorem~\ref{thm:gamma-vs-uc}, it follows by abstract nonsense that the structure isomorphisms of the strong symmetric monoidal functor $\nerve\colon\cat{$\bm{E\mathcal M}$-$\bm G$-Cat}^\tau\to\cat{$\bm{E\mathcal M}$-$\bm G$-SSet}^\tau$ assemble into the desired isomorphism.
\end{proof}
\end{cor}

\begin{proof}[Proof of Proposition~\ref{prop:K-G-gl-homotopical}]
By Corollary~\ref{cor:compatibility-digammas} we can factor $\cat{K}_{\textup{$G$-gl}}$ up to isomorphism as
\begin{align*}
\cat{$\bm G$-ParSumCat}&\xrightarrow{\nerve}\cat{$\bm G$-ParSumSSet}\xrightarrow{\digamma}\cat{$\bm\Gamma$-$\bm{E\mathcal M}$-$\bm G$-SSet}_*^\tau\hookrightarrow\cat{$\bm\Gamma$-$\bm{E\mathcal M}$-$\bm G$-SSet}_*\\
&\xrightarrow{(\blank)[\omega^\bullet]}\cat{$\bm\Gamma$-$\bm G$-$\bm{\mathcal I}$-SSet}_*\xrightarrow{\mathcal E^\otimes}\cat{$\bm G$-Spectra}.
\end{align*}
Of these the first arrow preserves $G$-global weak equivalences as $\nerve$ preserves limits; moreover, the second arrow sends these to $G$-global level weak equivalences by Proposition~\ref{prop:digamma-quillen-adjunction}, which are preserved by the inclusion by definition and by the penultimate arrow by Theorem~\ref{thm:gamma-ev-omega}. Finally, $\mathcal E^\otimes$ is homotopical by Corollary~\ref{cor:E-otimes-homotopical}.
\end{proof}

$G$-global algebraic $K$-theory is compatible with restrictions along group homomorphisms in the following sense:

\begin{prop}\label{prop:K-G-gl-parsum-functoriality}
Let $\phi\colon H\to G$ be any group homomorphism. Then the diagram
\begin{equation*}
\begin{tikzcd}
\cat{$\bm G$-ParSumCat}^\infty\arrow[d, "(\phi^*)^\infty"']\arrow[r, "\cat{K}_{\textup{$G$-gl}}^\infty"]&[1em]\cat{$\bm G$-Spectra}^\infty_{\textup{$G$-global}}\arrow[d,"(\phi^*)^\infty"]\\
\cat{$\bm H$-ParSumCat}^\infty\arrow[r, "\cat{K}_{\textup{$H$-gl}}^\infty"']&\cat{$\bm H$-Spectra}_{\textup{$H$-global}}^\infty
\end{tikzcd}
\end{equation*}
commutes up to canonical equivalence.
\begin{proof}
By construction, we even have an equality $\phi^*\circ\cat{K}_{\textup{$G$-gl}}=\cat{K}_{\textup{$H$-gl}}\circ\phi^*$ of homotopical functors.
\end{proof}
\end{prop}

\subsection{$\bm G$-global $\bm K$-theory of symmetric monoidal $\bm G$-categories} Recall that a \emph{permutative category}\index{permutative G-category@permutative $G$-category}\index{symmetric monoidal category|seeonly{permutative $G$-category}} is a symmetric monoidal category in which the unitality and associativity isomorphisms are required to be the respective identities. We write $\cat{PermCat}$\nomenclature[aPermCat]{$\cat{PermCat}$}{category of small permutative categories and strict symmetric monoidal functors} for the $1$-category of small permutative categories and \emph{strict} symmetric monoidal functors. In \cite[Constructions~11.1 and~11.6]{schwede-k-theory}, Schwede constructs a specific functor $\Phi\colon\cat{PermCat}\to\cat{ParSumCat}$\nomenclature[aPhi]{$\Phi$}{parsummable category associated to a small permutative category} as a variation of a `strictification' construction due to Schlichtkrull and Solberg \cite[Sections~4.14 and~7]{solberg-schlichtkrull}:

\begin{constr}
While we will be able to completely black box the definition of $\Phi$, let us give the basic idea of its construction for motivational purposes; for details, we refer the reader to Schwede's article.

Let $\mathscr C$ be a small permutative category. An object of $\Phi\mathscr C$ is a family $(X_i)_{i\in\omega}$ of objects of $\mathscr C$ such that $X_i=\bm1$ for all but finitely many $i\in\omega$. If $Y_\bullet$ is another object of $\Phi\mathscr C$, then
\begin{equation*}
\Hom_{\Phi\mathscr C}(X_\bullet,Y_\bullet)\mathrel{:=}\Hom_{\mathscr C}\left(\bigotimes_{i\in\omega} X_i,\bigotimes_{i\in\omega} Y_i\right);
\end{equation*}
note that these infinite tensor products are indeed well-defined as the tensor product of $\mathscr C$ is strictly unital. The composition in $\Phi\mathscr C$ is given by the composition in $\mathscr C$.

We have an action of $\mathcal M$ on $\Ob(\Phi\mathscr C)$ via `shuffling and extending by $\bm1$,' i.e.
\begin{equation*}
(u.X_\bullet)_i=\begin{cases}
X_j & \text{if $i=u(j)$}\\
\bm1 & \text{if $i\notin\im u$}
\end{cases}
\end{equation*}
for any $u\in\mathcal M$, $X_\bullet\in\Phi\mathscr C$. For any further $v\in\mathcal M$, the structure isomorphism $[v,u]$ is given by a suitable composition of the symmetry isomorphisms of the permutative structure on $\mathscr C$.

Clearly, $\supp X_\bullet=\{i\in\omega : X_i\not=\bm1\}$, so we can define the sum of two disjointly supported objects $X_\bullet,Y_\bullet$ by
\begin{equation*}
(X_\bullet+Y_\bullet)_i=\begin{cases}
X_i & \text{if $i\in\supp X_\bullet$}\\
Y_i & \text{otherwise}.
\end{cases}
\end{equation*}
The additive unit is the constant family at $\bm1$. The sum of morphisms $f\colon X_\bullet\to X_\bullet'$, $g\colon Y_\bullet\to Y_\bullet'$ with pairwise disjoint sources and pairwise disjoint targets is given by suitably conjugating the tensor product $f\otimes g$ in $\mathscr C$ by symmetry isomorphisms.

Finally, if $F\colon\mathscr C\to\mathscr D$ is a strict symmetric monoidal functor, then $\Phi(F)$ is given by pushforward along $F$.
\end{constr}

Again, we formally extend this to a functor $\cat{$\bm G$-PermCat}\to\cat{$\bm G$-ParSumCat}$ from the category of \emph{permutative $G$-categories}, i.e.~permutative categories with a strict $G$-action through strict symmetric monoidal functors, to the category of \emph{$G$-parsummable categories}.

While the underlying non-equivariant spectrum of $\cat{K}_\text{gl}(\Phi\mathscr C)$ recovers the usual $K$-theory of $\mathscr C$ \cite[Remark~4.11]{schwede-k-theory}, already for $G=1$ the actual \emph{global} homotopy type $\cat{K}_\text{gl}(\Phi\mathscr C)$ is not the `correct' definition of the global algebraic $K$-theory of $\mathscr C$---for example, if $\mathscr C$ is a small permutative replacement of the category $\mathscr P(\mathbb C)$ of finite dimensional $\mathbb C$-vector spaces and $\mathbb C$-linear isomorphisms under $\oplus$, then $\cat{K}_{\textup{gl}}(\Phi\mathscr C)$ does not agree with the global algebraic $K$-theory of $\mathbb C$ introduced in Example~\ref{ex:global-K-rings}, see~\cite[Proposition~11.9]{schwede-k-theory}. In fact, as we explain in \cite[Remark~3.27]{perm-parsum-categorical}, there is no (small) permutative category $\mathscr D$ at all such that $\cat{K}_{\textup{gl}}(\Phi\mathscr D)\simeq\cat{K}_{\textup{gl}}(\mathbb C)$.

\subsubsection{Saturation}
Ultimately, the above issue stems from the fact that while we have good control over the underlying category of $\Phi(\mathscr C)$ (which is equivalent to $\mathscr C$ itself), the equivariant information encoded in the $E\mathcal M$-action is not that natural. On the other hand, there is an interesting class of $G$-parsummable categories for which the categorical information is enough to describe the underlying $G$-global homotopy type:

\begin{defi}\index{G-parsummable category@$G$-parsummable category!saturated|textbf}\index{saturated|seeonly{$G$-parsummable category, saturated}}
An $E\mathcal M$-$G$ category $\mathcal C$ is called \emph{saturated} if the canonical map $\mathcal C^\phi\to\mathcal  C^{\myh\phi}\mathrel{:=}\Fun(EH,\mathcal C)^\phi$ from the honest fixed points to the categorical homotopy fixed points (induced by restricting along $EH\to *$) is an equivalence of categories for each universal subgroup $H\subset\mathcal M$ and each homomorphism $\phi\colon H\to G$; here $H$ acts on $EH$ from the right as usual.
\end{defi}

Here we use the notation $\myh$ instead of the usual $h$ in order to emphasize that while our notion of weak equivalences of $G$-parsummable categories is based on weak homotopy equivalences, the above are homotopy fixed points with respect to the \emph{underlying equivalences of categories}. In particular, $(\blank)^{\myh\phi}$ does not preserve $G$-global weak equivalences, i.e.~`homotopy' fixed points are not homotopy invariant.

\begin{rk}
The canonical map $\mathcal C^\phi\to\mathcal C^{\myh\phi}$ is always fully faithful as a limit of fully faithful functors.
\end{rk}

\begin{lemma}\label{lemma:cat-between-sat}
Let $f\colon\mathcal C\to\mathcal D$ be a map of saturated $E\mathcal M$-$G$-categories that is an equivalence of underlying categories. Then $f^\phi$ is an equivalence of categories for any universal subgroup $H\subset\mathcal M$ and any homomorphism $\phi\colon H\to G$; in particular, $f$ is a $G$-global weak equivalence if $\mathcal C$ and $\mathcal D$ are small.
\begin{proof}
See \cite[Remark~1.29]{sym-mon-global}.
\end{proof}
\end{lemma}

\begin{ex}\label{ex:module-saturated}
Let $R$ be any $G$-ring. Then the $G$-parsummable category $\mathcal P(R)$ of Example~\ref{ex:G-global-K-rings} is saturated. For $G=1$ this appears as \cite[Theorem~10.3-(i)]{schwede-k-theory}; we will now give a similar argument for the general case:

Let $H\subset\mathcal M$ be a universal subgroup, let $\phi\colon H\to G$ be any homomorphism, and let $F\in\mathcal P(R)^{\myh\phi}$ arbitrary; we want to show that $F$ is isomorphic to the image of some $N\in\mathcal P(R)^\phi$. To this end, we set $M\mathrel{:=}F(1)$; a straightforward diagram chase shows that $M$ admits a $\phi^*R$-semilinear $H$-action with $h\in H$ acting via
\begin{equation*}
M=F(1)\xrightarrow{(h,\phi(h)).\blank} (h,\phi(h)).F(1)=F(h)\xrightarrow{F(1,h)} F(1)=M,
\end{equation*}
where the first arrow is given by the action on $R^{(\omega)}$.

We now pick $i_1,i_2,\dots\in\omega$ such that the $H$-orbits $H{i_j}\subset\omega$ are free and pairwise disjoint, and we write $X$ for the $R$-linear span $R\langle e_{i_1},e_{i_2},\dots\rangle\subset R^{(\omega)}$. As $M$ is finitely generated, there exists an $R$-linear surjection $p\colon X\to M$, which then admits a section $s\colon M\to X$ as $M$ is projective. We now define $\tilde s\colon M\to R^{(\omega)}$ by
\begin{equation*}
\tilde s(m)\mathrel{:=}\sum_{h\in H}(h,\phi(h)).s(h^{-1}.m),
\end{equation*}
which one easily checks to be $R$-linear and $H$-equivariant (with $H$ acting on $R^{(\omega)}$ via $\mathcal M$ and $\phi$ as before). Moreover, we let $\tilde p\colon R^{(\omega)}\to M$ be the unique $R$-linear map with $\tilde p(e_i)=p(e_i)$ for $i\in\{i_1,i_2,\dots\}$ and $\tilde p(e_i)=0$ otherwise.

As $(h,\phi(h)).X\cap X=0$ for $h\not=1$, it follows that $\tilde p\tilde s=\id_M$, so that $N\mathrel{:=}\tilde s(M)\subset R^{(\omega)}$ is an object of $\mathcal P(R)$. Unravelling definitions and using $H$-equivariance of $\tilde s$, one then readily verifies that $N$ is $\phi$-fixed and that the maps $\big((h,\phi(h)).\blank\big)\circ\tilde s=\tilde s(h.\blank)$ for varying $h\in H$ define an isomorphism in $\mathcal P(R)^{\myh\phi}$ between the given functor $F$ and the functor constant at $N$.
\end{ex}

Apart from the above example, most $G$-parsummable categories that come up in practice are not saturated; in particular, the $G$-parsummable categories of the form $\Phi\mathscr C$ are usually not saturated, even for $G=1$. Despite this seeming scarcity, saturated $G$-parsummable categories actually exist in abundance: as we prove in \cite[Theorem~5.9]{sym-mon-global}, any $G$-parsummable category is $G$-globally weakly equivalent to a saturated one.

In our present situation, however, we instead want to replace the $G$-parsum\-mable category $\Phi\mathscr C$ by a saturated one without changing the \emph{underlying category}, though we will be happy to change its $G$-global weak homotopy type in the process if necessary. A universal way to achieve this is provided by our \emph{saturation construction}, which first appeared for $G=1$ as \cite[Construction~7.20]{schwede-k-theory} and for general $G$ as \cite[Construction~1.29]{sym-mon-global}:

\begin{constr}\index{saturation|seeonly{$G$-parsummable category, saturation}}\index{G-parsummable category@$G$-parsummable category!saturation|textbf}
Let $\mathcal C$ be a tame $E\mathcal M$-category and equip $\Fun(E\mathcal M,\mathcal C)$ with the diagonal of the $E\mathcal M$-action on $\mathcal C$ and the left $E\mathcal M$-action induced by the \emph{right} $E\mathcal M$-action on $E\mathcal M$ via precomposition. We define $\mathcal C^{\textup{sat}}\mathrel{:=}\Fun(E\mathcal M,\mathcal C)^\tau$\nomenclature[asat]{$(\blank)^{\textup{sat}}$}{saturation of an $E\mathcal M$-$G$-category or parsummable category} and call it the \emph{saturation} of $\mathcal C$. The natural functor $s\colon \mathcal C\to\Fun(E\mathcal M,\mathcal C)$ sending an object $c\in\mathcal  C$ to the constant functor at $c$ restricts to $\mathcal C\to\mathcal  C^{\textup{sat}}$; we omit the easy verification that this is natural with respect to the evident functoriality of $(\blank)^{\textup{sat}}$.

Finally, we lift $(\blank)^{\textup{sat}}$ to an endofunctor of $\cat{$\bm{E\mathcal M}$-$\bm G$-Cat}^\tau$ by pulling through the $G$-action; we observe that $s$ automatically defines a natural transformation $\id_{\cat{$\bm{E\mathcal M}$-$\bm G$-Cat}^\tau}\Rightarrow(\blank)^{\textup{sat}}$.
\end{constr}

\begin{thm}\label{thm:saturation}
The above functor $(\blank)^{\textup{sat}}$ takes values in the full subcategory $\cat{$\bm{E\mathcal M}$-$\bm G$-Cat}^{\tau,s}$ of \emph{saturated} tame $E\mathcal M$-$G$-categories. Moreover:
\begin{enumerate}
\item The map $s\colon\mathcal C\to\mathcal C^{\textup{sat}}$ is an underlying equivalence of categories for any $\mathcal C\in\cat{$\bm{E\mathcal M}$-$\bm G$-Cat}^\tau$.
\item The inclusion $\mathcal C^{\textup{sat}}\hookrightarrow\Fun(E\mathcal M,\mathcal C)$ induces \emph{equivalences of categories} on $\phi$-fixed points for every universal $H\subset\mathcal M$ and each $\phi\colon H\to G$; in particular, it is a $G$-global weak equivalence.
\end{enumerate}
\begin{proof}
See \cite[Theorem~1.30]{sym-mon-global} and its proof.
\end{proof}
\end{thm}

\begin{rk}\label{rk:saturated-alt}\index{G-parsummable category@$G$-parsummable category!saturated}
There is an alternative characterization of saturated $E\mathcal M$-categories that will be useful below: the forgetful functor $\cat{$\bm{E\mathcal M}$-Cat}\to\cat{Cat}$ admits a right adjoint $\Fun(E\mathcal M,\blank)$; we claim that an $E\mathcal M$-$G$-category is saturated if and only if the unit $\mathcal C\to\Fun(E\mathcal M,\forget\mathcal C)$ induces equivalences on $\phi$-fixed points for all universal $H\subset\mathcal M$ and all $\phi\colon H\to G$.

Indeed, the unit is always an equivalence of categories and its target is saturated by the previous theorem. Thus, if also the source is saturated, then $\eta^\phi$ is an equivalence by Lemma~\ref{lemma:cat-between-sat}. Conversely, saturated $E\mathcal M$-$G$-categories are clearly closed under functors inducing \emph{equivalences} on all $\phi$-fixed points for $\phi$ as above.
\end{rk}

As explained in \cite[Construction~7.25]{schwede-k-theory}, the usual lax symmetric monoidal structure on $(\blank)^{\textup{sat}}$ with respect to the cartesian product restricts to a lax symmetric monoidal structure with respect to the box product. In particular, $(\blank)^{\textup{sat}}$ canonically lifts to an endofunctor of $\cat{$\bm G$-ParSumCat}$. Using this we can now finally introduce:

\begin{defi}
We define the \emph{$G$-global algebraic $K$-theory}\index{G-global algebraic K-theory@$G$-global algebraic $K$-theory!of permutative G-categories@of permutative $G$-categories|textbf} of permutative $G$-categories as the composition
\begin{equation*}
\cat{$\bm G$-PermCat}\xrightarrow{\Phi}\cat{$\bm G$-ParSumCat}\xrightarrow{(\blank)^{\textup{sat}}}\cat{$\bm G$-ParSumCat}\xrightarrow{\cat{K}_{\textup{$G$-gl}}}\cat{$\bm G$-Spectra},
\end{equation*}
which we denote by $\cat{K}_{\textup{$G$-gl}}$ again.
\end{defi}

\begin{rk}
For $G=1$ the above agrees with Schwede's global definition implicit in \cite[discussion after Proposition~11.9]{schwede-k-theory}.
\end{rk}

It is not hard to show that the above composition sends maps in $\cat{$\bm G$-PermCat}$ that are equivalences of underlying categories to $G$-global weak equivalences. We will however be interested in the following coarser notion of weak equivalence:

\begin{defi}\index{G-global weak equivalence@$G$-global weak equivalence!in G-Cat@in $\cat{$\bm G$-Cat}$|textbf}
A $G$-equivariant functor $f\colon\mathscr C\to\mathscr D$ of small $G$-categories is called a \emph{$G$-global weak equivalence} if the induced functor $f^{\myh\phi}\colon\mathscr C^{\myh\phi}\to\mathscr D^{\myh\phi}$ is a weak homotopy equivalence for every finite group $H$ and every homomorphism $\phi\colon H\to G$.
\end{defi}

For $G=1$ this recovers the `global equivalences' of \cite[Definition~3.2]{schwede-cat}.

\begin{rk}\label{rk:G-global-we-Fun-EM}
Using that $EH\hookrightarrow E\mathcal M$ is an $H$-equivariant equivalence of right $H$-categories (i.e.~an equivalence in the $2$-category of right $H$-categories, equivariant functors, and equivariant natural transformations) for every $H\subset\mathcal M$, we easily see that $f$ is a $G$-global weak equivalence in the above sense if and only if $\nerve\Fun(E\mathcal M,f)$ is a $G$-global weak equivalence in $\cat{$\bm{E\mathcal M}$-$\bm G$-SSet}$.
\end{rk}

\begin{prop}
The functors $(\blank)^{\textup{sat}}\circ\Phi\colon\cat{$\bm G$-PermCat}\to\cat{$\bm G$-ParSumCat}$ and $\cat{K}_{\textup{$G$-gl}}\colon\cat{$\bm G$-PermCat}\to\cat{$\bm G$-Spectra}$ preserve $G$-global weak equivalences.
\begin{proof}
The first statement is proven as \cite[Lemma~6.11]{sym-mon-global}. The second one now follows immediately from Proposition~\ref{prop:K-G-gl-homotopical}.
\end{proof}
\end{prop}

\subsubsection{An alternative description} We now want to extend the above construction in order to define the $G$-global algebraic $K$-theory of a more general class of symmetric monoidal categories with $G$-actions.

We denote by $\cat{SymMonCat}$\nomenclature[aSymMonCat]{$\cat{SymMonCat}$}{category of small symmetric monoidal categories and strong symmetric monoidal functors} the category of small symmetric monoidal categories and \emph{strong} symmetric monoidal functors, and we define $\cat{SymMonCat}^0$\nomenclature[aSymMonCat]{$\cat{SymMonCat}^0$}{category of small symmetric monoidal categories and strictly unital strong symmetric monoidal functors} as the wide subcategory with morphisms those strong symmetric monoidal functors that are \emph{strictly unital}, i.e.~for which the unit isomorphism is the identity. As usual, we write $\cat{$\bm G$-SymMonCat}$ and $\cat{$\bm G$-SymMonCat}^0$ for the corresponding categories of $G$-objects; in particular, objects of $\cat{$\bm G$-SymMonCat}$ are symmetric monoidal categories with a \emph{strict} $G$-action through \emph{strong} symmetric monoidal functors, whereas objects of $\cat{$\bm G$-SymMonCat}^0$ have $G$-actions through \emph{strictly unital} strong symmetric monoidal functors.

The following lemma is well-known, and we also give a proof as \cite[Proposition~6.7]{sym-mon-global}:

\begin{lemma}\label{lemma:g-perm-0-g-sym-mon}
The inclusions
\begin{equation*}
\cat{$\bm G$-PermCat}\hookrightarrow\cat{$\bm G$-SymMonCat}^0\hookrightarrow\cat{$\bm G$-SymMonCat}
\end{equation*}
are homotopy equivalences (in the sense of Definition~\ref{defi:homotopy-equivalence}) with respect to the underlying equivalences of categories, hence in particular with respect to the $G$-global weak equivalences.\qed
\end{lemma}

The lemma already tells us that we can extend the definition of $G$-global algebraic $K$-theory from $\cat{$\bm G$-PermCat}$ to $\cat{$\bm G$-SymMonCat}$ by simply precomposing with a homotopy inverse to the inclusion. However, the usual homotopy inverse is based on Mac Lane's strictification construction \cite[Theorem~XI.3.1]{cat-working}, which is rather complicated.

Instead, we will now use the construction from Example~\ref{ex:shimada-shimakawa-g-global} to give an explicit construction on $\cat{$\bm G$-SymMonCat}^0$ more in the spirit of Shimakawa's approach to $G$-equivariant algebraic $K$-theory. This already covers many examples arising in practice, and if necessary the extension to $\cat{$\bm G$-SymMonCat}$ can then again be constructed by precomposing with a homotopy inverse (which is much simpler in this case and simply duplicates the tensor unit).

\begin{defi}\label{defi:K-G-gl-prime}\index{G-global algebraic K-theory@$G$-global algebraic $K$-theory!of symmetric monoidal G-categories@of symmetric monoidal $G$-categories|textbf}
We write $\cat{K}'_{\textup{$G$-gl}}$\nomenclature[aKGglprime]{$\cat{K}'_{\textup{$G$-gl}}$}{variant of $\cat{K}_{\textup{$G$-gl}}$ for small symmetric monoidal categories with strictly unital $G$-action} for the composition
\begin{align*}
\cat{$\bm G$-SymMonCat}^0&\xrightarrow{\Gamma}\cat{$\bm\Gamma$-$\bm G$-Cat}_*\xrightarrow{\Fun(E\mathcal M,\blank)}\cat{$\bm\Gamma$-$\bm{E\mathcal M}$-$\bm G$-Cat}_*\\
&\xrightarrow{\nerve}\cat{$\bm\Gamma$-$\bm{E\mathcal M}$-$\bm G$-SSet}_*\xrightarrow{(\blank)\langle\mathbb S\rangle}\cat{$\bm G$-Spectra}.
\end{align*}
Here $\Gamma$ is again obtained from the Segal-May-Shimada-Shimakawa functor (see Examples~\ref{ex:shimada-shimakawa} and~\ref{ex:shimada-shimakawa-equivariant}) by pulling through the $G$-action. Note moreover that while $\nerve\Gamma(\mathscr C)$ will typically not be $G$-equivariantly special, this is rectified by applying $\Fun(E\mathcal M,\blank)$, so that the composite $\cat{$\bm G$-SymMonCat}^0\to\cat{$\bm\Gamma$-$\bm{E\mathcal M}$-$\bm G$-SSet}_*$ actually lands in special $G$-global $\Gamma$-spaces, see Example~\ref{ex:shimada-shimakawa-g-global}.
\end{defi}

In order to relate the above to our previous definition of the $G$-global algebraic $K$-theory of permutative $G$-categories, it will be useful to also conversely express the $G$-global algebraic $K$-theory of saturated $G$-parsummable categories in terms of $\cat{K}_{\textup{$G$-gl}}'$. For this we let $\mu\colon\bm2\times\omega\to\omega$ be any injection. Schwede shows in \cite[Proposition~5.6]{schwede-k-theory} that if $\mathcal C$ is any parsummable category, then its underlying category admits a canonical symmetric monoidal structure with tensor product given on objects by $X\otimes Y=\mu(1,\blank)_*X+\mu(2,\blank)_*Y$; more precisely \cite[Constructions~5.1 and 5.5]{schwede-k-theory}:

\begin{constr}
Let $\mathcal C$ be a parsummable category; for any $n\ge 0$ and any injection $\psi\colon\bm n\times\omega\to\omega$, we write $\psi_*$ for the functor $\mathcal C^{\times n}\to\mathcal C$ given on objects by $\psi_*(X_1,\dots,X_n)\mathrel{:=}\psi(1,\blank)_*(X_1)+\cdots+\psi(n,\blank)_*(X_n)$ and similarly on morphisms; note that this sum is indeed well-defined by Lemma~\ref{lemma:support-vs-action-M}. If $\theta\colon\bm{n}\times\omega\to\omega$ is any other injection, then we write $[\theta,\psi]$ for the natural transformation $\psi_*\Rightarrow\theta_*$ given by the sum $[\theta(1,\blank),\psi(1,\blank)]+\cdots+[\theta(n,\blank),\psi(n,\blank)]$.

Now fix an injection $\mu\colon\bm2\times\omega\to\omega$; we define a symmetric monoidal category $\mu^*\mathcal C$ with underlying category $\mathcal C$, monoidal unit $0$, and tensor product $\mu_*$ as follows: for any $X\in\mathcal C$, the left unitality isomorphism $\bm1\otimes X=\mu_*(0,X)=\mu(2,\blank)_*(X)\to X$ is given by $[1,\mu(2,\blank)]_X$, and similarly for the right unitality isomorphism. If $Y\in\mathcal C$ is another object, then the symmetry isomorphism $\mu_*(X,Y)\to\mu_*(Y,X)$ is given by $[\mu\circ t,\mu]_{(X,Y)}$, where $t\colon\bm2\times\omega\cong\bm 2\times\omega$ exchanges the two copies of $\omega$. Finally, if $Z\in\mathcal C$ is yet another object, then the associativity isomorphism $\mu_*(\mu_*(X,Y),Z)\to\mu_*(X,\mu_*(Y,Z))$ is the natural map $[\mu(\id,\mu),\mu(\mu,\id)]_{X,Y,Z}$; here we write $\mu(f,g)$ for all $f\colon\bm m\times\omega\to\omega, g\colon\bm n\times\omega\to\omega$ for the map with $\mu(f,g)(i,x)=\mu(1,f(i,x))$ for $i\le m$ and $\mu(f,g)(i,x)=\mu(2,g(i-m,x))$ otherwise.
\end{constr}

Any morphism $\mathcal C\to\mathcal D$ of parsummable categories is actually strict symmetric monoidal when viewed as a functor $\mu^*\mathcal C\to\mu^*\mathcal D$ \cite[Remark~5.7]{schwede-k-theory}, in particular yielding a functor $\cat{ParSumCat}\to\cat{SymMonCat}^0$, which we as usual lift to $\cat{$\bm G$-ParSumCat}\to\cat{$\bm G$-SymMonCat}^0$. We will now prove:

\begin{thm}\label{thm:K-G-gl-comparison}
The functor $\cat{K}'_{\textup{$G$-gl}}$ preserves $G$-global weak equivalences and there exist natural equivalences of functors to $\cat{$\bm G$-Spectra}_{\textup{$G$-global}}^\infty$
\begin{align}\label{eq:K-prime-gl-on-mu}
\cat{K}'_{\textup{$G$-gl}}\circ\mu^*|_{\cat{$\bm G$-ParSumCat}^s}&\simeq\cat{K}_{\textup{$G$-gl}}|_{\cat{$\bm G$-ParSumCat}^s}\\
\intertext{(for any injection $\mu\colon\bm2\times\omega\to\omega$) and}
\label{eq:K-gl-on-Phi}
\cat{K}'_{\textup{$G$-gl}}|_{\cat{$\bm G$-PermCat}}&\simeq\cat{K}_{\textup{$G$-gl}}.
\end{align}
\end{thm}

For the proof we will use the following categorical comparison due to Schwede; here we call a map $f$ in $\cat{$\bm\Gamma$-$\bm G$-Cat}_*$ a \emph{categorical equivalence}\index{categorical equivalence} if each $f(S_+)$ is an underlying equivalence of categories.

\begin{prop}
For any injection $\mu\colon{\bm 2}\times\omega\to\omega$ there exists a natural levelwise categorical equivalence $\Psi\colon{\forget}\circ\digamma\Rightarrow\Gamma\circ\mu^*$ of functors $\cat{$\bm G$-ParSumCat}\to\cat{$\bm\Gamma$-$\bm G$-Cat}_*$.
\end{prop}

We will not need any information on the concrete construction of the transformation $\Psi$ below, except that it is given in degree $1^+$ by the canonical isomorphism $\iota\colon\mathscr C\to\Gamma(\mathscr C)(1^+)$ from Example~\ref{ex:shimada-shimakawa}. The curious reader can find the complete definition as \cite[Construction~5.13]{schwede-k-theory}.

\begin{proof}
As the $G$-actions are simply pulled through everywhere, it suffices to prove this for $G=1$, which is \cite[proof of Theorem~5.14-(ii)]{schwede-k-theory}.
\end{proof}

\begin{proof}[Proof of Theorem~\ref{thm:K-G-gl-comparison}]\index{G-global algebraic K-theory@$G$-global algebraic $K$-theory!these are compatible|textbf}
We first observe that $\nerve\circ\Fun(E\mathcal M,\blank)$ sends underlying equivalences of $G$-categories to $G$-global weak equivalences as a consequence of Remark~\ref{rk:G-global-we-Fun-EM}. Applying this with $G$ replaced by $G\times\Sigma_S$ for varying finite set $S$, we therefore see that $\nerve\circ\Fun(E\mathcal M,\blank)$ sends categorical equivalences to $G$-global level weak equivalences, and (using in addition that $\nerve$ and $\Fun(E\mathcal M,\blank)$ preserve products) that it sends categorically special $\Gamma$-$G$-categories to special $G$-global $\Gamma$-spaces.

In particular, as we already recalled in Definition~\ref{defi:K-G-gl-prime}, $\nerve\Fun(E\mathcal M,\Gamma(\blank))$ takes values in special $G$-global $\Gamma$-spaces since $\Gamma$ sends symmetric monoidal categories with strictly unital $G$-actions to categorically special $\Gamma$-$G$-categories. Thus, $\nerve\Fun(E\mathcal M,\Gamma(f))$ is a $G$-global level weak equivalence if and only if $\nerve\Fun(E\mathcal M,f)\cong\nerve\Fun(E\mathcal M,\Gamma(f))(1^+)$ is a $G$-global weak equivalence. The latter is equivalent to $f$ being a $G$-global weak equivalence in $\cat{$\bm G$-SymMonCat}^0$ by Remark~\ref{rk:G-global-we-Fun-EM}, in particular proving that $\cat{K}_{\textup{$G$-gl}}'$ is invariant under $G$-global weak equivalences.

On the other hand, applying this reasoning to the map $\Psi$ from above yields a natural levelwise $G$-global level weak equivalence
\begin{equation*}
\nerve\Fun(E\mathcal M,\Psi)\colon\nerve\circ\Fun(E\mathcal M,\blank)\circ{\forget}\circ\digamma\Rightarrow\nerve\circ\Fun(E\mathcal M,\blank)\circ\Gamma\circ\mu^*.
\end{equation*}
Moreover, we may conclude in the same way from \cite[proof of Theorem~2.33]{schwede-k-theory} that the $G$-global $\Gamma$-space $\nerve\Fun(E\mathcal M,\forget\digamma(\mathcal C))$ is $G$-globally special for any $G$-parsummable category $\mathcal C$.

Since also $\nerve\digamma(\mathcal C)\cong\digamma(\nerve\mathcal C)$ is special, the map
\begin{equation*}
\nerve\eta\colon\nerve\digamma(\mathcal C)\to\nerve\Fun(E\mathcal M,\forget\digamma(\mathcal C))
\end{equation*}
induced by the unit is therefore a $G$-global level weak equivalence if and only if $\eta\colon\mathcal C\to\Fun(E\mathcal M,\forget\mathcal C)$ is a $G$-global weak equivalence, which is in particular true for saturated $\mathcal C$ by Remark~\ref{rk:saturated-alt}. This provides the equivalence $(\ref{eq:K-prime-gl-on-mu})$.

For the equivalence $(\ref{eq:K-gl-on-Phi})$ it is then enough to construct a zig-zag of levelwise $G$-global weak equivalences between $\mu^*\circ(\blank)^\textup{sat}\circ\Phi$ and the inclusion $\cat{$\bm G$-PermCat}\hookrightarrow\cat{$\bm G$-SymMonCat}^0$. We will even give a zig-zag of underlying equivalences: by Theorem~\ref{thm:saturation}, $(\blank)^{\textup{sat}}\simeq\id$, so it suffices that $\mu^*\circ\Phi$ is equivalent to the inclusion. This has been sketched by Schwede in \cite[Remark~11.4]{schwede-k-theory}, and we give a full proof in \cite[Corollary~2.19]{perm-parsum-categorical}.
\end{proof}

\begin{rk}
Unravelling all definitions, one can actually check that the equivalence $(\ref{eq:K-gl-on-Phi})$ is induced by a single natural transformation on the point-set level, and that this point-set level map is independent of the choice of $\mu\colon\bm2\times\omega\to\omega$ in the strong sense that for any other injection $\nu$ the resulting map is even \emph{equal}.
\end{rk}

\begin{rk}\label{rk:homotopy-coherent-actions}\index{permutative G-category@permutative $G$-category!homotopy coherent version|textbf}
Instead of requiring strict actions, we can also canonically extend the above construction to accept small symmetric monoidal categories with \emph{homotopy coherent} $G$-actions. While this generalization will play no role here, let us briefly sketch what goes into this:

It is well-known (and an easy application of Proposition~\ref{prop:enrichment-vs-localization}) that the inclusion $\nerve(\cat{SymMonCat})\hookrightarrow\nerve_\Delta(\cat{SymMonCat}_{2,1})$ is a quasi-localization at the underlying equivalences of categories, where $\cat{SymMonCat}_{2,1}$ is the strict $(2,1)$-category of small symmetric monoidal categories, strong symmetric monoidal functors, and symmetric monoidal isomorphisms. By Lemma~\ref{lemma:g-perm-0-g-sym-mon}, the same is then true for $\gamma\colon\nerve(\cat{PermCat})\hookrightarrow\nerve_\Delta(\cat{SymMonCat}_{2,1})$.

On the other hand, the $1$-category $\cat{PermCat}$ inherits the structure of a homotopy $\cat{U}$-complete category with fibrations and weak equivalences in the sense of \cite[Definitions~7.4.12 and~7.7.2]{cisinski-book} (where $\cat{U}$ is our implicitly chosen Grothendieck universe in the background) from the canonical model structure on $\cat{Cat}$; the only non-formal statement is the existence of factorizations into weak equivalences followed by fibrations, for which one can use the usual mapping cocylinder construction. By Example~7.9.6${}^{\op}$ or alternatively Remark~7.9.7${}^\op$ of \emph{op.~cit.}, $\cat{PermCat}$ is moreover $\cat{U}$-hereditary in the sense of Definition~7.9.4${}^\op$ of \emph{op.~cit.}, so Cisinski's Theorem 7.9.8${}^\op$ shows that quasi-localizations of $\cat{PermCat}$ are stable under taking diagram categories, i.e.~for any small category $I$ the composite
\begin{equation*}
\nerve(\cat{PermCat}^I)\cong\nerve(\cat{PermCat})^{\nerve I}\xrightarrow{\smash{\gamma^{\nerve I}}}\nerve_\Delta(\cat{SymMonCat}_{2,1})^{\nerve I}
\end{equation*}
is a quasi-localization at the underlying equivalences of categories (and the same holds for $\nerve(\cat{SymMonCat}^I)\to\nerve_\Delta(\cat{SymMonCat}_{2,1})^{\nerve I}$ by another application of Lemma~\ref{lemma:g-perm-0-g-sym-mon}). Specializing to $I=BG$ then yields the desired claim as $\cat{K}_{\textup{$G$-gl}}$ is invariant under underlying equivalences of categories.

In more $2$-categorical terms, the objects of $\nerve_\Delta(\cat{SymMonCat}_{2,1})^{\nerve(BG)}$ can be identified with the symmetric monoidal categories with strictly unital pseudofunctorial $G$-action through strong symmetric monoidal functors, see e.g.~\cite[Tags~\href{https://kerodon.net/tag/00AV}{00AV} and~\href{https://kerodon.net/tag/00KY}{00KY}]{kerodon}. While I do not know of a satisfactory proof of this in the literature, it is well-known (and follows for example from the results announced on \cite[p.~204]{duskin}) that $\nerve_\Delta(\cat{SymMonCat}_{2,1})^{\nerve(BG)}$ is in fact equivalent to the homotopy coherent nerve of the strict $(2,1)$-category of symmetric monoidal categories with pseu\-do\-functorial $G$-actions, pseudoequivariant strong symmetric monoidal functors, and pseudoequivariant symmetric monoidal isomorphisms.
\end{rk}

\begin{rk}\label{rk:G-global-K-rings-symmetric-monoidal}\index{G-global algebraic K-theory@$G$-global algebraic $K$-theory!these are compatible|textbf}
Let $R$ be any $G$-ring, and fix an injection $\mu\colon\bm2\times\omega\to\omega$. Then applying Theorem~\ref{thm:K-G-gl-comparison} to the saturated $G$-parsummable category $\mathcal P(R)$ (Examples~\ref{ex:G-global-K-rings} and~\ref{ex:module-saturated}), shows that the $G$-global algebraic $K$-theory $\cat{K}_{\textup{$G$-gl}}(R)$ of $R$ is $G$-globally weakly equivalent to the $G$-global algebraic $K$-theory of the symmetric monoidal category $\mu^*\mathcal P(R)$.

Let us now write $\mathscr P(R)$ for the category of finitely generated projective $R$-modules and $R$-linear isomorphisms. To be entirely precise, $\mathscr P(R)$ is not small as there are too many representatives of each isomorphism class, so we should restrict to a small essentially wide subcategory; to simplify the exposition below, we do this as follows: we fix any set $\mathscr U$ of sets containg all subsets of $R^{(\omega)}$ and restrict to those projective $R$-modules $M$ whose underlying set belongs to $\mathscr U$.

The category $\mathscr P(R)$ inherits a symmetric monoidal structure from `the' cocartesian symmetric monoidal structure on the ambient category of $R$-modules. While the result will of course be independent of the chosen coproducts up to canonical isomorphism, we will make a very specific choice now that will trivialize several computations below: namely, we insist that the monoidal unit be the trivial submodule of $R^{(\omega)}$, that our choices of coproducts be obtained from fixed choices of coproducts in $\cat{Ab}$ (so that the underlying abelian groups of our chosen coproducts only depend on the underlying abelian groups of their inputs \emph{on the nose}), and that for all $M,N\in\mathcal P(R)$ their chosen coproduct be given by $\mu_*(M,N)$ with structure maps
\begin{equation*}
M\xrightarrow{\smash{[\mu(1,\blank),1]}}\mu(1,\blank)_*(M)\hookrightarrow \mu_*(M,N)
\text{ and }
N\xrightarrow{\smash{[\mu(2,\blank),1]}}\mu(2,\blank)_*(N)\hookrightarrow \mu_*(M,N).
\end{equation*}
Following \cite[4.3]{merling}, we can make $\mathscr P(R)$ into an object in $\cat{$\bm G$-SymMonCat}^0$ as follows: if $g\in G$ is arbitrary, and $M\in\mathscr P(R)$, then ${}^gM$ has the same underlying abelian group, but the scalar multiplication is instead given by $r\cdot_{g}m=(g^{-1}.r)\cdot m$ (note that Merling lets $r$ act by $g.r$ instead, but this would give a \emph{right} $G$-action); moreover, if $f\colon M\to N$ is any morphism, then ${}^gf\colon {}^gM\to{}^gN$ agrees with $f$ as a morphism of underlying abelian groups. We omit the easy verification that this is a well-defined $G$-action on the underlying category. Our specific choice of coproducts then ensures that each ${}^g(\blank)$ is actually strict symmetric monoidal, in particular yielding an object in $\cat{$\bm G$-SymMonCat}^0$.

Ignoring $G$-actions, the inclusion $\iota\colon\mu^*\mathcal P(R)\hookrightarrow\mathscr P(R)$ is again strict symmetric monoidal by our specific choice of coproducts, and it is clearly an equivalence of categories. However, $\iota$ is typically \emph{not} $G$-equivariant unless $G$ acts trivially on $R$: namely, if $M\in\mathcal P(R)$, then $g.M=(g.\blank)(M)$ will typically have a different underlying set than $M$ and hence ${}^gM$. However, acting with $g^{-1}$ on $R^{(\omega)}$ still gives us a preferred isomorphism $g.M \to {}^gM$, and one easily checks that for varying $g$ and $M$ these equip $\iota$ with the structure of a pseudoequivariant strong symmetric monoidal functor.

Thus, taking the previous remark for granted, we see that $\cat{K}'_{\textup{$G$-gl}}(\mu^*\mathcal P(R))$ and $\cat{K}'_{\textup{$G$-gl}}(\mathscr P(R))$ are $G$-globally weakly equivalent, and hence so are $\cat{K}_{\textup{$G$-gl}}(R)$ and $\cat{K}'_{\textup{$G$-gl}}(\mathscr P(R))$. However, we can actually also prove this equivalence unconditionally, using a trick of Merling's:

Consider $\Fun(EG,\mathscr P(R))$ with the levelwise symmetric monoidal structure and $G$-action via conjugation; here we view $EG$ as coming with the usual \emph{left} $G$-action as opposed to the right action considered before in order to keep our notation consistent with Merling's. We define a functor $\tilde\iota\colon\mu^*\mathcal P(R)\to\Fun(EG,\mathscr P(R))$ as follows: if $M\in\mathcal P(R)$ is arbitrary, then $\tilde\iota(M)(g)={}^g\iota(g^{-1}.M)$ for all $g\in G$, with structure maps $\tilde\iota(M)(g_2,g_1)=g_2^{-1}g_1.\blank\colon{}^{g_1}\iota(g_1^{-1}.M)\to{}^{g_2}\iota(g_2^{-1}.M)$; moreover, if $f\colon M\to N$ is any morphism in $\mathcal P(R)$, then $\tilde\iota(f)_g=\iota(g^{-1}.f)$. \cite[Proposition~3.3]{merling} (restricted to constant diagrams) shows that this is a well-defined and strictly $G$-equivariant functor. Moreover, it is once again strict symmetric monoidal, and it is an equivalence of underlying categories as it becomes one after evaluating at $1\in EG$. We therefore obtain a zig-zag
\begin{equation*}
\mu^*\mathcal P(R)\xrightarrow{\tilde\iota}\Fun(EG,\mathscr P(R))\xleftarrow{\const}\mathscr P(R)
\end{equation*}
of underlying equivalences in $\cat{$\bm G$-SymMonCat}^0$, which is then sent under $\cat{K}'_{\textup{$G$-gl}}$ to the desired zig-zag of $G$-global weak equivalences.
\end{rk}

Finally we note that we have analogously to Proposition~\ref{prop:K-G-gl-parsum-functoriality}:

\begin{prop}\label{prop:K-G-gl-symmon-functoriality}
Let $\phi\colon H\to G$ be any group homomorphism. Then the diagram
\begin{equation*}
\begin{tikzcd}[column sep=large]
(\cat{$\bm G$-SymMonCat}^0)^\infty\arrow[d, "(\phi^*)^\infty"']\arrow[r, "(\cat{K}'_{\textup{$G$-gl}})^\infty"]&[1em]\cat{$\bm G$-Spectra}^\infty_{\textup{$G$-global}}\arrow[d,"(\phi^*)^\infty"]\\
(\cat{$\bm H$-SymMonCat}^0)^\infty\arrow[r, "(\cat{K}'_{\textup{$H$-gl}})^\infty"']&\cat{$\bm H$-Spectra}_{\textup{$H$-global}}^\infty
\end{tikzcd}
\end{equation*}
commutes up to preferred equivalence.\qed
\end{prop}

\subsubsection{Comparison to $G$-equivariant algebraic $K$-theory} Let $G$ be finite; the alternative description of $G$-global algebraic $K$-theory provided by Theorem~\ref{thm:K-G-gl-comparison} above allows for an easy comparison to Shimakawa's construction of $G$-equivariant algebraic $K$-theory \cite[Theorem~$\text{A}'$]{shimakawa}:

\begin{defi}\index{algebraic K-theory@algebraic $K$-theory!G-equivariant@$G$-equivariant|seeonly{$G$-equivariant algebraic $K$-theory}}\index{G-equivariant algebraic K-theory@$G$-equivariant algebraic $K$-theory|textbf}
We write $\cat{K}_G$\nomenclature[aKG]{$\cat{K}_G$}{$G$-equivariant algebraic $K$-theory} for the composition
\begin{align*}
\cat{$\bm G$-SymMonCat}^0&\xrightarrow{\Gamma}\cat{$\bm\Gamma$-$\bm G$-Cat}_*\xrightarrow{\Fun(EG,\blank)}\cat{$\bm\Gamma$-$\bm G$-Cat}_*\\
&\xrightarrow\nerve\cat{$\bm\Gamma$-$\bm G$-SSet}_*\xrightarrow{\mathcal E_G}\cat{$\bm G$-Spectra};
\end{align*}
here $G$ acts on $EG$ from the right in the obvious way.

If $\mathscr C$ is any small symmetric monoidal category with strictly unital $G$-action, then $\cat{K}_G(\mathscr C)$ is called the \emph{$G$-equivariant algebraic $K$-theory} of $\mathscr C$.
\end{defi}

Strictly speaking, Shimakawa uses a bar construction instead of the prolongation functor $\mathcal E_G$; however, the two delooping constructions are equivalent, as explained in \cite[Sections~2--3]{may-merling-osorno}.

\begin{thm}\label{thm:G-global-K-vs-G-equiv}\index{G-equivariant algebraic K-theory@$G$-equivariant algebraic $K$-theory!vs G-global algebraic K-theory@vs.~$G$-global algebraic $K$-theory|textbf}\index{G-global algebraic K-theory@$G$-global algebraic $K$-theory!vs G-equivariant algebraic K-theory@vs $G$-equivariant algebraic $K$-theory|textbf}
There are natural equivalences of functors
\begin{align*}
(\cat{K}_G|_{\cat{$\bm G$-PermCat}})^\infty\simeq\und_G^\infty\circ(\cat{K}_{\textup{$G$-gl}})^\infty&\colon\cat{$\bm G$-PermCat}^\infty\to\cat{$\bm G$-Spectra}_{\textup{$G$-equiv.}}^\infty\\
\cat{K}_G^\infty\simeq\und_G^\infty\circ(\cat{K}'_{\textup{$G$-gl}})^\infty&\colon(\cat{$\bm G$-SymMonCat}^0)^\infty\to\cat{$\bm G$-Spectra}_{\textup{$G$-equiv.}}^\infty.
\end{align*}
\end{thm}

For the proof we will use:

\begin{lemma}\label{lemma:comparison-E-s-Gamma}
Let $i\colon G\to\mathcal M$ be an injective homomorphism with universal image. Then the diagram
\begin{equation*}
\begin{tikzcd}
\cat{$\bm\Gamma$-$\bm{E\mathcal M}$-$\bm G$-SSet}_*\arrow[d, "i^*"'] \arrow[r, "(\blank)\langle\mathbb S\rangle"] & \cat{$\bm G$-Spectra}_\textup{$G$-global}\arrow[d,"="]\\
\cat{$\bm\Gamma$-$\bm G$-SSet}_*\arrow[r, "\mathcal E_G"'] & \cat{$\bm G$-Spectra}_\textup{$G$-equivariant}
\end{tikzcd}
\end{equation*}
of homotopical functors commutes up to a zig-zag of levelwise natural weak equivalences.
\begin{proof}
Plugging in the definition of $(\blank)\langle\mathbb S\rangle$, the above diagram can be expanded to
\begin{equation}\label{eq:comparison-deloopings}
\begingroup\hfuzz=10in
\setbox0=\hbox{
\begin{tikzcd}[row sep=1.75em]
\cat{$\bm\Gamma$-$\bm{E\mathcal M}$-$\bm G$-SSet}_*\arrow[r, "{(\blank)[\omega^\bullet]}"]\arrow[rd, bend right=15pt, "="'] & \cat{$\bm\Gamma$-$\bm G$-$\bm{\mathcal I}$-SSet}_*\arrow[r, "\mathcal E^\otimes"]\arrow[d, "\ev_\omega"']\arrow[dd, xshift=25pt, bend left=45pt, "\ev_{i^*\omega}"] & \cat{$\bm G$-Spectra}_{\textup{$G$-global}}\arrow[dd, "="]\\
& \cat{$\bm\Gamma$-$\bm{E\mathcal M}$-$\bm G$-SSet}_*\arrow[d, "i^*"']\\
& \cat{$\bm\Gamma$-$\bm{G}$-SSet}_*\arrow[r, "\mathcal E_G"'] & \cat{$\bm G$-Spectra}_{\textup{$G$-equivariant}\rlap.}
\end{tikzcd}}
\hbox to 0pt{\hskip-.5\wd0\unhbox0}
\endgroup
\end{equation}
As $i^*\omega$ is a complete $G$-set universe, Lemma~\ref{lemma:und-E-vs-E-H} (for $\phi=\id$) shows that the rectangle on the right commutes up to a zig-zag of levelwise $G$-equivariant weak equivalences; explicitly, this is given for a $G$-global $\Gamma$-space $X$ in degree $A$ by the maps
\begin{equation*}
(\mathcal E^\otimes X)(A) = X(A)(S^A) \to X(A\amalg i^*\omega)(S^A)\gets X(i^*\omega)(S^A) =(\mathcal E_G X(i^*\omega))(A)
\end{equation*}
induced by the inclusions $A\hookrightarrow A\amalg i^*\omega\hookleftarrow i^*\omega$.

On the other hand, the triangle on the left commutes up to $G$-global level weak equivalence by Theorem~\ref{thm:gamma-ev-omega}, which finishes the proof.
\end{proof}
\end{lemma}

\begin{proof}[Proof of Theorem~\ref{thm:G-global-K-vs-G-equiv}]
By Theorem~\ref{thm:K-G-gl-comparison} it suffices to prove the second statement. For this we choose an injective homomorphism $i\colon G\to\mathcal M$ with universal image; by the above lemma, it will then be enough to show that the $\Gamma$-$G$-spaces $i^*\nerve\Fun(E\mathcal M,C)$ and $\nerve\Fun(EG,C)$ are naturally $G$-equivariantly weakly equivalent for any $\Gamma$-$G$-category $C$. But $Ei\colon EG\to E\mathcal M$ is an equivalence in the $2$-category of right $G$-categories, equivariant functors, and equivariant natural transformations; thus, restricting along $Ei$ produces the desired $G$-equivariant level weak equivalence.
\end{proof}

\begin{rk}
From a higher categorical point of view, the choice of $i$ can again be shown to be inessential: for example, the right hand portion of $(\ref{eq:comparison-deloopings})$ can be parameterized over the contractible $2$-groupoid with objects the injective homomorphisms $G\to\mathcal M$, $1$-morphisms $i\to j$ the invertible $\phi\in\mathcal M$ such that $j(g) = \phi i(g) \phi^{-1}$ for all $g\in G$ (cf.~Lemma~\ref{lemma:uniqueness-of-universal-homom}), and a unique $2$-cell between any pair of parallel arrows.
\end{rk}

\begin{rk}\label{rk:G-global-K-rings-vs-G-equivariant}\index{G-global algebraic K-theory@$G$-global algebraic $K$-theory!vs G-equivariant algebraic K-theory@vs.~$G$-equivariant algebraic $K$-theory|textbf}\index{G-equivariant algebraic K-theory@$G$-equivariant algebraic $K$-theory!vs G-global algebraic K-theory@vs.~$G$-global algebraic $K$-theory|textbf}
Let $R$ be any $G$-ring. By the comparison from Remark~\ref{rk:G-global-K-rings-symmetric-monoidal}, we in particular conclude from the above that the underlying $G$-equivariant spectrum of $\cat{K}_{\textup{$G$-gl}}(R)$ recovers $\cat{K}_G(\mathscr P(R))$. Thus, the $G$-global algebraic $K$-theory of $G$-rings as introduced above refines Merling's $G$-equivariant construction by \cite[Corollary~5.26]{merling}.
\end{rk}

Together with Proposition~\ref{prop:K-G-gl-symmon-functoriality} applied to the unique homomorphism ${G\to1}$, we in particular get the following comparison between global and equivariant algebraic $K$-theory:

\begin{cor}\index{global algebraic K-theory@global algebraic $K$-theory!vs G-equivariant algebraic K-theory@vs.~$G$-equivariant algebraic $K$-theory|textbf}
The diagram
\begin{equation*}
\begin{tikzcd}
\cat{PermCat}^\infty\arrow[d, "\triv^\infty"']\arrow[r,"\cat{K}_{\textup{gl}}^\infty"]& \cat{Spectra}^\infty_{\textup{global}}\arrow[d, "\und_G^\infty"]\\
(\cat{$\bm G$-SymMonCat}^0)^\infty\arrow[r, "\cat{K}_{G}^\infty"'] & \cat{$\bm G$-Spectra}^\infty_{\textup{$G$-equivariant}}
\end{tikzcd}
\end{equation*}
commutes up to natural equivalence.\qed
\end{cor}

\begin{rk}\label{rk:global-K-rings-vs-G-equivariant}
Analogously, we conclude from Remark~\ref{rk:G-global-K-rings-vs-G-equivariant} above that Schwede's global algebraic $K$-theory of rings forgets for any finite group $G$ to Merling's $G$-equivariant construction restricted to rings with trivial action, confirming the expectation expressed in \cite[Remark~10.4]{schwede-k-theory}.
\end{rk}

\section[$G$-global versions of the Barratt-Priddy-Quillen Theorem]{\texorpdfstring{$\except{toc}{\bm G}\for{toc}{G}$}{G}-global versions of the Barratt-Priddy-Quillen Theorem}
It is custom that any paper on infinite loop space machinery should prove some version of the Barratt-Priddy-Quillen Theorem, see e.g.~\cite[Propositions~3.5 and~3.6]{segal-gamma}, \cite[Theorems~6.1 and~6.2]{guillou-may}, and \cite[Theorems~8.7 and~9.7]{schwede-k-theory} for non-equivariant, equivariant, and global versions, respectively. In this section we will honour this tradition by providing several versions for the two different constructions of $G$-global algebraic $K$-theory given above.

\subsection{Versions for $\bm G$-parsummable simplicial sets} We begin by identifying the $G$-global spectra associated to certain free $G$-parsummable simplicial sets. For this we will use the following comparison of $G$-global $\Gamma$-spaces, which one can view as a `non-group-completed' $G$-global Barratt-Priddy-Quillen Theorem, also cf.~the `special' Barratt-Priddy-Quillen Theorem of \cite[Theorem~1.2]{boavida-moerdijk}:

\index{Barratt-Priddy-Quillen Theorem|seeonly{$G$-global Barratt-Priddy-Quillen Theorem}}
\begin{prop}\label{prop:BPQ-non-group}\index{G-global Barratt-Priddy-Quillen Theorem@$G$-global Barratt-Priddy-Quillen-Theorem!non-group-completed|textbf}
Let $X\in\cat{$\bm{E\mathcal M}$-$\bm G$-SSet}^\tau_*$ be cofibrant in the positive $G$-global model structure. Then the map
\begin{equation*}
\tilde\eta\colon\Gamma(1^+,\blank)\smashp X\to \digamma(\SP^\infty X)
\end{equation*}
adjoint to the unit $\eta\colon X\to\forget\SP^\infty X=\digamma(\SP^\infty X)(1^+)$ is a $G$-global special weak equivalence.
\begin{proof}
Obviously, $\ev_{1^+}\colon\cat{$\bm\Gamma$-$\bm{E\mathcal M}$-$\bm G$-SSet}_*^\tau\to\cat{$\bm{E\mathcal M}$-$\bm G$-SSet}_*^\tau$ is right Quillen for the $G$-global positive level model structure on the source and the $G$-global positive model structure on the target, hence also with respect to the $G$-global special model structure on the source. Thus, we have a Quillen adjunction
\begin{equation*}
i_!\mathrel{:=}\Gamma(1^+,\blank)\smashp\blank\colon(\cat{$\bm{E\mathcal M}$-$\bm G$-SSet}_*^\tau)_{\textup{positive}}\rightleftarrows(\cat{$\bm\Gamma$-$\bm{E\mathcal M}$-$\bm G$-SSet}_*^\tau)_{\textup{special}} :\!\ev_{1^+}
\end{equation*}
and one easily checks that the left adjoint is fully homotopical, while the right adjoint is homotopical in $G$-global level weak equivalences; in particular, $\cat{R}\ev_{1^+}$ can be computed by choosing for any $G$-global $\Gamma$-space a $G$-global special weak equivalence $\iota\colon X\to X'$ to some special $X'$; we agree that $\iota$ should be the identity whenever $X$ is already special. With this convention, the middle square in
\begingroup\hfuzz=2in
\begin{equation*}
\hskip-42pt
\begin{tikzcd}[cramped]
& \arrow[dl, bend right=15pt, "="'{name=A}]\hskip-2pt\Ho(\cat{$\bm\Gamma$-$\bm{E\mathcal M}$-$\bm G$-SSet}_*^\tau)\hskip-2pt\arrow[d, "\cat{R}\ev_{1^+}"'] &[-.45em] \arrow[l, "\;\Ho(\digamma)"'] \hskip-2pt\Ho(\cat{ParSumCat})\hskip-2pt\arrow[d, "\Ho(\forget)"] & \arrow[l, "\;\cat{L}\SP^\infty"'] \Ho(\cat{$\bm{E\mathcal M}$-$\bm G$-SSet}_*^\tau)\arrow[dl, bend left=15pt, "="{name=B}]\\
\Ho(\cat{$\bm\Gamma$-$\bm{E\mathcal M}$-$\bm G$-SSet}_*^\tau)\hskip-2pt& \twocell[to=A,"\scriptscriptstyle\epsilon"{yshift=5pt,xshift=2pt}]\arrow[l, "{\Ho(i_!)}"]\hskip-2pt\Ho(\cat{$\bm{E\mathcal M}$-$\bm G$-SSet}_*^\tau)\hskip-2pt &\arrow[l, "\,="] \hskip-2pt\Ho(\cat{$\bm{E\mathcal M}$-$\bm G$-SSet}_*^\tau)\hskip-2pt\twocell[ul, "\scriptscriptstyle\id"{yshift=5pt,xshift=2pt}]\twocell[from=B, u, "\scriptscriptstyle\eta"{yshift=5pt,xshift=2pt}]
\end{tikzcd}
\end{equation*}
\endgroup
commutes \emph{strictly}. On the other hand, also $\Ho(\forget)$ admits a left adjoint as the forgetful functor is right Quillen (with left adjoint given by $\SP^\infty$). Thus, we can pass to canonical mates with respect to the vertical arrows, yielding the above pasting. As the horizontal maps in the middle square are equivalences (the non-trivial case being Theorem~\ref{thm:gamma-vs-uc}), this pasting is an isomorphism. Chasing through $X$ therefore yields an isomorphism $\Gamma(1^+,\blank)\smashp X\cong \digamma(\SP^\infty X)$ in the homotopy category and it only remains to show that this is represented by $\tilde\eta$.

Indeed, as $\forget$ is fully homotopical, we can choose the derived unit of $\cat{L}\SP^\infty\dashv\forget$ to be represented by the ordinary unit on all cofibrant objects; similarly, we can take the derived counit of $\Ho(i_!)\dashv\cat{R}\,\ev_{1^+}$ to be represented by the ordinary counit on all special $G$-global $\Gamma$-spaces. With these choices, the above isomorphism in the homotopy category is then represented by the composite
\begin{align*}
\Gamma(1^+,\blank)\smashp X&\xrightarrow{\Gamma(1^+,\blank)\smashp\eta}\Gamma(1^+,\blank)\smashp\forget\SP^\infty X\\
&=\Gamma(1^+,\blank)\smashp\digamma(\SP^\infty X)(1^+)\xrightarrow{\epsilon}\digamma(\SP^\infty X),
\end{align*}
which is precisely $\tilde\eta$ as desired.
\end{proof}
\end{prop}

\begin{constr}
Define $\rho\colon\Sigma^\bullet Y(S^0)\to\mathcal E^\otimes Y$ for any $Y\in\cat{$\bm\Gamma$-$\bm G$-$\bm{\mathcal I}$-SSet}_*$\nomenclature[arho3]{$\rho$}{a specific map from the suspension spectrum of the underlying space of a global $\Gamma$-space to its delooping} as the map given in degree $A$ by
\begin{equation*}
S^A\smashp Y(S^0)(A)\xrightarrow{\asm} Y(S^A\smashp S^0)(A)\cong Y(S^A)(A)=(\mathcal E^\otimes Y)(A)
\end{equation*}
where the unlabelled isomorphism is induced by the unitality isomorphism. It is clear that $\rho$ is natural.
\end{constr}

\begin{ex}
Let $X\in\cat{$\bm G$-$\bm{\mathcal I}$-SSet}_*$ be arbitrary. Then
\begin{equation*}
\rho\colon\Sigma^\bullet(\Gamma(1^+,\blank)\smashp X)(S^0)\to\mathcal E^\otimes(\Gamma(1^+,\blank)\smashp X)
\end{equation*}
is an isomorphism: indeed, arguing levelwise we are reduced to the statement that the analogous map $\Sigma^\infty(\Gamma(1^+,\blank)\smashp X')(S^0)\to(\Gamma(1^+,\blank)\smashp X')(\mathbb S)$ is an isomorphism for every $X'\in\cat{SSet}_*$, which is well-known (and just an instance of the enriched co-Yoneda Lemma).

In particular, precomposing with the unit isomorphism yields an isomorphism $\varrho\colon \Sigma^\bullet X\to\mathcal E^\otimes(\Gamma(1^+,\blank)\smashp X)$.
\end{ex}

\begin{thm}\index{G-global Barratt-Priddy-Quillen Theorem@$G$-global Barratt-Priddy-Quillen-Theorem!for pointed G-parsummable simplicial sets@for pointed $G$-parsummable simplicial sets}\label{thm:BPQ-group-completed-tame-EM-pointed}
The composition
\begin{equation}\label{eq:BPQ-ParSumSSet-we}
\Sigma^\bullet X[\omega^\bullet]\xrightarrow{\Sigma^\bullet\eta[\omega^\bullet]}
\Sigma^\bullet (\SP^\infty X)[\omega^\bullet]\xrightarrow{\rho}
\mathcal E^\otimes\big(\digamma(\SP^\infty X)[\omega^\bullet]\big)=
\digamma(\SP^\infty X)\langle\mathbb S\rangle
\end{equation}
is a $G$-global weak equivalence for any positively cofibrant $X\in\cat{$\bm{E\mathcal M}$-$\bm G$-SSet}_*^\tau$.
\begin{proof}
We begin with the following observation:

\begin{claim*}
The functor $(\blank)\langle\mathbb S\rangle=\mathcal E^\otimes\circ(\blank)[\omega^\bullet]\colon\cat{$\bm\Gamma$-$\bm{E\mathcal M}$-$\bm G$-SSet}_*^\tau\to\cat{$\bm G$-Spectra}$ sends $G$-global special weak equivalences to $G$-global weak equivalences.
\begin{proof}
By Theorem~\ref{thm:gamma-ev-omega} and Corollary~\ref{cor:E-otimes-homotopical}, $(\blank)\langle\mathbb S\rangle$ descends to a functor
\begin{equation}\label{eq:eval-sphere-level-induced}
(\cat{$\bm\Gamma$-$\bm{E\mathcal M}$-$\bm G$-SSet}_*^\tau)^\infty_{\textup{$G$-global level}}\to\cat{$\bm G$-Spectra}^\infty_{\textup{$G$-global}},
\end{equation}
and it suffices that this inverts $G$-global special weak equivalences. However, $(\ref{eq:eval-sphere-level-induced})$ is left adjoint to $\ev_\omega^\infty\circ\cat{R}\Phi^\otimes$ by Theorem~\ref{thm:gamma-ev-omega} together with Theorem~\ref{thm:group-completion}, and it is then enough that this right adjoint takes values in the local objects for the special $G$-global weak equivalences, i.e.~the special $G$-global $\Gamma$-spaces. This is however immediate from Theorem~\ref{thm:group-completion} and Corollary~\ref{cor:comparison-special}.
\end{proof}
\end{claim*}

Applying this to Proposition~\ref{prop:BPQ-non-group}, we therefore get a $G$-global weak equivalence $\tilde\eta\langle\mathbb S\rangle\colon  (\Gamma(1^+,\blank)\smashp X)\langle\mathbb S\rangle\to\digamma(\SP^\infty X)\langle\mathbb S\rangle$ while the previous example tells us that $\varrho\colon\Sigma^\bullet X[\omega^\bullet]\to (\Gamma(1^+,\blank)\smashp X)\langle\mathbb S\rangle$ is an isomorphism. The claim follows as $(\ref{eq:BPQ-ParSumSSet-we})$ agrees with the composition $\tilde\eta\langle\mathbb S\rangle\circ\varrho$ by naturality of $\rho$ and the definition of $\tilde\eta$.
\end{proof}
\end{thm}

Next, let us give an unpointed version:

\begin{cor}\label{cor:BPQ-SSet-unpointed}\index{G-global Barratt-Priddy-Quillen Theorem@$G$-global Barratt-Priddy-Quillen-Theorem!for G-parsummable simplicial sets@for $G$-parsummable simplicial sets}
Let $X\in\cat{$\bm{E\mathcal M}$-$\bm G$-SSet}^\tau\!$, and assume that $X$ has no vertices of empty support. Then the composition
\begin{align*}
\Sigma^\bullet_+ X[\omega^\bullet]
=\Sigma^\bullet X_+[\omega^\bullet]
\xrightarrow{\Sigma^\bullet\tilde\eta[\omega^\bullet]}\Sigma^\bullet(\cat{P}X)[\omega^\bullet]\xrightarrow{\rho}\digamma(\cat{P}X)\langle\mathbb S\rangle
\end{align*}
is a $G$-global weak equivalence; here $\tilde\eta$ is the adjoint of the unit $X\to\forget\cat{P}X$.
\begin{proof}
If $X$ is cofibrant in the $G$-global positive model structure, the claim follows from the previous theorem using that there exists a (unique) isomorphism $\cat{P}X\cong\SP^\infty(X_+)$ compatible with the units.

The general case now follows by $2$-out-of-$3$ since both sides are homotopical in $G$-global weak equivalences between objects without vertices of empty support: the left hand side is even fully homotopical by Proposition~\ref{prop:suspension-loop-projective} together with Proposition~\ref{prop:omega-bullet-inverse}, and so are $\digamma$ and $(\blank)\langle\mathbb S\rangle$ as recalled in the proof of Proposition~\ref{prop:K-G-gl-homotopical}. Finally, $\cat P=\coprod_{n\ge0}\Sym^n$ preserves $G$-global weak equivalences between objects without vertices of empty support by Corollary~\ref{cor:sym-homotopical-EM} applied to each summand of the coproduct.
\end{proof}
\end{cor}

Let $H$ be a finite group, and let $A$ be a non-empty finite faithful $H$-set. For the next theorem, we recall the bijection $\tau\colon \omega\cong\omega^A$ chosen in Construction \ref{constr:omega-bullet}. The given $H$-action on $A$ then induces an $H$-action on $\omega^A$, which we transport to $\omega$ along $\tau$.

\begin{thm}\label{thm:BPQ-SSet-corep}
Let $A$ and $H$ be as above, let $\phi\colon H\to G$ be a homomorphism, and let $f\colon A\to\omega$ be an injection that is $H$-equivariant with respect to the action on $\omega$ described above. Then the composition
\begin{equation*}
\Sigma^\bullet_+(I(A,\blank)\times_\phi G)\xrightarrow{\Sigma^\bullet(\tilde\eta[\omega^\bullet]\circ\tilde f_+)}
\Sigma^\bullet\cat{P}(E\Inj(A,\omega)\times_\phi G)\xrightarrow\rho\digamma\big(\cat{P}(E\Inj(A,\omega)\times_\phi G)\big)\langle\mathbb S\rangle
\end{equation*}
is a $G$-global weak equivalence, where $\tilde f\colon I(A,\blank)\times_\phi G\to(E\Inj(A,\omega)\times_\phi G)[\omega^\bullet]$ is the map classifying $[f;1]$. Moreover, such an injection $f$ exists, and for any two such choices the resulting maps are simplicially homotopic.
\begin{proof}
One easily checks that $[f;1]$ is $\phi$-fixed (so that $\tilde f$ is well-defined), and that for any other $H$-equivariant injection $f'$ so is the edge $[f',f;1]$, so that $\tilde f$ and $\tilde f'$ are simplicially homotopic. By the previous corollary, it therefore suffices to give one such map $f$ and to show that the induced map $\Sigma^\bullet_+(I(A,\blank)\times_\phi G)\to \Sigma^\bullet_+\big((E\Inj(A,\omega)\times_\phi G)[\omega^\bullet]\big)$ is a $G$-global weak equivalence, for which it is in turn enough by Proposition~\ref{prop:suspension-loop-projective} that $\tilde f$ itself is a $G$-global weak equivalence.

For this we take $f=\tau\chi$, where $\chi\colon A\to\omega^A$ sends $a\in A$ to its characteristic function $1_a$. This is clearly $H$-equivariant, and moreover $\tilde f$ factors up to the canonical isomorphism $E\Inj(A,\omega)\cong\mathcal I(A,\blank)(\omega)$ as the composition of the inclusion $I(A,\blank)\times_\phi G\hookrightarrow\mathcal I(A,\blank)\times_\phi G$ (which is a $G$-global weak equivalence by Theorem~\ref{thm:strict-global-I-model-structure}) with the map $\theta_X$ from Remark~\ref{rk:theta} (which is a $G$-global weak equivalence by the proof of Proposition~\ref{prop:omega-bullet-inverse}), finishing the proof.
\end{proof}
\end{thm}

We may also think of $[f;1]$ as a vertex of $\cat{P}(E\Inj(A,\omega)\times_\phi G)$, and by slight abuse of notation we also denote the map ${(I(A,\blank)\hskip-1.5pt\times_\phi\hskip-1pt G)_+\hskip-2pt\to\hskip-1pt\cat{P}\hskip-.5pt(E\Inj(A,\omega)\hskip-1.5pt\times_\phi \hskip-1ptG)[\omega^A]}$ classifying $[f;1]$ by $\tilde f$. With this notation, the left hand map in the above composition will then simply be written as $\Sigma^\bullet\tilde f$.

\subsection{Versions for $\bm G$-parsummable categories} Also the forgetful functor $\cat{$\bm G$-ParSumCat}\to\cat{$\bm{E\mathcal M}$-$\bm G$-Cat}^\tau$ has a left adjoint $\cat{P}$. Using the above, we can now very easily describe the $G$-global algebraic $K$-theory of certain \emph{free $G$-parsummable categories} $\cat{P}\mathcal C$, generalizing Schwede's global Barratt-Priddy-Quillen Theorem \cite[Theorem~8.7]{schwede-k-theory}:

\begin{thm}\index{G-global Barratt-Priddy-Quillen Theorem@$G$-global Barratt-Priddy-Quillen-Theorem!for G-parsummable categories@for $G$-parsummable categories|textbf}
Let $\mathcal C$ be a tame $E\mathcal M$-$G$-category without objects of empty support. Then the composition
\begin{align*}
\Sigma^\bullet_+\nerve(\mathcal C)[\omega^\bullet]\xrightarrow{\Sigma^\bullet\widetilde{\nerve(\eta)}[\omega^\bullet]}\Sigma^\bullet\nerve(\forget\cat{P}\mathcal C)[\omega^\bullet]\xrightarrow{\rho}\nerve\big(\digamma(\cat{P}\mathcal C)\big)\langle\mathbb S\rangle=\cat{K}_{\textup{$G$-gl}}(\cat{P}\mathcal C)
\end{align*}
is a $G$-global weak equivalence.
\begin{proof}
By adjointness, there is a unique map $p\colon\cat{P}(\nerve\mathcal C)\to\nerve(\cat{P}\mathcal C)$ compatible with the units of the two free-forgetful adjunctions, and we claim that this map is an isomorphism. With this established, the theorem will follow from Corollary~\ref{cor:BPQ-SSet-unpointed}.

For the proof of the claim it suffices to observe that $p$ is the coproduct of the maps $\nerve(\mathcal C)^{\boxtimes n}/\Sigma_n\cong \nerve(\mathcal C^{\boxtimes n})/\Sigma_n\to\nerve(\mathcal C^{\boxtimes n}/\Sigma_n)$, where the first map is induced by the strong symmetric monoidal structure on $\nerve$, and the second one is the canonical comparison map. However, $\Sigma_n$ again acts freely on $\mathcal C^{\boxtimes n}$ as $\mathcal C$ has no objects of empty support, and the nerve preserves quotients by \emph{free} group actions.
\end{proof}
\end{thm}

Let $H$ be a finite group, $\phi\colon H\to G$ a homomorphism, and $A$ a non-empty finite faithful $H$-set. We write $\mathcal F_\phi$\nomenclature[aFphi3]{$\mathcal F_\phi$}{$G$-parsummable version of $\mathscr F_\phi$} for the $G$-parsummable category $\cat{P}(E\Inj(A,\omega)\times_\phi G)$. Analogously to the above, we can conclude from Theorem~\ref{thm:BPQ-SSet-corep}:

\begin{thm}\label{thm:sigma-corep-vs-F-phi-parsum}\index{G-global Barratt-Priddy-Quillen Theorem@$G$-global Barratt-Priddy-Quillen-Theorem!for G-parsummable categories@for $G$-parsummable categories|textbf}
Let $f\colon A\to\omega$ be an $H$-equivariant injection (with respect to the $H$-action transported from $\omega^A$ as before). Then the composite map
\begin{equation*}
\Sigma^\bullet_+(I(A,\blank)\times_\phi G)\xrightarrow{\Sigma^\bullet\tilde f}\Sigma^\bullet\nerve(\mathcal F_\phi)[\omega^\bullet]\xrightarrow{\rho}\cat{K}_{\textup{$G$-gl}}(\mathcal F_\phi)
\end{equation*}
is a $G$-global weak equivalence. Moreover, this map is independent of the choice of $f$ up to simplicial homotopy.\qed
\end{thm}

Together with Example~\ref{ex:suspension-IAphi}, we in particular see that $\cat{K}_{\textup{$G$-gl}}(\mathcal F_\phi)$ corepresents the true zeroth $\phi$-equivariant homotopy group.

\subsection{Versions for symmetric monoidal $\bm G$-categories}
Throughout, let $H$, $\phi$, and $A$ be as above. For our final version of the $G$-global Barratt-Priddy-Quillen Theorem, we want to express $\Sigma^\bullet_+ I(A,\blank)\times_\phi G$ as $G$-global algebraic $K$-theory of certain explicit symmetric monoidal $G$-categories. The crucial observation that allows reducing this to category theoretic considerations is the following:

\begin{prop}\label{prop:F-phi-saturated}
The $G$-parsummable category $\mathcal F_\phi$ is saturated.
\end{prop}

For the proof we will need:

\begin{lemma}\label{lemma:E-Inj-times-phi-sat}
Let $K$ be a finite group, let $B$ be a countable faithful $K$-set, and let $X$ be any $G$-$K$-biset. Then the $E\mathcal M$-$G$-category $E\Inj(B,\omega)\times_KX$ is saturated.
\begin{proof}
Let $L\subset\mathcal M$ universal and let $\theta\colon L\to G$ be any homomorphism. We want to show (see Remark~\ref{rk:saturated-alt}) that $\eta^\theta$ is an equivalence of categories, for which it is enough to show that it is essentially surjective. If we equip $L$ with the right $L$-action via $\ell'.\ell=\ell^{-1}\ell'$, then $L\to\mathcal M,\ell\mapsto\ell^{-1}$ is right-$L$-equivariant; as both sides are free, we conclude as before that $\Fun(E\mathcal M,\forget\mathcal C)^\theta\to\Fun(EL,\forget\mathcal C)^\theta=\Fun^L(EL,\theta^*\forget\mathcal C)$ is an equivalence of categories for all $E\mathcal M$-$G$-categories $\mathcal C$.

Taking $\mathcal C=E\Inj(B,\omega)\times_KX$, it is therefore enough to show that the resulting composition $\Theta\colon(E\Inj(B,\omega)\times_K X)^\theta\to\Fun^L(EL,\theta^*\forget(E\Inj(B,\omega)\times_K X))$ is essentially surjective; we can explicitly describe $\Theta$ on objects as follows: $[u;x]$ is sent to the left $L$-equivariant functor $F_{[u;x]}\colon EL\to \theta^*(\forget E\Inj(B,\omega)\times_K X)$ with $F_{[u;x]}(\ell_0,\dots,\ell_n)=[\ell_0^{-1}u,\dots,\ell_n^{-1}u;x]$ ($n=0,1$).

Now let $F\colon EL\to \theta^*\forget(E\Inj(B,\omega)\times_K X)$ be any $L$-equivariant functor and fix a representative $(u;x)$ of $F(1)$. By $L$-equivariance we then have $F(\ell)=\ell.(F(1))=[u;\theta(\ell).x]$ for any $\ell\in L$. Lemma~\ref{lemma:emg-equiv-hom-sets-general} therefore implies that there is a unique $\kappa(\ell)\in K$ such that $F(\ell,1)=[u.\kappa(\ell)^{-1},u;x]$ and
\begin{equation}\label{eq:emg-equiv-kappa-two-cell}
x.\kappa(\ell)=\theta(\ell).x;
\end{equation}
the usual argument then shows that $\kappa$ is a group homomorphism $L\to K$. Together with $(\ref{eq:emg-equiv-kappa-two-cell})$ we see that $\kappa$ satisfies the assumptions of Corollary~\ref{cor:emg-equiv-group-hom-realization-general}, so that we get a $v\in\Inj(B,\omega)$ with $\ell v=v.\kappa(\ell)$ and such that $[v;x]$ is a $\theta$-fixed point of $E\Inj(A,\omega)\times_K X$. To finish the proof, it suffices to show that the maps $\tau(\ell)\mathrel{:=}[u,v;\theta(\ell).x]$ assemble into an $L$-equivariant natural isomorphism $F_{[v;x]}\cong F$.

For this we first observe that $\tau$ is obviously $L$-equivariant and a levelwise isomorphism, so it only remains to check naturality. As the edges $(\ell,1)$ generate the morphisms of $EL$ (in the groupoid sense), it suffices to prove compatibility with respect to those, i.e.~to check commutativity of the squares
\begin{equation*}
\begin{tikzcd}[column sep=1.2in]
{[v;x]} \arrow[d, "{[\ell^{-1}v,v;x]}"'] \arrow[r, "{\tau(1)=[u,v;x]}"] & {[u;x]}\arrow[d, "{F(\ell,1)=[u.\kappa(\ell)^{-1},u;x]}"]\\
{[\ell^{-1}v;x]} \arrow[r, "{\tau(\ell)=[u,v;\theta(\ell).x]}"'] & {[u,\theta(\ell).x]}
\end{tikzcd}
\end{equation*}
in $E\Inj(B,\omega)\times_K X$. But indeed, the top path through this diagram is given by $[u.\kappa(\ell)^{-1},v;x]$ whereas the lower one evaluates to
\begin{equation*}
[u,v;\underbrace{\theta(\ell).x}_{{}=x.\kappa(\ell)}][\ell^{-1}v,v;x]=[u.\kappa(\ell)^{-1},\underbrace{v.\kappa(\ell)^{-1}}_{{}=\ell^{-1}v};x][\ell^{-1}v,v;x]=[u.\kappa(\ell)^{-1},v;x]
\end{equation*}
as desired. This completes the proof of the lemma.
\end{proof}
\end{lemma}

\begin{proof}[Proof of Proposition~\ref{prop:F-phi-saturated}]
By the same argument as in the simplicial setting, $\cat{P}(E\Inj(A,\omega)\times_\phi G)\cong\coprod_{n\ge 0} E\Inj(\bm n\times A,\omega)\times_{\Sigma_n\wr H} G^n$. As $\Sigma_n\wr H$ acts faithfully on $\bm n\times A$, the previous lemma implies that all summands are saturated. The claim follows as both $\Fun(E\mathcal M,\blank)$ and $(\blank)^\theta$ preserve coproducts.
\end{proof}

\begin{cor}\label{cor:K-G-gl-prime-mu-F-phi}\index{G-global Barratt-Priddy-Quillen Theorem@$G$-global Barratt-Priddy-Quillen-Theorem!for symmetric monoidal G-categories@for symmetric monoidal $G$-categories}
Let $f\colon A\to\omega$ be an $H$-equivariant injection, and set $F\mathrel{:=}\eta[f;1]\in\Fun(E\mathcal M,\forget\mathcal F_\phi)^\phi$. Then the composition
\begin{equation}\label{eq:equivalence-mu-star-F-phi}
\Sigma_+^\bullet\big(I(A,\blank)\times_\phi G\big)\xrightarrow{\Sigma^\bullet\widetilde{\iota\circ F}}
\Sigma^\bullet\nerve\Fun(E\mathcal M,\Gamma(\mu^*\mathcal F_\phi))[\omega^\bullet]\xrightarrow{\rho}\cat{K}_{\textup{$G$-gl}}'(\mu^*\mathcal F_\phi)
\end{equation}
is a $G$-global weak equivalence. Moreover, such an $f$ exists and the resulting map is independent of the choice of $f$ up to simplicial homotopy.
\begin{proof}
We have a commutative diagram
\begin{equation*}
\begin{tikzcd}[cramped]
\Sigma_+^\bullet\big(I(A,\blank)\times_\phi G\big)\arrow[d, equal]\arrow[r,"\Sigma_\bullet\tilde f"] &[8pt] \Sigma^\bullet\nerve(\mathcal F_\phi)[\omega^\bullet]\arrow[d, "\sim"]\arrow[r,"\rho"] &[-4pt] \cat{K}_{\textup{$G$-gl}}(\mathcal F_\phi)\arrow[d, "\sim"]\\
\Sigma_+^\bullet\big(I(A,\blank)\times_\phi G\big)\arrow[r, "\Sigma^\bullet\widetilde{\iota\circ F}\;\;"'] & \Sigma^\bullet\nerve\big(\Fun(E\mathcal M,\Gamma(\mu^*\mathcal F_\phi)(S^0))\big)[\omega^\bullet] \arrow[r, "\rho\;"'] & \cat{K}_{\textup{$G$-gl}}'(\mu^*\mathcal F_\phi)
\end{tikzcd}
\end{equation*}
where the vertical arrows on the right are induced by the $G$-global level weak equivalence $\nerve\digamma(\mathcal F_\phi)\simeq\nerve\Fun(E\mathcal M,\Gamma(\mu^*\mathcal F_\phi))$ from the proof of Theorem~\ref{thm:K-G-gl-comparison}; in particular, the middle vertical arrow is induced by $\nerve\Fun(E\mathcal M,\iota)\circ\eta$. The claim therefore follows from Theorem~\ref{thm:sigma-corep-vs-F-phi-parsum}.
\end{proof}
\end{cor}

\begin{constr}
Let $\mathscr C\in\cat{$\bm G$-SymMonCat}^0$ be arbitrary. Then we write $\rho'$ for the composition
\begin{equation*}
\Sigma^\bullet\nerve\Fun(E\mathcal M,\mathscr C)[\omega^\bullet]\hskip0pt minus 2pt\xrightarrow{\Sigma^\bullet\nerve\Fun(E\mathcal M,\iota)[\omega^\bullet]}\hskip0pt minus 2pt\Sigma^\bullet\nerve\Fun(E\mathcal M,\Gamma(\mathscr C))(S^0)[\omega^\bullet]\hskip0pt minus 2pt
\xrightarrow{\rho}\hskip0pt minus 2pt\cat{K}'_{\textup{$G$-gl}}(\mathscr C).
\end{equation*}
Simply by naturality, $(\ref{eq:equivalence-mu-star-F-phi})$ then agrees with
\begin{equation*}
\Sigma_+^\bullet\big(I(A,\blank)\times_\phi G\big)\xrightarrow{\Sigma^\bullet\widetilde{F}}
\Sigma^\bullet\nerve\Fun(E\mathcal M, \mu^*\mathcal F_\phi)[\omega^\bullet]\xrightarrow{\rho'}\cat{K}_{\textup{$G$-gl}}'(\mu^*\mathcal F_\phi).
\end{equation*}
\end{constr}

Corollary~\ref{cor:K-G-gl-prime-mu-F-phi} is not yet the result we are after as the symmetric monoidal $G$-category $\mu^*\mathcal F_\phi$ is a rather unnatural and also somewhat unwieldy object. In the remainder of this section, we therefore want to discuss two simpler but equivalent symmetric monoidal $G$-categories:

\begin{constr}\nomenclature[aFphi1]{$\mathscr F_\phi$}{a small symmetric monoidal $G$-category, twisted version of the symmetric monoidal category of finite free $H$-sets and $H$-equivariant isomorphisms}
We define a small category $\widehat{\mathscr F}_\phi$ as follows: an object of $\smash{\widehat{\mathscr F}_\phi}$ is a finite free right $H$-set $A$ together with a right $H$-equivariant map $f\colon A\to G$; here $H$ acts on $G$ from the right via $\phi$, i.e.~$g.h\mathrel{:=}g\phi(h)$. To be entirely precise, we again have to restrict to a suitable \emph{set} $\mathscr U$ of such $H$-sets, and for technical reasons that will become apparent soon, we take one that contains the sets $\bm n\times H$ for $n\in\mathbb N$.

A map $(A,f)\to (A',f')$ is a right $H$-equivariant map $\lambda\colon A\to A'$ such that $f=f'\lambda$; composition and identities are defined in the evident way. Put differently, $\widehat{\mathscr F}_\phi$ is the full subcategory of the slice $\cat{Set-$\bm H$}\downarrow G$ of the category of right $H$-sets over $G$, that is spanned by the objects whose source is finite free (and contained in $\mathscr U$). In particular, $\widehat{\mathscr F}_\phi$ has finite coproducts and these are created in $\cat{Set-$\bm H$}$.

The group $G$ acts on the objects of $\widehat{\mathscr F}_\phi$ from the left via postcomposition, i.e.~$g.(A,f)=(A,(g.\blank)\circ f)$, and this extends to a $G$-action on all of $\widehat{\mathscr F}_\phi$ by letting $g$ send a map $\lambda\colon (A,f)\to (A',f')$ to $\lambda\colon (A,(g.\blank)\circ f)\to (A',(g.\blank)\circ f')$.

Any choice of coproducts in $\cat{Set-$\bm H$}$ (that preserves $\mathscr U$) makes ${\widehat{\mathscr F}_\phi}$ into a symmetric monoidal category; by construction, the $G$-action is through strict symmetric monoidal functors, in particular yielding an object of $\cat{$\bm G$-SymMonCat}^0$, and this is independent of choices in the sense that for any two different choices the identity functor admits a unique $G$-equivariant strong symmetric monoidal structure; this is automatically strictly unital as in both cases the tensor unit is the \emph{unique} initial object. In particular, we get a symmetric monoidal $G$-structure on the maximal subgroupoid $\mathscr F_\phi\mathrel{:=}\core(\widehat{\mathscr F}_\phi)$, which is independent of choices in the sense that for any two different choices of coproducts in $\cat{Set-$\bm H$}$ the identity of $\mathscr F_\phi$ admits a canonical $G$-equivariant strictly unital strong symmetric monoidal structure.
\end{constr}

\begin{ex}
If $G=1$, then forgetting the reference map to $G$ defines an equivalence between $\widehat{\mathscr F}_\phi$ and the category of finite free right $H$-sets, and this functor admits a unique strictly unital strong symmetric monoidal structure. In particular, we get a canonical strictly unital strong symmetric monoidal equivalence between $\mathscr F_\phi$ and the category of finite free right $H$-sets and $H$-equivariant isomorphisms (with symmetric monoidal structure via disjoint union).
\end{ex}

In order to compare this to $\mu^*\mathcal F_\phi$, we will need the following more rigid version of the above construction:

\begin{constr}\nomenclature[aFphi2]{$\mathfrak F_\phi$}{a specific small permutative $G$-category equivalent to $\mathscr F_\phi$}
We write $\mathfrak F_\phi$ for the coproduct $\coprod_{n\ge 0}G^n\hq(\Sigma_n\wr H)$ of the action groupoids. Here the left action of $\Sigma_n\wr H$ on $G$ is as follows: $\Sigma_n$ acts naturally on $G^n$ from the right via permutation of the factors, whereas $H$ acts on $G$ from the right via $\phi$. These assemble into a right $\Sigma_n\wr H$-action, which we turn into a left action by passing to inverses, i.e.~$(\sigma;h_1,\dots,h_n).(g_1,\dots,g_n)=(g_{\sigma^{-1}(1)}\phi(h^{-1}_{\sigma^{-1}(1)}),\dots,g_{\sigma^{-1}(n)}\phi(h^{-1}_{\sigma^{-1}(n)}))$.

Explicitly, this means that an object of $\mathfrak F_\phi$ is an $n$-tuple $(g_1,\dots,g_n)$ of elements of $G$ ($n\ge0$), and a map $(g_1,\dots,g_n)\to (g'_1,\dots,g'_n)$ is a $(\sigma;h_1,\dots,h_n)\in\Sigma_n\wr H$ with $(g'_1,\dots,g'_n)=(\sigma;h_1,\dots,h_n).(g_1,\dots,g_n)$ while there are no maps between tuples of different lengths.

The group $G$ acts on $\mathfrak F_\phi$ from the left as follows: the action on objects is given by $g.(g_1,\dots,g_n)=(gg_1,\dots,gg_n')$, and a morphism $(\sigma;h_1,\dots,h_n)\colon (g_1,\dots,g_n)\to (g'_1,\dots,g'_n)$ is sent to the morphism $(\sigma;h_1,\dots,h_n)\colon g.(g_1,\dots,g_n)\to g.(g'_1,\dots,g'_n)$. We omit the easy verification that this is well-defined and functorial.

We make $\mathfrak F_\phi$ into a permutative category as follows: the tensor product is given on objects by $(g_1,\dots,g_n)\otimes (g'_1,\dots,g'_{n'})=(g_1,\dots,g_n,g'_1,\dots,g'_{n'})$, and on morphisms by $(\sigma;h_1,\dots,h_n)\otimes(\sigma';h'_1,\dots,h'_{n'})=(\sigma\times\sigma',h_1,\dots,h_n,h'_1,\dots,h'_{n'})$, where $\sigma\times\sigma'\in\Sigma_{n+n'}$ is the usual block sum. We omit the easy verification that this is well-defined, functorial, strictly associative, and strictly unital with unit the empty tuple $\epsilon$. It is moreover clear that the tensor product is strictly $G$-equivariant. Finally, we define the symmetry isomorphism as
\begin{equation*}
\tau\mathrel{:=}(\chi;1,\dots,1)\colon (g_1,\dots,g_n)\otimes(g'_1,\dots,g'_{n'})\to (g'_1,\dots,g'_{n'})\otimes(g_1,\dots,g_n)
\end{equation*}
where $\chi\in\Sigma_{n+n'}$ is the shuffle moving the first $n$ entries to the end. We omit the easy verification that $\tau$ satisfies the coherence conditions required to make $\mathfrak F_\phi$ into a permutative category (which will also follow automatically from the argument given in Lemma~\ref{lemma:mathfrak-vs-mathscr-F-phi} below). It is then clear, that $\tau$ strictly commutes with the $G$-action, i.e.~$G$ acts by strict symmetric monoidal functors, so that $\mathfrak F_\phi$ becomes an object of $\cat{$\bm G$-PermCat}$.
\end{constr}

\begin{constr}
We now construct a functor $i\colon\mathfrak F_\phi\to\mathscr F_\phi$ on objects as follows: $i(g_1,\dots,g_n)=(\bm n\times H,\tilde g)$, where $\tilde g\colon\bm n\times H\to G$ is the unique right $H$-equivariant map with $\tilde g(k,1)=g_k$. Moreover, if $(\sigma;h_1,\dots,h_n)$ is a morphism $(g_1,\dots,g_n)\to (g_1',\dots,g'_n)$, then $i(\sigma;h_1,\dots,h_n)$ is the isomorphism $\bm n\times H\to\bm n\times H$ given by acting with $(\sigma;h_1,\dots,h_n)$ from the left, i.e.~sending $(k,h)$ to $(\sigma(k),h_kh)$.
\end{constr}

\begin{lemma}\label{lemma:mathfrak-vs-mathscr-F-phi}
The above is well-defined, $G$-equivariant, and an equivalence of categories. Moreover, $i$ admits a preferred $G$-equivariant strictly unital strong symmetric monoidal structure, making it into a morphism in $\cat{$\bm G$-SymMonCat}^0$.
\begin{proof}
Let us first show that $i$ is well-defined and fully faithful. For this we recall that we have a bijection between $\Sigma_n\wr H$ and the right-$H$-equivariant automorphisms of $\bm n\times H$ given by sending $f\in\Sigma_n\wr H$ to the map $f.\blank\colon\bm n\times H\to\bm n\times H$ (for the usual left $(\Sigma_n\wr H)$-action). It therefore suffices to show that $(g_1',\dots,g_n')=(\sigma;h_1,\dots,h_n).(g_1,\dots,g_n)$ if and only if $\tilde g'\circ((\sigma;h_1,\dots,h_n).\blank)=\tilde g$. For `$\Leftarrow$,' plugging in $(k,1)$ shows $g'_{\sigma(k)}\phi(h_k)=\tilde g'(\sigma(k),h_k)=g_k$, i.e.~$(g_1,\dots,g_n)=(g'_1,\dots,g'_n).(\sigma;h_1,\dots,h_n)=(\sigma;h_1,\dots,h_n)^{-1}.(g_1',\dots,g_n')$, and the claim follows by acting on both sides with $(\sigma;h_1,\dots,h_n)$. Conversely, reading the same calculation backwards shows that $\tilde g'\circ((\sigma;h_1,\dots,h_n).\blank)$ and $\tilde g$ agree on the elements $(k,1)\in\bm{n}\times H$ if $(g_1',\dots,g_n')=(\sigma;h_1,\dots,h_n).(g_1,\dots,g_n)$; the claim follows by right $H$-equivariance of both sides.

Next, let us show that $i$ is essentially surjective: if $(A,f)$ is any object of $\mathscr F_\phi$, then we can choose an isomorphism $\lambda\colon\bm n\times H\cong A$ for a suitable $n\ge 0$; clearly, $\lambda$ defines an isomorphism $(\bm n\times H,f\lambda)\to(A,f)$. On the other hand, $f\lambda$ is uniquely described by its values at $(k,1)$ for $k=1,\dots,n$; thus, if we set $g_k\mathrel{:=} f\lambda(k,1)$, then $(\bm n\times H,f\lambda)=i(g_1,\dots,g_n)$. This completes the proof that $i$ is an equivalence of categories.

It is clear that $i$ is strictly $G$-equivariant. To prove that we can make it into a morphism in $\cat{$\bm G$-SymMonCat}^0$ in a preferred way, it suffices to consider one specific choice of coproducts on $\cat{Set-$\bm H$}$. We therefore agree that for all $m,n\ge 0$ the coproduct $(\bm m\times H)\amalg(\bm n\times H)$ should be taken to be $\bm{(m+n)}\times H$ with the inclusions induced by $\bm m\hookrightarrow\bm{m+n}$ and $\bm n\to\bm{m+n},k\mapsto k+m$. With this choices, one easily checks that the full symmetric monoidal $G$-subcategory $\mathscr F'_\phi\subset\mathscr F_\phi$ spanned by the objects of the form $(\bm n\times H,f)$ is actually a \emph{permutative} $G$-category, and that $i$ defines a strict symmetric monoidal functor $\mathfrak F_\phi\to\mathscr F'_\phi\subset\mathscr F_\phi$.
\end{proof}
\end{lemma}

\begin{constr}
Choose for each $n\ge 0$ an injection $f^{(n)}\colon\bm n\times A\to\omega$. We define a functor $j\colon\mathfrak F_\phi\to\mu^*\mathcal F_\phi$ as follows: an object $g_\bullet\mathrel{:=}(g_1,\dots,g_n)$ is sent to $j(g_\bullet)=[f^{(n)};g_\bullet]$ (where we confuse $\mathcal F_\phi$ with $\coprod_{n\ge 0}(E\Inj(\bm n\times A,\omega)\times G^n)/(\Sigma_n\wr H)$ as before); moreover, if $\alpha\in\Sigma_n\wr H$ defines a morphism $g_\bullet\to g'_\bullet$, then $j(\alpha)=[f^{(n)}.\alpha, f^{(n)};g_\bullet]$. Finally, we define for all $g_\bullet=(g_1,\dots,g_n),g'_\bullet=(g_1',\dots,g'_{n'})$
\begin{equation*}
\nabla_{g_\bullet,g'_\bullet}\mathrel{:=}[f^{(n+n')},\mu(f^{(n)},f^{(n')});g_\bullet\otimes g'_\bullet]\colon j(g_\bullet)\otimes j(g'_\bullet)\to j(g_\bullet\otimes g'_\bullet).
\end{equation*}
\end{constr}

\begin{prop}
This is a well-defined morphism in $\cat{$\bm G$-SymMonCat}^0$ and an underlying equivalence of categories.
\begin{proof}
Let $\alpha\in\Sigma_n\wr H$ define a map $g_\bullet\to g'_\bullet$; then Lemma~\ref{lemma:emg-equiv-hom-sets-general} implies that $j(\alpha)$ is indeed a map $j(g_\bullet)\to j(g'_\bullet)$, and that $j$ actually defines a bijection $\Hom_{\mathfrak F_\phi}(g_\bullet,g'_\bullet)\to\Hom_{\mathcal F_\phi}(j(g_\bullet),j(g'_\bullet))$. It is then easy to check that $j$ preserves compositions and identities, making it into a fully faithful functor. To see that it is an equivalence of categories, it is therefore enough to show that it is essentially surjective. But indeed, if $X\mathrel{:=}[f';g_\bullet]$ is any object of $\mathcal F_\phi$, then $[f',f^{(n)};g_\bullet]$ defines an isomorphism $j(g_\bullet)\cong X$.

It is clear that $j$ is strictly $G$-equivariant and that it strictly preserves the tensor unit, so it only remains to show that $\nabla$ is a strictly $G$-equivariant natural transformation compatible with associativity, unitality, and symmetry isomorphisms.

Equivariance of $\nabla$ is again obvious; moreover, one easily checks that it is natural and compatible with the unit and associativity isomorphisms on both sides (as in each of these cases the required equalities already hold on the level of the chosen representatives), so it only remains to show that $\nabla$ is compatible with the symmetry isomorphisms.

By definition, the symmetry isomorphism $j(g_\bullet)\otimes j(g'_\bullet)\to j(g'_\bullet)\otimes j(g_\bullet)$ is given by $[(\mu\circ t)(f^{(n)},f^{(n')}),\mu(f^{(n)},f^{(n')});g_\bullet\otimes g'_\bullet]$. If we write $\chi\in\Sigma_{n+n'}$ for the permutation shuffling the first $n$ entries to the end, then $(\mu\circ t)(f^{(n)},f^{(n')})=\mu(f^{(n')},f^{(n)})\circ(\chi\times\id)=\mu(f^{(n')},f^{(n)}).(\chi;1,\dots,1)$, hence the symmetry isomorphism agrees with
$[\mu(f^{(n')},f^{(n)}).(\chi;1,\dots,1),\mu(f^{(n)},f^{(n')});g_\bullet\otimes g'_\bullet]$. Moreover,
\begin{align*}
\nabla_{g_\bullet',g_\bullet}&=[f^{(n'+n)},\mu(f^{(n')},f^{(n)});g'_\bullet\otimes g_\bullet]\\
&=[f^{(n'+n)},\mu(f^{(n')},f^{(n)});(g_\bullet\otimes g'_\bullet).(\chi;1,\dots,1)^{-1}]\\
&=[f^{(n'+n)}.(\chi;1,\dots,1),\mu(f^{(n')},f^{(n)}).(\chi;1,\dots,1);g_\bullet\otimes g'_\bullet],
\end{align*}
hence $\nabla_{g'_\bullet,g_\bullet}\circ\tau=[f^{(n'+n)}.(\chi;1,\dots,1),\mu(f^{(n)},f^{(n')});g_\bullet\otimes g'_\bullet]$.

On the other hand, straight from the definitions
\begin{align*}
j(\tau)\hskip-.8pt\circ\hskip-.8pt\nabla_{g_\bullet,g'_\bullet}&=[f^{(n+n')}.(\chi;1,\dots,1),f^{(n+n')};g_\bullet\hskip-.8pt\otimes\hskip-.8pt g'_\bullet][f^{(n+n')},\mu(f^{(n)},f^{(n')});g_\bullet\hskip-.8pt\otimes\hskip-.8pt g'_\bullet]\\
&=[f^{(n+n')}.(\chi;1,\dots,1),\mu(f^{(n)},f^{(n')});g_\bullet\otimes g'_\bullet],
\end{align*}
finishing the proof that $\nabla$ makes $j$ into a morphism in $\cat{$\bm G$-SymMonCat}^0$.
\end{proof}
\end{prop}

\begin{rk}
While we will not need this, it is in fact not hard to show that $j$ is independent of the choices of the injections $f^{(n)}$ up to canonical isomorphism; more precisely, for any other such choice $\tilde f^{(n)}$, the maps $[f^{(n)},\tilde f^{(n)},g_\bullet]$ for $(g_1,\dots,g_n)\in\mathfrak F_\phi$ assemble into a $G$-equivariant symmetric monoidal isomorphism.
\end{rk}

Now we can finally prove:

\begin{thm}\index{G-global Barratt-Priddy-Quillen Theorem@$G$-global Barratt-Priddy-Quillen-Theorem!for symmetric monoidal G-categories@for symmetric monoidal $G$-categories|textbf}\label{thm:BPQ-sym-mon}
\begin{enumerate}
\item Let $F\in\Fun(E\mathcal M,\mathfrak F_\phi)^\phi$ such that
\begin{equation*}
F(1.h_2,1.h_1)=\big(h_2h_1^{-1}\colon (\phi(h_1^{-1})) \to (\phi(h_2^{-1}))\big)
\end{equation*}
for all $h_1,h_2\in H$ (where $H$ acts on $\mathcal M$ from the right via its action on $\omega$ as before). Then
\begin{equation}\label{eq:K-G-gl-mathfrak-F-we}
\Sigma^\bullet_+\big(I(A,\blank)\times_\phi G\big)\xrightarrow{\Sigma^\bullet\widetilde{F}}
\Sigma^\bullet\nerve\big(\Fun(E\mathcal M,\mathfrak F_\phi)\big)[\omega^\bullet]\xrightarrow{\rho'}\cat{K}'_{\textup{$G$-gl}}(\mathfrak F_\phi)
\end{equation}
is a $G$-global weak equivalence.
\item Let $F\in\Fun(E\mathcal M,\mathscr F_\phi)^\phi$ such that $F(1.h_2,1.h_1)$ is given by
\begin{equation}\label{eq:mathscr-F-phi-tautological-class-G-glob}
\begin{tikzcd}[column sep=small]
H\arrow[dr, bend right=10pt, "\phi(h_1^{-1}\cdot\blank)"']\arrow[rr, "h_2h_1^{-1}\cdot\blank"] && H\arrow[dl, bend left=10pt, "\phi(h_2^{-1}\cdot\blank)"]\\
&G
\end{tikzcd}
\end{equation}
for all $h_1,h_2\in H$. Then the composition
\begin{equation*}
\Sigma^\bullet_+\big(I(A,\blank)\times_\phi G\big)\xrightarrow{\Sigma^\bullet\widetilde{F}}
\Sigma^\bullet\nerve\big(\Fun(E\mathcal M,\mathscr F_\phi)\big)[\omega^\bullet]\xrightarrow{\rho'}\cat{K}'_{\textup{$G$-gl}}(\mathscr F_\phi)
\end{equation*}
is a $G$-global weak equivalence.
\end{enumerate}
Moreover, in each of these cases such a functor $F$ exists, and for any two such choices the resulting maps are canonically simplicially homotopic.
\begin{proof}
Let us prove the first statement. As $H$ acts faithfully on $\omega$, its right action on $\mathcal M$ is free; in particular, restricting along the map induced by $h\mapsto 1.h$ again induces an equivalence $\Fun(E\mathcal M,\mathscr C)^\phi\to\Fun(EH,\mathscr C)^\phi$ for any $G$-category $\mathscr C$. As the prescribed restriction of $F$ is obviously $\phi$-fixed, we conclude that such an $F$ exists and that any two choices are uniquely isomorphic relative $EH$, yielding a canonical simplicial homotopy between the induced maps.

To prove that $(\ref{eq:K-G-gl-mathfrak-F-we})$ is a $G$-global weak equivalence, we now pick any such $F$, and we choose an injection $\mu\colon\bm2\times\omega\to\omega$ as well as for each $n\ge 0$ an injection $f^{(n)}\colon\bm n\times A\to\omega$ such that $f^{(1)}$ is $H$-equivariant, yielding an underlying equivalence $j\colon\mathfrak F_\phi\to\mu^*\mathcal F_\phi$ in $\cat{$\bm G$-SymMonCat}^0$. If we write $f\colon A\to\omega$ for the composition of $f^{(1)}$ with the canonical bijection $A\cong\bm1\times A$, then one immediately checks (using $H$-equivariance of $f$) that $j\circ F$ agrees with $\eta[f;1]$ on $EH$ for $\eta\colon\mathcal F_\phi\to\Fun(E\mathcal M,\forget\mathcal F_\phi)$ as before. In particular, these are equivariantly isomorphic, yielding a simplicial homotopy filling the left hand square in
\begin{equation*}
\begin{tikzcd}
\Sigma_+^\bullet\big(I(A,\blank)\times_\phi G\big)\arrow[d,equal]\arrow[r, "\Sigma^\bullet\tilde F"] &[1.5em] \Sigma^\bullet\nerve\big(\Fun(E\mathcal M,\mathfrak F_\phi)\big)[\omega^\bullet]\arrow[d, "\simeq"]\arrow[r,"{\rho'}"]&[-.5em]\cat{K}_{\textup{$G$-gl}}'(\mathfrak F_\phi)\arrow[d, "\simeq"]\\
\Sigma_+^\bullet\big(I(A,\blank)\times_\phi G\big)\arrow[r, "\Sigma^\bullet\widetilde{\eta{[f;1]}}"'] & \Sigma^\bullet\nerve\big(\Fun(E\mathcal M,\mu^*\mathcal F_\phi)\big)[\omega^\bullet]\arrow[r, "{\rho'}"'] & \cat{K}'_{\textup{$G$-gl}}(\mu^*\mathcal F_\phi)
\end{tikzcd}
\end{equation*}
where the vertical arrows are induced by $j$. By the previous proposition, the right hand vertical map is a $G$-global weak equivalence, and so is the lower horizontal composite by Corollary~\ref{cor:K-G-gl-prime-mu-F-phi}. Thus, also $(\ref{eq:K-G-gl-mathfrak-F-we})$ is a $G$-global weak equivalence by $2$-out-of-$3$ as desired.

The second statement follows similarly from the first statement together with Lemma~\ref{lemma:mathfrak-vs-mathscr-F-phi}.
\end{proof}
\end{thm}

We also have a `non-group-completed' version of the above theorem:

\begin{thm}
Assume $H\subset\mathcal M$ is a universal subgroup.
\begin{enumerate}
\item Let $F\in\Fun(E\mathcal M,\mathfrak F_\phi)^\phi$ with $F(h_2,h_1)=\big(h_2h_1^{-1}\colon (\phi(h_1^{-1}))\to (\phi(h_2^{-1}))\big)$ for all $h_1,h_2\in H$. Then the map
\begin{equation*}
\Gamma(1^+,\blank)\smashp (E\mathcal M\times_{\phi} G)_+\to\nerve\Fun(E\mathcal M, \Gamma(\mathfrak F_\phi))
\end{equation*}
classifying $F$ is a $G$-global special weak equivalence.
\item Let $F\in\Fun(E\mathcal M,\mathscr F_\phi)^\phi$ such that $F(h_2,h_1)$ is given by the map $(\ref{eq:mathscr-F-phi-tautological-class-G-glob})$ for all $h_1,h_2\in H$. Then the map
\begin{equation*}
\Gamma(1^+,\blank)\smashp (E\mathcal M\times_{\phi} G)_+\to\nerve\Fun(E\mathcal M,\Gamma(\mathscr F_\phi))
\end{equation*}
classifying $F$ is a $G$-global special weak equivalence.
\end{enumerate}
\begin{proof}
We will prove the first statement; the second statement will then again follow via Lemma~\ref{lemma:mathfrak-vs-mathscr-F-phi}. 

Let $A\subset\omega$ be a non-empty finite faithful $H$-subset. By an analogous computation to the one in the previous theorem, it is enough to show that the map
\begin{equation*}
\Gamma(1^+,\blank)\smashp (E\mathcal M\times_{\phi} G)_+\to \digamma\big(\cat{P}(E\Inj(A,\omega)\times_\phi G)\big)\cong\digamma\big(\SP^\infty(E\Inj(A,\omega)\times_\phi G)_+\big)
\end{equation*}
classifying the element $[A\hookrightarrow\omega;1]$ is a $G$-global special weak equivalence. However, the restriction map $E\mathcal M\times_\phi G\to E\Inj(A,\omega)\times_\phi G$ is a $G$-global weak equivalence (Proposition~\ref{prop:restriction-to-faithful}), so the claim follows from Proposition~\ref{prop:BPQ-non-group}.
\end{proof}
\end{thm}

We can now use this to prove a $G$-equivariant Barratt-Priddy-Quillen Theorem for symmetric monoidal $G$-categories.

\begin{prop}\label{prop:equiv-BPQ-v1}
Let $G$ be a finite group, let $H$ be any subgroup, and write $\iota\colon H\hookrightarrow G$ for the inclusion.
\begin{enumerate}
\item Let $F\in\Fun(EG,\mathfrak F_{\iota})^H$ with $F(h_2,h_1)=\big(h_2h_1^{-1}\colon (h_1^{-1})\to (h_2^{-1})\big)$ for all $h_1,h_2\in H$. Then the map
\begin{equation}\label{eq:G-equiv-BPQ-map}
\Sigma^\infty_+G/H\to \cat{K}_G(\mathfrak F_{\iota})
\end{equation}
classifying $F\in \cat{K}_G(\mathfrak F_{\iota})(\varnothing)^H\cong\nerve\Fun(EG,\mathfrak F_\iota)^H$ is a $G$-equivariant weak equivalence.
\item Let $F\in\Fun(EG,\mathscr F_\iota)^H$ with $F(h_2,h_1)$ given by
\begin{equation*}
\begin{tikzcd}[column sep=small]
H\arrow[dr, bend right=10pt, "h_1^{-1}\cdot\blank"']\arrow[rr, "h_2h_1^{-1}\cdot\blank"] && H\arrow[dl, bend left=10pt, "h_2^{-1}\cdot\blank"]\\
& G
\end{tikzcd}
\end{equation*}
for all $h_1,h_2\in H$. Then the map $\Sigma^\infty_+G/H\to\cat{K}_G(\mathscr F_\iota)$ classifying $F$ is a $G$-equivariant weak equivalence.
\end{enumerate}
\begin{proof}
Again, it will be enough to prove the first statement, where we may assume without loss of generality that $G$ is a universal subgroup of $\mathcal M$.

Let $\bar F$ be an extension of $F$ to an $H$-equivariant functor $E\mathcal M\to \mathfrak F_\iota$ (obtained for example by picking an equivariant retraction $\mathcal M\to G$); by the previous theorem, the map $f\colon\Gamma(1^+,\blank)\smashp (E\mathcal M\times_HG)_+\to\nerve\Fun(E\mathcal M,\Gamma(\mathfrak F_\iota))$ classifying $\bar F$ is a $G$-global special weak equivalence, so it in particular induces a $G$-global weak equivalence of $G$-global spectra after applying $(\blank)\langle\mathbb S\rangle=\mathcal E^\otimes\circ(\blank)[\omega^\bullet]$, cf.~the proof of Theorem~\ref{thm:BPQ-group-completed-tame-EM-pointed}. Thus, Lemma~\ref{lemma:comparison-E-s-Gamma} implies that if we view both sides as $\Gamma$-$G$-spaces via the {diagonal $G$-action}, then $f$ induces a $G$-equivariant weak equivalence on $\mathcal E_G$.

On the other hand, the projection $p\colon E\mathcal M\times_HG\to G/H$ is a $G$-global weak equivalence (Proposition~\ref{prop:EM-inj-discrete}), hence in particular a $G$-equivariant weak equivalence with respect to the diagonal actions. It follows that also the diagonal embedding $G/H\to E\mathcal M\times_HG,[g]\mapsto[g;g]$ is a $G$-equivariant weak equivalence as it is right inverse to $p$, so that the induced map $\Gamma(1^+,\blank)\smashp G/H_+\to \Gamma(1^+,\blank)\smashp (E\mathcal M\times_HG)_+$ is a $G$-equivariant level weak equivalence.

Finally, as already observed in the proof of Theorem~\ref{thm:G-global-K-vs-G-equiv} the restriction map $\nerve\Fun(E\mathcal M,\Gamma(\mathfrak F_\iota))\to\nerve\Fun(EG,\Gamma(\mathfrak F_\iota))$ is a $G$-equivariant level weak equivalence. Composing all of the above, we see that applying $\mathcal E_G$ to the map $\Gamma(1^+,\blank)\smashp G/H_+\allowbreak\to\allowbreak\nerve\Fun(EG,\Gamma(\mathfrak F_\iota))$ classifying $\bar F|_{EG}=F$ yields a $G$-equivariant weak equivalence. The claim follows as this in turn agrees with $(\ref{eq:G-equiv-BPQ-map})$ up to isomorphism.
\end{proof}
\end{prop}

Guillou and May \cite[Theorem~9.9]{guillou-may} formulated their original equivariant Barratt-Priddy-Quillen Theorem in terms of a certain $G$-$E_\infty$-algebra whose underlying category is the groupoid core of the category of finite sets over $G/H$ (or more generally any $G$-space). We close this discussion by giving a version of the above result in the spirit of Guillou's and May's theorem:

\begin{constr}
Let $H\subset G$ be finite groups and write $\widehat{\mathscr E}_{G/H}$ for the category of finite sets over $G/H$, viewed as a $G$-category via postcomposition (strictly speaking, we again have to restrict to subsets of a suitable infinite set $\mathscr U$ in order to actually obtain a small category). Picking coproducts of finite sets makes this into an object of $\cat{$\bm G$-SymMonCat}^0$; in particular we can also view its core $\mathscr E_{G/H}$ as an object of $\cat{$\bm G$-SymMonCat}^0$.
\end{constr}

\begin{thm}
Let $H\subset G$ be finite and let $F\in\Fun(EG,\mathscr E_{G/H})^H$ such that its restriction to $EH$ is the constant functor at $[1]\colon *\to G/H$. Then the map $\Sigma^\infty_+G/H\to\cat{K}_G(\mathscr E_{G/H})$ classifying $F$ is a $G$-equivariant weak equivalence.
\begin{proof}
Let $\iota\colon H\hookrightarrow G$ denote the inclusion again. Then the functor $\hat\pi\colon\widehat{\mathscr F}_\iota\to \widehat{\mathscr E}_{G/H}$ given by quotiening out the right $H$-action is strictly $G$-equivariant, it preserves coproducts, and it strictly preserves the (unique) initial object, so it induces a map $\pi\colon \mathscr F_\iota\to \mathscr E_{G/H}$ in $\cat{$\bm G$-SymMonCat}^0$. We claim that $\hat\pi$ and hence also $\pi$ is an underlying equivalence of categories; with this established, the theorem will follow immediately from Proposition~\ref{prop:equiv-BPQ-v1} and the invariance of $\cat{K}_G$ under equivalences of categories.

For essential surjectivity we let $f\colon X\to G/H$ be any map, and we pick for every $x\in X$ a representative $\bar f(x)\in G$ of $f(x)$. Then the image of $X\times H\to G,\allowbreak (x,h)\mapsto \bar f(x)h$ under $\hat\pi$ is clearly isomorphic to $f$.

Next, we will show that for every $f\colon A\to G$ in $\widehat{\mathscr F}_\iota$ the square
\begin{equation}\label{diag:equiv-BPQ-pullback}
\begin{tikzcd}
A\arrow[d, "f"']\arrow[r, two heads] & A/H\arrow[d, "f/H"]\\
G\arrow[r, two heads] & G/H
\end{tikzcd}
\end{equation}
is a pullback in $\cat{Set-$\bm H$}$ (when we equip the sets on the right hand side with the trivial right $H$-action); this will then complete the proof as fullness and faithfulness for maps into $f$ directly correspond to the existence and uniqueness part in the universal property.

To prove that $(\ref{diag:equiv-BPQ-pullback})$ is a pullback, we have to show that there is for every $g\in G$ and $\alpha\in A/H$ with $[g]=(f/H)(\alpha)$ in $G/H$ a unique representative $a\in\alpha$ such that $f(a)=g$. For existence, we first take \emph{any} representative $a'\in\alpha$ and observe that $[g]=[f(a')]$ in $G/H$, i.e.~$g=f(a')h$ for some $h\in H$. Then $a\mathrel{:=}a'h$ is the desired element by $H$-equivariance of $f$. For uniqueness, we let $a_1,a_2\in\alpha$; then $a_2=a_1h$ for some $h\in H$, so $f(a_2)=f(a_1)h$. But if also $f(a_2)=g=f(a_1)$, then we necessarily have $h=1$, i.e.~$a_1=a_2$ as desired.
\end{proof}
\end{thm}

\section[$G$-global algebraic $K$-theory as a quasi-localization]{\texorpdfstring{$\except{toc}{\bm G}\for{toc}{G}$}{G}-global algebraic \texorpdfstring{$\except{toc}{\bm K}\for{toc}{K}$}{K}-theory as a quasi-localization}
In this final section, we will will use almost all of the theory developed above to prove, as our main results in this monograph, various $G$-equivariant, global, and $G$-global versions of Thomason's classical equivalence between symmetric monoidal categories and connective stable homotopy types.

We begin with the parsummable case of Theorem~\ref{introthm:thomason-general} from the introduction:
\index{Thomason Theorem!G-global@$G$-global|seeonly{$G$-global Thomason Theorem}}%
\index{Thomason Theorem!G-equivariant@$G$-equivariant|seeonly{$G$-equivariant Thomason Theorem}}%
\index{Thomason Theorem!global|seeonly{$G$-global Thomason Theorem}}%
\index{global Thomason Theorem|seeonly{$G$-global Thomason Theorem}}%

\begin{thm}\label{thm:thomason-parsummable-spectra}\index{G-global Thomason Theorem@$G$-global Thomason Theorem!for G-parsummable categories@for $G$-parsummable categories|textbf}
The $G$-global algebraic $K$-theory functor $\cat{K}_{\textup{$G$-gl}}$ defines a quasi-localization $\cat{$\bm G$-ParSumCat}\to(\cat{$\bm G$-Spectra}_{\textup{$G$-global}}^{\ge0})^\infty$.
\end{thm}

For $G=1$, this in particular shows that global algebraic $K$-theory defines a quasi-localization $\cat{K}_{\textup{gl}}\colon\cat{ParSumCat}\to(\cat{Spectra}_{\textup{global}}^{\ge0})^\infty$, proving a conjecture of Schwede \cite[p.~1331]{schwede-k-theory} and one half of Theorem~\ref{introthm:thomason-global} from the introduction.

The only missing ingredient for the proof of the above theorem is the following comparison, which we proved as \cite[Theorem~5.8]{sym-mon-global}:

\begin{thm}\index{G-parsummable category@$G$-parsummable category!vs G-parsummable simplicial sets@vs.~$G$-parsummable simplicial sets|textbf}\index{G-parsummable simplicial set@$G$-parsummable simplicial set!vs G-parsummable categories@vs.~$G$-parsummable categories|textbf}
The nerve $\nerve\colon\cat{$\bm G$-ParSumCat}\to\cat{$\bm G$-ParSumSSet}$ descends to an equivalence of the quasi-localizations at the $G$-global weak equivalences.\qed
\end{thm}

\begin{proof}[Proof of Theorem~\ref{thm:thomason-parsummable-spectra}]
At this point, this only amounts to collecting results we proved above. We first recall from the proof of Proposition~\ref{prop:K-G-gl-homotopical} that $\cat{K}_{\textup{$G$-gl}}$ agrees up to isomorphism with the composition
\begin{align*}
\cat{$\bm G$-ParSumCat}&\xrightarrow{\nerve}\cat{$\bm G$-ParSumSSet}\xrightarrow{\digamma}\cat{$\bm\Gamma$-$\bm{E\mathcal M}$-$\bm G$-SSet}_*^\tau\hookrightarrow\cat{$\bm\Gamma$-$\bm{E\mathcal M}$-$\bm G$-SSet}_*\\
&\xrightarrow{(\blank)[\omega^\bullet]}\cat{$\bm\Gamma$-$\bm G$-$\bm{\mathcal I}$-SSet}_*\xrightarrow{\mathcal E^\otimes}\cat{$\bm G$-Spectra}^{\ge0};
\end{align*}
moreover, if we equip these categories with the $G$-global weak equivalences or $G$-global level weak equivalences, then all of the above functors are homotopical.

By the previous theorem, the first of these induces an equivalence on quasi-localizations; moreover, we have seen in Theorem~\ref{thm:gamma-vs-uc} that $\digamma$ induces an equivalence $\cat{$\bm G$-ParSumSSet}^\infty\to(\cat{$\bm\Gamma$-$\bm{E\mathcal M}$-$\bm G$-SSet}_*^{\tau,\textup{special}})^\infty$, while the next two functors induce equivalences between the respective quasi-categories of special $G$-global $\Gamma$-spaces by Corollary~\ref{cor:comparison-special}. Finally, Corollary~\ref{cor:E-otimes-special} shows that $\mathcal E^\otimes$ induces a Bousfield localization $(\cat{$\bm\Gamma$-$\bm G$-$\bm{\mathcal I}$-SSet}_*^{\textup{special}})^\infty\to(\cat{$\bm G$-Spectra}_{\textup{$G$-global}}^{\ge0})^\infty$.
\end{proof}

In fact, the above argument also shows:

\begin{thm}\label{thm:G-global-Thomason-non-group-completed}\index{G-global Thomason Theorem@$G$-global Thomason Theorem!for G-parsummable categories@for $G$-parsummable categories!non-group-completed|textbf}
The functor $\cat{$\bm G$-ParSumCat}\to\cat{$\bm\Gamma$-$\bm{E\mathcal M}$-$\bm G$-SSet}_*$ from the definition of $\cat{K}_{\textup{$G$-gl}}$ descends to an equivalence between the quasi-categories of $G$-par\-sum\-mable categories (with respect to the $G$-global weak equivalences) and the \emph{special} $G$-global $\Gamma$-spaces (with respect to the $G$-global level weak equivalences).\qed
\end{thm}

We can view this as a `non-group-completed' version of the $G$-global Thomason Theorem in the spirit of Mandell's non-equivariant comparison between small permutative categories and special $\Gamma$-spaces \cite[Theorem~1.4]{mandell}.

\begin{rk}
In view of the other comparisons established in Chapter~\ref{chapter:coherent}, this also yields equivalences between $G$-parsummable categories on the one hand and $G$-ultra-commutative monoids or any of the other flavours of special $G$-global $\Gamma$-spaces discussed above on the other hand.
\end{rk}

Next, we will establish the corresponding results for the $G$-global algebraic $K$-theory of permutative (or symmetric monoidal) categories with $G$-action. For this we will need the following comparison which we proved as \cite[Theorem~6.7]{sym-mon-global}:

\begin{thm}\label{thm:Phi-sat-equivalence}\index{G-parsummable category@$G$-parsummable category!vs permutative G-categories@vs.~permutative $G$-categories|textbf}\index{permutative G-category@permutative $G$-category!vs G-parsummable categories@vs.~$G$-parsummable categories|textbf}
The composition
\begin{equation*}
\cat{$\bm G$-PermCat}\xrightarrow\Phi\cat{$\bm G$-ParSumCat}\xrightarrow{(\blank)^{\textup{sat}}}\cat{$\bm G$-ParSumCat}
\end{equation*}
descends to an equivalence of the quasi-localizations at the $G$-global weak equivalences.\qed
\end{thm}

We now immediately obtain the remaining half of Theorem~\ref{introthm:thomason-general} (and hence also of Theorem~\ref{introthm:thomason-global}) from the introduction:

\begin{thm}\label{thm:thomason-SymMon-Spectra}\index{G-global Thomason Theorem@$G$-global Thomason Theorem!for permutative G-categories@for permutative $G$-categories|textbf}
The functors $\cat{K}_{\textup{$G$-gl}}\colon\cat{$\bm G$-PermCat}\to(\cat{$\bm G$-Spectra}_{\textup{$G$-global}}^{\ge0})^\infty$ and $\cat{K}'_{\textup{$G$-gl}}\colon\cat{$\bm G$-SymMonCat}^0\to(\cat{$\bm G$-Spectra}_{\textup{$G$-global}}^{\ge0})^\infty$ are quasi-localizations.
\end{thm}

Together with Lemma~\ref{lemma:g-perm-0-g-sym-mon} this in particular allows us to express connective stable $G$-global homotopy theory as quasi-localization of $\cat{$\bm G$-SymMonCat}$.

\begin{proof}
The first statement follows from the previous theorem together with Theorem~\ref{thm:thomason-parsummable-spectra}. Moreover, as $\cat{$\bm G$-PermCat}\hookrightarrow\cat{$\bm G$-SymMonCat}^0$ is a homotopy equivalence with respect to the $G$-global weak equivalences, the second statement now follows from Theorem~\ref{thm:K-G-gl-comparison}.
\end{proof}


\begin{thm}\index{G-global Thomason Theorem@$G$-global Thomason Theorem!for permutative G-categories@for permutative $G$-categories!non-group-completed|textbf}
The composition
\begin{equation*}
\cat{$\bm G$-SymMonCat}^0\hskip0pt minus 1pt\xrightarrow{\Gamma}\hskip0pt minus 1pt\cat{$\bm\Gamma$\kern-.5pt-$\bm G$-Cat}_*\hskip0pt minus 1pt\xrightarrow{\!\Fun(E\mathcal M,\blank)}\hskip0pt minus 1pt\cat{$\bm\Gamma$\kern-.5pt-$\bm{E\mathcal M}$-$\bm G$-Cat}_*\hskip0pt minus 1pt\xrightarrow{\nerve}\hskip0pt minus 1pt\cat{$\bm\Gamma$\kern-.5pt-$\bm{E\mathcal M}$-$\bm G$-SSet}_*
\end{equation*}
defines an equivalence
$(\cat{$\bm G$-SymMonCat}^0)^\infty_{\textup{$G$-gl.}}\simeq (\cat{$\bm\Gamma$-$\bm{E\mathcal M}$-$\bm G$-SSet}_*^{\textup{special}})^\infty_{\textup{$G$-gl.~level}}$.
\begin{proof}
As before, it suffices to prove this after restricting to $\cat{$\bm G$-PermCat}$. As we have seen in the proof of Theorem~\ref{thm:K-G-gl-comparison}, the resulting functor factors up to homotopy as
\begin{equation*}
\cat{$\bm G$-PermCat}\xrightarrow{(\blank)^{\textup{sat}}\circ\Phi}\cat{$\bm G$-ParSumCat}\to \cat{$\bm\Gamma$-$\bm{E\mathcal M}$-$\bm G$-SSet}_*^{\textup{special}}
\end{equation*}
where the second functor is the composite discussed before. The claim therefore follows from Theorem~\ref{thm:Phi-sat-equivalence} together with Theorem~\ref{thm:G-global-Thomason-non-group-completed}.
\end{proof}
\end{thm}
\nobreak
Again, this in particular yields the corresponding statements for the functor arising in Schwede's construction of global algebraic $K$-theory.\goodbreak

We close this discussion by establishing the analogous results for $G$-equivariant algebraic $K$-theory, in particular proving Theorem~\ref{introthm:thomason-equivariant} from the introduction.

\begin{thm}\label{thm:G-equivariant-Thomason-non-group-completed}\index{G-equivariant Thomason Theorem@$G$-equivariant Thomason Theorem!for permutative G-categories@for permutative $G$-categories!non-group-completed|textbf}
Assume $G$ is finite. The composition
\begin{equation}\label{eq:thomason-equivariant-ngc}
\cat{$\bm G$-SymMonCat}^0\xrightarrow{\Gamma}\cat{$\bm\Gamma$-$\bm G$-Cat}_*\xrightarrow{\Fun(EG,\blank)}\cat{$\bm\Gamma$-$\bm G$-Cat}_*\xrightarrow{\nerve}\cat{$\bm\Gamma$-$\bm G$-SSet}_*
\end{equation}
yields a quasi-localization
$\cat{$\bm G$-SymMonCat}^0\to (\cat{$\bm\Gamma$-$\bm G$-SSet}_*^{\textup{special}})^\infty_{\textup{$G$-equiv.~level}}$.
\begin{proof}
By the proof of Theorem~\ref{thm:G-global-K-vs-G-equiv}, $(\ref{eq:thomason-equivariant-ngc})$ agrees up to $G$-equivariant level weak equivalence with the composition
\begin{align*}
\cat{$\bm G$-SymMonCat}^0&\to\cat{$\bm\Gamma$-$\bm{E\mathcal M}$-$\bm G$-SSet}_*^{\textup{special}}\xrightarrow{(\blank)[\omega^\bullet]}\cat{$\bm\Gamma$-$\bm G$-$\bm{\mathcal I}$-SSet}_*^{\textup{special}}\\
&\xrightarrow{\ev_{\mathcal U_G}}\cat{$\bm\Gamma$-$\bm G$-SSet}_*^{\textup{special}}
\end{align*}
where the unlabelled arrow is the composite from the previous theorem, and $\mathcal U_G$ is any complete $G$-set universe. The claim now follows from the previous theorem together with Corollary~\ref{cor:comparison-special} and Theorem~\ref{thm:special-G-global-Gamma-vs-Gamma-G}.
\end{proof}
\end{thm}

\begin{thm}\index{G-equivariant Thomason Theorem@$G$-equivariant Thomason Theorem!for permutative G-categories@for permutative $G$-categories|textbf}
The equivariant $K$-theory functor
\begin{equation*}
\cat{K}_G\colon\cat{$\bm G$-SymMonCat}^0\to(\cat{$\bm G$-Spectra}_{\textup{$G$-equivariant}}^{\ge0})^\infty
\end{equation*}
is a quasi-localization.
\begin{proof}
This follows from the previous theorem by the usual $G$-equivariant Delooping Theorem (which we recalled as Theorem~\ref{thm:equivariant-delooping}).
\end{proof}
\end{thm}

\begin{rk}\index{permutative G-category@permutative $G$-category!vs genuine permutative G-categories@vs.~genuine permutative $G$-categories}
We emphasize that we are working with the na\"ive notion of symmetric monoidal or permutative categories with $G$-action here, not with the \emph{genuine} permutative $G$-categories in the sense of Guillou and May \cite[Definition~4.5]{guillou-may}. Following Shimakawa, they observed in Proposition~4.6 of \emph{op.~cit.} that the endofunctor $\Fun(EG,\blank)$ of $\cat{$\bm G$-PermCat}$ lifts to a functor into the category of genuine permutative $G$-categories, and they further mentioned that (while expecting other examples to exist) they were not aware of any genuine permutative $G$-category not arising this way.

On the other hand, the theorem above tells us that na\"ive permutative $G$-categories actually model---via a similar construction---all (genuine) $G$-equivariantly coherently commutative monoids. This makes it plausible that the Guillou-May-Shimakawa construction should actually yield an equivalence of homotopy theories between $\cat{$\bm G$-PermCat}$ (with respect to the maps inducing \emph{weak} equivalences on categorical homotopy fixed points) and the genuine permutative $G$-categories (with respect to maps inducing \emph{weak} equivalences on honest fixed points). This is indeed the case, as we prove in \cite{gnpg} building on the above results. In particular, up to $G$-weak equivalence every genuine permutative $G$-category indeed arises via the aforementioned construction from a na\"ive one.
\end{rk}

\appendix
\chapter{Abstract homotopy theory}
In this appendix we recall for easy reference some general results about quasi-categories and model categories that we use in the main text.

\section{Quasi-localizations} \label{appendix:quasi-loc} While we in most cases use tools from homotopical algebra to prove our statements, we are ultimately interested in comparisons on the level of quasi-categories. The passage from the former to the latter is as usual provided by the following definition:

\begin{defi}
Let $\mathscr C$ be a quasi-category and let $W\subset\mathscr C_1$ be a collection of morphisms. A functor $\gamma\colon\mathscr C\to\mathscr D$ of quasi-categories is called a \emph{quasi-localization} at $W$ if it has the following universal property: for every quasi-category $\mathscr T$, restriction along $\gamma$ induces an equivalence
\begin{equation}\label{eq:defi-quasi-loc}
\Fun(\mathscr D,\mathscr T)\to\Fun^W(\mathscr C,\mathscr T),
\end{equation}
where $\Fun^W(\mathscr C,\mathscr T)\subset\Fun(\mathscr C,\mathscr T)$ is the full subcategory spanned by the \emph{homotopical functors}\index{homotopical functor|textbf}, i.e.~functors sending morphisms in $W$ to equivalences. (We remark that taking $\mathscr T=\mathscr D$ implies that $\gamma$ itself sends morphisms in $W$ to equivalences.)

By common abuse, we will often supress $W$ from notation and call $\gamma\colon \mathscr C\to\mathscr D$ (and in fact, by further abuse also simply $\mathscr D$ itself) `the' quasi-localization of $\mathscr C$.
\end{defi}

\begin{warn}
Lurie \cite[Definition~5.2.7.2]{htt} uses the term `localization' for a functor with a fully faithful right adjoint, for which we use the more classical name `Bousfield localization'\index{Bousfield localization} in this paper; the above terminology is used for example by Joyal in \cite{joyal-quasi-loc}.

While every Bousfield localization is a quasi-localization by \cite[Proposition~5.2.7.12]{htt}, the converse does not hold and in fact most of the quasi-localizations we are interested in are not Bousfield localizations.
\end{warn}

\begin{rk}
In the $1$-categorical situation one often adds the condition that $W$ is a wide subcategory (in which case the pair $(\mathscr C,W)$ is called a \emph{relative category})\index{relative category|textbf} and sometimes also that it satisfies the $2$-out-of-$3$ property (\emph{categories with weak equivalences})\index{category with weak equivalences|textbf} or even $2$-out-of-$6$ (\emph{homotopical categories}).\index{homotopical category|textbf} While in all of our applications $W$ is indeed a wide subcategory, the stronger properties need not always be satisfied.

However, assume $\gamma\colon\mathscr C\to\mathscr D$ is a quasi-localization at some collection $W$. Then by definition any functor sending $W$ to equivalences factors up to equivalence through $\gamma$ and hence it sends (by $2$-out-of-$3$ for equivalences in quasi-categories) more generally all morphisms $f\in\overline W$ to equivalences where $\overline W$ is the collection of morphisms sent to equivalences \emph{by $\gamma$}. It follows that $\gamma$ is also a quasi-localization at any collection $W'$ of morphisms such that $W\subset W'\subset \overline W$. As $\overline W$ is a wide subcategory satisfying the $2$-out-of-$6$ property, this will allow us in some proofs to restrict to this case.
\end{rk}

\begin{rk}\label{rk:quasi-loc-ordinary-loc}
Specializing $\mathscr T$ in $(\textup{\ref{eq:defi-quasi-loc}})$ to nerves of categories and using the enriched adjunction $\h\dashv\nerve$, we see that if $\mathscr C\to\mathscr D$ is a quasi-localization at some collection $W$ of morphisms, then the induced functor $\h\mathscr C\to\h\mathscr D$\nomenclature[ah]{$\h$}{homotopy category of a quasi-category (or general simplicial set)} is a localization at the same collection. In particular, if $\mathscr C$ is an ordinary relative category (which we confuse with its nerve) and $\mathscr C\to\mathscr D$ is a quasi-localization, then $\mathscr C\cong\h\nerve\mathscr C\to\h\mathscr D$ is an ordinary localization (where the isomorphism is the inverse of the counit).
\end{rk}

\subsection{Simplicial localization}
\index{localization!simplicial|seeonly{simplicial localization}}\index{simplicial localization|(}
While we use the above theory mostly as a conceptual way to pass from homotopical algebra to higher category theory, there is also a simplicial version of it (which historically predates the systematic study of higher categories), that is central to some arguments in Chapter~\ref{chapter:unstable}.

For this we recall that a simplicial category $\mathscr C$\index{simplicial category} can be equivalently viewed as a simplicial object in categories that is constant on objects. We will write $\mathscr C_n$ (an ordinary category) for the $n$-simplices of the corresponding simplicial object. For $n=0$ this yields what is usually called the \emph{underlying category} of $\mathscr C$\index{underlying category of a simplicial category}\index{simplicial category!underlying category|seeonly{underlying category of a simplicial category}} and accordingly we will write $\und\mathscr C\mathrel{:=}\mathscr C_0$.

\begin{defi}\label{defi:simplicial-localization}
Let $\mathscr C$ be a simplicial category that is fibrant in the Bergner model structure \cite{bergner-cat-delta},\index{simplicial category!Bergner model structure|seeonly{Bergner model structure!on simplicial categories}}\index{Bergner model structure on simplicial categories} i.e.~all its mapping spaces are Kan complexes. Moreover, let $W$ be any collection of morphisms in $\und\mathscr C$. A simplicially enriched functor $\gamma\colon\mathscr C\to\mathscr D$ is called a \emph{simplicial localization}\index{simplicial localization|textbf} at $W$ if for some (hence any) fibrant replacement $\mathscr D\to\mathscr E$ the induced map $\nerve_\Delta\mathscr C\to\nerve_\Delta\mathscr E$\nomenclature[aNDelta]{$\nerve_\Delta$}{homotopy coherent nerve} on homotopy coherent nerves\index{homotopy coherent nerve} is a quasi-localization at $W$ (considered as a collection of morphisms in $\nerve_\Delta\mathscr C$ in the obvious way).
\end{defi}

\begin{constr}
We refer the reader to \cite[2.1]{dwyer-kan-calculating} for the definition of the \emph{Hammock localization} $L^H$\nomenclature[aLH]{$L^H$}{Hammock localization}\index{relative category!Hammock localization|seeonly{Hammock localization}}\index{localization!Hammock|seeonly{Hammock localization}}\index{Hammock localization|textbf}\index{relative simplicial category!Hammock localization|seeonly{Hammock localization}} of a relative category $(\mathscr C,W)$. We will not need any details about the construction (except in the proof of Proposition~\ref{prop:localization-subcategory}, where the necessary properties will be recalled), but only that it is a strict $1$-functor and that there is a natural map $\mathscr C\to L^H(\mathscr C,W)$.

We also recall from \cite[Remark~2.5]{dwyer-kan-calculating} that this definition can be extended to a simplicial category $\mathscr C$ together with a wide simplicial subcategory $W$ as follows: we define $\bm{L^H}(\mathscr C,W)$ to be the diagonal of the bisimplicial object in categories obtained by applying $L^H$ levelwise to $(\mathscr C_n, W_n)$. We remark that $\bm{L^H}$ becomes a functor by employing functoriality of $L^H$ levelwise, and we get a natural map $\mathscr C\to \bm{L^H}(\mathscr C,W)$ induced from the unenriched situation.
\end{constr}

\begin{defi}
A \emph{relative simplicial category}\index{simplicial category!relative|seeonly{relative simplicial category}}\index{relative simplicial category|textbf} consists of a simplicial category $\mathscr C$ together with a wide simplicial subcategory $W$. We call the pair $(\mathscr C,W)$ \emph{fibrant}\index{relative simplicial category!fibrant|textbf} if both $\mathscr C$ and $W$ are fibrant in the Bergner model structure.
\end{defi}

\begin{thm}[Dwyer \& Kan, Hinich]\label{thm:lh-universal}
Let $(\mathscr C,W)$ be a fibrant relative simplicial category. Then the natural map $\mathscr C\to\bm{L^H}(\mathscr C,W)$ is a simplicial localization at the $1$-morphisms of $W$.
\begin{proof}
This is \cite[Proposition~1.2.1]{dwyer-kan-revisited}, also cf.~\cite[1.1.3]{dwyer-kan-revisited}.
\end{proof}
\end{thm}

Note that this theorem in particular tells us that if $(\mathscr C, W)$ is an ordinary relative category, then $\mathscr C\to L^H(\mathscr C,W)$ is a simplicial localization.

\begin{rk}\label{rk:associated-quasi}
Let $\mathscr C$ be a relative category. Then the above map $\mathscr C\to L^H(\mathscr C)$ is the identity on objects. A standard fibrant replacement in the Bergner model structure proceeds by applying Kan's $\Ex^\infty$-functor to morphism spaces and hence again is the identity on objects, i.e.~we get a simplicial localization $\gamma\colon\mathscr C\to\mathscr D$ with $\mathscr D$ fibrant such that $\gamma$ is the identity on objects.

As the objects of $\nerve\mathscr C$ and $\nerve_\Delta(\mathscr D)$ are canonically identified with the objects of $\mathscr C$ and $\mathscr D$, respectively, we conclude that any relative category $\mathscr C$ admits a quasi-localization $\nerve(\mathscr C)\to\mathscr E$ that is an isomorphism on objects. In this case we can of course just rename the $0$-simplices of $\mathscr E$ so that our quasi-localization is the \emph{identity} on objects. We will call the quasi-localization constructed this way (and by abuse of terminology, also the resulting quasi-category $\mathscr E$) the \emph{associated quasi-category}\index{associated quasi-category} of $\mathscr C$. We write $\mathscr E\mathrel{=:}\mathscr C^\infty$.\nomenclature[aCinfinity]{$\mathscr C^\infty$}{associated quasi-category (quasi-localization that is the identity on objects)}

Of course, insisting on a statement about equality (or already about isomorphism) of the class of objects is `evil' in the sense that it is not invariant under equivalences. However, we think that this is justified as this convention allows us to simplify several statements and it is moreover also close to the way we usually think about quasi-localizations (as higher versions of the ordinary homotopy category, whose standard construction is the identity on objects).
\end{rk}

Dwyer and Kan already proved that for a simplicial model category $\mathscr C$ its full subcategory $\mathscr C^\circ$ of cofibrant-fibrant objects is the simplicial localization of $\und\mathscr C^\circ$ at the weak homotopy equivalences; more precisely, the inclusion $\und\mathscr C^\circ\hookrightarrow\mathscr C^\circ$ is a quasi-localization, see \cite[Proposition~4.8]{dwyer-kan-function}. Their proof in fact gives a general criterion, which we exploit in the main text:

\begin{prop}\label{prop:enrichment-vs-localization}
Let $\mathscr C$ be a fibrant simplicial category and let $W\subset\und\mathscr C$ be a wide subcategory all of whose morphisms are homotopy equivalences in $\mathscr C$.

Assume moreover that for every $n\ge 0$ the map
\begin{equation*}
L^H(\mathscr C_0,W)\to L^H(\mathscr C_n,s^*W)
\end{equation*}
induced by the unique map $s\colon[n]\to[0]$ is a Dwyer-Kan equivalence (i.e.~it induces an equivalence of $\Ho(\cat{SSet})$-enriched homotopy categories). Then $\und\mathscr C\hookrightarrow\mathscr C$ is a simplicial localization at $W$.
\begin{proof}
This is implicit in \cite[Proof of Proposition 4.8]{dwyer-kan-function}, also cf.~\cite[1.4.3 and 1.4.4]{dwyer-kan-revisited}. We begin with the following observation:

\begin{claim*}
The canonical map $\mathscr C\to\bm{L^H}(\mathscr C,W)$ is a Dwyer-Kan equivalence.
\begin{proof}
Let $\bm{L^H}(\mathscr C,W)\to\mathscr E$ be a fibrant replacement. Then Theorem~\ref{thm:lh-universal} implies that the induced map $\nerve_\Delta(\mathscr C)\to\nerve_\Delta(\mathscr E)$ is a quasi-localization at a subclass of the equivalences of $\mathscr C$, hence it is an equivalence of quasi-categories. As $\nerve_\Delta$ reflects weak equivalences between fibrant simplicial categories (as the right half of a Quillen equivalence), we conclude that $\mathscr C\to\mathscr E$ is a Dywer-Kan equivalence. The claim follows from $2$-out-of-$3$.
\end{proof}
\end{claim*}

Looking at the naturality square
\begin{equation*}
\begin{tikzcd}
\und\mathscr C\arrow[d]\arrow[r,hook] & \mathscr C\arrow[d]\\
L^H(\und\mathscr C,W)\arrow[r] & \bm{L^H}(\mathscr C,W)
\end{tikzcd}
\end{equation*}
it therefore suffices by $2$-out-of-$3$ that $L^H(\und\mathscr C,W)\to \bm{L^H}(\mathscr C,W)$ is a Dwyer-Kan equivalence. This map is the identity on objects, so it suffices to prove that it is given on mapping spaces by weak equivalences of simplicial sets.

For this we observe that for objects $X,Y\in\mathscr C$ the mapping space on the right is by definition the diagonal of the bisimplicial set $\Maps_{L^H(\mathscr C_\bullet,W)}(X,Y)$, and the map $\Maps_{L^H(\und\mathscr C,W)}(X,Y)\to\Maps_{\bm{L^H}(\mathscr C,W)}(X,Y)$ in question is the diagonal of the map of bisimplicial sets
\begin{equation*}
\const\Maps_{L^H(\mathscr C_0,W)}(X,Y)\to\Maps_{L^H(\mathscr C_\bullet,W)}(X,Y)
\end{equation*}
induced in degree $n$ by the degeneracy $[n]\to[0]$. The assumption guarantees that this is a levelwise weak equivalence and hence its diagonal is a weak equivalence by the `Diagonal Lemma' from simplicial homotopy theory, see \cite[Proposition~IV.1.7]{goerss}, finishing the proof.
\end{proof}
\end{prop}

\begin{rk}
The same proof works in the case where $W$ is any wide fibrant simplicial subcategory and one replaces $s^*W$ by $W_n$. However, we will only need the above version.
\end{rk}

We recall one of the standard ways to produce Dwyer-Kan equivalences on simplicial localizations, also see~\cite[2.5]{dwyer-kan-function}.

\begin{defi}
Let $F,G\colon\mathscr C\to\mathscr D$ be homotopical simplicial functors of relative simplicial categories. Then $F$ and $G$ are called \emph{homotopic}\index{homotopic functors|textbf} if they can be connected by a zig-zag of simplicially enriched transformations that are at the same time levelwise weak equivalences.
\end{defi}

\begin{defi}\label{defi:homotopy-equivalence}
A homotopical functor $F\colon\mathscr C\to\mathscr D$ of relative simplicially enriched categories is called a \emph{homotopy equivalence}\index{homotopy equivalence!between relative (simplicial) categories|textbf} if there exists a homotopical functor $G\colon\mathscr D\to\mathscr C$ such that $FG$ is homotopic to the identity of $\mathscr D$ and $GF$ is homotopic to the identity of $\mathscr C$.
\end{defi}

\begin{cor}\label{cor:homotopy-equivalence}
Let $F\colon\mathscr C\to\mathscr D$ be a homotopy equivalence of fibrant simplicial categories equipped with fibrant wide subcategories of weak equivalences. Then the induced functor on $\bm{L^H}$ is a Dwyer-Kan equivalence.
\begin{proof}
This follows from the universal property in the same way as in the claim in the proof of Proposition~\ref{prop:enrichment-vs-localization}.

Alternatively, we reduce to the case of ordinary categories as in the proof of Proposition~\ref{prop:enrichment-vs-localization}. We can then enlarge the weak equivalences on both sides so that they satisfy $2$-out-of-$3$, so that all the intermediate functors in the zig-zags are themselves homotopical. The claim then follows by inductively applying \cite[Proposition~3.5]{dwyer-kan-calculating}.
\end{proof}
\end{cor}

Finally we note:

\begin{prop}\label{prop:localization-subcategory}
Let $\mathscr C$ be a model category with functorial factorizations and let $\mathscr B\subset\mathscr C$ be a full subcategory closed under weak equivalences. Then the inclusion $\mathscr B\hookrightarrow\mathscr C$ induces a fully faithful functor on quasi-localizations.
\begin{proof}
It suffices that $L^H(\mathscr B)\to L^H(\mathscr C)$ is a weak equivalence on morphism spaces. For this let us pick $X,Y\in\mathscr B$ arbitrary.

We consider the category $\mathscr M_{X,Y}$ with objects the zig-zags
\begin{equation*}
\begin{tikzcd}
X & \arrow[l, "\sim"'] A \arrow[r] & B & \arrow[l, "\sim"'] Y
\end{tikzcd}
\end{equation*}
and morphisms the commutative diagrams
\begin{equation*}
\begin{tikzcd}[row sep=small]
 & \arrow[dl, "\sim"'] A\arrow[dd, "\sim"] \arrow[r] & B\arrow[dd, "\sim"'] \\
X & & & \arrow[lu, "\sim"'] Y.\arrow[ld, "\sim"]\\
& \arrow[ul, "\sim"] A' \arrow[r] & B
\end{tikzcd}
\end{equation*}
The important observation is that it does not matter whether we form this in $\mathscr B$ or $\mathscr C$ as $\mathscr B$ is closed under weak equivalences. In the terminology of \cite{dwyer-kan-calculating}, the $k$-simplices of $\nerve(\mathscr M_{X,Y})$ are \emph{hammocks} of width $k$ and length $3$ (in $\mathscr B$ or equivalently in $\mathscr C$) between $X$ and $Y$, whereas the $k$-simplices of $\Maps_{L^H(\mathscr B)}(X,Y)$ and $\Maps_{L^H(\mathscr C)}(X,Y)$ are \emph{reduced} hammocks of width $k$ and arbitrary length between $X$ and $Y$ in $\mathscr B$ resp.~$\mathscr C$. We then have a commutative diagram
\begin{equation*}
\begin{tikzcd}
\nerve(\mathscr M_{X,Y})\arrow[d, equal] \arrow[r] & \Maps_{L^H(\mathscr B)}(X,Y)\arrow[d, hook]\\
\nerve(\mathscr M_{X,Y})\arrow[r] & \Maps_{L^H(\mathscr C)}(X,Y)
\end{tikzcd}
\end{equation*}
where the horizontal maps are given by \emph{reduction} (i.e.~iteratively removing identity columns and composing adjacent horizontal arrows in the same direction). The lower map is a weak homotopy equivalence by \cite[7.2]{dwyer-kan-function}; it therefore suffices that the top map is also a weak equivalence, for which we can use the same argument as Dwyer and Kan. Namely, it is enough by \cite[Proposition~6.2-(i)]{dwyer-kan-calculating} together with \cite[Proposition~8.2]{dwyer-kan-calculating} to exhibit subcategories $W_1,W_2$ of the weak homotopy equivalences in $\mathscr B$ satisfying the following conditions:
\begin{enumerate}
\item Any diagram
\begin{equation*}
\begin{tikzcd}
X \arrow[d, "i\in W_1"'] \arrow[r] & Y\\
X'
\end{tikzcd}
\end{equation*}
with $i\in W_1$ can be functorially completed to a square
\begin{equation*}
\begin{tikzcd}
X \arrow[d, "i\in W_1"'] \arrow[r] & Y\arrow[d, "j\in W_1"]\\
X'\arrow[r] & Y'
\end{tikzcd}
\end{equation*}
with $j\in W_1$.
\item Any diagram
\begin{equation*}
\begin{tikzcd}
& Y'\arrow[d, "p\in W_2"]\\
X \arrow[r] & Y
\end{tikzcd}
\end{equation*}
with $p\in W_2$ can be functorially completed to a square
\begin{equation*}
\begin{tikzcd}
X' \arrow[r]\arrow[d, "q\in W_2"'] & Y'\arrow[d, "p\in W_2"]\\
X \arrow[r] & Y
\end{tikzcd}
\end{equation*}
with $q\in W_2$.
\item Every weak equivalence $w$ admits a functorial factorization $w=w_2w_1$ with $w_i\in W_i$, $i=1,2$.
\end{enumerate}
All of these properties are actually formal consequences of the fact that $\mathscr B$ is closed under weak equivalences together with the respective properties for $\mathscr C$, cf.~\cite[Proposition 8.4]{dwyer-kan-calculating}, but as Dwyer and Kan omit most of the (simple) argument for these, we give a full proof here for completeness:

We take $W_1$ to be the acyclic cofibrations and $W_2$ to be the acyclic fibrations. Then invoking the functorial factorizations of $\mathscr C$ in, say, an acyclic cofibration followed by a fibration (which is automatically acyclic by $2$-out-of-$3$), we get for any weak equivalence $w\colon X\to Y$ in $\mathscr B$ a functorial diagram
\begin{equation*}
X\xrightarrow{w_1\in W_1} H\xrightarrow{w_2\in W_2} Y
\end{equation*}
in $\mathscr C$; as $\mathscr B$ is closed under weak equivalences, also $H\in\mathscr B$, and as $\mathscr B$ is full, this yields the desired functorial factorizations, proving Condition $(3)$.

It remains to verify Condition $(1)$, Condition $(2)$ will then follow from duality. For this we simply pass to pushouts in $\mathscr C$; then $j$ is again an acylic cofibration, and hence also $Y'\in\mathscr B$. Functoriality is induced (and in fact, uniquely determined) by the universal property. This finishes the proof.
\end{proof}
\end{prop}
\index{simplicial localization|)}

\subsection{Quasi-localizations and functor categories}
Let $\mathscr A$ and $\mathscr B$ be simplicial categories; it will be convenient to use the exponential notation $\mathscr A^{\mathscr B}\mathrel{:=}\FUN(\mathscr B,\mathscr A)$ for the enriched category of simplicially enriched functors below, and similarly for functor categories of quasi-categories. In this subsection, we want to prove the following model categorical manifestation of the universal property of quasi-localizations:

\begin{thm}[Dwyer \& Kan]\label{thm:simpl-localization-model-cat}\index{simplicial localization}
Let $A,B$ be small and fibrant simplicial categories, let $f\colon A\to B$ be a simplicial localization at some collection $W$ of morphisms, and let $\mathscr C$ be any combinatorial simplicial model category. Then the functor induced by $f^*\colon\mathscr C^B\to\mathscr C^A$ on associated quasi-categories is fully faithful and its essential image consists precisely of those elements of $\mathscr C^A$ that send morphisms in $W$ to weak homotopy equivalences in $\mathscr C$.
\end{thm}

The fibrancy assumptions on $A$ and $B$ are not necessary, but merely an artifact of our methods. In the special case that $\mathscr C=\cat{SSet}$ with the usual Kan-Quillen model structure (which is the only case we use in the main text), a proof without this assumption can be found as \cite[Theorem~2.2]{dwyer-kan-loc-diagrams}. However, the language and setup used in \emph{op.~cit.} differ from ours, and rigorously translating their statement would require quite a bit of additional terminology. Instead, we will give an alternative proof here using a `rigidification' result due to Lurie. This requires some preparations, throughout which we fix the combinatorial simplicial model category $\mathscr C$.

\begin{lemma}\label{lemma:localization-replacement}
Let $A$ be a small simplicial category and equip $\mathscr C^A$ with either the projective or injective model structure. Then each of the inclusions
\begin{equation*}
(\mathscr C^A)^\circ\hookrightarrow (\mathscr C^\circ)^A\hookrightarrow \mathscr C^A
\end{equation*}
induces an equivalence on the quasi-localizations of the underlying categories at the levelwise weak equivalences.
\begin{proof}
This is true for the composition (in even greater generality) by \cite[Proposition~1.3.8]{dwyer-kan-revisited}; we remark that in our case a simpler proof can be given analogously to \cite[7.1]{dwyer-kan-function}. We will use the same strategy to prove that the first map also has the desired property, which is enough to prove the lemma.

The model category $\mathscr C^A$ is combinatorial and hence we can in particular find functorial cofibrant and fibrant replacements, i.e.~a functor $Q\colon\mathscr C^A\to\mathscr C^A$ together with a natural transformation $\pi\colon Q\Rightarrow\id$ and a functor $P\colon\mathscr C^A\to\mathscr C^A$ together with a natural transformation $\iota\colon\id\Rightarrow P$ such that for each $X\in\mathscr C^A$, $\pi_X\colon QX\to X$ is an acyclic fibration with cofibrant source and $\iota_X\colon X\to PX$ is an acyclic cofibration with fibrant target; we remark that we do not make any claim about these functors or transformations being simplicially enriched. We observe that $P$ and $Q$ are homotopical, preserve the underlying categories of $(\mathscr C^A)^\circ$ and $(\mathscr C^\circ)^A$ (the latter because (co)fibrations in $\mathscr C^A$ are in particular levelwise (co)fibrations), and that moreover $PQ$ sends all of $\und(\mathscr C^A)$ to $\und((\mathscr C^A)^\circ)$.

We claim that (the restriction of) $PQ$ is homotopy inverse to the inclusion $\und\big((\mathscr C^A)^\circ\big)\hookrightarrow \und\big((\mathscr C^\circ)^A\big)$. Indeed, we have for each $X\in\mathscr C^A$ a natural zig-zag of levelwise weak equivalences
\begin{equation*}
PQ(X)\xleftarrow{(\iota_{Q\circ X})}Q(X)\xrightarrow{\pi_{X}} X
\end{equation*}
in $\mathscr C^A$. By the above remarks, $Q$ preserves $\und((\mathscr C^\circ)^A)$ (hence the above exhibits $PQ$ as right homotopy inverse to the inclusion) and $\und((\mathscr C^A)^\circ)$ (hence the above exhibits $PQ$ also as left homotopy inverse); this finishes the proof.
\end{proof}
\end{lemma}

\begin{prop}\label{prop:localization-of-diagrams}
Let $A$ be a small fibrant simplicial category. Then the composition
\begin{equation*}
\nerve\big({(\mathscr C^\circ)}^A\big)\hookrightarrow\nerve_\Delta\big({(\mathscr C^\circ)}^A\big) \to\nerve_\Delta(\mathscr C^\circ)^{\nerve_\Delta(A)},
\end{equation*}
where the second map is adjunct to
\begin{equation*}
\nerve_\Delta\big({(\mathscr C^\circ)}^A\big)\times\nerve_\Delta(A)\cong \nerve_\Delta\big({(\mathscr C^\circ)}^A\times A\big) \xrightarrow{\textup{eval}} \nerve_\Delta(\mathscr C^\circ),
\end{equation*}
is a quasi-localization at the levelwise weak equivalences.
\end{prop}

Here we abbreviate $\nerve(\mathscr D)\mathrel{:=}\nerve(\und\mathscr D)$ for any simplicial category $\mathscr D$.

\begin{proof}
Equip $\mathscr C^A$ with either the projective or injective model structure. By the previous lemma we may  restrict to $(\mathscr C^A)^\circ$, and this can then be factored as
\begin{equation*}
\nerve\big((\mathscr C^A)^\circ\big)\hookrightarrow\nerve_\Delta\big((\mathscr C^A)^\circ\big)
\hookrightarrow\nerve_\Delta\big({(\mathscr C^\circ)}^A\big)\to\nerve_\Delta(\mathscr C^\circ)^{\nerve_\Delta(A)}.
\end{equation*}
Here the first map is a quasi-localization at the levelwise weak equivalences by \cite[Proposition~4.8]{dwyer-kan-function}, while the composition of the latter two maps is an equivalence by a special case of \cite[Proposition~4.2.4.4]{htt}, also see \cite[proof of Corollary~4.2.4.7]{htt}. Thus, the whole composition is a quasi-localization at the levelwise weak equivalences, finishing the proof.
\end{proof}

\begin{proof}[Proof of Theorem~\ref{thm:simpl-localization-model-cat}]
Let us say that the essential image of a functor is \emph{created} by a collection $X$ if it consists precisely of the objects equivalent to elements of $X$. The above statement is then equivalent to demanding that the essential image of $f^*$ be created by those diagrams that send morphisms in $W$ to weak equivalences (as this collection is already closed under equivalences by the $2$-out-of-$3$-property for weak equivalences in $\mathscr C$ and as weak equivalences in $\mathscr C^A$ are saturated). This notion has the advantage that it is invariant under equivalences.

We observe that the subcollection of \emph{levelwise cofibrant-fibrant} diagrams inverting $W$ creates the same essential image: indeed, by the proof of the previous proposition, any $X\in\mathscr C^A$ admits a zig-zag of levelwise weak equivalences in $\mathscr C^A$ (and hence of equivalences in the quasi-localization) to one that is levelwise cofibrant-fibrant, and if $X$ inverts $W$, so does this replacement by $2$-out-of-$3$.

Now let us consider the diagram
\begin{equation*}
\begin{tikzcd}
\nerve\big((\mathscr C^\circ)^B\big)\arrow[d, "f^*"'] \arrow[r, hook] & \nerve(\mathscr C^B)\arrow[d, "f^*"] \arrow[r] & \mathscr D\arrow[d]\\
\nerve\big((\mathscr C^\circ)^A\big) \arrow[r, hook] & \nerve(\mathscr C^A) \arrow[r] & \mathscr E
\end{tikzcd}
\end{equation*}
where the right hand horizontal maps are quasi-localizations and the right hand vertical map is induced by $f^*$; the right hand square commutes up to equivalence while the left hand one commutes strictly.

By Lemma~\ref{lemma:localization-replacement} also the horizontal composites are quasi-localizations at the levelwise weak equivalences, and obviously the right hand vertical map still qualifies as induced map. Altogether we see that it suffices to prove: \emph{the functor induced by $f^*\colon(\mathscr C^\circ)^B\to(\mathscr C^\circ)^A$ on quasi-localizations is fully faithful and its essential image is created by those elements of $(\mathscr C^\circ)^A$ that invert $W$} (or, more precisely, their images under the quasi-localization functor).

For this we may pick any model of quasi-localizations, and we choose the one from the previous proposition. We then have a (strictly) commutative diagram
\begin{equation*}
\begin{tikzcd}
\nerve\big((\mathscr C^\circ)^B\big)\arrow[d, "f^*"'] \arrow[r, hook] & \nerve_\Delta(\mathscr C^\circ)^{\nerve_\Delta(B)}\arrow[d, "\nerve_\Delta(f)^*"]\\
\nerve\big((\mathscr C^\circ)^A\big) \arrow[r, hook] & \nerve_\Delta(\mathscr C^\circ)^{\nerve_\Delta(A)}
\end{tikzcd}
\end{equation*}
which allows us to identify the induced map with the obvious restriction. By assumption, this is fully faithful with essential image those $\nerve_\Delta(A)\to\nerve_\Delta(\mathscr C^\circ)$ that invert $W$, which obviously includes all images of diagrams in $(\mathscr C^\circ)^A$ sending $W$ to weak homotopy equivalences (as weak homotopy equivalences in $\mathscr C^\circ$ agree with honest homotopy equivalences). So it only remains to prove the following converse: each such diagram $\nerve_\Delta(A)\to\nerve_\Delta(\mathscr C^\circ)$ is equivalent to the image of an object of $(\mathscr C^\circ)^A$ inverting $W$. But indeed, as quasi-localizations are essentially surjective (see e.g.~Remark~\ref{rk:associated-quasi}), it is equivalent to some $X\in(\mathscr C^\circ)^A$, which by $2$-out-of-$3$ for equivalences in the quasi-category $\nerve_\Delta(\mathscr C^\circ)$ then sends $W$ to (weak) homotopy equivalences as desired.
\end{proof}

\subsection{Derived functors}
Assume we are given a Quillen adjunction
\begin{equation}\label{eq:quillen-adjunction}
F\colon \mathscr C\rightleftarrows\mathscr D :\!G.
\end{equation}
While the above associates quasi-categories to $\mathscr C$ and $\mathscr D$, in most cases $F$ and $G$ are not homotopical, so they a priori do not give rise to functors between them. In the classical situation, one instead has \emph{left} and \emph{right derived functors}, respectively, yielding an adjunction $\cat{L}F\colon\Ho(\mathscr C)\rightleftarrows\Ho(\mathscr D) :\!\cat{R}G$ on the level of homotopy categories. We are interested in the following higher categorical version of this:

\begin{thm}[Mazel-Gee]\label{thm:derived-adjunction}
In the above situation, there is an adjunction\index{derived functor|textbf}\index{left derived functor|seeonly{derived functor}}\index{right derived functor|seeonly{derived functor}}
\begin{equation}\label{eq:adjunction-induced}
\textbf{\textup L}F\colon \mathscr C^\infty\rightleftarrows\mathscr D^\infty :\!\textbf{\textup R}G
\end{equation}
\nomenclature[aLF]{$\cat{L}F$}{left derived functor (on homotopy categories or associated quasi-categories)}%
\nomenclature[aRG]{$\cat{R}G$}{right derived functor (on homotopy categories or associated quasi-categories)}%
of associated quasi-categories that is induced by $(\textup{\ref{eq:quillen-adjunction}})$ in the following sense: the restriction of $\textbf{\textup L}F$ to $\mathscr C^c$ is equivalent (in a preferred way) to the composition
\begin{equation*}
\mathscr C^c\xrightarrow{F}\mathscr D\to \mathscr D^\infty,
\end{equation*}
and dually for the restriction of $\textbf{\textup R}G$ to $\mathscr D^f$.
\begin{proof}
This is \cite[Theorem~2.1]{mazel-gee} and its proof.
\end{proof}
\end{thm}

In analogy with the classical situation, we call $\textbf{\textup L}F$ the \emph{left derived functor} of $F$ and $\textbf{\textup R}G$ the \emph{right derived functor} of $G$. The adjunction $(\textup{\ref{eq:adjunction-induced}})$ is called the \emph{derived adjunction}.

\begin{rk}
By \cite[Proposition~5.2]{dwyer-kan-function} the inclusions $\mathscr C^c\hookrightarrow\mathscr C$ and $\mathscr D^f\hookrightarrow\mathscr D$ induce equivalences on quasi-localizations, also see \cite[Lemma 2.8]{mazel-gee} or \cite[Proposition~1.3.4]{dwyer-kan-revisited}. It follows that the above functors $\textbf{L}F$ and $\textbf{R}G$ are determined up to (canonical) equivalence by their restrictions to $\mathscr C^c$ or $\mathscr D^f$, respectively, prescribed above.

Now assume $F$ is actually homotopical. Then it induces $F^\infty\colon\mathscr C^\infty\to\mathscr D^\infty$ by the universal property of quasi-localizations. But the restriction of $F^\infty$ to $\mathscr C^c$ is by definition equivalent to the restriction of $\textbf{L}F$, so we conclude from the above that $\textbf{L}F\simeq F^\infty$ (in a preferred way) in this case. Dually, we see that $\textbf{R}G\simeq G^\infty$ whenever $G$ is homotopical.
\end{rk}

\begin{rk}
The above characterization together with Remark~\ref{rk:quasi-loc-ordinary-loc} already implies that the functors induced by $\textbf{L}F$ and $\textbf{R}G$ on the homotopy categories are canonically equivalent to the classical derived functors.

However, it is a bit subtle to identify the unit or counit transformation of the quasi-categorical adjunction $(\textup{\ref{eq:adjunction-induced}})$. In the case that $\mathscr C$ and $\mathscr D$ admit functorial factorizations, there are obvious candidates of this, mimicking the construction of the derived adjunction on homotopy categories, and Mazel-Gee sketches a proof that one can use this to get a unit transformation in \cite[A.3.1]{mazel-gee}.

While all the model categories appearing in this monograph admit functorial functorizations, we do not need this identification. In fact, we only care about the adjunction data to exist and whether unit and/or counit are equivalences. The latter can be checked on the level of homotopy categories, where they are conjugate to the unit and counit, respectively, of the classical derived adjunction for purely formal reasons.
\end{rk}

\section{Some homotopical algebra}
\subsection{Bousfield localizations of model categories}
Let $\mathscr C$ be any model category. We recall that a \emph{(left) Bousfield localization}\index{localization!Bousfield|seeonly{Bousfield localization}}\index{Bousfield localization!for model categories|(}\index{left Bousfield localization|seeonly{Bousfield localization}}\index{right Bousfield localization|seeonly{Bousfield localization}} of $\mathscr C$ is a model structure $\mathscr C_{\textup{loc}}$ on the same underlying category with the same cofibrations as $\mathscr C$ and such that each weak equivalence of $\mathscr C$ is also a weak equivalence in $\mathscr C_{\textup{loc}}$. If $W_{\textup{loc}}$ is the collection of weak equivalences in $\mathscr C_{\textup{loc}}$, we also call $\mathscr C_{\textup{loc}}$ the \emph{Bousfield localization of $\mathscr C$ with respect to $W_{\textup{loc}}$}.

It follows directly from the definitions that
\begin{equation}\label{eq:bousfield-characteristic-adjunction}
\id\colon\mathscr C\rightleftarrows\mathscr C_{\textup{loc}} :\!\id
\end{equation}
is a Quillen adjunction with homotopical left adjoint. Note that $(\textup{\ref{eq:bousfield-characteristic-adjunction}})$ indeed induces a Bousfield localization on associated quasi-categories as any adjoint of a quasi-localization is fully faithful \cite[Proposition~7.1.17]{cisinski-book}.

We will now recall two criteria for the existence of Bousfield localizations for combinatorial simplicial model categories. These are based on the following notions:

\begin{defi}
Let $\mathscr C$ be a simplicial model category, and write $\langle X,Y\rangle$ for the mapping spaces in the $\Ho(\cat{SSet})$-enriched homotopy category of $\mathscr C$; explicitly, $\langle X,Y\rangle$ can be computed as the mapping space in $\mathscr C$ between a cofibrant replacement of $X$ and a fibrant replacement of $Y$.

Now let $S$ be any set of maps in $\mathscr C$.
\begin{enumerate}
\item An object $Z\in\mathscr C$ is called \emph{$S$-local}\index{local object|textbf} if $\langle f,Z\rangle$ is an isomorphism in $\Ho(\cat{SSet})$ for all $f\in S$.
\item A map $g\colon X\to Y$ in $\mathscr C$ is called an \emph{$S$-weak equivalence} if $\langle g,Z\rangle$ is an isomorphism in $\Ho(\cat{SSet})$ for all $S$-local $Z\in\mathscr C$.
\end{enumerate}
\end{defi}

We remark that in the second condition we could equivalently ask that the induced map $[g,Z]$ of hom-sets in the \emph{unenriched} homotopy category be bijective (as the $S$-local objects are closed under the cotensoring over $\cat{SSet}$). However, in the definition of $S$-locality the use of the mapping spaces is indeed essential.

\begin{thm}\label{thm:local-model-structure}\index{Bousfield localization!for model categories!at a set of maps|textbf}
Let $\mathscr C$ be a left proper, combinatorial, and simplicial model category, and let $S$ be a set of maps in $\mathscr C$. Then there exists a unique model structure $\mathscr C_{\textup{loc}}$ on $\mathscr C$ with the same cofibrations as $\mathscr C$ and the $S$-weak equivalences as weak equivalences. This model structure is left proper, simplicial, and combinatorial. Moreover, its fibrant objects consist precisely of the fibrant objects of $\mathscr C$ that are in addition $S$-local.
\begin{proof}
The notion of $S$-locality is obviously invariant under changing the elements of $S$ by conjugation with weak equivalences, and hence so is the notion of $S$-weak equivalences. We may therefore assume without loss of generality that $S$ consists of cofibrations of $\mathscr C$, for which the above appears as \cite[Proposition~A.3.7.3]{htt}.
\end{proof}
\end{thm}

\begin{thm}[Lurie]\label{thm:lurie-localization-criterion}\index{Bousfield localization!for model categories!along a Quillen adjunction|textbf}
Let
\begin{equation}\label{eq:lurie-localization-basic-adjunction}
F\colon \mathscr C\rightleftarrows\mathscr D :\!G
\end{equation}
be a simplicial Quillen adjunction of left proper simplicial combinatorial model categories, and assume that the right derived functor $\textbf{\textup R}G$ is fully faithful (say, as a functor between homotopy categories). Then there exists a Bousfield localization $\mathscr C_{\textup{loc}}$ of $\mathscr C$ whose weak equivalences are precisely those maps $f$ such that $(\textbf{\textup L}F)(f)$ is an isomorphism in $\Ho(\mathscr D)$. The resulting model structure is again left proper, and an object is fibrant in $\mathscr C_{\textup{loc}}$ if and only if it is fibrant in $\mathscr C$ and contained in the essential image of $\textbf{\textup R}G$.

Moreover, $(\textup{\ref{eq:lurie-localization-basic-adjunction}})$ defines a simplicial Quillen equivalence
\begin{equation}\label{eq:lurie-localization-induced-equivalence}
F\colon \mathscr C_{\textup{loc}}\rightleftarrows\mathscr D :\!G.
\end{equation}
\begin{proof}
\cite[Corollary~A.3.7.10]{htt} proves all of this except for the characterization of the fibrant objects. We record that the proof of the aforementioned corollary uses \cite[Proposition~A.3.7.3]{htt} (i.e.~the previous theorem) to construct the model structure, so the fibrant objects of $\mathscr C_{\textup{loc}}$ are closed inside the fibrant objects of $\mathscr C$ under the weak equivalences of $\mathscr C$ (this is actually true for {any} Bousfield localization, but this fact does not seem to appear in the literature).

It is obvious that the fibrant objects of $\mathscr C_{\textup{loc}}$ are fibrant in $\mathscr C$. Let us now consider the composition
\begin{equation}\label{eq:right-derived-explicit}
\Ho(\mathscr D)\xrightarrow{\text{fibrant replacement}}\Ho(\mathscr D^f)\xrightarrow{\Ho(G)}\Ho(\mathscr C_{\text{loc}}^f).
\end{equation}
Postcomposing with the equivalence $\Ho(\mathscr C_{\text{loc}}^f)\hookrightarrow\Ho(\mathscr C_{\textup{loc}})$ yields the standard construction of the right derived functor $\textbf{R}G$ associated to the Quillen equivalence $(\textup{\ref{eq:lurie-localization-induced-equivalence}})$, hence we conclude that $(\textup{\ref{eq:right-derived-explicit}})$ is an equivalence of categories. In particular, $\Ho(G)\colon\Ho(\mathscr D^f)\to \Ho(\mathscr C_{\text{loc}}^f)$ is essentially surjective. By Ken Brown's Lemma applied to the Quillen adjunction $\mathscr C\rightleftarrows\mathscr C_{\textup{loc}}$, the weak equivalences between the fibrant objects of $\mathscr C_{\text{loc}}$ are the same as the weak equivalences in $\mathscr C$, so we can find for any fibrant object $X$ of $\mathscr{C}_{\textup{loc}}$ a zig-zag of $\mathscr C$-weak equivalences to some $GY$ with $G\in\mathscr D^f$, i.e.~$X$ is contained in the essential image of $\textbf{R}G\colon\Ho(\mathscr D)\to\Ho(\mathscr C)$.

It remains to prove the converse, i.e.~that any $\mathscr C$-fibrant object $X$ in the essential image of $\textbf{R}G$ is also fibrant in $\mathscr C_{\textup{loc}}$. For this we observe that, since $\Ho(\mathscr C^f)\to\Ho(\mathscr C)$ and $\Ho(\mathscr D^f)\to\Ho(\mathscr D)$ are equivalences, $X$ can be connected by a zig-zag of weak equivalences {inside $\mathscr C^f$} to $GY$ for some $Y\in\mathscr D^f$. The latter is fibrant in $\mathscr C_{\textup{loc}}$ as $(\textup{\ref{eq:lurie-localization-induced-equivalence}})$ is a Quillen adjunction. But $\mathscr C_{\textup{loc}}^f\subset\mathscr C^f$ is closed under $\mathscr C$-weak equivalences as remarked above, hence also $X$ is fibrant in $\mathscr C_{\textup{loc}}$, finishing the proof.
\end{proof}
\end{thm}

In the main text we are usually concerned with model categories $\mathscr C$ in which filtered colimits are homotopical, i.e.~for any small filtered category $I$ the colimit functor $\mathscr C^I\to\mathscr C$ sends levelwise weak equivalences to weak equivalences. The following lemma in particular tells us that this property is preserved under the localization process of the above theorems:

\begin{lemma}\label{lemma:filtered-still-homotopical}
Let $\mathscr C$ be a model category in which filtered colimits are homotopical, and let $\mathscr C_{\textup{loc}}$ be a Bousfield localization that is moreover cofibrantly generated as a model category. Then filtered colimits in $\mathscr C_{\textup{loc}}$ are also homotopical.
\begin{proof}
Let $I$ be a small filtered category. As $\mathscr C_{\textup{loc}}$ is cofibrantly generated, the projective model structure on $(\mathscr C_{\textup{loc}})^I$ exists, and with respect to this the colimit functor is left Quillen.

If now $f\colon X\to Y$ is any weak equivalence in $(\mathscr C_{\textup{loc}})^I$, then we can factor it into an acyclic cofibration $i$ followed by a fibration $p$, which is again acyclic by $2$-out-of-$3$. As the acyclic fibrations in $(\mathscr C_{\textup{loc}})^I$ are defined levelwise, and since the acyclic fibrations of $\mathscr C_{\textup{loc}}$ and $\mathscr C$ agree, $p$ is then in particular levelwise a weak equivalence of $\mathscr C$. The assumption therefore guarantees that $\colim_I p$ is a weak equivalence (in $\mathscr C$ and hence also in $\mathscr C_{\textup{loc}}$). On the other hand, $\colim_I$ being left Quillen implies that $\colim_I i$ is an acyclic cofibration and hence in particular a weak equivalence in $\mathscr C_{\textup{loc}}$. The claim then follows by another application of $2$-out-of-$3$.
\end{proof}
\end{lemma}
\index{Bousfield localization!for model categories|)}

The localization criterion given in Theorem~\ref{thm:lurie-localization-criterion} has the disadvantage that the weak equivalences of $\mathscr C_{\textup{loc}}$ can be hard to grasp as soon as the left adjoint $F$ in $(\textup{\ref{eq:lurie-localization-basic-adjunction}})$ is not homotopical, and of course the situation is even worse for Theorem~\ref{thm:local-model-structure}. On the other hand, in both cases we usually have good control over the fibrant objects, so it will be useful to know how much this tells us about the model structure.

\begin{prop}\label{prop:cofibrations-fibrant-determine}
Let $\mathscr C_0,\mathscr C_1$ be model structures on the same underlying category that have the same cofibrations as well as the same fibrant objects. Then the two model structures actually agree.
\begin{proof}
A proof can be found as \cite[Theorem~15.3.1]{cathtpy}, where this result is attributed to Joyal.
\end{proof}
\end{prop}

There is also a `relative' version of this statement that needs some additional assumptions:

\begin{prop}\label{prop:cofibrations-fibrant-qa}
Let $\mathscr C,\mathscr D$ be simplicial model categories such that $\mathscr D$ is left proper, and let
\begin{equation}\label{eq:wannabe-quillen}
F\colon\mathscr C\rightleftarrows\mathscr D :\!G
\end{equation}
be a simplicial adjunction such that $F$ preserves cofibrations and $G$ preserves fibrant objects. Then $(\textup{\ref{eq:wannabe-quillen}})$ is already a Quillen adjunction.
\begin{proof}
See~\cite[Corollary~A.3.7.2]{htt}.
\end{proof}
\end{prop}

\subsection{Homotopy pushouts}\index{homotopy pushout|(}
While we freely use several basic properties of homotopy pushouts in the main text, this appendix is devoted to proving the following closure property:

\begin{prop}\label{prop:U-sharp-closure}
Let $\mathscr C$ be a left proper model category, let $\mathscr D$ be any cocomplete category, and let $U\colon\mathscr D\to\mathscr C$ be any functor preserving filtered colimits up to weak equivalence, i.e.~for any small filtered category $I$ and any diagram $X_\bullet\colon I\to\mathscr D$, the canonical map $\colim_{i\in I} U(X_i)\to U(\colim_{i\in I}X_i)$ is a weak equivalence. Assume moreover that filtered colimits in $\mathscr C$ are homotopical.

Then the class $\mathscr H$ of morphisms $i\colon A\to B$ in $\mathscr D$ such that $U$ sends all pushouts along $i$ to homotopy pushouts is closed under pushouts, transfinite compositions, and retracts.
\end{prop}

\begin{rk}
If $U$ is the identity, then the above says that in a left proper model category in which filtered colimits are homotopical, the class of maps $i$ such that pushouts along $i$ are homotopy pushouts is stable under pushout, transfinite compositions, and retracts. Because of the left properness assumption these agree with the so-called \emph{flat} maps \cite[Proposition~1.6]{batanin-berger}, i.e.~those maps such that pushouts along pushouts of them preserve weak equivalences. The flat maps are always (i.e.~without the above assumptions) stable under pushout, \emph{finite} composition, and retracts, as appears for example without proof in \cite[Proposition~B.11]{hhr} or \cite[Lemma~1.3]{batanin-berger}.
\end{rk}

\begin{rk}
The corresponding closure properties when one considers ordinary pushouts as opposed to homotopy pushouts in $\mathscr C$ are proven (by a very similar argument to the one we will give below) as \cite[Lemma~3.1]{gcat}.
\end{rk}

To prove the proposition, we will need some preliminary work.

\begin{lemma}\label{lemma:pushout-filtered}
Let $\mathscr C$ be a cocomplete category, let $I$ be a small filtered category with an initial object $\varnothing$, let $X_\bullet,Y_\bullet\colon I\to\mathscr C$ be functors, and let $\tau\colon X_\bullet\Rightarrow Y_\bullet$ be a natural transformation such that for each $i\in I$ the naturality square
\begin{equation*}
\begin{tikzcd}
X_\varnothing\arrow[d, "\tau_\varnothing"']\arrow[r] & X_i\arrow[d, "\tau_i"]\\
Y_\varnothing\arrow[r] & Y_i
\end{tikzcd}
\end{equation*}
is a pushout. Then the induced square
\begin{equation}\label{diag:pushout-filtered}
\begin{tikzcd}
X_\varnothing\arrow[d, "\tau_\varnothing"']\arrow[r] & \colim_{i\in I}X_i\arrow[d, "\colim_{i\in I}\tau_i"]\\
Y_\varnothing\arrow[r] & \colim_{i\in I}Y_i
\end{tikzcd}
\end{equation}
is also a pushout.
\begin{proof}
We decompose $(\textup{\ref{diag:pushout-filtered}})$ as
\begin{equation*}
\begin{tikzcd}
X_\varnothing\arrow[d, "\tau_\varnothing"']\arrow[r] & \colim_{i\in I}X_\varnothing \arrow[r]\arrow[d, "\colim_{i\in I}\tau_\varnothing"] & \colim_{i\in I}X_i\arrow[d, "\colim_{i\in I}\tau_i"]\\
Y_\varnothing\arrow[r] & \colim_{i\in I} Y_\varnothing \arrow[r] & \colim_{i\in I}Y_i.
\end{tikzcd}
\end{equation*}
The right hand square is a pushout as a colimit of pushout squares. But the two left hand horizontal maps are isomorphisms as filtered categories have connected nerve, so also the total rectangle is a pushout as desired.
\end{proof}
\end{lemma}

Here is the homotopy theoretic analogue of the previous result:

\begin{lemma}\label{lemma:homotopy-pushout-filtered}
Let $\mathscr C$ be a left proper model category such that filtered colimits in $\mathscr C$ are homotopical, let $I$ be a small filtered category with an initial object $\varnothing$, let $X_\bullet,Y_\bullet\colon I\to\mathscr C$ be functors, and let $\tau\colon X_\bullet\Rightarrow Y_\bullet$ be a natural transformation such that for each $i\in I$ the naturality square
\begin{equation}\label{diag:homotopy-pushout-base}
\begin{tikzcd}
X_\varnothing\arrow[d, "\tau_\varnothing"']\arrow[r] & X_i\arrow[d, "\tau_i"]\\
Y_\varnothing\arrow[r] & Y_i
\end{tikzcd}
\end{equation}
is a homotopy pushout. Then the induced square
\begin{equation}\label{diag:homotopy-pushout-filtered}
\begin{tikzcd}
X_\varnothing\arrow[d, "\tau_\varnothing"']\arrow[r] & \colim_{i\in I}X_i\arrow[d, "\colim_{i\in I}\tau_i"]\\
Y_\varnothing\arrow[r] & \colim_{i\in I}Y_i
\end{tikzcd}
\end{equation}
is also a homotopy pushout.
\begin{proof}
We choose a factorization $\smash{X_\varnothing\xrightarrow{\kappa_\varnothing} P_\varnothing\xrightarrow{\phi_\varnothing} Y_\varnothing}$ of $\tau_\varnothing\colon X_\varnothing\to Y_\varnothing$ into a cofibration followed by a weak equivalence. Next, we choose for each $i\in I$ an honest pushout
\begin{equation}\label{diag:honest-pushout}
\begin{tikzcd}
X_\varnothing\arrow[d, "\kappa_\varnothing"']\arrow[r] & X_i\arrow[d, "\kappa_i"]\\
P_\varnothing\arrow[r] & P_i
\end{tikzcd}
\end{equation}
in such a way that for $i=\varnothing$ the lower horizontal map $P_\varnothing\to P_\varnothing$ is the identity. By the universal property of pushouts there is a unique way to extend the chosen maps $P_\varnothing\to P_i$ to a functor $P_\bullet\colon I\to\mathscr C$ in such a way that the $\kappa_i$ assemble into a natural transformation $\kappa\colon X_\bullet\Rightarrow P_\bullet$. Moreover, we define $\phi_i\colon P_i\to Y_i$ to be the map induced by the square $(\textup{\ref{diag:homotopy-pushout-base}})$ via the universal property of the pushout $(\textup{\ref{diag:honest-pushout}})$. Then each $\phi_i$ is a weak equivalence by definition of homotopy pushouts in left proper model categories. Moreover, the $\phi_i$ obviously assemble into a natural transformation $\phi\colon P_\bullet\Rightarrow Y_\bullet$ such that $\tau=\phi\kappa$.

Thus we can factor $(\textup{\ref{diag:homotopy-pushout-filtered}})$ as
\begin{equation*}
\begin{tikzcd}
X_\varnothing\arrow[d, "\kappa_\varnothing"']\arrow[r] & \colim_{i\in I}X_i\arrow[d, "\colim_{i\in I}\kappa_i"]\\
P_\varnothing\arrow[d, "\phi_\varnothing"']\arrow[r] & \colim_{i\in I}P_i\arrow[d, "\colim_{i\in I}\phi_i"]\\
Y_\varnothing\arrow[r] & \colim_{i\in I}Y_i.
\end{tikzcd}
\end{equation*}
The top square is a pushout by the previous lemma and the left hand vertical map is a cofibration by construction, so it is a homotopy pushout by left properness of $\mathscr C$. On the other hand, the lower left vertical map is a weak equivalence by construction, whereas the lower right vertical map is a weak equivalence as a filtered colimit of weak equivalences. Thus the total square is a homotopy pushout, finishing the proof.
\end{proof}
\end{lemma}

\begin{proof}[Proof of Proposition~\ref{prop:U-sharp-closure}]
Consider any diagram
\begin{equation*}
\begin{tikzcd}
A\arrow[d, "f"'] \arrow[r] & C\arrow[d]\arrow[r] & E\arrow[d]\\
B\arrow[r] & D\arrow[r] & F
\end{tikzcd}
\end{equation*}
in $\mathscr D$ with $f\in\mathscr H$ and such that both squares are pushouts. Applying $U$ to it then yields by assumption a diagram in $\mathscr C$ such that both the left square as well as the big rectangle are homotopy pushouts. It follows that also the right square is a homotopy pushout, and hence we conclude that $\mathscr H$ is indeed closed under pushouts.

Next, assume that we have a diagram
\begin{equation*}
\begin{tikzcd}
A'\arrow[d, "{f'}"'] \arrow[r, "i"] & A\arrow[d, "f"] \arrow[r, "r"] & A'\arrow[d, "{f'}"]\\
B'\arrow[r, "j"'] & B \arrow[r, "s"'] & B'
\end{tikzcd}
\end{equation*}
exhibiting $f'$ as a retract of some $f\in\mathscr H$ and consider any pushout
\begin{equation}\label{diag:pushout-of-retract}
\begin{tikzcd}
A'\arrow[d, "{\alpha'}"'] \arrow[r, "{f'}"] & B'\arrow[d, "{\beta'}"]\\
C'\arrow[r] & D'.
\end{tikzcd}
\end{equation}
We now construct a commutative cube
\begin{equation*}
\begin{tikzcd}[row sep=small, column sep=small]
& A\arrow[dd]\arrow[rr, "f"] && B\arrow[dd]\\
A'\arrow[dd]\arrow[ur, "i"] \arrow[rr, "{f'}" near end, crossing over] && B'\arrow[ur, "j"']\\
& C\arrow[rr] && D\\
C'\arrow[rr]\arrow[ur, "k"] && D'\arrow[from=uu, crossing over]\arrow[ur, "\ell"']
\end{tikzcd}
\end{equation*}
as follows: starting from
\begin{equation*}
\begin{tikzcd}[row sep=small, column sep=small]
& A\arrow[rr, "f"] && B\\
A'\arrow[dd, "{\alpha'}"']\arrow[ur, "i"] \arrow[rr, "{f'}" near end] && B'\arrow[dd, "{\beta'}"]\arrow[ur, "j"']\\
\phantom{A}\\
C'\arrow[rr] && D'
\end{tikzcd}
\end{equation*}
we form the pushout $C$ of $C'\gets A'\to A$, yielding the left hand square. With this constructed, we can also form the pushout $D$ of $C\gets A\to B$, yielding the back square. Finally, the map $\ell\colon D'\to D$ is induced by $i,j,k$ via the functoriality of pushouts.

On the other hand we also have a commutative cube
\begin{equation*}
\begin{tikzcd}[row sep=small, column sep=small]
& A'\arrow[dd]\arrow[rr, "{f'}"] && B'\arrow[dd]\\
A\arrow[dd]\arrow[ur, "r"] \arrow[rr, "f" near end, crossing over] && B\arrow[ur, "s"']\\
& C'\arrow[rr] && D'.\\
C\arrow[rr]\arrow[ur, "t"] && D\arrow[from=uu, crossing over]\arrow[ur, "u"']
\end{tikzcd}
\end{equation*}
Here $t$ is induced via the universal property of the pushout $C$ by $\alpha'r\colon A\to C'$ and $\id\colon C'\to C'$; in particular, it makes the left hand square commute by construction. Finally, we let $u$ be the map induced by the remaining front-to-back maps via the functoriality of pushouts again.

By assumption we have $ri=\id$ and $sj=\id$ while $tk=\id$ by construction; the universal property of the pushout $D'$ then implies that also $u\ell=\id$. Altogether, we have exhibited the pushout $(\textup{\ref{diag:pushout-of-retract}})$ as retract of a pushout along $f\in\mathscr H$. In particular we see that applying $U$ to $(\textup{\ref{diag:pushout-of-retract}})$ yields a square that is a retract of a homotopy pushout, hence itself a homotopy pushout. We conclude that $\mathscr H$ is indeed closed under retracts.

Finally, we have to prove that $\mathscr H$ is closed under transfinite compositions. Analogously to the argument in the proof of Lemma~\ref{lemma:saturated-objects} this amounts to the following: given any ordinal $\alpha$ and a functor $X_\bullet\colon\{\beta\le\alpha\}\to\mathscr D$ such that
\begin{enumerate}
\item[(A)] For each $\beta<\alpha$, the map $X_\beta\to X_{\beta+1}$ is in $\mathscr H$,
\item[(B)] If $\beta\le\alpha$ is a limit ordinal, then the maps $X_\gamma\to X_\beta$ for $\gamma<\beta$ express $X_\beta$ as (filtered) colimit,
\end{enumerate}
then the induced map $X_0\to X_\alpha$ is in $\mathscr H$ again. So let us consider any pushout
\begin{equation}\label{diag:pushout-transfinite-comp}
\begin{tikzcd}
X_0\arrow[r]\arrow[d, "f"'] & X_\alpha\arrow[d, "g"]\\
C\arrow[r] & D
\end{tikzcd}
\end{equation}
in $\mathscr D$; we have to prove that $U$ sends this to a homotopy pushout.

For this we define $C_0=C$, $\tau_0=f$, and then choose for each $\beta\le\alpha$ a pushout
\begin{equation}\label{diag:transfinite-pushout-intermediate}
\begin{tikzcd}
X_0\arrow[d, "\tau_0"'] \arrow[r] & X_\beta\arrow[d, "\tau_\beta"]\\
C_0\arrow[r] & C_\beta
\end{tikzcd}
\end{equation}
in such a way that for $\beta=0$ the lower horizontal map is the identity and for $\beta=\alpha$ the chosen pushout is given by $(\textup{\ref{diag:pushout-transfinite-comp}})$. By the universal property of pushouts, there is a unique way to extend the chosen maps $C_0\to C_\beta$ as to yield a functor $C_\bullet\colon\{\beta\le\alpha\}\to\mathscr C$ so that the $\tau_\beta$ assemble into a natural transformation $X_\bullet\Rightarrow C_\bullet$. We observe that by Lemma~\ref{lemma:pushout-filtered} together with the uniqueness of colimits, the maps $C_\gamma\to C_\beta$, $\gamma<\beta$ exhibit $C_\beta$ as colimit for every limit ordinal $\beta\le\alpha$.

We will now prove by transfinite induction that the square $(\textup{\ref{diag:transfinite-pushout-intermediate}})$ is sent by $U$ to a homotopy pushout for each $\beta\le\alpha$. Specializing to $\beta=\alpha$ then precisely proves the claim.

For $\beta=0$ both horizontal maps are identity arrows, hence the same holds after applying $U$. In particular, the resulting square is a homotopy pushout. Now assume $\beta>0$ and we have already proven the claim for all $\gamma<\beta$.

If $\beta$ is a successor ordinal, $\beta=\gamma+1$, then we decompose our square into
\begin{equation*}
\begin{tikzcd}
X_0\arrow[d, "\tau_0"'] \arrow[r] & X_\gamma\arrow[d, "\tau_\gamma"] \arrow[r] & X_\beta\arrow[d, "\tau_\beta"]\\
C_0\arrow[r] & C_\gamma \arrow[r] & C_\beta.
\end{tikzcd}
\end{equation*}
Both the left hand square as well as the total rectangle are pushouts by construction, hence so is the right hand square. The left hand square is sent by $U$ to a homotopy pushout by induction hypothesis; on the other hand $(X_\gamma\to X_\beta)\in\mathscr H$ by assumption, hence also the right hand square is sent to a homotopy pushout. But then the whole square is sent to a homotopy pushout, just as desired.

On the other hand, if $\beta$ is a limit ordinal, then we decompose our square into
\begin{equation*}
\begin{tikzcd}
X_0\arrow[d, "f"'] \arrow[r] & \colim_{\gamma<\beta}X_\gamma\arrow[d, "\colim_{\gamma<\beta}\tau_\gamma"] \arrow[r] & X_\beta\arrow[d, "\tau_\beta"]\\
C_0\arrow[r] & \colim_{\gamma<\beta}C_\gamma \arrow[r] & C_\beta.
\end{tikzcd}
\end{equation*}
We have seen above that the two right hand horizontal maps are isomorphisms. Together with the assumption that $U$ preserve filtered colimits up to weak equivalence, we conclude that in the diagram
\begin{equation*}
\begin{tikzcd}
U(X_0)\arrow[d, "U(f)"'] \arrow[r] & \colim_{\gamma<\beta} U(X_\gamma)\arrow[d, "\colim_{\gamma<\beta}U(\tau_\gamma)"] \arrow[r] & U(X_\beta)\arrow[d, "U(\tau_\beta)"]\\
U(C_0)\arrow[r] & \colim_{\gamma<\beta}U(C_\gamma) \arrow[r] & U(C_\beta).
\end{tikzcd}
\end{equation*}
the right hand horizontal maps are weak equivalences. Thus it suffices to prove that the left hand square is a homotopy pushout, which is just an instance of the previous lemma. We conclude that $\mathscr H$ is indeed stable under transfinite composition, which finishes the proof.
\end{proof}
\index{homotopy pushout|)}

\subsection{Transferring model structures}\label{appendix:transfer} \index{transferred model structure|(}
Here we collect several basic facts about transferred model structures that are used in the main text. We begin by recalling the definition once more for convenience:

\begin{defi}\index{transferred model structure|textbf}
Let $\mathscr D$ be a complete and cocomplete category, let $\mathscr C$ be a model category, and let $F\colon\mathscr C\rightleftarrows\mathscr D :\!U$ be an (ordinary) adjunction. The \emph{model structure transferred along $F\dashv U$} on $\mathscr D$ is the (unique if it exists) model structure where a morphism $f$ is a weak equivalence or fibration if and only if $Uf$ is.
\end{defi}

Crans \cite{transferred} and Kan gave useful criteria on when we can transfer a cofibrantly generated model structure in the above sense:

\begin{prop}\label{prop:transfer-criterion}
Let $\mathscr D$ be a complete and cocomplete category, let $\mathscr C$ be a cofibrantly generated model category, $I$ a set of generating cofibrations, $J$ a set of generating acyclic cofibrations, and let
\begin{equation*}
F\colon\mathscr C\rightleftarrows\mathscr D :\!U
\end{equation*}
be an (ordinary) adjunction. Assume the following:
\begin{enumerate}
\item The sets $FI$ and $FJ$ permit the small object argument.
\item Any relative $FJ$-cell complex is a weak equivalence (i.e.~sent by $U$ to a weak equivalence).
\end{enumerate}
Then the transferred model structure exists. Moreover, it is cofibrantly generated with sets of generating cofibrations $FI$ and generating acylic cofibrations $FJ$.
\end{prop}

Here we as usual say that a set $K$ of maps \emph{permits the small object argument}\index{permit the small object argument}\index{transferred model structure!permit the small object argument|seeonly{permit the small object argument}} if the sources of all maps in $K$ are small relative to relative $K$-cell complexes.

\begin{proof}
Crans \cite[Theorem~3.3]{transferred} proved a slightly weaker statement; the above appears as \cite[Theorem~11.3.2]{hirschhorn}, where this is attributed to Kan.
\end{proof}

The following folklore result tells us that several convenient model categorical properties are preserved by the above transfer process:

\begin{lemma}\label{lemma:transferred-properties}
Let $F\colon\mathscr C\rightleftarrows\mathscr D :\!U$ be as above, and assume the transferred model structure on $\mathscr D$ exists.
\begin{enumerate}
\item If $\mathscr C$ is right proper, then so is $\mathscr D$.\label{item:tpr-proper}
\item Assume that $\mathscr C$ is a simplicial model category, $\mathscr D$ is enriched, tensored, and cotensored over $\cat{SSet}$, and that $F\dashv U$ can be enriched to a simplicial adjunction. Then $\mathscr D$ is a simplicial model category.\label{item:tpr-simplicial}
\end{enumerate}
\begin{proof}
We prove the second statement, the argument for the first one being similar but easier. For this we have to verify that for any cofibration $i\colon K\hookrightarrow L$ in $\cat{SSet}$ and any fibration $p\colon X\to Y$ in $\mathscr D$ the map
\begin{equation*}
(i^*,p_*)\colon X^L \to X^K\times_{Y^K}Y^L
\end{equation*}
is a fibration, and that this is acyclic if at least one of $i$ and $p$ is. By definition, this can be detected after applying $U$, and as $U$ is a simplicial right adjoint the resulting map agrees up to conjugation by isomorphisms with
\begin{equation*}
(i^*,U(p)_*)\colon U(X)^L \to U(X)^K\times_{U(Y)^K} U(Y)^L.
\end{equation*}
As $U$ is right Quillen, $U(p)$ is a fibration, acyclic if $p$ is. Hence the claim follows from the fact that $\mathscr C$ is simplicial.
\end{proof}
\end{lemma}

The question whether the transferred model structure is also left proper is more subtle. The following criterion covers all cases relevant to us:

\begin{lemma}\label{lemma:U-pushout-preserve-reflect}\index{homotopy pushout}
Let $\mathscr C$ be a left proper model category, let $\mathscr D$ be any model category, and let $U\colon\mathscr D\to\mathscr C$ be a functor of their underlying categories. Assume that $U$ preserves and reflects weak equivalences, and that it sends pushouts along cofibrations to homotopy pushouts. Then $\mathscr D$ is left proper, and $U$ preserves and reflects homotopy pushouts.
\begin{proof}
Consider any commutative square
\begin{equation}\label{diag:homotopy-po-candidate}
\begin{tikzcd}
A\arrow[d, "f"'] \arrow[r, "i"] & B\arrow[d, "g"]\\
C \arrow[r] & D
\end{tikzcd}
\end{equation}
in $\mathscr D$. If this is a pushout and $i$ is a cofibration, then its image under $U$ is a homotopy pushout by assumption. In particular, if $f$ is in addition a weak equivalence, then so is $Uf$ and hence also $Ug$. As $U$ creates weak equivalence, this means that $g$ is a weak equivalence, proving left properness of $\mathscr D$.

For a general commutative square $(\text{\ref{diag:homotopy-po-candidate}})$, we now factor $i$ into a cofibration $A\to H$ followed by a weak equivalence $H\to B$, yielding a commutative cube
\begin{equation*}
\begin{tikzcd}[row sep=small, column sep=small]
& A\arrow[dd] \arrow[rr] && B\arrow[dd]\\
A\arrow[ur, "="] \arrow[rr, crossing over] && H\arrow[ur, "\sim"']\\
& C\arrow[rr] && D\\
C\arrow[from=uu]\arrow[ur, "="] \arrow[rr] && P\arrow[from=uu, crossing over]\arrow[ur]
\end{tikzcd}
\end{equation*}
where the front face is a pushout. All front-to-back maps are weak equivalences except possibly the lower right one, which is (as $\mathscr D$ is left proper) a weak equivalence if and only if $(\text{\ref{diag:homotopy-po-candidate}})$ is a homotopy pushout.

Applying $U$ to this diagram, we get a commutative cube in $\mathscr C$ whose front face is a homotopy pushout (by assumption) and all of whose front-to-back maps except possibly the lower right one are weak equivalences (as $U$ is homotopical). We conclude that the lower right front-to-back map of this cube is a weak equivalence if and only if also the back square is a homotopy pushout.

If now $(\text{\ref{diag:homotopy-po-candidate}})$ is a homotopy pushout, then $P\to D$ and hence also $U(P\to D)$ is a weak equivalence. We conclude from the above that $U$ preserves homotopy pushouts. On the other hand, if the back square of the cube obtained by applying $U$, i.e.~the result of applying $U$ to $(\text{\ref{diag:homotopy-po-candidate}})$, is a homotopy pushout, then the map $U(P\to D)$ is also a weak equivalence. As $U$ reflects weak equivalences, so is $P\to D$, and we conclude that $U$ reflects homotopy pushouts, finishing the proof.\index{transferred model structure|)}
\end{proof}
\end{lemma}

\subsection{Enlarging the class of cofibrations} At several points, we want to enlarge the cofibrations of a given model structure, in particular to construct several \emph{injective model structures}\index{injective model structure}. In this section we devise two criteria for this, both of them relying on the following characterization of certain combinatorial model structures, which is due to Lurie building on previous work of Jeff Smith. Here the notion of a \emph{perfect} class of morphisms occurs as a black box; the curious reader can find the definition as \cite[Definition~A.2.6.10]{htt}.

\begin{thm}[Lurie]\label{thm:combinatorial-perfect-criterion}
Let $\mathscr C$ be a locally presentable category, let $W$ be a class of morphisms in $\mathscr C$, and let $I$ be a \emph{set} of morphisms in $\mathscr C$. Then the following are equivalent:
\begin{enumerate}
\item There exists a (necessarily unique) left proper combinatorial model structure on $\mathscr C$ with weak equivalences $W$ and generating cofibrations $I$. Moreover, filtered colimits are homotopical with respect to this model structure.
\item All of the following conditions are satisfied:
\begin{enumerate}
\item The class $W$ is perfect.
\item $W$ is closed under pushouts along pushouts of morphisms of $I$.
\item If a morphism $f$ in $\mathscr C$ has the right lifting property against $I$, then $f\in W$.
\end{enumerate}
\end{enumerate}
\begin{proof}
The implication $(1)\Rightarrow(2)$ is \cite[Remark~A.2.6.14]{htt}, whereas $(2)\Rightarrow(1)$ is \cite[Proposition A.2.6.13]{htt}.
\end{proof}
\end{thm}

\begin{cor}\label{cor:enlarge-generating-cof}
Let $\mathscr C$ be a left proper combinatorial model category in which filtered colimits are homotopical. Let $I_1$ be a set of generating cofibrations and let $I_2$ be any other set of morphisms in $\mathscr C$ such that weak equivalences in $\mathscr C$ are closed under pushouts along pushouts of maps in $I_2$.

Then there exists a unique combinatorial model structure on $\mathscr C$ with generating cofibrations $I\mathrel{:=}I_1\cup I_2$ and weak equivalences the weak equivalences of $\mathscr C$. This model structure is left proper, and filtered colimits in it are homotopical. It is right proper if the original one was.
\begin{proof}
Let us prove that the sets $W$ and $I$ satisfy the assumptions for the implication `$\Leftarrow$' of the previous theorem.

We first observe that $W$ is perfect by the other direction of the theorem, verifying Condition~(a). Moreover, weak equivalences are closed under pushouts along pushouts of maps in $I_1$ (by left properness of the original model structure) and closed under pushouts along pushouts of maps in $I_2$ (by assumption), which proves Condition~(b). Finally, as $I\supset I_1$, any map with the right lifting property against $I$ also has the right lifting property against $I_1$ and hence is a weak equivalence, verifying Condition~(c). This proves that the model structure exists, is left proper, and that filtered colimits in it are homotopical.

It only remains to show that the new model structure is right proper provided that the original one was. It is a general (and somewhat surprising) fact, that being right proper only depends on the class of weak equivalences, see~\cite[Proposition~2.5]{rezk-proper}. However, in the present situation there is an easier argument: both model structures have the same weak equivalences and any cofibration in the old model structure is also a cofibration in the new one, so that any new fibration is already a fibration in the old model structure. The claim thus follows immediately from the definition of right properness.
\end{proof}
\end{cor}

As a special case of this we can now often `borrow' the cofibrations of another model structure on the same underlying category:

\begin{cor}\label{cor:mix-model-structures}
Let $\mathscr C_1,\mathscr C_2$ be combinatorial model structures on the same category, such that every cofibration of $\mathscr C_1$ is also a cofibration in $\mathscr C_2$. Assume moreover that the weak equivalences of $\mathscr C_1$ are closed under filtered colimits and under pushouts along cofibrations of $\mathscr C_2$ (!).

Then there exists a unique model structure $\mathscr C_{12}$ on the same category whose weak equivalences are the weak equivalences of $\mathscr C_1$ and whose cofibrations are the cofibrations of $\mathscr C_2$. This model structure is left proper, combinatorial, and filtered colimits in it are homotopical. If $\mathscr C_1$ is right proper, then so is $\mathscr C_{12}$.
\end{cor}

Note that this is not an instance of Cole's Mixing Theorem \cite[Theorem~2.1]{cole-mix} nor of its dual---in particular, we do not require any containment relation between the two classes of weak equivalences. In fact, in our applications in the main text, the weak equivalences of $\mathscr C_1$ (and hence also of $\mathscr C_{12}$) will be strictly finer than the ones of $\mathscr C_2$, whereas Cole's result produces a model structure for the coarser notion of weak equivalence.

\begin{proof}
As $\mathscr C_1$-cofibrations are in particular $\mathscr C_2$-cofibrations, the assumptions guarantee that $\mathscr C_1$ is left proper; moreover, we have explicitly assumed that filtered colimits in it are homotopical.

Now let $I_k$ be any set of generating cofibrations of $\mathscr C_k$, $k=1,2$. As the weak equivalences of $\mathscr C_1$ are stable under pushouts along all cofibrations of $\mathscr C_2$, they are in particular stable under pushouts along pushouts of $I_2$. We can therefore apply the previous corollary to get a left proper combinatorial model structure (right proper if $\mathscr C_1$ is) with the same weak equivalences as $\mathscr C_1$ and with generating cofibrations $I_1\cup I_2$. We omit the easy verification that this is the desired model structure.
\end{proof}

\backmatter
\begingroup
\frenchspacing
\bibliographystyle{amsalpha}
\bibliography{literature}

\newcommand{\etalchar}[1]{$^{#1}$}
\providecommand{\bysame}{\leavevmode\hbox to3em{\hrulefill}\thinspace}
\providecommand{\MR}{\relax\ifhmode\unskip\space\fi MR }
\providecommand{\MRhref}[2]{%
  \href{http://www.ams.org/mathscinet-getitem?mr=#1}{#2}
}
\providecommand{\href}[2]{#2}
\begin{thebibliography}{BMO{\etalchar{+}}15}

\bibitem[AS69]{atiyah-segal}
Michael~F. {Atiyah} and Graeme~B. {Segal}, \emph{{Equivariant $K$-Theory and
  Completion}}, {J. Differ. Geom.} \textbf{3} (1969), 1--18.

\bibitem[{Ati}61]{atiyah-old}
Michael~F. {Atiyah}, \emph{{Characters and Cohomology of Finite Groups}},
  {Publ. Math., Inst. Hautes \'Etud. Sci.} \textbf{9} (1961), 247--288.

\bibitem[{Bar}10]{barwick-tractable}
Clark {Barwick}, \emph{{On Left and Right Model Categories and Left and Right
  Bousfield Localizations}}, {Homology Homotopy Appl.} \textbf{12} (2010),
  no.~2, 245--320.

\bibitem[{Bar}20]{barrero}
Miguel {Barrero Santamaría}, \emph{{Operads in Unstable Global Homotopy
  Theory}}, {Master's Thesis, Rheinische Friedrich-Wilhelms-Universität Bonn},
  2020.

\bibitem[BB17]{batanin-berger}
Michael~A. {Batanin} and Clemens {Berger}, \emph{{Homotopy Theory for Algebras
  Over Polynomial Monads}}, {Theory Appl. Categ.} \textbf{32} (2017), 148--253.

\bibitem[{Ber}07]{bergner-cat-delta}
Julia~E. {Bergner}, \emph{{A Model Category Structure on the Category of
  Simplicial Categories}}, {Trans. Am. Math. Soc.} \textbf{359} (2007), no.~5,
  2043--2058.

\bibitem[BF78]{bousfield-friedlander}
Aldridge~K. {Bousfield} and Eric~M. {Friedlander}, \emph{{Homotopy Theory of
  \(\Gamma\)-Spaces, Spectra, and Bisimplicial Sets}}, {Geom. Appl. Homotopy
  Theory, II, Proc. Conf., Evanston 1977}, {Lect. Notes Math.}, vol. 658, 1978,
  pp.~80--130.

\bibitem[{Blu}06]{blumberg-infinite-loop-spaces}
Andrew~J. {Blumberg}, \emph{{Continuous Functors as a Model for the Equivariant
  Stable Homotopy Category}}, {Algebr. Geom. Topol.} \textbf{6} (2006),
  2257--2295.

\bibitem[BM20]{boavida-moerdijk}
Pedro {Boavida de Brito} and Ieke {Moerdijk}, \emph{{Dendroidal Spaces, \(
  \Gamma \)-Spaces and the Special Barratt-Priddy-Quillen Theorem}}, {J. Reine
  Angew. Math.} \textbf{760} (2020), 229--265.

\bibitem[BMO{\etalchar{+}}15]{gcat}
Anna~Marie {Bohmann}, Kristen {Mazur}, Ang\'elica~M. {Osorno}, Viktoriya
  {Ozornova}, Kate {Ponto}, and Carolyn {Yarnall}, \emph{{A Model Structure on
  \(G{\mathcal C}at\)}}, {Women in Topology: Collaborations in Homotopy Theory.
  WIT: Women in Topology Workshop, Banff International Research Station, Banff,
  Alberta, Canada, August 18--23, 2013}, Providence, RI: American Mathematical
  Society (AMS), 2015, pp.~123--134.

\bibitem[Bro73]{brown-factorization}
Kenneth~S. Brown, \emph{{Abstract Homotopy Theory and Generalized Sheaf
  Cohomology}}, {Trans. Am. Math. Soc.} \textbf{186} (1973), 419--458.

\bibitem[{Car}84]{carlsson-segal}
Gunnar {Carlsson}, \emph{{Equivariant Stable Homotopy and Segal's Burnside Ring
  Conjecture}}, {Ann. Math. (2)} \textbf{120} (1984), 189--224.

\bibitem[{Car}91]{carlsson-sullivan}
\bysame, \emph{{Equivariant Stable Homotopy and Sullivan's Conjecture}},
  {Invent. Math.} \textbf{103} (1991), no.~3, 497--525.

\bibitem[{Car}92]{carlsson-survey}
\bysame, \emph{{A Survey of Equivariant Stable Homotopy Theory}}, {Topology}
  \textbf{31} (1992), no.~1, 1--27.

\bibitem[{Cis}19]{cisinski-book}
Denis-Charles {Cisinski}, \emph{{Higher Categories and Homotopical Algebra}},
  {Camb. Stud. Adv. Math.}, vol. 180, Cambridge: Cambridge University Press,
  2019.

\bibitem[{Col}06]{cole-mix}
Michael {Cole}, \emph{{Mixing Model Structures}}, {Topology Appl.} \textbf{153}
  (2006), no.~7, 1016--1032.

\bibitem[{Cra}95]{transferred}
Sjoerd~E. {Crans}, \emph{{Quillen Closed Model Structures for Sheaves}}, {J.
  Pure Appl. Algebra} \textbf{101} (1995), no.~1, 35--57.

\bibitem[{Day}70]{day-convolution}
Brian {Day}, \emph{{On Closed Categories of Functors}}, {Reports of the Midwest
  Category Seminar IV} (Saunders {MacLane}, ed.), {Lect. Notes Math.}, vol.
  137, Springer Berlin Heidelberg, 1970, pp.~1--38.

\bibitem[DHL{\etalchar{+}}19]{proper-equivariant}
Dieter {Degrijse}, Markus {Hausmann}, Wolfang {L{\"u}ck}, Irakli {Patchkoria},
  and Stefan {Schwede}, \emph{{Proper Equivariant Stable Homotopy Theory}},
  {available as \texttt{arxiv:1908.00779}, to appear in
  \textit{Mem.~Am.~Math.~Soc.}}, 2019.

\bibitem[DK80a]{dwyer-kan-calculating}
William~G. {Dwyer} and Daniel~M. {Kan}, \emph{{Calculating Simplicial
  Localizations}}, {J. Pure Appl. Algebra} \textbf{18} (1980), 17--35.

\bibitem[DK80b]{dwyer-kan-function}
\bysame, \emph{{Function Complexes in Homotopical Algebra}}, {Topology}
  \textbf{19} (1980), 427--440.

\bibitem[DK84]{dwyer-kan-equivariant}
\bysame, \emph{{Singular Functors and Realization Functors}}, {Indag. Math.}
  \textbf{46} (1984), 147--153.

\bibitem[DK87]{dwyer-kan-loc-diagrams}
\bysame, \emph{{Equivalences Between Homotopy Theories of Diagrams}},
  {Algebraic Topology and Algebraic $K$-Theory, Proc. Conf., Princeton, NJ
  (USA)}, {Ann. Math. Stud.}, vol. 113, 1987, pp.~180--205.

\bibitem[{Dus}01]{duskin}
John~W. {Duskin}, \emph{{Simplicial Matrices and the Nerves of Weak
  \(n\)-Categories I: Nerves of Bicategories}}, {Theory Appl. Categ.}
  \textbf{9} (2001), 198--308.

\bibitem[{Elm}83]{elmendorf}
Anthony~D. {Elmendorf}, \emph{{Systems of Fixed Point Sets}}, {Trans. Am. Math.
  Soc.} \textbf{277} (1983), 275--284.

\bibitem[FHM82]{fhm}
Zbigniew {Fiedorowicz}, Henning {Hauschild}, and J.~Peter {May},
  \emph{{Equivariant Algebraic \(K\)-Theory}}, {Algebraic \(K\)-Theory, Proc.
  Conf., Oberwolfach 1980, Part I} (R.~Keith {Dennis}, ed.), {Lect. Notes
  Math.}, vol. 967, Springer Berlin Heidelberg, 1982, pp.~23--80.

\bibitem[FW69]{froehlich-wall}
Albrecht {Fr{\"o}hlich} and Charles~T.C. {Wall}, \emph{{Foundations of
  Equivariant Algebraic $K$-Theory}}, {Algebraic $K$-Theory and its Geometric
  Applications} (Robert~M.F. {Moss} and Charles~B. {Thomas}, eds.), Springer
  Berlin Heidelberg, 1969, pp.~12--27.

\bibitem[GG16]{sym-powers}
Sergey {Gorchinskiy} and Vladimir {Guletski\u{\i}}, \emph{{Symmetric Powers in
  Abstract Homotopy Categories}}, {Adv. Math.} \textbf{292} (2016), 707--754.

\bibitem[GJ99]{goerss}
Paul~G. {Goerss} and John~F. {Jardine}, \emph{{Simplicial Homotopy Theory}},
  {Prog. Math.}, vol. 174, Basel: Birkh\"auser, 1999.

\bibitem[GM17]{guillou-may}
Bertrand~J. {Guillou} and J.~Peter {May}, \emph{{Equivariant Iterated Loop
  Space Theory and Permutative \(G\)-Categories}}, {Algebr. Geom. Topol.}
  \textbf{17} (2017), no.~6, 3259--3339.

\bibitem[{Hau}17]{hausmann-equivariant}
Markus {Hausmann}, \emph{{\(G\)-Symmetric Spectra, Semistability and the
  Multiplicative Norm}}, {J. Pure Appl. Algebra} \textbf{221} (2017), no.~10,
  2582--2632.

\bibitem[{Hau}19]{hausmann-global}
\bysame, \emph{{Symmetric Spectra Model Global Homotopy Theory of Finite
  Groups}}, {Algebr. Geom. Topol.} \textbf{19} (2019), no.~3, 1413--1452.

\bibitem[{Hau}22]{hausmann-formal}
\bysame, \emph{{Global Group Laws and Equivariant Bordism Rings}}, {Ann. Math.
  (2)} \textbf{195} (2022), 841--910.

\bibitem[HHR16]{hhr}
Michael~A. {Hill}, Michael~J. {Hopkins}, and Douglas~C. {Ravenel}, \emph{{On
  the Nonexistence of Elements of Kervaire Invariant One}}, {Ann. Math. (2)}
  \textbf{184} (2016), no.~1, 1--262.

\bibitem[{Hin}16]{dwyer-kan-revisited}
Vladimir {Hinich}, \emph{{Dwyer-Kan Localization Revisited}}, {Homology
  Homotopy Appl.} \textbf{18} (2016), no.~1, 27--48.

\bibitem[{Hir}03]{hirschhorn}
Philip~S. {Hirschhorn}, \emph{{Model Categories and Their Localizations}},
  {Math. Surv. Monogr.}, vol.~99, Providence, RI: American Mathematical Society
  (AMS), 2003.

\bibitem[{Hir}15]{hirschhorn-slice}
\bysame, \emph{{Overcategories and Undercategories of Model Categories}},
  {available as \texttt{arXiv:1507.01642}}, 2015.

\bibitem[{Hov}99]{hovey}
Mark {Hovey}, \emph{{Model Categories}}, {Math. Surv. Monogr.}, vol.~63,
  Providence, RI: American Mathematical Society, 1999.

\bibitem[HSS00]{hss}
Mark {Hovey}, Brooke {Shipley}, and Jeff {Smith}, \emph{{Symmetric Spectra}},
  {J. Am. Math. Soc.} \textbf{13} (2000), no.~1, 149--208.

\bibitem[{Joy}08]{joyal-quasi-loc}
Andr\'e {Joyal}, \emph{{Notes on Quasi-Categories}}, {lecture notes, available
  at \url{www.math.uchicago.edu/~may/IMA/Joyal.pdf}}, 2008.

\bibitem[{Kel}05]{enriched-presheaves}
G.~Maxwell {Kelly}, \emph{{Basic Concepts of Enriched Category Theory}}, {Repr.
  Theory Appl. Categ.} \textbf{2005} (2005), no.~10, 1--136.

\bibitem[{Len}20]{sym-mon-global}
Tobias {Lenz}, \emph{{On the Global Homotopy Theory of Symmetric Monoidal
  Categories}}, {preprint, available as \texttt{arXiv:2009.07004}}, 2020.

\bibitem[{Len}21]{perm-parsum-categorical}
\bysame, \emph{{Parsummable Categories as a Strictification of Symmetric
  Monoidal Categories}}, {Theory Appl. Categ.} \textbf{37} (2021), no.~17,
  482--529.

\bibitem[{Len}22]{gnpg}
\bysame, \emph{{Genuine vs.~Na\"ive Permutative \(G\)-Categories}}, {preprint,
  available as \url{arXiv:2203.02277}}, 2022.

\bibitem[{Lin}13]{lind}
John~A. {Lind}, \emph{{Diagram Spaces, Diagram Spectra and Spectra of Units}},
  {Algebr. Geom. Topol.} \textbf{13} (2013), no.~4, 1857--1935.

\bibitem[LMS86]{lms}
L.~Gaunce {Lewis, Jr.}, J.~Peter {May}, and Mark {Steinberger},
  \emph{{Equivariant Stable Homotopy Theory. With contributions by J. E.
  McClure}}, {Lect. Notes Math.}, vol. 1213, Springer, Berlin Heidelberg New
  York London Paris Tokyo, 1986.

\bibitem[{Lur}09]{htt}
Jacob {Lurie}, \emph{{Higher Topos Theory}}, Princeton, NJ: Princeton
  University Press, 2009, {updated version available from the author's homepage
  \url{www.math.ias.edu/~lurie/papers/HTT.pdf}}.

\bibitem[{Lur}18]{kerodon}
\bysame, \emph{{Kerodon}}, \url{https://kerodon.net}, 2018.

\bibitem[{Mac}98]{cat-working}
Saunders {Mac Lane}, \emph{{Categories for the Working Mathematician}}, 2nd
  ed., {Grad. Texts Math.}, vol.~5, New York, NY: Springer, 1998.

\bibitem[{Man}10]{mandell}
Michael~A. {Mandell}, \emph{{An Inverse \(K\)-Theory Functor}}, {Doc. Math.}
  \textbf{15} (2010), 765--791.

\bibitem[{May}72]{may-loop-spaces}
J.~Peter {May}, \emph{{The Geometry of Iterated Loop Spaces}}, {Lect. Notes
  Math.}, vol. 271, Springer Berlin Heidelberg, 1972.

\bibitem[{May}74]{may-permutative}
\bysame, \emph{{E\(_\infty\) Spaces, Group Completions, and Permutative
  Categories}}, {New Developments in Topology. The edited and revised
  proceedings of the symposium on algebraic topology, Oxford, June 1972}
  (Graeme~B. {Segal}, ed.), {London Mathematical Society Lecture Note Series},
  vol.~11, Cambridge University Press, 1974, pp.~61--94.

\bibitem[{May}78]{may-unique}
\bysame, \emph{{The Spectra Associated to Permutative Categories}}, {Topology}
  \textbf{17} (1978), 225--228.

\bibitem[{Maz}16]{mazel-gee}
Aaron {Mazel-Gee}, \emph{{Quillen Adjunctions Induce Adjunctions of
  Quasicategories}}, {New York J. Math.} \textbf{22} (2016), 57--93.

\bibitem[{McC}69]{mccord}
Michael~C. {McCord}, \emph{{Classifying Spaces and Infinite Symmetric
  Products}}, {Trans. Am. Math. Soc.} \textbf{146} (1969), 273--298.

\bibitem[{Mer}17]{merling}
Mona {Merling}, \emph{{Equivariant Algebraic \(K\)-Theory of \(G\)-Rings}},
  {Math. Z.} \textbf{285} (2017), no.~3-4, 1205--1248.

\bibitem[{Mil}84]{miller}
Haynes {Miller}, \emph{{The Sullivan Conjecture on Maps From Classifying
  Spaces}}, {Ann. Math. (2)} \textbf{120} (1984), 39--87.

\bibitem[MM02]{mandell-may}
Michael~A. {Mandell} and J.~Peter {May}, \emph{{Equivariant Orthogonal Spectra
  and \(S\)-Modules}}, vol. 159, {Mem. Am. Math. Soc.}, no. 755, Providence,
  RI: American Mathematical Society (AMS), 2002.

\bibitem[MMO17]{may-merling-osorno}
J.~Peter {May}, Mona {Merling}, and Ang{\'e}lica~M. {Osorno},
  \emph{{Equivariant Infinite Loop Space Theory, the Space Level Story}},
  {preprint, available as \texttt{arXiv:1704.03413}}, 2017.

\bibitem[{Ost}16]{equivariant-gamma}
Dominik {Ostermayr}, \emph{{Equivariant \(\Gamma\)-Spaces}}, {Homology Homotopy
  Appl.} \textbf{18} (2016), no.~1, 295--324.

\bibitem[{Qui}67]{quillen}
Daniel~G. {Quillen}, \emph{{Homotopical Algebra}}, Lect. Notes Math., vol.~43,
  Springer, Berlin Heidelberg New York, 1967.

\bibitem[{Rez}02]{rezk-proper}
Charles {Rezk}, \emph{{Every Homotopy Theory of Simplicial Algebras Admits a
  Proper Model}}, {Topology Appl.} \textbf{119} (2002), no.~1, 65--94.

\bibitem[{Rie}14]{cathtpy}
Emily {Riehl}, \emph{{Categorical Homotopy Theory}}, {New Math. Monogr.},
  vol.~24, Cambridge: Cambridge University Press, 2014.

\bibitem[{Sch}99]{schwede-gamma}
Stefan {Schwede}, \emph{{Stable Homotopical Algebra and \(\Gamma\)-Spaces}},
  {Math. Proc. Camb. Philos. Soc.} \textbf{126} (1999), no.~2, 329--356.

\bibitem[{Sch}07]{schwede-symmetric-book}
\bysame, \emph{{An untitled book project about symmetric spectra}}, {available
  at \url{www.math.uni-bonn.de/people/schwede/SymSpec.pdf}}, 2007, {version
  dated July 12th, 2007}.

\bibitem[{Sch}08]{schwede-semistable}
\bysame, \emph{{On the Homotopy Groups of Symmetric Spectra}}, {Geom. Topol.}
  \textbf{12} (2008), no.~3, 1313--1344.

\bibitem[{Sch}17]{schwede-sym-prod}
\bysame, \emph{{Equivariant Properties of Symmetric Products}}, {J. Am. Math.
  Soc.} \textbf{30} (2017), no.~3, 673--711.

\bibitem[{Sch}18]{schwede-book}
\bysame, \emph{{Global Homotopy Theory}}, {New Math. Monogr.}, vol.~34,
  Cambridge: Cambridge University Press, 2018.

\bibitem[{Sch}19]{schwede-cat}
\bysame, \emph{{Categories and Orbispaces}}, {Algebr.~Geom.~Top.} \textbf{19}
  (2019), no.~6, 3171--3215.

\bibitem[{Sch}20a]{schwede-equivariant}
\bysame, \emph{{Lectures on Equivariant Stable Homotopy Theory}}, {lecture
  notes, available at
  \url{www.math.uni-bonn.de/people/schwede/equivariant.pdf}}, 2020, {version
  dated February 12th 2020}.

\bibitem[{Sch}20b]{schwede-orbi}
\bysame, \emph{{Orbispaces, Orthogonal Spaces, and the Universal Compact Lie
  Group}}, {Math.~Z.} \textbf{294} (2020), 71--107.

\bibitem[{Sch}22]{schwede-k-theory}
\bysame, \emph{{Global Algebraic $K$-Theory}}, {J. Topol.} \textbf{15} (2022),
  1325--1454.

\bibitem[{Seg}74]{segal-gamma}
Graeme~B. {Segal}, \emph{{Categories and Cohomology Theories}}, {Topology}
  \textbf{13} (1974), 293--312.

\bibitem[{Shi}89]{shimakawa}
Kazuhisa {Shimakawa}, \emph{{Infinite Loop $G$-Spaces Associated to Monoidal
  $G$-Graded Categories}}, {Publ. Res. Inst. Math. Sci.} \textbf{25} (1989),
  no.~2, 239--262.

\bibitem[{Shi}91]{shimakawa-simplify}
\bysame, \emph{{A Note on $\Gamma_G$-Spaces}}, Osaka J. Math. \textbf{28}
  (1991), no.~2, 223--228.

\bibitem[{Shi}04]{shipley-mixed}
Brooke {Shipley}, \emph{{A Convenient Model Category for Commutative Ring
  Spectra}}, {Homotopy Theory: Relations with Algebraic Geometry, Group
  Cohomology, and Algebraic \(K\)-Theory. Papers from the international
  conference on algebraic topology, Northwestern University, Evanston, IL, USA,
  March 24--28, 2002}, {Contemp. Math.}, vol. 346, Providence, RI: American
  Mathematical Society (AMS), 2004, pp.~473--483.

\bibitem[SS79]{shimada-shimakawa}
Nobuo {Shimada} and Kazuhisa {Shimakawa}, \emph{{Delooping Symmetric Monoidal
  Categories}}, {Hiroshima Math. J.} \textbf{9} (1979), 627--645.

\bibitem[SS00]{schwede-shipley-monoidal}
Stefan {Schwede} and Brooke~E. {Shipley}, \emph{{Algebras and Modules in
  Monoidal Model Categories}}, {Proc. Lond. Math. Soc. (3)} \textbf{80} (2000),
  no.~2, 491--511.

\bibitem[SS12]{sagave-schlichtkrull}
Steffen {Sagave} and Christian {Schlichtkrull}, \emph{{Diagram Spaces and
  Symmetric Spectra}}, {Adv. Math.} \textbf{231} (2012), no.~3-4, 2116--2193.

\bibitem[SS16]{solberg-schlichtkrull}
Christian {Schlichtkrull} and Mirjam {Solberg}, \emph{{Braided Injections and
  Double Loop Spaces}}, {Trans. Am. Math. Soc.} \textbf{368} (2016), no.~10,
  7305--7338.

\bibitem[SS20]{I-vs-M-1-cat}
Steffen {Sagave} and Stefan {Schwede}, \emph{{Homotopy Invariance of
  Convolution Products}}, {Int. Math. Res. Not.} \textbf{2021} (2020), no.~8,
  6246--6292.

\bibitem[{Ste}16]{cellular}
Marc {Stephan}, \emph{{On Equivariant Homotopy Theory for Model Categories}},
  {Homology Homotopy Appl.} \textbf{18} (2016), no.~2, 183--208.

\bibitem[{Sul}05]{sullivan}
Dennis~P. {Sullivan}, \emph{{Geometric Topology. Localization, Periodicity and
  Galois Symmetry. The 1970 MIT Notes. Edited by Andrew Ranicki}}, Dordrecht:
  Springer, 2005.

\bibitem[{Tho}95]{thomason}
Robert~W. {Thomason}, \emph{{Symmetric Monoidal Categories Model All Connective
  Spectra}}, {Theory Appl. Categ.} \textbf{1} (1995), 78--118.

\bibitem[{Whi}17]{white-cmon}
David {White}, \emph{{Model Structures on Commutative Monoids in General Model
  Categories}}, {J. Pure Appl. Algebra} \textbf{221} (2017), no.~12,
  3124--3168.

\end{thebibliography}
\endgroup

\ifidx
\def\pagedeclaration#1{,~\textbf{#1}}
\nomlabelwidth=0pt
\renewcommand*{\nompreamble}{\markboth{\textsc{\MakeUppercase{\nomname}}}{\textsc{\MakeUppercase{\nomname}}}\begin{multicols*}{2}\footnotesize\begingroup\baselineskip=\the\baselineskip plus 1pt\labelsep=1em\raggedright}
\renewcommand*{\nompostamble}{\endgroup
\end{multicols*}}
\printnomenclature
\printindex
\fi
\end{document}